\documentclass[10pt]{amsart}
\usepackage{xypic,amscd,amssymb,latexsym}
\usepackage{amsmath}	
\usepackage{graphicx}
\usepackage{pdfsync}
\usepackage{color}
\usepackage{cancel}
\usepackage{fancyhdr}

\usepackage[pdftex]{hyperref}

\usepackage{soul}

\usepackage{mathrsfs}

\usepackage{url}
\usepackage{pifont}

\DeclareMathAlphabet{\mathpzc}{OT1}{pzc}{m}{it}

\usepackage{setspace}

\usepackage{multirow}

\usepackage{changepage}
\usepackage{xcolor}

\usepackage{tikz}

\begin{document}
\pagenumbering{arabic}

\renewcommand{\circled}[1]{\scalebox{0.8}{\textcircled{\fontsize{8}{9}\selectfont \raisebox{-0.1pt}{#1}}}}

\newtheorem{theorem}{Theorem}[section]
\newtheorem{proposition}[theorem]{Proposition}
\newtheorem{lemma}[theorem]{Lemma}
\newtheorem{corollary}[theorem]{Corollary}
\newtheorem{remark}[theorem]{Remark}
\newtheorem{definition}[theorem]{Definition}
\newtheorem{question}[theorem]{Question}
\newtheorem{claim}[theorem]{Claim}
\newtheorem{conjecture}[theorem]{Conjecture}
\newtheorem{defprop}[theorem]{Definition and Proposition}
\newtheorem{example}[theorem]{Example}
\newtheorem{deflem}[theorem]{Definition and Lemma}

\def\qed{{\quad \vrule height 8pt width 8pt depth 0pt}}

\newcommand{\vs}[0]{\vspace{2mm}}

\newcommand{\til}[1]{\widetilde{#1}}

\newcommand{\mcal}[1]{\mathcal{#1}}

\renewcommand{\ul}[1]{\underline{#1}}

\newcommand{\ol}[1]{\overline{#1}}

\newcommand{\wh}[1]{\widehat{#1}}

\newcommand{\mattwo}[4]{
\left(
\begin{array}{cc}
#1 & #2 \\
#3 & #4 \\
\end{array}
\right)
}

\newcommand{\matthree}[9]{
\left(
\begin{array}{ccc}
#1 & #2 & #3 \\
#4 & #5 & #6 \\
#7 & #8 & #9
\end{array}
\right)
}

\newcommand{\smallmatthree}[9]{
\left(
\begin{smallmatrix}
#1 & #2 & #3 \\
#4 & #5 & #6 \\
#7 & #8 & #9
\end{smallmatrix}
\right)
}

\definecolor{pinky}{rgb}{1.0, 0, 1.0}

\newcommand{\red}[1]{{\color{red}#1}}
\newcommand{\redfix}[1]{{\color{black}#1}}
\newcommand{\blue}[1]{{\color{blue}#1}}
\newcommand{\bluefix}[1]{{\color{black}#1}}
\newcommand{\purple}[1]{{\color{pinky}#1}}
\newcommand{\green}[1]{{\color{green}#1}}


\author[H. Kim]{Hyun Kyu Kim}
\email{hkim@kias.re.kr; hyunkyu87@gmail.com}

\address{School of Mathematics, Korea Institute for Advanced Study (KIAS), 85 Hoegiro Dongdaemun-gu, Seoul 02455, Republic of Korea}

\numberwithin{equation}{section}

\title[${\rm SL}_3$-laminations as bases for ${\rm PGL}_3$ cluster varieties for surfaces]{\resizebox{150mm}{!}{${\rm SL}_3$-laminations as bases for ${\rm PGL}_3$ cluster varieties for surfaces}}

\begin{abstract}
In this paper we \redfix{partially} settle Fock-Goncharov's duality conjecture for cluster varieties associated to their moduli spaces of ${\rm G}$-local systems on a punctured surface $\frak{S}$ with boundary data, when ${\rm G}$ is a group of type $A_2$, namely ${\rm SL}_3$ and ${\rm PGL}_3$. Based on Kuperberg's ${\rm SL}_3$-webs, we introduce the notion of ${\rm SL}_3$-laminations on $\frak{S}$ defined as certain ${\rm SL}_3$-webs with integer weights. We introduce coordinate systems for ${\rm SL}_3$-laminations, and show that ${\rm SL}_3$-laminations satisfying a congruence property are geometric realizations of the tropical integer points of the cluster $\mathscr{A}$-moduli space $\mathscr{A}_{{\rm SL}_3,\frak{S}}$. Per each such ${\rm SL}_3$-lamination, we construct a regular function on the cluster $\mathscr{X}$-moduli space $\mathscr{X}_{{\rm PGL}_3,\frak{S}}$. We show that these functions form a basis of the ring of all regular functions. For a proof, we develop ${\rm SL}_3$ quantum and classical trace maps for any triangulated bordered surface with marked points, and state-sum formulas for them. We construct quantum versions of the basic regular functions on $\mathscr{X}_{{\rm PGL}_3,\frak{S}}$. \redfix{The bases constructed in this paper are built from non-elliptic webs, hence could be viewed as higher `bangles' bases, and the corresponding `bracelets' versions can also be considered as direct analogs of Fock-Goncharov's and Allegretti-Kim's bases for the ${\rm SL}_2$-${\rm PGL}_2$ case.}
\end{abstract}

\maketitle

\vspace{-7mm}

\tableofcontents

\section{Introduction}

\subsection{Background on Fock-Goncharov's duality conjecture}

Let $\frak{S}$ be a \ul{\em punctured surface}, i.e. a compact oriented surface of genus $g\ge 0$ minus $n\ge 1$ punctures. We say $\frak{S}$ is \ul{\em triangulable} if $g=0$ and $n\ge 4$, or if $g\ge 1$ and $n\ge 1$ (so we exclude the three-punctured sphere; see \S\ref{subsec:surface}). Let ${\rm G}$ be a split reductive algebraic group over $\mathbb{Q}$. A {\em ${\rm G}$-local system} $\mathcal{L}$ on $\frak{S}$ can be thought of as a right principal ${\rm G}$-bundle on $\frak{S}$ together with a flat ${\rm G}$-connection. An isomorphism class of $\mathcal{L}$ is captured by the monodromy representation $\pi_1(\frak{S}) \to {\rm G}$ which is a group homomorphism defined up to conjugation in ${\rm G}$. Hence the moduli space $\mathscr{L}_{{\rm G},\frak{S}}$ of ${\rm G}$-local systems on $\frak{S}$ is identified with the \redfix{${\rm G}$-character stack for $\frak{S}$, namely}
$$
\mathscr{L}_{{\rm G},\frak{S}} = {\rm Hom}(\pi_1(\frak{S}),{\rm G})/{\rm G}.
$$
Fock and Goncharov defined \cite{FG06} two related moduli stacks
$$
\mathscr{A}_{{\rm G},\frak{S}} \quad\mbox{and}\quad 
\mathscr{X}_{{\rm G},\frak{S}},
$$
where $\mathscr{A}_{{\rm G},\frak{S}}$ parametrizes the {\em decorated} ${\rm G}$-local systems, while $\mathscr{X}_{{\rm G},\frak{S}}$ parametrizes the {\em framed} ${\rm G}$-local systems. To briefly recall the definitions, consider $\frak{S}$ as being given by a compact oriented surface minus $n$ open discs, so that punctures now become boundary circle components. Choose a Borel subgroup ${\rm B}$ of ${\rm G}$, and let $\mathcal{B} := {\rm G} / {\rm B}$ be the flag variety. For a ${\rm G}$-local system $\mathcal{L}$ on $\frak{S}$, let $\mathcal{L}_{\mathcal{B}} := \mathcal{L} \times_{\rm G} \mathcal{B}$ be the associated flag bundle. A {\em framing} on $\mathcal{L}$ is a choice $\beta$ of a flat section of the restriction $\mathcal{L}_\mathcal{B}|_{\partial \frak{S}}$ of $\mathcal{L}_\mathcal{B}$ to the boundary of $\frak{S}$.
\begin{definition}[\cite{FG06}]
A \ul{\em framed ${\rm G}$-local system} on $\frak{S}$ is a pair $(\mathcal{L},\beta)$ of a ${\rm G}$-local system $\mathcal{L}$ on $\frak{S}$ together with a framing on $\mathcal{L}$. Let $\mathscr{X}_{{\rm G},\frak{S}}$ be the moduli stack parametrizing framed ${\rm G}$-local systems on $\frak{S}$.
\end{definition}
For our case, define a decorated ${\rm G}$-local system and its moduli space $\mathscr{A}_{{\rm G},\frak{S}}$ analogously, with ${\rm B}$ being replaced by the maximal unipotent subgroup ${\rm U} := [{\rm B},{\rm B}]$. For the case when ${\rm G}$ is of type $A_1$, Fock and Goncharov showed \cite{FG06} that $\mathscr{A}_{{\rm SL}_2,\frak{S}}$ and $\mathscr{X}_{{\rm PGL}_2,\frak{S}}$ recover the decorated Teichm\"uller space and the enhanced Teichm\"uller space of the surface $\frak{S}$ respectively\redfix{, as sets of positive real points of the real loci}. For higher rank groups ${\rm G}$, \redfix{the sets of positive real points of} these spaces \redfix{$\mathscr{A}_{{\rm G},\frak{S}}$ and $\mathscr{X}_{{\rm G},\frak{S}}$} can be viewed as providing models for \redfix{Fock-Goncharov's versions of} {\em higher Teichm\"uller spaces}.

\vs

\begin{figure}[htbp!]
\vspace{-4mm}
\begin{center}
\hspace*{-0mm}
\raisebox{-0.3\height}{\scalebox{0.6}{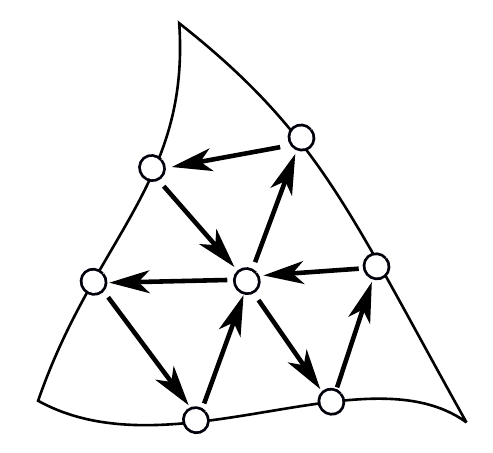}} 
\end{center}
\vspace{-5mm}
\caption{3-triangulation quiver, for one triangle}
\vspace{-2mm}
\label{fig:3-triangulation}
\end{figure}

The present paper concerns the case when ${\rm G}$ is of type $A_2$, or more precisely the spaces $\mathscr{A}_{{\rm SL}_3,\frak{S}}$ and $\mathscr{X}_{{\rm PGL}_3,\frak{S}}$. Pivotal in the study of these spaces are Fock-Goncharov's special coordinate systems \cite{FG06}
\begin{align}
\label{eq:FG_charts}
\mathscr{A}_{{\rm SL}_3,\frak{S}} \dashrightarrow (\mathbb{G}_m)^N \quad\mbox{and}\quad
\mathscr{X}_{{\rm PGL}_3,\frak{S}} \dashrightarrow (\mathbb{G}_m)^N
\end{align}
which are birational maps, associated to each choice of an \ul{\em ideal triangulation} $\Delta$ of $\frak{S}$, which is a maximal collection of mutually disjoint simple arcs in $\frak{S}$ running between punctures of $\frak{S}$, where $\frak{S}$ is viewed as a punctured surface again, dividing $\frak{S}$ into ideal triangles; we assume the valence of $\Delta$ at each puncture is at least $3$ (see Rem.\ref{rem:non-regular_triangulations} for a comment on this condition). A remarkable fact is that for two ideal triangulations, the coordinate change maps are {\em positive} rational, not involving any subtraction, and moreover, they follow patterns called the {\em cluster mutations} appearing in the theory of cluster algebras and cluster varieties \cite{FZ07} \cite{FG06}. To elaborate, for an ideal triangulation $\Delta$, define the \ul{\em 3-triangulation} $Q_\Delta$ of $\Delta$ as the quiver obtained by gluing the quivers associated to triangles of $\Delta$, as depicted in Fig.\ref{fig:3-triangulation}. So $Q_\Delta$ will have two nodes lying on each arc of $\Delta$ and one node lying in the interior of each triangle of $\Delta$, while the arrows are as in Fig.\ref{fig:3-triangulation} for each triangle. Let
$$
\mathcal{V}(Q_\Delta) = \mbox{the set of all nodes of the 3-triangulation quiver $Q_\Delta$ of the triangulation $\Delta$}.
$$
The Fock-Goncharov coordinates (eq.\eqref{eq:FG_charts}), for each of $\mathscr{A}_{{\rm SL}_3,\frak{S}}$ and $\mathscr{X}_{{\rm PGL}_3,\frak{S}}$, are enumerated by $\mathcal{V}(Q_\Delta)$.

In the theory of cluster algebras and cluster varieties, there is a certain combinatorial process called the \ul{\em quiver mutation} at a node of a quiver, transforming a quiver into another one according to some rule (Def.\ref{def:quiver_mutation}). When ideal triangulations $\Delta$ and $\Delta'$ differ only by one arc $e$, we say that they are related by a {\em flip}, and it is known that the quiver $Q_{\Delta'}$ can be obtained from $Q_\Delta$ by a sequence of four quiver mutations, first at the two nodes lying in the arc $e$ of $\Delta$, then at the two nodes lying in the interiors of the two triangles of $\Delta$ having $e$ as a side (Lem.\ref{lem:flip_as_four_quiver_mutations}). Per each quiver mutation, the \ul{\em cluster $\mathscr{A}$-mutation} and the \ul{\em cluster $\mathscr{X}$-mutation} are associated, which are certain coordinate change formulas for the coordinate functions associated to the nodes (Def.\ref{def:A}, \ref{def:X}). Fock and Goncharov \cite{FG06} showed that their coordinates on $\mathscr{A}_{{\rm SL}_3,\frak{S}}$ and those on $\mathscr{X}_{{\rm PGL}_3,\frak{S}}$ associated to $\Delta$ and $\Delta'$ indeed transform according to the corresponding compositions of the four cluster mutations (Prop.\ref{prop:cluster_atlases_of_FG}). One can summarize this result as having constructed birational equivalences
$$
\mathscr{A}_{{\rm SL}_3,\frak{S}} \dashrightarrow \mathscr{A}_{|Q_{\Delta}|} \quad\mbox{and}\quad
\mathscr{X}_{{\rm PGL}_3,\frak{S}} \dashrightarrow \mathscr{X}_{|Q_{\Delta}|}
$$
from the moduli spaces $\mathscr{A}_{{\rm SL}_3,\frak{S}}$ and $\mathscr{X}_{{\rm PGL}_3,\frak{S}}$ to the {\em cluster varieties} $\mathscr{A}_{|Q_{\Delta}|}$ and $\mathscr{X}_{|Q_{\Delta}|}$, which are abstract schemes constructed by gluing split algebraic tori $(\mathbb{G}_m)^N = (\mathbb{G}_m)^{\mathcal{V}(Q_\Delta)}$ along the cluster mutation maps, where $|Q_\Delta|$ denotes the quiver mutation equivalence class of $Q_\Delta$. One of the important conjectures set out by Fock and Goncharov in \cite{FG06} is the following, which they proposed as a generalization of their ${\rm SL}_2$-${\rm PGL}_2$ result. We first define the necessary notions.
\begin{definition}
The cluster coordinate charts of $\mathscr{A}_{{\rm SL}_3,\frak{S}}$ in eq.\eqref{eq:FG_charts} associated to ideal triangulations $\Delta$ of $\frak{S}$ are related by positive rational maps, so it makes sense to consider the set $\mathscr{A}_{{\rm SL}_3,\frak{S}}(\mathbb{P})$ of points valued in a {\em semi-field} $\mathbb{P}$, which is a set with an addition and a multiplication, where the multiplication makes $\mathbb{P}$ an abelian group, while the addition is merely associative, commutative, and distributive for the multiplication. In particular, this applies for the \ul{\em semi-field of tropical integers} $\mathbb{Z}^{\redfix{T}}$, which is $\mathbb{Z}$ as a set, with the tropical addition of $a,b$ being $\max(a,b)$ and the tropical multiplication of $a,b$ being $a+b$.
\end{definition}
It is known that, as a set, $\mathscr{A}_{{\rm SL}_3,\frak{S}}(\mathbb{P})$ is in bijection with $\mathbb{P}^N$ (or $\mathbb{P}^{\mathcal{V}(Q_\Delta)}$).
\begin{definition}
\label{def:intro_L}
Let ${\bf L}(\mathscr{X}_{{\rm PGL}_3,\frak{S}})$ be the ring of all rational functions on $\mathscr{X}_{{\rm PGL}_3,\frak{S}}$ that are regular with respect to all the cluster coordinate charts of $\mathscr{X}_{{\rm PGL}_3,\frak{S}}$ in eq.\eqref{eq:FG_charts} associated to ideal triangulations $\Delta$ of $\frak{S}$, i.e. the ring of all {\em universally Laurent polynomial} functions for {\em all ideal triangulations}.
\end{definition}
\begin{conjecture}[Fock-Goncharov's duality conjecture for ${\rm SL}_3$-${\rm PGL}_3$; \cite{FG06}]
\label{conj:duality}
There exists a canonical map
$$
\mathbb{I} : \mathscr{A}_{{\rm SL}_3,\frak{S}}(\mathbb{Z}^{\redfix{T}}) \to {\bf L}(\mathscr{X}_{{\rm PGL}_3,\frak{S}})
$$
satisfying favorable properties\redfix{, such as: $\mathbb{I}$ is injective and its image forms a basis of ${\bf L}(\mathscr{X}_{{\rm PGL}_3,\frak{S}})$; the structure constants for this basis are positive integers; and for each $\ell \in \mathscr{A}_{{\rm SL}_3,\frak{S}}(\mathbb{Z}^{\redfix{T}})$, the function $\mathbb{I}(\ell)$ is a Laurent polynomial for each ideal triangulation with positive integer coefficients}.
\end{conjecture}
A kind of implicit prerequisite conjecture is:
\begin{conjecture}[Conjectural geometric model of tropical integer points of $\mathscr{A}_{{\rm SL}_3,\frak{S}}$; \cite{FG06}]
\label{conj:geometric_model_of_tropical_integer_points}
There is a natural geometric model for the set $\mathscr{A}_{{\rm SL}_3,\frak{S}}(\mathbb{Z}^{\redfix{T}})$.
\end{conjecture}
There have been attempts for Conjecture \ref{conj:geometric_model_of_tropical_integer_points}, e.g. by Ian Le \cite{Ian} (`higher' laminations) and by Goncharov and Shen \cite{GS15} (top-dimensional components of `surface affine Grassmannian' stack, for `positive' integral tropical points), but these are not as direct and intuitive as Fock-Goncharov's answer \cite{FG06} for the ${\rm SL}_2$ case $\mathscr{A}_{{\rm SL}_2,\frak{S}}(\mathbb{Z}^{\redfix{T}})$ by the {\em integral laminations} on $\frak{S}$, which are certain collections of simple closed curves in $\frak{S}$ with integer weights. More importantly, a good answer must also immediately help answering Conjecture \ref{conj:duality}, but Conjecture \ref{conj:duality} has been elusive.

\vs

In the present paper we settle the above two conjectures to a large extent. In particular, we provide a geometrically intuitive model for $\mathscr{A}_{{\rm SL}_3,\frak{S}}(\mathbb{Z}^{\redfix{T}})$, and for each $\ell \in \mathscr{A}_{{\rm SL}_3,\frak{S}}(\mathbb{Z}^{\redfix{T}})$ we explicitly construct a universally Laurent function on $\mathscr{X}_{{\rm PGL}_3,\frak{S}}$, and prove some important properties. In fact, ${\bf L}(\mathscr{X}_{{\rm PGL}_3,\frak{S}})$ is replaced by $\mathscr{O}(\mathscr{X}_{{\rm PGL}_3,\frak{S}})$ and 
$\mathscr{O}_{\rm cl}(\mathscr{X}_{{\rm PGL}_3,\frak{S}})$ which we believe are more correct target spaces, where
$$
\mathscr{O}_{\rm cl}(\mathscr{X}_{{\rm PGL}_3,\frak{S}}) := \mbox{the ring of all functions on $\mathscr{X}_{{\rm PGL}_3,\frak{S}}$ that are regular for {\em all} cluster $\mathscr{X}$-charts},
$$
which a priori may be different from ${\bf L}(\mathscr{X}_{{\rm PGL}_3,\frak{S}})$ as there are cluster $\mathscr{X}$-charts that are not associated to ideal triangulations $\Delta$ (see \S\ref{subsec:cluster_atlases}); this better suits the theory of cluster varieties, e.g. we have $\mathscr{O}_{\rm cl}(\mathscr{X}_{{\rm PGL}_3,\frak{S}}) \cong \mathscr{O}(\mathscr{X}_{|Q_\Delta|})$.

\vs

Meanwhile, there is a more general version of Conjecture \ref{conj:duality}, for cluster varieties for any quivers $Q$ (or more generally for any skew-symmetrizable integer matrices $\varepsilon$), which was solved in the celebrated paper by Gross, Hacking, Keel and Kontsevich \cite{GHKK}. They showed that if a quiver $Q$ satisfies a certain combinatorial condition, then there exists a canonical map
$$
\mathbb{I} : \mathscr{A}_{|Q|}(\mathbb{Z}^t) \to \mathscr{O}(\mathscr{X}_{|Q|})
$$
for the cluster varieties $\mathscr{A}_{|Q|}$ and $\mathscr{X}_{|Q|}$, satisfying favorable properties\redfix{, where $\mathbb{Z}^t$ is the version of the tropical integer semi-field whose tropical addition is defined as $\min(\cdot,\cdot)$ (we note that $\mathbb{Z}^T$ and $\mathbb{Z}^t$ are isomorphic)}. It was shown by Goncharov and Shen \cite{GS18} that the 3-triangulation quiver $Q = Q_\Delta$ satisfies the condition of \cite{GHKK} if $\frak{S}$ has at least two punctures, hence proving the existence of an ${\rm SL}_3$-${\rm PGL}_3$ duality map as being sought for in Conjecture \ref{conj:duality}. However, the construction in \cite{GHKK} is quite combinatorial and uniform for all quivers $Q$, not giving special geometric meaning for the quivers coming from surfaces. Actual computations of their functions $\mathbb{I}(\ell) \in \mathscr{O}(\mathscr{X}_{|Q|})$ require enormous amount of combinatorics in large dimensional Euclidean spaces, and to find a geometric meaning of the resulting (universally) Laurent polynomials is a big important challenge. \redfix{As of now, a direct computation of the functions $\mathbb{I}(\ell)$ obtained in \cite{GHKK} is in fact almost impossible, not just being difficult, even for the simplest surfaces like once-punctured torus, because a crucial ingredient called the `consistent scattering diagram' is only known to exist but has never been constructed in a manner that is explicit enough for the purpose of this computation.} 

\vs

On the other hand, our \redfix{${\rm SL}_3$-${\rm PGL}_3$} duality map in the present paper is down to earth\redfix{, explicitly constructive, and geometric}. Of course a very natural conjecture would then be to compare our map with Gross-Hacking-Keel-Kontsevich's, but we expect that it will be quite difficult. \redfix{We notice that even for the ${\rm SL}_2$-${\rm PGL}_2$ case,  i.e. when $Q$ is the {\em 2-triangulation} \cite{FG06} of an ideal triangulation $\Delta$ of $\frak{S}$ (which is a quiver whose set of nodes is in bijection with $\Delta$), Gross-Hacking-Keel-Kontsevich's duality map $\mathbb{I}$ \cite{GHKK} has not been computed nor proved to match Fock-Goncharov's constructive duality map $\mathbb{I}$ \cite{FG06}; only recently, Mandel and Qin announced that they proved the equality of these two $\mathbb{I}$'s in their upcoming paper \cite{MQ}.} Besides, our \redfix{${\rm SL}_3$-${\rm PGL}_3$} duality map do\redfix{es} not exclude \redfix{punctured} surfaces $\frak{S}$ with only one puncture\redfix{, unlike the one coming from \cite{GHKK} \cite{GS18}}.

\subsection{${\rm SL}_3$-webs and ${\rm SL}_3$-laminations}

One of the major objects to tackle is $\mathscr{O}(\mathscr{X}_{{\rm PGL}_3,\frak{S}})$, the ring of regular functions on the moduli stack $\mathscr{X}_{{\rm PGL}_3,\frak{S}}$. We will see step by step in the present paper how this is closely related to $\mathscr{O}(\mathscr{L}_{{\rm SL}_3,\frak{S}})$, which has been studied in relation to higher rank versions of the so-called {\em skein algebras} of the surface $\frak{S}$. We first briefly recall the ${\rm SL}_2$-${\rm PGL}_2$ story. For each {\em unoriented} closed curve $\gamma$ in $\frak{S}$, there is a natural regular function $f_\gamma$ on $\mathscr{L}_{{\rm SL}_2,\frak{S}}$ given by the \ul{\em trace of monodromy} along $\gamma$, namely, whose value at the point of $\mathscr{L}_{{\rm SL}_2,\frak{S}}$ represented by a monodromy homomorphism $\rho : \pi_1(\frak{S}) \to {\rm SL}_2$ is defined as 
\begin{align}
\label{eq:intro_f_gamma}
f_\gamma([\rho]) = {\rm tr}(\rho(\gamma)),
\end{align}
where in the right hand side $\gamma$ is given an arbitrary orientation. Then $f_\gamma$ is well-defined because the trace is invariant under conjugation, and under taking the inverse in ${\rm SL}_2$. Due to the matrix identity $({\rm tr} A)({\rm tr}B) = {\rm tr}(AB) + {\rm tr}(AB^{-1})$ in ${\rm SL}_2$, these trace-of-monodromy functions $f_\gamma$ satisfy the relation

\vspace{-5mm}

$$
f_\gamma f_{\gamma'} = f_{\gamma_1} + f_{\gamma_2}, \quad \mbox{when $\gamma,\gamma',\gamma_1,\gamma_2$ look like \raisebox{-0.5\height}{
\begingroup%
  \makeatletter%
  \providecommand\color[2][]{%
    \errmessage{(Inkscape) Color is used for the text in Inkscape, but the package 'color.sty' is not loaded}%
    \renewcommand\color[2][]{}%
  }%
  \providecommand\transparent[1]{%
    \errmessage{(Inkscape) Transparency is used (non-zero) for the text in Inkscape, but the package 'transparent.sty' is not loaded}%
    \renewcommand\transparent[1]{}%
  }%
  \providecommand\rotatebox[2]{#2}%
  \newcommand*\fsize{\dimexpr\f@size pt\relax}%
  \newcommand*\lineheight[1]{\fontsize{\fsize}{#1\fsize}\selectfont}%
  \ifx\svgwidth\undefined%
    \setlength{\unitlength}{107.71653543bp}%
    \ifx\svgscale\undefined%
      \relax%
    \else%
      \setlength{\unitlength}{\unitlength * \real{\svgscale}}%
    \fi%
  \else%
    \setlength{\unitlength}{\svgwidth}%
  \fi%
  \global\let\svgwidth\undefined%
  \global\let\svgscale\undefined%
  \makeatother%
  \begin{picture}(1,0.39473684)%
    \lineheight{1}%
    \setlength\tabcolsep{0pt}%
    \put(0,0){\includegraphics[width=\unitlength,page=1]{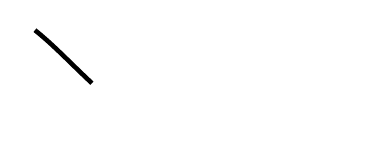}}%
    \put(0.03024523,0.10823152){\color[rgb]{0,0,0}\makebox(0,0)[lt]{\lineheight{1.25}\smash{\begin{tabular}[t]{l}$\gamma$\end{tabular}}}}%
    \put(0,0){\includegraphics[width=\unitlength,page=2]{intro_skein.pdf}}%
    \put(0.25059693,0.11081698){\color[rgb]{0,0,0}\makebox(0,0)[lt]{\lineheight{1.25}\smash{\begin{tabular}[t]{l}$\gamma'$\end{tabular}}}}%
    \put(0,0){\includegraphics[width=\unitlength,page=3]{intro_skein.pdf}}%
    \put(0.50240355,0.10112851){\color[rgb]{0,0,0}\makebox(0,0)[lt]{\lineheight{1.25}\smash{\begin{tabular}[t]{l}$\gamma_1$\end{tabular}}}}%
    \put(0,0){\includegraphics[width=\unitlength,page=4]{intro_skein.pdf}}%
    \put(0.83661403,0.10112851){\color[rgb]{0,0,0}\makebox(0,0)[lt]{\lineheight{1.25}\smash{\begin{tabular}[t]{l}$\gamma_2$\end{tabular}}}}%
    \put(0,0){\includegraphics[width=\unitlength,page=5]{intro_skein.pdf}}%
  \end{picture}%
\endgroup%
} in a small disc.}
$$

\vspace{-4mm}

Since a {\em commutative Kauffman bracket skein algebra} is a free associative algebra generated by closed curves in $\frak{S}$ up to isotopies mod out by the {\em skein relations} which model the above relation, one obtains an algebra homomorphism from the skein algebra to $\mathscr{O}(\mathscr{L}_{{\rm SL}_2,\frak{S}})$, which had been known to be an isomorphism \cite{P76} \cite{S01}. By using the skein relations repeatedly, any element $\mathscr{O}(\mathscr{L}_{{\rm SL}_2,\frak{S}})$ can be expressed as a linear combination of products $f_{\gamma_1} \cdots f_{\gamma_m}$, where $\gamma_1,\ldots,\gamma_m$ are mutually disjoint simple loops, hence forming a multicurve or an example of a {\em lamination}. For our ${\rm SL}_3$-${\rm PGL}_3$ situation, let $\gamma$ be an {\em oriented} closed curve in $\frak{S}$. Still, the trace\redfix{-}of\redfix{-}monodromy function $f_\gamma$ on $\mathscr{L}_{{\rm SL}_3,\frak{S}}$ is defined by the formula eq.\eqref{eq:intro_f_gamma}, and it is known \cite{P76} \cite{S01} that they generate $\mathscr{O}(\mathscr{L}_{{\rm SL}_3,\frak{S}})$, but the algebraic relations among them are different from the case of ${\rm SL}_2$. In particular, oriented multicurves do not yield a basis $\mathscr{O}(\mathscr{L}_{{\rm SL}_3,\frak{S}})$, and one needs to consider certain tri-valent oriented graphs on the surface, called the {\em webs}, first studied by Kuperberg \cite{Kuperberg} in the context of invariant theory, which is kind of directly relevant to our situation which can be viewed as a surface version of invariant theory. Webs for groups of other Dynkin types can be considered ($A_1$ type yielding the unoriented curves), and the corresponding (higher) versions of skein algebras have been studied, notably by Sikora and collaborators \cite{S01} \cite{S05} \cite{SW}. For our purposes, we take the {\em ${\rm SL}_3$-webs} (which are of type $A_2$) and the corresponding {\em ${\rm SL}_3$-skein algebras}, which are extensively investigated recently in Frohman-Sikora \cite{FS} and Higgins \cite{Higgins}. First, we only use the commutative versions, which we somewhat simplified for this section.

\begin{figure}[htbp!]

\vspace{-5mm}

\begin{center}
\hspace*{-5mm}\begin{tabular}{|c|c|c|c|}
\hline
\raisebox{-0.4\height}{
\begingroup%
  \makeatletter%
  \providecommand\color[2][]{%
    \errmessage{(Inkscape) Color is used for the text in Inkscape, but the package 'color.sty' is not loaded}%
    \renewcommand\color[2][]{}%
  }%
  \providecommand\transparent[1]{%
    \errmessage{(Inkscape) Transparency is used (non-zero) for the text in Inkscape, but the package 'transparent.sty' is not loaded}%
    \renewcommand\transparent[1]{}%
  }%
  \providecommand\rotatebox[2]{#2}%
  \newcommand*\fsize{\dimexpr\f@size pt\relax}%
  \newcommand*\lineheight[1]{\fontsize{\fsize}{#1\fsize}\selectfont}%
  \ifx\svgwidth\undefined%
    \setlength{\unitlength}{79.37007874bp}%
    \ifx\svgscale\undefined%
      \relax%
    \else%
      \setlength{\unitlength}{\unitlength * \real{\svgscale}}%
    \fi%
  \else%
    \setlength{\unitlength}{\svgwidth}%
  \fi%
  \global\let\svgwidth\undefined%
  \global\let\svgscale\undefined%
  \makeatother%
  \begin{picture}(1,0.32142857)%
    \lineheight{1}%
    \setlength\tabcolsep{0pt}%
    \put(0.27117991,0.13588895){\color[rgb]{0,0,0}\makebox(0,0)[lt]{\lineheight{1.25}\smash{\begin{tabular}[t]{l}$=3{\O}=$\end{tabular}}}}%
    \put(0,0){\includegraphics[width=\unitlength,page=1]{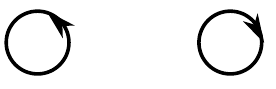}}%
  \end{picture}%
\endgroup%
} & \raisebox{-0.4\height}{
\begingroup%
  \makeatletter%
  \providecommand\color[2][]{%
    \errmessage{(Inkscape) Color is used for the text in Inkscape, but the package 'color.sty' is not loaded}%
    \renewcommand\color[2][]{}%
  }%
  \providecommand\transparent[1]{%
    \errmessage{(Inkscape) Transparency is used (non-zero) for the text in Inkscape, but the package 'transparent.sty' is not loaded}%
    \renewcommand\transparent[1]{}%
  }%
  \providecommand\rotatebox[2]{#2}%
  \newcommand*\fsize{\dimexpr\f@size pt\relax}%
  \newcommand*\lineheight[1]{\fontsize{\fsize}{#1\fsize}\selectfont}%
  \ifx\svgwidth\undefined%
    \setlength{\unitlength}{102.04724409bp}%
    \ifx\svgscale\undefined%
      \relax%
    \else%
      \setlength{\unitlength}{\unitlength * \real{\svgscale}}%
    \fi%
  \else%
    \setlength{\unitlength}{\svgwidth}%
  \fi%
  \global\let\svgwidth\undefined%
  \global\let\svgscale\undefined%
  \makeatother%
  \begin{picture}(1,0.25)%
    \lineheight{1}%
    \setlength\tabcolsep{0pt}%
    \put(0.46080193,0.09917433){\color[rgb]{0,0,0}\makebox(0,0)[lt]{\lineheight{1.25}\smash{\begin{tabular}[t]{l}$=-2$\end{tabular}}}}%
    \put(0,0){\includegraphics[width=\unitlength,page=1]{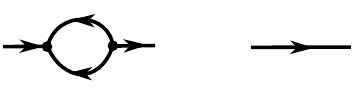}}%
  \end{picture}%
\endgroup%
} & \raisebox{-0.4\height}{
\begingroup%
  \makeatletter%
  \providecommand\color[2][]{%
    \errmessage{(Inkscape) Color is used for the text in Inkscape, but the package 'color.sty' is not loaded}%
    \renewcommand\color[2][]{}%
  }%
  \providecommand\transparent[1]{%
    \errmessage{(Inkscape) Transparency is used (non-zero) for the text in Inkscape, but the package 'transparent.sty' is not loaded}%
    \renewcommand\transparent[1]{}%
  }%
  \providecommand\rotatebox[2]{#2}%
  \newcommand*\fsize{\dimexpr\f@size pt\relax}%
  \newcommand*\lineheight[1]{\fontsize{\fsize}{#1\fsize}\selectfont}%
  \ifx\svgwidth\undefined%
    \setlength{\unitlength}{130.39370079bp}%
    \ifx\svgscale\undefined%
      \relax%
    \else%
      \setlength{\unitlength}{\unitlength * \real{\svgscale}}%
    \fi%
  \else%
    \setlength{\unitlength}{\svgwidth}%
  \fi%
  \global\let\svgwidth\undefined%
  \global\let\svgscale\undefined%
  \makeatother%
  \begin{picture}(1,0.36956522)%
    \lineheight{1}%
    \setlength\tabcolsep{0pt}%
    \put(0.3100311,0.18879898){\color[rgb]{0,0,0}\makebox(0,0)[lt]{\lineheight{1.25}\smash{\begin{tabular}[t]{l}$=$\end{tabular}}}}%
    \put(0,0){\includegraphics[width=\unitlength,page=1]{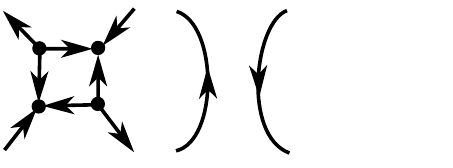}}%
    \put(0.66664339,0.18879898){\color[rgb]{0,0,0}\makebox(0,0)[lt]{\lineheight{1.25}\smash{\begin{tabular}[t]{l}$+$\end{tabular}}}}%
    \put(0,0){\includegraphics[width=\unitlength,page=2]{skein_rel3.pdf}}%
  \end{picture}%
\endgroup%
}
 & \raisebox{-0.4\height}{
\begingroup%
  \makeatletter%
  \providecommand\color[2][]{%
    \errmessage{(Inkscape) Color is used for the text in Inkscape, but the package 'color.sty' is not loaded}%
    \renewcommand\color[2][]{}%
  }%
  \providecommand\transparent[1]{%
    \errmessage{(Inkscape) Transparency is used (non-zero) for the text in Inkscape, but the package 'transparent.sty' is not loaded}%
    \renewcommand\transparent[1]{}%
  }%
  \providecommand\rotatebox[2]{#2}%
  \newcommand*\fsize{\dimexpr\f@size pt\relax}%
  \newcommand*\lineheight[1]{\fontsize{\fsize}{#1\fsize}\selectfont}%
  \ifx\svgwidth\undefined%
    \setlength{\unitlength}{124.72440945bp}%
    \ifx\svgscale\undefined%
      \relax%
    \else%
      \setlength{\unitlength}{\unitlength * \real{\svgscale}}%
    \fi%
  \else%
    \setlength{\unitlength}{\svgwidth}%
  \fi%
  \global\let\svgwidth\undefined%
  \global\let\svgscale\undefined%
  \makeatother%
  \begin{picture}(1,0.34090909)%
    \lineheight{1}%
    \setlength\tabcolsep{0pt}%
    \put(0.28804389,0.1519262){\color[rgb]{0,0,0}\makebox(0,0)[lt]{\lineheight{1.25}\smash{\begin{tabular}[t]{l}$=$\end{tabular}}}}%
    \put(0,0){\includegraphics[width=\unitlength,page=1]{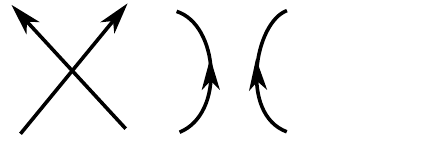}}%
    \put(0.69694537,0.1519262){\color[rgb]{0,0,0}\makebox(0,0)[lt]{\lineheight{1.25}\smash{\begin{tabular}[t]{l}$+$\end{tabular}}}}%
    \put(0,0){\includegraphics[width=\unitlength,page=2]{skein_rel4.pdf}}%
  \end{picture}%
\endgroup%
} \\ \hline
{\rm (S1)} & {\rm (S2)} & {\rm (S3)} & {\rm (S4)} \\ \hline
\end{tabular}
\end{center}
\vspace{-2mm}
\caption{${\rm SL}_3$-skein relations, drawn locally (${\O}$ means empty); the regions bounded by a loop, a $2$-gon, or a $4$-gon in (S1), (S2), (S3) are contractible.}
\vspace{-2mm}
\label{fig:A2-skein_relations}
\end{figure}

\begin{definition}[\cite{Kuperberg} \cite{S01} \cite{SW} \cite{FS}; Def.\ref{def:A2-web}, \ref{def:A2-skein_algebra}]
An \ul{\em ${\rm SL}_3$-web} $W$ in a punctured surface $\frak{S}$ is a union of any finite number of isotopy classes of oriented loops in $\frak{S}$ and/or oriented tri-valent connected graphs such that each tri-valent vertex is either a sink or a source. Let $\mathcal{R}$ be a ring with $1$.

The \ul{\em (commutative) ${\rm SL}_3$-skein algebra} $\mathcal{S}(\frak{S}; \mathcal{R})$ is a free $\mathcal{R}$-module generated by ${\rm SL}_3$-webs mod out by the \ul{\em ${\rm SL}_3$-skein relations} in Fig.\ref{fig:A2-skein_relations}, where the product of elements of $\mathcal{S}(\frak{S}; \mathcal{R})$ is given by the union.
\end{definition}

It is known \cite{S01} that there is a ring isomorphism
\begin{align}
\label{eq:intro_Phi}
\Phi : \mathcal{S}(\frak{S};\mathbb{Q}) \to \mathscr{O}(\mathscr{L}_{{\rm SL}_3,\frak{S}})
\end{align}
sending each oriented loop $\gamma$ to the trace-of-monodromy function $f_\gamma$. So a basis of $\mathscr{O}(\mathscr{L}_{{\rm SL}_3,\frak{S}})$ can be obtained from a basis of $\mathcal{S}(\frak{S};\mathbb{Q})$; there is a nice basis consisting of the so-called {\em non-elliptic} ${\rm SL}_3$-webs.
\begin{definition}[\cite{Kuperberg} \cite{SW}; Def.\ref{def:non-elliptic}]
An ${\rm SL}_3$-web in a punctured surface $\frak{S}$ is \ul{\em non-elliptic} if it has no self-intersection other than (possibly) the tri-valent vertices and does not bound a contractible region bounded by a loop or by two or four oriented edges like in (S1), (S2), (S3). 
\end{definition}
The first main definition of the present paper is the following.
\begin{definition}[${\rm SL}_3$-laminations; Def.\ref{def:A2-lamination}]
An \ul{\em ${\rm SL}_3$-lamination} $\ell$ in a punctured surface $\frak{S}$ is a non-elliptic ${\rm SL}_3$-web $W(\ell)$ in $\frak{S}$ together with integer weights on the (connected) components of $W(\ell)$, subject to the following conditions and equivalence relations:
\begin{enumerate}
\itemsep0em
\item[\rm (1)] the weight of a component containing a tri-valent vertex is $1$;

\item[\rm (2)]  the weight of a component is non-negative unless the component is a \ul{\em peripheral loop}, i.e. a small loop surrounding a puncture of $\frak{S}$;

\item[\rm (3)] an ${\rm SL}_3$-lamination with one of the components having weight $0$ is      equivalent to the ${\rm SL}_3$-lamination with this component removed;

\item[\rm (4)] an ${\rm SL}_3$-lamination with two of the components being isotopic and with weights $a,b$ is equivalent to the ${\rm SL}_3$-lamination with one of these two components removed and the weight $a+b$ given on the other.
\end{enumerate}
Let $\mathscr{A}_{\rm L}(\frak{S};\mathbb{Z})$ be the set of all ${\rm SL}_3$-laminations. 
\end{definition}

\begin{figure}[htbp!]
\vspace{-5mm}
\raisebox{-0.5\height}{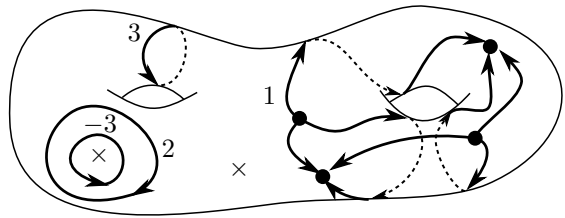}
\vspace{-2mm}
\caption{Example of an ${\rm SL}_3$-lamination ($\times$ are punctures)}
\vspace{-3mm}
\label{fig:example_of_an_A2-lamination}
\end{figure}
We note that in the first version of the present paper \cite{Kim}, ${\rm SL}_3$-webs and ${\rm SL}_3$-laminations were called $A_2$-webs and $A_2$-laminations; we believe ${\rm SL}_3$ is a better label for these objects (see e.g. \redfix{\S\ref{sec:conjectures}} \redfix{of the second and third versions of the present paper}).

\vs

The negative weights for peripheral loops will be used to compensate the difference between $\mathscr{O}(\mathscr{X})$ and $\mathscr{O}(\mathscr{L})$, like in Fock-Goncharov's ${\rm SL}_2$-${\rm PGL}_2$ duality. Crucial in the study of ${\rm SL}_3$-laminations is a special coordinate system for them, associated to an ideal triangulation $\Delta$ of $\frak{S}$. We make use of Frohman-Sikora's coordinates \cite{FS} and Douglas-Sun's coordinates \cite{DS1} for non-elliptic ${\rm SL}_3$-webs. Here we only recall Douglas-Sun's, multiplied by $\frac{1}{3}$. Let $\Delta$ be an ideal triangulation of $\frak{S}$, and let $\wh{\Delta}$ be a \ul{\em split ideal triangulation} of $\frak{S}$ for $\Delta$, obtained from $\Delta$ by adding one more arc $e'$ in $\frak{S}$ per each arc $e$ of $\Delta$, where $e'$ is isotopic to \redfix{$e$ rel endpoints} but disjoint \redfix{from} $e$ and such that $\Delta'$ is a mutually disjoint collection \cite{BW}. Then $\wh{\Delta}$ divides $\frak{S}$ into ideal triangles and ideal {\em biangles}, where a biangle of $\wh{\Delta}$ is associated to each arc of $\Delta$, and a triangle $\wh{t}$ of $\wh{\Delta}$ to each triangle $t$ of $\Delta$.

\begin{figure}[htbp!]
\begin{center}
{\scalebox{0.8}{
\begingroup%
  \makeatletter%
  \providecommand\color[2][]{%
    \errmessage{(Inkscape) Color is used for the text in Inkscape, but the package 'color.sty' is not loaded}%
    \renewcommand\color[2][]{}%
  }%
  \providecommand\transparent[1]{%
    \errmessage{(Inkscape) Transparency is used (non-zero) for the text in Inkscape, but the package 'transparent.sty' is not loaded}%
    \renewcommand\transparent[1]{}%
  }%
  \providecommand\rotatebox[2]{#2}%
  \newcommand*\fsize{\dimexpr\f@size pt\relax}%
  \newcommand*\lineheight[1]{\fontsize{\fsize}{#1\fsize}\selectfont}%
  \ifx\svgwidth\undefined%
    \setlength{\unitlength}{496.06299213bp}%
    \ifx\svgscale\undefined%
      \relax%
    \else%
      \setlength{\unitlength}{\unitlength * \real{\svgscale}}%
    \fi%
  \else%
    \setlength{\unitlength}{\svgwidth}%
  \fi%
  \global\let\svgwidth\undefined%
  \global\let\svgscale\undefined%
  \makeatother%
  \begin{picture}(1,0.22285714)%
    \lineheight{1}%
    \setlength\tabcolsep{0pt}%
    \put(0,0){\includegraphics[width=\unitlength,page=1]{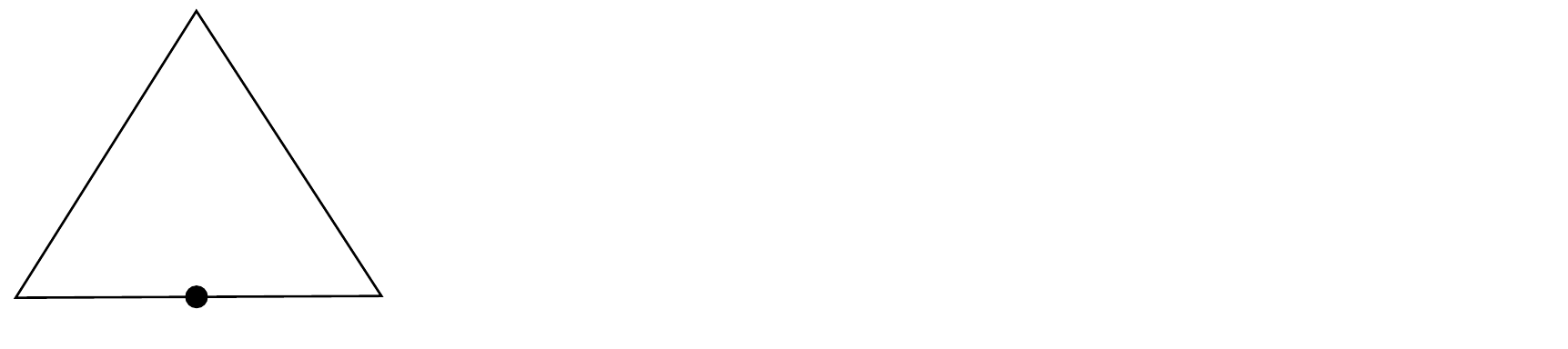}}%
    \put(0.14787992,0.01165401){\color[rgb]{0,0,0}\makebox(0,0)[lt]{\lineheight{1.25}\smash{\begin{tabular}[t]{l}$H_1$\end{tabular}}}}%
    \put(0,0){\includegraphics[width=\unitlength,page=2]{pyramid1.pdf}}%
    \put(0.35954647,0.01165401){\color[rgb]{0,0,0}\makebox(0,0)[lt]{\lineheight{1.25}\smash{\begin{tabular}[t]{l}$H_{-2}$\end{tabular}}}}%
    \put(0,0){\includegraphics[width=\unitlength,page=3]{pyramid1.pdf}}%
    \put(0.60447507,0.01165401){\color[rgb]{0,0,0}\makebox(0,0)[lt]{\lineheight{1.25}\smash{\begin{tabular}[t]{l}$H_{-3}$\end{tabular}}}}%
    \put(0,0){\includegraphics[width=\unitlength,page=4]{pyramid1.pdf}}%
    \put(0.85242778,0.01165401){\color[rgb]{0,0,0}\makebox(0,0)[lt]{\lineheight{1.25}\smash{\begin{tabular}[t]{l}$H_4$\end{tabular}}}}%
    \put(0,0){\includegraphics[width=\unitlength,page=5]{pyramid1.pdf}}%
  \end{picture}%
\endgroup%
}} 
\end{center}
\vspace{-7mm}
\caption{Pyramids $H_d$ in a triangle}
\vspace{-2mm}
\label{fig:pyramids}
\end{figure}

\begin{figure}[htbp!]
\vspace{-5mm}
\begin{center}
\hspace*{-12mm}  {\scalebox{1.0}{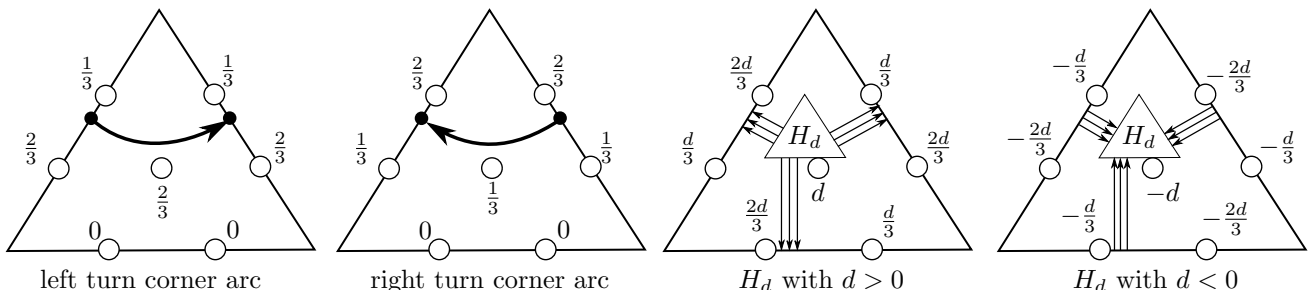}} 
\end{center}
\vspace{-5mm}
\caption{The tropical coordinates for elementary ${\rm SL}_3$-webs in a triangle}
\vspace{-4mm}
\label{fig:DS_coordinates}
\end{figure}

\begin{definition}[\cite{FS}; Def.\ref{def:canonical_wrt_split_ideal_triangulation}]
A non-elliptic ${\rm SL}_3$-web $W$ in a triangulable punctured surface $\frak{S}$ is in a \ul{\em canonical position} with respect to a split ideal triangulation $\wh{\Delta}$ of $\frak{S}$ for a triangulation $\Delta$, if
\begin{enumerate}
\itemsep0em
\item[\rm (1)] for each triangle $t$ of $\Delta$, $W\cap \wh{t}$ is {\em canonical} (Def.\ref{def:canonical_web_in_a_triangle}), i.e. consists of some finite number of left turn or right turn \ul{\em corner arcs}, each of which connects two distinct sides of $\wh{t}$, and a single \ul{\em degree $d$ pyramid $H_d$} for some $d\in \mathbb{Z}$ defined as in Def.\ref{def:canonical_web_in_a_triangle}\redfix{;} in particular $H_0 = {\O}$ and some examples of $H_d$ for nonzero $d$'s are as depicted in Fig.\ref{fig:pyramids};

\item[\rm (2)] for each biangle $B$ of $\wh{\Delta}$, $W\cap B$ is a {\em minimal crossbar ${\rm SL}_3$-web} in $B$ (Def.\ref{def:crossbar}), i.e. when the orientations are forgotten, is homeomorphic to the union of some finite number of simple arcs connecting two distinct sides of $B$, called {\em strands}, and some finite number of simple arcs connecting two adjacent strands, called {\em crossbars}, where the intersections of strands and crossbars are transverse double and lie in the interior of $B$, such that in between any two adjacent strands there is no consecutive crossbars that form a contractible 4-gon as in (S3) of Fig.\ref{fig:A2-skein_relations}, and that under a homeomorphism of the biangle $B$ to $\mathbb{R} \times [0,1]$ (the two sides going to $\mathbb{R}\times \{0\}$ and $\mathbb{R}\times \{1\}$) each strand is of the form $\{c\} \times [0,1]$ (i.e. vertical) and each crossbar is of the form $[c_1,c_2] \times \{a\}$ (i.e. horizontal).
\end{enumerate}
\end{definition}
For example, \raisebox{-0.5\height}{\scalebox{0.8}{
\begingroup%
  \makeatletter%
  \providecommand\color[2][]{%
    \errmessage{(Inkscape) Color is used for the text in Inkscape, but the package 'color.sty' is not loaded}%
    \renewcommand\color[2][]{}%
  }%
  \providecommand\transparent[1]{%
    \errmessage{(Inkscape) Transparency is used (non-zero) for the text in Inkscape, but the package 'transparent.sty' is not loaded}%
    \renewcommand\transparent[1]{}%
  }%
  \providecommand\rotatebox[2]{#2}%
  \newcommand*\fsize{\dimexpr\f@size pt\relax}%
  \newcommand*\lineheight[1]{\fontsize{\fsize}{#1\fsize}\selectfont}%
  \ifx\svgwidth\undefined%
    \setlength{\unitlength}{73.7007874bp}%
    \ifx\svgscale\undefined%
      \relax%
    \else%
      \setlength{\unitlength}{\unitlength * \real{\svgscale}}%
    \fi%
  \else%
    \setlength{\unitlength}{\svgwidth}%
  \fi%
  \global\let\svgwidth\undefined%
  \global\let\svgscale\undefined%
  \makeatother%
  \begin{picture}(1,0.53846154)%
    \lineheight{1}%
    \setlength\tabcolsep{0pt}%
    \put(0,0){\includegraphics[width=\unitlength,page=1]{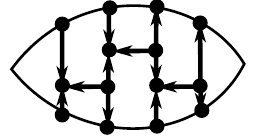}}%
  \end{picture}%
\endgroup%
}} is a minimal crossbar ${\rm SL}_3$-web in a biangle, with 4 strands and 3 crossbars; it has 14 vertices, where 6 among them are tri-valent, and 13 edges.
\begin{definition}[Douglas-Sun coordinates of non-elliptic ${\rm SL}_3$-webs; \cite{Douglas} \cite{DS1} \cite{DS2}; Def.\ref{def:Douglas-Sun}]
Let $W$ be a non-elliptic ${\rm SL}_3$-web in a triangulable punctured surface $\frak{S}$, with a chosen split ideal triangulation $\wh{\Delta}$ (for $\Delta$). Put $W$ into a canonical position with respect to $\wh{\Delta}$ by isotopy. For each triangle $t$ of $\Delta$, for the nodes $v$ of the $3$-triangulation quiver $Q_\Delta$ living in $t$, define ${\rm a}_v(W)$ as the sum of the coordinates at $v$ of the components of $W\cap \wh{t}$, as defined in Fig.\ref{fig:DS_coordinates}.
\end{definition}

\begin{definition}[the tropical coordinates for ${\rm SL}_3$-laminations; Def.\ref{def:tropical_coordinates}] Let $\ell$ be an ${\rm SL}_3$-lamination in a triangulable punctured surface $\frak{S}$, represented by an ${\rm SL}_3$-web $W$. Let $W_1,\ldots,W_m$ be the components of $W$, with weights $k_1,\ldots,k_m$ respectively. For each node $v\in \mathcal{V}(Q_\Delta)$, define the \ul{\em tropical coordinate} of the ${\rm SL}_3$-lamination $\ell$ at $v$ as the weighted sum of Douglas-Sun coordinates ${\rm a}_v(\ell) := \sum_{i=1}^m k_i {\rm a}_v(W_i)$.
\end{definition}
We also give explicit formulas for ${\rm a}_v(\ell)$ in terms of the Frohman-Sikora coordinates of ${\rm SL}_3$-webs (see Def.\ref{def:tropical_coordinates}). We then prove:
\begin{proposition}[well-definedness, and the image of the tropical coordinates of ${\rm SL}_3$-laminations; Prop.\ref{prop:tropical_coordinate_is_well-defined}, Prop.\ref{prop:tropical_coordinates_are_coordinates}]
\label{prop:intro_well-definedness_and_image_of_tropical_coordinates}
For an ideal triangulation $\Delta$ of a punctured surface $\frak{S}$, the tropical coordinate map
$$
{\bf a}_\Delta ~:~\{\mbox{${\rm SL}_3$-laminations in $\frak{S}$}\} \longrightarrow (\textstyle \frac{1}{3}\mathbb{Z})^{\mathcal{V}(Q_\Delta)}
$$
is well-defined and is a bijection onto the set of all \ul{\em balanced} elements of $(\frac{1}{3}\mathbb{Z})^{\mathcal{V}(Q_\Delta)}$, where an element $({\rm a}_v)_v \in (\frac{1}{3}\mathbb{Z})^{\mathcal{V}(Q_\Delta)}$ is said to be balanced if for each triangle $t$ of $\Delta$, if we denote by $e_1,e_2,e_3$ the sides of $t$ by $v_{e_1,1},v_{e_1,2},v_{e_2,1},v_{e_2,2},v_{e_3,1}\redfix{,}v_{e_3,2}$ the nodes of $Q_\Delta$ lying on the sides of $t$ appearing clockwise in this order (with $v_{e_\alpha,*}$ on $e_\alpha$), \redfix{and by $v_t$ the node of $Q_\Delta$ in the interior of $t$,} then \redfix{$\sum_{\alpha=1}^3 {\rm a}_{v_{e_\alpha,1}}$, $\sum_{\alpha=1}^3 {\rm a}_{v_{e_\alpha,2}}$, \, ${\rm a}_{v_{e_\beta,1}} + {\rm a}_{v_{e_\beta,2}}$ and $-{\rm a}_{v_t} + {\rm a}_{v_{e_\beta,2}} + {\rm a}_{v_{e_{\beta+1},1}}$ (for $\beta=1,2,3$) all belong to $\mathbb{Z}$.}
\end{proposition}

\redfix{For a combinatorial and representation-theoretic meaning of the balancedness condition in terms of the so-called Knutson-Tao rhombi \cite{KT99} \cite{GS15}, we refer the readers to \cite{DS2} (and also \cite{DS1}).}

\begin{proposition}[tropical coordinates \redfix{transform by} tropical $\mathscr{A}$-mutations; Prop.\ref{prop:coordinate_change_formulas}]
\label{prop:intro_tropical_change}
Let $\Delta,\Delta'$ be ideal triangulations of a punctured surface $\frak{S}$ related by a flip. Then the coordinates ${\bf a}_\Delta$ and ${\bf a}_{\Delta'}$ are related by the sequence of tropicalized versions of the cluster $\mathscr{A}$-mutations relating the cluster $\mathscr{A}$-charts of $\mathscr{A}_{{\rm SL}_3,\frak{S}}$ associated to $\Delta$ \redfix{and} $\Delta'$.
\end{proposition}
Prop.\ref{prop:intro_tropical_change} is a consequence of the corresponding statement for the Douglas-Sun coordinates of non-elliptic ${\rm SL}_3$-webs (Prop.\ref{prop:DS_coordinate_change}), which is proved \redfix{as a main result} in \cite{DS2}.

\begin{definition}[Def.\ref{def:Delta-congruent}]
For an ideal triangulation $\Delta$ of a punctured surface $\frak{S}$, an ${\rm SL}_3$-lamination $\ell$ in $\frak{S}$ is said to be \ul{\em $\Delta$-congruent} if the tropical coordinates ${\rm a}_v(\ell)$, $v\in \mathcal{V}(Q_\Delta)$, are all integers.

\vspace{0,5mm}

We say $\ell$ is \ul{\em congruent} if it is $\Delta$-congruent for all $\Delta$.
\end{definition}
\begin{proposition}[congruence is independent on triangulations; Prop.\ref{prop:congruence_condition_is_indepdent_on_triangulation}]
\label{prop:intro_congruence_independence}
For any two ideal triangulations $\Delta$ and $\Delta'$ of a punctured surface $\frak{S}$, an ${\rm SL}_3$-lamination $\ell$ in $\frak{S}$ is $\Delta$-congruent iff it is $\Delta'$-congruent.
\end{proposition}
Prop.\ref{prop:intro_congruence_independence} would be an easy corollary of Prop.\ref{prop:intro_tropical_change}, which in turn is a consequence of the corresponding result of \cite{DS2}. In the present paper, we present a proof of Prop.\ref{prop:intro_congruence_independence} independent of Prop.\ref{prop:intro_tropical_change}. Note that, although the statement of Prop.\ref{prop:intro_congruence_independence} itself is of topological and combinatorial nature, our proof of it heavily uses the proof of one of our main results (Thm.\ref{thm:intro_main}, \ref{thm:main}), which in turn is heavily algebraic. Anyways, consequently we have
\begin{theorem}[geometric model of tropical integer points of $\mathscr{A}_{{\rm SL}_3,\frak{S}}$; Thm.\ref{thm:geometric_model}]
\label{thm:intro_geometric_model}
Let $\frak{S}$ be a triangulable punctured surface. The tropical coordinate maps ${\bf a}_\Delta$ for ideal triangulations $\Delta$ of $\frak{S}$ provide bijections
$$
\{\mbox{congruent ${\rm SL}_3$-laminations in $\frak{S}$} \} ~\longrightarrow~ \mathbb{Z}^{\mathcal{V}(Q_\Delta)}
$$
which, under changes of ideal triangulations, transform by sequences of tropical versions of the cluster $\mathscr{A}$-mutations for the corresponding cluster $\mathscr{A}$-charts of $\mathscr{A}_{{\rm SL}_3,\frak{S}}$. Therefore we have the identification
$$
\{\mbox{congruent ${\rm SL}_3$-laminations in $\frak{S}$} \} ~ \longleftrightarrow ~  \mathscr{A}_{{\rm SL}_3,\frak{S}}(\mathbb{Z}^{\redfix{T}}).
$$
\end{theorem}
This is our solution to Conjecture \ref{conj:geometric_model_of_tropical_integer_points}. The full content of Thm.\ref{thm:intro_geometric_model} depends on Prop.\ref{prop:intro_tropical_change}; however, even without \cite{DS2}, we have a weaker version of Thm.\ref{thm:intro_geometric_model}, due to our proof of Prop.\ref{prop:intro_congruence_independence}. 

\vs

We note that, in the main text, the constructions and \redfix{the} statements in the present subsection are extended to any surface $\frak{S}$ having a boundary with marked points, called a {\em generalized marked surface} (Def.\ref{def:generalized_marked_surface}).

\subsection{The ${\rm SL}_3$-${\rm PGL}_3$ duality map: the first main theorem}

Let $\frak{S}$ be a triangulable punctured surface (i.e. without boundary). For each congruent ${\rm SL}_3$-lamination $\ell$ in $\frak{S}$, we should now describe our regular function $\mathbb{I}(\ell)$ on $\mathscr{X}_{{\rm PGL}_3,\frak{S}}$. We do this through several steps. First, we let
\begin{align*}
\mathscr{A}_{\rm L}^0(\frak{S};\mathbb{Z}) := \mbox{the set of all ${\rm SL}_3$-laminations in $\frak{S}$ with non-negative weights}.
\end{align*}
Since $\mathscr{A}_{\rm L}^0(\frak{S};\mathbb{Z})$ is in bijection with the set of all non-elliptic ${\rm SL}_3$-webs in $\frak{S}$, it embeds into the ${\rm SL}_3$-skein algebra $\mathcal{S}(\frak{S})$ as a basis \cite{SW} (see Cor.\ref{cor:A2-bangles_basis_for_L_SL3} of the present paper), hence from eq.\eqref{eq:intro_Phi} we get a map
$$
\mathbb{I}^0_{{\rm SL}_3} : \mathscr{A}^0_{\rm L}(\frak{S};\mathbb{Z}) \longrightarrow \mathscr{O}(\mathscr{L}_{{\rm SL}_3,\frak{S}})
$$
which is injective and whose image forms a basis of $\mathscr{O}(\mathscr{L}_{{\rm SL}_3,\frak{S}})$. Pullback of the natural frame-forgetting regular map
$$
F : \mathscr{X}_{{\rm SL}_3,\frak{S}}\longrightarrow \mathscr{L}_{{\rm SL}_3,\frak{S}}
$$
yields a map 
$$
F^* : \mathscr{O}(\mathscr{L}_{{\rm SL}_3,\frak{S}}) \longrightarrow \mathscr{O}(\mathscr{X}_{{\rm SL}_3,\frak{S}}).
$$
The gap between (the image under $F^*$ of) $\mathscr{O}(\mathscr{L}_{{\rm SL}_3,\frak{S}})$ and $\mathscr{O}(\mathscr{X}_{{\rm SL}_3,\frak{S}})$ is filled in by the regular functions on $\mathscr{X}_{{\rm SL}_3,\frak{S}}$ which read the framing data at punctures as follows (as in \cite{FG06}). Namely, for a framed ${\rm SL}_3$-local system on $\frak{S}$, the monodromy along a peripheral loop surrounding a puncture $p$ yields an element of a Borel subgroup \redfix{${\rm B}$} of ${\rm SL}_3$, and by reading the semi-simple part one gets an element of the {\em Cartan group} ${\rm H} := {\rm B}/{\rm U}$ of ${\rm SL}_3$. Choosing ${\rm B}$ to be the subgroup of all upper triangular matrices, ${\rm H}$ is isomorphic to the subgroup of all diagonal matrices. A monodromy is defined only up to conjugation, so from the monodromy alone we really get an element of ${\rm H}/{\rm W}$ where ${\rm W}$ is the Weyl group. However, the framing data pins down an element of ${\rm H}$ indeed, and we get a well-defined regular map
$$
\pi_p ~:~ \mathscr{X}_{{\rm SL}_3,\frak{S}} \longrightarrow {\rm H}
$$
for each puncture $p$; we elaborate more on this process in the main text (\S\ref{subsec:bases_of_rings}). Hence one obtains
$$
\mathscr{O}(\mathscr{X}_{{\rm SL}_3,\frak{S}}) \cong \mathscr{O}(\mathscr{L}_{{\rm SL}_3,\frak{S}}) \otimes_{\mathscr{O}( ({\rm H}/{\rm W})^\mathcal{P})} \mathscr{O}({\rm H}^{\mathcal{P}})
$$
as done in \cite{FG06}, where $\mathcal{P}$ is the set of all punctures of $\frak{S}$.  In the map $\pi_p$, by further reading the three diagonal entries of ${\rm H}$, we get three monomial regular functions $(\pi_p)_i : \mathscr{X}_{{\rm SL}_3,\frak{S}}  \to \mathbb{G}_m$, $i=1,2,3$, which fill in the gap between $\mathscr{L}_{{\rm SL}_3,\frak{S}}$ and the ${\rm SL}_3$ $\mathscr{X}$-moduli space $\mathscr{X}_{{\rm SL}_3,\frak{S}}$.
\begin{definition}[Def.\ref{def:I_SL3}]
\label{def:intro_I_SL3}
For a (triangulable) punctured surface $\frak{S}$, define the map
\begin{align}
\label{eq:intro_I_SL3}
\mathbb{I}_{{\rm SL}_3} ~:~ \mathscr{A}_{\rm L}(\frak{S};\mathbb{Z}) \longrightarrow \mathscr{O}(\mathscr{X}_{{\rm SL}_3,\frak{S}})
\end{align}
as follows. For $\ell \in \mathscr{A}_{\rm L}(\frak{S};\mathbb{Z})$, write $\ell = \ell_1 \cup \cdots \cup \ell_m$ as a disjoint union of ${\rm SL}_3$-laminations with weighted single-component ${\rm SL}_3$-webs. Let
$$
\mathbb{I}_{{\rm SL}_3}(\ell_i) := \left\{
\begin{array}{cl}
( (\pi_p)_1 )^{k_i} & \mbox{if $\ell_i$ is a peripheral loop going counterclockwise around $p$ with weight $k_i$,} \\
( (\pi_p)_3 )^{-k_i} & \mbox{if $\ell_i$ is a peripheral loop going clockwise around $p$ with weight $k_i$,} \\
F^*(\mathbb{I}^0_{{\rm SL}_3}(\ell_i)) & \mbox{otherwise},
\end{array}
\right.
$$
and let $\mathbb{I}_{{\rm SL}_3}(\ell) := \mathbb{I}_{{\rm SL}_3}(\ell_1) \cdot \cdots \cdot \mathbb{I}_{{\rm SL}_3}(\ell_m)$.
\end{definition}
We combine the results and the arguments above to prove:
\begin{proposition}[Prop.\ref{prop:I_SL3_is_a_basis}]
$\mathbb{I}_{{\rm SL}_3}$ is injective, and its image forms a basis of $\mathscr{O}(\mathscr{X}_{{\rm SL}_3,\frak{S}})$.
\end{proposition}
Our original interest is the space $\mathscr{X}_{{\rm PGL}_3,\frak{S}}$, for which the cluster $\mathscr{X}$-coordinate variables $X_v$, $v\in \mathcal{V}(Q_\Delta)$ are defined. In our final answer, for each $\ell \in \mathscr{A}_{{\rm SL}_3,\frak{S}}(\mathbb{Z}^{\redfix{T}})$ and $\Delta$, we shall construct a Laurent polynomial function in the variables $X_v$, $v\in \mathcal{V}(Q_\Delta)$. As a tool to relate the spaces $\mathscr{X}_{{\rm PGL}_3,\frak{S}}$ and $\mathscr{X}_{{\rm SL}_3,\frak{S}}$, we make use of the evaluations at the positive-real semi-field and the field of real numbers: let
\begin{align*}
\mathscr{X}^+_{{\rm PGL}_3, \frak{S}} := \mathscr{X}_{{\rm PGL}_3,\frak{S}}(\mathbb{R}_{>0}) ~~(~ \subset \mathscr{X}_{{\rm PGL}_3,\frak{S}}(\mathbb{R})~) ~ = ~  \mbox{the \ul{\em higher Teichm\"uller space}},
\end{align*}
so that $\mathscr{X}^+_{{\rm PGL}_3, \frak{S}}$ is given the structure of a smooth manifold \cite{FG06}. We construct a map
\begin{align}
\label{eq:intro_Psi}
\Psi ~:~ \mathscr{X}^+_{{\rm PGL}_3,\frak{S}} \longrightarrow \mathscr{X}_{{\rm SL}_3,\frak{S}}(\mathbb{R})
\end{align}
by using Fock-Goncharov's basic {\em monodromy matrices} \cite[\S9]{FG06}, suitably normalized (see \S\ref{subsec:lifting_PGL3_to_SL3} of the present paper). Namely, given positive real numbers $(X_v)_v \in (\mathbb{R}_{>0})^{\mathcal{V}(Q_\Delta)}$, we construct a point of $\mathscr{X}_{{\rm SL}_3,\frak{S}}(\mathbb{R})$, i.e. a framed ${\rm SL}_3(\mathbb{R})$-local system of $\frak{S}$ (\S\ref{subsec:lifting_PGL3_to_SL3}). The monodromy along each loop in $\frak{S}$ is constructed as the product of the basic monodromy matrices associated to small elementary pieces of this loop, where each basic monodromy matrix is an element of ${\rm SL}_3(\mathbb{R})$ whose entries depend on $(X_v)_v$; also the framing data can be constructed from this process, as done in \cite[\S9]{FG06}. In particular, by composing eq.\eqref{eq:intro_I_SL3} and the pullback of eq.\eqref{eq:intro_Psi} we obtain a map
$$
\mathbb{I}^+_{{\rm PGL}_3} ~:~ \mathscr{A}_{\rm L}(\frak{S}; \mathbb{Z}) \longrightarrow C^\infty(\mathscr{X}^+_{{\rm PGL}_3\redfix{, \frak{S}}}).
$$
For each $\ell \in \mathscr{A}_{\rm L}(\frak{S};\mathbb{Z})$, the smooth function $\mathbb{I}^+_{{\rm PGL}_3}\redfix{(\ell)}$ can be expressed, per each ideal triangulation $\Delta$ of $\frak{S}$, as a Laurent polynomial in the cube-root coordinate functions $\{X_v^{1/3} \, | \, v\in \mathcal{V}(Q_\Delta)\}$. It is this cube-root Laurent polynomial function $\mathbb{I}^+_{{\rm PGL}_3}(\ell)$ for some easy oriented loops $\ell$ that had been computed by some authors before; e.g. by Xie \cite{Xie} (see also \cite{CGT}), inspiring the tropical $\mathscr{A}$-coordinates for some elementary examples of ${\rm SL}_3$-webs \cite{Xie} \cite{Douglas} \cite{DS1}. If an ${\rm SL}_3$-web or an ${\rm SL}_3$-lamination $\ell$ has many tri-valent vertices, then the computation of $\mathbb{I}^+_{{\rm PGL}_3}(\ell)$ is quite difficult. By performing such computations, we prove:
\begin{proposition}[the highest term; Prop.\ref{prop:highest_term}]
\label{prop:intro_highest_term}
For $\ell \in \mathscr{A}_{\rm L}(\frak{S};\mathbb{Z})$, the expression of $\mathbb{I}^+_{{\rm PGL}_3}(\ell)$ as a Laurent polynomial in $\{X_v^{1/3}\,|\,v\in \mathcal{V}(Q_\Delta)\}$ for any chosen ideal triangulation $\Delta$ has the unique highest Laurent monomial in the natural partial ordering, and it is $\prod_{v\in \mathcal{V}(Q_\Delta)} X_v^{{\rm a}_v(\ell)}$, with coefficient $1$.
\end{proposition}
\begin{proposition}[congruence of terms; Prop.\ref{prop:congruence-compatibility_of_terms_of_each_basic_regular_function}]

For $\ell \in \mathscr{A}_{\rm L}(\frak{S};\mathbb{Z})$, we have \\
\hspace*{75mm} $\mathbb{I}^+_{{\rm PGL}_3}(\ell) \in (\prod_{v\in \mathcal{V}(Q_\Delta)} X_v^{{\rm a}_v(\ell)}) \cdot \mathbb{Z}[\{X_v^{\pm 1}\,|\, v\in \mathcal{V}(Q_\Delta)\}]$.
\end{proposition}
In particular, for a congruent ${\rm SL}_3$-lamination $\ell$, $\mathbb{I}^+_{{\rm PGL}_3}(\ell)$ is a Laurent polynomial in the (positive-real evaluations of the) usual cluster $\mathscr{X}$-variables $X_v$'s, $v\in \mathcal{V}(Q_\Delta)$. As a matter of fact, proving these innocent-looking assertions was much more of a challenge than it looked. For this we developed a whole set of new machinery throughout the entire section \S\ref{sec:SL3_trace}, which we call the \ul{\em ${\rm SL}_3$ classical trace maps}. A full treatment of them forces us to consider the surfaces with boundary with marked points.
\begin{definition}
Let $\frak{S}$ be a \ul{\em generalized marked surface} (Def.\ref{def:generalized_marked_surface}), i.e. a compact oriented surface with (possibly empty) boundary minus a finite set of points called {\em marked points}, such that each boundary component of $\frak{S}$ is homeomorphic to an open interval. For a ring $\mathcal{R}$, define the \ul{\em (commutative) stated} \ul{\em ${\rm SL}_3$-skein algebra} $\mathcal{S}_{\rm s}(\frak{S};\mathcal{R})$ as in Def.\ref{def:commutative_stated_A2-skein_algebra}, using ${\rm SL}_3$-webs $W$ in $\frak{S}$ that can have endpoints in $\partial \frak{S}$ (Def.\ref{def:A2-web}), together with a state $s:\partial W\to \{1,2,3\}$. For an ideal triangulation $\Delta$ of $\frak{S}$, let $\mathcal{Z}_\Delta := \mathbb{Z}[\{Z_v^{\pm 1}\,|\, v\in \mathcal{V}(Q_\Delta)\}]$ be the {\em classical cube-root Fock-Goncharov algebra} (Def.\ref{def:Fock-Goncharov_algebra_quantum}).
\end{definition}
\begin{proposition}[the ${\rm SL}_3$ classical trace map; Prop.\ref{cor:SL3_classical_trace_map}]
There exists a family of ring homomorphisms
$$
{\rm Tr}_{\Delta;\frak{S}} ~:~ \mathcal{S}_{\rm s}(\frak{S};\mathbb{Z}) \longrightarrow \mathcal{Z}_\Delta
$$
for each triangulable generalized marked surface $\frak{S}$ and its ideal triangulation $\Delta$, satisfying favorable properties, e.g. restricting to $\mathbb{I}^+_{{\rm PGL}_3}$ on the ${\rm SL}_3$-webs not containing peripheral loops, with $X_v = Z_v^3$.
\end{proposition}
The ${\rm SL}_3$ classical trace maps behave well under cutting along an arc of $\Delta$, and this cutting property yields a {\em state-sum} type formula for computing the values of ${\rm Tr}_{\Delta;\frak{S}}$; we give a bit more detail about these maps in the next subsection. What makes the computation of $\mathbb{I}^+_{{\rm PGL}_3}(\ell)$ difficult is the tri-valent vertices of ${\rm SL}_3$-webs, and the state-sum formula for the ${\rm SL}_3$ classical trace maps that we develop in \S\ref{sec:SL3_trace} allows us to push these tri-valent vertices to biangles of a split ideal triangulation $\wh{\Delta}$, so that the computation boils down essentially to computations of Reshetikhin-Turaev type invariants \cite{RT} in biangles, which we interpret as values of the counit of a commutative Hopf algebra $\mathcal{O}({\rm SL}_3)$, using the results of Higgins on the stated ${\rm SL}_3$-skein algebras \cite{Higgins}. As a useful by-product, we obtain expressions for the basic `monodromy matrices' for small tri-valent oriented graphs, as certain $3\times 3\times 3$ matrices.

\vs

Moreover, we also perform the computation of the effect on $\mathbb{I}^+_{{\rm PGL}_3}(\ell)$ of a single mutation at every possible node of $Q_\Delta$, which too \redfix{requires} quite heavy a calculation (\S\ref{subsec:mutation_of_basic_regular_functions})\redfix{, as well as usage of the machinery develop in \S\ref{sec:SL3_trace} on the ${\rm SL}_3$ classical trace maps}.
\begin{proposition}[the effect of a single mutation; Prop.\ref{prop:mutation_of_basic_semi-regular_function_at_interior_node_of_triangle}--\ref{prop:mutation_of_basic_semi-regular_function_at_edge_node_of_triangle}, Cor.\ref{cor:mutation_of_congruent_ell}]
For any ideal triangulation $\Delta$ of a punctured surface $\frak{S}$, any $\Delta$-congruent ${\rm SL}_3$-lamination $\ell$, and any node $v$ of the $3$-triangulation quiver $Q_\Delta$ of $\Delta$, if we write $\mathbb{I}^+_{{\rm PGL}_3}(\ell)$ as a function in the (cube-root) coordinate functions $\{X'_{v'} \, | \, v' \in \mathcal{V}(Q')\}$ for the cluster $\mathscr{X}$-chart of $\mathscr{X}_{{\rm PGL}_3,\frak{S}}$ (with quiver $Q' = \mu_v(Q_\Delta)$) obtained from the cluster $\mathscr{X}$-chart for $\Delta$ by mutating at $v$, it is a Laurent polynomial in $\{X'_{v'} \, | \, v' \in \mathcal{V}(Q')\}$.
\end{proposition}
Using these propositions, together with Shen's result $\mathscr{O}(\mathscr{X}_{{\rm PGL}_3,\frak{S}}) = \mathscr{O}_{\rm cl}(\mathscr{X}_{{\rm PGL}_3,\frak{S}})$ \cite{Shen} (Prop.\ref{prop:cluster_regular_is_regular}), and a theorem of Gross-Hacking-Keel \cite{GHK} (Prop.\ref{prop:one-level_mutations_are_enought}) stating that if a regular function on a cluster $\mathscr{X}$-chart stays regular after all possible single mutations then it is regular for {\em all} cluster $\mathscr{X}$-charts, we are able to prove the first main theorem of the present paper:
\begin{theorem}[the first main theorem: the ${\rm SL}_3$-${\rm PGL}_3$ duality map; Thm.\ref{thm:main}]
\label{thm:intro_main}
Let $\frak{S}$ be a triangulable punctured surface. The above process using $\mathbb{I}^+_{{\rm PGL}_3}$ yields a map
$$
\mathbb{I} ~:~ \mathscr{A}_{{\rm SL}_3,\frak{S}}(\mathbb{Z}^{\redfix{T}})  \longrightarrow \mathscr{O}_{\rm cl}(\mathscr{X}_{{\rm PGL}_3,\frak{S}})
$$
satisfying
\begin{enumerate}
\itemsep0em
\item[\rm (1)] $\mathbb{I}$ is injective, and the image of $\mathbb{I}$ forms a basis of $\mathscr{O}_{\rm cl}(\mathscr{X}_{{\rm PGL}_3,\frak{S}}) = \mathscr{O}(\mathscr{X}_{{\rm PGL}_3,\frak{S}})$, which we refer to as an \ul{\em $A_2$-bangles basis};

\item[\rm (2)] for each congruent ${\rm SL}_3$-lamination $\ell \in \mathscr{A}_{{\rm SL}_3,\frak{S}}(\mathbb{Z}^{\redfix{T}}) \subset \mathscr{A}_{\rm L}(\frak{S};\mathbb{Z})$ in $\frak{S}$ and any ideal triangulation $\Delta$, $\mathbb{I}(\ell)$ is a Laurent polynomial in the cluster $\mathscr{X}$-variables $\{X_v\,|\,v\in \mathcal{V}(Q_\Delta)\}$ with integer coefficients, with the unique highest term being $\prod_{v\in \mathcal{V}(Q_\Delta)} X_v^{{\rm a}_v(\ell)}$ with coefficient $1$, where $({\rm a}_v(\ell))_v\in \mathbb{Z}^{\mathcal{V}(Q_\Delta)}$ are the tropical coordinates constructed in the present paper;

\item[\rm (3)] if $\ell$ consists of peripheral loops only, then for each ideal triangulation $\Delta$, $\mathbb{I}(\ell) = \prod_{v\in \mathcal{V}(Q_\Delta)} X_v^{{\rm a}_v(\ell)}$;

\item[\rm (4)] the structure constants of this basis are integers, i.e. for each $\ell,\ell' \in \mathscr{A}_{{\rm SL}_3,\frak{S}}(\mathbb{Z}^{\redfix{T}})$, we have
\begin{align}
\nonumber
\mathbb{I}(\ell) \, \mathbb{I}(\ell') = \underset{\ell'' \in \mathscr{A}_{{\rm SL}_3,\frak{S}}(\mathbb{Z}^{\redfix{T}})}{\textstyle \sum}  c(\ell,\ell';\ell'') \, \mathbb{I}(\ell'')
\end{align}
where $c(\ell,\ell';\ell'')\in \mathbb{Z}$ and $c(\ell,\ell';\ell'')$ are zero for all but at most finitely many $\ell''$.
\end{enumerate}
\end{theorem}
As mentioned, we obtain a proof of Prop.\ref{prop:intro_congruence_independence} during our proof of Thm.\ref{thm:intro_main}. This first main theorem, Thm.\ref{thm:intro_main}, settles Conjecture \ref{conj:duality}, as promised. Note that one may write the domain and the codomain of this map $\mathbb{I}$ as $\mathscr{A}_{|Q_\Delta|}(\mathbb{Z}^{\redfix{T}})$ and $\mathscr{O}(\mathscr{X}_{|Q_\Delta|})$ respectively in terms of the cluster varieties. 

\subsection{The ${\rm SL}_3$ quantum trace map: the second main theorem}
\label{subsec:intro_quantum}

After finishing the first version of the present paper \cite{Kim}, we realized that the same framework of \cite[\S5]{Kim} for the constructions and proofs for the ${\rm SL}_3$ classical trace can be adapted to obtain a quantum version, and decided that it would be more economical to incorporate the quantum construction in the paper too. The domain of the sought-for \ul{\em ${\rm SL}_3$ quantum trace map} is a version of a non-commutative stated ${\rm SL}_3$-skein algebra modeled on stated ${\rm SL}_3$-webs living in the {\em thickened surface} $\frak{S} \times {\bf I}$ (Def.\ref{def:webs_in_3d}), which is a 3-manifold\redfix{; here, ${\bf I} = (-1,1)$}. So, a crossing of an ${\rm SL}_3$-web now carries an overpassing/underpassing information, and the defining relations should be written with coefficients in $\mathbb{Z}[\omega^{\pm 1/2}]$, where $\omega^{1/2} = q^{1/18}$ is a formal quantum parameter.
\begin{definition}[\cite{S05} \cite{FS} \cite{Higgins}; Def.\ref{def:stated_A2-skein_algebra_quantum}]
Let $\frak{S}$ be a generalized marked surface. \\ Define the \ul{\em (non-commutative) stated ${\rm SL}_3$-skein algebra} $\mathcal{S}^\omega_{\rm s}(\frak{S})$ as in Def.\ref{def:stated_A2-skein_algebra_quantum}, as the free $\mathbb{Z}[\omega^{\pm 1/2}]$-module generated by isotopy classes of ${\rm SL}_3$-webs in $\frak{S}\times {\bf I}$ together with states, mod out by the non-commutative ${\rm SL}_3$-skein relations in Fig.\ref{fig:A2-skein_relations_quantum}, with the product given by superposition. Define the \ul{\em reduced} stated ${\rm SL}_3$-skein algebra $\mathcal{S}^\omega_{\rm s}(\frak{S})_{\rm red}$ as the quotient of $\mathcal{S}^\omega_{\rm s}(\frak{S})$ by the boundary relations in Fig.\ref{fig:stated_boundary_relations}.
\end{definition}
The codomain of the ${\rm SL}_3$ quantum trace map is the following generalized quantum torus algebra.
\begin{definition}[\cite{FG06} \cite{GS19} \cite{Douglas} \cite{Douglas21}; Def.\ref{def:Fock-Goncharov_algebra_quantum}]
Let $\Delta$ be an ideal triangulation of a generalized marked surface $\frak{S}$. 

Define the 3-triangulation quiver $Q_\Delta$ as in the case of punctured surfaces, by gluing the quivers in Fig.\ref{fig:3-triangulation} for triangles of $\Delta$. Let $\varepsilon=(\varepsilon_{vw})_{v,w\in \mathcal{V}(Q_\Delta)}$ be the signed adjacency matrix for $Q_\Delta$ (Def.\ref{def:quiver_and_epsilon_matrix}), i.e.
$$
\varepsilon_{vw} = \mbox{(number of arrows from $v$ to $w$)} -  \mbox{(number of arrows from $w$ to $v$)}.
$$
For $v,w \in \mathcal{V}(Q_\Delta)$, let
\begin{align}
\nonumber
\wh{\varepsilon}_{vw} = \left\{
\begin{array}{cl}
\frac{1}{2} & \mbox{if $v,w$ lie in a same boundary arc of $\frak{S}$ and $\overrightarrow{vw}$ matches the boundary-orientation (Def.\ref{def:holed_surface})}, \\
-\frac{1}{2} & \mbox{if $v,w$ lie in a same boundary arc of $\frak{S}$ and $\overrightarrow{wv}$ matches the boundary-orientation (Def.\ref{def:holed_surface})}, \\
\varepsilon_{vw} & \mbox{otherwise}.
\end{array}
\right.
\end{align}
Let $\mathcal{Z}^\omega_\Delta$ be the \ul{\em (quantum non-commutative) cube-root Fock-Goncharov algebra} as the free associative $\mathbb{Z}[\omega^{\pm 1/2}]$-algebra generated by $\{\wh{Z}_v^{\pm 1} \, | \, v \in \mathcal{V}(Q_\Delta)\}$ with relations $\wh{Z}_v \wh{Z}_w = \omega^{2\wh{\varepsilon}_{vw}} \wh{Z}_w \wh{Z}_v$, $\forall v,w\in \mathcal{V}(Q_\Delta)$.
\end{definition}
When $\frak{S}$ is a punctured surface, the matrix $(\wh{\varepsilon}_{vw})_{v,w} = (\varepsilon_{vw})_{v,w}$ encodes the Poisson structure on the moduli space $\mathscr{X}_{{\rm PGL}_3,\frak{S}}$ given by $\{X_v, X_w\} = \varepsilon_{vw} X_v X_w$, or $\{Z_v, Z_w\} = \varepsilon_{vw} Z_v Z_w$ \cite{FG06}.

\vs

In \S\ref{sec:SL3_trace}, we prove the following, which is now the second main theorem of the paper. 
\begin{theorem}[the second main theorem: the ${\rm SL}_3$ quantum trace map; Thm.\ref{thm:SL3_quantum_trace_map}]
There exists a family of $\mathbb{Z}[\omega^{\pm 1/2}]$-algebra homomorphisms
$$
{\rm Tr}_{\Delta;\frak{S}}^\omega ~:~ \mathcal{S}^\omega_{\rm s}(\frak{S})_{\rm red} \longrightarrow \mathcal{Z}^\omega_\Delta
$$
for each triangulable generalized marked surface $\frak{S}$ and its ideal triangulation $\Delta$, satisfying favorable properties, e.g. the cutting/gluing property, restricting to the ${\rm SL}_3$ classical trace ${\rm Tr}_{\Delta;\frak{S}}$ when $\omega^{1/2}=1$.
\end{theorem}
This ${\rm SL}_3$ quantum trace map ${\rm Tr}^\omega_{\Delta;\frak{S}}$ can be viewed as a surface generalization of the Reshetikhin-Turaev  operator invariant for $\mathcal{U}_q(\frak{sl}_3)$ \redfix{(for biangles)} \cite{RT}, and is an ${\rm SL}_3$ version of Bonahon-Wong's celebrated ${\rm SL}_2$ quantum trace map \cite{BW}. We note that the ${\rm SL}_3$ quantum trace map is partially dealt with in \cite{Douglas} \cite{Douglas21} \cite{GLM} and \cite{Gabella} only for oriented curves and in \cite{CGT} a bit more generally for small surfaces, and here we provide a full version for the first time; we expect that our ${\rm SL}_3$ quantum trace map will be essentially obtained as a special case of the ${\rm SL}_n$ quantum trace in an upcoming work \cite{Lcoming}, although constructed by a different method than ours. In \S\ref{sec:SL3_trace} of the present paper, we construct a state-sum formula for ${\rm Tr}^\omega_{\Delta;\frak{S}}$, with the help of its cutting/gluing property. A crucial ingredient in this formula is the {\em biangle ${\rm SL}_3$ quantum trace} (Prop.\ref{prop:biangle_SL3_quantum_trace})
$$
{\rm Tr}^\omega_B ~:~ \mathcal{S}^\omega_{\rm s}(B)_{\rm red} \longrightarrow \mathbb{Z}[\omega^{\pm 1/2}],
$$
whose values can be viewed as the matrix elements of some Reshetikhin-Turaev type operator invariants of $\mathcal{U}_q(\frak{sl}_3)$. We show its existence by using the results of Higgins \cite{Higgins}, who proved that his version of the stated ${\rm SL}_3$-skein algebra of a biangle is isomorphic to the quantum group $\mathcal{O}_q({\rm SL}_3)$, which is a Hopf algebra; we take the counit of this Hopf algebra following the idea of \cite{CL}, and show that it satisfies the desired properties of ${\rm Tr}^\omega_B$. In practice, one actual difficult part was to come up with a suitable twist (eq.\eqref{eq:isomorphism_from_ours_to_Higgins}) of the stated ${\rm SL}_3$-skein algebra used by Higgins, so that the resulting stated ${\rm SL}_3$-skein algebra and the biangle ${\rm SL}_3$ quantum trace fit into our framework for the state-sum construction of the surface ${\rm SL}_3$ quantum trace maps, with all the desired properties satisfied. In particular, we show that our state-sum construction indeed provides a well-defined family of homomorphisms ${\rm Tr}^\omega_{\Delta;\frak{S}}$, by checking the invariance of the values under isotopy (\S\ref{subsec:isotopy_invariance}). \redfix{We note that such isotopy invariance for oriented loops was checked more or less in \cite{Douglas} \cite{Douglas21} \cite{CS}; here we provide hands-on proofs (without computer), and we deal with general ${\rm SL}_3$-webs}. 

\vs

The ${\rm SL}_3$ quantum trace map is itself an interesting and important object of study\redfix{; for example, it can be used to develop a representation theory of the ${\rm SL}_3$-skein algebras, in the style of Bonahon-Wong's work on the ${\rm SL}_2$ case \cite{BW16}}. \redfix{Meanwhile,} in the context of the present paper, it has two specific and significant roles. One is the role which its classical version played for our proof of the first main theorem (Thm.\ref{thm:intro_main}). The other is that it yields a quantum version of the duality map $\mathbb{I}$ of Thm.\ref{thm:intro_main}, via the arguments (\S\ref{subsec:quantum_duality_map}) similar to those in \cite{AK} for ${\rm SL}_2$-${\rm PGL}_2$.
\begin{theorem}[the ${\rm SL}_3$-${\rm PGL}_3$ quantum duality map; Thm.\ref{thm:quantum_duality_map}]
Let $\frak{S}$ be a triangulable punctured surface. For each ideal triangulation $\Delta$, define the {\em Fock-Goncharov algebra} $\mathcal{X}^q_\Delta$ as the free associative $\mathbb{Z}[\redfix{q^{\pm 1/18}}]$-algebra generated by $\{\wh{X}^{\pm 1}_v \, | \,v\in \mathcal{V}(Q_\Delta)\}$ with relations $\wh{X}_v \wh{X}_v = q^{2\varepsilon_{vw}} \wh{X}_w \wh{X}_v$, $\forall v,w\in \mathcal{V}(Q_\Delta)$. For each $\Delta$ there is a quantum duality map
$$
\mathbb{I}^q_\Delta ~:~ \mathscr{A}_{{\rm SL}_3,\frak{S}}(\mathbb{Z}^{\redfix{T}}) \longrightarrow \mathcal{X}^q_\Delta
$$
that recovers $\mathbb{I}$ as \redfix{$q^{1/18}=1$} with $\wh{X}_v\mapsto X_v$, and satisfies favorable properties analogous to those in Thm.\ref{thm:intro_main}.
\end{theorem}
In particular, for each $\Delta$ one constructs a map
$
\mathscr{O}(\mathscr{X}_{{\rm PGL}_3,\frak{S}}) \to \mathcal{X}^q_\Delta
$
sending $\mathbb{I}(\ell)$ to $\mathbb{I}^q_\Delta(\ell)$, for each $\ell \in \mathscr{A}_{{\rm SL}_3,\frak{S}}(\mathbb{Z}^\redfix{T})$. \redfix{It is shown in a follow-up work \cite[Thm.1.1]{Kim21} to the present paper that} these maps for different $\Delta$'s are compatible (\redfix{Prop.\ref{prop:compatibility_of_I_q}}) under the quantum coordinate change maps $\redfix{\Phi^q_{\Delta,\Delta'}} : {\rm Frac}(\mathcal{X}^q_{\Delta'}) \to {\rm Frac}(\mathcal{X}^q_\Delta)$ \redfix{(e.g. of \cite{FG09})} between the skew-fields of fractions of the Fock-Goncharov algebras, so that they collectively constitute a {\em deformation quantization map}
$$
\mathscr{O}(\mathscr{X}_{{\rm PGL}_3,\frak{S}}) 
\to \mathscr{O}^q_{\redfix{\rm tri}}(\mathscr{X}_{{\rm PGL}_3,\frak{S}}),
$$
where $\mathscr{O}^q_{\redfix{\rm tri}}(\mathscr{X}_{{\rm PGL}_3,\frak{S}}):= \bigcap_{\Delta'} \redfix{\Phi}^q_{\Delta,\Delta'}(\mathcal{X}^q_{\Delta'})$ is the algebra of quantum universally Laurent elements \redfix{for ideal triangulations}.

\vs

We discuss further research topics in \S\ref{sec:conjectures}.

\vs

\noindent{\bf Acknowledgments.} This work was supported by the National Research Foundation of Korea(NRF) grant
funded by the Korea government(MSIT) (No. 2020R1C1C1A01011151). The entire work grew out of discussions with Sangjib Kim on representation theory and invariant theory, and H. Kim thanks him for all the help and encouragement. H. Kim thanks Dylan Allegretti and Linhui Shen for answering massive questions about the topic, and thanks Vijay Higgins and Daniel Douglas for kindly explaining their works. H. Kim thanks Seung-Jo Jung and Thang L\^e for helpful discussions. \redfix{H. Kim thanks anonymous referees for helpful suggestions.}


\section{Moduli spaces of $A_2$ local systems on surfaces}

\subsection{Generalized marked surfaces and ideal triangulations}
\label{subsec:surface}

We first recall some basic definitions about the surfaces \redfix{used in the definition of} moduli spaces. We mostly adapt conventions used by L\^e \cite{Le17} \cite{Le18}.

\begin{definition}[\cite{Le17} \cite{Le18}]
\label{def:generalized_marked_surface}
$\bullet$ A \ul{\em generalized marked surface} $(\Sigma,\mathcal{P})$ consists of an oriented compact smooth surface $\Sigma$ with possibly-empty boundary $\partial \Sigma$ and a non-empty finite subset $\mathcal{P}$ of $\Sigma$. We always require that each component of $\partial \Sigma$ contains at least one point of $\mathcal{P}$. Two generalized marked surfaces $(\Sigma_1,\mathcal{P}_1)$ and $(\Sigma_2,\mathcal{P}_2)$ are \ul{\em isomorphic} if there exists an orientation-preserving diffeomorphism $\Sigma_1 \to \Sigma_2$ mapping $\mathcal{P}_1$ onto $\mathcal{P}_2$.

\vs

$\bullet$ The elements of $\mathcal{P}$ are called \ul{\em marked points}. A marked point that lies in the interior
$$
\mathring{\Sigma} := \Sigma\setminus\partial\Sigma
$$
of $\Sigma$ is called a \ul{\em puncture} of $(\Sigma,\mathcal{P})$.

\vs

$\bullet$ A generalized marked surface $(\Sigma,\mathcal{P})$ is called a \ul{\em punctured surface} if $\partial \Sigma = {\O}$.

\vs

$\bullet$ Each component of $(\partial \Sigma)\setminus \mathcal{P}$ is called a \ul{\em boundary arc} of $(\Sigma,\mathcal{P})$.
\end{definition}
For a generalized marked surface $(\Sigma,\mathcal{P})$, we often let
$$
\frak{S} := \Sigma\setminus\mathcal{P},
$$
and refer to $\frak{S}$ as the generalized marked surface, by identifying it with the data $(\Sigma,\mathcal{P})$, by abuse of notation. Note
$$
\partial \frak{S} = (\partial \Sigma) \setminus \mathcal{P}, \qquad
\mathring{\frak{S}} = \mathring{\Sigma}\setminus\mathcal{P}.
$$
A crucial ingredient is an {\em ideal triangulation} of the surface $\frak{S}$.
\begin{definition}[\cite{Le17} \cite{Le18}]
\label{def:triangulation}
Let $(\Sigma,\mathcal{P})$ be a generalized marked surface, and $\frak{S} = \Sigma\setminus\mathcal{P}$.

\vs

$\bullet$ A \ul{\em $\mathcal{P}$-arc in $\Sigma$}, or an \ul{\em ideal arc in $\frak{S}$}, is the image of an immersion $\alpha: [0,1] \to \Sigma$ such that $\alpha(\{0,1\}) \subset \mathcal{P}$ and $\alpha|_{(0,1)}$ is an embedding \redfix{into $\frak{S}$}. The elements of $\alpha(\{0,1\})$ are called \ul{\em endpoints}, and the subset $\alpha(0,1)$ is called the \ul{\em interior} of this $\mathcal{P}$-arc. Two $\mathcal{P}$-arcs are said to be \ul{\em isotopic} if they are isotopic within the class of $\mathcal{P}$-arcs. We often identify $\alpha$ with its image $\alpha([0,1])$ in $\Sigma$, or even with $\alpha((0,1))$ in $\frak{S}$.

\vs

$\bullet$ We say that $(\Sigma,\mathcal{P})$, or $\frak{S}$, is \ul{\em triangulable} \bluefix{if it} is none of the following:

-- \ul{\em monogon}, i.e. $\Sigma$ is diffeomorphic to a closed disc, and $\mathcal{P}$ consists of a single point on $\partial \Sigma$;

-- \ul{\em biangle}, i.e. $\Sigma$ is diffeomorphic to a closed disc, and $\mathcal{P}$ consists of two points on $\partial \Sigma$;

-- sphere with one or two punctures, i.e. $\Sigma$ is diffeomorphic to the sphere $S^2$, and $|\mathcal{P}|\le 2$.

\vs

$\bullet$ When $(\Sigma,\mathcal{P})$ is triangulable, a \ul{\em $\mathcal{P}$-triangulation of $\Sigma$}, or \ul{\em an ideal triangulation of $\frak{S}$}, is defined as a collection $\Delta$ of $\mathcal{P}$-arcs in $\Sigma$ such that
\begin{enumerate}
\itemsep0em
\item[\rm (IT1)] no arc in $\Delta$ bounds a disk whose interior is in $\Sigma\setminus\mathcal{P}$;

\item[\rm (IT2)] no two members of $\Delta$ are isotopic or intersect in $\Sigma\setminus \mathcal{P}$;

\item[\rm (IT3)] $\Delta$ is maximal among the collections satisfying (IT1) and (IT2).
\end{enumerate}
Two ideal triangulations are \ul{\em isotopic} if they are related by a simultaneous isotopy of their members, within the class of ideal triangulations.

\vs

We assume (by applying an isotopy if necessary) that each constituent arc of $\Delta$ that is isotopic to a boundary arc of $\frak{S}$ is indeed a boundary arc of $\frak{S}$. The constituent arcs of $\Delta$ that are not boundary arcs are called \ul{\em internal arcs} of $\Delta$. Constituent arcs of $\Delta$ are often called \ul{\em edges} of $\Delta$.
\end{definition}
\begin{remark}
The monogon and biangle will play crucial roles later in the present paper.
\end{remark}
Let $\Delta$ be an ideal triangulation of a triangulable generalized marked surface $(\Sigma,\mathcal{P})$. Let $\mathring{t}$ be a connected component of the complement $\Sigma\setminus (\bigcup_{e\in \Delta} e)$. The closure $t$ of $\mathring{t}$ in $\Sigma$ is called an \ul{\em ideal triangle} of $\Delta$. Let
$$
\mathcal{F}(\Delta) := \mbox{the set of all ideal triangles of $\Delta$}.
$$
Note that $t\setminus \mathring{t}$ is union of two or three ideal arcs in $\Delta$, which are called \ul{\em sides} of $t$. In case $t$ has only two sides, $t$ is called \ul{\em self-folded}. In the present paper, we do not allow ideal triangulation having a self-folded ideal triangle. In fact, we only use ideal triangulations satisfying a somewhat stronger condition.
\begin{definition}[from \cite{FG06}]
\label{def:regular_triangulation}
An ideal triangulation $\Delta$ of a triangulable generalized marked surface $(\Sigma,\mathcal{P})$ is \ul{\em regular} if for each puncture $p$ of $(\Sigma,\mathcal{P})$, the valence of $\Delta$ at $p$ is at least $3$, where the valence of $\Delta$ at $p$ is the total number of arcs of $\Delta$ meeting $p$ counted with multiplicity, where the multiplicity of an arc is $1$ if $p$ is exactly one of the two distinct endpoints of the arc and is $2$ if both endpoints of the arc coincide with $p$. 
\end{definition}
We require that $\frak{S}$ admits at least one regular ideal triangulation. Especially, our main theorems will be on punctured surfaces, hence we need the following observation:
\begin{lemma}
Every triangulable punctured surface except for the sphere with three punctures admits at least one regular ideal triangulation.
\end{lemma}
{\it Proof.} Let's use induction. Suppose that a triangulable punctured surface $(\Sigma,\mathcal{P)}$, where $\Sigma$ is a compact oriented surface of genus $g$ and $|\mathcal{P}|=n$, admits a regular ideal triangulation $\Delta$. Let $\mathcal{P}' := \mathcal{P} \cup \{x\}$ where $x$ is a point lying in the interior of some ideal triangle $t$ of $\Delta$. Then, by adding to $\Delta$ three $\mathcal{P}'$-arcs, each connecting $x$ and a vertex marked point $\in\mathcal{P}$ of $t$, one obtains an ideal triangulation $\Delta'$ of $(\Sigma,\mathcal{P}')$. The valence of $\Delta'$ at $x$ is $3$, and the valence of $\Delta'$ at each $p\in \mathcal{P}$ is at least the valence of $\Delta$ at $p$ hence is at least $3$. So $\Delta'$ is regular. Hence, for each genus $g$ surface $\Sigma$, it suffices to prove the statement for $\mathcal{P}$ with minimal possible $|\mathcal{P}|$. For genus $0$ surface, i.e. sphere, when $|\mathcal{P}|=4$, one can easily find a regular triangulation. For genus $g\ge 1$ surface $\Sigma$ with $|\mathcal{P}|=1$, it is well known that $(\Sigma,\mathcal{P})$ admits an ideal triangulation $\Delta$, and that any ideal triangulation of it has $6g-6+3 = 6g-3$ arcs. So the unique $\mathcal{P}$ has valence $2(6g-3) \ge 3$, hence $\Delta$ is regular. \qed

\vs

Basic constructions of Fock-Goncharov's higher Teichm\"uller theory \cite{FG06} make use of the choice of a regular ideal triangulation $\Delta$ of a generalized marked surface $\frak{S}$. A key point is then to assure certain compatibility under changes of ideal triangulations. One common strategy is to focus on certain elementary changes called {\em flips}, which change an ideal triangulation by only one edge at a time.

\begin{definition}
We say that two ideal triangulations of a generalized marked surface $\frak{S}$, defined up to isotopy, are related by a \ul{\em flip} if they differ by exactly one edge. 
\end{definition}
If $\Delta$ and $\Delta'$ are related by a flip, then we have a natural bijection between $\Delta$ and $\Delta'$ as sets, through which we identify the two sets. If the only differing edge is labeled by $e$ (both in $\Delta$ and $\Delta'$), we say that this flip is a \ul{\em flip at the edge $e$}.

\vs

One thing to keep in mind is:
$$
\mbox{Throughout the paper, by an \ul{\em ideal triangulation} we mean a regular ideal triangulation.}
$$
In particular, when we refer to a \ul{\em triangulable} generalized marked surface, we mean one that admits at least one regular ideal triangulation, and we will only consider flips between regular ideal triangulations.

\begin{remark}
\label{rem:non-regular_triangulations}
Most, but perhaps not all, of the constructions and proofs of the present paper immediately apply also for ideal triangulations $\Delta$ that have punctures of valence $2$, i.e. that have `eyes' in the 3-valent (fat) graph dual to $\Delta$ in the language of \cite[\S3.8]{FG06}. To be safe, we only focus on the regular ideal triangulations throughout the paper. However, after proving the main theorems, one will be able to obtain the statements for non-regular ideal triangulations in a suitable sense, even when there are punctures of valence $1$, i.e. $\Delta$ has self-folded triangles or `viruses'. But we will not elaborate on this \redfix{in the present paper. We refer the readers to the subsequent paper \cite{JK} on non-regular triangulations}.
\end{remark}

\subsection{Moduli spaces of ${\rm G}$-local systems on surfaces}

We recall Fock-Goncharov's versions \cite{FG06} of moduli spaces of ${\rm G}$-local systems on a triangulable generalized marked surface $\frak{S}$, where ${\rm G}$ is a split reductive algebraic group. A {\em ${\rm G}$-local system} $\mathcal{L}$ on $\frak{S}$ may be understood as a right principal ${\rm G}$-bundle on $\frak{S}$ together with a flat ${\rm G}$-connection on it. 
\begin{definition}[$\mathscr{L}$-moduli space; \cite{FG06}]
Let $\mathscr{L}_{{\rm G},\frak{S}}$ be the moduli stack parametrizing all isomorphism classes of ${\rm G}$-local systems on $\frak{S}$.
\end{definition}
A ${\rm G}$-local system $\mathcal{L}$ induces a group homomorphism $\pi_1(\frak{S}) \to {\rm G}$ defined up to conjugation by an element of ${\rm G}$, which is called a {\em monodromy} representation of $\mathcal{L}$, and which in fact determines the isomorphism class of $\mathcal{L}$. Hence we have a natural identification
$$
\mathscr{L}_{{\rm G},\frak{S}} \cong {\rm Hom}(\pi_1(\frak{S}), {\rm G})/{\rm G}
$$
where the right hand side is the quotient of ${\rm Hom}(\pi_1(\frak{S}), {\rm G})$ by the action of ${\rm G}$ by conjugation. 

\vs

Choose any Borel subgroup ${\rm B}$ of ${\rm G}$, and let ${\rm U}:=[{\rm B}, {\rm B}]$ be a maximal unipotent subgroup of ${\rm G}$. Let $\mathcal{B} = {\rm G} / {\rm B}$ be the {\em flag variety} for ${\rm G}$.  For a ${\rm G}$-local system $\mathcal{L}$ on $\frak{S}$, denote by
$$
\mathcal{L}_\mathcal{B} := \mathcal{L} \times_{\rm G} \mathcal{B} \qquad\mbox{and}\qquad
\mathcal{L}_\mathcal{A} := \mathcal{L}/{\rm U}
$$
the associated {\em flag bundle} and {\em principal affine bundle} on $\frak{S}$, respectively; each of these associated bundles is also naturally equipped with a flat connection induced by $\mathcal{L}$. 

\vs

The present paper concerns ${\rm G} = {\rm PGL}_3$ and ${\rm SL}_3$. For ${\rm G} = {\rm PGL}_3$ we may choose ${\rm B}$ to be the subgroup of all upper triangular matrices in ${\rm PGL}_3$, and then for ${\rm G}={\rm SL}_3$ we choose ${\rm U}$ to be the subgroup of all upper triangular matrices with all diagonal entries being $1$.

\vs

To describe an extra boundary data, it is more convenient to deal with a {\em holed} surface, instead of a punctured surface. We fix a notation, which is similar but slightly different from what is in \cite{FG06}.
\begin{definition}[holed surface]
\label{def:holed_surface}
Let $(\Sigma,\mathcal{P})$ be a generalized marked surface, with $\frak{S} = \Sigma\setminus\mathcal{P}$. For each puncture $p$ of $(\Sigma,\mathcal{P})$, choose a neighborhood $N_p$ of $p$ in $\Sigma$ diffeomorphic to an open disc. Let
$$
\til{\frak{S}} := \Sigma \setminus \textstyle \bigcup_{p \, : \, {\rm puncture ~of}~ (\Sigma,\mathcal{P})} N_p
$$
be a \ul{\em holed surface} for $\frak{S}$. Each boundary component of $\til{\frak{S}}$ that is diffeomorphic to a circle is called a \ul{\em hole} of $\til{\frak{S}}$, and other boundary components of $\til{\frak{S}}$ are called \ul{\em boundary arcs} of $\til{\frak{S}}$.

\vs

The \ul{\em boundary-orientation} on a hole of $\til{\frak{S}}$ is the orientation induced by the counterclockwise orientation along the boundary of the disc $N_p$ (hence is `clockwise' when viewed from the interior of $\til{\frak{S}}$).
\end{definition}
Holes of $\til{\frak{S}}$ correspond to punctures of $\frak{S}$. One can view $\til{\frak{S}}$ as being embedded as a subspace of $\frak{S}$, onto which $\frak{S}$ deformation retracts. In particular, the inclusion $\til{\frak{S}} \to \frak{S}$ naturally induces isomorphism $\pi_1(\til{\frak{S}}) \to \pi_1(\frak{S})$ and hence an identification of $\mathscr{L}_{{\rm G},\til{\frak{S}}}$ with $\mathscr{L}_{{\rm G},\frak{S}}$. We will use $\frak{S}$ and $\til{\frak{S}}$ interchangeably. In particular, a hole of $\til{\frak{S}}$ is sometimes regarded as the oriented loop in $\frak{S}$, oriented according to the boundary-orientation as in Def.\ref{def:holed_surface}. The boundary arcs of $\til{\frak{S}}$ are naturally identified with boundary arcs of $\frak{S}$.

\begin{definition}[$\mathscr{X}$-moduli space; \cite{FG06} \cite{A19}]
Let ${\rm G} = {\rm PGL}_3$. A \ul{\em framing} for a ${\rm G}$-local system $\mathcal{L}$ on $\frak{S}$ is a flat section $\beta$ of the restriction of $\mathcal{L}_\mathcal{B}$ to $\partial \til{\frak{S}}$. A pair $(\mathcal{L},\beta)$ is called a \ul{\em framed ${\rm G}$-local system} on $\frak{S}$. Two framed ${\rm G}$-local systems $(\mathcal{L}_1,\beta_1)$ and $(\mathcal{L}_2,\beta_2)$ are \ul{\em isomorphic} if there is an isomorphism $\mathcal{L}_1 \to \mathcal{L}_2$ of ${\rm G}$-local systems whose induced map $(\mathcal{L}_1)_\mathcal{B} \to (\mathcal{L}_2)_\mathcal{B}$ sends $\beta_1$ to $\beta_2$. Let $\mathscr{X}_{{\rm G},\frak{S}}$ be the moduli stack parametrizing all isomorphism classes of framed ${\rm G}$-local systems on $\frak{S}$.
\end{definition}
\begin{definition}[$\mathscr{A}$-moduli space; \cite{FG06} \cite{A19}]
Let ${\rm G}={\rm SL}_3$. A \ul{\em decoration} for a ${\rm G}$-local system $\mathcal{L}$ on $\frak{S}$ is a flat section $\alpha$ of the restriction of $\mathcal{L}_\mathcal{A}$ to $\partial \til{\frak{S}}$. A pair $(\mathcal{L},\alpha)$ is called a \ul{\em decorated ${\rm G}$-local system} on $\frak{S}$. Two decorated ${\rm G}$-local systems $(\mathcal{L}_1,\alpha_1)$ and $(\mathcal{L}_2,\alpha_2)$ are \ul{\em isomorphic} if there is an isomorphism $\mathcal{L}_1\to\mathcal{L}_2$ of ${\rm G}$-local systems whose induced map $(\mathcal{L}_1)_{\redfix{\mathcal{A}}} \to (\mathcal{L}_2)_{\redfix{\mathcal{A}}}$ sends $\alpha_1$ to $\alpha_2$. Let $\mathscr{A}_{{\rm G},\frak{S}}$ be the moduli stack parametrizing all isomorphism classes of decorated ${\rm G}$-local systems on $\frak{S}$.
\end{definition}
The moduli spaces $\mathscr{X}_{{\rm G},\frak{S}}$ and $\mathscr{A}_{{\rm G},\frak{S}}$ can be defined for other groups ${\rm G}$. For general cases, the definition of $\mathscr{A}_{{\rm G},\frak{S}}$ is more complicated; see \cite{FG06} \cite{A19}. 

\subsection{Cluster atlases}
\label{subsec:cluster_atlases}

Let $\frak{S}$ be a triangulable generalized marked surface. Fock and Goncharov \cite{FG06} constructed special coordinate systems for $\mathscr{A}_{{\rm SL}_3,\frak{S}}$ and $\mathscr{X}_{{\rm PGL}_3,\frak{S}}$ respectively, per each choice of an ideal triangulation $\Delta$ of $\frak{S}$. They showed that, upon each change of ideal triangulations, the coordinates transform according to the {\em mutation} formulas appearing in the theory of cluster varieties. We first recall and define some basic notions needed in this theory.

\begin{definition}
\label{def:quiver_and_epsilon_matrix}
By a \ul{\em quiver} $Q$ we mean a directed graph without cycles of length $1$ or $2$. Its vertices are called \ul{\em nodes} of $Q$ and depicted as hollow circles $\circ$, while its oriented edges are called \ul{\em arrows} of $Q$ and depicted for example as $\overset{v}{\circ} \to \overset{w}{\circ}$. Denote by $\mathcal{V}(Q)$ and $\mathcal{E}(Q)$ the set of all nodes and the set of all arrows of $Q$.

\vs

The \ul{\em signed adjacency matrix} of a quiver $Q$ is the $\mathcal{V}(Q) \times \mathcal{V}(Q)$ matrix $\varepsilon_Q = \varepsilon$ whose $(v,w)$-th entry is
$$
\varepsilon_{vw} = \varepsilon_{v,w} = \mbox{(number of arrows from $v$ to $w$)} -  \mbox{(number of arrows from $w$ to $v$)}.
$$
\end{definition}
In the present paper, the \ul{\em $(\alpha,\beta)$-th entry} of a matrix refers to the entry in the $\alpha$-th row and $\beta$-th column. 

\vs

\redfix{Once we fix a surface $\frak{S}$, we will be dealing only with quivers $Q$ with a fixed set of nodes $\mathcal{V}$. The set $\mathcal{E}(Q)$ of arrows should be understood as a multiset of elements of $\mathcal{V} \times \mathcal{V}$. We identify two quivers $Q$ and $Q'$ if and only if they have the same sets of arrows, which is equivalent to $\varepsilon_Q = \varepsilon_{Q'}$. In particular, even when $Q$ and $Q'$ are isomorphic as quivers, we do not necessarily identify them. The correspondence $Q\leftrightarrow \varepsilon_Q$ is a bijection between the set of all quivers having $\mathcal{V}$ as the sets of nodes and the set of all $\mathcal{V} \times \mathcal{V}$ skew-symmetric integer matrices.}

\vs

We need to recall a certain transformation rule for quivers.
\begin{definition}
\label{def:quiver_mutation}
Let $Q$ be a quiver with the set of nodes $\mathcal{V}$ and the signed adjacency matrix $\varepsilon$. Let $k\in \mathcal{V}$. The \ul{\em quiver mutation} $\mu_k$ at the node $k$ transforms the quiver $Q$ into another quiver $Q' = \mu_k(Q)$ whose set of nodes is $\mathcal{V}$ and the signed adjacency matrix $\varepsilon'$ is defined as
$$
\varepsilon'_{ij} = \left\{
\begin{array}{ll}
- \varepsilon_{ij} & \mbox{if $k\in \{i,j\}$,} \\
\varepsilon_{ij} + \frac{1}{2}(\varepsilon_{ik}|\varepsilon_{kj}| + |\varepsilon_{ik}| \varepsilon_{kj}) & \mbox{if $k\notin\{i,j\}$.}
\end{array}
\right.
$$
\end{definition}
As is well-known, the quiver mutation can be described combinatorially as follows. First, from $Q$, reverse the orientations of all arrows starting or ending at $k$. Second, for each pair of an incoming arrow $\overset{i}{\circ} \to \overset{k}{\circ}$ and an outgoing arrow $\overset{k}{\circ} \to \overset{j}{\circ}$ at the node $k$, add the arrow $\overset{j}{\circ} \to \overset{i}{\circ}$ (i.e. `complete the $3$-cycle through $k$'). Finally, remove cycles of length $2$, until there are none. Then the resulting quiver is $\mu_k(Q)$.

\vs

To characterize Fock-Goncharov's special coordinate systems on the moduli spaces $\mathscr{A}_{{\rm SL}_3,\frak{S}}$ and $\mathscr{X}_{{\rm PGL}_3,\frak{S}}$, we establish some terminology, based on \cite{FG06}.
\begin{definition}
\label{def:A}
Let $\mathscr{S}$ be an irreducible stack or a scheme. A \ul{\em cluster $\mathscr{A}$-chart} on $\mathscr{S}$ is a pair $(Q,\psi)$, where $Q$ is a (labeled) quiver and
$$
\psi : \mathscr{S} \dashrightarrow (\mathbb{G}_m)^{\mathcal{V}(Q)}
$$
is a birational map, providing a rational coordinate system for $\mathscr{S}$. Denote by $A_i$ the coordinate function for the node $i \in \mathcal{V}(Q)$, which is called a \ul{\em cluster $\mathscr{A}$-variable} for this chart.

\vs

We say that a cluster $\mathscr{A}$-chart $(Q,\psi)$ is related to another cluster $\mathscr{A}$-chart $\mu_k(Q,\psi) = (Q',\psi')$ by the \ul{\em cluster $\mathscr{A}$-mutation} at the node $k\in \mathcal{V}(Q)$ if $Q' = \mu_k(Q)$ holds, so that we have an identification of $\mathcal{V}(Q)$ and $\mathcal{V}(Q')$, and the coordinate functions $A_i'$ for $\psi'$ are related to those $A_i$ of $\psi$ as
$$
A_i' = \left\{
\begin{array}{ll}
A_i & \mbox{if $i\neq k$,} \\
A_k^{-1}(\prod_j A_j^{[\varepsilon_{jk}]_+} + \prod_j A_j^{[-\varepsilon_{jk}]_+}) & \mbox{if $i=k$},
\end{array}
\right.
$$
where $[a]_+$ is the {\em positive part} of a real number $a$:
\begin{align}
\label{eq:positive_part}
[a]_+ := \left\{
\begin{array}{ll}
a & \mbox{if $a\ge 0$,} \\
0 & \mbox{if $a<0$}.
\end{array}
\right.
\end{align}
\end{definition}
Two cluster $\mathscr{A}$-charts $(Q,\psi)$ and $(Q',\psi')$ are identified if each coordinate function $A_i$ for the former coincides with the corresponding coordinate function $A_i'$ for the latter, on an open dense subset. 
\begin{definition}
A \ul{\em cluster $\mathscr{A}$-atlas} on an irreducible stack or a scheme $\mathscr{S}$ is a collection $\mathcal{C}$ of cluster $\mathscr{A}$-charts on $\mathscr{S}$ such that each two members of $\mathcal{C}$ are related by a finite sequence of members of $\mathcal{C}$ such that each pair of consecutive membe\redfix{r}s are related by a cluster $\mathscr{A}$-mutation.
\end{definition}
Note that cluster $\mathscr{A}$-mutation can also be used as a tool to construct from a cluster $\mathscr{A}$-chart another cluster $\mathscr{A}$-chart. So, starting from any cluster $\mathscr{A}$-atlas, by adding all cluster $\mathscr{A}$-charts obtained by applying cluster $\mathscr{A}$-mutations, one gets a uniquely determined maximal cluster $\mathscr{A}$-atlas.

\begin{definition}
\label{def:X}
Let $\mathscr{S}$ be an irreducible stack or a scheme. A \ul{\em cluster $\mathscr{X}$-chart} on $\mathscr{S}$ is a pair $(Q,\psi)$, where $Q$ is a (labeled) quiver and
$$
\psi : \mathscr{S} \dashrightarrow (\mathbb{G}_m)^{\mathcal{V}(Q)}
$$
is a birational map, providing a rational coordinate system for $\mathscr{S}$. Denote by $X_i$ the coordinate function for the node $i \in \mathcal{V}(Q)$, which is called a \ul{\em cluster $\mathscr{X}$-variable} for this chart.

\vs

We say that a cluster $\mathscr{X}$-chart $(Q,\psi)$ is related to another cluster $\mathscr{X}$-chart $\mu_k(Q,\psi) = (Q',\psi')$ by the \ul{\em cluster $\mathscr{X}$-mutation} at the node $k\in \mathcal{V}(Q)$ if $Q' = \mu_k(Q)$ holds, so that we have an identification of $\mathcal{V}(Q)$ and $\mathcal{V}(Q')$, and the coordinate functions $X_i'$ for $\psi'$ are related to those $X_i$ of $\psi$ as
\begin{align}
\label{eq:X-mutation_formula}
X_i' = \left\{
\begin{array}{ll}
X_k^{-1} & \mbox{if $i=k$} \\
X_i(1+X_k^{-{\rm sgn}(\varepsilon_{ik})})^{-\varepsilon_{ik}} & \mbox{if $i\neq k$},
\end{array}
\right.
\end{align}
where ${\rm sgn}(a)$ means the sign of a real number:
\begin{align}
\label{eq:sgn}
{\rm sgn}(a) = \left\{
\begin{array}{rl}
1 & \mbox{if $a>0$,} \\
0 & \mbox{if $a=0$,} \\
-1 & \mbox{if $a<0$.}
\end{array}
\right.
\end{align}
\end{definition}
\begin{definition}
Define the notion of cluster $\mathscr{X}$-atlas similarly as for cluster $\mathscr{A}$-atlas.
\end{definition}
One of the major objects of study in the theory of cluster varieties is the following.
\begin{definition}[Ring of cluster $\mathscr{X}$-regular functions; \cite{GS18} \cite{Shen}]
\label{def:O_cl}
Let $\mathscr{S}$ be an irreducible stack or a scheme, equipped with a chosen cluster $\mathscr{X}$-atlas. Let $\mathscr{O}_{\rm cl}(\mathscr{S})$ the ring of all \ul{\em cluster $\mathscr{X}$-regular} functions on $\mathscr{S}$, i.e. the rational functions on $\mathscr{S}$ that are regular on each cluster $\mathscr{X}$-chart belonging to the maximal cluster $\mathscr{X}$-atlas determined by the given cluster $\mathscr{X}$-atlas.
\end{definition}
Observe that a rational function defined on a cluster $\mathscr{X}$-chart $(Q,\psi)$ is regular on this chart iff it can be written as a Laurent polynomial in the coordinate functions $X_i$, $i \in \mathcal{V}(Q)$, for this chart. So a cluster $\mathscr{X}$-regular function is often referred to as \ul{\em universally ($\mathscr{X}$-)Laurent (polynomial)} functions. Note that the similarly defined ring of all cluster $\mathscr{A}$-regular functions coincides with the notion of the so-called {\em upper cluster algebra}. However, in the present paper, only the cluster $\mathscr{X}$-regular functions will be dealt with.

\vs

Fock and Goncharov showed \cite{FG06} that the moduli spaces $\mathscr{A}_{{\rm SL}_3,\frak{S}}$ and $\mathscr{X}_{{\rm PGL}_3,\frak{S}}$ exhibit a cluster $\mathscr{A}$-atlas and a cluster $\mathscr{X}$-atlas respectively, with some special quivers $Q = Q_\Delta$ associated to ideal triangulations $\Delta$ of $\frak{S}$.

\begin{definition}[\cite{FG06}]
\label{def:3-triangulation}
Let $\Delta$ be an ideal triangulation of a triangulable generalized marked surface $\frak{S}$. Let $Q_\Delta$ be the quiver, which may be drawn on the surface $\frak{S}$, whose nodes and arrows are defined as follows.

\vs

For each edge $e$ of $\Delta$, there are two nodes of $Q_\Delta$ lying in the interior of $e$. For each triangle $t$ of $\Delta$, there is one node of $Q_\Delta$ lying in the interior of $t$.

\vs

For each triangle of $\Delta$, consider a quiver as depicted in Fig.\ref{fig:3-triangulation}, consisting of three small counterclockwise cycles of length $3$. The quiver $Q_\Delta$ is obtained by gluing these quivers for all triangles of $\Delta$.

\vs

The quiver $Q_\Delta$ is called the \ul{\em $3$-triangulation} of $\frak{S}$ associated to the ideal triangulation $\Delta$.
\end{definition}
Although the construction of $Q_\Delta$ made use of the surface $\frak{S}$, the resulting quiver $Q_\Delta$ can be considered as an abstract quiver. One property of these quivers is that when we flip an ideal triangulation to another one, the corresponding quivers are related by a sequence of quiver mutations.

\begin{figure}[htbp!]
\hspace*{-5mm}
\raisebox{-0.5\height}{\scalebox{0.7}{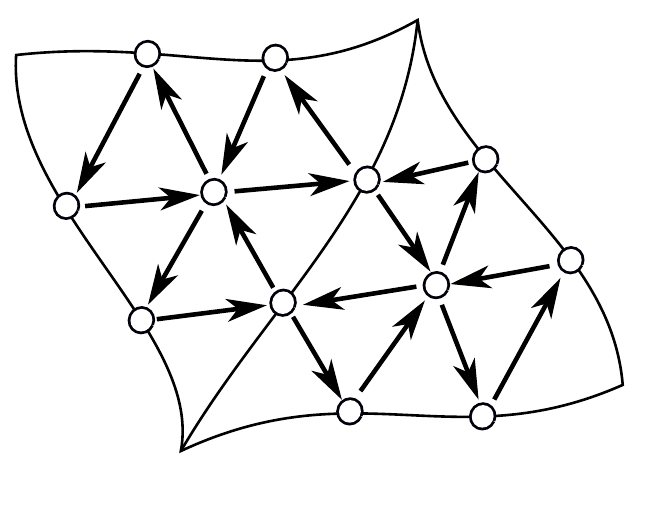}} \hspace{-4mm} $\overset{\mu_{\redfix{v}_3} \mu_{\redfix{v}_4}}{\longrightarrow}$ \hspace{-1mm}\raisebox{-0.5\height}{\scalebox{0.7}{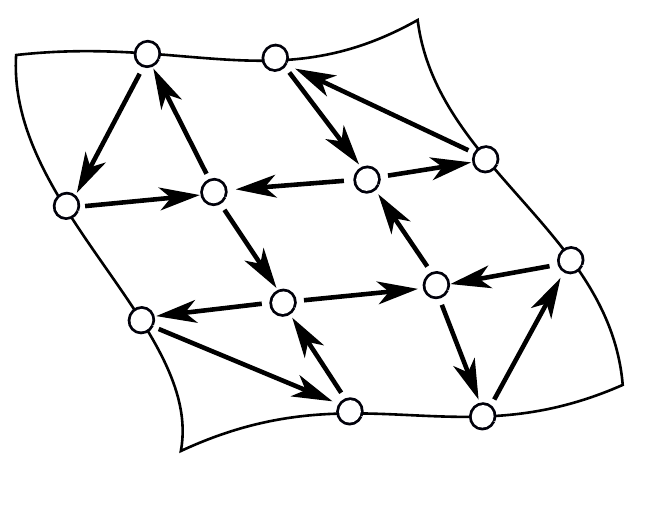}} \hspace{-4mm} $\overset{\mu_{\redfix{v}_7} \mu_{\redfix{v}_{12}}}{\longrightarrow}$  \hspace{-2mm} \raisebox{-0.5\height}{\scalebox{0.7}{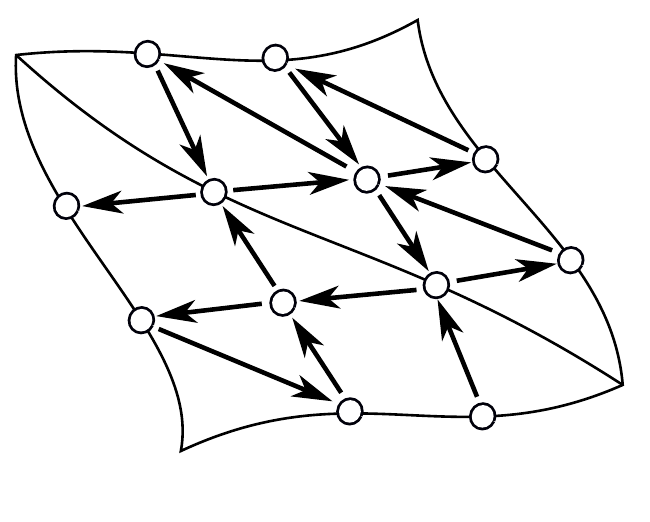}}
\hspace{-2mm}
\vspace{-2mm}
\caption{Sequence of mutations for a flip at an edge, transforming $Q_\Delta$ to $Q_{\Delta'}$}
\vspace{-2mm}
\label{fig:mutations_for_a_flip}
\end{figure}

\begin{lemma}[flip as four quiver mutations; \cite{FG06}]
\label{lem:flip_as_four_quiver_mutations}
Let $\Delta$ and $\Delta'$ be ideal triangulations of a triangulable generalized marked surface $\frak{S}$ that are related by the flip at an edge labeled by $e$. Let $k_1,k_2$ be the two nodes of $Q_\Delta$ lying on the edge $e$ of $\Delta$, and $k_3,k_4$ be the nodes of $Q_\Delta$ lying in the interiors of the two ideal triangles of $\Delta$ having $e$ as one of their sides. Then \redfix{$\mu_{k_1}$ commutes with $\mu_{k_2}$, and $\mu_{k_3}$ commutes with $\mu_{k_4}$. The quivers $Q_\Delta$ and $Q_{\Delta'}$ are related as follows:}
\begin{align}
\label{eq:quiver_mutation_sequence}
Q_{\Delta'} = \redfix{\mu_{k_4} \mu_{k_3} \mu_{k_2} \mu_{k_1}} (Q_\Delta).
\end{align}
\end{lemma}
This lemma is straightforward to check, and is partly depicted in Fig.\ref{fig:mutations_for_a_flip}.

\vs

We now state Fock-Goncharov's special atlases.
\begin{proposition}[Cluster atlases for Fock-Goncharov moduli spaces; \cite{FG06}]
\label{prop:cluster_atlases_of_FG}
Let $\frak{S}$ be a triangulable punctured surface. 

\vs

$\bullet$ For each ideal triangulation $\Delta$ of $\frak{S}$, there exists a cluster $\mathscr{A}$-chart $(Q_\Delta,\psi_\Delta)$ of $\mathscr{A}_{{\rm SL}_3,\frak{S}}$, such that these charts are contained a cluster $\mathscr{A}$-atlas, so that if two triangulations $\Delta,\Delta'$ are related by a flip, then the corresponding cluster $\mathscr{A}$-charts are related by the sequence $\redfix{\mu_{k_4} \mu_{k_3} \mu_{k_2} \mu_{k_1}}$ of four cluster $\mathscr{A}$-mutations as appearing in eq.\eqref{eq:quiver_mutation_sequence}.

\vs

$\bullet$ For each ideal triangulation $\Delta$ of $\frak{S}$, there exists a cluster $\mathscr{X}$-chart $(Q_\Delta,\psi_\Delta)$ of $\mathscr{X}_{{\rm PGL}_3,\frak{S}}$, such that these charts are contained a cluster $\mathscr{X}$-atlas, so that if two triangulations $\Delta,\Delta'$ are related by a flip, then the corresponding cluster $\mathscr{X}$-charts are related by the sequence $\redfix{\mu_{k_4} \mu_{k_3} \mu_{k_2} \mu_{k_1}}$ of four cluster $\mathscr{X}$-mutations as appearing in eq.\eqref{eq:quiver_mutation_sequence}.
\end{proposition}
For explicit construction of these charts, see \cite{FG06} \cite{Douglas}. What we will do recall later is the reconstruction map for the above cluster $\mathscr{X}$-charts of $\mathscr{X}_{{\rm PGL}_3,\frak{S}}$; namely, given the cluster $\mathscr{X}$-coordinates for an ideal triangulation, we will reconstruct a framed ${\rm PGL}_3$-local system on $\frak{S}$.

\vs

Before going to the next section, we give a couple more remarks on a stack $\mathscr{S}$ equipped with a cluster $\mathscr{A}$- or a cluster $\mathscr{X}$-atlas. When one focuses only on the coordinate change formulas, this cluster structure on $\mathscr{S}$ is completely determined by the quiver $Q$ for any of the cluster chart chosen. As a matter of fact, Given any abstract quiver $Q$ with $N$ nodes, or any skew-symmetric exchange matrix $\varepsilon = \varepsilon_Q$, one may start from the affine scheme $(\mathbb{G}_m)^N$ associated to $Q$, package it as the data of a \ul{\em seed torus} $(Q,(\mathbb{G}_m)^N)$, then by mutating at a node $k$ construct another seed torus $(\mu_k(Q), (\mathbb{G}_m)^N)$ which is glued to the original seed along the cluster $\mathscr{A}$- or cluster $\mathscr{X}$-mutation map. Staring from one seed torus, one can mutate in $N$ directions to get $N$ seed tori to glue, then mutate at nodes of these new seed tori, etc. Gluing all such seed tori, one obtains the so-called \ul{\em cluster $\mathscr{A}$-variety} $\mathscr{A}_Q = \mathscr{A}_{|Q|}$, and the \ul{\em cluster $\mathscr{X}$-variety} $\mathscr{X}_Q = \mathscr{X}_{|Q|}$, where $|Q|$ means the mutation equivalence class of a quiver $Q$. These cluster varieties are schemes defined abstractly and combinatorially, and having a cluster atlas of a stack $\mathscr{S}$ provides a birational map from $\mathscr{S}$ to the corresponding cluster variety. Many properties and questions on $\mathscr{S}$ related to the chosen cluster structure can be translated to those on the corresponding abstract cluster varieties. For example, $\mathscr{O}_{\rm cl}(\mathscr{S})$ is isomorphic to the ring of all regular functions on the corresponding cluster variety. 

\vspace{1mm}

Observe that cluster mutation coordinate change formulas involve only multiplication, division, and addition, but not subtraction. This allows us to consider the set of points valued in a {\em semi-field}, not just in a field. A semi-field $(\mathbb{P},\oplus, \odot)$ means a set $\mathbb{P}$ equipped with two binary operations $\oplus$ and $\odot$, where $(\mathbb{P},\odot)$ is an abelian group, where $\oplus$ is required only to be associative and commutative, so that $\oplus$ and $\odot$ satisfy the distributive law. Of our particular interest are two examples of semi-fields:
\begin{align*}
& \mathbb{R}_{>0} ~:~ \mbox{semi-field of positive-real numbers, with usual addition and usual multiplication of real numbers} \\
& \mathbb{Z}^{\redfix{T}} ~:~ \mbox{semi-field of \ul{\em tropical integers}, where $\mathbb{Z}^{\redfix{T}}=\mathbb{Z}$ as a set, and $a\oplus b := \max(a,b)$ and $a\odot b := a+b$.}
\end{align*}
\redfix{In the previous versions of the present paper, $\mathbb{Z}^T$ was denoted by $\mathbb{Z}^t$, which is fixed now as suggested by Linhui Shen.} 
What will play important roles are $\mathscr{A}_{{\rm SL}_3,\frak{S}}(\mathbb{Z}^{\redfix{T}})$ and $\mathscr{X}_{{\rm PGL}_3,\frak{S}}(\mathbb{R}_{>0})$. We note that, in general, for a stack $\mathscr{S}$ equipped with a cluster atlas, or for a cluster variety $\mathscr{S}$, and for a semi-field $\mathbb{P}$, the set $\mathscr{S}(\mathbb{P})$ can be understood as being obtained by gluing the sets $\mathbb{P}^N$ along the {\em tropicalized} version of the cluster mutation formulas, namely by replacing the operations $+,\cdot, \div$ by $\oplus, \odot, \circled{\mbox{$\div$}}$ (where $\circled{\mbox{$\div$}}$ is the inverse operation of the tropical multiplication $\odot$). Note also that, unlike the general case of fields, these gluing maps between $\mathbb{P}^N$ are bijections, so that $\mathscr{S}(\mathbb{P})$ is in fact $\mathbb{P}^N$ as a set.

\vspace{0mm}

\section{${\rm SL}_3$-webs and ${\rm SL}_3$-laminations on a surface}

\vspace{-1mm}

\subsection{${\rm SL}_3$-webs and ${\rm SL}_3$-skein algebras}

A {\em web} on a surface is a certain oriented graph on the surface $\frak{S}$. Since Kuperberg \cite{Kuperberg} introduced it for the case when the surface is a disc, for a diagrammatic calculus for representation theory of the (quantized) Lie algebra of rank up to two, it has been extensively studied by many authors, being generalized in several directions. In particular, the $A_2$-type webs are generalized to corresponding objects living in a thickened surface $\frak{S} \times (-1,1)$, leading to the definition of a\redfix{n} `${\rm SU}_3$-skein algebra' \cite{S05} \cite{FS}. We start from the following definition, taken from \cite{Kuperberg} \cite{S05} \cite{FS} and modified to fit our purpose. 
\begin{definition}
\label{def:A2-web}
Let $(\Sigma,\mathcal{P})$ be a generalized marked surface, and let $\frak{S} = \Sigma\setminus\mathcal{P}$. An \ul{\em ${\rm SL}_3$-web} $W$ in $(\Sigma,\mathcal{P})$ (or in $\frak{S}$) consists of

$\bullet$ a finite subset of $\partial \frak{S}$, whose elements are called \ul{\em external vertices}; when there is no confusion, we refer to them as \ul{\em endpoints} of $W$, and let $\partial W$ be the set of all endpoints of $W$;

$\bullet$ a finite subset of $\mathring{\frak{S}}$, whose elements are called \ul{\em internal vertices};

$\bullet$ a finite set of non-closed oriented smooth curves in $\redfix{\frak{S}}$ ending at external or internal vertices, whose elements are called \ul{\em (oriented) edges} of $W$;

$\bullet$ a finite set of closed oriented smooth curves in $\mathring{\frak{S}}$, whose elements are called \ul{\em (oriented) loops} of $W$,

\vs

subject to the following conditions\redfix{, where we often let the symbol $W$ to denote the subset of $\frak{S}$ given by the union of all edges and loops of $W$}:
\begin{enumerate}
\itemsep0em
\item[\rm (W1)] each external vertex $v$ is $1$-valent, i.e. exactly one edge of $W$ ends at $v$, \redfix{and $W$ meets a boundary arc transversally at an external vertex (if $W$ has endpoints)};

\item[\rm (W2)] each internal vertex $v$ is either a $3$-valent sink or a $3$-valent source, i.e. exactly three edges of $W$ end at $v$, and the orientations of them are either all toward $v$ or all outgoing from $v$\redfix{;}

\item[\rm (W3)] each self-intersection of \redfix{$W$ that is not an internal vertex} is a transverse double intersection lying \redfix{$\mathring{\frak{S}}$}\redfix{, and is called a \ul{\em crossing} of $W$};

\item[\rm (W4)] there are at most finitely many crossings.
\end{enumerate}
\end{definition}

We depict the external and the internal vertices of $W$ by bullets $\bullet$. We allow the empty ${\rm SL}_3$-web ${\O}$. 

\begin{definition}
\label{def:A2-skein_algebra}
Let $\frak{S}$ be a generalized marked surface.

$\bullet$ An \ul{\em isotopy} of ${\rm SL}_3$-webs in $\frak{S}$ is an isotopy within the class of ${\rm SL}_3$-webs in $\frak{S}$. Two ${\rm SL}_3$-webs in $\frak{S}$ are said to be \ul{\em equivalent} if they are related by a sequence of isotopies and the following moves:
\begin{enumerate}
\itemsep0em
\item[\rm (M1)] Reidemeister moves I \raisebox{-0.4\height}{\scalebox{0.8}{
\begingroup%
  \makeatletter%
  \providecommand\color[2][]{%
    \errmessage{(Inkscape) Color is used for the text in Inkscape, but the package 'color.sty' is not loaded}%
    \renewcommand\color[2][]{}%
  }%
  \providecommand\transparent[1]{%
    \errmessage{(Inkscape) Transparency is used (non-zero) for the text in Inkscape, but the package 'transparent.sty' is not loaded}%
    \renewcommand\transparent[1]{}%
  }%
  \providecommand\rotatebox[2]{#2}%
  \newcommand*\fsize{\dimexpr\f@size pt\relax}%
  \newcommand*\lineheight[1]{\fontsize{\fsize}{#1\fsize}\selectfont}%
  \ifx\svgwidth\undefined%
    \setlength{\unitlength}{31.18110236bp}%
    \ifx\svgscale\undefined%
      \relax%
    \else%
      \setlength{\unitlength}{\unitlength * \real{\svgscale}}%
    \fi%
  \else%
    \setlength{\unitlength}{\svgwidth}%
  \fi%
  \global\let\svgwidth\undefined%
  \global\let\svgscale\undefined%
  \makeatother%
  \begin{picture}(1,0.90909091)%
    \lineheight{1}%
    \setlength\tabcolsep{0pt}%
    \put(0,0){\includegraphics[width=\unitlength,page=1]{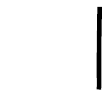}}%
    \put(0.48908514,0.44621243){\color[rgb]{0,0,0}\makebox(0,0)[lt]{\lineheight{1.25}\smash{\begin{tabular}[t]{l}$\leftrightarrow$\end{tabular}}}}%
    \put(0,0){\includegraphics[width=\unitlength,page=2]{RM1.pdf}}%
  \end{picture}%
\endgroup%
}}, II \raisebox{-0.4\height}{\scalebox{0.8}{
\begingroup%
  \makeatletter%
  \providecommand\color[2][]{%
    \errmessage{(Inkscape) Color is used for the text in Inkscape, but the package 'color.sty' is not loaded}%
    \renewcommand\color[2][]{}%
  }%
  \providecommand\transparent[1]{%
    \errmessage{(Inkscape) Transparency is used (non-zero) for the text in Inkscape, but the package 'transparent.sty' is not loaded}%
    \renewcommand\transparent[1]{}%
  }%
  \providecommand\rotatebox[2]{#2}%
  \newcommand*\fsize{\dimexpr\f@size pt\relax}%
  \newcommand*\lineheight[1]{\fontsize{\fsize}{#1\fsize}\selectfont}%
  \ifx\svgwidth\undefined%
    \setlength{\unitlength}{34.01574803bp}%
    \ifx\svgscale\undefined%
      \relax%
    \else%
      \setlength{\unitlength}{\unitlength * \real{\svgscale}}%
    \fi%
  \else%
    \setlength{\unitlength}{\svgwidth}%
  \fi%
  \global\let\svgwidth\undefined%
  \global\let\svgscale\undefined%
  \makeatother%
  \begin{picture}(1,0.83333333)%
    \lineheight{1}%
    \setlength\tabcolsep{0pt}%
    \put(0,0){\includegraphics[width=\unitlength,page=1]{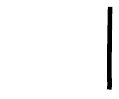}}%
    \put(0.40423083,0.40902805){\color[rgb]{0,0,0}\makebox(0,0)[lt]{\lineheight{1.25}\smash{\begin{tabular}[t]{l}$\leftrightarrow$\end{tabular}}}}%
    \put(0,0){\includegraphics[width=\unitlength,page=2]{RM2.pdf}}%
  \end{picture}%
\endgroup%
}} and III \raisebox{-0.4\height}{\scalebox{0.8}{
\begingroup%
  \makeatletter%
  \providecommand\color[2][]{%
    \errmessage{(Inkscape) Color is used for the text in Inkscape, but the package 'color.sty' is not loaded}%
    \renewcommand\color[2][]{}%
  }%
  \providecommand\transparent[1]{%
    \errmessage{(Inkscape) Transparency is used (non-zero) for the text in Inkscape, but the package 'transparent.sty' is not loaded}%
    \renewcommand\transparent[1]{}%
  }%
  \providecommand\rotatebox[2]{#2}%
  \newcommand*\fsize{\dimexpr\f@size pt\relax}%
  \newcommand*\lineheight[1]{\fontsize{\fsize}{#1\fsize}\selectfont}%
  \ifx\svgwidth\undefined%
    \setlength{\unitlength}{56.69291339bp}%
    \ifx\svgscale\undefined%
      \relax%
    \else%
      \setlength{\unitlength}{\unitlength * \real{\svgscale}}%
    \fi%
  \else%
    \setlength{\unitlength}{\svgwidth}%
  \fi%
  \global\let\svgwidth\undefined%
  \global\let\svgscale\undefined%
  \makeatother%
  \begin{picture}(1,0.5)%
    \lineheight{1}%
    \setlength\tabcolsep{0pt}%
    \put(0.4012885,0.24541683){\color[rgb]{0,0,0}\makebox(0,0)[lt]{\lineheight{1.25}\smash{\begin{tabular}[t]{l}$\leftrightarrow$\end{tabular}}}}%
    \put(0,0){\includegraphics[width=\unitlength,page=1]{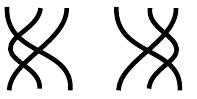}}%
  \end{picture}%
\endgroup%
}}, with all possible orientations;

\item[\rm (M2)] the web Reidemeister move \raisebox{-0.4\height}{\scalebox{0.9}{
\begingroup%
  \makeatletter%
  \providecommand\color[2][]{%
    \errmessage{(Inkscape) Color is used for the text in Inkscape, but the package 'color.sty' is not loaded}%
    \renewcommand\color[2][]{}%
  }%
  \providecommand\transparent[1]{%
    \errmessage{(Inkscape) Transparency is used (non-zero) for the text in Inkscape, but the package 'transparent.sty' is not loaded}%
    \renewcommand\transparent[1]{}%
  }%
  \providecommand\rotatebox[2]{#2}%
  \newcommand*\fsize{\dimexpr\f@size pt\relax}%
  \newcommand*\lineheight[1]{\fontsize{\fsize}{#1\fsize}\selectfont}%
  \ifx\svgwidth\undefined%
    \setlength{\unitlength}{51.02362205bp}%
    \ifx\svgscale\undefined%
      \relax%
    \else%
      \setlength{\unitlength}{\unitlength * \real{\svgscale}}%
    \fi%
  \else%
    \setlength{\unitlength}{\svgwidth}%
  \fi%
  \global\let\svgwidth\undefined%
  \global\let\svgscale\undefined%
  \makeatother%
  \begin{picture}(1,0.38888889)%
    \lineheight{1}%
    \setlength\tabcolsep{0pt}%
    \put(0.41647795,0.19421307){\color[rgb]{0,0,0}\makebox(0,0)[lt]{\lineheight{1.25}\smash{\begin{tabular}[t]{l}$\leftrightarrow$\end{tabular}}}}%
    \put(0,0){\includegraphics[width=\unitlength,page=1]{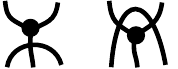}}%
  \end{picture}%
\endgroup%
}}, with all possible orientations (according to Def.\ref{def:A2-web});

\item[\rm (M3)] the boundary exchange move \raisebox{-0.4\height}{\scalebox{0.7}{
\begingroup%
  \makeatletter%
  \providecommand\color[2][]{%
    \errmessage{(Inkscape) Color is used for the text in Inkscape, but the package 'color.sty' is not loaded}%
    \renewcommand\color[2][]{}%
  }%
  \providecommand\transparent[1]{%
    \errmessage{(Inkscape) Transparency is used (non-zero) for the text in Inkscape, but the package 'transparent.sty' is not loaded}%
    \renewcommand\transparent[1]{}%
  }%
  \providecommand\rotatebox[2]{#2}%
  \newcommand*\fsize{\dimexpr\f@size pt\relax}%
  \newcommand*\lineheight[1]{\fontsize{\fsize}{#1\fsize}\selectfont}%
  \ifx\svgwidth\undefined%
    \setlength{\unitlength}{42.51968504bp}%
    \ifx\svgscale\undefined%
      \relax%
    \else%
      \setlength{\unitlength}{\unitlength * \real{\svgscale}}%
    \fi%
  \else%
    \setlength{\unitlength}{\svgwidth}%
  \fi%
  \global\let\svgwidth\undefined%
  \global\let\svgscale\undefined%
  \makeatother%
  \begin{picture}(1,0.7)%
    \lineheight{1}%
    \setlength\tabcolsep{0pt}%
    \put(0,0){\includegraphics[width=\unitlength,page=1]{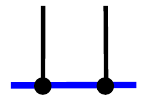}}%
  \end{picture}%
\endgroup%
}} $\leftrightarrow$ \raisebox{-0.4\height}{\scalebox{0.7}{
\begingroup%
  \makeatletter%
  \providecommand\color[2][]{%
    \errmessage{(Inkscape) Color is used for the text in Inkscape, but the package 'color.sty' is not loaded}%
    \renewcommand\color[2][]{}%
  }%
  \providecommand\transparent[1]{%
    \errmessage{(Inkscape) Transparency is used (non-zero) for the text in Inkscape, but the package 'transparent.sty' is not loaded}%
    \renewcommand\transparent[1]{}%
  }%
  \providecommand\rotatebox[2]{#2}%
  \newcommand*\fsize{\dimexpr\f@size pt\relax}%
  \newcommand*\lineheight[1]{\fontsize{\fsize}{#1\fsize}\selectfont}%
  \ifx\svgwidth\undefined%
    \setlength{\unitlength}{42.51968504bp}%
    \ifx\svgscale\undefined%
      \relax%
    \else%
      \setlength{\unitlength}{\unitlength * \real{\svgscale}}%
    \fi%
  \else%
    \setlength{\unitlength}{\svgwidth}%
  \fi%
  \global\let\svgwidth\undefined%
  \global\let\svgscale\undefined%
  \makeatother%
  \begin{picture}(1,0.7)%
    \lineheight{1}%
    \setlength\tabcolsep{0pt}%
    \put(0,0){\includegraphics[width=\unitlength,page=1]{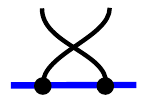}}%
  \end{picture}%
\endgroup%
}}, with all possible orientations, where the horizontal blue line is boundary.
\end{enumerate}

$\bullet$ Let $\mathcal{R}$ be a commutative ring with unity $1$. The \ul{\em (commutative) ${\rm SL}_3$-skein algebra} $\mathcal{S}(\frak{S};\mathcal{R})$ is the free $\mathcal{R}$-module with the set of all equivalence classes of ${\rm SL}_3$-webs in $\frak{S}$ as a free basis, mod out by the \ul{\em ${\rm SL}_3$-skein relations} (S1), (S2), (S3) and (S4) in Fig.\ref{fig:A2-skein_relations}.

\vs

$\bullet$ For an equivalence class of ${\rm SL}_3$-webs $W$ in $\frak{S}$, the corresponding element of the ${\rm SL}_3$-skein algebra $\mathcal{S}(\frak{S};\mathcal{R})$ is denoted by $[W]$ and is called an \ul{\em ${\rm SL}_3$-skein}.
\end{definition}
Note that (M1) and (M2) make (S2) and (S3) redundant. The following special class of ${\rm SL}_3$-webs are important.
\begin{definition}[\cite{Kuperberg} \cite{SW} \cite{FS}]
\label{def:non-elliptic}
Let $\frak{S}$ be a generalized marked surface. 

$\bullet$ An ${\rm SL}_3$-web $W$ in $\frak{S}$ is said to be \ul{\em non-elliptic} if all of the following hold:
\begin{enumerate}
\itemsep0em
\item[\rm (NE1)] $W$ has no crossings;

\item[\rm (NE2)] none of the loops of $W$ is a contractible loop in $\frak{S}$;

\item[\rm (NE3)] none of the components of the complement in $\frak{S}$ of \redfix{$W$} is a contractible region bounded by either two or four edges of $W$ (as appearing in the first term of (S2) or (S3)).
\end{enumerate}

$\bullet$ A non-elliptic ${\rm SL}_3$-web $W$ is \ul{\em weakly reduced} if it contains none of \raisebox{-0.3\height}{\scalebox{0.7}{
\begingroup%
  \makeatletter%
  \providecommand\color[2][]{%
    \errmessage{(Inkscape) Color is used for the text in Inkscape, but the package 'color.sty' is not loaded}%
    \renewcommand\color[2][]{}%
  }%
  \providecommand\transparent[1]{%
    \errmessage{(Inkscape) Transparency is used (non-zero) for the text in Inkscape, but the package 'transparent.sty' is not loaded}%
    \renewcommand\transparent[1]{}%
  }%
  \providecommand\rotatebox[2]{#2}%
  \newcommand*\fsize{\dimexpr\f@size pt\relax}%
  \newcommand*\lineheight[1]{\fontsize{\fsize}{#1\fsize}\selectfont}%
  \ifx\svgwidth\undefined%
    \setlength{\unitlength}{42.51968504bp}%
    \ifx\svgscale\undefined%
      \relax%
    \else%
      \setlength{\unitlength}{\unitlength * \real{\svgscale}}%
    \fi%
  \else%
    \setlength{\unitlength}{\svgwidth}%
  \fi%
  \global\let\svgwidth\undefined%
  \global\let\svgscale\undefined%
  \makeatother%
  \begin{picture}(1,0.56666667)%
    \lineheight{1}%
    \setlength\tabcolsep{0pt}%
    \put(0,0){\includegraphics[width=\unitlength,page=1]{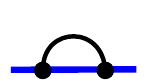}}%
  \end{picture}%
\endgroup%
}} and \raisebox{-0.3\height}{\scalebox{0.7}{
\begingroup%
  \makeatletter%
  \providecommand\color[2][]{%
    \errmessage{(Inkscape) Color is used for the text in Inkscape, but the package 'color.sty' is not loaded}%
    \renewcommand\color[2][]{}%
  }%
  \providecommand\transparent[1]{%
    \errmessage{(Inkscape) Transparency is used (non-zero) for the text in Inkscape, but the package 'transparent.sty' is not loaded}%
    \renewcommand\transparent[1]{}%
  }%
  \providecommand\rotatebox[2]{#2}%
  \newcommand*\fsize{\dimexpr\f@size pt\relax}%
  \newcommand*\lineheight[1]{\fontsize{\fsize}{#1\fsize}\selectfont}%
  \ifx\svgwidth\undefined%
    \setlength{\unitlength}{42.51968504bp}%
    \ifx\svgscale\undefined%
      \relax%
    \else%
      \setlength{\unitlength}{\unitlength * \real{\svgscale}}%
    \fi%
  \else%
    \setlength{\unitlength}{\svgwidth}%
  \fi%
  \global\let\svgwidth\undefined%
  \global\let\svgscale\undefined%
  \makeatother%
  \begin{picture}(1,0.56666667)%
    \lineheight{1}%
    \setlength\tabcolsep{0pt}%
    \put(0,0){\includegraphics[width=\unitlength,page=1]{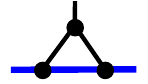}}%
  \end{picture}%
\endgroup%
}}, and is \ul{\em reduced} if furthermore it contains none of \raisebox{-0.3\height}{\scalebox{0.7}{
\begingroup%
  \makeatletter%
  \providecommand\color[2][]{%
    \errmessage{(Inkscape) Color is used for the text in Inkscape, but the package 'color.sty' is not loaded}%
    \renewcommand\color[2][]{}%
  }%
  \providecommand\transparent[1]{%
    \errmessage{(Inkscape) Transparency is used (non-zero) for the text in Inkscape, but the package 'transparent.sty' is not loaded}%
    \renewcommand\transparent[1]{}%
  }%
  \providecommand\rotatebox[2]{#2}%
  \newcommand*\fsize{\dimexpr\f@size pt\relax}%
  \newcommand*\lineheight[1]{\fontsize{\fsize}{#1\fsize}\selectfont}%
  \ifx\svgwidth\undefined%
    \setlength{\unitlength}{42.51968504bp}%
    \ifx\svgscale\undefined%
      \relax%
    \else%
      \setlength{\unitlength}{\unitlength * \real{\svgscale}}%
    \fi%
  \else%
    \setlength{\unitlength}{\svgwidth}%
  \fi%
  \global\let\svgwidth\undefined%
  \global\let\svgscale\undefined%
  \makeatother%
  \begin{picture}(1,0.56666667)%
    \lineheight{1}%
    \setlength\tabcolsep{0pt}%
    \put(0,0){\includegraphics[width=\unitlength,page=1]{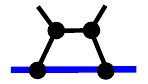}}%
  \end{picture}%
\endgroup%
}}; in these pictures, the blue line is boundary, the edges can be given all possible orientations (according to Def.\ref{def:A2-web}), and the boundary 2-gon, 3-gon and 4-gon are contractible.
\end{definition}

\begin{remark}
\label{rmk:punctured_surface_non-elliptic}
For a punctured surface $\frak{S}$, any non-elliptic ${\rm SL}_3$-web is reduced.
\end{remark}

\begin{proposition}[\cite{SW} {\cite[Thm.2]{FS}}]
\label{prop:basis_of_A2-skein_algebra}
Let $\frak{S}$ be a punctured surface. The set of all ${\rm SL}_3$-skeins for non-elliptic ${\rm SL}_3$-webs form a basis of $\mathcal{S}(\frak{S};\mathcal{R})$.
\end{proposition}
This can be viewed as an $A_2$ analog of the so-called bangles basis of the usual (\redfix{Kauffman} bracket) skein algebra (of type $A_1$); we will elaborate in the next section \redfix{(\S\ref{sec:regular_functions_on_moduli_spaces})}.

\subsection{${\rm SL}_3$-laminations}

For $A_1$-type theory, Fock and Goncharov \cite{FG06} introduced certain versions of {\em laminations} on a surface $\frak{S}$, where a lamination is defined as a collection of mutually non-intersecting simple unoriented curves equipped with weights, where a constituent curve is either closed or ends at components of $\partial \frak{S}$, and weights are rational or integer numbers. For $A_2$-type theory, here based on reduced non-elliptic ${\rm SL}_3$-webs we propose a generalization of Fock-Goncharov's integral $\mathcal{A}$-laminations which were based on just simple curves, or $A_1$-type webs. Basic idea is to consider reduced non-elliptic ${\rm SL}_3$-webs with weights given on its (connected) components. The weights are required to be non-negative integers, except for the {\em special} curves as in \cite{FG07} for the $A_1$-type webs, which we call {\em peripheral}.
\begin{definition} 
\label{def:peripheral}
Let $(\Sigma,\mathcal{P})$ be a generalized marked surface, and $\frak{S} := \Sigma\setminus\mathcal{P}$.

\vs

$\bullet$ A simple loop in $\frak{S}$ is called  a \ul{\em peripheral loop} if it bounds a region \redfix{diffeo}morphic to a disc with one puncture in the interior. If the corresponding puncture is $p\in \mathcal{P}$, we say that this peripheral loop \ul{\em surrounds} $p$.

\vs

$\bullet$ A \ul{\em peripheral arc} in $\frak{S}$ is a simple curve $e$ in $\frak{S}$ that ends at points of $\partial \frak{S}$ and that bounds a region \redfix{diffeo}morphic to an upper-half disc with one puncture on the boundary, i.e. $e$ is homotopic in $\Sigma$ (rel endpoints) to a simple arc $e'$ lying in $\partial \Sigma$ such that $e'$ contains exactly one point of $\mathcal{P}$, say $p$; in this case, we say this peripheral arc \ul{\em surrounds} $p$.

\vs

$\bullet$ Peripheral loops and peripheral arcs are referred to as \ul{\em peripheral curves}.
\end{definition}

\begin{definition}
Let $W$ be an ${\rm SL}_3$-web in a generalized marked surface $\frak{S}$. We define \ul{\em components} of $W$ as follows. First, each loop of $W$ is a component of $W$, and each edge of $W$ whose two endpoints are both external vertices of $W$ is a component. A union of a collection $C$ of at least two edges of $W$ is called a component if
\begin{enumerate}
\itemsep0em
\item[\rm (C1)] for any two distinct edges $e$ and $e'$ of $C$, there is a sequence of edges $e_1,\ldots,e_n$ \redfix{of $W$} such that $e_1=e$, $e_n=e'$, and $e_i$ meets $e_{i+1}$ at an internal vertex of $W$ for each $i=1,\ldots,n-1$;

\item[\rm (C2)] $C$ is maximal among the collections satisfying (C1).
\end{enumerate}
\end{definition}
\begin{remark}
Each component of $W$ is an ${\rm SL}_3$-web on its own.
\end{remark}
The following is the first main definition of the present paper.
\begin{definition}[${\rm SL}_3$-laminations]
\label{def:A2-lamination}
Let $\frak{S}$ be a generalized marked surface.

\vs

$\bullet$ An \ul{\em (integral) (bounded) ${\rm SL}_3$-lamination} $\ell$ in $\frak{S}$ consists of the equivalence class of a reduced non-elliptic ${\rm SL}_3$-web $W=W(\ell)$ in $\frak{S}$ and the assignment of an integer weight to each component of $W(\ell)$, subject to the following conditions and equivalence relation:
\begin{enumerate}
\itemsep0em
\item[\rm (L1)] the weight for each component of $W(\ell)$ containing an internal vertex is $1$;

\item[\rm (L2)] the weight for each component of $W(\ell)$ that is not a peripheral curve is non-negative;

\item[\rm (L3)] an ${\rm SL}_3$-lamination containing a component of weight zero is equivalent to the ${\rm SL}_3$-lamination with this component removed;

\item[\rm (L4)] an ${\rm SL}_3$-lamination with two of its components being homotopic with weights $a$ and $b$ is equivalent to the ${\rm SL}_3$-lamination with one of these components removed and the other having weight $a+b$.
\end{enumerate}

Let $\mathscr{A}_{\rm L}(\frak{S};\mathbb{Z})$ be the set of all (integral) ${\rm SL}_3$-laminations in $\frak{S}$.

\vs

Let $\mathscr{A}_{\rm L}^0(\frak{S};\mathbb{Z})$ be the set of all (integral) ${\rm SL}_3$-laminations in $\frak{S}$ with no negative weights.
\end{definition}
\begin{lemma}
\label{lem:non-negative_A2-laminations_and_non-elliptic_A2-webs}
An ${\rm SL}_3$-lamination in $\mathscr{A}_{\rm L}^0(\frak{S};\mathbb{Z})$ can be represented by a \redfix{reduced non-elliptic} ${\rm SL}_3$-\redfix{web} whose weights are all $1$. This gives a bijection
$$
\mathscr{A}_{\rm L}^0(\frak{S};\mathbb{Z}) \leftrightarrow \{\mbox{\redfix{equivalence classes of} reduced non-elliptic ${\rm SL}_3$-webs in $\frak{S}$}\}. \qed
$$
\end{lemma}

Crucial in the study of ${\rm SL}_3$-webs and ${\rm SL}_3$-laminations is a coordinate system for them. The coordinate system which we will construct requires the choice of an ideal triangulation $\Delta$ of the surface $\frak{S}$. We first isotope an ${\rm SL}_3$-web to be in a minimal position with respect to $\Delta$.
\begin{definition}[\cite{FS}]
Let $\frak{S}$ be a generalized marked surface, and let $\Delta$ be a collection of mutually non-intersecting ideal arcs in $\frak{S}$. A non-elliptic ${\rm SL}_3$-web $W$ in a triangulable punctured surface $\frak{S}$ is said to be in a \ul{\em minimal position} with respect to $\Delta$ of $\frak{S}$ if the cardinality of the intersection $W\cap \Delta$ equals the minimum of the cardinality of $W'\cap \Delta$ among all non-elliptic ${\rm SL}_3$-webs $W'$ in $\frak{S}$ isotopic to $W$.
\end{definition}
For our convenience, we may assume that a non-elliptic ${\rm SL}_3$-web in a minimal position with respect to \redfix{an ideal triangulation} $\Delta$ meets edges of $\Delta$ transversally. In fact, putting into a minimal position with respect to $\Delta$ is not sufficient for the purpose of constructing our coordinates, and we need a further tidying-up process; we use a result obtained in \cite{FS}. In the end, we would like our ${\rm SL}_3$-web in each \redfix{ideal} triangle $t$ of $\Delta$ to be a disjoint union of elementary pieces; namely, peripheral arcs of $t$ \redfix{which we call} {\em corner arcs} in $t$, \redfix{and} special webs having internal vertices called {\em pyramids} $H_d$ for $d \in \mathbb{Z}\setminus\{0\}$, some examples of which are depicted in Fig.\ref{fig:pyramids}. In particular, $H_d$ has $|d|$ external vertices on each \redfix{of the three} side\redfix{s} of $t$, and $H_{-d}$ can be obtained from $H_d$ by reversing the orientation of all edges of $H_d$. 
%
%
By looking at these pictures, we believe that the readers can deduce the definition of $H_d$ for each $d\in \mathbb{Z}\setminus\{0\}$; see \cite[\S10]{FS} for a precise recipe for constructing $H_d$, and also see \cite{DS1} where $H_d$ is called a honeycomb-web.  
\begin{definition}[\cite{FS}]
\label{def:canonical_web_in_a_triangle}
Let $t$ be a triangle, viewed as a generalized marked surface diffeomorphic to a closed disc with three marked points on the boundary. 

\vs

$\bullet$ For $d\in \mathbb{Z}\setminus\{0\}$, the ${\rm SL}_3$-web $H_d$ \redfix{in $t$} described above is a \ul{\em degree $d$ pyramid} in $t$.  Let $H_0 = {\O}$.

\vs

$\bullet$ A single-component non-elliptic ${\rm SL}_3$-web in $t$  consisting of an edge connecting the $1$-valent vertices lying in two distinct sides of $t$ (i.e. a peripheral arc in $t$) is called a \ul{\em corner arc} in $t$.

\vs

$\bullet$ An ${\rm SL}_3$-web in $t$ is \ul{\em canonical} if it is a disjoint union of one $H_d$ for some $d\in \mathbb{Z}$ and some number of (possibly none of) corner arcs.
\end{definition}
We note that it is not always possible to isotope a reduced non-elliptic ${\rm SL}_3$-web in a triangul\redfix{able} surface so that it is canonical in each triangle of \redfix{a chosen ideal triangulation} $\Delta$. One remedy is to fatten each edge of $\Delta$ to a {\em biangle}, and push some $3$-valent vertices into the biangles.
\begin{definition}[\cite{BW}]
\label{def:split_ideal_triangulation}
Let $\Delta$ be a\redfix{n ideal} triangulation of a triangulable generalized marked surface $\frak{S} = \Sigma\setminus\mathcal{P}$. For each edge $e$ of $\Delta$, choose a $\mathcal{P}$-arc $e'$ in $\Sigma$ isotopic to $e$ as $\mathcal{P}$-arcs, such that $\wh{\Delta} := \Delta \cup \{e' : e \in \Delta\}$ is a mutually non-intersecting collection of $\mathcal{P}$-arcs. We call $\wh{\Delta}$ a \ul{\em split ideal triangulation} of $\frak{S}$ for the triangulation $\Delta$. The closure (in $\Sigma$ or $\frak{S}$) of a connected component of the complement of (the union of members of) $\wh{\Delta}$ in $\Sigma$ is called a \ul{\em (ideal) triangle} of $\wh{\Delta}$ if it is bounded by three edges of $\wh{\Delta}$, and a \ul{\em (ideal) biangle} of $\wh{\Delta}$ if it is bounded by two edges of $\wh{\Delta}$.
\end{definition}
Edges of $\Delta$ are in bijection with the biangles of $\wh{\Delta}$, and triangles of $\Delta$ are in bijection with the triangles of $\wh{\Delta}$. Each triangle and biangle of $\wh{\Delta}$ may be viewed as a generalized marked surface on its own, in a natural way. In particular, a {\em biangle} $B$ (Def.\ref{def:triangulation}) can be considered as being a generalized marked surface diffeomorphic to a closed disc, with two marked points on the boundary, with no punctures.

\vs

We now recall from \cite{FS} some special classes of ${\rm SL}_3$-webs in a biangle. To match our convention used in a later section, we change some notations and definitions, and re-interpret a little bit.
\begin{definition}[{\cite[\S9]{FS}}]
\label{def:crossbar}
$\bullet$ Let $B$ be a biangle, viewed as a generalized marked surface. For integer $n\ge 1$, consider a mutually disjoint collection $n$ \ul{\em strands}, each of which is a simple curve connecting the two sides of $B$. Choose a possibly-empty mutually disjoint finite collection of \ul{\em crossbars}, each of which is a simple curve connecting two adjacent strands, such that each intersection of crossbars and strands are transverse double intersection in the interior of $B$, and that under a homeomorphism of the biangle $B$ to $\mathbb{R} \times [0,1]$ (the two sides going to $\mathbb{R}\times \{0\}$ and $\mathbb{R}\times \{1\}$) each strand is of the form $\{c\} \times [0,1]$ (i.e. vertical) and each crossbar is of the form $[c_1,c_2] \times \{a\}$ (i.e. horizontal). Union of all $n$ strands together with all of these chosen crossbars is called a \ul{\em crossbar graph} of \ul{\em index} $n$. 

\vs

$\bullet$ A crossing-less ${\rm SL}_3$-web $W$ in $B$ is called a \ul{\em crossbar ${\rm SL}_3$-web} (in $B$) of \ul{\em index} $n$ if the union of all its components, with orientations forgotten, is a crossbar graph of index $n$. 

\vs

$\bullet$ The \ul{\em signature} of a crossbar ${\rm SL}_3$-web $W$ is the map $\partial W \to \{{\rm sink},{\rm source}\}$, which records the sink-source information of $W$ at each external vertex of $W$.

\vs

$\bullet$ A \ul{\em minimal crossbar ${\rm SL}_3$-web} in $B$ (of index $n$) is a crossbar ${\rm SL}_3$-web in $B$ (of index $n$) that is non-elliptic (i.e. does not contain a contractible 4-gon as in (S3) of Fig.\ref{fig:A2-skein_relations}).
\end{definition}
\begin{remark}
Not every crossbar graph admits a crossbar ${\rm SL}_3$-web structure.
\end{remark}
By convention, a crossbar ${\rm SL}_3$-web of index $0$ means the empty ${\rm SL}_3$-web, which is a minimal crossbar ${\rm SL}_3$-web. In the above definition, one observes that a crossbar ${\rm SL}_3$-web is minimal if there are no two `consecutive' crossbars. It is easy to observe that a minimal crossbar ${\rm SL}_3$-web is weakly reduced.

\begin{lemma}[minimal crossbar ${\rm SL}_3$-web is determined by signature; \cite{FS}]
\label{lem:minimal_crossbar_A2-web_determined_by_signature}
In a biangle $B$, one has:
\begin{enumerate}
\itemsep0em
\item[\rm (MC1)] Each minimal crossbar ${\rm SL}_3$-web in $B$ is completely determined, up to isotopy, by its signature.

\item[\rm (MC2)] For any crossbar ${\rm SL}_3$-web \redfix{$W$} in $B$, its signature $\partial W \to \{{\rm sink},{\rm source}\}$ is \ul{\em sign-preserving} in the sense that, the number of sources of $W$ on one side of $B$ coincides with the number of sinks of $W$ on the other side of $B$.

\item[\rm (MC3)] For $n\ge 1$, pick a subset $\mathcal{V}$ of $B$ consisting of $n$ points in the interior of one side of $B$ and $n$ points in the interior of the other side of $B$. Choose any function $\mathcal{V} \to \{{\rm sink},{\rm source}\}$ that is sign-preserving in the above sense. Then there exists a minimal crossbar ${\rm SL}_3$-web $W$ in $B$ with $\partial W = \mathcal{V}$ such that the signature of $W$ coincides with this function.
\end{enumerate}
\end{lemma}

\begin{definition}[\cite{FS}]
\label{def:canonical_wrt_split_ideal_triangulation}
A weakly reduced non-elliptic ${\rm SL}_3$-web $W$ in a triangulable generalized marked surface $\frak{S}$ is said to be \ul{\em canonical} with respect to a split ideal triangulation $\wh{\Delta}$ of $\frak{S}$ if:
\begin{enumerate}
\itemsep0em
\item[\rm (CW1)] for each triangle $\wh{t}$ of $\wh{\Delta}$, the intersection $W\cap \wh{t}$ is a canonical ${\rm SL}_3$-web in $\wh{t}$ (Def.\ref{def:canonical_web_in_a_triangle});

\item[\rm (CW2)] for each biangle $B$ of $\wh{\Delta}$, the intersection $W\cap B$ is a minimal crossbar ${\rm SL}_3$-web in $B$ (Def.\ref{def:crossbar}).
\end{enumerate}

\end{definition}

\begin{lemma}[\cite{FS}]
Let $\wh{\Delta}$ be a split ideal triangulation of a triangulable generalized marked surface $\frak{S}$. Any weakly reduced non-elliptic ${\rm SL}_3$-web $W$ in $\frak{S}$ is isotopic to a weakly reduced non-elliptic ${\rm SL}_3$-web in $\frak{S}$ that is canonical with respect to $\wh{\Delta}$.
\end{lemma}

A useful observation:
\begin{lemma}
Any canonical ${\rm SL}_3$-web in a triangle is a reduced non-elliptic ${\rm SL}_3$-web in that triangle, when the triangle is viewed as a generalized marked surface.
\end{lemma}
\begin{corollary}
\label{cor:intersection_in_wh_t_is_A2-lamination}
Let $\ell$ be an ${\rm SL}_3$-lamination in a triangulable generalized marked surface $\frak{S}$ that is represented by a weighted reduced non-elliptic ${\rm SL}_3$-web $W$ in $\frak{S}$ that is canonical with respect to a split ideal triangulation $\wh{\Delta}$. For each triangle $\wh{t}$ of $\wh{\Delta}$, $\ell \cap \wh{t}$ is an ${\rm SL}_3$-lamination in $\wh{t}$ (represented by weighted ${\rm SL}_3$-web $W\cap \wh{t}$), when $\wh{t}$ is viewed as a generalized marked surface on its own.
\end{corollary}

\subsection{Frohman-Sikora coordinates and Douglas-Sun coordinates for ${\rm SL}_3$-webs}

In the present subsection we recall  two coordinate systems for non-elliptic ${\rm SL}_3$-webs in $\frak{S}$ with respect to an ideal triangulation $\Delta$ of $\frak{S}$, one by Frohman-Sikora \cite{FS} and the other by Douglas-Sun \cite{Douglas} \cite{DS1} \cite{DS2}. We modified the notations to fit our purpose.

\begin{definition}[\cite{FS}]
Let $t$ be a triangle, viewed as a generalized marked surface. Let $e_1,e_2,e_3$ be the sides of $t$, appearing clockwise in this order in $\partial t$. Let $W$ be a canonical ${\rm SL}_3$-web in $t$ (Def.\ref{def:canonical_web_in_a_triangle}). 

\vs

For each side $e_\alpha$, let ${\rm e}_{{\rm out},\alpha}(W)$ be the number of external vertices of $W$ that are sinks and lie on $e_\alpha$, and ${\rm e}_{{\rm in},\alpha}(W)$ be the number of external vertices of $W$ that are sources and lie on $e_\alpha$. These six numbers are called \ul{\em intersection coordinates} of $W$.

\vs

A corner arc of $W$ is said to be \ul{\em left turn} if it starts at a vertex in $e_\alpha$ and terminates at a vertex in $e_{\alpha+1}$ (where $e_4=e_1$), and \ul{\em right turn} otherwise. Let
\begin{align*}
{\rm r}_t(W) & = \mbox{(number of left turn corner arcs of $W$)} - \mbox{(number of right turn corner arcs of $W$)},
\end{align*}
which is called the \ul{\em rotation number} of $W$. 
\end{definition}

\begin{definition}[\cite{FS}]
\label{def:FS_coordinates}
Let $\wh{\Delta}$ be a split ideal triangulation of a triangulable generalized marked surface $\frak{S}$. Let $W$ be a weakly reduced non-elliptic ${\rm SL}_3$-web in $\frak{S}$ that is canonical with respect to $\wh{\Delta}$. For each triangle $\wh{t}$ of $\wh{\Delta}$, consider the intersection coordinates and the rotation number for the ${\rm SL}_3$-web $W\cap \wh{t}$ in the triangle $\wh{t}$. These numbers are the \ul{\em Frohman-Sikora coordinates} of $W$ with respect to $\wh{\Delta}$.
\end{definition}
Each edge $e$ of $\wh{\Delta}$ is a side of a unique triangle of $\wh{\Delta}$, say $\wh{t}$; denote by ${\rm e}_{{\rm out},e}(W)$ and ${\rm e}_{\redfix{\rm in},e}(W)$ the intersection coordinates of $W\cap \wh{t}$ at this side $e$. Write ${\rm r}_{\wh{t}}(W) := {\rm r}_{\wh{t}}(W\cap \wh{t}\,\,)$. If $e$ and $e'$ are edges of $\wh{\Delta}$ forming a biangle, then it is easy to see from Lem.\ref{lem:minimal_crossbar_A2-web_determined_by_signature}(MC2) that ${\rm e}_{{\rm out},e}(W) = {\rm e}_{{\rm in},e'}(W)$ and ${\rm e}_{{\rm in},e}(W) = {\rm e}_{{\rm out},e'}(W)$. So, one can think of the intersection coordinates \redfix{as being} assigned to edges of $\Delta$, instead of edges of $\wh{\Delta}$, and hence one may also say that the Frohman-Sikora coordinates are defined with respect to the ideal triangulation $\Delta$ instead of $\wh{\Delta}$; we might be using $\Delta$ and $\wh{\Delta}$ interchangeably in this respect. The following asserts that these coordinates indeed form a coordinate system, and is one of the two main results of \cite{FS}. 
\begin{proposition}[{\cite{FS}}]
\label{prop:FS_coordinates}
Let $\wh{\Delta}$ be a split ideal triangulation of a triangulable generalized marked surface $\frak{S}$. Let $W$ be a weakly reduced non-elliptic ${\rm SL}_3$-web in $\frak{S}$, not necessarily canonical with respect to $\wh{\Delta}$.
\begin{enumerate}
\itemsep0em
\item[\rm (FS1)] Define Frohman-Sikora coordinates of $W$ with respect to $\wh{\Delta}$ by using any weakly reduced non-elliptic ${\rm SL}_3$-web $W'$ in $\frak{S}$ that is equivalent to $W$ and is canonical with respect to $\wh{\Delta}$. Then these coordinates are well-defined, i.e. do not depend on the choice of $W'$;

\item[\rm (FS2)] If $W$ is {\em reduced}, the Frohman-Sikora coordinates of $W$ with respect to $\wh{\Delta}$ completely determine $W$ up to equivalence, i.e. two {\em reduced} non-elliptic ${\rm SL}_3$-webs with same Frohman-Sikora coordinates are equivalent.
\end{enumerate}
\end{proposition}
This coordinate system is geometrically intuitive, and gives an injection
$$
\{\mbox{equivalence clases of reduced non-elliptic ${\rm SL}_3$-webs in $\frak{S}$}\} \longrightarrow (\mathbb{Z}_{\ge 0})^{\Delta} \times (\mathbb{Z}_{\ge 0})^{\Delta} \times \mathbb{Z}^{\mathcal{F}(\Delta)},
$$
where $\mathcal{F}(\Delta)$ is the set of all ideal triangles of $\Delta$. This coordinate map is not surjective, so one may want to study the structure of the image set; see \cite{FS} for a discussion.

\vs

We now recall another set of coordinates studied by Douglas-Sun \cite{Douglas} \cite{DS1} \cite{DS2}, which better suits our purposes. Their coordinates are parametrized by the nodes of the $3$-triangulation quiver $Q_\Delta$ (Def.\ref{def:3-triangulation}).
\begin{definition}[Douglas-Sun \cite{Douglas} \cite{DS1} \cite{DS2}]
\label{def:Douglas-Sun}
Let $\frak{S}$ be a triangulable generalized marked surface, $\Delta$ an ideal triangulation of $\frak{S}$ and $\wh{\Delta}$ a split ideal triangulation for $\Delta$. Let $Q_\Delta$ be the 3-triangulation quiver for $\Delta$ (Def.\ref{def:3-triangulation}). Define the {\rm Douglas-Sun} coordinate map
\begin{align}
\label{eq:DS}
\{\mbox{equivalence clases of reduced non-elliptic ${\rm SL}_3$-webs in $\frak{S}$}\} \longrightarrow (\textstyle \frac{1}{3} \mathbb{Z}_{\ge 0})^{\mathcal{V}(Q_\Delta)},
\end{align}
given, for a reduced non-elliptic ${\rm SL}_3$-web $W$ in $\frak{S}$, by the following number per each node of $Q_\Delta$.

\vs

Let $W'$ be any reduced non-elliptic ${\rm SL}_3$-web in $\frak{S}$ that is equivalent to $W$ and is canonical with respect to $\wh{\Delta}$. Let $t$ be a triangle of $\Delta$, and $\wh{t}$ be the triangle of $\wh{\Delta}$ corresponding to $t$, so that $W' \cap \wh{t}$ is a canonical ${\rm SL}_3$-web in $\wh{t}$. The coordinates of $W$ for these nodes $\mathcal{V}(Q_\Delta) \cap t$ are defined as the coordinates of the ${\rm SL}_3$-web $\redfix{W}' \cap \wh{t}$ in $\wh{t}$ for these nodes, given as follows. 

\vs

We require that the coordinates are additive for \redfix{each} $t$, in the sense that for $v\in \mathcal{V}(Q_\Delta)\cap t$, if $W_1,W_2$ are disjoint ${\rm SL}_3$-webs in $\wh{t}$, then the coordinate of $W_1\cap W_2$ for $v$ equals the sum of the coordinates of $W_1$ and $W_2$ for $v$. Then it suffices to define the coordinates for a corner arc ${\rm SL}_3$-web in $\wh{t}$ and for a pyramid ${\rm SL}_3$-web $H_d$ in $\wh{t}$, which are given in Fig.\ref{fig:DS_coordinates}.
\end{definition}

In fact, the original Douglas-Sun coordinates \cite{Douglas} \cite{DS1} \cite{DS2} are 3 times the ones depicted in Fig.\ref{fig:DS_coordinates}, hence are integers. The reason why we use the $\frac{1}{3}$-scaled version will be justified by one of our main results (Thm.\ref{thm:main}).
\begin{proposition}[\cite{DS1} \cite{DS2}]
\label{prop:DS_injection}
The above coordinate system yields a well-defined injection as in eq.\eqref{eq:DS}.
\end{proposition}
Def.\ref{def:Douglas-Sun} and Prop.\ref{prop:DS_injection} were stated for punctured surfaces in \cite{DS1}, and for generalized marked surfaces in \cite{DS2}.

\vs

In \cite{DS1} \cite{DS2}, the image of eq.\eqref{eq:DS} is studied in detail. As mentioned in \cite{DS1}, this coordinate system is inspired by the degrees of the highest term of a (sought-for) canonical regular function on $\mathscr{X}_{{\rm PGL}_3}$ associated to each ${\rm SL}_3$-web $W$, and this idea goes back to Xie \cite{Xie}. In a sense, the results of the present paper will fully justify this idea. Even without the result of the present paper, one can study some remarkable properties of the Douglas-Sun coordinate systems, a crucial one being the behavior under flip of an ideal triangulation.
\begin{proposition}[coordinate change formula for Douglas-Sun coordinates; {\cite[Thm.4.4]{DS2}}]
\label{prop:DS_coordinate_change}
Let $\Delta$ and $\Delta'$ be ideal triangulations of a triangulable generalized marked surface $\frak{S}$ related to each other by a flip at an edge. Let $W$ be a reduced non-elliptic ${\rm SL}_3$-web in $\frak{S}$. The Douglas-Sun coordinates $({\rm a}_v)_v \in (\frac{1}{3} \mathbb{Z}_{\ge 0})^{\mathcal{V}(Q_\Delta)}$ and $({\rm a}_{v'}')_{v'} \in (\frac{1}{3} \mathbb{Z}_{\ge 0})^{\mathcal{V}(Q_{\Delta'})}$ of $W$ with respect to $\Delta$ and $\Delta'$ are related by the sequence of tropical cluster $\mathscr{A}$-mutations with respect to the sequence of mutations associated to a flip. To be more precise, if we label the nodes of $Q_\Delta$ and $Q_{\Delta'}$ for triangles having the flipped arc as a side as in Fig.\ref{fig:mutations_for_a_flip}, then
\begin{align*}
& {\rm a}_{v_3'}' = -{\rm a}_{v_3} +  \max({\rm a}_{v_2} + {\rm a}_{v_{12}}, {\rm a}_{v_7}+{\rm a}_{v_8}), \quad
{\rm a}_{v_4'}' = -{\rm a}_{v_4} + \max({\rm a}_{v_7}+{\rm a}_{v_{11}}, {\rm a}_{v_5}+{\rm a}_{v_{12}}) \\
& {\rm a}_{v_7'}' = -{\rm a}_{v_7} + \max({\rm a}_{v_1}+{\rm a}_{v_4'}', {\rm a}_{v_6}+{\rm a}_{v_3'}'), \quad
{\rm a}_{v_{12}'}' = -{\rm a}_{v_{12}} + \max({\rm a}_{v_3'}' +{\rm a}_{v_{10}}, {\rm a}_{v_4'}' + {\rm a}_{v_9}).
\end{align*}
Nodes $v$ in $Q_\Delta$ other than $v_3,v_4,v_7,v_{12}$ in Fig.\ref{fig:mutations_for_a_flip} are naturally in bijection with nodes $v'$ in $Q_{\Delta'}$ other than $v_3',v_4',v_7',v_{12}'$, and ${\rm a}_v = {\rm a}_{v'}'$ holds for them.
\end{proposition}
Some simple cases can be checked by hand easily, but to prove it fully seems not so easy; see \cite{DS2}.

\subsection{Tropical coordinates for ${\rm SL}_3$-laminations}

We introduce a coordinate system for ${\rm SL}_3$-laminations, by extending the coordinate systems for ${\rm SL}_3$-webs, especially Douglas-Sun's. One notable aspect is that our coordinate map will eventually map bijectively onto $\mathbb{Z}^\Delta \times \mathbb{Z}^\Delta \times \mathbb{Z}^{\mathcal{F}(\Delta)}$. 
\begin{definition}[tropical coordinates for ${\rm SL}_3$-laminations]
\label{def:tropical_coordinates}
Let $\Delta$ be an ideal triangulation of a triangulable generalized marked surface $\frak{S}$. Let $Q_\Delta$ be the $3$-triangulation of $\Delta$, and let $\wh{\Delta}$ be a split ideal triangulation for $\Delta$. Let $\ell$ be an ${\rm SL}_3$-lamination in $\frak{S}$, represented by a weighted reduced non-elliptic ${\rm SL}_3$-web $W(\ell)$ that is canonical with respect to $\wh{\Delta}$. For each node $v$ of $Q_\Delta$, define the integer ${\rm a}_v(\ell)$ as follows.
\begin{enumerate}
\itemsep0em
\item[\rm (TC1)] (edge coordinates) Let $e\in \wh{\Delta}$. Let $\wh{t}$ be the unique triangle of $\wh{\Delta}$ having $e$ as a side, and let $W(\ell) \cap \wh{t}$ be naturally given the structure of a weighted weakly reduced non-elliptic ${\rm SL}_3$-web in $\wh{t}$. Let $e$ denote also the corresponding edge of $\Delta$ by abuse of notation, and let $t$ be the triangle of $\Delta$ corresponding to $\wh{t}$. Let $v_{e,1}$ and $v_{e,2}$ be the nodes of $Q_\Delta$ lying in $e\in \Delta$, such that the direction $v_{e,1} \to v_{e,2}$ matches the clockwise orientation on $\partial t$. Let the \ul{\em intersection weights} of $\ell$ at edge $e$ be
\begin{align*}
{\rm e}_{{\rm out},e}(\ell)&:=\mbox{sum of weights of edges of $W(\ell) \cap \wh{t}$ whose terminal endpoints lie in $e\in \wh{\Delta}$}, \\
{\rm e}_{{\rm in},e}(\ell)&:=\mbox{sum of weights of edges of $W(\ell) \cap \wh{t}$ whose initial endpoints lie in $e\in \wh{\Delta}$}.
\end{align*}
Define the \ul{\em edge coordinates} of $\ell$ for the edge $e$ of the ideal triangulation $\Delta$ as
\begin{align}
\label{eq:edge_coordinates}
\textstyle {\rm a}_{v_{e,1}}(\ell) :=  \frac{1}{3} {\rm e}_{{\rm out},e}(\ell) + \frac{2}{3} {\rm e}_{{\rm in},e}(\ell), \qquad
{\rm a}_{v_{e,2}}(\ell) :=  \frac{2}{3} {\rm e}_{{\rm out},e}(\ell) + \frac{1}{3} {\rm e}_{{\rm in},e}(\ell).
\end{align}

\item[\rm (TC2)] (triangle coordinates) Let $t$ be a triangle of $\Delta$, $v_t$ be the node of $Q_\Delta$ lying in the interior of $t$, and let $\wh{t}$ be the triangle of $\wh{\Delta}$ corresponding to $t$. Let $e_1,e_2,e_3$ be the sides of $\wh{t}$, appearing clockwise in this order along $\partial \wh{t}$. Let the \ul{\em rotation weight} of $\ell$ at triangle $t$ be
\begin{align*}
{\rm r}_t(\ell) & := \mbox{(sum of weights of left turn corner arcs of $W(\ell) \cap \wh{t}$)} \\
& \quad - \mbox{(sum of weights of right turn corner arcs of $W(\ell) \cap \wh{t}$)}.
\end{align*}
Define the \ul{\em degree} ${\rm d}_t(\ell)$ of $\ell$ for $t$ as
\begin{align}
\label{eq:degree}
{\rm d}_t(\ell) := \textstyle \frac{1}{3} \sum_{\alpha=1}^3 {\rm e}_{{\rm out},e_{\alpha}}(\ell) - \frac{1}{3} \sum_{\alpha=1}^3 {\rm e}_{{\rm in},e_{\alpha}}(\ell).
\end{align}
Define the \ul{\em triangle coordinate} of $\ell$ for the triangle $t$ of the triangulation $\Delta$ as
\begin{align}
\label{eq:triangle_coordinate}
{\rm a}_{v_t}(\ell) & := \left\{
{\renewcommand{\arraystretch}{1.4}
\begin{array}{ll}
\frac{1}{6}( {\rm r}_t(\ell) + 3{\textstyle \sum}_{\alpha=1}^3 {\rm a}_{v_{e_\alpha,2}}(\ell)) & \mbox{if ${\rm d}_t(\ell)\ge 0$,} \\
\frac{1}{6}( {\rm r}_t(\ell) + 3{\textstyle \sum}_{\alpha=1}^3 {\rm a}_{v_{e_\alpha,1}}(\ell) ) & \mbox{if ${\rm d}_t(\ell)\le 0$.} \\
\end{array}} \right.
\end{align}
\end{enumerate}
The numbers ${\rm a}_v(\ell)$ are called \ul{\em tropical coordinates} for $\ell$.
\end{definition}
\begin{remark}
\label{rem:degree_is_generalized}
The degree ${\rm d}_t(\ell)$ is a generalization of a corresponding concept defined for ${\rm SL}_3$-webs in \cite[\S12]{FS}, which detects the degree of the pyramid in triangles of $\wh{\Delta}$.
\end{remark}

It is easy to see that the degree for a triangle $t$ can be expressed using the edge coordinates for the sides $e_1,e_2,e_3$ of that triangle: 
\begin{align}
\label{eq:d_t_in_terms_of_a}
{\rm d}_t(\ell) = \textstyle \sum_{\alpha=1}^3 {\rm a}_{v_{e_\alpha,2}}(\ell) - \sum_{\alpha=1}^3 {\rm a}_{v_{e_\alpha,1}}(\ell).
\end{align}
The reason why we consider the specific \redfix{(piecewise-)}linear combinations of analogs of Frohman-Sikora coordinates as in eq.\eqref{eq:edge_coordinates} and eq.\eqref{eq:triangle_coordinate}, as well as the word {\em tropical}, is related to the coordinate change formula under flips of triangulations which we will soon discuss, as seen for Douglas-Sun coordinates of ${\rm SL}_3$-webs. Indeed, one can verify that our coordinates agree with Douglas-Sun's on reduced non-elliptic ${\rm SL}_3$-webs which can naturally be viewed as ${\rm SL}_3$-laminations (with all weights being $1$).
\begin{lemma}
\label{lem:ours_are_compatible_with_DS}
Let $\ell$ be represented by a reduced non-elliptic ${\rm SL}_3$-web $W$ with \redfix{all} weight\redfix{s being} $1$, in a triangulable generalized marked surface $\frak{S}$. For an ideal triangulation $\Delta$ and for each node $v$ of $Q_\Delta$, our coordinate ${\rm a}_v(\ell)$ coincides with Douglas-Sun's coordinate of $W$ at $v$.
\end{lemma}
This lemma may be useful already, \redfix{because in \cite{DS1}} Douglas-Sun's coordinates are defined \redfix{in the style of} Def.\ref{def:Douglas-Sun} \redfix{only}, \redfix{while} explicit formulas for them in terms of Frohman-Sikora coordinates are not given. We postpone a proof of this lemma until a little bit later. 

\vs

\redfix{For now, we begin with some basic observations. Let $\frak{S}$ be a triangulable generalized marked surface, $\Delta$ an ideal triangulation of $\frak{S}$, $\wh{\Delta}$ a split ideal triangulation for $\Delta$, and $\ell$ an ${\rm SL}_3$-lamination in $\frak{S}$.} \bluefix{Represent $\ell$ as a weighted reduced non-elliptic ${\rm SL}_3$-web $W(\ell)$ that is canonical with respect to $\wh{\Delta}$. Let $\wh{t}$ be a triangle of $\wh{\Delta}$ corresponding to a triangle $t$ of $\Delta$, and let $e_1,e_2,e_3$ be the sides of $\wh{t}$ appearing clockwise in this order along $\partial \wh{t}$. For $\alpha,\beta \in \{1,2,3\}$, let
$$
c_{\alpha,\beta} = c_{\alpha,\beta;t}(\ell)
$$
be the sum of weights of corner arcs of the ${\rm SL}_3$-web $W(\ell)\cap \wh{t}$ in $\wh{t}$ going from edge $e_\alpha$ to edge $e_\beta$. By Def.\ref{def:canonical_wrt_split_ideal_triangulation}(CW1) and Def.\ref{def:canonical_web_in_a_triangle}, it follows that $W(\ell) \cap \wh{t}$ is a union of a single pyramid $H_{d_t}$ of some degree $d_t \in \mathbb{Z}$ and \redfix{some number of} corner arcs. The lamination \redfix{$\ell_t := \ell\cap \wh{t}$ in $\wh{t}$} is completely determined by these numbers $c_{\alpha,\beta}$ and $d_t$. By Def.\ref{def:A2-lamination}(L1), $H_{d_t}$ has weight $1$. So we have
\begin{align}
\label{eq:FS_intersection_coordinates_of_lamination}
{\rm e}_{{\rm out},e_\alpha}(\ell) = c_{\alpha+1,\alpha}+c_{\alpha-1,\alpha}+[d_t]_+, 
\qquad 
{\rm e}_{{\rm in},e_\alpha}(\ell) = c_{\alpha,\alpha+1}+c_{\alpha,\alpha-1}+[-d_t]_+,
\qquad
\alpha=1,2,3,
\end{align}
where the subscript indices in $c_{\alpha,\beta}$ are considered modulo $3$ (e.g. $c_{42} = c_{12}$), and $[\sim]_+$ is as in eq.\eqref{eq:positive_part}. In particular, note $[a]_+ = (a+|a|)/2$, and hence $[a]_+-[-a]_+=a$.  Thus we observe
\begin{align}
\label{eq:FS_degree_of_lamination}
\textstyle {\rm d}_t(\ell) ~\stackrel{{\rm eq}.\eqref{eq:degree}}{=}~
\frac{1}{3} \sum_{\alpha=1}^3 {\rm e}_{{\rm out},e_\alpha}(\ell) - \frac{1}{3} \sum_{\alpha=1}^3 {\rm e}_{{\rm in},e_\alpha}(\ell) = [d_t]_+ - [-d_t]_+ = d_t,
\end{align}
justifying Rem.\ref{rem:degree_is_generalized}.}

\vs

\bluefix{
In \cite[Lem.23]{FS} it is observed that the intersection coordinates and the rotation numbers of an ${\rm SL}_3$-web completely determine the number of each kind of corner arcs. Likewise, we show that the weights $c_{\alpha,\beta}$ of corner arcs are completely determined by the intersection weights and the rotation weight; we will give explicit reconstruction formulas. For convenience, let $\til{\ell}_t$ be the ${\rm SL}_3$-lamination in $\wh{t}$ obtained from $\ell_t$ by removing the pyramid \redfix{$H_{d_t}$}. Then $\ell_t$ and $\til{\ell}_t$ have same corner weights $c_{\alpha,\beta}$, hence the same rotation weight ${\rm r}_t$. The intersection weights of $\til{\ell}_t$, denoted by $\til{\rm e}_{{\rm out},e_\alpha}$ and $\til{\rm e}_{{\rm in},e_\alpha}$, are obtained from those of $\ell_t$ by subtracting $[d_t]_+$ and $[-d_t]_+$:
$$
\til{\rm e}_{{\rm out},e_\alpha} = {\rm e}_{{\rm out},e_\alpha}(\ell) - [d_t]_+ = c_{\alpha+1,\alpha} + c_{\alpha-1,\alpha}, \qquad
\til{\rm e}_{{\rm in},e_\alpha} = {\rm e}_{{\rm in},e_\alpha}(\ell) - [-d_t]_+ = c_{\alpha,\alpha+1}+c_{\alpha,\alpha-1}.
$$
Note
$$
{\rm r}_t = \textstyle \sum_{\alpha=1}^3 c_{\alpha,\alpha+1} - \sum_{\alpha=1}^3 c_{\alpha+1,\alpha}.
$$
\redfix{Define} the left and the right rotation weights \redfix{as}
$$
\textstyle {\rm r}_{\rm left} := \sum_{\alpha=1}^3 c_{\alpha,\alpha+1}, \qquad {\rm r}_{\rm right} := \sum_{\alpha=1}^3 c_{\alpha+1,\alpha}.
$$
Note $\sum_{\alpha=1}^3 \til{\rm e}_{{\rm out},e_\alpha} = \sum_{\alpha=1}^3 \til{\rm e}_{{\rm in},e_\alpha} = {\rm r}_{\rm left} + {\rm r}_{\rm right}$, while ${\rm r}_t = {\rm r}_{\rm left} - {\rm r}_{\rm right}$. Thus we can express ${\rm r}_{\rm left}$ and ${\rm r}_{\rm right}$ in terms of the intersection weights and the rotation weight as
$$
\textstyle {\rm r}_{\rm left} = \frac{1}{2} (\sum_{\alpha=1}^3 \til{\rm e}_{{\rm out},e_\alpha}) + \frac{1}{2} {\rm r}_t, \qquad
{\rm r}_{\rm right} = \frac{1}{2} (\sum_{\alpha=1}^3 \til{\rm e}_{{\rm out},e_\alpha}) - \frac{1}{2} {\rm r}_t.
$$
Observe now
\begin{align*}
& {\rm r}_{\rm left} + \til{\rm e}_{{\rm in},e_1} + \til{\rm e}_{{\rm out},e_2} - \til{\rm e}_{{\rm out},e_3} - \til{\rm e}_{{\rm in},e_3}  \\
& = (c_{1,2}+c_{2,3}+c_{3,1})+(c_{1,2}+c_{1,3})+(c_{3,2}+c_{1,2})-(c_{1,3}+c_{2,3})-(c_{3,1}+c_{3,2}) = 3 c_{1,2}.
\end{align*}
Exchanging ${\rm r}_{\rm left}$ with ${\rm r}_{\rm right}$, and each $\til{\rm e}_{{\rm in},e_\alpha}$ with $\til{\rm e}_{{\rm out},e_\alpha}$ and vice versa results in exchanging the order of subscripts of $c_{\alpha,\beta}$, so we obtain ${\rm r}_{\rm right} + \til{\rm e}_{{\rm out},e_1} + \til{\rm e}_{{\rm in},e_2} - \til{\rm e}_{{\rm in},e_3} - \til{\rm e}_{{\rm out},e_3}=3c_{2,1}$. By the cyclicity of the subscript indices $1,2,3$, we thus get
\begin{align*}
& \textstyle c_{\alpha,\alpha+1} = \frac{1}{3}({\rm r}_{\rm left} + \til{\rm e}_{{\rm in},e_\alpha} + \til{\rm e}_{{\rm out},e_{\alpha+1}} - \til{\rm e}_{{\rm out},e_{\alpha-1}} - \til{\rm e}_{{\rm in},e_{\alpha-1}}), \\
& \textstyle c_{\alpha+1,\alpha} = \frac{1}{3}({\rm r}_{\rm right} + \til{\rm e}_{{\rm out},e_\alpha} + \til{\rm e}_{{\rm in},e_{\alpha+1}} - \til{\rm e}_{{\rm in},e_{\alpha-1}} - \til{\rm e}_{{\rm out},e_{\alpha-1}}).
\end{align*}
\redfix{So} we expressed all corner weights $c_{\alpha,\beta}$ in terms of the intersection weights and the rotation weight.
}

\vs

\redfix{One can also express these structural numbers $c_{\alpha,\beta}$ and $d_t$ for $\ell_t = \ell \cap \wh{t}$ completely in terms of our tropical coordinates ${\rm a}_v(\ell)$ for the nodes $v$ of $Q_\Delta$ living in the triangle $t$. The formula for $d_t = {\rm d}_t(\ell)$ is already in eq.\eqref{eq:d_t_in_terms_of_a}, which can be used to rewrite eq.\eqref{eq:triangle_coordinate} as
$$
\textstyle {\rm a}_{v_t}(\ell) = \frac{1}{6} {\rm r}_t(\ell) + \frac{1}{2} \sum_{\alpha=1}^3 {\rm a}_{v_{e_\alpha,1}}(\ell) + \frac{1}{2} [{\rm d}_t(\ell)]_+ = \frac{1}{6} {\rm r}_t(\ell) + \frac{1}{2} \sum_{\alpha=1}^3 {\rm a}_{v_{e_\alpha,2}}(\ell) + \frac{1}{2} [-{\rm d}_t(\ell)]_+ ,
$$
which in turn enables us to express the rotation weight ${\rm r}_t(\ell) = {\rm r}_t$ as
$$
\textstyle {\rm r}_t(\ell) = 6 \, {\rm a}_{v_t}(\ell) - 3\sum_{\alpha=1}^3 {\rm a}_{v_{e_\alpha,1}}(\ell) - 3 [{\rm d}_t(\ell)]_+
= 6 \, {\rm a}_{v_t}(\ell) - 3\sum_{\alpha=1}^3 {\rm a}_{v_{e_\alpha,2}}(\ell) - 3 [-{\rm d}_t(\ell)]_+
$$
in terms of the tropical coordinates; use eq.\eqref{eq:d_t_in_terms_of_a} for $d_t = {\rm d}_t(\ell)$. Then, using the arguments above for $c_{\alpha,\beta}$ and the following easy observation from eq.\eqref{eq:edge_coordinates}
$$
{\rm e}_{{\rm out},e_\alpha} (\ell) = 2 {\rm a}_{v_{e_\alpha,2}} (\ell) - {\rm a}_{v_{e_\alpha,1}} (\ell), \qquad
{\rm e}_{{\rm in},e_\alpha} (\ell) = 2 {\rm a}_{v_{e_\alpha,1}} (\ell) - {\rm a}_{v_{e_\alpha,2}} (\ell),
$$
one can compute an explicit expression for each $c_{\alpha,\beta} = c_{\alpha,\beta;t}(\ell)$ in terms of the tropical coordinates:
\begin{align*}
\textstyle c_{\alpha,\alpha+1;t}(\ell) 
& = {\rm a}_{v_t}(\ell) + {\rm a}_{v_{e_{\alpha+1},2}}(\ell) - {\rm a}_{v_{e_{\alpha+1},1}}(\ell) - {\rm a}_{v_{e_{\alpha-1},1}}(\ell)  - [{\rm d}_t(\ell)]_+ , \\ 
\textstyle c_{\alpha+1,\alpha;t}(\ell) & = - {\rm a}_{v_t}(\ell) + {\rm a}_{v_{e_{\alpha},2}}(\ell)+ {\rm a}_{v_{e_{\alpha+1},1}} (\ell). 
\end{align*}

}

\vs

\redfix{We are ready to state} the first major assertion about our coordinates for ${\rm SL}_3$-laminations.
\begin{proposition}
\label{prop:tropical_coordinate_is_well-defined}
The coordinates of Def.\ref{def:tropical_coordinates} provide a well-defined map
\begin{align*}
{\bf a}_\Delta : \{\mbox{${\rm SL}_3$-laminations in $\frak{S}$}\} & \longrightarrow \mathcal{B}_\Delta\subset ({\textstyle \frac{1}{3}}\mathbb{Z})^{\mathcal{V}(Q_\Delta)} \\
\ell & \longmapsto ({\rm a}_v(\ell))_{v\in \mathcal{V}(Q_\Delta)}
\end{align*}
where $\mathcal{B}_\Delta$ is the set of all \ul{\em balanced} elements of $(\frac{1}{3}\mathbb{Z})^{\redfix{\mathcal{V}}(Q_\Delta)}$, where an element $({\rm a}_v)_v \in (\frac{1}{3}\mathbb{Z})^{\redfix{\mathcal{V}}(Q_\Delta)}$ is said to be balanced if, for each triangle $t$ of $\Delta$, \redfix{with its sides denoted by $e_1,e_2,e_3$ in the clockwise order and the nodes of $Q_\Delta$ living in $t$ denoted as in Def.\ref{def:tropical_coordinates},} \redfix{
\begin{enumerate}
\item[\rm (BE1)] the numbers $\textstyle \sum_{\alpha=1}^3 {\rm a}_{v_{e_\alpha,1}}$ and $\textstyle \sum_{\alpha=1}^3{\rm a}_{v_{e_\alpha,2}}$ both belong to $\mathbb{Z}$;

\item[\rm (BE2)] for each $\alpha=1,2,3$, the number ${\rm a}_{v_{e_\alpha,1}} + {\rm a}_{v_{e_\alpha,2}}$ belongs to $\mathbb{Z}$;

\item[\rm (BE3)] for each $\alpha=1,2,3$, the number $-{\rm a}_{v_t} + {\rm a}_{v_{e_\alpha,2}} + {\rm a}_{v_{e_{\alpha+1},1}}$ (or the number ${\rm a}_{v_t} + {\rm a}_{v_{e_\alpha,1}} + {\rm a}_{v_{e_{\alpha+1},2}}$) belongs to $\mathbb{Z}$.
\end{enumerate}}
\end{proposition}

{\it Proof.} \redfix{The} well-definedness of the above coordinates follows from that of Frohman-Sikora coordinates for (weakly) reduced non-elliptic ${\rm SL}_3$-webs. Next, we should check whether the coordinates have values in $\frac{1}{3}\mathbb{Z}$; this is clear for the edge coordinates. Let's show that the triangle coordinates also have values in $\frac{1}{3}\mathbb{Z}$. \redfix{Pick any ${\rm SL}_3$-lamination $\ell$ in $\frak{S}$.} Let $\wh{t}$ be a triangle of $\wh{\Delta}$ corresponding to a triangle $t$ of $\Delta$, \redfix{with all the notations as before.} When ${\rm d}_t(\ell)=d_t\ge 0$, note
\begin{align*}
& \textstyle \redfix{\rm r}_t(\ell) + 3\sum_{\alpha=1}^3 {\rm a}_{v_{e_\alpha,2}}(\ell)
= (\sum_{\alpha=1}^3 c_{\alpha,\alpha+1} - \sum_{\alpha=1}^3 c_{\alpha+1,\alpha}) + (2\sum_{\alpha=1}^3 {\rm e}_{{\rm out},e_\alpha}
+ \sum_{\alpha=1}^3 {\rm e}_{{\rm in},e_\alpha}) \\
& \textstyle = 4\sum_{\alpha=1}^3 c_{\alpha,\alpha+1} + 2 \sum_{\alpha=1}^3 c_{\alpha+1,\alpha} + \redfix{6} [d_t]_+ + \cancel{\redfix{3} [-d_t]_+} ~\in~ 2\mathbb{Z}
\end{align*}
so ${\rm a}_{v_t}(\ell) = \frac{1}{6}(r_t(\ell) + 3\sum_{\alpha=1}^3 {\rm a}_{v_{e_\alpha,2}}(\ell)) \in \frac{1}{3} \mathbb{Z}$. Similarly, when ${\rm d}_t(\ell)=d_t \le 0$, one observes ${\rm a}_{v_t}(\ell) = \frac{1}{6}(\redfix{\rm r}_t(\ell)+3\sum_{\alpha=1}^3 {\rm a}_{v_{e_\alpha},1}(\ell))=\frac{1}{6}(4\sum_{\alpha=1}^3 c_{\alpha,\alpha+1} + 2\sum_{\alpha=1}^3 c_{\alpha+1,\alpha} + \cancel{\redfix{3}[d_t]_+} + \redfix{6}[-d_t]_+) \in \frac{1}{3}\mathbb{Z}$. So indeed, all coordinate values lie in $\frac{1}{3}\mathbb{Z}$. \redfix{By the previous arguments, the image $({\rm a}_v)_v = ({\rm a}_v(\ell))_v$ satisfies 
\begin{enumerate}
\item[\rm (BE1')] the number ${\rm d}_t := \textstyle \sum_{\alpha=1}^3 {\rm a}_{v_{e_\alpha,2}} - \sum_{\alpha=1}^3{\rm a}_{v_{e_\alpha,1}}$ belongs to $\mathbb{Z}$;

\item[\rm (BE3')] for each $\alpha=1,2,3$, the numbers $c_{\alpha,\alpha+1;t}:={\rm a}_{v_t} + {\rm a}_{v_{e_{\alpha+1},2}} - {\rm a}_{v_{e_{\alpha+1},1}} - {\rm a}_{v_{e_{\alpha-1},1}} - [{\rm d}_t]_+$ and $c_{\alpha+1,\alpha;t}:=-{\rm a}_{v_t} + {\rm a}_{v_{e_\alpha,2}} + {\rm a}_{v_{e_{\alpha+1},1}}$ belong to $\mathbb{Z}$,
\end{enumerate}
and (BE2). Then $c_{\alpha,\alpha+1;t} + c_{\alpha+1,\alpha;t} + [{\rm d}_t]_+ +({\rm a}_{v_{e_{\alpha-1},1}} + {\rm a}_{v_{e_{\alpha-1},2}}) = \sum_{\beta=1}^3 {\rm a}_{v_{e_{\beta},2}}$ belongs to $\mathbb{Z}$, and hence so does ${\rm d}_t - \sum_{\beta=1}^3 {\rm a}_{v_{e_{\beta},2}} = \sum_{\beta=1}^3 {\rm a}_{v_{e_{\beta},1}}$; therefore (BE1) holds. Thus the image $({\rm a}_v)_v = ({\rm a}_v(\ell))_v$ of the coordinate map is balanced. For later use, one should also remark that having (BE1), (BE2), (BE3) is equivalent to having (BE1'), (BE2'), (BE3'). \qed}

\vs

We shall prove that the image of the coordinate map coincides with $\mathcal{B}_\Delta$.  We first establish one useful lemma, which is straightforward to see.
\begin{definition}
\label{def:disjoint}
We say that ${\rm SL}_3$-laminations $\ell_1,\ldots,\ell_n$ in a generalized marked surface are \ul{\em disjoint} if they can be represented by weighted reduced non-elliptic ${\rm SL}_3$-webs that are mutually disjoint. We denote by $\ell_1 \cup \cdots \cup \ell_n$ the ${\rm SL}_3$-lamination obtained by taking the union of them.
\end{definition}
\begin{lemma}[additivity of coordinates under disjoint union]
\label{lem:additivity_of_coordinates}
Suppose that $\ell_1,\cdots,\ell_n$ are disjoint ${\rm SL}_3$-laminations in a triangulable generalized marked  surface $\frak{S}$. For any triangulation $\Delta$ of $\frak{S}$, we have
$$
{\bf a}_\Delta(\ell_1\cup \cdots \cup \ell_n) = {\bf a}_\Delta(\ell_1) + \cdots + {\bf a}_\Delta(\ell_n),
$$
i.e. ${\rm a}_v(\ell_1 \cup \cdots\cup \ell_n) = {\rm a}_v(\ell_1) + \cdots + {\rm a}_v(\ell_n)$ holds for every node $v$ of $Q_\Delta$.
\end{lemma}
{\it Proof.} It suffices to prove the assertion when $\frak{S}$ is a triangle, and when each $\ell_i$ can be represented by a single-component weakly reduced non-elliptic ${\rm SL}_3$-web in a triangle. One can see that at most one of $\ell_1,\ldots,\ell_n$ can contain an internal vertex, i.e. can be a pyramid $H_d$ with $d\neq 0$, and others are all corner arcs. If $d>0$ or there is no pyramid, then ${\rm d}_t(\ell_i)\ge 0$ for all $i$ and ${\rm d}_t(\ell_1\cup \cdots \cup \ell_n)\ge 0$, hence ${\rm a}_{v_t}(\ell_i)$ as well as ${\rm a}_{v_t}(\ell_1\cup \cdots \cup \ell_n)$ is given by the first line formula of eq.\eqref{eq:triangle_coordinate}. So, all of ${\rm a}_{v_{e,1}}(\cdot)$, ${\rm a}_{v_{e,2}}(\cdot)$, ${\rm r}_t(\cdot)$, and hence also ${\rm a}_{v_t}(\cdot)$, are additive for $\ell_1,\ldots,\ell_n$. Likewise, if $d<0$ or there is no pyramid, the second line formula of eq.\eqref{eq:triangle_coordinate} applies to all $\ell_1,\ldots,\ell_n$ and $\ell_1\cup\cdots \cup \ell_n$, so the coordinates are additive. \qed

\vs

Before proceeding further, we use this lemma to prove the promised easy lemma, Lem.\ref{lem:ours_are_compatible_with_DS}.

\vs

{\it Proof of Lem.\ref{lem:ours_are_compatible_with_DS}.} It suffices to check this for \redfix{each} node $v$ living in a triangle $t$ of $\Delta$. Let $\wh{t}$ be the triangle of the split ideal triangulation $\wh{\Delta}$. When $W$ is canonical with respect to $\wh{\Delta}$, note that $W\cap \wh{t}$ is a (weakly) reduced ${\rm SL}_3$-web in $\wh{t}$, and the tropical coordinates of $W\cap \wh{t}$ coincide with the tropical coordinates of $W$ for the nodes of $Q_\Delta$ living in $t$. Since our coordinates are additive (Lem.\ref{lem:additivity_of_coordinates}) and so are Douglas-Sun's by construction, it suffices to show the equality for each component of $W\cap \wh{t}$, which is a corner arc or a pyramid. For these elementary cases, the Douglas-Sun coordinates are as in Fig.\ref{fig:DS_coordinates}, which we verify to be same as ours as follows.

\vs

Let the side names $e_\alpha$, node names $v_{e_\alpha,1}$, $v_{e_\alpha,2}$, $v_t$ be as in Def.\ref{def:tropical_coordinates}. Let $W_{\alpha,\alpha+1}$ be a left turn corner arc ${\rm SL}_3$-web in $\wh{t}$, going from side $e_\alpha$ to $e_{\alpha+1}$ (where $e_4=e_1$). The Frohman-Sikora coordinates are $0 = {\rm a}_{v_{e_\alpha},{\rm out}}(W_{\alpha,\alpha+1}) = {\rm a}_{v_{e_{\alpha+1}},{\rm in}}(W_{\alpha,\alpha+1})={\rm a}_{v_{e_{\alpha+2}},{\rm in}}(W_{\alpha,\alpha+1}) = {\rm a}_{v_{e_{\alpha+2}},{\rm out}}(W_{\alpha,\alpha+1})$, $1 = {\rm a}_{v_{e_\alpha},{\rm in}}(W_{\alpha,\alpha+1})={\rm a}_{v_{e_{\alpha+1}},{\rm out}}(W_{\alpha,\alpha+1})$, ${\rm r}_t(W_{\alpha,\alpha+1}) = 1$, and ${\rm d}_t(W_{\alpha,\alpha+1}) = 0$, so
\begin{align}
\label{eq:tropical_coordinates_for_a_left_turn_W}
\left\{ \begin{array}{l}
{\rm a}_{v_{e_\alpha},1}(W_{\alpha,\alpha+1}) = \frac{2}{3} = {\rm a}_{v_{e_{\alpha+1}},2}(W_{\alpha,\alpha+1}), \quad
{\rm a}_{v_{e_\alpha},2}(W_{\alpha,\alpha+1}) = \frac{1}{3} = {\rm a}_{v_{e_{\alpha+1}},1}(W_{\alpha,\alpha+1}), \\
{\rm a}_{v_{e_{\alpha+2}},1}(W_{\alpha,\alpha+1}) = {\rm a}_{v_{e_{\alpha+2}},2}(W_{\alpha,\alpha+1})=0, \quad
{\rm a}_{v_t}(W_{\alpha,\alpha+1}) = \frac{1}{6}(1+3)=\frac{2}{3},
\end{array}
\right.
\end{align}
which matches Fig.\ref{fig:DS_coordinates}. Let $W_{\alpha+1,\alpha}$ be a right turn corner arc ${\rm SL}_3$-web in $\wh{t}$, going from side $e_{\alpha+1}$ to $e_\alpha$. The Frohman-Sikora coordinates are $0 = {\rm a}_{v_{e_\alpha},{\rm in}}(W_{\alpha+1,\alpha}) = {\rm a}_{v_{e_{\alpha+1}},{\rm out}}(W_{\alpha+1,\alpha})={\rm a}_{v_{e_{\alpha+2}},{\rm in}}(W_{\alpha+1,\alpha}) = {\rm a}_{v_{e_{\alpha+2}},{\rm out}}(W_{\alpha+1,\alpha})$, $1 = {\rm a}_{v_{e_\alpha},{\rm out}}(W_{\alpha+1,\alpha})={\rm a}_{v_{e_{\alpha+1}},{\rm in}}(W_{\alpha+1,\alpha})$, ${\rm r}_t(W_{\alpha+1,\alpha}) = -1$, and ${\rm d}_t(W_{\alpha+1,\alpha}) = 0$, so
\begin{align}
\label{eq:tropical_coordinates_for_a_right_turn_W}
\left\{ \begin{array}{l}
{\rm a}_{v_{e_\alpha},1}(W_{\alpha+1,\alpha}) = \frac{1}{3} = {\rm a}_{v_{e_{\alpha+1}},2}(W_{\alpha+1,\alpha}), \quad
{\rm a}_{v_{e_\alpha},2}(W_{\alpha+1,\alpha}) = \frac{2}{3} = {\rm a}_{v_{e_{\alpha+1}},1}(W_{\alpha+1,\alpha}), \\
{\rm a}_{v_{e_{\alpha+2}},1}(W_{\alpha+1,\alpha}) = {\rm a}_{v_{e_{\alpha+2}},2}(W_{\alpha+1,\alpha})=0, \quad
{\rm a}_{v_t}(W_{\alpha+1,\alpha}) = \frac{1}{6}(-1+3)=\frac{1}{3},
\end{array}
\right.
\end{align}
which matches Fig.\ref{fig:DS_coordinates}. Now let $H_d$ be a pyramid with $d>0$. The Frohman-Sikora coordinates are $0={\rm a}_{v_{e_\alpha,{\rm in}}}(H_d)$, $d = {\rm a}_{v_{e_\alpha,{\rm out}}}(H_d)$, $\alpha=1,2,3$, ${\rm r}_t(H_d)=0$, and ${\rm d}_t(H_d) = d>0$, so 
\begin{align}
\label{eq:tropical_coordinates_for_Hd_positive}
\textstyle {\rm a}_{v_{e_\alpha},1}(H_d) = \frac{d}{3}, \quad
 {\rm a}_{v_{e_{\redfix{\alpha}}},2}(H_d) = \frac{2d}{3}, \quad
\alpha=1,2,3, \quad
{\rm a}_{v_t}(H_d) = \frac{1}{6}(0+6d)=d.
\end{align}
which matches Fig.\ref{fig:DS_coordinates}. Finally, let $H_d$ be a pyramid with $d<0$. The Frohman-Sikora coordinates are $-d={\rm a}_{v_{e_\alpha,{\rm in}}}(H_d)$, $0 = {\rm a}_{v_{e_\alpha,{\rm out}}}(H_d)$, $\alpha=1,2,3$, ${\rm r}_t(H_d)=0$, and ${\rm d}_t(H_d) = d<0$, so 
\begin{align}
\nonumber
\textstyle {\rm a}_{v_{e_\alpha},1}(H_d) = -\frac{2d}{3}, \quad
 {\rm a}_{v_{e_{\redfix{\alpha}}},2}(H_d) = -\frac{d}{3}, \quad
\alpha=1,2,3, \quad
{\rm a}_{v_t}(H_d) = \frac{1}{6}(0-6d)=-d.
\end{align}
which matches Fig.\ref{fig:DS_coordinates}. \qed

\vs

What we will use right away is another easy observation.
\begin{lemma}
\label{lem:disjoint_union_with_peripheral_curve}
\redfix{Let} $\ell_0$ be a\redfix{n} ${\rm SL}_3$-lamination in $\frak{S}$ represented by a weighted ${\rm SL}_3$-web \redfix{consisting only of} oriented peripheral curves. Then for any ${\rm SL}_3$-lamination $\ell$ in $\frak{S}$,
\begin{enumerate}
\itemsep0em
\item[\rm (1)] $\ell_0$ is disjoint from $\ell$;

\item[\rm (2)] If we denote by $-\ell_0$ the ${\rm SL}_3$-lamination obtained from $\ell_0$ by multiplying the weight on each constituent peripheral curve by $-1$, then $\ell \cup \ell_0 \cup (-\ell_0) = \ell$ as ${\rm SL}_3$-laminations. \qed
\end{enumerate}
\end{lemma}

We now prove:
\begin{proposition}
\label{prop:tropical_coordinates_are_coordinates}
The coordinate map in Prop.\ref{prop:tropical_coordinate_is_well-defined} is a bijection onto $\mathcal{B}_\Delta$.
\end{proposition}

{\it Proof.} We construct an inverse map to the coordinate map. Let $\vec{\rm a} = ({\rm a}_v)_{v\in\mathcal{V}(Q_\Delta)}$ be any balanced element of $(\frac{1}{3}\mathbb{Z})^{\redfix{\mathcal{V}}(Q_\Delta)}$. We will construct an ${\rm SL}_3$-lamination $\ell$ having these as its coordinates. We shall construct a weighted ${\rm SL}_3$-web in each triangle of $\wh{t}$, `fill in' the biangles, then remove the boundary 4-gons, to construct a sought-for ${\rm SL}_3$-lamination $\ell$ in $\frak{S}$. Let $t$ be a triangle of $\Delta$, and let $\wh{t}$ be the corresponding triangle of $\wh{\Delta}$. \bluefix{Let's show that there exists a unique \redfix{canonical ${\rm SL}_3$-lamination} $\ell_t$ \redfix{in $\wh{t}$} whose tropical coordinates coincide with those assigned by the element $\vec{\rm a}$.} \redfix{Let's use the notations and arguments for investigation of the structure of $\ell_t$, as we presented right before Prop.\ref{prop:tropical_coordinate_is_well-defined}. As we saw already, $\ell_t$ is completely determined by its corner weights $c_{\alpha,\beta}$ and the degree $d_t$, and these numbers are in turn determined by the coordinate numbers $\vec{a}$. By the balancedness condition of $\vec{a}$ as written in (BE1') and (BE3'), the numbers $c_{\alpha,\beta}$ and $d_t$ are integers. Hence we indeed obtain a well-defined unique canonical ${\rm SL}_3$-lamination $\ell_t$ in $\wh{t}$ whose tropical coordinates equal to those assigned by $\vec{\rm a}$.}

\vs

We now modify $\vec{\rm a}$ before proceeding to biangles. Let $\vec{k} \in (\mathbb{Z}\times \mathbb{Z})^\mathcal{P}$ be the choice of two integers $k_{{\rm CW};p}$ and $k_{{\rm CCW};p}$ for each marked point $p \in \mathcal{P}$ of $\frak{S}$. Let $\ell_{\vec{k}}$ be the ${\rm SL}_3$-lamination in $\frak{S}$ consisting of two peripheral curves per marked point $p\in\mathcal{P}$ surrounding $p$ (so having $2|\mathcal{P}|$ components), where one of them has clockwise orientation around $p$ \redfix{(i.e. the orientation opposite to the boundary-orientation of the corresponding hole of $\til{\frak{S}}$ (Def.\ref{def:holed_surface}))} and has weight $k_{{\rm CW};p}$, and the other is counterclockwise \redfix{(i.e. the orientation matching the boundary-orientation of the corresponding hole of $\til{\frak{S}}$)} with weight $k_{{\rm CCW};p}$. Let $\vec{\rm a}_{\vec{k}} := {\bf a}_\Delta(\ell_{\vec{k}})$ be the element of $(\frac{1}{3}\mathbb{Z})^{\mathcal{V}(Q_\Delta)}$ for the tropical coordinates of $\ell_{\vec{k}}$. Let
$$
\textstyle \vec{\rm a}' := \vec{\rm a} + \vec{\rm a}_{\vec{k}} \in (\frac{1}{3}\mathbb{Z})^{\mathcal{V}(Q_\Delta)}.
$$
Now, repeat the previous process for this new element $\vec{\rm a}'$, to get an ${\rm SL}_3$-lamination $\ell_t'$ in each triangle $\wh{t}$. We claim that we can choose $\vec{k}$ so that all corner weights of $\ell'_t$ for all triangles $t$ are non-negative. For example, for a fixed positive integer $k$, let $\vec{k}$ be such that $k_{{\rm CW};p} = k_{{\rm CCW};p} = k$ for all $p\in \mathcal{P}$. Each corner of a triangle $\wh{t}$ of $\wh{\Delta}$ is attached at some unique marked point of $\mathcal{P}$, and hence in this corner there are $k$ corner arcs of $\ell_{\vec{k}} \cap \wh{t}$ in one direction and $k$ corner arcs of $\ell_{\vec{k}} \cap \wh{t}$ in the opposite direction. Meanwhile, the degree ${\rm d}_t(\ell_{\vec{k}})$ is zero for all triangles $t$. Hence, each corner weight $c_{\alpha,\beta}'$ of $\ell'_t$ constructed from $\vec{\rm a}'$ equals $c_{\alpha,\beta}+k$, while the degree ${\rm d}_t(\ell_t')$ equals ${\rm d}_t(\ell_t) = d_t$. So, for a sufficiently large $k$, we see that all corner weights of $\ell_t'$ for each triangle $t$ are non-negative.

\vs

Now, represent $\ell_t'$ by a weighted ${\rm SL}_3$-web $W_t'$ such that all weights are $1$. So, for each ordered pair $(\alpha,\beta)$ of distinct indices in $\{1,2,3\}$, we draw $c_{\alpha,\beta}'$ number of corner arcs going from side $e_\alpha$ to side $e_\beta$, and a pyramid $H_{d_t}$, so that these are all disjoint. This way, in each side $e_i$ of the triangle $\wh{t}$, there are some ${\rm e}'_{{\rm in},e_\alpha}$ number of source external vertices, and some ${\rm e}'_{{\rm out},e_\alpha}$ number of sink external vertices. Pick one side $e_\alpha$ of $\wh{t}$. Let $f_\beta$ be the edge of $\wh{\Delta}$ parallel to $e_\alpha$ hence forming a biangle with $e_\alpha$, where $f_1,f_2,f_3$ are edges of the triangle $\wh{r}$ of $\wh{\Delta}$ (corresponding to triangle $r$ of $\Delta$ adjacent to $t$), where the ${\rm SL}_3$-lamination $\ell'_r$ in $\wh{r}$ is drawn as a weighted ${\rm SL}_3$-web $W_r'$ with all weights being $1$. So, on $f_\beta$, there are ${\rm e}_{{\rm in},f_\beta}'$ source external vertices of $\ell'_r$ and ${\rm e}_{{\rm out},f_\beta}'$ sink external vertices of $\ell'_r$. By construction of the ${\rm SL}_3$-laminations $\ell'_r$ and $\ell'_t$, one has the compatibility ${\rm e}_{{\rm in},e_\alpha}' = {\rm e}_{{\rm out},f_\beta}'$ and ${\rm e}_{{\rm out},e_\alpha}' = {\rm e}_{{\rm in},f_\beta}'$ at the common edge of $t$ and $r$. Then, by Lem.\ref{lem:minimal_crossbar_A2-web_determined_by_signature}(MC3), we can fill in the biangle formed by $e_\alpha$ and $f_\beta$ by a (uniquely determined) minimal crossbar web. 

\vs

Gluing the canonical ${\rm SL}_3$-webs $W_t'$ in $\wh{t}$ for all triangles $\wh{t}$ of $\wh{\Delta}$ and the minimal crossbar webs for all biangles of $\wh{\Delta}$, we obtain a crossingless ${\rm SL}_3$-web $W'$, without boundary $1$-gon or $2$-gon. Replace each \raisebox{-0.4\height}{\scalebox{0.8}{
\begingroup%
  \makeatletter%
  \providecommand\color[2][]{%
    \errmessage{(Inkscape) Color is used for the text in Inkscape, but the package 'color.sty' is not loaded}%
    \renewcommand\color[2][]{}%
  }%
  \providecommand\transparent[1]{%
    \errmessage{(Inkscape) Transparency is used (non-zero) for the text in Inkscape, but the package 'transparent.sty' is not loaded}%
    \renewcommand\transparent[1]{}%
  }%
  \providecommand\rotatebox[2]{#2}%
  \newcommand*\fsize{\dimexpr\f@size pt\relax}%
  \newcommand*\lineheight[1]{\fontsize{\fsize}{#1\fsize}\selectfont}%
  \ifx\svgwidth\undefined%
    \setlength{\unitlength}{36.8503937bp}%
    \ifx\svgscale\undefined%
      \relax%
    \else%
      \setlength{\unitlength}{\unitlength * \real{\svgscale}}%
    \fi%
  \else%
    \setlength{\unitlength}{\svgwidth}%
  \fi%
  \global\let\svgwidth\undefined%
  \global\let\svgscale\undefined%
  \makeatother%
  \begin{picture}(1,1.03846154)%
    \lineheight{1}%
    \setlength\tabcolsep{0pt}%
    \put(0,0){\includegraphics[width=\unitlength,page=1]{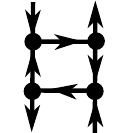}}%
  \end{picture}%
\endgroup%
}} by \raisebox{-0.4\height}{\scalebox{0.8}{
\begingroup%
  \makeatletter%
  \providecommand\color[2][]{%
    \errmessage{(Inkscape) Color is used for the text in Inkscape, but the package 'color.sty' is not loaded}%
    \renewcommand\color[2][]{}%
  }%
  \providecommand\transparent[1]{%
    \errmessage{(Inkscape) Transparency is used (non-zero) for the text in Inkscape, but the package 'transparent.sty' is not loaded}%
    \renewcommand\transparent[1]{}%
  }%
  \providecommand\rotatebox[2]{#2}%
  \newcommand*\fsize{\dimexpr\f@size pt\relax}%
  \newcommand*\lineheight[1]{\fontsize{\fsize}{#1\fsize}\selectfont}%
  \ifx\svgwidth\undefined%
    \setlength{\unitlength}{36.8503937bp}%
    \ifx\svgscale\undefined%
      \relax%
    \else%
      \setlength{\unitlength}{\unitlength * \real{\svgscale}}%
    \fi%
  \else%
    \setlength{\unitlength}{\svgwidth}%
  \fi%
  \global\let\svgwidth\undefined%
  \global\let\svgscale\undefined%
  \makeatother%
  \begin{picture}(1,1.03846154)%
    \lineheight{1}%
    \setlength\tabcolsep{0pt}%
    \put(0,0){\includegraphics[width=\unitlength,page=1]{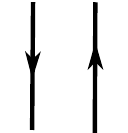}}%
  \end{picture}%
\endgroup%
}} and each \raisebox{-0.4\height}{\scalebox{0.8}{
\begingroup%
  \makeatletter%
  \providecommand\color[2][]{%
    \errmessage{(Inkscape) Color is used for the text in Inkscape, but the package 'color.sty' is not loaded}%
    \renewcommand\color[2][]{}%
  }%
  \providecommand\transparent[1]{%
    \errmessage{(Inkscape) Transparency is used (non-zero) for the text in Inkscape, but the package 'transparent.sty' is not loaded}%
    \renewcommand\transparent[1]{}%
  }%
  \providecommand\rotatebox[2]{#2}%
  \newcommand*\fsize{\dimexpr\f@size pt\relax}%
  \newcommand*\lineheight[1]{\fontsize{\fsize}{#1\fsize}\selectfont}%
  \ifx\svgwidth\undefined%
    \setlength{\unitlength}{42.51968504bp}%
    \ifx\svgscale\undefined%
      \relax%
    \else%
      \setlength{\unitlength}{\unitlength * \real{\svgscale}}%
    \fi%
  \else%
    \setlength{\unitlength}{\svgwidth}%
  \fi%
  \global\let\svgwidth\undefined%
  \global\let\svgscale\undefined%
  \makeatother%
  \begin{picture}(1,0.7)%
    \lineheight{1}%
    \setlength\tabcolsep{0pt}%
    \put(0,0){\includegraphics[width=\unitlength,page=1]{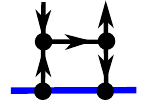}}%
  \end{picture}%
\endgroup%
}} by \raisebox{-0.4\height}{\scalebox{0.8}{
\begingroup%
  \makeatletter%
  \providecommand\color[2][]{%
    \errmessage{(Inkscape) Color is used for the text in Inkscape, but the package 'color.sty' is not loaded}%
    \renewcommand\color[2][]{}%
  }%
  \providecommand\transparent[1]{%
    \errmessage{(Inkscape) Transparency is used (non-zero) for the text in Inkscape, but the package 'transparent.sty' is not loaded}%
    \renewcommand\transparent[1]{}%
  }%
  \providecommand\rotatebox[2]{#2}%
  \newcommand*\fsize{\dimexpr\f@size pt\relax}%
  \newcommand*\lineheight[1]{\fontsize{\fsize}{#1\fsize}\selectfont}%
  \ifx\svgwidth\undefined%
    \setlength{\unitlength}{42.51968504bp}%
    \ifx\svgscale\undefined%
      \relax%
    \else%
      \setlength{\unitlength}{\unitlength * \real{\svgscale}}%
    \fi%
  \else%
    \setlength{\unitlength}{\svgwidth}%
  \fi%
  \global\let\svgwidth\undefined%
  \global\let\svgscale\undefined%
  \makeatother%
  \begin{picture}(1,0.7)%
    \lineheight{1}%
    \setlength\tabcolsep{0pt}%
    \put(0,0){\includegraphics[width=\unitlength,page=1]{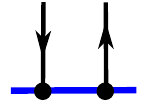}}%
  \end{picture}%
\endgroup%
}}, and do the same for the cases with all orientations reversed. Repeat until no more replacing is possible, so that the resulting ${\rm SL}_3$-web $W''$ is a reduced non-elliptic ${\rm SL}_3$-web in $\frak{S}$. One can observe that the internal 3-valent vertices appearing in this process must be in biangles, so the process removes some crossbars. Thus, in each triangle $\wh{t}$ of $\wh{\Delta}$, each step of such a process results only in exchanging positions of some corner arcs, hence the part in $\wh{t}$ is still canonical. For each biangle $B$ of $\wh{\Delta}$, one such process removes an `outermost' crossbar, and one can observe that the resulting picture is still a non-elliptic crossbar ${\rm SL}_3$-web in $B$, and hence is a minimal crossbar ${\rm SL}_3$-web. Thus, the final  ${\rm SL}_3$-web $W''$ is non-elliptic, reduced, and is canonical with respect to $\wh{\Delta}$, and for each triangle $\wh{t}$ of $\wh{\Delta}$, the ${\rm SL}_3$-webs $W'' \cap \wh{t}$ and $W'\cap \wh{t} = W_t' = \ell_t'$ in $\wh{t}$ have the same numbers of each kind of corner arcs and the same degree of pyramid. Thus the tropical coordinates of $W''$ form the vector $\vec{\rm a}'$. Let $\ell''$ be the ${\rm SL}_3$-lamination represented by the reduced non-elliptic ${\rm SL}_3$-web $W''$ with weight $1$. Now let $\ell$ be the ${\rm SL}_3$-lamination defined as $\ell := \ell'' \cup \ell_{-\vec{k}}$, which makes sense by Lem.\ref{lem:disjoint_union_with_peripheral_curve}(1). Note that the tropical coordinates of $\ell$ form the vector $\vec{\rm a}' - \vec{\rm a}_{\vec{k}} = \vec{\rm a}$, as desired. This shows the surjectivity of the coordinate map.

\vs

Let's now show the injectivity of the coordinate map. Let $\ell_1$ and $\ell_2$ be ${\rm SL}_3$-laminations in $\frak{S}$ having same tropical coordinates. As discussed above, we can find some $\vec{k}$ such that the ${\rm SL}_3$-laminations $\ell'_1 := \ell_1 \cup \ell_{\vec{k}}$ and $\ell'_2 := \ell_2 \cup \ell_{\vec{k}}$ have non-negative corner weights, so that they can be represented as weighted ${\rm SL}_3$-webs with all weights being $1$; they can be viewed as reduced non-elliptic ${\rm SL}_3$-webs. The intersection weights and the rotation weights, which can be easily seen to be determined by the tropical coordinates, then coincide with Frohman-Sikora's intersection coordinates and rotation numbers. Thus from Prop.\ref{prop:FS_coordinates}(FS2) it follows that $\ell_1'$ and $\ell_2'$ are equivalent as ${\rm SL}_3$-webs. Hence $\ell_1' = \ell_2'$ as ${\rm SL}_3$-laminations. Thus by Lem.\ref{lem:disjoint_union_with_peripheral_curve}(2) we get $\ell_1 = \ell_1' \cup \ell_{-\vec{k}} = \ell_2' \cup \ell_{-\vec{k}} = \ell_2$ as ${\rm SL}_3$-laminations, finishing the proof of injectivity. \qed

\vs

As mentioned already, one of the favorable properties of our coordinates is the compatibility formulas under change of ideal triangulations.
\begin{proposition}[coordinate change formula for tropical coordinates]
\label{prop:coordinate_change_formulas}
Let $\Delta$ and $\Delta'$ be ideal triangulations of a triangulable generalized marked surface $\frak{S}$ related to each other by a flip at an edge. Let $\ell$ be an ${\rm SL}_3$-lamination in $\frak{S}$. The tropical coordinates $({\rm a}_v(\ell))_v \in (\mathbb{Z}_{\ge 0})^{\mathcal{V}(Q_\Delta)}$ and $({\rm a}_{v'}'(\ell))_{v'} \in (\mathbb{Z}_{\ge 0})^{\mathcal{V}(Q_{\Delta'})}$ of $W$ with respect to $\Delta$ and $\Delta'$ are related by the sequence of tropical $\mathscr{A}$-mutations with respect to the sequence of mutations associated to a flip, i.e. by the same coordinate change formulas as described in Prop.\ref{prop:DS_coordinate_change}.
\end{proposition}
\begin{remark}
Prop.\ref{prop:coordinate_change_formulas} is also stated in \cite[Cor.4.5]{DS2}, following the previous version of the present paper \cite[Prop.3.35]{Kim}.
\end{remark}

{\it Proof.} Core of a proof of this proposition is just the corresponding statement for the Douglas-Sun coordinates of ${\rm SL}_3$-webs, i.e. Prop.\ref{prop:DS_coordinate_change}. Indeed, by additivity (Lem.\ref{lem:additivity_of_coordinates}), it suffices to show the statement for ${\rm SL}_3$-laminations $\ell$ that \redfix{can} be represented as a single-component ${\rm SL}_3$-web $W$, with some weight. Let $W$ be a single-component reduced non-elliptic ${\rm SL}_3$-web, and for each integer $k$, define $k\ell$ \redfix{to} be the ${\rm SL}_3$-lamination represented by $W$ with weight $k$, whenever it can be defined. Then, in view of the definition of the tropical coordinates, it is easy to observe ${\rm a}_v(k\ell) = k {\rm a}_v(\ell)$ for all $v\in \mathcal{V}(Q_\Delta)$ and ${\rm a}'_{v'}(k\ell) = k {\rm a}'_{v'}(\ell)$ for all $a' \in \mathcal{V}(Q_{\Delta'})$. Note ${\rm a}_v(\ell)$ and ${\rm a}'_{v'}(\ell)$ coincide with the Douglas-Sun coordinates of the ${\rm SL}_3$-web $W$ (Lem.\ref{lem:ours_are_compatible_with_DS}), and they transform as asserted, by Prop.\ref{prop:DS_coordinate_change}. The transformation formulas as presented in Prop.\ref{prop:DS_coordinate_change} are equivariant under the multiplication action by $\mathbb{Z}$, hence ${\rm a}_v(k\ell)$'s and ${\rm a}'_{v'}(k\ell)$'s also transform as wanted. \qed

\begin{remark}
This compatibility with tropical $\mathscr{A}$-mutation formulas is how we found the definition of our coordinates, up to scalar. Namely, we verified that those particular linear combinations of Frohman-Sikora coordinates enjoy these coordinate change formulas, at least for some simple cases. 
\end{remark}

Note
$$
\textstyle \mathbb{Z}^{\redfix{\mathcal{V}}(Q_\Delta)} \subset \mathcal{B}_\Delta \subset (\frac{1}{3}\mathbb{Z})^{\redfix{\mathcal{V}}(Q_\Delta)}.
$$
What will eventually play a major role are the ${\rm SL}_3$-laminations whose tropical coordinates lie in $\mathbb{Z}^{\redfix{\mathcal{V}}(Q_\Delta)}$. One consequence of Prop.\ref{prop:coordinate_change_formulas} is that, if all tropical coordinates of an ${\rm SL}_3$-lamination $\ell$ with respect to some ideal triangulation $\Delta$ are integers, then so are those of $\ell$ with respect to {\em any} ideal triangulation.
\begin{definition}
\label{def:congruent}
An ${\rm SL}_3$-lamination in a triangulable generalized marked surface $\frak{S}$ is said to be \ul{\em congruent} if for some, hence for every, ideal triangulation $\Delta$ of $\frak{S}$, its tropical coordinates are all integers.
\end{definition}
So we have a bijection
$$
{\bf a}_\Delta : \{ \mbox{congruent ${\rm SL}_3$-laminations in $\frak{S}$}\} \to \mathbb{Z}^{\redfix{\mathcal{V}}(Q_\Delta)} \cong \mathscr{A}_{|Q_\Delta|}(\mathbb{Z}^{\redfix{T}}).
$$
which is compatible under the tropical $\mathscr{A}$-mutations; see \S\ref{subsec:cluster_atlases} for $\mathscr{A}_{|Q_\Delta|}(\mathbb{Z}^{\redfix{T}})$. Hence the set of all congruent ${\rm SL}_3$-laminations in $\frak{S}$ works as a geometric model of $\mathscr{A}_{{\rm SL}_3,\frak{S}}(\mathbb{Z}^{\redfix{T}})$, the set of tropical integer points of the moduli space $\mathscr{A}_{{\rm SL}_3,\frak{S}}$, or that of the corresponding cluster $\mathscr{A}$-variety. 
\begin{theorem}
\label{thm:geometric_model}
For a triangulable generalized marked surface $\frak{S}$, we have a geometric model of $\mathscr{A}_{{\rm SL}_3,\frak{S}}(\mathbb{Z}^{\redfix{T}})$, the set of tropical integer points of the moduli space $\mathscr{A}_{{\rm SL}_3,\frak{S}}$, or that of the corresponding cluster $\mathscr{A}$-variety:
\begin{align}
\label{eq:geometric_model}
\mathscr{A}_{{\rm SL}_3, \frak{S}}(\mathbb{Z}^{\redfix{T}}) \leftrightarrow \{ \mbox{congruent ${\rm SL}_3$-laminations in $\frak{S}$}\}.
\end{align}
\end{theorem}

\vs

We suggest the readers to compare our model with previously proposed models of Le \cite{Ian} (`higher' laminations) and Goncharov-Shen \cite{GS15} (top-dimensional components of `surface affine Grassmannian' stack). 

\vs

Note that the assertion that our congruent ${\rm SL}_3$-laminations indeed provides a model of $\mathscr{A}_{{\rm SL}_3,\frak{S}}(\mathbb{Z}^{\redfix{T}})$ depends on Prop.\ref{prop:coordinate_change_formulas}, which in turn heavily depends on Prop.\ref{prop:DS_coordinate_change}, which is a result in \cite{DS2}. However, as a corollary of one of our main results, Thm.\ref{thm:main} (which is algebraic), whose proof does not depend on the validity of Prop.\ref{prop:coordinate_change_formulas} or Prop.\ref{prop:DS_coordinate_change}, we will provide a self-contained proof of a weaker version of Thm.\ref{thm:geometric_model}.
\begin{definition}
\label{def:Delta-congruent}
Let $\frak{S}$ be a triangulable generalized marked surface. For an ideal triangulation $\Delta$ of $\frak{S}$, we say that an ${\rm SL}_3$-lamination $\ell \in \mathscr{A}_{\rm L}(\frak{S};\mathbb{Z})$ in $\frak{S}$ is \ul{\em $\Delta$-congruent} if all tropical coordinates of $\ell$ for $\Delta$ are integers, i.e. ${\rm a}_v(\ell) \in \mathbb{Z}$, $\forall v\in \mathcal{V}(Q_\Delta)$. Let
$$
\mathscr{A}_\Delta(\mathbb{Z}^{\redfix{T}}) := \{ \mbox{$\Delta$-congruent ${\rm SL}_3$-laminations in $\frak{S}$} \} ~\subset~\mathscr{A}_{\rm L}(\frak{S};\mathbb{Z}).
$$
\end{definition}
\begin{proposition}[congruence condition is independent on triangulation]
\label{prop:congruence_condition_is_indepdent_on_triangulation}
For any ideal triangulations $\Delta$ and $\Delta'$ of a triangulable punctured surface $\frak{S}$, we have $\mathscr{A}_\Delta(\mathbb{Z}^{\redfix{T}}) = \mathscr{A}_{\Delta'}(\mathbb{Z}^{\redfix{T}})$.
\end{proposition}
In particular, this Prop.\ref{prop:congruence_condition_is_indepdent_on_triangulation} would justify Def.\ref{def:congruent}, and also eq.\eqref{eq:geometric_model} of Thm.\ref{thm:geometric_model}, for punctured surfaces. In the next section, Prop.\ref{prop:congruence_condition_is_indepdent_on_triangulation} will be proved (without Prop.\ref{prop:coordinate_change_formulas}) only at the end, so until then, we will mostly use the notion $\mathscr{A}_\Delta(\mathbb{Z}^{\redfix{T}})$ instead of $\mathscr{A}_{{\rm SL}_3,\frak{S}}(\mathbb{Z}^{\redfix{T}})$, to make the paper more independent.

\section{Regular functions on moduli spaces}
\label{sec:regular_functions_on_moduli_spaces}

One of the original Fock-Goncharov's duality conjectures \cite{FG06} is on the existence of a basis of the ring ${\bf L}(\mathscr{X}_{{\rm PGL}_3,\frak{S}})$ (Def.\ref{def:intro_L}) enumerated by $\mathscr{A}_{{\rm SL}_3,\frak{S}}(\mathbb{Z}^{\redfix{T}})$. We will construct a map $\mathscr{A}_{{\rm SL}_3,\frak{S}}(\mathbb{Z}^{\redfix{T}}) \to {\bf L}(\mathscr{X}_{{\rm PGL}_3,\frak{S}})$. By mimicking Fock-Goncharov's argument \cite{FG06} for ${\rm SL}_2$ and ${\rm PGL}_2$, we show that the image of this map is a basis of $\mathscr{O}(\mathscr{X}_{{\rm PGL}_3,\frak{S}})$. To do that, we investigate the ring $\mathscr{O}(\mathscr{L}_{{\rm PGL}_3,\frak{S}})$, and its relationship with $\mathscr{O}(\mathscr{X}_{{\rm PGL}_3,\frak{S}})$. We will observe that $\mathscr{O}(\mathscr{X}_{{\rm PGL}_3,\frak{S}})$ coincides with $\mathscr{O}_{\rm cl}(\mathscr{X}_{{\rm PGL}_3,\frak{S}}) =  \mathscr{O}(\mathscr{X}_{|Q_\Delta|})$.  We also investigate favorable properties of our duality map.  Throughout this section, $\frak{S} = \Sigma\setminus\mathcal{P}$ is a triangulable punctured surface. 

For a stack or a scheme $\mathscr{S}$, we denote by $\mathscr{O}(\mathscr{S})$ the ring of all its regular functions.

\subsection{Bases of rings of regular functions on ${\rm SL}_3$-moduli spaces}
\label{subsec:bases_of_rings}

We begin with $\mathscr{O}(\mathscr{L}_{{\rm SL}_3,\frak{S}})$, the ring of regular functions on the moduli space $\mathscr{L}_{{\rm SL}_3,\frak{S}}$ of ${\rm SL}_3$-local systems on the punctured surface $\frak{S}$, where ${\rm SL}_3$ is viewed as a scheme over $\mathbb{Q}$. \redfix{Note that $\mathscr{O}(\mathscr{L}_{{\rm SL}_3,\frak{S}})$ is the ring of invariants $A^{{\rm SL}_3}$, where $A$ is the coordinate ring of the affine variety ${\rm Hom}(\pi_1(\frak{S}),{\rm SL}_3)$, on which ${\rm SL}_3$ acts by conjugation.} Here is a standard element \redfix{of} $\mathscr{O}(\mathscr{L}_{{\rm SL}_3,\frak{S}})$\redfix{:}
\begin{definition}[trace of monodromy on $\mathscr{L}_{{\rm SL}_3,\frak{S}}$]
\label{def:trace-of-monodromy}
Let $\gamma$ be an oriented loop in $\frak{S}$. Denote by $f_\gamma$ the function on $\mathscr{L}_{{\rm SL}_3,\frak{S}}$ given by the trace of monodromy along $\gamma$. That is, for an ${\rm SL}_3$-local system $\mathcal{L}$ on $\frak{S}$, if $\rho:\pi_1(\frak{S}) \to {\rm SL}_3$ is the monodromy representation of $\mathcal{L}$ defined up to conjugation, define 
$$
f_\gamma(\mathcal{L}) := {\rm tr}(\rho([\gamma])).
$$
\end{definition}
It is easy to see that $f_\gamma$ is a well-defined regular function on $\mathscr{L}_{{\rm SL}_3,\frak{S}}$. It is known from Procesi \cite{P76} that $\mathscr{O}(\mathscr{L}_{{\rm SL}_3,\frak{S}})$ is generated by these trace-of-monodromy functions along loops. Sikora \cite{S01} found a complete set of relations among the trace-of-monodromy functions, and thus obtained an algebra isomorphism between the ${\rm SL}_3$-skein algebra $\mathcal{S}(\frak{S};\mathbb{Q})$ (Def.\ref{def:A2-skein_algebra}) and $\mathscr{O}(\mathscr{L}_{{\rm SL}_3,\frak{S}})$. 
Note that a single-component ${\rm SL}_3$-web $W$ in $\frak{S}$ with no internal or external vertices (Def.\ref{def:A2-web}) is an oriented loop $\gamma$ in $\frak{S}$; let's denote this $W$ by $W_\gamma$.
\begin{proposition}[\cite{S01}] 
\label{prop:skein_algebra_isomorphism}
For a punctured surface $\frak{S}$, there is a unique isomorphism 
$$
\Phi : \mathcal{S}(\frak{S};\mathbb{Q}) \to \mathscr{O}(\mathscr{L}_{{\rm SL}_3,\frak{S}})
$$
that sends each ${\rm SL}_3$-skein $[W_\gamma]$ consisting of one oriented loop $\gamma$ to the trace-of-monodromy function $f_\gamma$.
\end{proposition}
Recall Prop.\ref{prop:basis_of_A2-skein_algebra}, the result of Sikora and Westbury \cite{SW}, saying that the non-elliptic ${\rm SL}_3$-webs form a basis of $\mathcal{S}(\frak{S};\mathbb{Q})$, and recall that the set of all non-elliptic ${\rm SL}_3$-webs is in bijection with the set $\mathscr{A}_{\rm L}^0(\frak{S};\mathbb{Z})$ of all (integral) ${\rm SL}_3$-laminations with non-negative weights (Lem.\ref{lem:non-negative_A2-laminations_and_non-elliptic_A2-webs}); recall also \redfix{Rem}.\ref{rmk:punctured_surface_non-elliptic}.
\begin{corollary}[$A_2$-bangles basis of $\mathscr{O}(\mathscr{L}_{{\rm SL}_3,\frak{S}})$]
\label{cor:A2-bangles_basis_for_L_SL3}
The above construction yields an injective map
$$
\mathbb{I}^0_{{\rm SL}_3} : \mathscr{A}_{\rm L}^0(\frak{S};\mathbb{Z}) \to \mathscr{O}(\mathscr{L}_{{\rm SL}_3,\frak{S}})
$$
whose image set forms a basis of $\mathscr{O}(\mathscr{L}_{{\rm SL}_3,\frak{S}})$, which we call the \ul{\em $A_2$-bangles basis} of $\mathscr{O}(\mathscr{L}_{{\rm SL}_3,\frak{S}})$.
\end{corollary}
As mentioned earlier for the ${\rm SL}_3$-skein algebra $\mathcal{S}(\frak{S};\mathbb{Q})$, we may think of this basis of $\mathscr{O}(\mathscr{L}_{{\rm SL}_3,\frak{S}})$ as an $A_2$ version of a bangles basis for the well known $A_1$ theory. We will also discuss the $A_2$ version of the so-called bracelets basis. We recall the meaning of bangles and bracelets.

\begin{definition}
\label{def:bangles_and_bracelets}
Let $\gamma$ be an oriented simple loop in $\frak{S}$, hence forming a single-component non-elliptic ${\rm SL}_3$-web $W_\gamma$, and therefore a single-component ${\rm SL}_3$-lamination with weight $1$. Let $k \in \mathbb{Z}_{>0}$.
\begin{enumerate}
\item[\rm (1)] Define a \ul{\em $k$-bangle} $W_\gamma^k$ of $W_\gamma$ as a non-elliptic ${\rm SL}_3$-web consisting of $k$ copies of mutually disjoint oriented loops isotopic to $\gamma$. 
\item[\rm (2)] Define a \ul{\em $k$-bracelet} $W_\gamma^{(k)}$ as a single-component ${\rm SL}_3$-web obtained from the loop $\gamma^k = \gamma.\gamma.\cdots.\gamma$ by deforming it by a homotopy so that its self-intersections are transverse double.
\end{enumerate}
\end{definition}
When $W_\gamma$ is given as an element of the ${\rm SL}_3$-skein algebra $\mathcal{S}(\frak{S};\mathbb{Q})$, note that $W_\gamma^k$ and $W_\gamma^{(k)}$ yield well-defined elements $[W^k_\gamma]$ and $[W^{(k)}_\gamma]$ of $\mathcal{S}(\frak{S};\mathbb{Q})$. The notation for $k$-bangle is instructive, since $[W_\gamma^k] = [W_\gamma]^k$, with respect to the product structure of $\mathcal{S}(\frak{S};\mathbb{Q})$. By construction, (for an oriented simple loop $\gamma$) the $k$-bangle $W_\gamma^k$ can be viewed as an ${\rm SL}_3$-lamination, and we have
$$
\mathbb{I}^0_{{\rm SL}_3}(W^k_\gamma) = (\mathbb{I}^0_{{\rm SL}_3}(W_\gamma))^k = (f_\gamma)^k,
$$
where $f_\gamma$ is the trace-of-monodromy function along $\gamma$ (Def.\ref{def:trace-of-monodromy}); note that the ${\rm SL}_3$-lamination $W_\gamma^k$ can be represented by one component $W_\gamma$ with weight $k$. We will consider the bracelets version in \redfix{\S\ref{sec:conjectures}}. 

\vs

On the other hand, from the defining relations and the product structure of $\mathcal{S}(\frak{S};\mathbb{Q})$, we immediately get the following useful result.
\begin{corollary}
\label{cor:structure_constants_SL3_0}
The structure constants of the basis $\mathbb{I}^0_{{\rm SL}_3}(\mathscr{A}^0_{\rm L}(\frak{S};\mathbb{Z}))$ of $\mathscr{O}(\mathscr{L}_{{\rm SL}_3,\frak{S}})$ are integers. That is, for any $\ell,\ell' \in \mathscr{A}^0_{\rm L}(\frak{S};\mathbb{Z})$, we have
$$
\mathbb{I}^0_{{\rm SL}_3}(\ell) \, \mathbb{I}^0_{{\rm SL}_3}(\ell') = \underset{\ell'' \in \mathscr{A}_{\rm L}^0(\frak{S};\mathbb{Z})}{\textstyle \sum} c_{{\rm SL}_3}^0(\ell,\ell';\ell'') \, \mathbb{I}^0_{{\rm SL}_3}(\ell'')
$$
where $c_{{\rm SL}_3}^0(\ell,\ell';\ell'')\in \mathbb{Z}$ and $c_{{\rm SL}_3}^0(\ell,\ell';\ell'')$ are zero for all but at most finitely many $\ell''$. \qed 
\end{corollary}

\vs

Next step is to consider $\mathscr{O}(\mathscr{X}_{{\rm SL}_3,\frak{S}})$. We recall the result from \cite[\S12.5]{FG06} showing that it is a free module over $\mathscr{O}(\mathscr{L}_{{\rm SL}_3,\frak{S}})$. By abuse of notation, let ${\rm B}$ be a Borel subgroup of ${\rm SL}_3$, say the subgroup of all upper triangular matrices in ${\rm SL}_3$, and let ${\rm U} := [{\rm B},{\rm B}]$ be the corresponding maximal unipotent subgroup, which would be the subgroup of all upper triangular matrices with diagonal entries being all $1$. Define the \ul{\em Cartan group} of ${\rm SL}_3$ as ${\rm H} := {\rm B}/{\rm U}$, which for our case is canonically isomorphic to the subgroup of ${\rm SL}_3$ of diagonal matrices, and the quotient map ${\rm B}\to{\rm B}/{\rm U}$ just reads the diagonal entries. Let ${\rm W}$ be the corresponding Weyl group. For each puncture $p\in \mathcal{P}$, there is a canonical map
\begin{align}
\label{eq:pi_p}
\pi_p : \mathscr{X}_{{\rm SL}_3,\frak{S}} \to {\rm H}
\end{align}
provided by the framing and the semi-simple part of the monodromy along a peripheral loop surrounding $p$. In fact, we need to choose the orientation of the loop carefully.
\begin{definition}
\label{def:positive_orientation}
A peripheral loop in $\frak{S}$ (Def.\ref{def:A2-lamination}) around a puncture $p\in \mathcal{P}$ is \ul{\em positively oriented} if it is isotopic to the hole of $\til{\frak{S}}$ corresponding to $p$, given the boundary-orientation (Def.\ref{def:holed_surface}). We say it is \ul{\em negatively oriented} otherwise.
\end{definition}
Now, given a framed ${\rm SL}_3$-local system $(\mathcal{L},\beta)$ on $\frak{S}$, for each puncture $p\in \mathcal{P}$, consider the monodromy along a positively oriented peripheral loop $\gamma$ around $p$. This monodromy is defined only up to conjugation in ${\rm SL}_3$, and lives in some Borel subgroup of ${\rm SL}_3$, hence can be thought of as living in our fixed choice ${\rm B}$. The semi-simple part of this element of ${\rm B}$ can be obtained as the image of the quotient map ${\rm B} \to {\rm B}/{\rm U} = {\rm H}$. As said in \cite{FG06}, the semi-simple part of the monodromy alone yields an element of ${\rm H}/{\rm W}$, giving a map
$$
\mathscr{L}_{{\rm SL}_3,\frak{S}} \dashrightarrow {\rm H}/{\rm W}
$$
\redfix{defined on the locus of ${\rm SL}_3$-local systems with semi-simple monodromy around punctures,} and together with the framing data we get the map $\pi_p : \mathscr{X}_{{\rm SL}_3,\frak{S}} \to {\rm H}$. Let us give a more precise explanation of $\pi_p$ as it is important in the present paper, but is not described in \cite{FG06} in detail. Consider $p$ as a point of the hole of $\til{\frak{S}}$, hence in particular $p$ can be though\redfix{t} of as a point of $\frak{S}$, as in Def.\ref{def:holed_surface} and the discussion after that. So the framing $\beta$ yields a distinguished point of the fiber $(\mathcal{L}_\mathcal{B})_p$ of the associated flag bundle $\mathcal{L}_\mathcal{B}$; recall $(\mathcal{L}_\mathcal{B})_p = \mathcal{L}_p \times_G {\rm G}/{\rm B} = \{ [ v, g{\rm B}] \, | \, v\in \mathcal{L}_p, g\in {\rm G}\}$, where $[v,g{\rm B}] = [v',g'{\rm B}]$ iff $v'g' = vg b$ for some $b\in {\rm B}$. In particular, $[v,g{\rm B}] = [vg, {\rm B}] = [vgb,{\rm B}]$. Hence we can write the distinguished element of $(\mathcal{L}_\mathcal{B})_p$ assigned by $\beta$ as $\beta(p) = [v_0,{\rm B}]$ for some $v_0\in \mathcal{L}_p$ that is  uniquely determined up to right action of ${\rm B}$. The parallel transport map of $\mathcal{L}$ along $\gamma$ gives the monodromy map $\Pi_\gamma : \mathcal{L}_p \to \mathcal{L}_p$ that is equivariant under the right ${\rm G}$-actions. The induced monodromy $(\Pi_\mathcal{B})_\gamma : (\mathcal{L}_\beta)_p \to (\mathcal{L}_\beta)_p$ for $\mathcal{L}_\mathcal{B}$ then sends $[v,g{\rm B}]$ to $[\Pi_\gamma(v), g{\rm B}]$\redfix{. Since} $\beta$ is a flat\redfix{, or covariantly-constant,} section\redfix{, the values of $\beta$ at points are related by the parallel transport maps; in particular}, we have $[v_0,{\rm B}] = [\Pi_\gamma(v_0), {\rm B}]$. This means $\Pi_\gamma(v_0) =v_0 b_0$ for some $b_0\in {\rm B}$ which is unique determined by $v_0$. If $\beta(p)=[v_0',{\rm B}]$, then $v_0' = v_0b$ for some $b\in {\rm B}$, then $\Pi_\gamma(v_0')=\Pi_\gamma(v_0b) = \Pi_\gamma(v_0) b = v_0 b_0 b = v_0b(b^{-1} b_0 b) = v_0'(b^{-1} b_0 b)$, so $b_0' = b^{-1} b_0 b$ makes $\Pi_\gamma(v_0') = v_0'b_0'$. This means that out of the monodomy of $(\mathcal{L},\beta)$ along $\gamma$ we get an element of ${\rm B}$ uniquely determined up to conjugation by an element of ${\rm B}$. But $b_0$ and $b^{-1} b_0 b$ have same semi-simple parts (i.e. they have same diagonal entries), i.e. we get a well-defined element of ${\rm B}/{\rm U} = {\rm H}$. Thus we get the well-defined map $\pi_p : \mathscr{X}_{{\rm SL}_3,\frak{S}} \to {\rm H}$. 

\vs

The maps $\pi_p$ for all $p\in \mathcal{P}$ constitute the map
$$
\pi : \mathscr{X}_{{\rm SL}_3,\frak{S}} \to {\rm H}^\mathcal{P}
$$
where $\mathcal{P}$ is the set of all punctures of $\frak{S}$; here ${\rm H}^{\mathcal{P}}$ may be understood as ${\rm H}^{|\mathcal{P}|}$.  Likewise, the maps $\mathscr{L}_{{\rm SL}_3,\frak{S}} \dashrightarrow {\rm H}/{\rm W}$ for punctures $p\in \mathcal{P}$ constitute the map $\mathscr{L}_{{\rm SL}_3,\frak{S}} \dashrightarrow ({\rm H}/{\rm W})^\mathcal{P}$. Fock and Goncharov \cite[\S12.5]{FG06} state that we have a Cartesian square of stacks
$$
\xymatrix@R-3mm{
\mathscr{X}_{{\rm SL}_3,\frak{S}} \ar[r] \ar[d]_F & {\rm H}^\mathcal{P} \ar[d] \\
\mathscr{L}_{{\rm SL}_3,\frak{S}} \ar[r] & ({\rm H}/{\rm W})^\mathcal{P}
}
$$
where the left vertical arrow is the forgetting-of-framing map
$$
F : \mathscr{X}_{{\rm SL}_3,\frak{S}} \to \mathscr{L}_{{\rm SL}_3,\frak{S}}
$$
and that therefore
\begin{align}
\label{eq:O_of_X_SL3_and_O_of_L_SL3}
\mathscr{O}(\mathscr{X}_{{\rm SL}_3,\frak{S}}) = \mathscr{O}(\mathscr{L}_{{\rm SL}_3,\frak{S}}) \otimes_{\mathscr{O}(({\rm H}/{\rm W})^\mathcal{P})} \mathscr{O}({\rm H}^\mathcal{P}).
\end{align}
holds. They also state that, by Chevalley's theorem $\mathbb{Q}[{\rm H}]$ is a free $\mathbb{Q}[{\rm H}]^{\rm W}$-module with $|{\rm W}|$ generators, and hence $\mathscr{O}(\mathscr{X}_{{\rm SL}_3,\frak{S}})$ is a free $\mathscr{O}(\mathscr{L}_{{\rm SL}_3,\frak{S}})$-module of rank $|{\rm W}|^{|\mathcal{P}|}$, where a set of generators is obtained by pullbacks of generators of the $\mathbb{Q}[({\rm H}/{\rm W})^\mathcal{P}]$-module $\mathbb{Q}[{\rm H}]^\mathcal{P}$.

\vs

For $i=1,2,3$, composing the projection ${\rm H} \to \mathbb{G}_m$ to the $i$-th diagonal entry with the map $\pi_p : \mathscr{X}_{{\rm SL}_3,\frak{S}} \to {\rm H}$ yields a regular map
\begin{align}
\label{eq:pi_p_i}
(\pi_p)_i : \mathscr{X}_{{\rm SL}_3,\frak{S}} \to \mathbb{G}_m,
\end{align}
that is, we get $(\pi_p)_i \in \mathscr{O}(\mathscr{X}_{{\rm SL}_3,\frak{S}})$.

\vs

Note that the map $F : \mathscr{X}_{{\rm SL}_3,\frak{S}} \to \mathscr{L}_{{\rm SL}_3,\frak{S}}$ is regular, hence induces
\begin{align}
\label{eq:F_star}
F^* : \mathscr{O}(\mathscr{L}_{{\rm SL}_3,\frak{S}}) \to \mathscr{O}(\mathscr{X}_{{\rm SL}_3,\frak{S}}).
\end{align}
Since the image of $F$ is open and dense in $\mathscr{L}_{{\rm SL}_3,\frak{S}}$, it follows that $F^*$ is injective.

\vs

We now arrive at a canonical basis of $\mathscr{O}(\mathscr{X}_{{\rm SL}_3,\frak{S}})$.
\begin{definition}[$A_2$-bangles basis of $\mathscr{O}(\mathscr{X}_{{\rm SL}_3,\frak{S}})$]
\label{def:I_SL3}
For a punctured surface $\frak{S}$, define a map 
$$
\mathbb{I}_{{\rm SL}_3} : \mathscr{A}_{\rm L}(\frak{S};\mathbb{Z}) \to \mathscr{O}(\mathscr{X}_{{\rm SL}_3,\frak{S}})
$$
as follows. Let $\ell \in \mathscr{A}_{\rm L}(\frak{S};\mathbb{Z})$. Represent $\ell$ as disjoint union $\ell = \ell_1 \cup \ell_2 \cup \cdots \cup \ell_n$ (Def.\ref{def:disjoint}) of single-component ${\rm SL}_3$-laminations $\ell_1,\ldots,\ell_n$, whose underlying non-elliptic ${\rm SL}_3$-webs are mutually non-isotopic. Define $\mathbb{I}_{{\rm SL}_3}(\ell_i)$ as:
\begin{enumerate}
\itemsep0em
\item[\rm (CB1)] If $\ell_i$ consists of a peripheral  loop $\gamma_i$ with weight $k_i \in \mathbb{Z}$, surrounding a puncture $p \in \mathcal{P}$,
\begin{enumerate}
\itemsep0em
\item[\rm (CB1-1)] if $\gamma_i$ is positively oriented (Def.\ref{def:positive_orientation}), then
$$
\mathbb{I}_{{\rm SL}_3}(\ell_i) := ( (\pi_p)_1 )^{k_i};
$$

\item[\rm (CB1-2)] if $\gamma_i$ is negatively oriented (Def.\ref{def:positive_orientation}), then
$$
\mathbb{I}_{{\rm SL}_3}(\ell_i) := ( (\pi_p)_3 )^{-k_i};
$$
\end{enumerate}

\item[\rm (CB2)] Otherwise, define
$$
\mathbb{I}_{{\rm SL}_3}(\ell_i) = F^* \mathbb{I}_{{\rm SL}_3}^0(\ell_i).
$$
\end{enumerate}
Define
$$
\mathbb{I}_{{\rm SL}_3}(\ell) := \mathbb{I}_{{\rm SL}_3}(\ell_1) \, \mathbb{I}_{{\rm SL}_3}(\ell_2) \,\cdot \cdots \cdot\, \mathbb{I}_{{\rm SL}_3}(\ell_n).
$$
By convention, we set $\mathbb{I}_{{\rm SL}_3}({\O}):=1$.

\vs

The image set $\mathbb{I}_{{\rm SL}_3}(\mathscr{A}_{\rm L}(\frak{S};\mathbb{Z}))$ is called the \ul{\em $A_2$-bangles basis} of $\mathscr{O}(\mathscr{X}_{{\rm SL}_3,\frak{S}})$, by a slight abuse of notation.
\end{definition}

\begin{proposition}
\label{prop:I_SL3_is_a_basis}
For a punctured surface $\frak{S}$, one has:
\begin{enumerate}
\itemsep0em
\item[\rm (1)] The map $\mathbb{I}_{{\rm SL}_3}$ is injective, and the image set of $\mathbb{I}_{{\rm SL}_3}$ is indeed a basis of $\mathscr{O}(\mathscr{X}_{{\rm SL}_3,\frak{S}})$.

\vs

\item[\rm (2)] The structure constants of this $A_2$-bangles basis of $\mathscr{O}(\mathscr{X}_{{\rm SL}_3,\frak{S}})$ are integers. That is, for any $\ell,\ell' \in \mathscr{A}_{\rm L}(\frak{S};\mathbb{Z})$, we have
\begin{align}
\label{eq:I_SL3_structure_constants}
\mathbb{I}_{{\rm SL}_3}(\ell) \, \mathbb{I}_{{\rm SL}_3}(\ell') = \underset{\ell'' \in \mathscr{A}_{\rm L}(\frak{S};\mathbb{Z})}{\textstyle \sum} c_{{\rm SL}_3}(\ell,\ell';\ell'') \, \mathbb{I}_{{\rm SL}_3}(\ell'')
\end{align}
where $c_{{\rm SL}_3}(\ell,\ell';\ell'')\in \mathbb{Z}$ and $c_{{\rm SL}_3}(\ell,\ell';\ell'')$ are zero for all but at most finitely many $\ell''$.
\end{enumerate}
\end{proposition}

{\it Proof.} (1) We first consider the restriction of $\mathbb{I}_{{\rm SL}_3}$ to the subset $\mathscr{A}^0_{\rm L}(\frak{S};\mathbb{Z})$ of $\mathscr{A}_{\rm L}(\frak{S};\mathbb{Z})$ consisting of ${\rm SL}_3$-laminations with non-negative weights. This restricted map doesn't exactly equal $F^* \circ \mathbb{I}^0_{{\rm SL}_3}$ because of the peripheral loops. Let $\ell_i$ be a single-component ${\rm SL}_3$-lamination consisting of a single peripheral (simple) loop $\gamma_i$ around a puncture $p\in \mathcal{P}$, with weight $k_i \in \mathbb{Z}$. One may view $\gamma_i$ as the hole of $\til{\frak{S}}$ corresponding to $p$; assume that the orientation matches the boundary-orientation of the hole (Def.\ref{def:holed_surface}). As discussed before, the monodromy around $\gamma_i$ can be thought of as living in ${\rm B}$, i.e. being an upper triangular matrix. Recall that the map $\pi_p : \mathscr{X}_{{\rm SL}_3,\frak{S}} \to {\rm H}$ reads the diagonal part. Hence it follows that
$$
F^*\mathbb{I}^0_{{\rm SL}_3}(\ell_i) = F^* f_{\gamma_i^{k_i}} = ((\pi_p)_1)^{k_i} + ((\pi_p)_2)^{k_i} + ((\pi_p)_3)^{k_i}.
$$
By definition, one notes $\mathbb{I}_{{\rm SL}_3}(\ell_i) = ((\pi_p)_1)^{k_i}$. As in Lem.\ref{lem:disjoint_union_with_peripheral_curve}, $-\ell_i$ would denote a single-component ${\rm SL}_3$-lamination consisting of $\gamma_i$ with weight $-k_i$; then $\mathbb{I}_{{\rm SL}_3}(-\ell_i) =((\pi_p)_1)^{-k_i}$. On the other hand, denote by $\ol{\ell_i}$ the single-component ${\rm SL}_3$-lamination consisting of the peripheral loop $\ol{\gamma_i}$ with weight $k_i$, where $\ol{\gamma_i}$ is same as $\gamma_i$ with the orientation reversed. So, by definition, $\mathbb{I}_{{\rm SL}_3}(\ol{\ell_i}) = ((\pi_p)_3)^{-k_i}$, and $\mathbb{I}_{{\rm SL}_3}(-\ol{\ell_i}) = ((\pi_p)_3)^{k_i}$. From $(\pi_p)_1 (\pi_p)_2 (\pi_p)_3 = 1$ it follows that $((\pi_p)_2)^{k_i} = ((\pi_p)_1)^{-k_i}((\pi_p)_3)^{-k_i} = \mathbb{I}_{{\rm SL}_3}(-\ell_i) \mathbb{I}_{{\rm SL}_3}(\ol{\ell_i})$, which, in turn, by Lem.\ref{lem:disjoint_union_with_peripheral_curve}, equals $\mathbb{I}_{{\rm SL}_3}((-\ell_i)\cup \ol{\ell_i})$. To summarize,
\begin{align}
\label{eq:difference_for_peripheral_loop}
F^*\mathbb{I}^0_{{\rm SL}_3}(\ell_i)
= \mathbb{I}_{{\rm SL}_3}(\ell_i) + \mathbb{I}_{{\rm SL}_3}( (-\ell_i)\cup \ol{\ell_i}) + \mathbb{I}_{{\rm SL}_3}(-\ol{\ell_i})
\end{align}
when $\ell_i$ is a single-component ${\rm SL}_3$-lamination consisting of a peripheral loop, oriented according to the bounary-orientation along the corresponding hole of $\til{\frak{S}}$. Now suppose that $\ell_i$ is a single peripheral loop $\gamma_i$ with weight $k_i$, but $\gamma_i$ is negatively oriented (Def.\ref{def:positive_orientation}). Then, by definition of $\pi_p$, we have $F^*\mathbb{I}^0_{{\rm SL}_3}(\ell_i) = f_{\gamma_i^{k_i}} = ((\pi_p)_1)^{-k_i} + ((\pi_p)_2)^{-k_i} + ((\pi_p)_3)^{-k_i}$. This time, we can observe that $\mathbb{I}_{{\rm SL}_3}(\ell_i) = ((\pi_p)_3)^{-k_i}$ and $\mathbb{I}_{{\rm SL}_3}(-\ol{\ell_i}) = ((\pi_p)_1)^{-k_i}$, hence eq.\eqref{eq:difference_for_peripheral_loop} still holds.

\vs

On the other hand, if $\ell \in \mathscr{A}_{\rm L}(\frak{S},\mathbb{Z})$ does not contain any peripheral loop, then $\ell \in \mathscr{A}_{\rm L}^0(\frak{S};\mathbb{Z})$, and
\begin{align}
\label{eq:same_for_non-peripheral_loop}
F^*\mathbb{I}^0_{{\rm SL}_3}(\ell) = \mathbb{I}_{{\rm SL}_3}(\ell).
\end{align}
Note that $\mathbb{I}^0_{{\rm SL}_3}(\mathscr{A}^0_{\rm L}(\frak{S};\mathbb{Z}))$ spans $\mathscr{O}(\mathscr{L}_{{\rm SL}_3,\frak{S}})$ (Cor.\ref{cor:A2-bangles_basis_for_L_SL3}). We just saw that the set $F^*(\mathbb{I}^0_{{\rm SL}_3}(\mathscr{A}^0_{\rm L}(\frak{S};\mathbb{Z}))) \subset F^*(\mathscr{O}(\mathscr{L}_{{\rm SL}_3,\frak{S}})) \subset \mathscr{O}(\mathscr{X}_{{\rm SL}_3,\frak{S}})$ lies in the span of $\mathbb{I}_{{\rm SL}_3}(\mathscr{A}_{\rm L}(\frak{S};\mathbb{Z}))$. In view of eq.\eqref{eq:O_of_X_SL3_and_O_of_L_SL3}, elements of $F^*(\mathscr{O}(\mathscr{L}_{{\rm SL}_3,\frak{S}}))$ tensored with elements of $\mathscr{O}({\rm H}^\mathcal{P}) \cong \mathscr{O}({\rm H})^{\otimes |\mathcal{P}|}$ span $\mathscr{O}(\mathscr{X}_{{\rm SL}_3,\frak{S}})$. One copy of $\mathscr{O}({\rm H}) \cong \mathbb{Q}[a^{\pm 1}, b^{\pm 1}, c^{\pm 1}]/(abc-1)$ (where $a,b,c$ are coordinate functions of ${\rm H}$ for the diagonal entries) is associated to each puncture $p\in \mathcal{P}$, and by definition of $\pi_p$ and $(\pi_p)_i$, one can observe that the functions $(\pi_p)_i^k$, $i=1,2,3$, $k \in \mathbb{Z}$, span this copy of $\mathscr{O}({\rm H})$. By (CB1-1)--(CB1-2), $(\pi_p)_1^k$ and $(\pi_p)_3^k$ (for each $k\in \mathbb{Z}$) belong to $\mathbb{I}_{{\rm SL}_3}(\mathscr{A}_{\rm L}(\frak{S};\mathbb{Z}))$, and we saw above that $(\pi_p)_2^k$ also belongs to $\mathbb{I}_{{\rm SL}_3}(\mathscr{A}_{\rm L}(\frak{S};\mathbb{Z}))$. This shows that $\mathbb{I}_{{\rm SL}_3}(\mathscr{A}_{\rm L}(\frak{S};\mathbb{Z}))$ spans $\mathscr{O}(\mathscr{X}_{{\rm SL}_3,\frak{S}})$. We only sketch a proof for the linear independence of this set and the injectivity of $\mathbb{I}_{{\rm SL}_3}$, as we will not really use these facts; but we will definitely be using the spanning property. From the injectivity of $F^*$, the injectivity of $F^* \mathbb{I}^0_{{\rm SL}_3}$ and the linear independence of the set $F^*(\mathbb{I}^0_{{\rm SL}_3}(\mathscr{A}^0_{\rm L}(\frak{S};\mathbb{Z})))$ follow. One can explicitly write down this much result in terms of $\mathbb{I}_{{\rm SL}_3}$. The remaining is essentially the investigation of a basis of $\mathscr{O}({\rm H}) \cong \mathbb{Q}[a^{\pm 1}, b^{\pm 1}, c^{\pm 1}]/(abc-1)$; a non-redundant set of (all possible) Laurent monomials in $a,b,c$ will be a basis. For each puncture $p$, such a set is in bijection with the set of all distinct ${\rm SL}_3$-laminations consisting only of peripheral loops surrounding $p$; the redundancy relation $a^kb^kc^k=1$ is exactly captured by the fact that for a ${\rm SL}_3$-lamination $\ell$ with a single peripheral loop with weight $k$, the ${\rm SL}_3$-lamination $\ell \cup ( (-\ell) \cup \ol{\ell} ) \cup (-\ol{\ell})$ equals the empty ${\rm SL}_3$-lamination (Lem.\ref{lem:disjoint_union_with_peripheral_curve}(2)).

\vs

(2) Let's first establish a lemma, which is easily observed (with the help of Lem.\ref{lem:disjoint_union_with_peripheral_curve}):
\begin{lemma}
\label{lem:product_factors}
Let $\ell,\ell' \in \mathscr{A}_{\rm L}(\frak{S},\mathbb{Z})$. If $\ell$ and $\ell'$ are disjoint (Def.\ref{def:disjoint}), so that $\ell\cup \ell'$ makes sense as an ${\rm SL}_3$-lamination, then
$$
\mathbb{I}_{{\rm SL}_3}(\ell\cup \ell') = \mathbb{I}_{{\rm SL}_3}(\ell) \, \mathbb{I}_{{\rm SL}_3}(\ell').
$$
If furthermore $\ell,\ell'\in \mathscr{A}_{\rm L}^0(\frak{S},\mathbb{Z})$, then
$$
\mathbb{I}^0_{{\rm SL}_3}(\ell\cup \ell') = \mathbb{I}^0_{{\rm SL}_3}(\ell) \, \mathbb{I}^0_{{\rm SL}_3}(\ell').
$$
For example, if $\ell$ or $\ell'$ (both belonging to $\mathscr{A}_{\rm L}(\frak{S},\mathbb{Z})$ or \redfix{both} to $\mathscr{A}^0_{\rm L}(\frak{S},\mathbb{Z})$, respectively) consists only of peripheral loops, then $\ell$ and $\ell'$ are disjoint, and the above holds. \qed
\end{lemma}

Now let $\ell,\ell' \in \mathscr{A}_{\rm L}(\frak{S};\mathbb{Z})$. We can decompose them into disjoint unions as $\ell = \ell_1 \cup \ell_2$ and $\ell' = \ell'_1 \cup \ell_2'$, where each of $\ell_2$ and $\ell_2'$ is either empty or consists only of peripheral loops, while each of $\ell_1$ and $\ell_1'$ is either empty or does not contain any peripheral loop. In particular, $\ell_1,\ell_1' \in \mathscr{A}^0_{\rm L}(\frak{S};\mathbb{Z})$, and $F^*\mathbb{I}^0_{{\rm SL}_3}(\ell_1) = \mathbb{I}_{{\rm SL}_3}(\ell_1)$, $F^*\mathbb{I}^0_{{\rm SL}_3}(\ell_1') = \mathbb{I}_{{\rm SL}_3}(\ell_1')$. So
\begin{align}
\nonumber
\mathbb{I}_{{\rm SL}_3}(\ell) \, \mathbb{I}_{{\rm SL}_3}(\ell')
& = F^*( \mathbb{I}^0_{{\rm SL}_3}(\ell_1) \,\mathbb{I}^0_{{\rm SL}_3}(\ell_1') )\,\mathbb{I}_{{\rm SL}_3}(\ell_2) \, \mathbb{I}_{{\rm SL}_3}(\ell_2') \\
\nonumber
& = F^* \left( \textstyle \underset{\ell''\in \mathscr{A}^0_{\rm L}(\frak{S};\mathbb{Z})}{\sum} c^0_{{\rm SL}_3}(\ell_1,\ell_1';\ell'') \, \mathbb{I}^0_{{\rm SL}_3}(\ell'')  \,\right) \mathbb{I}_{{\rm SL}_3}(\ell_2) \, \mathbb{I}_{{\rm SL}_3}(\ell_2') \\
\label{eq:structure_constant_proof_equation}
& \textstyle = \underset{\ell''\in \mathscr{A}^0_{\rm L}(\frak{S};\mathbb{Z})}{\sum} c^0_{{\rm SL}_3}(\ell_1,\ell_1';\ell'')\, (F^*\mathbb{I}^0_{{\rm SL}_3})(\ell'') \, \mathbb{I}_{{\rm SL}_3}(\ell_2) \, \mathbb{I}_{{\rm SL}_3}(\ell_2').
\end{align}
Decompose $\ell'' \in \mathscr{A}_{\rm L}^0(\frak{S};\mathbb{Z})$ into disjoint union $\ell''_1 \cup \ell''_2$, where $\ell''_2$ consists only of peripheral loops and $\ell''_1$ has no peripheral loop. Then, as seen,
\begin{align*}
F^*\mathbb{I}^0_{{\rm SL}_3}(\ell'')
& = F^*\mathbb{I}^0_{{\rm SL}_3}(\ell''_1) \,
F^*\mathbb{I}^0_{{\rm SL}_3}(\ell''_2) \quad (\mbox{$\because$ Lem.\ref{lem:product_factors}}) \\
& = \mathbb{I}_{{\rm SL}_3}(\ell_1'') \, ( \mathbb{I}_{{\rm SL}_3}(\ell_2'') + \mathbb{I}_{{\rm SL}_3}((-\ell_2'')\cup \ol{\ell''_2}) + \mathbb{I}_{{\rm SL}_3}(-\ol{\ell_2''})) \quad (\because \mbox{eq.\eqref{eq:difference_for_peripheral_loop}--\eqref{eq:same_for_non-peripheral_loop}}) \\
& = \mathbb{I}_{{\rm SL}_3}(\ell_1''\cup \ell''_2) + \mathbb{I}_{{\rm SL}_3}(\ell_1'' \cup(-\ell_2'') \cup \ol{\ell''_2}) + \mathbb{I}_{{\rm SL}_3}(\ell_1'' \cup (-\ol{\ell_2''})) \quad (\mbox{$\because$ Lem.\ref{lem:product_factors}}).
\end{align*}
Putting into eq.\eqref{eq:structure_constant_proof_equation} and using Lem.\ref{lem:product_factors}, we obtain the desired statement for item (2). \qed

\subsection{Lifting ${\rm PGL}_3$ monodromies to ${\rm SL}_3$}
\label{subsec:lifting_PGL3_to_SL3}

Consider the natural regular map
$$
P : \mathscr{X}_{{\rm SL}_3,\frak{S}} \to \mathscr{X}_{{\rm PGL}_3,\frak{S}}
$$
induced by the natural quotient ${\rm SL}_3 \to {\rm PGL}_3$, yielding a map
\begin{align}
\label{eq:P_star}
P^* : \mathscr{O}(\mathscr{X}_{{\rm PGL}_3,\frak{S}}) \to \mathscr{O}(\mathscr{X}_{{\rm SL}_3,\frak{S}})\redfix{.}
\end{align}
In the previous subsection, we obtained a basis of $\mathscr{O}(\mathscr{X}_{{\rm SL}_3,\frak{S}})$. Now we have to figure out which elements of $\mathscr{O}(\mathscr{X}_{{\rm SL}_3,\frak{S}})$ belong to the image of $P^*$. Or, going in the other direction, given an element of $\mathscr{O}(\mathscr{X}_{{\rm PGL}_3,\frak{S}})$, the image of it under $P^*$ would be an element of $\mathscr{O}(\mathscr{X}_{{\rm SL}_3,\frak{S}})$, hence can be written as linear combination of elements of the $A_2$-bangles basis $\mathbb{I}_{{\rm SL}_3}(\mathscr{A}_{\rm L}(\frak{S},\mathbb{Z}))$ we obtained. Each $A_2$-bangles basis vector is a product of trace-of-monodromy functions along loops and certain functions associated to punctures. Now, for example, what kind of function on $\mathscr{X}_{{\rm PGL}_3,\frak{S}}$ should correspond to the trace-of-monodromy function on $\mathscr{X}_{{\rm SL}_3,\frak{S}}$? The monodromy for a point of $\mathscr{X}_{{\rm PGL}_3,\frak{S}}$ gives only a homomorphism $\pi_1(\frak{S}) \to {\rm PGL}_3$ (defined up to conjugation), hence the naive trace-of-monodromy along a loop is not well-defined (or, its value is defined in $\mathbb{A}^1$ only up to $\mathbb{G}_m$, which is not useful).

\vs

As an auxiliary device, we will make use of the set of positive real points, i.e. 
$$
\mathscr{X}_{{\rm PGL}_3,\frak{S}}^+ := \mathscr{X}_{{\rm PGL}_3,\frak{S}}(\mathbb{R}_{>0})
$$
which was studied by Fock and Goncharov \cite{FG06} and called a {\em higher Teichm\"uller space}. It is topologized e.g. as a subspace of $\mathscr{X}_{{\rm PGL}_3,\frak{S}}(\mathbb{R})$. Or, one could think of it as being obtained by gluing $(\mathbb{R}_{>0})^{|Q_\Delta|}$ associated to each cluster $\mathscr{X}$-chart (not just the cluster charts for ideal triangulations $\Delta$), along the mutation gluing maps. In this case, the gluing maps are diffeomorphisms, and so $\mathscr{X}_{{\rm PGL}_3,\frak{S}}^+$ is a smooth manifold diffeomorphic to $(\mathbb{R}_{>0})^{|Q_\Delta|}$. Given a regular function on $\mathscr{X}_{{\rm PGL}_3,\frak{S}}$, i.e. an element of $\mathscr{O}(\mathscr{X}_{{\rm PGL}_3,\frak{S}})$, for the cluster $\mathscr{X}$-chart associated to any ideal triangulation $\Delta$, this function can be written as a Laurent polynomial in the coordinate functions $X_v$'s, $v\in \mathcal{V}(Q_\Delta)$. This Laurent polynomial expression can be thought of as a smooth function on the manifold $\mathscr{X}_{{\rm PGL}_3,\frak{S}}^+$. In particular, each $X_v$ is a positive real valued smooth function on $\mathscr{X}_{{\rm PGL}_3,\frak{S}}^+$.

\vs

We will observe that there is an embedding
\begin{align}
\label{eq:section_of_PGL3_to_SL3}
\Psi : \mathscr{X}^+_{{\rm PGL}_3,\frak{S}} \to \mathscr{X}_{{\rm SL}_3,\frak{S}}(\mathbb{R}),
\end{align}
whose inverse map on the image coincides with the map $P$. Then we use this to translate the functions $\mathbb{I}_{{\rm SL}_3}(\ell) \in \mathscr{O}(\mathscr{X}_{{\rm SL}_3,\frak{S}})$ (for $\ell \in \mathscr{A}_{\rm L}(\frak{S};\mathbb{Z})$) to functions on the manifold $\mathscr{X}^+_{{\rm PGL}_3,\frak{S}}$.
\begin{definition}[translation of ${\rm SL}_3$ regular functions to ${\rm PGL}_3$]
\label{def:translation_to_PGL3}
For each $\ell \in \mathscr{A}_{\rm L}(\frak{S};\mathbb{Z})$, denote by $\mathbb{I}^+_{{\rm PGL}_3}(\ell)$ the function on $\mathscr{X}^+_{{\rm PGL}_3,\frak{S}}$ obtained as the pullback under the map eq.\eqref{eq:section_of_PGL3_to_SL3} of the function $\mathbb{I}_{{\rm SL}_3}(\ell) \in\mathscr{O}(\mathscr{X}_{{\rm SL}_3,\frak{S}})$ (evaluated at $\mathbb{R}$). Call this
$$
\mathbb{I}^+_{{\rm PGL}_3}(\ell) := \Psi^*(\mathbb{I}_{{\rm SL}_3}(\ell)(\mathbb{R})) ~\in~ C^\infty(\mathscr{X}^+_{{\rm PGL}_3,\frak{S}})
$$
a \ul{\em basic semi-regular function} on $\mathscr{X}^+_{{\rm PGL}_3,\frak{S}}$. They can be viewed as forming a map
$$
\mathbb{I}^+_{{\rm PGL}_3} : \mathscr{A}_{\rm L}(\frak{S};\mathbb{Z}) \to C^\infty(\mathscr{X}^+_{{\rm PGL}_3,\frak{S}}).
$$
\end{definition}

\vs

To construct the map $\Psi$ in eq.\eqref{eq:section_of_PGL3_to_SL3}, we partially recall Fock-Goncharov's \cite{FG06} reconstruction of a framed ${\rm PGL}_3$-local system on $\frak{S}$ out of the cluster $\mathscr{X}$-coordinates; see also \cite{Douglas}. Let $\Delta$ be an ideal triangulation of a punctured surface $\frak{S}$.  Let $\gamma$ be an oriented  loop in $\frak{S}$, not necessarily simple. Deform $\gamma$ by an isotopy if necessary, so that $\gamma$ meets $\Delta$ transversally, at finitely many points\redfix{, and use Reidemeister move I \raisebox{-0.4\height}{\scalebox{0.8}{}} to remove all kinks}. We call the elements of $\gamma \cap \Delta$ the \ul{\em $\Delta$-junctures} of $\gamma$. For each juncture of $\gamma$, choose a small neighborhood of it in $\gamma$, which is an oriented path meeting $\Delta$ exactly once; call this a \ul{\em juncture segment} of $\gamma$ corresponding to this $\Delta$-juncture. Each \redfix{maximal curve segment in} the complement in $\gamma$ of the union of all juncture segments is called a \ul{\em triangle segment} of $\gamma$. \redfix{Now, by a} \ul{\em segment} of $\gamma$ \redfix{we mean} either a juncture segment or a triangle segment. Then, by choosing a starting segment of $\gamma$, one can express $\gamma$ as a concatenation (or, path product) of a sequence of segments
\begin{align}
\label{eq:gamma_concatenation}
\gamma = \gamma_1 . \gamma_2 . \cdots . \gamma_N;
\end{align}
here $\gamma_1$ is the initial segment, and as one travels on $\gamma$ along its orientation, one then meets $\gamma_2$, and then $\gamma_3$, etc. Notice that this sequence alternates between juncture segments and triangle segments, and that $N$ is even. So, if $\gamma_1$ is a triangle segment, then $\gamma_2$ is a juncture segment, $\gamma_3$ is a triangle segment, and so on, and the last $\gamma_N$ is a juncture segment. Note that a triangle segment is exactly one of a \ul{\em left turn}, a \ul{\em right turn}, or a \ul{\em U-turn}. Examples are shown below.

\begin{center}
\scalebox{1.0}{
\begingroup%
  \makeatletter%
  \providecommand\color[2][]{%
    \errmessage{(Inkscape) Color is used for the text in Inkscape, but the package 'color.sty' is not loaded}%
    \renewcommand\color[2][]{}%
  }%
  \providecommand\transparent[1]{%
    \errmessage{(Inkscape) Transparency is used (non-zero) for the text in Inkscape, but the package 'transparent.sty' is not loaded}%
    \renewcommand\transparent[1]{}%
  }%
  \providecommand\rotatebox[2]{#2}%
  \newcommand*\fsize{\dimexpr\f@size pt\relax}%
  \newcommand*\lineheight[1]{\fontsize{\fsize}{#1\fsize}\selectfont}%
  \ifx\svgwidth\undefined%
    \setlength{\unitlength}{447.87401575bp}%
    \ifx\svgscale\undefined%
      \relax%
    \else%
      \setlength{\unitlength}{\unitlength * \real{\svgscale}}%
    \fi%
  \else%
    \setlength{\unitlength}{\svgwidth}%
  \fi%
  \global\let\svgwidth\undefined%
  \global\let\svgscale\undefined%
  \makeatother%
  \begin{picture}(1,0.3164557)%
    \lineheight{1}%
    \setlength\tabcolsep{0pt}%
    \put(0,0){\includegraphics[width=\unitlength,page=1]{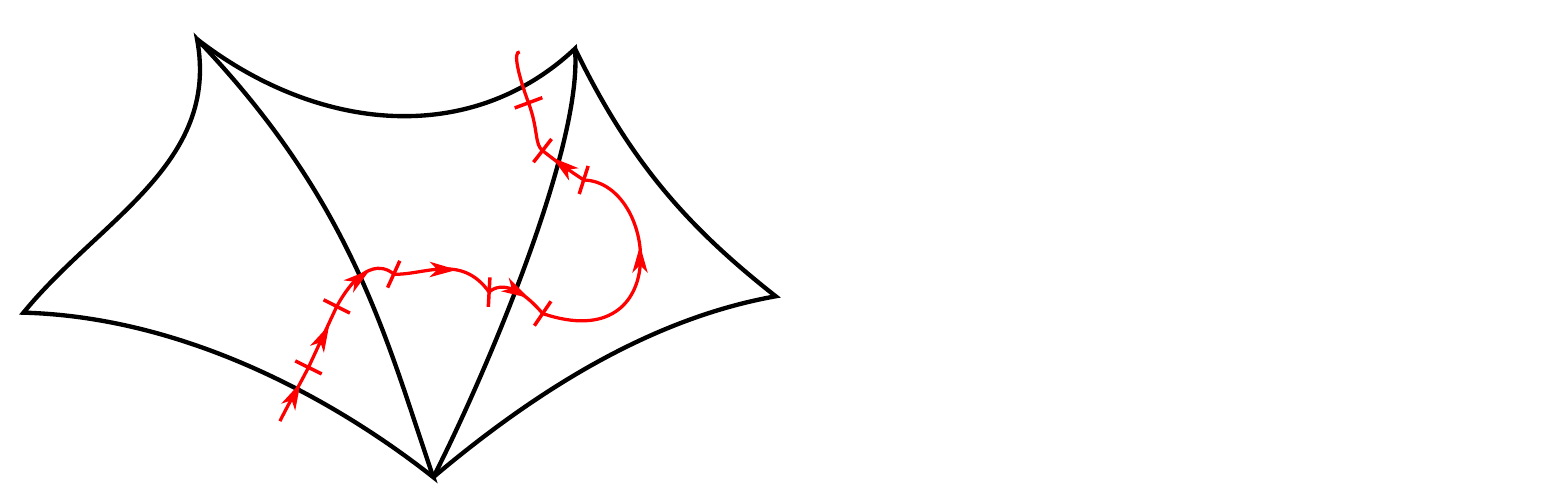}}%
    \put(0.18137654,0.10512986){\color[rgb]{0,0,0}\makebox(0,0)[lt]{\lineheight{1.25}\smash{\begin{tabular}[t]{l}$\red{\gamma_1}$\end{tabular}}}}%
    \put(0.19708522,0.14272192){\color[rgb]{0,0,0}\makebox(0,0)[lt]{\lineheight{1.25}\smash{\begin{tabular}[t]{l}$\red{\gamma_2}$\end{tabular}}}}%
    \put(0.27687879,0.15262337){\color[rgb]{0,0,0}\makebox(0,0)[lt]{\lineheight{1.25}\smash{\begin{tabular}[t]{l}$\red{\gamma_3}$\end{tabular}}}}%
    \put(0.33819446,0.1410829){\color[rgb]{0,0,0}\makebox(0,0)[lt]{\lineheight{1.25}\smash{\begin{tabular}[t]{l}$\red{\gamma_4}$\end{tabular}}}}%
    \put(0.41593805,0.13145534){\color[rgb]{0,0,0}\makebox(0,0)[lt]{\lineheight{1.25}\smash{\begin{tabular}[t]{l}$\red{\gamma_5}$\end{tabular}}}}%
    \put(0.36674518,0.2197074){\color[rgb]{0,0,0}\makebox(0,0)[lt]{\lineheight{1.25}\smash{\begin{tabular}[t]{l}$\red{\gamma_6}$\end{tabular}}}}%
    \put(0.15286716,0.06269594){\color[rgb]{0,0,0}\makebox(0,0)[lt]{\lineheight{1.25}\smash{\begin{tabular}[t]{l}$\red{\gamma_N}$\end{tabular}}}}%
    \put(0,0){\includegraphics[width=\unitlength,page=2]{segments_of_gamma.pdf}}%
    \put(0.23311407,0.25603466){\color[rgb]{0,0,0}\makebox(0,0)[lt]{\lineheight{1.25}\smash{\begin{tabular}[t]{l}$\red{\gamma_i}$\end{tabular}}}}%
    \put(0.25118654,0.19919116){\color[rgb]{0,0,0}\makebox(0,0)[lt]{\lineheight{1.25}\smash{\begin{tabular}[t]{l}$\red{\gamma_{i+1}}$\end{tabular}}}}%
    \put(0.18905753,0.20319822){\color[rgb]{0,0,0}\makebox(0,0)[lt]{\lineheight{1.25}\smash{\begin{tabular}[t]{l}$\red{\gamma_{i+2}}$\end{tabular}}}}%
    \put(0.11513543,0.17123351){\color[rgb]{0,0,0}\makebox(0,0)[lt]{\lineheight{1.25}\smash{\begin{tabular}[t]{l}$\red{\gamma_{i+3}}$\end{tabular}}}}%
    \put(0.0715426,0.07725836){\color[rgb]{0,0,0}\makebox(0,0)[lt]{\lineheight{1.25}\smash{\begin{tabular}[t]{l}$\red{\gamma_{i+4}}$\end{tabular}}}}%
    \put(0,0){\includegraphics[width=\unitlength,page=3]{segments_of_gamma.pdf}}%
    \put(0.31406919,0.2296885){\color[rgb]{0,0,0}\makebox(0,0)[lt]{\lineheight{1.25}\smash{\begin{tabular}[t]{l}$\red{\gamma_7}$\end{tabular}}}}%
    \put(0.31575596,0.29261097){\color[rgb]{0,0,0}\makebox(0,0)[lt]{\lineheight{1.25}\smash{\begin{tabular}[t]{l}$\red{\gamma_8}$\end{tabular}}}}%
    \put(0.44898311,0.2369948){\color[rgb]{0,0,0}\makebox(0,0)[lt]{\lineheight{1.25}\smash{\begin{tabular}[t]{l}triangle segments : $\gamma_1$, $\gamma_3$, $\gamma_5$, \ldots, $\gamma_{i+1}$, $\gamma_{i+3}$, \ldots, $\gamma_{N-1}$\end{tabular}}}}%
    \put(0.44598538,0.27734658){\color[rgb]{0,0,0}\makebox(0,0)[lt]{\lineheight{1.25}\smash{\begin{tabular}[t]{l}juncture segments : $\gamma_2$, $\gamma_4$, $\gamma_6$, \ldots, $\gamma_i$, $\gamma_{i+2}$, \ldots, $\gamma_N$\end{tabular}}}}%
    \put(0.52256355,0.18968258){\color[rgb]{0,0,0}\makebox(0,0)[lt]{\lineheight{1.25}\smash{\begin{tabular}[t]{l}left turns : $\gamma_{i+3}$, \ldots\end{tabular}}}}%
    \put(0.5226617,0.15287372){\color[rgb]{0,0,0}\makebox(0,0)[lt]{\lineheight{1.25}\smash{\begin{tabular}[t]{l}right turns : $\gamma_1$, $\gamma_3$, $\gamma_7$, $\gamma_{i+1}$, \ldots\end{tabular}}}}%
    \put(0.52294925,0.11599803){\color[rgb]{0,0,0}\makebox(0,0)[lt]{\lineheight{1.25}\smash{\begin{tabular}[t]{l}U-turns : $\gamma_5$, \ldots\end{tabular}}}}%
  \end{picture}%
\endgroup%
}
\end{center}

To each segment $\gamma_i$, we assign a \ul{\em monodromy matrix}
$$
{\bf M}_{\gamma_i} \in {\rm SL}_3 ( \mathbb{Z}[\{X_v^{\pm \frac{1}{3}} \, | \, v\in \mathcal{V}(Q_\Delta) \}])
$$
as follows, where $X_v$'s are the Fock-Goncharov $\mathscr{X}$-coordinates of the space $\mathscr{X}_{{\rm PGL}_3,\frak{S}}$ associated to nodes $v$ of the quiver $Q_\Delta$, i.e. the coordinates for the cluster $\mathscr{X}$-chart for $\Delta$. One can view the symbol $X_v^{\frac{1}{3}}$ as a generator of a formally defined Laurent polynomial ring $\mathbb{Z}[\{X_v^{\pm \frac{1}{3}} \, | \, v\in \mathcal{V}(Q_\Delta) \}]$, in which $\mathbb{Z}[\{X_v^{\pm 1} \, | \, v\in \mathcal{V}(Q_\Delta) \}]$ embeds into as $X_v \mapsto (X_v^{\frac{1}{3}})^3$. Or, as in \cite{FG06}, we can also view $X_v^{\pm\frac{1}{3}}$ as functions on a covering space $\wh{\mathscr{X}_{{\rm PGL}_3,\frak{S}}}$ of $\mathscr{X}_{{\rm PGL}_3,\frak{S}}$. Our approach here will be to view each $X_v$ as a positive-real valued smooth function on the manifold $\mathscr{X}_{{\rm PGL}_3,\frak{S}}^+$; then $X_v^{\pm \frac{1}{3}}$ is well defined as a positive real valued smooth function on $\mathscr{X}_{{\rm PGL}_3,\frak{S}}^+$.

\vs

\begin{figure}[htbp!]
\vspace{-2mm}
\begin{center}
\scalebox{1.0}{
\begingroup%
  \makeatletter%
  \providecommand\color[2][]{%
    \errmessage{(Inkscape) Color is used for the text in Inkscape, but the package 'color.sty' is not loaded}%
    \renewcommand\color[2][]{}%
  }%
  \providecommand\transparent[1]{%
    \errmessage{(Inkscape) Transparency is used (non-zero) for the text in Inkscape, but the package 'transparent.sty' is not loaded}%
    \renewcommand\transparent[1]{}%
  }%
  \providecommand\rotatebox[2]{#2}%
  \newcommand*\fsize{\dimexpr\f@size pt\relax}%
  \newcommand*\lineheight[1]{\fontsize{\fsize}{#1\fsize}\selectfont}%
  \ifx\svgwidth\undefined%
    \setlength{\unitlength}{141.73228346bp}%
    \ifx\svgscale\undefined%
      \relax%
    \else%
      \setlength{\unitlength}{\unitlength * \real{\svgscale}}%
    \fi%
  \else%
    \setlength{\unitlength}{\svgwidth}%
  \fi%
  \global\let\svgwidth\undefined%
  \global\let\svgscale\undefined%
  \makeatother%
  \begin{picture}(1,0.42)%
    \lineheight{1}%
    \setlength\tabcolsep{0pt}%
    \put(0,0){\includegraphics[width=\unitlength,page=1]{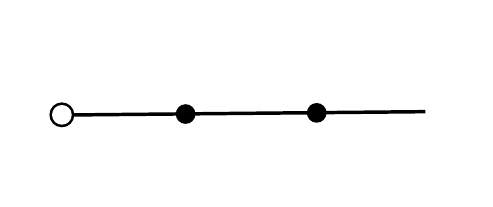}}%
    \put(0.23476515,0.33308462){\color[rgb]{0,0,0}\makebox(0,0)[lt]{\lineheight{1.25}\smash{\begin{tabular}[t]{l}\red{juncture segment $\gamma_i$}\end{tabular}}}}%
    \put(0,0){\includegraphics[width=\unitlength,page=2]{edge_crossing2.pdf}}%
    \put(0.06732883,0.23319851){\color[rgb]{0,0,0}\makebox(0,0)[lt]{\lineheight{1.25}\smash{\begin{tabular}[t]{l}edge of triangulation\end{tabular}}}}%
    \put(0.32758923,0.11758736){\color[rgb]{0,0,0}\makebox(0,0)[lt]{\lineheight{1.25}\smash{\begin{tabular}[t]{l}$X_1$\end{tabular}}}}%
    \put(0.60358231,0.11870907){\color[rgb]{0,0,0}\makebox(0,0)[lt]{\lineheight{1.25}\smash{\begin{tabular}[t]{l}$X_2$\end{tabular}}}}%
  \end{picture}%
\endgroup%
}
\end{center}
\vspace{-9mm}
\caption{Juncture segment $\gamma_i$, intersecting an edge of $\Delta$}
\vspace{-2mm}
\label{fig:juncture_segment}
\end{figure}

\begin{enumerate}
\item[\rm (MM1)] (edge matrix) If $\gamma_i$ is a juncture segment that meets the edge $e$ of $T$, then if the Fock-Goncharov $\mathscr{X}$-coordinates at the two nodes of $Q_\Delta$ lying in $e$ are $X_1$ and $X_2$, where $X_1$ is located at the left and $X_2$ at the right with respect to the orientation $\gamma_i$, as depicted in Fig.\ref{fig:juncture_segment}, we have
$$
{\bf M}_{\gamma_i} = \smallmatthree{X_1^{1/3} X_2^{2/3}}{0}{0}{0}{X_1^{1/3} X_2^{-1/3}}{0}{0}{0}{X_1^{-2/3} X_2^{-1/3}};
$$

\item[\rm (MM2)] (left turn matrix) If $\gamma_i$ is a triangle segment in a triangle $t$ and turns to left, then if $X$ is the Fock-Goncharov $\mathscr{X}$-coordinate at the node of $Q_\Delta$ lying in the interior of $t$, we have
$$
{\bf M}_{\gamma_i} = \smallmatthree{X^{2/3}}{~X^{2/3}+X^{-1/3}~}{X^{-1/3}}{0}{X^{-1/3}}{X^{-1/3}}{0}{0}{X^{-1/3}};
$$

\item[\rm (MM3)] (right turn matrix) If $\gamma_i$ is a triangle segment in a triangle $t$ and turns to right, then if $X$ is the Fock-Goncharov $\mathscr{X}$-coordinate at the node of $Q_\Delta$ lying in the interior of $t$, we have
$$
{\bf M}_{\gamma_i} = \smallmatthree{X^{1/3}}{0}{0}{X^{1/3}}{X^{1/3}}{0}{X^{1/3}}{~X^{1/3}+X^{-2/3}~}{X^{-2/3}};
$$

\item[\rm (MM4)] (U-turn matrix) If $\gamma_i$ is a triangle segment in a triangle $t$ taking a U-turn, we have
$$
{\bf M}_{\gamma_i} = \smallmatthree{0}{0}{1}{0}{-1}{0}{1}{0}{0}.
$$
\end{enumerate}
Define then the \ul{\em monodromy matrix} along $\gamma$ to be the product
\begin{align}
\label{eq:bf_M_gamma}
{\bf M}_\gamma := {\bf M}_{\gamma_1} {\bf M}_{\gamma_2} \cdots {\bf M}_{\gamma_N} \quad \in \quad {\rm SL}_3 ( \mathbb{Z}[\{X_v^{\pm \frac{1}{3}} \, | \, v\in \mathcal{V}(Q_\Delta) \}]).
\end{align}
One can always isotope $\gamma$ so that there is no U-turn, so the entries of ${\bf M}_\gamma$ actually lies in $\mathbb{Z}_{\ge 0}[\{X_v^{\pm \frac{1}{3}} \, | \, v\in \mathcal{V}(Q_\Delta) \}]$\redfix{, with emphasis on the non-negativity}. Hence its trace ${\rm tr}({\bf M}_\gamma)$ is an element of $\mathbb{Z}_{\ge 0}[\{X_v^{\pm \frac{1}{3}} \, | \, v\in \mathcal{V}(Q_\Delta) \}]$:
$$
{\rm tr}({\bf M}_\gamma) \quad \in \quad \mathbb{Z}_{\ge 0}[\{X_v^{\pm \frac{1}{3}} \, | \, v\in \mathcal{V}(Q_\Delta) \}].
$$

\vs

Observe that for each segment $\gamma_i$, the monodromy matrix ${\bf M}_{\gamma_i}$ times a monomial $\prod_v X_v^{k_v/3}$ for some integer $k_v$'s lies in ${\rm GL}_3 (\mathbb{Z}[\{X_v^{\pm 1} \, | \, v\in \mathcal{V}(Q_\Delta) \}])$. Let $\til{\bf M}_{\gamma_i} = (\prod_v X_v^{k_v/3}) {\bf M}_{\gamma_i}$ for each $i$, and
$$
\til{\bf M}_\gamma := \til{\bf M}_{\gamma_1} \cdots \til{\bf M}_{\gamma_N} \quad \in \quad {\rm GL}_3 (\mathbb{Z}[\{X_v^{\pm 1} \, | \, v\in \mathcal{V}(Q_\Delta) \}]).
$$
Since the determinant of each $\til{\bf M}_{\gamma_i}$ is a Laurent monomial in $X_v$'s, so is that of $\til{\bf M}_\gamma$. However, note that such a normalization for $\til{\bf M}_\gamma$ is not unique, and it is defined only up to a Laurent monomial in $X_v$'s. Later, we will use a specific choice of normalization (which makes the $(1,1)$-th entry to be $1$). \redfix{We note that what are used in \cite{FG06} are some normalizations $\til{\bf M}_{\gamma_i}$ and $\til{\bf M}_\gamma$, instead of ${\bf M}_{\gamma_i}$ and ${\bf M}_\gamma$, where the latter ones play crucial roles in the present paper.}

\vs

An important example is a peripheral loop surrounding a puncture. An easy observation:
\begin{lemma}
\label{lem:peripheral_is_same_turns}
An oriented loop $\gamma$ \redfix{without a kink} is a peripheral loop if and only if it can be isotoped so that the triangle segments are either all left turns, or all right turns. 

\vs

When the orientation of $\gamma$ matches the boundary-orientation of the corresponding hole of $\til{\frak{S}}$ (Def.\ref{def:holed_surface}), the triangle segments are all left turns. Otherwise, all right turns. \qed
\end{lemma}
Therefore we get:
\begin{lemma}
If $\gamma$ is a peripheral loop surrounding a puncture, then ${\bf M}_\gamma$ (hence also $\til{\bf M}_\gamma$) is either an upper triangular matrix or a lower triangular matrix. It is upper triangular if and only if $\gamma$ matches the boundary-orientation of the corresponding hole of $\til{\frak{S}}$.
\end{lemma}
So, given a tuple of (nonzero) coordinates $(X_v)_{v\in \mathcal{V}(Q_\Delta)}$, i.e. given a point of $(\mathbb{G}_m)^{\mathcal{V}(Q_\Delta)}$, one can reconstruct a monodromy representation $\pi_1(\frak{S}) \to {\rm PGL}_3$ defined up to conjugation, by setting the image of $[\gamma]$ for a loop $\gamma$ under the sought-for map $\pi_1(\frak{S}) \to {\rm PGL}_3$ to be the image of the matrix $\til{\bf M}_\gamma = \til{\bf M}_\gamma( (X_v)_{v\in \mathcal{V}(Q_\Delta)} ) \in {\rm GL}_3$ under the projection ${\rm GL}_3 \to {\rm PGL}_3$. More precisely:
\begin{proposition}[\cite{FG06}]
\label{prop:reconstruction}
Let $\gamma$ be an oriented simple loop in $\frak{S}$ and $\Delta$ an ideal triangulation of $\frak{S}$. Then, for a framed ${\rm PGL}_3$-local system on $\frak{S}$, the image under the underlying monodromy $\pi_1(\frak{S}) \to {\rm PGL}_3$ of the equivalence class $[\gamma] \in \pi_1(\frak{S})$ of $\gamma$ coincides with the image under the projection ${\rm GL}_3 \to {\rm PGL}_3$ of the matrix $\til{\bf M}_\gamma = \til{\bf M}_\gamma( (X_v)_{v\in \mathcal{V}(Q_\Delta)} ) \in {\rm GL}_3$ constructed above.
\end{proposition}

In fact, the monodromy matrices ${\bf M}_{\gamma_i}$\redfix{, or $\til{\bf M}_{\gamma_i}$,} let us completely reconstruct a point of $\mathscr{X}_{{\rm PGL}_3,\frak{S}}$; namely, Fock Goncharov \cite[\S9]{FG06} considered a certain graph on $\frak{S}$ and assigned these matrices to its graph, and constructed a ${\rm PGL}_3$-local system explicitly (not just its monodromy), together with a framing. 

\vs

Coming back to our strategy, let us construct the promised map $\Psi$ of eq.\eqref{eq:section_of_PGL3_to_SL3}. Given a point of the domain $\mathscr{X}^+_{{\rm PGL}_3,\frak{S}}$, we can record it by its positive real coordinates $X_v$'s, for $v\in \mathcal{V}(Q_\Delta)$, for any chosen ideal triangulation $\Delta$. Consider the monodromy $\rho : \pi_1(\frak{S}) \to {\rm PGL}_3(\mathbb{R})$ for this point, which is a homomorphism defined up to conjugation. Above, we saw explicitly how $\rho([\gamma]) \in {\rm PGL}_3(\mathbb{R})$ is given in terms of the coordinates $X_v$'s, for each $\gamma \in \pi_1(\frak{S})$. In fact, we can lift it to a\redfix{n} ${\rm SL}_3(\mathbb{R})$ monodromy $\til{\rho} : \pi_1(\frak{S}) \to {\rm SL}_3(\mathbb{R})$. Pick any basepoint $x\in \frak{S}$, and let $\gamma$ be a loop based at $x$. Define $\til{\rho}([\gamma]) := {\bf M}_\gamma \in {\rm SL}_3(\mathbb{R})$, which maps to $\rho([\gamma]) \in {\rm PGL}_3(\mathbb{R})$ under the projection ${\rm SL}_3(\mathbb{R}) \to {\rm PGL}_3(\mathbb{R})$. Since ${\rm SL}_3(\mathbb{R}) \to {\rm PGL}_3(\mathbb{R})$ is a bijection, it follows that $\til{\rho}:\pi_1(\frak{S}) \to {\rm SL}_3(\mathbb{R})$ is a genuine (as opposed to projective) group homomorphism, defined up to conjugation. As mentioned above, the monodromy matrices ${\bf M}_\gamma$ also yield the specific choice of the framing data at punctures, so that one indeed obtains a point of $\mathscr{X}_{{\rm SL}_3,\frak{S}}(\mathbb{R})$. Using the fact that the projection ${\rm SL}_3(\mathbb{R}) \to {\rm PGL}_3(\mathbb{R})$ is bijective, one can observe that the resulting point of $\mathscr{X}_{{\rm SL}_3,\frak{S}}(\mathbb{R})$ does not depend on $\Delta$. For us, what we need to know about the framing data are the functions $\pi_p : \mathscr{X}_{{\rm SL}_3,\frak{S}}(\mathbb{R}) \to {\rm H}(\mathbb{R})$ at punctures $p$ (eq.\eqref{eq:pi_p}), and corresponding regular functions $(\pi_p)_i$ (eq.\eqref{eq:pi_p_i}). As observed above, the ${\rm SL}_3(\mathbb{R})$ monodromy ${\bf M}_{\gamma_p}$ along a \redfix{positively oriented (Def.\ref{def:positive_orientation})} peripheral loop $\gamma_p$ surrounding a puncture $p$ is upper triangular, hence belongs to our choice of the Borel subgroup ${\rm B}(\mathbb{R})$ of ${\rm G}(\mathbb{R}) = {\rm SL}_3(\mathbb{R})$. Composing with the quotient map ${\rm B}(\mathbb{R}) \to {\rm B}(\mathbb{R})/{\rm U}(\mathbb{R}) = {\rm H}(\mathbb{R})$ which extracts the semi-simple, i.e. the diagonal, part, we obtain an element of ${\rm H}(\mathbb{R})$, yielding the value of the function $\pi_p : \mathscr{X}_{{\rm SL}_3,\frak{S}}(\mathbb{R}) \to {\rm H}(\mathbb{R})$ at this point of $\mathscr{X}_{{\rm SL}_3,\frak{S}}(\mathbb{R})$. Namely, for this $\gamma_p$ we have:
\begin{align}
\label{eq:pi_p_i_investigated}
(\pi_p)_i = \mbox{the $i$-th diagonal entry of the upper triangular monodromy matrix ${\bf M}_{\gamma_p}$}.
\end{align}

\vs

Using the map $\Psi$ just constructed, we apply Def.\ref{def:translation_to_PGL3} and get the basic semi-regular functions $\mathbb{I}^+_{{\rm PGL}_3}(\ell)$ on $\mathscr{X}^+_{{\rm PGL}_3,\frak{S}}$ by pulling back the $A_2$-bangles basis functions $\mathbb{I}_{{\rm SL}_3}(\ell) \in \mathscr{O}(\mathscr{X}_{{\rm SL}_3,\frak{S}})$ for $\ell \in \mathscr{A}_{\rm L}(\frak{S};\mathbb{Z})$. A basic example is when $\ell$ is a single oriented simple loop. More generally, we consider the pullback of trace-of-monodromy function along any oriented loop, which we also call a trace-of-monodromy.
\begin{definition}
\label{def:f_plus_gamma}
For an oriented loop $\gamma$ in $\frak{S}$, we define the trace-of-monodromy function $f^+_\gamma$ on $\mathscr{X}^+_{{\rm PGL}_3,\frak{S}}$ as
$$
f^+_\gamma := {\rm tr}({\bf M}_\gamma)
$$
which is a smooth function on the manifold $\mathscr{X}^+_{{\rm PGL}_3,\frak{S}}$. 
\end{definition}
For any oriented loop $\gamma$, for each triangulation $\Delta$, this function $f^+_\gamma$ can be written as a Laurent polynomial in $\{ X_v^{1/3} \, | \, v\in \mathcal{V}(Q_\Delta)\}$, with non-negative integer coefficients. Changing the basepoint $x$ results in a new matrix ${\bf M}_\gamma$ related to the previous one by a conjugation, hence the trace doesn't change.

\vs

If $\ell$ is a single-component ${\rm SL}_3$-lamination consisting of an oriented {\em simple} \redfix{non-peripheral} loop $\gamma$ with weight $1$, then
\begin{align}
\label{eq:I_plus_PGL3_simple_loop}
\mathbb{I}^+_{{\rm PGL}_3}(\ell) = f^+_\gamma.
\end{align}
We will be dealing with more general ${\rm SL}_3$-laminations in the coming subsections.

\subsection{Functions for punctures}

In order to investigate $\mathbb{I}^+_{{\rm PGL}_3}(\ell)$ in case $\ell$ consists only of a peripheral loop, we study the trace-of-monodromy function $f^+_\gamma = {\rm tr}({\bf M}_\gamma)$ on $\mathscr{X}^+_{{\rm PGL}_3,\frak{S}}$ for a peripheral loop $\gamma$. 
\begin{proposition}[peripheral monodromy]
\label{prop:peripheral_monodromy}
Let $\gamma$ be an oriented peripheral loop. Denote by $\ell$ the ${\rm SL}_3$-lamination consisting just of $\gamma$ with weight $1$, and by $\ol{\ell}$ the ${\rm SL}_3$-lamination consisting of just of the orientation-reversed loop $\ol{\gamma}$ with weight $1$.

\vs

If $\gamma$ is positively oriented (Def.\ref{def:positive_orientation}), then for any choice of basepoint of $\gamma$, the monodromy matrix ${\bf M}_\gamma$ is given by the upper triangular matrix
\begin{align}
\label{eq:peripheral_monodromy_matrix}
{\bf M}_\gamma = \matthree{\prod_v X_v^{{\rm a}_v(\ell)}}{*}{*}{0}{\prod_v X_v^{-{\rm a}_v(\ell)+{\rm a}_v(\ol{\ell})}}{*}{0}{0}{\prod_v X_v^{-{\rm a}_v(\ol{\ell})}}
\end{align}
in terms of the coordinate functions of any triangulation $\Delta$, where the products $\prod_v$ are taken over all $v$ in $\mathcal{V}(Q_\Delta)$. Hence the trace-of-monodromy function $f^+_\gamma$ on $\mathscr{X}^+_{{\rm PGL}_3,\frak{S}}$ is
\begin{align}
\label{eq:peripheral_monodromy_trace}
f^+_\gamma = \underset{v\in \mathcal{V}(Q_\Delta)}{\textstyle \prod} X_v^{{\rm a}_v(\ell)}
+ \underset{v\in \mathcal{V}(Q_\Delta)}{\textstyle \prod} X_v^{-{\rm a}_v(\ell)+{\rm a}_v(\ol{\ell})}
+ \underset{v\in \mathcal{V}(Q_\Delta)}{\textstyle \prod} X_v^{-{\rm a}_v(\ol{\ell})}
\end{align}
where ${\rm a}_v(\ell), {\rm a}_v(\ol{\ell}) \in \frac{1}{3}\mathbb{Z}$ are the tropical coordinates of the ${\rm SL}_3$-laminations $\ell$ and $\ol{\ell}$ at the node $v$ of the quiver $Q_\Delta$; in particular, ${\rm a}_v(\ell)$ and ${\rm a}_v(\ol{\ell})$ are non-negative. The trace formula eq.\eqref{eq:peripheral_monodromy_trace} holds also when $\gamma$ is negatively oriented.
\end{proposition}

{\it Proof.} Pick any basepoint of $\gamma$, and express $\gamma$ as a concatenation of segments as in eq.\eqref{eq:gamma_concatenation}. The triangle segments can be assumed to be all right turns, or all left turns (Lem.\ref{lem:peripheral_is_same_turns}). We assume that they are all left turns, i.e. $\gamma$ is positively oriented. Proof for the other case would be similar. Different choice of a basepoint results in cyclically shifting the concatenation expression in eq.\eqref{eq:gamma_concatenation}. Since each ${\bf M}_{\gamma_i}$ is upper triangular, it follows that ${\bf M}_\gamma = {\bf M}_{\gamma_1} {\bf M}_{\gamma_2} \cdots {\bf M}_{\gamma_N}$ is also upper triangular, and cyclic shift of the product order yields an upper triangular matrix ${\bf M}_\gamma'$ with same diagonal entries as ${\bf M}_\gamma$.

\vs

Let $p$ be the puncture that $\gamma$ is surrounding. For an edge $e$ of $\Delta$, $\gamma$ meets $e$ once if only one of the two endpoints of $e$ is $p$, twice if both endpoints of $e$ are $p$, and does not meet $e$ if none of the endpoints of $e$ is $p$. When $\gamma$ meets $e$ twice, they meet in different configuration of orientations as follows; given an arbitrary orientation on $e$, at each of the two intersection points $x$ of $\gamma$ and $e$, the velocity vectors of $\gamma$ and $e$ (in this order) form a positively oriented basis of $T_x\frak{S}$ (according to the orientation of the surface $\frak{S}$) at one $x$ and a negatively oriented basis for the other $x$. Now, let $t$ be any ideal triangle of $\Delta$, in which there is at least one triangle segment of $\gamma$. The triangle segments of $\gamma$ in $t$ are all left turns, and by the above discussion, each corner of $t$ can have at most one such triangle segment; if there were two, then $\gamma$ would meet some edge of $\Delta$ twice with same configuration of orientations. 

\vs

To investigate the tropical coordinates of $\ell = \gamma$, consider a split ideal triangulation $\wh{\Delta}$ for $\Delta$. For convenience, one can isotope so that the intersection points of $\gamma$ with $\wh{\Delta}$ are exactly the breaking points of the concatenation decomposition of $\gamma$ as in eq.\eqref{eq:gamma_concatenation}. That is, the intersection points $\gamma \cap \wh{\Delta}$ divide $\gamma$ into the pieces, where a piece in a biangle is a juncture segment, and a piece in a triangle is a triangle segment. In particular, now a triangle segment is what we called a corner arc before. Observe that as of now, $\ell=\gamma$ is an ${\rm SL}_3$-lamination that is canonical with respect to $\wh{\Delta}$ (Def.\ref{def:canonical_wrt_split_ideal_triangulation}), so we can read the tropical coordinates as in Def.\ref{def:tropical_coordinates}.

\vs

Let $t$ be an ideal triangle of $\Delta$, and $\wh{t}$ be the corresponding triangle of $\wh{\Delta}$. Let $e_1,e_2,e_3$ be the sides of $\wh{t}$, appearing clockwise in this order along $\partial \wh{t}$. On each $e_\alpha$, there are two nodes $v_{e_\alpha,1}$ and $v_{e_\alpha,2}$ of $Q_\Delta$ so that the direction $v_{e_\alpha,1} \to v_{e_\alpha,2}$ matches the clockwise orientation of $\partial \wh{t}$; in fact, these nodes should be viewed as living on an edge of $\Delta$ \redfix{(instead of $\wh{\Delta}$)}, but now \redfix{we} are focusing on only one triangle, so we can be ambiguous. Let $v_t$ be the node of $Q_\Delta$ lying in the interior of $\wh{t}$. So, in total, we are considering seven nodes of $Q_\Delta$ in $\wh{t}$ (or in $t$). Let $\gamma_j$ be a triangle segment of $\gamma$ in $\wh{t}$, which is a left turn segment and hence a left turn corner arc in $\wh{t}$. Say, the initial endpoint of $\gamma_j$ lies in the side $e_\alpha$; then the terminal endpoint of $\gamma_j$ lies in $e_{\alpha+1}$ (where $e_4 = e_1$). The tropical coordinates of this $\gamma_j$ are given as in eq.\eqref{eq:tropical_coordinates_for_a_left_turn_W} with $W_{\alpha,\alpha+1}=\gamma_j$. Denoting by $\ol{\gamma_j}$ the triangle segment of the orientation-reversed loop $\ol{\gamma}$ corresponding to $\gamma_j$, by viewing it as an ${\rm SL}_3$-lamination in $\wh{t}$ that is a right turn corner arc in $\wh{t}$, its tropical coordinates are as given in eq.\eqref{eq:tropical_coordinates_for_a_right_turn_W} with $W_{\alpha+1,\alpha} = \ol{\gamma_j}$.

\vs

On the other hand, let's now consider the monodromy matrix contribution, from the three segments $\gamma_{j-1},\gamma_j,\gamma_{j+1}$. We claim that, for a fixed triangle $t$, the basepoint of $\gamma$ could have been chosen in the beginning such that for each triangle segment $\gamma_j$ in $\wh{t}$ we have $1<j<N$. Indeed, since there are at least two triangles meeting $\gamma$, one could choose the basepoint of $\gamma$ such that the initial segment $\gamma_1$ is a triangle segment not living in $t$; thus $1<j$ for any triangle segment $\gamma_j$ living in $t$. Meanwhile, the concatenation sequence $\gamma_1,\ldots,\gamma_N$ must end with a juncture segment, hence it follows that $j<N$ for any triangle segment $\gamma_j$ living in $t$, as desired.

\vs

Note that the triples $(\gamma_{j-1},\gamma_j,\gamma_{j+1})$ associated to different triangle segments $\gamma_j$ living in $\wh{t}$ (or $t$) are disjoint with each other.  From (MM1)--(MM2), it follows that the corresponding product of monodromy matrices ${\bf M}_{\gamma_{j-1}} {\bf M}_{\gamma_j} {\bf M}_{\gamma_{j+1}}$ equals
\begin{align*}
 & \hspace{-10mm} \smallmatthree{ X_{v_{e_\alpha,2}}^{1/3} X_{v_{e_\alpha,1}}^{2/3} }{ 0 }{ 0 }{ 0 }{ X_{v_{e_\alpha,2}}^{1/3} X_{v_{e_\alpha,1}}^{-1/3} }{0}{0}{0}{X_{v_{e_\alpha,2}}^{-2/3} X_{v_{e_\alpha,1}}^{-1/3}}
\smallmatthree{X_{v_t}^{2/3}}{*}{*}{0}{X_{v_t}^{-1/3}}{*}{0}{0}{X_{v_t}^{-1/3}}
\smallmatthree{ X_{v_{e_{\alpha+1},1}}^{1/3} X_{v_{e_{\alpha+1},2}}^{2/3} }{ 0 }{ 0 }{ 0 }{ X_{v_{e_{\alpha+1},1}}^{1/3} X_{v_{e_{\alpha+1},2}}^{-1/3} }{0}{0}{0}{X_{v_{e_{\alpha+1},1}}^{-2/3} X_{v_{e_{\alpha+1},2}}^{-1/3}} \\
&\hspace{-10mm} = \smallmatthree{\prod_v X_v^{{\rm a}_v(\gamma_j)}}{*}{*}{0}{\prod_v X_v^{-{\rm a}_v(\gamma_j)+{\rm a}_v(\ol{\gamma_j})}}{*}{0}{0}{\prod_v X_v^{-{\rm a}_v(\ol{\gamma_j})}}
\end{align*}
with the last equality holding in view of the tropical coordinate values as in eq.\eqref{eq:tropical_coordinates_for_a_left_turn_W} with $W_{\alpha,\alpha+1}=\gamma_j$ and eq.\eqref{eq:tropical_coordinates_for_a_right_turn_W} with $W_{\alpha+1,\alpha} = \ol{\gamma_j}$, where $\prod_v$ is taken over seven nodes of $Q_\Delta$ living in $\wh{t}$ (or $t$).

\vs

Note ${\bf M}_\gamma = {\bf M}_{\gamma_1} {\bf M}_{\gamma_2} \cdots {\bf M}_{\gamma_N}$, where each factor ${\bf M}_{\gamma_i}$ is upper triangular with diagonal entries being Laurent monomials in $X_v^{1/3}$, $v\in \mathcal{V}(Q_\Delta)$. For each of the three diagaonal entries of ${\bf M}_\gamma$, we need to know the power of $X_v^{1/3}$ for each $v\in \mathcal{V}(Q_\Delta)$. Let's read the powers of $X_v^{1/3}$ for nodes $v$ living in $t$ (or in $\wh{t}$\,\,). Note that for each $\gamma_i$ that is not part of a triple $(\gamma_{j_1},\gamma_j,\gamma_{j+1})$ for a triangle segment $\gamma_j$ living in $t$, the monodromy matrix ${\bf M}_{\gamma_i}$ does not involve any $X_v^{1/3}$ for nodes $v$ living in $t$. So we should focus on the product of ${\bf M}_{\gamma_{j-1}} {\bf M}_{\gamma_j} {\bf M}_{\gamma_{j+1}}$ over all triples $(\gamma_{j-1},\gamma_j,\gamma_{j+1})$ associated to triangle segments $\gamma_j$ living in $t$. The diagonal entries of this product are $\prod_v X_v^{\sum_j {\rm a}_v(\gamma_j)}$, $\prod_v X_v^{\sum_j (-{\rm a}_v(\gamma_j)+{\rm a}_v(\ol{\gamma_j}))}$, and $\prod_v X_v^{\sum_j (-{\rm a}_v(\ol{\gamma_j}))}$, in this order, where $\prod_v$ is over all nodes $v$ living in $t$, and the sum $\sum_j$ is over all $j$'s such that $\gamma_j$ is a triangle segment in $t$. By Lem.\ref{lem:disjoint_union_with_peripheral_curve} we have $\sum_j {\rm a}_v(\gamma_j) = {\rm a}_v(\cup_j \gamma_j) = {\rm a}_v(\gamma \cap \wh{t}\,\,)$ and $\sum_j {\rm a}_v(\ol{\gamma_j}) = {\rm a}_v(\cup_j \ol{\gamma_j}) = {\rm a}_v(\ol{\gamma}\cap \wh{t}\,\,)$. Meanwhile, in view of \redfix{the} definition of the tropical coordinates, we can see that ${\rm a}_v(\gamma \cap \wh{t}\,\,) = {\rm a}_v(\gamma)$ and ${\rm a}_v(\ol{\gamma}\cap \wh{t}\,\,) = {\rm a}_v(\ol{\gamma})$ for these $v$'s. Thus, we showed that, for each node $v$ of $Q_\Delta$ living in each triangle $t$ of $\Delta$, hence for each node $v$ in $Q_\Delta$, the powers of $X_v^{1/3}$ in the monomials appearing as the three diagonal entries of ${\bf M}_\gamma$ are ${\rm a}_v(\ell)$, $-{\rm a}_v(\ell)+{\rm a}_v(\ol{\ell})$, and $-{\rm a}_v(\ol{\ell})$, in this order, as desired in eq.\eqref{eq:peripheral_monodromy_matrix}. We showed this statement for any chosen triangle $t$ of $\Delta$. For any other triangle $t'$, one might have to choose a different basepoint of $\gamma$ for the above arguments to work, so that in the new resulting monodromy matrix ${\bf M}_\gamma'$, the diagonal entries have correct powers for $X_v^{1/3}$ for all nodes $v$ living in $t'$. As mentioned in the beginning of the proof, ${\bf M}_\gamma$ and ${\bf M}_\gamma'$ have same diagonal entries. This finishes the proof for the case when $\gamma$ is a positively oriented peripheral loop. 

\vs

When $\gamma$ is a negatively oriented peripheral loop surrounding $p$, the proof goes similarly, using the triples $(\gamma_{j-1},\gamma_j,\gamma_{j+1})$ for triangle segments $\gamma_j$ living in $\wh{t}$. Now $\gamma_j$ is a right turn, so we can assume it goes from the side $e_{\alpha+1}$ to $e_\alpha$ of $\wh{t}$. By (MM1) and (MM3), the product ${\bf M}_{\gamma_{j-1}} {\bf M}_{\gamma_j} {\bf M}_{\gamma_{j+1}}$ now looks
\begin{align*}
 & \hspace{-10mm} \smallmatthree{ X_{v_{e_{\alpha+1},2}}^{1/3} X_{v_{e_{\alpha+1},1}}^{2/3} }{ 0 }{ 0 }{ 0 }{ X_{v_{e_{\alpha+1},2}}^{1/3} X_{v_{e_{\alpha+1},1}}^{-1/3} }{0}{0}{0}{X_{v_{e_{\alpha+1},2}}^{-2/3} X_{v_{e_{\alpha+1},1}}^{-1/3}}
\smallmatthree{X_{v_t}^{1/3}}{0}{0}{*}{X_{v_t}^{1/3}}{0}{*}{*}{X_{v_t}^{-2/3}}
\smallmatthree{ X_{v_{e_{\alpha},1}}^{1/3} X_{v_{e_{\alpha},2}}^{2/3} }{ 0 }{ 0 }{ 0 }{ X_{v_{e_{\alpha},1}}^{1/3} X_{v_{e_{\alpha},2}}^{-1/3} }{0}{0}{0}{X_{v_{e_{\alpha},1}}^{-2/3} X_{v_{e_{\alpha},2}}^{-1/3}} \\
&\hspace{-10mm} = \smallmatthree{\prod_v X_v^{{\rm a}_v(\gamma_j)}}{0}{0}{*}{\prod_v X_v^{-{\rm a}_v(\gamma_j)+{\rm a}_v(\ol{\gamma_j})}}{0}{*}{*}{\prod_v X_v^{-{\rm a}_v(\ol{\gamma_j})}}
\end{align*}
with the last equality holding in view of the tropical coordinate values as in eq.\eqref{eq:tropical_coordinates_for_a_right_turn_W} with $W_{\alpha+1,\alpha} = \gamma_j$ and eq.\eqref{eq:tropical_coordinates_for_a_left_turn_W} with $W_{\alpha,\alpha+1}=\ol{\gamma_j}$, where $\prod_v$ is taken over seven nodes of $Q_\Delta$ living in the triangle $\wh{t}$ (or $t$). The rest of the arguments goes the same.  \qed

\vs

As seen in eq.\eqref{eq:pi_p_i_investigated}, the three diagonal entries of ${\bf M}_\gamma$ for a peripheral loop $\gamma$ are the sought-for puncture functions on $\mathscr{X}^+_{{\rm PGL}_3,\frak{S}}$, corresponding to the regular functions $(\pi_p)_i$ on $\mathscr{X}_{{\rm SL}_3,\frak{S}}$: for a single-component ${\rm SL}_3$-lamination $\ell_p$ (resp. $\ol{\ell_p}$) consisting of a positively oriented (resp. negatively oriented) peripheral loop surrounding $p$ with weight $1$, we let
$$
(\pi_p)^+_1 := \underset{v\in \mathcal{V}(Q_\Delta)}{\textstyle \prod} X_v^{{\rm a}_v(\ell_p)}, \quad
(\pi_p)^+_2 := \underset{v\in \mathcal{V}(Q_\Delta)}{\textstyle \prod} X_v^{-{\rm a}_v(\ell_p)+{\rm a}_v(\ol{\ell_p})}, \quad
(\pi_p)^+_3 := \underset{v\in \mathcal{V}(Q_\Delta)}{\textstyle \prod} X_v^{-{\rm a}_v(\ol{\ell_p})},
$$
defined as smooth functions on the smooth manifold $\mathscr{X}^+_{{\rm PGL}_3,\frak{S}}$. The following statement is not trivial, but is immediate from definitions.
\begin{lemma}
Each of these functions $(\pi_p)^+_i$ on $\mathscr{X}^+_{{\rm PGL}_3,\frak{S}}$ does not depend on the choice of an ideal triangulation $\Delta$. \qed
\end{lemma}
For example, if $\Delta'$ is any other ideal triangulation, then $\prod_{v\in \mathcal{V}(Q_\Delta)} X_v^{{\rm a}_v(\ell_p)} = \prod_{v' \in \mathcal{V}(Q_{\Delta'})} {X_{v'}'}^{{\rm a}_{v'}'(\ell_p)}$.

\vs

Recall the $A_2$-bangles basis of $\mathscr{O}(\mathscr{X}_{{\rm SL}_3,\frak{S}})$ constructed by the map $\mathbb{I}_{{\rm SL}_3} : \mathscr{A}_{\rm L}(\frak{S};\mathbb{Z}) \to \mathscr{O}(\mathscr{X}_{{\rm SL}_3,\frak{S}})$ in Def.\ref{def:I_SL3}. For each $\ell \in \mathscr{A}_{\rm L}(\frak{S};\mathbb{Z})$, the function $\mathbb{I}_{{\rm SL}_3}(\ell)$ is constructed by products (and powers) and $\mathbb{Z}$-linear combinations of the trace-of-monodromy functions along loops $\gamma$ and the puncture functions $(\pi_p)_i$. Hence, now this function can be translated as a smooth function on the manifold $\mathscr{X}^+_{{\rm PGL}_3,\frak{S}}$ using $f^+_\gamma$ and $(\pi_p)^+_i$'s, which is given for each ideal triangulation $\Delta$ as a Laurent polynomial in $\{X_v^{1/3} \, | \,v\in \mathcal{V}(Q_\Delta)\}$ with integer coefficients. Denote this function by
$$
\mathbb{I}_{{\rm PGL}_3}^+(\ell) \in C^\infty(\mathscr{X}_{{\rm PGL}_3,\frak{S}}^+),
$$
\redfix{which we refer to as a \ul{\em basic semi-regular function}.} In particular, if $\ell$ consists only of peripheral loops with arbitrary integer weights, then we have
\begin{align}
\label{eq:I_plus_peripheral}
\mathbb{I}^+_{{\rm PGL}_3}(\ell) = \underset{v \in \mathcal{V}(Q_\Delta)}{\textstyle \prod} X_v^{{\rm a}_v(\ell)}.
\end{align}

\subsection{A basis of the ring of regular functions on $\mathscr{X}_{{\rm PGL}_3,\frak{S}}$: the first main theorem}
\label{subsec:basis}

We go back to the strategy set out in \S\ref{subsec:lifting_PGL3_to_SL3}. Let $f\in \mathscr{O}(\mathscr{X}_{{\rm PGL}_3,\frak{S}})$. By eq.\eqref{eq:P_star} we get $P^* f \in \mathscr{O}(\mathscr{X}_{{\rm SL}_3,\frak{S}})$. By Prop.\ref{prop:I_SL3_is_a_basis}(1), we have
\begin{align}
\label{eq:P_star_f}
P^*f = \underset{\ell \in \mathscr{A}_{\rm L}(\frak{S};\mathbb{Z})}{\textstyle \sum} c_\ell(f)  \, \mathbb{I}_{{\rm SL}_3}(\ell)
\end{align}
for some $c_\ell(f) \in \mathbb{Q}$, which are zero for all but finitely many $\ell \in \mathscr{A}_{\rm L}(\frak{S};\mathbb{Z})$. Evaluating at the field $\mathbb{R}$, we view $P^*f$ and each $\mathbb{I}_{{\rm SL}_3}(\ell)$ as  functions on $\mathscr{X}_{{\rm SL}_3,\frak{S}}(\mathbb{R})$. Pulling back by the map in eq.\eqref{eq:section_of_PGL3_to_SL3}, these can be viewed as functions on $\mathscr{X}^+_{{\rm PGL}_3,\frak{S}}$. The pullback of $P^*f$ on $\mathscr{X}^+_{{\rm PGL}_3,\frak{S}}$ is just $f$ evaluated at the semi-field $\mathbb{R}_{>0}$, and the pullback of each $\mathbb{I}_{{\rm SL}_3}(\ell)$ is what we denoted by $\mathbb{I}^+_{{\rm PGL}_3}(\ell)$. For any ideal triangulation $\Delta$, since $f$ is regular on the cluster $\mathscr{X}$-chart of $\mathscr{X}_{{\rm PGL}_3,\frak{S}}$ for $\Delta$, it can be written as a Laurent polynomial in the variables $\{X_v \, | \, v\in \mathcal{V}(Q_\Delta)\}$ with integer coefficients. By evaluating at $\mathbb{R}_{>0}$, this Laurent polynomial expression can be viewed as a function on $\mathscr{X}^+_{{\rm PGL}_3,\frak{S}}$. On the other hand, this Laurent polynomial function on $\mathscr{X}^+_{{\rm PGL}_3,\frak{S}}$ must equal the function $\sum_\ell c_\ell(f) \, \mathbb{I}^+_{{\rm PGL}_3}(\ell)$, which is a priori a Laurent polynomial in $\{X_v^{1/3} \, | \, v\in \mathcal{V}(Q_\Delta)\}$ with integer coefficients. In our investigation of when this becomes a Laurent polynomial in $\{X_v \,|\, v\in\mathcal{V}(Q_\Delta)\}$, what play crucial roles are the highest term of each basic semi-regular function $\mathbb{I}^+_{{\rm PGL}_3}(\ell)$, and the {\em congruence} property of all (cube root) Laurent monomial terms for $\mathbb{I}^+_{{\rm PGL}_3}(\ell)$.

\begin{definition}[partial ordering and congruence on Laurent monomials]
Let $\Delta$ be an ideal triangulation of a punctured surface $\frak{S}$.

\vs

$\bullet$ On the set of all Laurent monomials in $\{X_v^{1/3} \, | \, v\in \mathcal{V}(Q_\Delta)\}$, define the partial ordering as follows: for $(a_v)_{v\in \mathcal{V}(Q_\Delta)}, (b_v)_{v\in \mathcal{V}(Q_\Delta)} \in (\frac{1}{3}\mathbb{Z})^{\mathcal{V}(Q_\Delta)}$,
$$
\textstyle \prod_v X_v^{a_v} \succ \prod_v X_v^{b_v} \qquad \overset{\rm def.}{\Longleftrightarrow} \qquad a_v \ge b_v, ~ \forall v\in \mathcal{V}(Q_\Delta).
$$
By convention, the zero monomial is set to be of the lowest ordering, i.e. $\prod_v X_v^{a_v} \succ 0$.

\vs

$\bullet$ For $(a_v)_{v\in \mathcal{V}(Q_\Delta)}, (b_v)_{v\in \mathcal{V}(Q_\Delta)} \in (\frac{1}{3}\mathbb{Z})^{\mathcal{V}(Q_\Delta)}$, we say 
$$
\mbox{$\prod_v X_v^{a_v}$ and $\prod_v X_v^{b_v}$ are \ul{\em congruent} to each other}
\qquad \overset{\rm def.}{\Longleftrightarrow} \qquad
\mbox{$a_v - b_v \in \mathbb{Z}$, $\forall v\in \mathcal{V}(Q_\Delta)$.}
$$
\end{definition}

\begin{proposition}[highest term of basic semi-regular function]
\label{prop:highest_term}
Let $\Delta$ be an ideal triangulation of a punctured surface $\frak{S}$.  For each $\ell \in \mathscr{A}_{\rm L}(\frak{S};\mathbb{Z})$, the basic semi-regular function $\mathbb{I}^+_{{\rm PGL}_3}(\ell) \in C^\infty(\mathscr{X}^+_{{\rm PGL}_3,\frak{S}})$ can be written as a Laurent polynomial in $\{X_v^{1/3} \, | \, v\in\mathcal{V}(Q_\Delta)\}$ with integer coefficients such that the monomial $\prod_{v\in \mathcal{V}(Q_\Delta)} X_v^{{\rm a}_v(\ell)}$ appears with coefficient $1$ and is the unique Laurent monomial having the highest partial order among all Laurent monomials appearing in this expression.
\end{proposition}

\begin{proposition}[congruence of terms of a basic regular function]
\label{prop:congruence-compatibility_of_terms_of_each_basic_regular_function}
Let $\frak{S}$ be a triangulable punctured surface. For each $\ell \in \mathscr{A}_{\rm L}(\frak{S};\mathbb{Z})$, the basic semi-regular function $\mathbb{I}^+_{{\rm PGL}_3}(\ell)$ on $\mathscr{X}^+_{{\rm PGL}_3,\frak{S}}$ satisfies the following, for each ideal triangulation $\Delta$ of $\frak{S}$:
\begin{align*}
\mathbb{I}^+_{{\rm PGL}_3}(\ell) ~\in~ ({\textstyle \prod}_{v\in \mathcal{V}(Q_\Delta)} X_v^{{\rm a}_v(\ell)}) \cdot \mathbb{Z}[\{X_v^{\pm 1} \, | \, v\in \mathcal{V}(Q_\Delta)\}]\redfix{.}
\end{align*}
That is, $\mathbb{I}^+_{{\rm PGL}_3}(\ell)$ can be written as a Laurent polynomial in $\{X_v^{1/3}\,|\, v\in\mathcal{V}(Q_\Delta)\}$ with integer coefficients such that all Laurent monomials appearing are congruent to each other.
\end{proposition}
In fact, proofs of these two propositions are much more involved than it might look at the first glance, so we postpone them until the next section. In the present section, let's assume them.

\begin{corollary}[congruence and integrality of powers]
\label{cor:congruence_and_integrality_of_powers}
Let $\frak{S}$ be a punctured surface. 
Let $(c_\ell)_{\ell \in \mathscr{A}_{\rm L}(\frak{S};\mathbb{Z})} \in \mathbb{Z}^{\mathscr{A}_{\rm L}(\frak{S};\mathbb{Z})}$, where $c_\ell$'s are zero for all but finitely many $\ell$'s. For any ideal triangulation $\Delta$ of $\frak{S}$, \, $\sum_\ell c_\ell \, \mathbb{I}^+_{{\rm PGL}_3}(\ell) \in \mathbb{Z}[\{X_v^{\pm 1/3} \, | \, v\in \mathcal{V}(Q_\Delta)\}]$ belongs to $\mathbb{Z}[\{X_v^{\pm 1} \, | \, v\in \mathcal{V}(Q_\Delta)\}]$ if and only if $c_\ell=0$ for all $\ell \in \mathscr{A}_{\rm L}(\frak{S};\mathbb{Z})$ not belonging to $\mathscr{A}_\Delta(\mathbb{Z}^{\redfix{T}})$ (Def.\ref{def:Delta-congruent}), i.e. $c_\ell=0$ for all $\ell$ such that ${\rm a}_v(\ell)\in \frac{1}{3}\mathbb{Z}$ does not belong to $\mathbb{Z}$ for at least one $v\in \mathcal{V}(Q_\Delta)$.
\end{corollary}

{\it Proof of Cor.\ref{cor:congruence_and_integrality_of_powers}.} Let $f^+ := \sum_{\ell\in\mathscr{A}_{\rm L}(\frak{S};\mathbb{Z})} c_\ell \, \mathbb{I}_{{\rm PGL}_3}^+(\ell)$ be a function on $\mathscr{X}^+_{{\rm PGL}_3,\frak{S}}$, with $c_\ell \in \mathbb{Z}$, which are zero for all but finitely many $\ell$'s. One direction is easy. Suppose $c_\ell = 0$ whenever $\ell$ is not in $\mathscr{A}_{\Delta}(\mathbb{Z}^{\redfix{T}})$, so we can write $f^+ := \sum_{\ell\in\mathscr{A}_{{\rm SL}_3,\frak{S}} (\mathbb{Z}^{\redfix{T}})} c_\ell \, \mathbb{I}_{{\rm PGL}_3}^+(\ell)$. By Prop.\ref{prop:congruence-compatibility_of_terms_of_each_basic_regular_function}, $\mathbb{I}^+_{{\rm PGL}_3}(\ell)$ belongs to $\mathbb{Z}[\{X_v^{\pm 1} \, | \, v\in \mathcal{V}(Q_\Delta)\}]$ for each $\ell \in \mathscr{A}_{\Delta}(\mathbb{Z}^{\redfix{T}})$. Hence $f^+ \in \mathbb{Z}[\{X_v^{\pm 1} \, | \, v\in \mathcal{V}(Q_\Delta)\}]$.

\vs

Now, for the converse, suppose $f^+ \in \mathbb{Z}[\{X_v^{\pm 1} \, | \, v\in \mathcal{V}(Q_\Delta)\}]$. Recall the partial ordering on the set of all Laurent monomials in $\{X_v^{1/3} \, | \, v\in \mathcal{V}(Q_\Delta)\}$. Choose any ordering on the set $\mathcal{V}(Q_\Delta)$, and consider the induced lexicographic total ordering on the set of all Laurent monomials in $\{X_v^{1/3} \, | \, v\in \mathcal{V}(Q_\Delta)\}$, which is compatible with the previous partial ordering. We expressed each $\mathbb{I}^+_{{\rm PGL}_3}(\ell)$ so that it has the unique Laurent monomial term of highest partial order (Prop.\ref{prop:highest_term}). Among all these highest Laurent monomials appearing in the sum $f^+ = \sum_{\ell\in\mathscr{A}_{\rm L}(\frak{S};\mathbb{Z})} c_\ell \, \mathbb{I}_{{\rm PGL}_3}^+(\ell)$, there must be one with the highest lexicographic ordering; in view of Prop.\ref{prop:highest_term} it is $\prod_v X_v^{{\rm a}_v(\ell_0)}$ (which is the highest term of $\mathbb{I}^+_{{\rm PGL}_3}(\ell_0)$)  for some $\ell_0$ contributing to the sum. This is in fact the unique term of highest lexicographic order, because of the injectivity of the coordinate-system map ${\bf a}_\Delta : \ell \mapsto ({\rm a}_v(\ell))_{v\in \mathcal{V}(Q_\Delta)}$ (Prop.\ref{prop:tropical_coordinates_are_coordinates}). Therefore, in order for $f^+$ to be a function that can be written as a Laurent polynomial in $\{X_v \, | \, v\in \mathcal{V}(Q_\Delta)\}$, it follows that the term $\prod_v X_v^{{\rm a}_v(\ell_0)}$ of the highest lexicographic order must be a Laurent monomial in $\{X_v \, | \, v\in \mathcal{V}(Q_\Delta)\}$, so ${\rm a}_v(\ell_0) \in \mathbb{Z}$ for all $v\in \mathcal{V}(Q_\Delta)$, or equivalently, $\ell_0 \in \mathscr{A}_{\Delta}(\mathbb{Z}^{\redfix{T}})$. By Prop.\ref{prop:congruence-compatibility_of_terms_of_each_basic_regular_function} we know $\mathbb{I}^+_{{\rm PGL}_3}(\ell_0) \in \mathbb{Z}[\{X_v^{\pm 1} \, | \, v\in \mathcal{V}(Q_\Delta)\}]$. Now $f^+ - c_{\ell_0} \, \mathbb{I}^+_{{\rm PGL}_3}(\ell_0)$ equals $ \sum_{\ell\in\mathscr{A}_{\rm L}(\frak{S};\mathbb{Z})\setminus\{\ell_0\}} c_\ell \, \mathbb{I}_{{\rm PGL}_3}^+(\ell)$, and therefore it has fewer summands than $f^+$ (i.e. fewer $\ell$'s contributing to the sum) and it belongs to $\mathbb{Z}[\{X_v^{\pm 1} \, | \, v\in \mathcal{V}(Q_\Delta)\}]$ again. By induction, we get that all $\ell \in \mathscr{A}_{\rm L}(\frak{S};\mathbb{Z})$ contributing to the sum $\sum_{\ell\in\mathscr{A}_{\rm L}(\frak{S};\mathbb{Z})} c_\ell \, \mathbb{I}_{{\rm PGL}_3}^+(\ell)$ must belong to $\mathscr{A}_{\Delta}(\mathbb{Z}^{\redfix{T}})$. \qed

\begin{corollary}[congruent ${\rm SL}_3$-laminations give genuinely regular functions]
\label{cor:congruent_A2-lamination_gives_genuinely_regular_function}
Let $\frak{S}$ be a triangulable punctured surface. Let $\Delta$ be any ideal triangulation of $\frak{S}$. For $\ell \in \mathscr{A}_{\rm L}(\frak{S};\mathbb{Z})$, the function $\mathbb{I}^+_{{\rm PGL}_3}(\ell)$ on $\mathscr{X}^+_{{\rm PGL}_3,\frak{S}}$ can be written as a Laurent polynomial in $\{X_v \, | \, v\in \mathcal{V}(Q_\Delta)\}$ with integer coefficients if and only if $\ell \in\mathscr{A}_{\Delta}(\mathbb{Z}^{\redfix{T}})$.  \qed
\end{corollary}
So, for $\ell \in \mathscr{A}_{\Delta}(\mathbb{Z}^{\redfix{T}})$, $\mathbb{I}^+_{{\rm PGL}_3}(\ell)$ comes from a rational function on $\mathscr{X}_{{\rm PGL}_3,\frak{S}}$ that is regular on the cluster $\mathscr{X}$-chart associated to each ideal triangulation $\Delta$. In fact, this rational function on $\mathscr{X}_{{\rm PGL}_3,\frak{S}}$ is a regular function on the entire moduli space $\mathscr{X}_{{\rm PGL}_3,\frak{S}}$.
\begin{proposition}
\label{prop:congruent_lamination_gives_regular_function}
Let $\Delta$ be an ideal triangulation of a punctured surface $\frak{S}$. For $\ell \in \mathscr{A}_{\Delta}(\mathbb{Z}^{\redfix{T}})$, the basic semi-regular function $\mathbb{I}^+_{{\rm PGL}_3}(\ell) \in C^\infty(\mathscr{X}^+_{{\rm PGL}_3,\frak{S}})$ comes from a regular function on $\mathscr{X}_{{\rm PGL}_3,\frak{S}}$.
\end{proposition}
Prop.\ref{prop:congruent_lamination_gives_regular_function} will be proved in the next subsection \redfix{\S\ref{subsec:mutation_of_basic_regular_functions} and the next section} through several steps. For now, let's assume it. 

\vs

Combining the results so far, we arrive at the first main theorem of the paper.
\begin{theorem}[the first main theorem; the ${\rm SL}_3$-${\rm PGL}_3$ duality map and the $A_2$-bangles basis of $\mathscr{O}(\mathscr{X}_{{\rm PGL}_3,\frak{S}})$]
\label{thm:main}
Let $\frak{S}$ be a triangulable punctured surface. Then the sets $\mathscr{A}_\Delta(\mathbb{Z}^{\redfix{T}}) \subset \mathscr{A}_{\rm L}(\frak{S};\mathbb{Z})$ (Def.\ref{def:Delta-congruent}) for all ideal triangulations $\Delta$ of $\frak{S}$ coincide with each other (i.e. Prop.\ref{prop:congruence_condition_is_indepdent_on_triangulation} holds); denote any one of them by $\mathscr{A}_{{\rm SL}_3,\frak{S}}(\mathbb{Z}^{\redfix{T}})$. Then, there exists a map
$$
\mathbb{I} : \mathscr{A}_{{\rm SL}_3,\frak{S}}(\mathbb{Z}^{\redfix{T}}) \to \mathscr{O}(\mathscr{X}_{{\rm PGL}_3,\frak{S}})
$$
such that
\begin{enumerate}
\item[\rm (1)] $\mathbb{I}$ is injective and the image set $\mathbb{I}(\mathscr{A}_{{\rm SL}_3,\frak{S}}(\mathbb{Z}^{\redfix{T}}))$ forms a basis of $\mathscr{O}(\mathscr{X}_{{\rm PGL}_3,\frak{S}})$, which we call \redfix{the} \ul{\em $A_2$-bangles basis} of $\mathscr{O}(\mathscr{X}_{{\rm PGL}_3,\frak{S}})$.

\item[\rm (2)] For $\ell \in \mathscr{A}_{{\rm SL}_3,\frak{S}}(\mathbb{Z}^{\redfix{T}})$, for any ideal triangulation $\Delta$ of $\frak{S}$, $\mathbb{I}(\ell)$ can be written as a Laurent polynomial in $\{X_v \, | \, v\in \mathcal{V}(Q_\Delta)\}$ with  integer coefficients, with the unique highest Laurent monomial being $\prod_{v\in \mathcal{V}(Q_\Delta)} X_v^{{\rm a}_v(\ell)}$, with coefficient $1$.

\item[\rm (3)] If $\ell \in \mathscr{A}_{{\rm SL}_3,\frak{S}}(\mathbb{Z}^{\redfix{T}})$ consists only of peripheral loops, then for each ideal triangulation $\Delta$, we have $\mathbb{I}(\ell)= \prod_{v\in \mathcal{V}(Q_\Delta)} X_v^{{\rm a}_v(\ell)}$.

\item[\rm (4)] The structure constants of this $A_2$-bangles basis of $\mathscr{O}(\mathscr{X}_{{\rm PGL}_3,\frak{S}})$ are integers. That is, for any $\ell,\ell' \in \mathscr{A}_{{\rm SL}_3,\frak{S}}(\mathbb{Z}^{\redfix{T}})$, we have
\begin{align}
\label{eq:I_structure_constants}
\mathbb{I}(\ell) \, \mathbb{I}(\ell') = \underset{\ell'' \in \mathscr{A}_{{\rm SL}_3,\frak{S}}(\mathbb{Z}^{\redfix{T}})}{\textstyle \sum}  c(\ell,\ell';\ell'') \, \mathbb{I}(\ell'')
\end{align}
where $c(\ell,\ell';\ell'')\in \mathbb{Z}$ and $c(\ell,\ell';\ell'')$ are zero for all but at most finitely many $\ell''$.

\end{enumerate}
\end{theorem}

The rest of this section is devoted to proof of Prop.\ref{prop:congruent_lamination_gives_regular_function} and Thm.\ref{thm:main}; as said, Prop.\ref{prop:highest_term} and Prop.\ref{prop:congruence-compatibility_of_terms_of_each_basic_regular_function} will be proved in the next section.

\subsection{Mutations of basic regular functions}
\label{subsec:mutation_of_basic_regular_functions}

In this subsection we prove Prop.\ref{prop:congruent_lamination_gives_regular_function}. First, recall from Def.\ref{def:O_cl} the notion $\mathscr{O}_{\rm cl}(\mathscr{X}_{{\rm PGL}_3,\frak{S}})$, the ring of all rational functions on $\mathscr{X}_{{\rm PGL}_3,\frak{S}}$ that are regular on {\em all} cluster $\mathscr{X}$-charts. Any element of $\mathscr{O}_{\rm cl}(\mathscr{X}_{{\rm PGL}_3,\frak{S}})$ is {\em universally Laurent} for all cluster $\mathscr{X}$-charts, hence in particular is universally Laurent in the Fock-Goncharov's weaker sense that it is a Laurent polynomial in the cluster $\mathscr{X}$-chart associated to every ideal triangulation $\Delta$, i.e. belongs to ${\bf L}(\mathscr{X}_{{\rm PGL}_3,\frak{S}})$ (Def.\ref{def:intro_L}). We recall the result of Shen:
\begin{proposition}[{\cite[Thm1.1]{Shen}}]
\label{prop:cluster_regular_is_regular}
$\mathscr{O}_{\rm cl}(\mathscr{X}_{{\rm PGL}_3,\frak{S}}) = \mathscr{O} (\mathscr{X}_{{\rm PGL}_3,\frak{S}})$.
\end{proposition}
Shen's result is written in terms of a slightly different moduli space $\mathscr{P}_{{\rm G},\mathbb{S}}$ introduced in \cite{GS19} \cite{GS15}, for a generalized marked surface $\mathbb{S}$. Putting ${\rm G} = {\rm PGL}_3$ and when $\mathbb{S}$ is a punctured surface $\frak{S}$, this moduli space is same as $\mathscr{X}_{{\rm PGL}_3,\frak{S}}$.

\vs

Next, we need the following statement, which follows from \redfix{the result of} Gross, Hacking and Keel \cite{GHK}. It tells us that, to check the universally Laurent condition, it suffices to check it for one cluster chart and for all charts obtained by applying a single mutation to this chart.
\begin{proposition}[{\cite[Thm.3.9]{GHK}},  {\cite[Lem.2.2]{Shen}}]
\label{prop:one-level_mutations_are_enought}
Let $f$ be a rational function on $\mathscr{X}_{{\rm PGL}_3,\frak{S}}$. Let $\Delta$ be an ideal triangulation of a puntured surface $\frak{S}$, and suppose that $f$ is regular on the cluster $\mathscr{X}$-chart for $\Delta$; that is, $f$ is a Laurent polynomial in the cluster $\mathscr{X}$-variables for this chart. If, for every node $v$ of $Q_\Delta$, $f$ is regular on the cluster $\mathscr{X}$-chart obtained from the cluster $\mathscr{X}$-chart for $\Delta$ by applying the mutation at this node, then $f$ belongs to $\mathscr{O}_{\rm cl}(\mathscr{X}_{{\rm PGL}_3,\frak{S}})$.
\end{proposition}

Our strategy to prove Prop.\ref{prop:congruent_lamination_gives_regular_function} is as follows. For $\ell \in \mathscr{A}_{\Delta}(\mathbb{Z}^{\redfix{T}})$, we know that $\mathbb{I}^+_{{\rm PGL}_3}(\ell)$ comes from a rational function on $\mathscr{X}_{{\rm PGL}_3,\frak{S}}$, say $\mathbb{I}_\Delta(\ell)$, that is regular on the cluster $\mathscr{X}$-chart for each ideal triangulation $\Delta$. We fix any triangulation $\Delta$, and will show that if we mutate at any node of $Q_\Delta$, the result is still a Laurent polynomial in the new cluster $\mathscr{X}$-variables. Then by Prop.\ref{prop:one-level_mutations_are_enought} it follows $\mathbb{I}_\Delta(\ell) \in \mathscr{O}_{\rm cl}(\mathscr{X}_{{\rm PGL}_3,\frak{S}})$, and in turn by Prop.\ref{prop:cluster_regular_is_regular} we get $\mathbb{I}_\Delta(\ell) \in \mathscr{O}(\mathscr{X}_{{\rm PGL}_3,\frak{S}})$, as desired in Prop.\ref{prop:congruent_lamination_gives_regular_function}.

\vs

In order to study the effect of mutation, we study the basic semi-regular functions $\mathbb{I}^+_{{\rm PGL}_3}(\ell)$ for $\ell \in \mathscr{A}_{\rm L}(\frak{S};\mathbb{Z})$, which are functions on the manifold $\mathscr{X}^+_{{\rm PGL}_3,\frak{S}}$ that can be written as Laurent polynomials in the cube roots of (positive real evaluations of) cluster $\mathscr{X}$-coordinate functions. We investigate the effect of mutation for these functions explicitly. We begin by mutating at a node lying in the interior of a triangle.

\begin{figure}[htbp!]
\vspace{-1mm}
\begin{center}
{\scalebox{1.0}{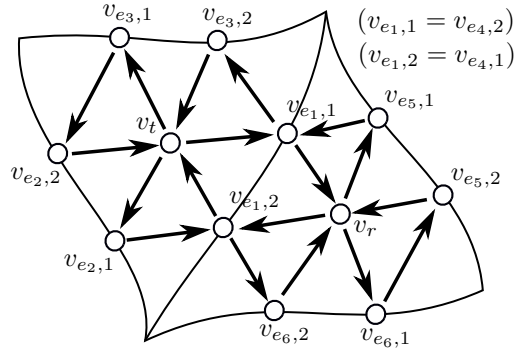}} 
\end{center}
\vspace{-7mm}
\caption{Quiver $Q_\Delta$ for a triangulation $\Delta$, drawn for two adjacent triangles $t$ and $r$}
\vspace{-4mm}
\label{fig:mutate_two_triangles}
\end{figure}

\begin{proposition}[mutation of basic semi-regular function at interior node of triangle]
\label{prop:mutation_of_basic_semi-regular_function_at_interior_node_of_triangle}
Let $\Delta$ be any ideal triangulation of a triangulable punctured surface $\frak{S}$. Consider the cluster $\mathscr{X}$-chart associated to $\Delta$, and mutate it at a node of $Q_\Delta$ lying in the interior of some triangle of $\Delta$. Denote the resulting quiver by $Q'$, while we naturally identify the sets of nodes $\mathcal{V}(Q_\Delta)$ and $\mathcal{V}(Q')$. Denote by $X_v'$ the $\mathscr{X}$-coordinate for the node $v$ of $Q'$ for this new chart obtained as the result of mutation. Then for any $\ell \in \mathscr{A}_{\rm L}(\frak{S};\mathbb{Z})$, we have
$$
\mathbb{I}^+_{{\rm PGL}_3}(\ell) ~\in~ {X'}_{\hspace{-1,5mm}v_t}^{-{\rm a}_{v_t}(\ell)} \, ({\textstyle \prod}_{v\in \mathcal{V}(Q')\setminus\{v_t\}} {X'}_{\hspace{-1,5mm}v}^{{\rm a}_v(\ell)}) \cdot \mathbb{Z}[\{{X'}_{\hspace{-1,5mm}v}^{\pm 1} \, | \, v\in \mathcal{V}(Q')\}].
$$
\end{proposition}

{\it \redfix{Partial proof.}} \redfix{Here we prove the statement only for the case when $\ell$ can be represented by an ${\rm SL}_3$-web without $3$-valent vertices. The general cases will be treated in the next section (\S\ref{subsec:on_the_effect_of_a_single_mutation}), using the machinery we develop in \S\ref{sec:SL3_trace} in order to prove Propositions \ref{prop:highest_term} and \ref{prop:congruence-compatibility_of_terms_of_each_basic_regular_function}.}

\vs

Consider mutation at $v_t$ of some triangle $t$ of $\Delta$. Let $e_1,e_2,e_3$ be the sides of $t$ appearing in this order clockwise along $\partial t$. Since we assumed that the triangulation $\Delta$ is regular (Def.\ref{def:regular_triangulation}), no two of these three sides are identified with each other. For each $e_\alpha$, let $v_{e_\alpha,1}$ and $v_{e_\alpha,2}$ be the nodes of $Q_\Delta$ so that $v_{e_\alpha,1} \to v_{e_\alpha,2}$ matches the clockwise orientation of $\partial t$; see the triangle on the left of Fig.\ref{fig:mutate_two_triangles}. The cluster $\mathscr{X}$-variables change under the mutation at $v_t$, by the formulas (eq.\eqref{eq:X-mutation_formula})
$$
X_{v_t}' = X_{v_t}^{-1}, \qquad
X_{v_{e_\alpha},1}' = X_{v_{e_\alpha},1}(1+X_{v_t}), \qquad
X_{v_{e_\alpha},2}' = X_{v_{e_\alpha},2}(1+X_{v_t}^{-1})^{-1}.
$$
for $\alpha=1,2,3$, with $X_v' = X_v$ for all $v\in \mathcal{V}(Q') = \mathcal{V}(Q_\Delta)$ not appearing in $t$; so, seven variables change. Writing the old variables in new variables:
\begin{align}
\label{eq:old_variable_as_new_new_variables1}
\hspace*{-2mm} \begin{array}{l}
X_{v_t} = {X'}_{\hspace{-1,5mm}v_t}^{-1}, \quad
X_{v_{e_\alpha},1} = X_{v_{e_\alpha},1}' {X'}_{\hspace{-1,5mm}v_t} (1+{X'}_{\hspace{-1,5mm}v_t})^{-1}, \quad
X_{v_{e_\alpha},2} = X_{v_{e_\alpha},2}'(1+X_{v_t}'), \quad \mbox{for $\alpha=1,2,3$}, \\
X_v = X_v', \quad \mbox{for all $v\in \mathcal{V}(Q_\Delta)$ not appearing in $t$.}
\end{array} 
\end{align}

\begin{figure}[htbp!]
\begin{center}
{\scalebox{1.0}{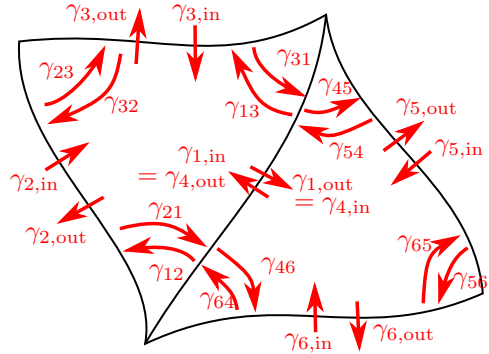}} 
\end{center}
\vspace{-5mm}
\caption{Segments for two triangles}
\vspace{-3mm}
\label{fig:mutate_two_triangles_segments}
\end{figure}

Let $\gamma$ be an oriented loop, decomposed into concatenation $\gamma_1.\gamma_2.\cdots.\gamma_N$ of triangle segments and juncture segments, as in eq.\eqref{eq:gamma_concatenation}. \redfix{For our situation, we may assume that there is no U-turn; in particular, $\gamma$ is non-contractible.} We study the monodromy matrix ${\bf M}_\gamma = {\bf M}_{\gamma_1} \cdots {\bf M}_{\gamma_N}$. A triangle segment in $t$, going from edge $e_{\alpha_1}$ to $e_{\alpha_2}$ is denoted by $\gamma_{\alpha_1\alpha_2}$. A juncture segment at the side $e_\alpha$ of $t$ coming out of this triangle $t$ is denoted by $\gamma_{\alpha,{\rm out}}$, and that going into $t$ by $\gamma_{\alpha,{\rm in}}$. See the left triangle $t$ of Fig.\ref{fig:mutate_two_triangles_segments}.

\vs

We should consider all possibilities of concatenations of segments in $t$ forming a `complete' concatenation in this \redfix{triangle}:
\begin{align}
\label{eq:cases_of_complete_concatenations1}
{\gamma_{\alpha_1,{\rm in}}} . {\gamma_{\alpha_1 \alpha_2}} . {\gamma_{\alpha_2,{\rm out}}}, \qquad \alpha_1,\alpha_2 \in \{1,2,3\}, \quad \alpha_1 \neq \alpha_2.
\end{align}
For each case of a complete concatenation, we should compute the effect of mutation on the product of corresponding monodromy matrices ${\bf M}_\cdot$ defined in (MM1)--(MM3) of \S\ref{subsec:lifting_PGL3_to_SL3}. We use the normalized matrices $\til{\bf M}_\cdot$, defined as follows. For a juncture segment $\gamma_i$ as in Fig.\ref{fig:juncture_segment}, let
\begin{align}
\label{eq:til_juncture_segment_matrix}
\til{\bf M}_{\gamma_i} := {\rm diag}(1, X_2^{-1}, X_1^{-1}X_2^{-1}).
\end{align}
Define the left and the right turn matrices for triangle $t$ as
\begin{align}
\label{eq:til_triangle_segment_matrix_t}
\til{\bf M}^{\rm left}_t = \smallmatthree{1}{~1+X_{v_t}^{-1}~}{X_{v_t}^{-1}}{0}{X_{v_t}^{-1}}{X_{v_t}^{-1}}{0}{0}{X_{v_t}^{-1}}, \qquad
\til{\bf M}^{\rm right}_t = \smallmatthree{1}{0}{0}{1}{1}{0}{1}{~1+X_{v_t}^{-1}~}{X_{v_t}^{-1}},
\end{align}
So these $\til{\bf M}_\cdot$ matrices are obtained by dividing by the $(1,1)$-th entry of the corresponding original matrix ${\bf M}_\cdot$.  A good way to keep track of the $(1,1)$-th entries of the original monodromy matrices is using ${\rm SL}_3$-webs and their tropical coordinates. Let $W_{\alpha_1\alpha_2;t}$ be an ${\rm SL}_3$-web in $t$ consisting just of one corner arc, from edge $e_{\alpha_1}$ to $e_{\alpha_2}$. In terms of the segments of $\gamma$, this can be viewed as a concatenation of (part of) $\gamma_{\alpha_1,{\rm in}}$, then $\gamma_{\alpha_1\alpha_2}$, then (part of) $\gamma_{\alpha_2,{\rm out}}$. 
When writing the corresponding product of matrices ${\bf M}_{W_{\alpha_1\alpha_2;t}} := {\bf M}_{\gamma_{\alpha_1,{\rm in}}} {\bf M}_{\gamma_{\alpha_1 \alpha_2}} {\bf M}_{\gamma_{\alpha_2,{\rm out}}}$ as the product of normalized matrices $\til{\bf M}_{W_{\alpha_1\alpha_2;t}} := \til{\bf M}_{\gamma_{\alpha_1,{\rm in}}} \til{\bf M}_{\gamma_{\alpha_1 \alpha_2}} \til{\bf M}_{\gamma_{\alpha_2,{\rm out}}}$ times some factor, this factor is the product of $(1,1)$-th entries of ${\bf M}_{\gamma_{\alpha_1,{\rm in}}}$, ${\bf M}_{\gamma_{\alpha_1 \alpha_2}}$, ${\bf M}_{\gamma_{\alpha_2,{\rm out}}}$, and one can observe that the power of each generator $X_v$ in this factor equals the tropical coordinate ${\rm a}_v(W_{\alpha_1 \alpha_2;t})$ of the ${\rm SL}_3$-web $W_{\alpha_1 \alpha_2;t}$ (see Fig.\ref{fig:DS_coordinates}), i.e. this factor equals $\prod_{v\in \mathcal{V}(Q_\Delta)\cap t} X_v^{{\rm a}_v(W_{\alpha_1 \alpha_2;t})}$; this was already seen in the proof of Prop.\ref{prop:peripheral_monodromy}.

\vs

Now, we will investigate the effect of mutation on the $(1,1)$-entry-factor and on the (products of) normalized matrices $\til{\bf M}_{W_{\alpha_1\alpha_2;t}}$. Note that, using this language of ${\rm SL}_3$-webs in $t$, the cases to be checked are $W_{\alpha_1\alpha_2;t}$ with $\alpha_1,\alpha_2\in\{1,2,3\}$, $\alpha_1\neq \alpha_2$. For convenience when studying the effect of mutation, we let
$$
{\bf X}_t := 1 + X'_{v_t},
$$
By eq.\eqref{eq:old_variable_as_new_new_variables1}, the effect of mutation on a monomial $\prod_{v\in\mathcal{V}(Q_\Delta)} X_v^{k_v/3}$, for $(k_v)_v \in \mathbb{Z}^{\mathcal{V}(Q_\Delta)}$, is
\begin{align*}
\underset{v\in\mathcal{V}(Q_\Delta)}{\textstyle \prod} \, X_v^{k_v/3} = 
{X'}_{\hspace{-1,5mm}v_t}^{(-k_{v_t} + \sum_{\alpha=1}^3 k_{v_{e_\alpha,1}})/3} {\bf X}_t^{\sum_{\alpha=1}^3 (-k_{v_{e_\alpha},1} + k_{v_{e_\alpha},2}) /3} \underset{v\in\mathcal{V}(Q')\setminus\{v_t\}}{\textstyle \prod} {X'}_{\hspace{-1,5mm}v}^{k_v/3}.
\end{align*}
For all the cases of ${\rm SL}_3$-webs $W = W_{\alpha_1\alpha_2;t}$ to be checked, we let $k_v = 3{\rm a}_v(W)$ for nodes $v$ of $Q_\Delta$ living in $t$, let $k_v=0$ for other $v\in \mathcal{V}(Q_\Delta)$. Note from eq.\eqref{eq:d_t_in_terms_of_a} that $\sum_{\alpha=1}^3 (-k_{v_{e_\alpha},1} + k_{v_{e_\alpha},2})\redfix{/3}$ equals ${\rm d}_t(W)$. 

\vs

In fact, by cyclic symmetry, it suffices to check only two cases $W_{12;t}$ and $W_{13;t}$. As can be seen in eq.\eqref{eq:tropical_coordinates_for_a_left_turn_W} and eq.\eqref{eq:tropical_coordinates_for_a_right_turn_W}, we have ${\rm d}_t(W)=0$ for these ${\rm SL}_3$-webs, as well as $\sum_{\alpha=1}^3 k_{v_{e_\alpha,1}}/3 = \sum_{\alpha=1}^3 {\rm a}_{v_{e_\alpha,1}}(W) \in \mathbb{Z}$ (which appears in the power of $X'_{v_t}$). 

\vs

We now investigate the effect of mutation on the normalized monodromy matrices. By eq.\eqref{eq:old_variable_as_new_new_variables1}, the left and the right turn matrices mutate as:
\begin{align*}
\til{\bf M}^{\rm left}_t = \smallmatthree{1}{~1+X'_{v_t}~}{X'_{v_t}}{0}{X'_{v_t}}{X'_{v_t}}{0}{0}{X'_{v_t}}
= \smallmatthree{1}{~{\bf X}_t~}{X'_{v_t}}{0}{X'_{v_t}}{X'_{v_t}}{0}{0}{X'_{v_t}}, \qquad
\til{\bf M}^{\rm right}_t = \smallmatthree{1}{0}{0}{1}{1}{0}{1}{~1+X'_{v_t}~}{X'_{v_t}}
= \smallmatthree{1}{0}{0}{1}{1}{0}{1}{~{\bf X}_t~}{X'_{v_t}}.
\end{align*}
The edge matrices mutate as:
\begin{align*}
& \til{\bf M}_{\gamma_{\alpha,{\rm in}}} = \smallmatthree{1}{0}{0}{0}{X_{v_{e_\alpha,1}}^{-1} }{0}{0}{0}{ X_{v_{e_\alpha,2}}^{-1} X_{v_{e_\alpha,1}}^{-1} }
= \smallmatthree{1}{0}{0}{0}{ {X'}_{\hspace{-1,5mm}v_{e_\alpha,1}}^{-1} {X'}_{\hspace{-1,5mm}v_t}^{-1} {\bf X}_t }{0}{0}{0}{ {X'}_{\hspace{-1,5mm}v_{e_\alpha,2}}^{-1} \cancel{ {\bf X}_t^{-1} } {X'}_{\hspace{-1,5mm}v_{e_\alpha,1}}^{-1} {X'}_{\hspace{-1,5mm}v_t}^{-1} \cancel{ {\bf X}_t } }, \\
& \til{\bf M}_{\gamma_{\alpha,{\rm out}}} = \smallmatthree{1}{0}{0}{0}{X_{v_{e_\alpha,2}}^{-1} }{0}{0}{0}{ X_{v_{e_\alpha,1}}^{-1} X_{v_{e_\alpha,2}}^{-1} }
= \smallmatthree{1}{0}{0}{0}{ {X'}_{\hspace{-1,5mm}v_{e_\alpha,2}}^{-1} {\bf X}_t^{-1} }{0}{0}{0}{ {X'}_{\hspace{-1,5mm}v_{e_\alpha,2}}^{-1} {X'}_{\hspace{-1,5mm}v_{e_\alpha,1}}^{-1} {X'}_{\hspace{-1,5mm}v_t}^{-1} }.
\end{align*}

\vs

What we would like to check is, for each ${\rm SL}_3$-web $W=W_{12;t}$ and $W_{13;t}$, that the corresponding product of normalized monodromy matrices lives in ${\rm GL}_3(\mathbb{Z}[\{{X'}_{\hspace{-1,5mm}v}^{\pm 1}\,|\,v\in\mathcal{V}(Q')\}])$, i.e. the entries are \ul{\em $X'$-Laurent}, i.e. Laurent polynomials in $\{{X'}_{\hspace{-1,5mm}v}\,|\,v\in\mathcal{V}(Q')\}$ with integer coefficients. The point is to make sure that there is no negative powers of ${\bf X}_t$. For $W_{12;t}$, the corresponding product of normalized matrices is
\begin{align*}
\hspace*{-3mm} \til{\bf M}_{\gamma_{1,{\rm in}}} \til{\bf M}^{\rm left}_t \til{\bf M}_{\gamma_{2,{\rm out}}}
= 
\smallmatthree{1}{0}{0}{0}{ {X'}_{\hspace{-1,5mm}v_{e_1,1}}^{-1} {X'}_{\hspace{-1,5mm}v_t}^{-1} {\bf X}_t }{0}{0}{0}{ {X'}_{\hspace{-1,5mm}v_{e_1,2}}^{-1} {X'}_{\hspace{-1,5mm}v_{e_1,1}}^{-1} {X'}_{\hspace{-1,5mm}v_t}^{-1} }
\smallmatthree{1}{~{\bf X}_t~}{X'_{v_t}}{0}{X'_{v_t}}{X'_{v_t}}{0}{0}{X'_{v_t}}
\smallmatthree{1}{0}{0}{0}{ {X'}_{\hspace{-1,5mm}v_{e_2,2}}^{-1} {\bf X}_t^{-1} }{0}{0}{0}{ {X'}_{\hspace{-1,5mm}v_{e_2,2}}^{-1} {X'}_{\hspace{-1,5mm}v_{e_2,1}}^{-1} {X'}_{\hspace{-1,5mm}v_t}^{-1} },
\end{align*}
and when we multiply these matrices, it is easy to see that in each entry there is no ${\bf X}_t^{-1}$ left, so that it is $X'$-Laurent. For $W_{13;t}$, the corresponding product is
\begin{align*}
\hspace*{-3mm}
\til{\bf M}_{\gamma_{1,{\rm in}}} \til{\bf M}^{\rm right}_t \til{\bf M}_{\gamma_{3,{\rm out}}}
= 
\smallmatthree{1}{0}{0}{0}{ {X'}_{\hspace{-1,5mm}v_{e_1,1}}^{-1} {X'}_{\hspace{-1,5mm}v_t}^{-1} {\bf X}_t }{0}{0}{0}{ {X'}_{\hspace{-1,5mm}v_{e_1,2}}^{-1} {X'}_{\hspace{-1,5mm}v_{e_1,1}}^{-1} {X'}_{\hspace{-1,5mm}v_t}^{-1} }
\smallmatthree{1}{0}{0}{1}{1}{0}{1}{~{\bf X}_t~}{X'_{v_t}}
\smallmatthree{1}{0}{0}{0}{ {X'}_{\hspace{-1,5mm}v_{e_3,2}}^{-1} {\bf X}_t^{-1} }{0}{0}{0}{ {X'}_{\hspace{-1,5mm}v_{e_3,2}}^{-1} {X'}_{\hspace{-1,5mm}v_{e_3,1}}^{-1} {X'}_{\hspace{-1,5mm}v_t}^{-1} },
\end{align*}
and again, when we multiply these matrices, we see that ${\bf X}_t^{-1}$ is cancelled, so that the entries are $X'$-Laurent.

\vs

Let's summarize the results so far. Writing the trace-of-monodromy $f^+_\gamma = {\rm tr}({\bf M}_{\gamma_1} \cdots {\bf M}_{\gamma_N})$ along oriented \redfix{non-contractible} simple loop $\gamma$ as a Laurent polynomial in the (cube-root) old variables $\{ X_v^{1/3} \, | \, v\in \mathcal{V}(Q_\Delta) \}$, by Prop.\ref{prop:congruence-compatibility_of_terms_of_each_basic_regular_function} we know 
$$
f^+_\gamma \in ( {\textstyle \prod}_{v\in \mathcal{V}(Q_\Delta)} X_v^{{\rm a}_v(\redfix{\gamma})} ) \cdot \mathbb{Z}[\redfix{\{}X_v^{\pm 1} \, | \, v\in \mathcal{V}(Q_\Delta)\}\redfix{]},
$$
\redfix{where $\gamma$ is being viewed as an ${\rm SL}_3$-lamination in $\frak{S}$ (with weight $1$).} 
We investigated the monodromy matrices ${\bf M}_{\gamma_i}$ in terms of new variables $\{ {X'}_{\hspace{-1,5mm}v}^{1/3} \, | \, v\in \mathcal{V}(Q')\}$, and found out that for the entries of the product matrix ${\bf M}_{\gamma_1} \cdots {\bf M}_{\gamma_N}$, the discrepancy between the power of an old variable $X_v$ and the corresponding new variable ${X'}_{\hspace{-1,5mm}v}$ (via natural identification $\mathcal{V}(Q_\Delta) \leftrightarrow \mathcal{V}(Q')$), considered up to integers, occurs only for the node $v_t$ which we are mutating at, where the previous power of $X_{v_t}$ is ${\rm a}_{v_t}(\redfix{\gamma})$ while the new power of ${X'}_{\hspace{-1,5mm}v_t}$ is $-{\rm a}_{v_t}(\redfix{\gamma}) $ (modulo $\mathbb{Z}$). So
$$
f^+_\gamma ~\in~ {X'}_{\hspace{-1,5mm}v_t}^{-{\rm a}_{v_t}(\redfix{\gamma})} \, ({\textstyle \prod}_{v\in \mathcal{V}(Q')\setminus\{v_t\}} {X'}_{\hspace{-1,5mm}v}^{{\rm a}_v(\redfix{\gamma})}) \cdot \mathbb{Z}[\{{X'}_{\hspace{-1,5mm}v}^{\pm 1} \, | \, v\in \mathcal{V}(Q')\}],
$$
\redfix{proving the desired statement for $\mathbb{I}^+_{{\rm PGL}_3}(\ell)$ in case when $\ell$ is represented by a non-peripheral non-contractible simple loop $\gamma$.}

\vs

In fact, when we apply the above investigation of monodromy matrices to a peripheral loop $\gamma$ around a puncture $p$, by looking at the diagonal entries, we obtain the following: if $\ell$ is a single-component ${\rm SL}_3$-lamination consisting only of this peripheral loop with an integer weight, while we already know from eq.\eqref{eq:I_plus_peripheral} that
$$
\mathbb{I}^+_{{\rm PGL}_3}(\ell) = X_{v_t}^{{\rm a}_{v_t}(\ell)} \, ({\textstyle \prod}_{v\in \mathcal{V}(Q_\Delta)\setminus\{v_t\}} X_v^{{\rm a}_v(\ell)}),
$$
we now know
$$
\mathbb{I}^+_{{\rm PGL}_3}(\ell) 
\redfix{~\in~} {X'}_{\hspace{-1,5mm}v_t}^{-{\rm a}_{v_t}(\ell)} \, ({\textstyle \prod}_{v\in \mathcal{V}(Q')\setminus\{v_t\}} {X'}_{\hspace{-1,5mm}v}^{{\rm a}_v(\ell)}) \redfix{\cdot \mathbb{Z}[\{{X'}_{\hspace{-1,5mm}v}^{\pm 1} \, | \, v\in \mathcal{V}(Q')\}]}.
$$

\vs

\redfix{When $\ell$ can be represented as an ${\rm SL}_3$-web without $3$-valent vertices, it can be represented by disjoint union of simple loops. From the above obtained results for single-loop ${\rm SL}_3$-laminations and by the additivity of tropical coordinates (Lem.\ref{lem:additivity_of_coordinates}), the desired result for $\mathbb{I}^+_{{\rm PGL}_3}(\ell)$ follows. \qed}

\vs

\redfix{Saying again, the statement for a general ${\rm SL}_3$-lamination will be dealt with later in \S\ref{subsec:on_the_effect_of_a_single_mutation}. One observation is that one can rewrite the statement of Prop.\ref{prop:mutation_of_basic_semi-regular_function_at_interior_node_of_triangle} in a slightly different way, as
\begin{align}
\label{eq:single_mutation_rewritten}
\mathbb{I}^+_{{\rm PGL}_3}(\ell) ~\in~ ({\textstyle \prod}_{v\in \mathcal{V}(Q')} {X'}_{\hspace{-1,5mm}v}^{\, {\rm a}'_v(\ell)}) \cdot \mathbb{Z}[\{{X'}_{\hspace{-1,5mm}v}^{\pm 1} \, | \, v\in \mathcal{V}(Q')\}],
\end{align}
where
$$
{\rm a}'_v(\ell) = \left\{
\begin{array}{ll}
- {\rm a}_{v_t}(\ell) + \max( \sum_{v\in \mathcal{V}(Q_\Delta)} [\varepsilon_{v,v_t}]_+ {\rm a}_v(\ell), \, \sum_{v\in \mathcal{V}(Q_\Delta)} [-\varepsilon_{v,v_t}]_+ {\rm a}_v(\ell)) & \mbox{if $v=v_t$}, \\
{\rm a}_v(\ell) & \mbox{if $v\neq v_t$},
\end{array}
\right.
$$
which is the tropical version of the cluster $\mathscr{A}$-mutation at $v_t$. This new choice of exponents of $X'_v$'s may be conceptually preferred, but we chose to work with our version for computational convenience, and because our choice of exponents are manifestly additive with respect to disjoint unions of ${\rm SL}_3$-laminations. To justify that eq.\eqref{eq:single_mutation_rewritten} is equivalent to our version of Prop.\ref{prop:mutation_of_basic_semi-regular_function_at_interior_node_of_triangle}, it suffices to show that ${\rm a}'_{v_t}(\ell) \equiv -{\rm a}_{v_t}(\ell)$ modulo $\mathbb{Z}$. This holds because $\sum_{v\in \mathcal{V}(Q_\Delta)} [\varepsilon_{v,v_t}]_+ {\rm a}_v(\ell) = \sum_{\alpha=1}^3 {\rm a}_{v_{e_\alpha,2}}(\ell)$ and $\sum_{v\in \mathcal{V}(Q_\Delta)} [-\varepsilon_{v,v_t}]_+ {\rm a}_v(\ell)) = \sum_{\alpha=1}^3 {\rm a}_{v_{e_\alpha,1}}(\ell)$ (in the notation of Fig.\ref{fig:mutate_two_triangles}) are both integers, due to the balancedness result in Prop.\ref{prop:tropical_coordinate_is_well-defined}(BE1).
}

\vs

\redfix{Now we turn to the effect of a single mutation at a node lying on an arc of $\Delta$.}

\begin{proposition}[mutations of a basic semi-regular function at edge nodes of a triangle]
\label{prop:mutation_of_basic_semi-regular_function_at_edge_node_of_triangle}
Let $\Delta$ be any ideal triangulation of a punctured surface $\frak{S}$. Consider the cluster $\mathscr{X}$-chart associated to $\Delta$, and mutate it at a node $v_0$ of $Q_\Delta$ lying in an edge of $\Delta$. Denote the resulting quiver by $Q''$, naturally identifying $\mathcal{V}(Q_\Delta)$ and $\mathcal{V}(Q'')$. Denote by $X_v''$ the $\mathscr{X}$-coordinate for the node $v$ of $Q''$ for this new chart after mutation. Then for any $\ell \in \mathscr{A}_{\rm L}(\frak{S},\mathbb{Z})$ we have
$$
\mathbb{I}^+_{{\rm PGL}_3}(\ell) ~\in~ {X''}_{\hspace{-2mm}v_0}^{-{\rm a}_{v_0}(\ell) + \sum_{v\in \mathcal{V}(Q_\Delta)} [\varepsilon_{v_0, v}]_+ {\rm a}_v(\ell)} \, (\underset{v\in \mathcal{V}(Q'')\setminus\{v_0\}}{\textstyle \prod} {X''}_{\hspace{-2mm}v}^{{\rm a}_v(\ell)}) \cdot \mathbb{Z}[\{{X''}_{\hspace{-2mm}v}^{\pm 1} \, | \, v\in \mathcal{V}(Q'')\}].
$$
\end{proposition}

\redfix{\it Partial proof.} \redfix{Like in Prop.\ref{prop:mutation_of_basic_semi-regular_function_at_interior_node_of_triangle}, here we prove the statement only for the case when $\ell$ can be represented by an ${\rm SL}_3$-web without $3$-valent vertices, while the general cases will be treated in the next section (\S\ref{subsec:on_the_effect_of_a_single_mutation}).} We use same notations as in the proof of Prop.\ref{prop:mutation_of_basic_semi-regular_function_at_interior_node_of_triangle} for $t,e_1,e_2,e_3$ and for nodes $v_t,v_{e_\alpha,1},v_{e_\alpha,2}$ $(\alpha=1,2,3)$ of $Q_\Delta$ appearing in $t$. It suffices to investigate the mutation at the node $v_{e_1,1}$; we do not lose generality. Let $r$ be the other triangle of $\Delta$ sharing $e_1$ as a side. Label the sides of $r$ as $e_4,e_5,e_6$ clockwise in $\partial r$, so that $e_4$ coincides with $e_1$. For each $\beta=4,5,6$, let $v_{e_\beta,1}$ and $v_{e_\beta,2}$ be the nodes of $Q_\Delta$ lying in $e_\beta$. In particular, we have $v_{e_1,1} = v_{e_4,2}$, $v_{e_1,2} = v_{e_4,1}$; see Fig.\ref{fig:mutate_two_triangles}. First, we assume that none of the external edges $e_2$, $e_3$, $e_5$, $e_6$ of the quadrilateral formed by $t$, $r$ are identified with each other. Under the mutation at the node $v_{e_1,1}$, the cluster $\mathscr{X}$-variables change as (eq.\eqref{eq:X-mutation_formula})
\begin{align*}
& X_{v_{e_1,1}}'' = X_{v_{e_1,1}}^{-1}, \quad
X_{v_{e_3,2}}'' = X_{v_{e_3,2}}(1+X_{v_{e_1,1}}), \quad
X_{v_r}'' = X_{v_r}(1+X_{v_{e_1,1}}), \\
& X_{v_{e_5,1}}'' = X_{v_{e_5,1}}(1+X_{v_{e_1,1}}^{-1})^{-1}, \qquad
X_{v_t}'' = X_{v_t}(1+X_{v_{e_1,1}}^{-1})^{-1},
\end{align*}
and $X''_v = X_v$ for all other nodes $v$ of $Q_\Delta$. Writing the old variables as new ones,
\begin{align}
\label{eq:old_variable_as_new_new_variables}
\left\{ \begin{array}{l}
X_{v_{e_1,1}} = {X''}_{\hspace{-2mm}v_{e_1,1}}^{-1}, \quad
X_{v_{e_3,2}} = {X''}_{\hspace{-2mm}v_{e_3,2}}(1+{X''}_{\hspace{-2mm}v_{e_1,1}}^{-1})^{-1} = {X''}_{\hspace{-2mm}v_{e_3,2}}{X''}_{\hspace{-2mm}v_{e_1,1}} (1+{X''}_{\hspace{-2mm}v_{e_1,1}})^{-1}, \\
X_{v_r} = X''_{v_r}(1+{X''}_{\hspace{-2mm}v_{e_1,1}}^{-1})^{-1} = X''_{v_r} {X''}_{\hspace{-2mm}v_{e_1,1}} (1+{X''}_{\hspace{-2mm}v_{e_1,1}})^{-1} , \\
X_{v_{e_5,1}} = X''_{v_{e_5,1}}(1+{X''}_{\hspace{-2mm}v_{e_1,1}}), \qquad
X_{v_t} = X_{v_t}'' (1+{X''}_{\hspace{-2mm}v_{e_1,1}}), \\
X_v = X_v'' \quad \mbox{for all other nodes $v$}.
\end{array} \right.
\end{align}
We then proceed as in the proof of Prop.\ref{prop:mutation_of_basic_semi-regular_function_at_interior_node_of_triangle}, to study the monodromy matrices of triangle and juncture segments. Triangle segments in $t$ are denoted by $\gamma_{\alpha_1\alpha_2}$ and juncture segments in $t$ by $\gamma_{\alpha,{\rm in}}$ and $\gamma_{\alpha,{\rm out}}$, as in the proof of Prop.\ref{prop:mutation_of_basic_semi-regular_function_at_interior_node_of_triangle}. Define triangle segments $\gamma_{\beta_1\beta_2}$ and juncture segments $\gamma_{\beta,{\rm in}}$ and $\gamma_{\beta,{\rm out}}$ for triangle $r$ analogously. In particular, $\gamma_{1,{\rm out}} = \gamma_{4,{\rm in}}$ and $\gamma_{1,{\rm in}} = \gamma_{4,{\rm out}}$ under this notation; see Fig.\ref{fig:mutate_two_triangles_segments}.

\vs

Using a similar argument used near the end of proof of Prop.\ref{prop:mutation_of_basic_semi-regular_function_at_interior_node_of_triangle}, it suffices to verify the sought-for statement only in the case when $\ell$ is a single oriented \redfix{non-contractible} simple loop $\gamma$, where we replace $\mathbb{I}^+_{{\rm PGL}_3}(\ell)$ in the statement by the trace-of-monodromy $f^+_\gamma$. We express $\gamma$ as concatenation $\gamma = \gamma_1.\gamma_2.\cdots.\gamma_N$ of triangle segments and juncture segments, and make use of $f^+_\gamma = {\rm tr}({\bf M}_{\gamma_1} \cdots {\bf M}_{\gamma_N})$. This time we are considering two adjacent triangles forming a quadrilateral. So we should consider all possibilities of concatenations of segments in Fig.\ref{fig:mutate_two_triangles_segments} forming a `complete' concatenation in this quadrilateral:
\begin{align*}
& {\gamma_{\alpha_1,{\rm in}}} . {\gamma_{\alpha_1 \alpha_2}} . {\gamma_{\alpha_2,{\rm out}}}, \quad
{\gamma_{\beta_1,{\rm in}}} . {\gamma_{\beta_1 \beta_2}} . {\gamma_{\beta_2,{\rm out}}}
\qquad \{\alpha_1,\alpha_2\} = \{2,3\}, ~\{\beta_1,\beta_2\} = \{5,6\}, \\
& {\gamma_{\alpha,{\rm in}}} . {\gamma_{\alpha 1}} . {\gamma_{1,{\rm out}}} . {\gamma_{4\beta}} . {\gamma_{\beta,{\rm out}}}, \qquad \alpha \in \{2,3\}, \quad \beta \in \{5,6\}, \\
& {\gamma_{\beta,{\rm in}}} . {\gamma_{\beta 4}} . {\gamma_{1,{\rm in}}} . {\gamma_{1 \alpha}} . {\gamma_{\alpha,{\rm out}}}, \qquad \alpha \in \{2,3\}, \quad \beta \in \{5,6\}.
\end{align*}
For each case of a complete concatenation, we should compute the effect of mutation on the product of corresponding monodromy matrices ${\bf M}_\cdot$ defined in (MM1)--(MM3) of \S\ref{subsec:lifting_PGL3_to_SL3}. We use the normalized matrices $\til{\bf M}_\cdot$ as in the proof of Prop.\ref{prop:mutation_of_basic_semi-regular_function_at_interior_node_of_triangle}. For a juncture segment $\gamma_i$ as in Fig.\ref{fig:juncture_segment}, let $\til{\bf M}_{\gamma_i}$ be as in eq.\eqref{eq:til_juncture_segment_matrix}. For triangle $t$, the left and the right turn matrices $\til{\bf M}^{\rm left}_t$ and $\til{\bf M}^{\rm right}_t$ are as in eq.\eqref{eq:til_triangle_segment_matrix_t}; define the corresponding matrices $\til{\bf M}^{\rm left}_r$ and $\til{\bf M}^{\rm right}_r$ for triangle $r$ by replacing each $X_{v_t}^{-1}$ by $X_{v_r}^{-1}$. So these $\til{\bf M}_\cdot$ matrices are obtained by dividing by the $(1,1)$-th entry of the corresponding original matrix ${\bf M}_\cdot$.

\vs

Let $W_{\alpha_1\alpha_2; t}$ be an ${\rm SL}_3$-web in triangle $t$ consisting just of one corner arc, from edge $e_{\alpha_1}$ to $e_{\alpha_2}$. In terms of the segments of $\gamma$, this can be viewed as concatenation of (part of) $\gamma_{\alpha_1,{\rm in}}$, then $\gamma_{\alpha_1 \alpha_2}$, then (part of) $\gamma_{\alpha_2,{\rm out}}$. We saw in the proof of Prop.\ref{prop:mutation_of_basic_semi-regular_function_at_interior_node_of_triangle} that the corresponding product of original monodromy matrices ${\bf M}_{W_{\alpha_1\alpha_2;t}} := {\bf M}_{\gamma_{\alpha_1,{\rm in}}} {\bf M}_{\gamma_{\alpha_1 \alpha_2}} {\bf M}_{\gamma_{\alpha_2,{\rm out}}}$ equals the product of normalized matrices $\til{\bf M}_{W_{\alpha_1\alpha_2;t}} := \til{\bf M}_{\gamma_{\alpha_1,{\rm in}}} \til{\bf M}_{\gamma_{\alpha_1 \alpha_2}} \til{\bf M}_{\gamma_{\alpha_2,{\rm out}}}$ times $\prod_{v\in \mathcal{V}(Q_\Delta)\cap t} X_v^{{\rm a}_v(W_{\alpha\beta;t})}$. Likewise for a corner arc ${\rm SL}_3$-web $W_{\beta_1\beta_2;r}$ in triangle $r$. Also, for an ${\rm SL}_3$-web in the union $t\cup r$ of two triangles given as union of a corner arc in $t$ and a corner arc in $r$, similar statement holds. Such web can be either in the form $W_{\alpha\beta;tr} := W_{\alpha 1;t} \cup W_{4\beta;r}$ or $W_{\beta\alpha;rt} := W_{\beta 4;r} \redfix{\cup} W_{1 \alpha;t}$. In the former case $W_{\alpha\beta;tr}$, the corresponding product of monodromy matrices is ${\bf M}_{W_{\alpha\beta;tr}} := {\bf M}_{\gamma_{\alpha,{\rm in}}} {\bf M}_{\gamma_{\alpha 1}} {\bf M}_{\gamma_{1,{\rm out}}} {\bf M}_{\gamma_{4\beta}} {\bf M}_{\gamma_{\beta,{\rm out}}}$ (note ${\bf M}_{\gamma_{1,{\rm out}}} = {\bf M}_{\gamma_{4,{\rm in}}}$), and its $(1,1)$-th entry can be seen to be $\prod_{v\redfix{\in \mathcal{V}(Q_\Delta)\cap(t\cup r)}} X_v^{{\rm a}_v(W_{\alpha\beta;tr})}$. In the latter case $W_{\beta\alpha;rt}$, the corresponding product of matrices is ${\bf M}_{W_{\beta\alpha;rt}} := {\bf M}_{\gamma_{\beta,{\rm in}}} {\bf M}_{\gamma_{\beta 4}} {\bf M}_{\gamma_{1,{\rm in}}} {\bf M}_{\gamma_{1\alpha}} {\bf M}_{\gamma_{\alpha,{\rm out}}}$ (note ${\bf M}_{\gamma_{1,{\rm in}}} = {\bf M}_{\gamma_{4,{\rm out}}}$), and its $(1,1)$-th entry can be seen to be $\prod_{v\in \mathcal{V}(Q_\Delta)\cap(t\cup r)} X_v^{{\rm a}_v(W_{\beta\alpha;rt})}$. Now, we will investigate the effect of mutation on the $(1,1)$-entry-factor and on the (products of) normalized matrices $\til{\bf M}_{W_{\alpha\beta;tr}}$ and $\til{\bf M}_{W_{\beta\alpha;rt}}$. Note that, using this language of ${\rm SL}_3$-webs in $t$ and $r$, the cases to be checked are
\begin{align}
\label{eq:W_cases_to_be_checked}
\left\{ \begin{array}{l}
W_{\alpha_1 \alpha_2;t}, \quad W_{\beta_1\beta_2;r} \qquad \mbox{for $\{\alpha_1,\alpha_2\} = \{2,3\}$},
 \quad \mbox{$\{\beta_1,\beta_2\}=\{5,6\}$}, \\
W_{\alpha\beta;tr} \qquad\qquad\qquad \mbox{for $\alpha \in \{2,3\}$, $\beta\in\{5,6\}$}, \\
W_{\beta\alpha;rt} \qquad \qquad\qquad \mbox{for $\alpha \in \{2,3\}$, $\beta\in\{5,6\}$}.
\end{array} \right.
\end{align}

\vs

For convenience when studying the effect of mutation, we let
$$
{\bf X}_1 := 1+{X''}_{\hspace{-2mm}v_{e_1,1}}.
$$
By eq.\eqref{eq:old_variable_as_new_new_variables}, the effect of mutation on a monomial $\prod_{v\in\mathcal{V}(Q_\Delta)} X_v^{k_v/3}$, for $(k_v)_v \in \mathbb{Z}^{\mathcal{V}(Q_\Delta)}$, is
\begin{align*}
\underset{v\in\mathcal{V}(Q_\Delta)}{\textstyle \prod} \, X_v^{k_v/3} = 
{X''}_{\hspace{-2mm}v_{e_1,1}}^{(-k_{v_{e_1,1}} + k_{v_{e_3,2}} + k_{v_r})/3} {\bf X}_1^{(-k_{v_{e_3,2}} - k_{v_r} + k_{v_{e_5,1}} + k_{v_t})/3} \underset{v\in\mathcal{V}(Q'')\setminus\{v_{e_1,1}\}}{\textstyle \prod} {X''}_{\hspace{-2mm}v}^{k_v/3}
\end{align*}
For all the cases of ${\rm SL}_3$-webs $W$ in eq.\eqref{eq:W_cases_to_be_checked} to be checked, we let $k_v = 3{\rm a}_v(W)$ for nodes $v$ of $Q_\Delta$ living in triangles containing part of $W$ and let $k_v=0$ for other $v\in \mathcal{V}(Q_\Delta)$, and let $k(W)$ be $k_{v_{e_3,2}} - k_{v_r} + k_{v_{e_5,1}} + k_{v_t}$ for this $W$, i.e.
$$
k(W):= 3( -{\rm a}_{v_{e_3,2}}(W) - {\rm a}_{v_r}(W) + {\rm a}_{v_{e_5,1}}(W) + {\rm a}_{v_t}(W) ).
$$
To compute this $k(W)$ for each of our ${\rm SL}_3$-web $W$, we first compute the following numbers for ${\rm SL}_3$-webs $W_{\alpha_1 \alpha_2;t}$ and $W_{\beta_1\beta_2;r}$ living in one triangle:
\begin{align*}
k_t(W_{\alpha_1\alpha_2;t}) := 3(- {\rm a}_{v_{e_3,2}}(W_{\alpha_1\alpha_2;t}) + {\rm a}_{v_t}(W_{\alpha_1\alpha_2;t})), \quad
k_r(W_{\beta_1\beta_2;r}) := 3(-{\rm a}_{v_r}(W_{\beta_1\beta_2;r}) + {\rm a}_{v_{e_5,1}}(W_{\beta_1\beta_2;r})).
\end{align*}
We list the results, which is easily verified from eq.\eqref{eq:tropical_coordinates_for_a_left_turn_W} and eq.\eqref{eq:tropical_coordinates_for_a_right_turn_W}:
\begin{align*}
& k_t(W_{12;t}) = 2, ~
k_t(W_{21;t}) = 1, ~
k_t(W_{13;t}) = -1, ~
k_t(W_{31;t}) = 1,~
k_t(W_{23;t}) = 0, ~
k_t(W_{32;t}) = 0, \\
& k_r(W_{45;r}) = -1, ~
k_r(W_{54;r}) = 1, ~
k_r(W_{46;r}) = -1,~
k_r(W_{64;r}) = -2, ~
k_r(W_{56;r}) = 0, ~
k_r(W_{65;r})=0.
\end{align*}
Then, for $\alpha,\alpha_1,\alpha_2 \in\{2,3\}$ and $\beta,\beta_1,\beta_2 \in\{5,6\}$ one can compute the $k(W)$ using
\begin{align*}
& k(W_{\alpha_1,\alpha_2;t}) = k_t(W_{\alpha_1,\alpha_2;t}), \qquad\qquad
k(W_{\beta_1,\beta_2;r}) = k_r(W_{\beta_1,\beta_2;r}), \\
& k(W_{\alpha\beta;tr}) = k_t(W_{\alpha 1;t}) + k_r(W_{4 \beta;r}), \qquad
k(W_{\beta\alpha;rt}) = k_t(W_{1\alpha;t}) + k_r(W_{\beta 4;r}).
\end{align*}

\vs

Now, we study the effect of mutation on the normalized monodromy matrices $\til{\bf M}_\cdot$ using eq.\eqref{eq:old_variable_as_new_new_variables}. The edge matrices for edges $e_2$ and $e_6$ are easy.
\begin{align*}
\hspace*{-3mm} \til{\bf M}_{\gamma_{2,{\rm in}}} = {\rm diag}(1,X_{v_{e_2,1}}^{-1}, X_{v_{e_2,2}}^{-1} X_{v_{e_2,1}}^{-1})
= {\rm diag}(1,{X''}_{\hspace{-2mm}v_{e_2,1}}^{-1}, {X''}_{\hspace{-2mm}v_{e_2,2}}^{-1} {X''}_{\hspace{-2mm}v_{e_2,1}}^{-1}), \quad
\til{\bf M}_{\gamma_{2,{\rm out}}} 
= {\rm diag}(1,{X''}_{\hspace{-2mm}v_{e_2,2}}^{-1}, {X''}_{\hspace{-2mm}v_{e_2,1}}^{-1} {X''}_{\hspace{-2mm}v_{e_2,2}}^{-1}), \\
\til{\bf M}_{\gamma_{6,{\rm in}}} 
= {\rm diag}(1,{X''}_{\hspace{-2mm}v_{e_6,1}}^{-1}, {X''}_{\hspace{-2mm}v_{e_6,2}}^{-1} {X''}_{\hspace{-2mm}v_{e_6,1}}^{-1}), \quad
\til{\bf M}_{\gamma_{6,{\rm out}}} 
= {\rm diag}(1,{X''}_{\hspace{-2mm}v_{e_6,2}}^{-1}, {X''}_{\hspace{-2mm}v_{e_6,1}}^{-1} {X''}_{\hspace{-2mm}v_{e_6,2}}^{-1}).
\end{align*}
The remaining edge matrices are
\begin{align*}
\til{\bf M}_{\gamma_{3,{\rm in}}} & = {\rm diag}(1,X_{v_{e_3,1}}^{-1}, X_{v_{e_3,2}}^{-1} X_{v_{e_3,1}}^{-1})
= {\rm diag}(1, {X''}_{\hspace{-2mm}v_{e_3,1}}^{-1}, 
{X''}_{\hspace{-2mm}v_{e_1,1}}^{-1} {\bf X}_1
{X''}_{\hspace{-2mm}v_{e_3,2}}^{-1}
{X''}_{\hspace{-2mm}v_{e_3,1}}^{-1}) \\
\til{\bf M}_{\gamma_{3,{\rm out}}} 
& = {\rm diag}(1, {X''}_{\hspace{-2mm}v_{e_1,1}}^{-1} {\bf X}_1
{X''}_{\hspace{-2mm}v_{e_3,2}}^{-1}, 
{X''}_{\hspace{-2mm}v_{e_1,1}}^{-1} {\bf X}_1 
{X''}_{\hspace{-2mm}v_{e_3,2}}^{-1}
{X''}_{\hspace{-2mm}v_{e_3,1}}^{-1}) \\
\til{\bf M}_{\gamma_{5,{\rm in}}} & = {\rm diag}(1,X_{v_{e_5,1}}^{-1}, X_{v_{e_5,2}}^{-1} X_{v_{e_5,1}}^{-1})
= {\rm diag}(1, {\bf X}_1^{-1} {X''}_{\hspace{-2mm}v_{e_5,1}}^{-1}, 
{\bf X}_1^{-1}
{X''}_{\hspace{-2mm}v_{e_5,2}}^{-1}
{X''}_{\hspace{-2mm}v_{e_5,1}}^{-1}) \\
\til{\bf M}_{\gamma_{5,{\rm out}}} & 
= {\rm diag}(1, {X''}_{\hspace{-2mm}v_{e_5,2}}^{-1}, 
{\bf X}_1^{-1}
{X''}_{\hspace{-2mm}v_{e_5,2}}^{-1}
{X''}_{\hspace{-2mm}v_{e_5,1}}^{-1}) \\
\til{\bf M}_{\gamma_{1,{\rm in}}} & = \til{\bf M}_{\gamma_{4,{\rm out}}} = {\rm diag}(1,X_{v_{e_1,1}}^{-1}, X_{v_{e_1,2}}^{-1} X_{v_{e_1,1}}^{-1})
= {\rm diag}(1, {X''}_{\hspace{-2mm}v_{e_1,1}}, 
{X''}_{\hspace{-2mm}v_{e_1,2}}^{-1}
{X''}_{\hspace{-2mm}v_{e_1,1}}
) \\
\til{\bf M}_{\gamma_{1,{\rm out}}} & = \til{\bf M}_{\gamma_{4,{\rm in}}}  
= {\rm diag}(1, {X''}_{\hspace{-2mm}v_{e_1,2}}^{-1}, 
{X''}_{\hspace{-2mm}v_{e_1,2}}^{-1}
{X''}_{\hspace{-2mm}v_{e_1,1}}
) 
\end{align*}
The left and the right turn matrices are:
\begin{align*}
&  \til{\bf M}^{\rm left}_t = 
\smallmatthree{1}{~1+{\bf X}_1^{-1} {X''}_{\hspace{-2mm}v_t}^{-1}~}{{\bf X}_1^{-1} {X''}_{\hspace{-2mm}v_t}^{-1}}{0}{{\bf X}_1^{-1} {X''}_{\hspace{-2mm}v_t}^{-1}}{{\bf X}_1^{-1} {X''}_{\hspace{-2mm}v_t}^{-1}}{0}{0}{{\bf X}_1^{-1} {X''}_{\hspace{-2mm}v_t}^{-1}},
\qquad\qquad
\til{\bf M}^{\rm right}_t = 
\smallmatthree{1}{0}{0}{1}{1}{0}{1}{~1+{\bf X}_1^{-1} {X''}_{\hspace{-2mm}v_t}^{-1}~}{{\bf X}_1^{-1} {X''}_{\hspace{-2mm}v_t}^{-1}}, \\
& \til{\bf M}^{\rm left}_r = 
\smallmatthree{1}{~1+{X''}_{\hspace{-2mm}v_{e_1,1}}^{-1} {\bf X}_1 {X''}_{\hspace{-2mm}v_r}^{-1}~}{{X''}_{\hspace{-2mm}v_{e_1,1}}^{-1} {\bf X}_1 {X''}_{\hspace{-2mm}v_r}^{-1}}{0}{{X''}_{\hspace{-2mm}v_{e_1,1}}^{-1} {\bf X}_1 {X''}_{\hspace{-2mm}v_r}^{-1}}{{X''}_{\hspace{-2mm}v_{e_1,1}}^{-1} {\bf X}_1 {X''}_{\hspace{-2mm}v_r}^{-1}}{0}{0}{{X''}_{\hspace{-2mm}v_{e_1,1}}^{-1} {\bf X}_1 {X''}_{\hspace{-2mm}v_r}^{-1}}, \quad
\til{\bf M}^{\rm right}_r = 
\smallmatthree{1}{0}{0}{1}{1}{0}{1}{~1+{X''}_{\hspace{-2mm}v_{e_1,1}}^{-1} {\bf X}_1 {X''}_{\hspace{-2mm}v_r}^{-1}~}{{X''}_{\hspace{-2mm}v_{e_1,1}}^{-1} {\bf X}_1 {X''}_{\hspace{-2mm}v_r}^{-1}}.
\end{align*}

\vs

We should check, for each ${\rm SL}_3$-web $W$ in eq.\eqref{eq:W_cases_to_be_checked}, that ${\bf X}_1^{k(W)/3}$ times the corresponding product of normalized monodromy matrices lives in ${\rm GL}_3(\mathbb{Z}[\{{X''}_{\hspace{-2mm}v}^{\pm 1}\,|\,v\in\mathcal{V}(Q'')\}])$, i.e. the entries are \ul{\em $X''$-Laurent}, i.e. Laurent polynomials in $\{{X''}_{\hspace{-2mm}v}\,|\,v\in\mathcal{V}(Q'')\}$ with integer coefficients. There are 12 cases to check in total. The point is to check that in the entries of the final matrices, we see no negative powers of ${\bf X}_1$. Note that the only basic normalized monodromy matrices that are not $X''$-Laurent are $\til{\bf M}_{\gamma_{5,*}}$, $\til{\bf M}^{\rm left}_t$, $\til{\bf M}^{\rm right}_t$ because they involve ${\bf X}_1^{-1}$, so we should keep an eye on them; on the other hand, keep track of $\til{\bf M}_{\gamma_{3,*}}$, $\til{\bf M}^{\rm left}_t$, $\til{\bf M}^{\rm right}_r$ as they involve ${\bf X}_1$. We cannot assume much symmetry, so we deal with all 12 cases explicitly. Case $1$ is $W_{23;t}$, where we have $k(W_{23;t})=0$, and the product of normalized matrices is $\til{\bf M}_{W_{23;t}}
= \til{\bf M}_{\gamma_{2,{\rm in}}} \til{\bf M}^{\rm left}_t \til{\bf M}_{\gamma_{3,{\rm out}}}$, which is manifestly $X''$-Laurent because each factor already is. In Case 2 we have $k(W_{32;t})=0$, and
\begin{align*}
\hspace*{-3mm}
\til{\bf M}_{W_{32;t}}
= \til{\bf M}_{\gamma_{3,{\rm in}}} \til{\bf M}^{\rm right}_t \til{\bf M}_{\gamma_{2,{\rm out}}}
= \smallmatthree{1}{0}{0}{0}{{X''}_{\hspace{-2mm}v_{e_3,1}}^{-1}}{0}{0}{0}{{X''}_{\hspace{-2mm}v_{e_1,1}}^{-1} {\bf X}_1
{X''}_{\hspace{-2mm}v_{e_3,2}}^{-1}
{X''}_{\hspace{-2mm}v_{e_3,1}}^{-1}}
\smallmatthree{1}{0}{0}{1}{1}{0}{1}{~1+{\bf X}_1^{-1} {X''}_{\hspace{-2mm}v_t}^{-1}~}{{\bf X}_1^{-1} {X''}_{\hspace{-2mm}v_t}^{-1}} \til{\bf M}_{\gamma_{2;{\rm out}}};
\end{align*}
when the first two matrices are multiplied, ${\bf X}_1^{-1}$ are cancelled, so the resulting entries are $X''$-Laurent. In Case 3 we have $k(W_{56;r})=0$, and
\begin{align*}
\hspace*{-3mm}
\til{\bf M}_{W_{56};r} = \til{\bf M}_{\gamma_{5,{\rm in}}} \til{\bf M}^{\rm left}_r \til{\bf M}_{\gamma_{6,{\rm out}}}
= 
\smallmatthree{1}{0}{0}{0}{{\bf X}_1^{-1} {X''}_{\hspace{-2mm}v_{e_5,1}}^{-1}}{0}{0}{0}{{\bf X}_1^{-1}
{X''}_{\hspace{-2mm}v_{e_5,2}}^{-1}
{X''}_{\hspace{-2mm}v_{e_5,1}}^{-1}}
\smallmatthree{1}{~1+{X''}_{\hspace{-2mm}v_{e_1,1}}^{-1} {\bf X}_1 {X''}_{\hspace{-2mm}v_r}^{-1}~}{{X''}_{\hspace{-2mm}v_{e_1,1}}^{-1} {\bf X}_1 {X''}_{\hspace{-2mm}v_r}^{-1}}{0}{{X''}_{\hspace{-2mm}v_{e_1,1}}^{-1} {\bf X}_1 {X''}_{\hspace{-2mm}v_r}^{-1}}{{X''}_{\hspace{-2mm}v_{e_1,1}}^{-1} {\bf X}_1 {X''}_{\hspace{-2mm}v_r}^{-1}}{0}{0}{{X''}_{\hspace{-2mm}v_{e_1,1}}^{-1} {\bf X}_1 {X''}_{\hspace{-2mm}v_r}^{-1}} \til{\bf M}_{\gamma_{6,{\rm out}}};
\end{align*}
when the first two matrices are multiplied, ${\bf X}_1^{-1}$'s are cancelled. In Case 4, $k(W_{65;r})=0$, and
\begin{align*}
\hspace*{-3mm}
\til{\bf M}_{W_{65};r} = \til{\bf M}_{\gamma_{6,{\rm in}}} \til{\bf M}^{\rm right}_r \til{\bf M}_{\gamma_{5,{\rm out}}}
= \til{\bf M}_{\gamma_{6,{\rm in}}}
\smallmatthree{1}{0}{0}{1}{1}{0}{1}{~1+{X''}_{\hspace{-2mm}v_{e_1,1}}^{-1} {\bf X}_1 {X''}_{\hspace{-2mm}v_r}^{-1}~}{{X''}_{\hspace{-2mm}v_{e_1,1}}^{-1} {\bf X}_1 {X''}_{\hspace{-2mm}v_r}^{-1}}
\smallmatthree{1}{0}{0}{0}{{X''}_{\hspace{-2mm}v_{e_5,2}}^{-1}}{0}{0}{0}{{\bf X}_1^{-1}
{X''}_{\hspace{-2mm}v_{e_5,2}}^{-1}
{X''}_{\hspace{-2mm}v_{e_5,1}}^{-1}};
\end{align*}
when the latter two matrices are multiplied, ${\bf X}_1^{-1}$ is cancelled. In Case 5, $k(W_{25;tr}) = 1-1=0$,
\begin{align*}
& \hspace*{-3mm} \til{\bf M}_{W_{25;tr}} = \til{\bf M}_{\gamma_{2,{\rm in}}} \til{\bf M}^{\rm right}_t \til{\bf M}_{\gamma_{1,{\rm out}}} \til{\bf M}^{\rm left}_r \til{\bf M}_{\gamma_{5,{\rm out}}} \\
& \hspace*{-3mm}  = \til{\bf M}_{\gamma_{2,{\rm in}}}
\smallmatthree{1}{0}{0}{1}{1}{0}{1}{~1+{\bf X}_1^{-1} {X''}_{\hspace{-2mm}v_t}^{-1}~}{{\bf X}_1^{-1} {X''}_{\hspace{-2mm}v_t}^{-1}}
\til{\bf M}_{\gamma_{1,{\rm out}}}
\smallmatthree{1}{~1+{X''}_{\hspace{-2mm}v_{e_1,1}}^{-1} {\bf X}_1 {X''}_{\hspace{-2mm}v_r}^{-1}~}{{X''}_{\hspace{-2mm}v_{e_1,1}}^{-1} {\bf X}_1 {X''}_{\hspace{-2mm}v_r}^{-1}}{0}{{X''}_{\hspace{-2mm}v_{e_1,1}}^{-1} {\bf X}_1 {X''}_{\hspace{-2mm}v_r}^{-1}}{{X''}_{\hspace{-2mm}v_{e_1,1}}^{-1} {\bf X}_1 {X''}_{\hspace{-2mm}v_r}^{-1}}{0}{0}{{X''}_{\hspace{-2mm}v_{e_1,1}}^{-1} {\bf X}_1 {X''}_{\hspace{-2mm}v_r}^{-1}}
\smallmatthree{1}{0}{0}{0}{{X''}_{\hspace{-2mm}v_{e_5,2}}^{-1}}{0}{0}{0}{{\bf X}_1^{-1}
{X''}_{\hspace{-2mm}v_{e_5,2}}^{-1}
{X''}_{\hspace{-2mm}v_{e_5,1}}^{-1}} \\
& \hspace*{-3mm} = \til{\bf M}_{\gamma_{2,{\rm in}}}
\smallmatthree{1}{0}{0}{1}{{X''}_{\hspace{-2mm}v_{e_1,2}}^{-1}}{0}{1}{~(1+{\bf X}_1^{-1} {X''}_{\hspace{-2mm}v_t}^{-1}){X''}_{\hspace{-2mm}v_{e_1,2}}^{-1} ~}{{\bf X}_1^{-1} {X''}_{\hspace{-2mm}v_t}^{-1} {X''}_{\hspace{-2mm}v_{e_1,2}}^{-1}
{X''}_{\hspace{-2mm}v_{e_1,1}} }
\smallmatthree{1}{~(1+{X''}_{\hspace{-2mm}v_{e_1,1}}^{-1} {\bf X}_1 {X''}_{\hspace{-2mm}v_r}^{-1}){X''}_{\hspace{-2mm}v_{e_5,2}}^{-1} ~}{{X''}_{\hspace{-2mm}v_{e_1,1}}^{-1}  {X''}_{\hspace{-2mm}v_r}^{-1} {X''}_{\hspace{-2mm}v_{e_5,2}}^{-1}
{X''}_{\hspace{-2mm}v_{e_5,1}}^{-1}}{0}{{X''}_{\hspace{-2mm}v_{e_1,1}}^{-1} {\bf X}_1 {X''}_{\hspace{-2mm}v_r}^{-1} {X''}_{\hspace{-2mm}v_{e_5,2}}^{-1} }{{X''}_{\hspace{-2mm}v_{e_1,1}}^{-1} {X''}_{\hspace{-2mm}v_r}^{-1} {X''}_{\hspace{-2mm}v_{e_5,2}}^{-1}
{X''}_{\hspace{-2mm}v_{e_5,1}}^{-1} }{0}{0}{{X''}_{\hspace{-2mm}v_{e_1,1}}^{-1}  {X''}_{\hspace{-2mm}v_r}^{-1} {X''}_{\hspace{-2mm}v_{e_5,2}}^{-1}
{X''}_{\hspace{-2mm}v_{e_5,1}}^{-1}}
\end{align*}
When taking product of the latter two matrices, the only entry that is not manifestly $X''$-Laurent is the $(3,3)$-th entry, which equals ${X''}_{\hspace{-2mm}v_{e_1,1}}^{-1}  {X''}_{\hspace{-2mm}v_r}^{-1} {X''}_{\hspace{-2mm}v_{e_5,2}}^{-1}
{X''}_{\hspace{-2mm}v_{e_5,1}}^{-1}$ times
\begin{align}
\label{eq:bf_X_1_canceling1}
\hspace*{-1mm} 1 + \ul{ (1+{\bf X}_1^{-1} {X''}_{\hspace{-2mm}v_t}^{-1} ) {X''}_{\hspace{-2mm}v_{e_1,2}}^{-1}  + {\bf X}_1^{-1} {X''}_{\hspace{-2mm}v_t} ^{-1} {X''}_{\hspace{-2mm}v_{e_1,2}}^{-1}  {X''}_{\hspace{-2mm}v_{e_1,1}} }
= 1 + \ul{ {X''}_{\hspace{-2mm}v_{e_1,2}}^{-1} 
+ \cancel{ ( 1+ {X''}_{\hspace{-2mm}v_{e_1,1}}) {\bf X}_1^{-1} } {X''}_{\hspace{-2mm}v_t}^{-1}  {X''}_{\hspace{-2mm}v_{e_1,2}}^{-1} }
\end{align}
hence is $X''$-Laurent. In Case 6, $k(W_{26;tr}) = 1-1=0$, and
\begin{align*}
& \hspace*{-3mm} \til{\bf M}_{W_{26;tr}} = \til{\bf M}_{\gamma_{2,{\rm in}}} \til{\bf M}^{\rm right}_t \til{\bf M}_{\gamma_{1,{\rm out}}} \til{\bf M}^{\rm right}_r \til{\bf M}_{\gamma_{6,{\rm out}}} \\
& \hspace*{-3mm}  = \til{\bf M}_{\gamma_{2,{\rm in}}}
\smallmatthree{1}{0}{0}{1}{{X''}_{\hspace{-2mm}v_{e_1,2}}^{-1}}{0}{1}{~(1+{\bf X}_1^{-1} {X''}_{\hspace{-2mm}v_t}^{-1}){X''}_{\hspace{-2mm}v_{e_1,2}}^{-1} ~}{ {\bf X}_1^{-1} {X''}_{\hspace{-2mm}v_t}^{-1} {X''}_{\hspace{-2mm}v_{e_1,2}}^{-1}
{X''}_{\hspace{-2mm}v_{e_1,1}} }
\smallmatthree{1}{0}{0}{1}{1}{0}{1}{~1+{X''}_{\hspace{-2mm}v_{e_1,1}}^{-1} {\bf X}_1 {X''}_{\hspace{-2mm}v_r}^{-1}~}{{X''}_{\hspace{-2mm}v_{e_1,1}}^{-1} {\bf X}_1 {X''}_{\hspace{-2mm}v_r}^{-1}}
\til{\bf M}_{\gamma_{6,{\rm out}}}.
\end{align*}
When taking the product of the middle two matrices, the only entries that are not manifestly $X''$-Laurent are $(3,1)$-th and $(3,2)$-th entries. The $(3,1)$-th entry is $X''$-Laurent due to the computation in eq.\eqref{eq:bf_X_1_canceling1}. The $(3,2)$-th entry is
\begin{align}
\label{eq:bf_X_1_canceling2}
\begin{array}{l}
(1+{\bf X}_1^{-1} {X''}_{\hspace{-2mm}v_t}^{-1}){X''}_{\hspace{-2mm}v_{e_1,2}}^{-1}
+ {\bf X}_1^{-1} {X''}_{\hspace{-2mm}v_t}^{-1} {X''}_{\hspace{-2mm}v_{e_1,2}}^{-1}
{X''}_{\hspace{-2mm}v_{e_1,1}} (1+{X''}_{\hspace{-2mm}v_{e_1,1}}^{-1} {\bf X}_1 {X''}_{\hspace{-2mm}v_r}^{-1}) \\
= \mbox{( underlined part in eq.\eqref{eq:bf_X_1_canceling1} )}
+ \cancel{ {\bf X}_1^{-1} } {X''}_{\hspace{-2mm}v_t}^{-1} {X''}_{\hspace{-2mm}v_{e_1,2}}^{-1}
{X''}_{\hspace{-2mm}v_{e_1,1}} {X''}_{\hspace{-2mm}v_{e_1,1}}^{-1} \cancel{ {\bf X}_1 } {X''}_{\hspace{-2mm}v_r}^{-1}
\end{array}
\end{align}
hence is $X''$-Laurent. In Case 7, we have $k(W_{35;tr}) = 1-1=0$, and
\begin{align*}
& \hspace*{-3mm} \til{\bf M}_{W_{35;tr}} = \til{\bf M}_{\gamma_{3,{\rm in}}} \til{\bf M}^{\rm left}_t \til{\bf M}_{\gamma_{1,{\rm out}}} \til{\bf M}^{\rm left}_r \til{\bf M}_{\gamma_{5,{\rm out}}} \\
&\hspace*{-3mm} = \til{\bf M}_{\gamma_{3,{\rm in}}} 
\smallmatthree{1}{~1+{\bf X}_1^{-1} {X''}_{\hspace{-2mm}v_t}^{-1}~}{{\bf X}_1^{-1} {X''}_{\hspace{-2mm}v_t}^{-1}}{0}{{\bf X}_1^{-1} {X''}_{\hspace{-2mm}v_t}^{-1}}{{\bf X}_1^{-1} {X''}_{\hspace{-2mm}v_t}^{-1}}{0}{0}{{\bf X}_1^{-1} {X''}_{\hspace{-2mm}v_t}^{-1}}
\til{\bf M}_{\gamma_{1,{\rm out}}} 
\smallmatthree{1}{~(1+{X''}_{\hspace{-2mm}v_{e_1,1}}^{-1} {\bf X}_1 {X''}_{\hspace{-2mm}v_r}^{-1}){X''}_{\hspace{-2mm}v_{e_5,2}}^{-1} ~}{{X''}_{\hspace{-2mm}v_{e_1,1}}^{-1}  {X''}_{\hspace{-2mm}v_r}^{-1} {X''}_{\hspace{-2mm}v_{e_5,2}}^{-1}
{X''}_{\hspace{-2mm}v_{e_5,1}}^{-1}}{0}{{X''}_{\hspace{-2mm}v_{e_1,1}}^{-1} {\bf X}_1 {X''}_{\hspace{-2mm}v_r}^{-1} {X''}_{\hspace{-2mm}v_{e_5,2}}^{-1} }{{X''}_{\hspace{-2mm}v_{e_1,1}}^{-1} {X''}_{\hspace{-2mm}v_r}^{-1} {X''}_{\hspace{-2mm}v_{e_5,2}}^{-1}
{X''}_{\hspace{-2mm}v_{e_5,1}}^{-1} }{0}{0}{{X''}_{\hspace{-2mm}v_{e_1,1}}^{-1}  {X''}_{\hspace{-2mm}v_r}^{-1} {X''}_{\hspace{-2mm}v_{e_5,2}}^{-1}
{X''}_{\hspace{-2mm}v_{e_5,1}}^{-1}} \\
&\hspace*{-3mm}  = \til{\bf M}_{\gamma_{3,{\rm in}}}  \smallmatthree{1}{~(1+{\bf X}_1^{-1} {X''}_{\hspace{-2mm}v_t}^{-1}){X''}_{\hspace{-2mm}v_{e_1,2}}^{-1}~}{ {\bf X}_1^{-1} {X''}_{\hspace{-2mm}v_t}^{-1} {X''}_{\hspace{-2mm}v_{e_1,2}}^{-1}
{X''}_{\hspace{-2mm}v_{e_1,1}} }{0}{{\bf X}_1^{-1} {X''}_{\hspace{-2mm}v_t}^{-1} {X''}_{\hspace{-2mm}v_{e_1,2}}^{-1} }{ {\bf X}_1^{-1} {X''}_{\hspace{-2mm}v_t}^{-1} {X''}_{\hspace{-2mm}v_{e_1,2}}^{-1}
{X''}_{\hspace{-2mm}v_{e_1,1}} }{0}{0}{ {\bf X}_1^{-1} {X''}_{\hspace{-2mm}v_t}^{-1} {X''}_{\hspace{-2mm}v_{e_1,2}}^{-1}
{X''}_{\hspace{-2mm}v_{e_1,1}} }
\smallmatthree{1}{~(1+{X''}_{\hspace{-2mm}v_{e_1,1}}^{-1} {\bf X}_1 {X''}_{\hspace{-2mm}v_r}^{-1}){X''}_{\hspace{-2mm}v_{e_5,2}}^{-1} ~}{{X''}_{\hspace{-2mm}v_{e_1,1}}^{-1}  {X''}_{\hspace{-2mm}v_r}^{-1} {X''}_{\hspace{-2mm}v_{e_5,2}}^{-1}
{X''}_{\hspace{-2mm}v_{e_5,1}}^{-1}}{0}{{X''}_{\hspace{-2mm}v_{e_1,1}}^{-1} {\bf X}_1 {X''}_{\hspace{-2mm}v_r}^{-1} {X''}_{\hspace{-2mm}v_{e_5,2}}^{-1} }{{X''}_{\hspace{-2mm}v_{e_1,1}}^{-1} {X''}_{\hspace{-2mm}v_r}^{-1} {X''}_{\hspace{-2mm}v_{e_5,2}}^{-1}
{X''}_{\hspace{-2mm}v_{e_5,1}}^{-1} }{0}{0}{{X''}_{\hspace{-2mm}v_{e_1,1}}^{-1}  {X''}_{\hspace{-2mm}v_r}^{-1} {X''}_{\hspace{-2mm}v_{e_5,2}}^{-1}
{X''}_{\hspace{-2mm}v_{e_5,1}}^{-1}}
\end{align*}
When taking the product of latter two matrices, the only entries that are not manifestly $X''$-Laurent are $(1,3)$-th, $(2,3)$-th, and $(3,3)$-th entries. The $(1,3)$-th entry is $X''$-Laurent due to eq.\eqref{eq:bf_X_1_canceling1}, and the $(3,3)$-th entry become $X''$-Laurent when we also multiply the matrix $\til{\bf M}_{\gamma_{3,{\rm in}}}$ from left because its $(3,3)$-th entry is divisible by ${\bf X}_1$. The $(2,3)$-th entry is ${X''}_{\hspace{-2mm}v_{e_1,1}}^{-1}  {X''}_{\hspace{-2mm}v_r}^{-1} {X''}_{\hspace{-2mm}v_{e_5,2}}^{-1}
{X''}_{\hspace{-2mm}v_{e_5,1}}^{-1}$ times
\begin{align}
\label{eq:bf_X_1_canceling3}
{\bf X}_1^{-1} {X''}_{\hspace{-2mm}v_t}^{-1} {X''}_{\hspace{-2mm}v_{e_1,2}}^{-1} + {\bf X}_1^{-1} {X''}_{\hspace{-2mm}v_t}^{-1} {X''}_{\hspace{-2mm}v_{e_1,2}}^{-1}
{X''}_{\hspace{-2mm}v_{e_1,1}}
= \cancel{ (1+{X''}_{\hspace{-2mm}v_{e_1,1}}) {\bf X}_1^{-1} } {X''}_{\hspace{-2mm}v_t}^{-1} {X''}_{\hspace{-2mm}v_{e_1,2}}^{-1}
\end{align}
which is $X''$-Laurent. In Case 8, we have $k(W_{36;tr}) = 1-1=0$, and
\begin{align*}
& \hspace*{-4mm} \til{\bf M}_{W_{36;tr}} = \til{\bf M}_{\gamma_{3,{\rm in}}} \til{\bf M}^{\rm left}_t \til{\bf M}_{\gamma_{1,{\rm out}}} \til{\bf M}^{\rm right}_r \til{\bf M}_{\gamma_{6,{\rm out}}} \\
& \hspace*{-4mm} = \til{\bf M}_{\gamma_{3,{\rm in}}}  \smallmatthree{1}{~(1+{\bf X}_1^{-1} {X''}_{\hspace{-2mm}v_t}^{-1}){X''}_{\hspace{-2mm}v_{e_1,2}}^{-1}~}{ {\bf X}_1^{-1} {X''}_{\hspace{-2mm}v_t}^{-1} {X''}_{\hspace{-2mm}v_{e_1,2}}^{-1}
{X''}_{\hspace{-2mm}v_{e_1,1}} }{0}{{\bf X}_1^{-1} {X''}_{\hspace{-2mm}v_t}^{-1} {X''}_{\hspace{-2mm}v_{e_1,2}}^{-1} }{ {\bf X}_1^{-1} {X''}_{\hspace{-2mm}v_t}^{-1} {X''}_{\hspace{-2mm}v_{e_1,2}}^{-1}
{X''}_{\hspace{-2mm}v_{e_1,1}} }{0}{0}{ {\bf X}_1^{-1} {X''}_{\hspace{-2mm}v_t}^{-1} {X''}_{\hspace{-2mm}v_{e_1,2}}^{-1}
{X''}_{\hspace{-2mm}v_{e_1,1}} }
\smallmatthree{1}{0}{0}{1}{1}{0}{1}{~1+{X''}_{\hspace{-2mm}v_{e_1,1}}^{-1} {\bf X}_1 {X''}_{\hspace{-2mm}v_r}^{-1}~}{{X''}_{\hspace{-2mm}v_{e_1,1}}^{-1} {\bf X}_1 {X''}_{\hspace{-2mm}v_r}^{-1}}
\til{\bf M}_{\gamma_{6,{\rm out}}}
\end{align*}
In the product of middle two matrices, we consider the following entries not manifestly $X''$-Laurent. The $(1,1)$-th, the $(1,2)$-th and the $(2,1)$-th entries are $X''$-Laurent due to eq.\eqref{eq:bf_X_1_canceling1}, eq.\eqref{eq:bf_X_1_canceling2}, and eq.\eqref{eq:bf_X_1_canceling3} respectively. So is $(2,2)$-th entry essentially due to eq.\eqref{eq:bf_X_1_canceling3}. The $(3,1$)-th, the $(3,2)$-th and the $(3,3)$-th become $X''$-Laurent when multiplied by $\til{\bf M}_{\gamma_{3,{\rm in}}}$ from the left, cancelling the ${\bf X}_1^{-1}$ factor. In Case 9, we have $k(W_{52;rt}) = 1+2=3$, and
\begin{align*}
& \hspace*{-3mm} \til{\bf M}_{W_{52};rt} = \til{\bf M}_{\gamma_{5,{\rm in}}} \til{\bf M}^{\rm right}_r \til{\bf M}_{\gamma_{1,{\rm in}}} \til{\bf M}^{\rm left}_t \til{\bf M}_{\gamma_{2,{\rm out}}} \\
& \hspace*{-3mm}  = 
\til{\bf M}_{\gamma_{5,{\rm in}}}
\smallmatthree{1}{0}{0}{1}{{X''}_{\hspace{-2mm}v_{e_1,1}}}{0}{1}{~(1+{X''}_{\hspace{-2mm}v_{e_1,1}}^{-1} {\bf X}_1 {X''}_{\hspace{-2mm}v_r}^{-1}){X''}_{\hspace{-2mm}v_{e_1,1}}~}{ \cancel{ {X''}_{\hspace{-2mm}v_{e_1,1}}^{-1} } {\bf X}_1 {X''}_{\hspace{-2mm}v_r}^{-1} {X''}_{\hspace{-2mm}v_{e_1,2}}^{-1}
\cancel{ {X''}_{\hspace{-2mm}v_{e_1,1}}} }
\smallmatthree{1}{~1+{\bf X}_1^{-1} {X''}_{\hspace{-2mm}v_t}^{-1}~}{{\bf X}_1^{-1} {X''}_{\hspace{-2mm}v_t}^{-1}}{0}{{\bf X}_1^{-1} {X''}_{\hspace{-2mm}v_t}^{-1}}{{\bf X}_1^{-1} {X''}_{\hspace{-2mm}v_t}^{-1}}{0}{0}{{\bf X}_1^{-1} {X''}_{\hspace{-2mm}v_t}^{-1}}
\til{\bf M}_{\gamma_{2,{\rm out}}}
\end{align*}
When we take the product of the middle two matrices, the $(3,3)$-th entry is ${\bf X}_1^{-1} {X''}_{\hspace{-2mm}v_t}^{-1}$ times
\begin{align}
\label{eq:bf_X_1_canceling4}
\begin{array}{l}
1
+
(1+{X''}_{\hspace{-2mm}v_{e_1,1}}^{-1} {\bf X}_1 {X''}_{\hspace{-2mm}v_r}^{-1}){X''}_{\hspace{-2mm}v_{e_1,1}}
+
{\bf X}_1 {X''}_{\hspace{-2mm}v_r}^{-1} {X''}_{\hspace{-2mm}v_{e_1,2}}^{-1} \\
= 1 + {X''}_{\hspace{-2mm}v_{e_1,1}} +  {\bf X}_1 {X''}_{\hspace{-2mm}v_r}^{-1}  ( 1 + {X''}_{\hspace{-2mm}v_{e_1,2}}^{-1})
= {\bf X}_1 ( 1 + {X''}_{\hspace{-2mm}v_r}^{-1}  ( 1 + {X''}_{\hspace{-2mm}v_{e_1,2}}^{-1})),
\end{array}
\end{align}
the $(3,2)$-th entry equals $1 + {\bf X}_1^{-1} {X''}_{\hspace{-2mm}v_t}^{-1} + {X''}_{\hspace{-2mm}v_{e_1,1}} {\bf X}_1^{-1} {X''}_{\hspace{-2mm}v_t}^{-1} + {X''}_{\hspace{-2mm}v_r}^{-1} {X''}_{\hspace{-2mm}v_t}^{-1} = 1 + {X''}_{\hspace{-2mm}v_t}^{-1} + {X''}_{\hspace{-2mm}v_r}^{-1} {X''}_{\hspace{-2mm}v_t}^{-1}$, so
\begin{align*}
\til{\bf M}_{W_{52};rt}
= \smallmatthree{1}{0}{0}{0}{{\bf X}_1^{-1} {X''}_{\hspace{-2mm}v_{e_5,1}}^{-1}}{0}{0}{0}{{\bf X}_1^{-1}
{X''}_{\hspace{-2mm}v_{e_5,2}}^{-1}
{X''}_{\hspace{-2mm}v_{e_5,1}}^{-1}}
\smallmatthree{1}{ ~(1+{\bf X}_1^{-1} {X''}_{\hspace{-2mm}v_t}^{-1})~ }{ {\bf X}_1^{-1} {X''}_{\hspace{-2mm}v_{\redfix{t}}}^{-1}  }{ 1 }{ 1 + {X''}_{\hspace{-2mm}v_t}^{-1} }{ {X''}_{\hspace{-2mm}v_t}^{-1} }{ 1 }{~( 1 + {X''}_{\hspace{-2mm}v_t}^{-1} + {X''}_{\hspace{-2mm}v_r}^{-1} {X''}_{\hspace{-2mm}v_t}^{-1} )~ }{ {X''}_{\hspace{-2mm}v_t}^{-1} ( 1 + {X''}_{\hspace{-2mm}v_r}^{-1}  ( 1 + {X''}_{\hspace{-2mm}v_{e_1,2}}^{-1})) }
\til{\bf M}_{\gamma_{2,{\rm out}}}
\end{align*}
If we multiply the first two matrices, in the entries we see some ${\bf X}_1^{-1}$ (but not higher powers of ${\bf X}_1^{-1}$); hence, multiplying $\til{\bf M}_{W_{52};rt}$ by ${\bf X}_1^{k(W_{52;rt})/3} = {\bf X}_1$ yields \redfix{an} $X''$-Laurent matrix. In Case 10, we have $k(W_{53;rt}) = 1-1=0$, and
\begin{align*}
& \hspace*{-3mm} \til{\bf M}_{W_{53};rt} = \til{\bf M}_{\gamma_{5,{\rm in}}} \til{\bf M}^{\rm right}_r \til{\bf M}_{\gamma_{1,{\rm in}}} \til{\bf M}^{\rm right}_t \til{\bf M}_{\gamma_{3,{\rm out}}} \\
& \hspace*{-3mm} = \til{\bf M}_{\gamma_{5,{\rm in}}}
\smallmatthree{1}{0}{0}{1}{{X''}_{\hspace{-2mm}v_{e_1,1}}}{0}{1}{~(1+{X''}_{\hspace{-2mm}v_{e_1,1}}^{-1} {\bf X}_1 {X''}_{\hspace{-2mm}v_r}^{-1}){X''}_{\hspace{-2mm}v_{e_1,1}}~}{ {\bf X}_1 {X''}_{\hspace{-2mm}v_r}^{-1} {X''}_{\hspace{-2mm}v_{e_1,2}}^{-1}}
\smallmatthree{1}{0}{0}{1}{1}{0}{1}{~1+{\bf X}_1^{-1} {X''}_{\hspace{-2mm}v_t}^{-1}~}{{\bf X}_1^{-1} {X''}_{\hspace{-2mm}v_t}^{-1}}
\til{\bf M}_{\gamma_{3,{\rm out}}}
\end{align*}
In the product of the middle two matrices, the (3,1)-th entry is as in eq.\eqref{eq:bf_X_1_canceling4}, the $(3,2)$-th entry is
\begin{align}
\label{eq:bf_X_1_canceling5}
\begin{array}{l}
(1+{X''}_{\hspace{-2mm}v_{e_1,1}}^{-1} {\bf X}_1 {X''}_{\hspace{-2mm}v_r}^{-1}) {X''}_{\hspace{-2mm}v_{e_1,1}} + {\bf X}_1 {X''}_{\hspace{-2mm}v_r}^{-1} {X''}_{\hspace{-2mm}v_{e_1,2}}^{-1}
 (1+{\bf X}_1^{-1} {X''}_{\hspace{-2mm}v_t}^{-1}) \\
= {X''}_{\hspace{-2mm}v_{e_1,1}} + {\bf X}_1 {X''}_{\hspace{-2mm}v_r}^{-1} (1 + {X''}_{\hspace{-2mm}v_{e_1,2}}^{-1})
+ {X''}_{\hspace{-2mm}v_r}^{-1}
{X''}_{\hspace{-2mm}v_{e_1,2}}^{-1}
{X''}_{\hspace{-2mm}v_t}^{-1}.
\end{array}
\end{align}
So $\til{\bf M}_{W_{53};rt}$ equals
\begin{align*}
& \hspace*{-8mm}  
\smallmatthree{1}{0}{0}{0}{{\bf X}_1^{-1} {X''}_{\hspace{-2mm}v_{e_5,1}}^{-1}}{0}{0}{0}{{\bf X}_1^{-1}
{X''}_{\hspace{-2mm}v_{e_5,2}}^{-1}
{X''}_{\hspace{-2mm}v_{e_5,1}}^{-1}} 
\smallmatthree{1}{0}{0}{{\bf X}_1}{{X''}_{\hspace{-2mm}v_{e_1,1}}}{0}{{\bf X}_1 ( 1 + {X''}_{\hspace{-2mm}v_r}^{-1}  ( 1 + {X''}_{\hspace{-2mm}v_{e_1,2}}^{-1}))~}{~\redfix{{X''}_{\hspace{-2mm}v_{e_1,1}} + {\bf X}_1 {X''}_{\hspace{-2mm}v_r}^{-1} (1 + {X''}_{\hspace{-2mm}v_{e_1,2}}^{-1})
+ {X''}_{\hspace{-2mm}v_r}^{-1}
{X''}_{\hspace{-2mm}v_{e_1,2}}^{-1}
{X''}_{\hspace{-2mm}v_t}^{-1}}
~}{ {X''}_{\hspace{-2mm}v_r}^{-1} {X''}_{\hspace{-2mm}v_{e_1,2}}^{-1} {X''}_{\hspace{-2mm}v_t}^{-1} } \\
& \cdot
\smallmatthree{1}{0}{0}{0}{{X''}_{\hspace{-2mm}v_{e_1,1}}^{-1} {\bf X}_1
{X''}_{\hspace{-2mm}v_{e_3,2}}^{-1}}{0}{0}{0}{{X''}_{\hspace{-2mm}v_{e_1,1}}^{-1} {\bf X}_1 
{X''}_{\hspace{-2mm}v_{e_3,2}}^{-1}
{X''}_{\hspace{-2mm}v_{e_3,1}}^{-1}}
\end{align*}
Concentrating on the ${\bf X}_1$'s and ${\bf X}_1^{-1}$'s, one sees that when one multiplies these three matrices, there is no ${\bf X}_1^{-1}$ left, hence the entries are $X''$-Laurent. In Case 11, we have $k(W_{62;rt}) = -2+2=0$, and
\begin{align*}
& \hspace*{-3mm} \til{\bf M}_{W_{62};rt} = \til{\bf M}_{\gamma_{6,{\rm in}}} \til{\bf M}^{\rm left}_r \til{\bf M}_{\gamma_{1,{\rm in}}} \til{\bf M}^{\rm left}_t \til{\bf M}_{\gamma_{2,{\rm out}}} \\
& \hspace*{-3mm}  = 
\til{\bf M}_{\gamma_{6,{\rm in}}}
\smallmatthree{1}{~(1+{X''}_{\hspace{-2mm}v_{e_1,1}}^{-1} {\bf X}_1 {X''}_{\hspace{-2mm}v_r}^{-1}){X''}_{\hspace{-2mm}v_{e_1,1}} ~}{\cancel{{X''}_{\hspace{-2mm}v_{e_1,1}}^{-1}} {\bf X}_1 {X''}_{\hspace{-2mm}v_r}^{-1} {X''}_{\hspace{-2mm}v_{e_1,2}}^{-1}
\cancel{{X''}_{\hspace{-2mm}v_{e_1,1}}} }{0}{ \cancel{ {X''}_{\hspace{-2mm}v_{e_1,1}}^{-1}} {\bf X}_1 {X''}_{\hspace{-2mm}v_r}^{-1} \cancel{ {X''}_{\hspace{-2mm}v_{e_1,1}}} }{\cancel{ {X''}_{\hspace{-2mm}v_{e_1,1}}^{-1}} {\bf X}_1 {X''}_{\hspace{-2mm}v_r}^{-1} {X''}_{\hspace{-2mm}v_{e_1,2}}^{-1}
\cancel{{X''}_{\hspace{-2mm}v_{e_1,1}}} }{0}{0}{\cancel{{X''}_{\hspace{-2mm}v_{e_1,1}}^{-1}} {\bf X}_1 {X''}_{\hspace{-2mm}v_r}^{-1} {X''}_{\hspace{-2mm}v_{e_1,2}}^{-1}
\cancel{{X''}_{\hspace{-2mm}v_{e_1,1}}}}
\smallmatthree{1}{~1+{\bf X}_1^{-1} {X''}_{\hspace{-2mm}v_t}^{-1}~}{{\bf X}_1^{-1} {X''}_{\hspace{-2mm}v_t}^{-1}}{0}{{\bf X}_1^{-1} {X''}_{\hspace{-2mm}v_t}^{-1}}{{\bf X}_1^{-1} {X''}_{\hspace{-2mm}v_t}^{-1}}{0}{0}{{\bf X}_1^{-1} {X''}_{\hspace{-2mm}v_t}^{-1}}
\til{\bf M}_{\gamma_{2,{\rm out}}}
\end{align*}
In the product of the middle two matrices, there are two entries that are not manifestly $X''$-Laurent are: the $(1,2)$-th entry which is
\begin{align*}
1+{\bf X}_1^{-1} {X''}_{\hspace{-2mm}v_t}^{-1}
+ (1+{X''}_{\hspace{-2mm}v_{e_1,1}}^{-1} {\bf X}_1 {X''}_{\hspace{-2mm}v_r}^{-1}){X''}_{\hspace{-2mm}v_{e_1,1}} {\bf X}_1^{-1} {X''}_{\hspace{-2mm}v_t}^{-1} = 1 + \cancel{ {\bf X}_1^{-1} } {X''}_{\hspace{-2mm}v_t}^{-1} \cancel{ (1 + {X''}_{\hspace{-2mm}v_{e_1,1}}) }
+ {X''}_{\hspace{-2mm}v_r}^{-1} {X''}_{\hspace{-2mm}v_t}^{-1}
\end{align*}
which is $X''$-Laurent, and the $(1,3)$-th entry which is ${\bf X}_1^{-1} {X''}_{\hspace{-2mm}v_t}^{-1}$ times eq.\eqref{eq:bf_X_1_canceling4}, hence is $X''$-Laurent. In Case 12, we have $k(W_{63;rt}) = -2-1=-3$, and
\begin{align*}
& \hspace*{-3mm} \til{\bf M}_{W_{63};rt} = \til{\bf M}_{\gamma_{6,{\rm in}}} \til{\bf M}^{\rm left}_r \til{\bf M}_{\gamma_{1,{\rm in}}} \til{\bf M}^{\rm right}_t \til{\bf M}_{\gamma_{3,{\rm out}}} \\
& = 
\til{\bf M}_{\gamma_{6,{\rm in}}}
\smallmatthree{1}{~(1+{X''}_{\hspace{-2mm}v_{e_1,1}}^{-1} {\bf X}_1 {X''}_{\hspace{-2mm}v_r}^{-1}){X''}_{\hspace{-2mm}v_{e_1,1}} ~}{ {\bf X}_1 {X''}_{\hspace{-2mm}v_r}^{-1} {X''}_{\hspace{-2mm}v_{e_1,2}}^{-1} }{0}{ {\bf X}_1 {X''}_{\hspace{-2mm}v_r}^{-1}  }{ {\bf X}_1 {X''}_{\hspace{-2mm}v_r}^{-1} {X''}_{\hspace{-2mm}v_{e_1,2}}^{-1}
 }{0}{0}{ {\bf X}_1 {X''}_{\hspace{-2mm}v_r}^{-1} {X''}_{\hspace{-2mm}v_{e_1,2}}^{-1} }
 \smallmatthree{1}{0}{0}{1}{1}{0}{1}{~1+{\bf X}_1^{-1} {X''}_{\hspace{-2mm}v_t}^{-1}~}{{\bf X}_1^{-1} {X''}_{\hspace{-2mm}v_t}^{-1}}
\til{\bf M}_{\gamma_{3,{\rm out}}}
\end{align*}
In the product of the middle two matrices, the $(1,1)$-th entry is eq.\eqref{eq:bf_X_1_canceling4}, the (1,2)-th entry is eq.\eqref{eq:bf_X_1_canceling5}; hence all entries are $X''$-Laurent, and the entries of the first column are divisible by ${\bf X}_1$. Multiplying by $\til{\bf M}_{\gamma_{3,{\rm out}}} 
= {\rm diag}(1, {X''}_{\hspace{-2mm}v_{e_1,1}}^{-1} {\bf X}_1
{X''}_{\hspace{-2mm}v_{e_3,2}}^{-1}, 
{X''}_{\hspace{-2mm}v_{e_1,1}}^{-1} {\bf X}_1 
{X''}_{\hspace{-2mm}v_{e_3,2}}^{-1}
{X''}_{\hspace{-2mm}v_{e_3,1}}^{-1})$ from right, one observes that all entries are $X''$-Laurent and are divisible by ${\bf X}_1$. Therefore ${\bf X}_1^{k(W_{63;rt})} \til{\bf M}_{W_{63};rt} = {\bf X}_1^{-1} \til{\bf M}_{W_{63};rt}$ is a $X''$-Laurent matrix. 

\vs

Now, we assume that some of the external edges $e_2$, $e_3$, $e_5$, $e_6$ of the quadrilateral formed by $t$, $r$ are identified with each other. We shall see that we can use the computational results obtained so far, and we just have to arrange each situation to fit one of the previous situations. Since we assumed that the triangulation $\Delta$ is regular (Def.\ref{def:regular_triangulation}), the only possibility is that $e_2$ is identified with $e_5$, and/or $e_3$ is identified with $e_6$. In particular, when both identifications happen, it means that the surface $\frak{S}$ is a once-punctured torus and $t,r$ are all the triangles of $\Delta$. Let's begin with the case when $e_2=e_5$ but $e_3 \neq e_6$. Then $v_{e_2,2} = v_{e_5,1}$ and $v_{e_2,1} = v_{e_5,2}$; the monodromy matrices $\til{\bf M}_{\gamma_{2,{\rm in}}}$ and $\til{\bf M}_{\gamma_{2,{\rm out}}}$ should be replaced by $\til{\bf M}_{\gamma_{5,{\rm out}}}$ and $\til{\bf M}_{\gamma_{5,{\rm in}}}$ respectively. The list of complete concatenations of segments in the `quadrilateral' formed by $t,r$ now changes to:
\begin{align*}
& \gamma_{3,{\rm in}} . \gamma_{31} . \gamma_{1,{\rm out}} . (\gamma_{45} . \gamma_{5,{\rm out}} . \gamma_{21} . \gamma_{1,{\rm out}})^n . \gamma_{46} . \gamma_{6,{\rm out}}, \quad n \in \mathbb{Z}_{\ge 0}, \\
& \gamma_{3,{\rm in}} . \gamma_{32} . \gamma_{5,{\rm in}} . ( \gamma_{54} . \gamma_{1,{\rm in}} . \gamma_{12} . \gamma_{5,{\rm in}} )^n . \gamma_{56} . \gamma_{6,{\rm out}}, \quad n\in \mathbb{Z}_{\ge 0}, \\
& \gamma_{6,{\rm in}} . \gamma_{64} . \gamma_{1,{\rm in}}. (\gamma_{12} . \gamma_{5,{\rm in}}. \gamma_{54} . \gamma_{1,{\rm in}})^n . \gamma_{13}. \gamma_{3,{\rm out}}, \quad n \in \mathbb{Z}_{\ge 0}, \\
& \gamma_{6,{\rm in}}. \gamma_{65} . \gamma_{5,{\rm out}} . (\gamma_{21}. \gamma_{1,{\rm out}} . \gamma_{45}. \gamma_{5,{\rm out}})^n . \gamma_{23}. \gamma_{3,{\rm out}}, \quad n \in \mathbb{Z}_{\ge 0}, \\
& \gamma_{5,{\rm in}} . \gamma_{54} . \gamma_{1,{\rm in}} . \gamma_{12}, \qquad 
\gamma_{5,{\rm out}} . \gamma_{21}. \gamma_{1,{\rm out}}. \gamma_{45}.
\end{align*}
In particular, the last two cases are loops by themselves.  One can observe that each of these complete concatenations is a concatenation of the following basic concatenations, which themselves are not necessarily complete concatenations:
\begin{align*}
& \gamma_{3,{\rm in}} . \gamma_{31} . \gamma_{1,{\rm out}} . \gamma_{45} . \gamma_{5,{\rm out}}, \quad
\gamma_{5,{\rm in}} . \gamma_{54} . \gamma_{1,{\rm in}} . \gamma_{13} . \gamma_{3,{\rm out}}, \\
& \gamma_{6,{\rm in}} . \gamma_{64} . \gamma_{1,{\rm in}} . \gamma_{12}, \quad
\gamma_{21} . \gamma_{1,{\rm out}} . \gamma_{46} . \gamma_{6,{\rm out}}
\\
& \gamma_{21} . \gamma_{1,{\rm out}} . \gamma_{45} . \gamma_{5,{\rm out}}, \quad
\gamma_{5,{\rm in}} . \gamma_{54} . \gamma_{1,{\rm in}} . \gamma_{12}, \\
& \gamma_{3,{\rm in}} . \gamma_{32}, \quad
\gamma_{23} . \gamma_{3,{\rm out}}, \quad
\gamma_{6,{\rm in}}. \gamma_{65} . \gamma_{5,{\rm out}}, \quad
\gamma_{5,{\rm in}} . \gamma_{56} . \gamma_{6,{\rm out}}
\end{align*}
Now, one can observe that for each of these basic concatenations, the effect of the mutation on the entries of the corresponding product of non-normalized monodromy matrices is as being asserted in this proof. More precisely, the first one was dealt with in Case 7 in our previous investigation, and the second one in Case 10. The third one $\gamma_{6,{\rm in}} . \gamma_{64} . \gamma_{1,{\rm in}} . \gamma_{12}$ was dealt with in Case 11, the fourth one in Case 6, the fifth one in Case 5, the sixth one in Case 9, the seventh one in Case 2, the eighth one in Case 1, the ninth one in Case 4, and the tenth one in Case 3; just ignore the factors $\til{\bf M}_{\gamma_{2,{\rm in}}}$ and $\til{\bf M}_{\gamma_{2,{\rm out}}}$ in our previous investigations, which were not playing any roles. Similarly, for the case when $e_2\neq e_5$ and $e_3=e_6$, we have $v_{e_6,1} = v_{e_3,2}$ and $v_{e_6,2}=v_{e_3,1}$, so that the monodromy matrices $\til{\bf M}_{\gamma_{6,{\rm in}}}$ and $\til{\bf M}_{\gamma_{6,{\rm out}}}$ should be replaced by $\til{\bf M}_{\gamma_{3,{\rm out}}}$ and $\til{\bf M}_{\gamma_{3,{\rm in}}}$, and each complete concatenation is a concatenation of the following basic concatenations:
\begin{align*}
& \gamma_{2,{\rm in}} . \gamma_{23} . \gamma_{3,{\rm out}}, \quad
\gamma_{3,{\rm in}} . \gamma_{32} . \gamma_{2,{\rm out}}, \quad
\gamma_{2,{\rm in}} . \gamma_{21} . \gamma_{1,{\rm out}} . \gamma_{46}, \quad
\gamma_{64}. \gamma_{1,{\rm in}} . \gamma_{12} . \gamma_{2,{\rm out}}, \\
& \gamma_{5,{\rm in}}. \gamma_{56}, \quad
\gamma_{65} . \gamma_{5,{\rm out}}, \quad
\gamma_{5,{\rm in}} . \gamma_{54} . \gamma_{1,{\rm in}} . \gamma_{13} . \gamma_{3,{\rm out}}, \quad
\gamma_{3,{\rm in}} . \gamma_{31} . \gamma_{1,{\rm out}} . \gamma_{45} . \gamma_{5,{\rm out}}, \\
& \gamma_{3,{\rm in}} . \gamma_{31} . \gamma_{1,{\rm out}} . \gamma_{46}, \quad
\gamma_{64} . \gamma_{1,{\rm in}} . \gamma_{13} . \gamma_{3,{\rm out}}.
\end{align*}
The effect of mutation on the entries of the product of non-normalized monodromy matrices for these basic concatenations were computed in Cases 1, 2, 6, 11, 3, 4, 10, 7, 8, 12, respectively in our previous investigation, where we need to just ignore the factors $\til{\bf M}_{\gamma_{6,{\rm in}}}$ and $\til{\bf M}_{\gamma_{6,{\rm out}}}$ which were not playing any roles before. Lastly, if $e_2=e_5$ and $e_3 = e_6$, then $v_{e_2,2} = v_{e_5,1}$, $v_{e_2,1} = v_{e_5,2}$, $v_{e_6,1} = v_{e_3,2}$, $v_{e_6,2} = v_{e_3,1}$, and each complete concatenation is a concatenation of the following basic concatenations:
\begin{align*}
& \gamma_{3,{\rm in}}. \gamma_{31} . \gamma_{1,{\rm out}} . \gamma_{45} . \gamma_{5,{\rm out}}, \quad
\gamma_{5,{\rm in}} . \gamma_{54} . \gamma_{1,{\rm in}} . \gamma_{13} . \gamma_{3,{\rm out}}, \\
& \gamma_{23} . \gamma_{3,{\rm out}}, \quad
\gamma_{3,{\rm in}} . \gamma_{32}, \quad
\gamma_{21} . \gamma_{1,{\rm out}} . \gamma_{46}, \quad
\gamma_{64} . \gamma_{1,{\rm in}} . \gamma_{12}, \quad
\gamma_{65} . \gamma_{5,{\rm out}}, \quad
\gamma_{5,{\rm in}} . \gamma_{56}, \\
& \gamma_{3,{\rm in}} . \gamma_{31} . \gamma_{1,{\rm out}} . \gamma_{46}, \quad
\gamma_{64} . \gamma_{1,{\rm in}} . \gamma_{13} . \gamma_{3,{\rm out}}, \quad
\gamma_{5,{\rm in}} . \gamma_{54} . \gamma_{1,{\rm in}} . \gamma_{12}, \quad
\gamma_{21} \gamma_{1,{\rm out}} . \gamma_{45} . \gamma_{5,{\rm out}},
\end{align*}
which were dealt with in Cases 7, 10, 1, 2, 6, 11, 4, 3, 8, 12, 9, 5 in our previous investigations, where we need to just ignore the factors $\til{\bf M}_{\gamma_{2,*}}$ and $\til{\bf M}_{\gamma_{6,*}}$ which were not playing any roles \redfix{before}.

\vs

Let's summarize the results so far. Writing the trace-of-monodromy $f^+_\gamma = {\rm tr}({\bf M}_{\gamma_1} \cdots {\bf M}_{\gamma_N})$ along oriented \redfix{non-contractible} simple loop $\gamma$ as a Laurent polynomial in (cube-root) old variables $\{ X_v^{1/3} \, | \, v\in \mathcal{V}(Q_\Delta) \}$, by Prop.\ref{prop:congruence-compatibility_of_terms_of_each_basic_regular_function} we know 
$$
f^+_\gamma \in ( {\textstyle \prod}_{v\in \mathcal{V}(Q_\Delta)} X_v^{{\rm a}_v(\ell)} ) \cdot \mathbb{Z}[X_v^{\pm 1} \, | \, v\in \mathcal{V}(Q_\Delta)\},
$$
where $\ell = \gamma$. We investigated the monodromy matrices ${\bf M}_{\gamma_i}$ in terms of new variables $\{ {X''}_{\hspace{-2mm}v}^{1/3} \, | \, v\in \mathcal{V}(Q'')\}$, and found out that for the entries of the product matrix ${\bf M}_{\gamma_1} \cdots {\bf M}_{\gamma_N}$, the discrepancy between the power of an old variable $X_v$ and the corresponding new variable ${X''}_{\hspace{-2mm}v}$, considered up to integers, occurs only for the node $v_{e_1,1}$ which we are mutating at, where the previous power of $X_{v_{e_1,1}}$ is ${\rm a}_{v_{e_1,1}}(\ell)$ while the new power of ${X''}_{\hspace{-2mm}v_{e_1,1}}$ is $-{\rm a}_{v_{e_1,1}}(\ell) + {\rm a}_{v_{e_3,2}}(\ell) + {\rm a}_{v_r}(\ell)$, which can be written as $-{\rm a}_{v_{e_1,1}}(\ell) + \sum_v [\varepsilon_{v_{e_1,1}, v}]_+ {\rm a}_v(\ell)$.

\vs

The rest of the argument goes as in the proof of Prop.\ref{prop:mutation_of_basic_semi-regular_function_at_interior_node_of_triangle}, as already mentioned. \qed

\vs

\redfix{If one prefers, one can rewrite the statement of Prop.\ref{prop:mutation_of_basic_semi-regular_function_at_edge_node_of_triangle} in the style of eq.\eqref{eq:single_mutation_rewritten}. In order to justify this, this time we need to show that our exponent $-{\rm a}_{v_0}(\ell) + \sum_{v\in \mathcal{V}(Q_\Delta)} [\varepsilon_{v_0, v}]_+ {\rm a}_v(\ell)$ appearing in Prop.\ref{prop:mutation_of_basic_semi-regular_function_at_edge_node_of_triangle} coincides with $- {\rm a}_{v_0}(\ell) + \max( \sum_{v\in \mathcal{V}(Q_\Delta)} [\varepsilon_{v,v_0}]_+ {\rm a}_v(\ell), \, \sum_{v\in \mathcal{V}(Q_\Delta)} [-\varepsilon_{v,v_0}]_+ {\rm a}_v(\ell))$ as appearing in eq.\eqref{eq:single_mutation_rewritten}. The difference between these exponents is either $0$ or equals $\sum_{v\in \mathcal{V}(Q_\Delta)} ([\varepsilon_{v,v_0}]_+-[-\varepsilon_{v,v_0}]_+) {\rm a}_v(\ell) = \sum_{v\in \mathcal{V}(Q_\Delta)} \varepsilon_{v,v_0} {\rm a}_v(\ell)$, which we claim to be an integer. In the notation of Fig.\ref{fig:mutate_two_triangles} with $v_0=v_{e_1,1}$, we have $\sum_v \varepsilon_{v,v_0} {\rm a}_v(\ell) = {\rm a}_{v_{e_5,1}}(\ell) + {\rm a}_{v_t}(\ell) - {\rm a}_{v_r}(\ell) - {\rm a}_{v_{e_3,2}}(\ell)$, which is an integer because it equals the difference between $-{\rm a}_{v_r}(\ell) + {\rm a}_{v_{e_1,1}}(\ell) + {\rm a}_{v_{e_5,1}}(\ell)$ and $-{\rm a}_{v_t}(\ell) + {\rm a}_{v_{e_3,2}}(\ell)+{\rm a}_{v_{e_1,1}}(\ell))$, both of which are integers due to the balancedness result Prop.\ref{prop:tropical_coordinate_is_well-defined}(BE3).}

\vs

\redfix{We now apply} Propositions \ref{prop:mutation_of_basic_semi-regular_function_at_interior_node_of_triangle} and \ref{prop:mutation_of_basic_semi-regular_function_at_edge_node_of_triangle} to $\ell \in \mathscr{A}_{\rm L}(\frak{S};\mathbb{Z})$ lying in $\mathscr{A}_{\Delta}(\mathbb{Z}^{\redfix{T}})$, i.e. when ${\rm a}_v(\ell) \in \mathbb{Z}$ for all nodes $v$ of $Q_\Delta$.
\begin{corollary}
\label{cor:mutation_of_congruent_ell}
Let $\Delta$ be any ideal triangulation \redfix{of a punctured surface $\frak{S}$}. Consider the cluster $\mathscr{X}$-chart associated to $\Delta$, and mutate it at a node of $Q_\Delta$. Denote the resulting quiver by $Q'$, naturally identifying $\mathcal{V}(Q_\Delta)$ and $\mathcal{V}(Q')$. Denoting by $X'_v$ the $\mathscr{X}$-coordinate for the node $v$ of $Q'$ for this new chart after mutation. Then for any $\ell \in \mathscr{A}_{\Delta}(\mathbb{Z}^{\redfix{T}})$, we have
$$
\mathbb{I}^+_{{\rm PGL}_3}(\ell) \in \mathbb{Z}[\{{X'}_{\hspace{-1,5mm}v}^{\pm 1} \, | \, v\in \mathcal{V}(Q')\}]. \qed
$$
\end{corollary}

Thus the argument we gave right after Prop.\ref{prop:one-level_mutations_are_enought} works, and this proves Prop.\ref{prop:congruent_lamination_gives_regular_function} as promised.

\subsection{Proof of the first main theorem}

We prove the first main theorem, Thm.\ref{thm:main}. First, choose any ideal triangulation $\Delta$ of $\frak{S}$. Define the map
$$
\mathbb{I}_\Delta : \mathscr{A}_\Delta(\mathbb{Z}^{\redfix{T}}) \to \mathscr{O}(\mathscr{X}_{{\rm PGL}_3,\frak{S}})
$$
as follows: for each $\ell \in \mathscr{A}_{\Delta}(\mathbb{Z}^{\redfix{T}})$, let $\mathbb{I}_\Delta(\ell) \in \mathscr{O}(\mathscr{X}_{{\rm PGL}_3,\frak{S}})$ be the regular function on $\mathscr{X}_{{\rm PGL}_3,\frak{S}}$ yielding $\mathbb{I}^+_{{\rm PGL}_3}(\ell)$ when evaluated at the semi-field $\mathbb{R}_{>0}$; such $\mathbb{I}_\Delta(\ell)$ exists by Prop.\ref{prop:congruent_lamination_gives_regular_function}, which we have explicitly constructed during the proof. We recall the arguments in \S\ref{subsec:basis}, to prove items Thm.\ref{thm:main}(1)--(4) for $\mathbb{I}_\Delta$.

\vs

(1) To show that $\mathbb{I}_\Delta(\mathscr{A}_{\Delta}(\mathbb{Z}^{\redfix{T}}))$ spans $\mathscr{O}(\mathscr{X}_{{\rm PGL}_3,\frak{S}})$, start from any $f\in \mathscr{O}(\mathscr{X}_{{\rm PGL}_3,\frak{S}})$. Bring it to $P^* f = \sum_{\ell\in \mathscr{A}_{\rm L}(\frak{S};\mathbb{Z})} c_\ell(f) \mathbb{I}_{{\rm SL}_3}(\ell) \in \mathscr{O}(\mathscr{X}_{{\rm SL}_3,\frak{S}})$ as in eq.\eqref{eq:P_star_f}, then pullback by eq.\eqref{eq:section_of_PGL3_to_SL3} to the function $f^+ := \sum_{\ell\in\mathscr{A}_{\rm L}(\frak{S};\mathbb{Z})} c_\ell(f) \mathbb{I}_{{\rm PGL}_3}^+(\ell)$ on $\mathscr{X}^+_{{\rm PGL}_3,\frak{S}}$, which should coincide with the evaluation of $f$ at the semi-field $\mathbb{R}_{>0}$. Since $f$ is regular on $\mathscr{X}_{{\rm PGL}_3,\frak{S}}$, it must be \redfix{a} Laurent polynomial in $\{X_v\,|\,v\in\mathcal{V}(Q_\Delta)\}$ for every ideal triangulation $\Delta$. So, $f^+$ must equal to a Laurent polynomial in $\{X_v\,|\,v\in\mathcal{V}(Q_\Delta)\}$, as a function on $\mathscr{X}^+_{{\rm PGL}_3,\frak{S}}$. By Cor.\ref{cor:congruence_and_integrality_of_powers} (after multiplying an integer to $f^+$), we see that each $\ell \in \mathscr{A}_{\rm L}(\frak{S};\mathbb{Z})$ contributing to the sum all belongs to $\mathscr{A}_{\Delta}(\mathbb{Z}^{\redfix{T}})$. For each $\ell \in\mathscr{A}_{\Delta}(\mathbb{Z}^{\redfix{T}})$, recall that $\mathbb{I}_\Delta(\ell)$ is the element of $\mathscr{O}(\mathscr{X}_{{\rm PGL}_3,\frak{S}})$ yielding $\mathbb{I}^+_{{\rm PGL}_3}(\ell)$ by evaluation at the semi-field $\mathbb{R}_{>0}$. So we have the equality
$$
f = \underset{\ell \in \mathscr{A}_{\Delta}(\mathbb{Z}^{\redfix{T}})}{\textstyle \sum} c_\ell(f) \, \mathbb{I}_\Delta (\ell)
$$
of elements of $\mathscr{O}(\mathscr{X}_{{\rm PGL}_3,\frak{S}})$, where the right hand side is a finite sum. This proves that $\mathbb{I}_\Delta(\mathscr{A}_{\Delta}(\mathbb{Z}^{\redfix{T}}))$ spans $\mathscr{O}(\mathscr{X}_{{\rm PGL}_3,\frak{S}})$. To show the linear independence, suppose the right hand side is zero, as a regular function on $\mathscr{X}_{{\rm PGL}_3,\frak{S}}$. Evaluate at the semi-field $\mathbb{R}_{>0}$, and consider the corresponding sum $\sum_{\ell \in \mathscr{A}_{\Delta}(\mathbb{Z}^{\redfix{T}})} c_\ell(f) \, \mathbb{I}^+_{{\rm PGL}_3}(\ell)$. By using a similar argument as in the proof of Cor.\ref{cor:congruence_and_integrality_of_powers} employing the lexicographic total ordering on the set of all Laurent monomials \redfix{(and taking advantage of the highest-term statement for $\mathbb{I}^+_{{\rm PGL}_3}(\ell)$ in Prop.\ref{prop:highest_term})}, one can show by induction that the coefficients $c_\ell(f)$ must be all zero. Hence the linear independence. This shows that $\mathbb{I}_\Delta(\mathscr{A}_{\Delta}(\mathbb{Z}^{\redfix{T}}))$ is a basis of $\mathscr{O}(\mathscr{X}_{{\rm PGL}_3,\frak{S}})$, and also shows the injectivity of the map $\mathbb{I}_\Delta$.

\vs

(2) This is immediate from the definition of $\mathbb{I}_\Delta$ and Prop.\ref{prop:highest_term}. 

\vs

(3) This is immediate from the definition of $\mathbb{I}_\Delta$ and eq.\eqref{eq:I_plus_peripheral}.

\vs

(4) Let $\ell,\ell' \in \mathscr{A}_{\Delta}(\mathbb{Z}^{\redfix{T}}) \subset \mathscr{A}_{\rm L}(\frak{S};\mathbb{Z})$. By Prop.\ref{prop:I_SL3_is_a_basis}(2) we get the product-to-sum decomposition as in eq.\eqref{eq:I_SL3_structure_constants} for some $c_{{\rm SL}_3}(\ell,\ell';\ell'') \in\mathbb{Z}$; this is an equality of elements of $\mathscr{O}(\mathscr{X}_{{\rm SL}_3,\frak{S}})$. Pulling back by the map eq.\eqref{eq:section_of_PGL3_to_SL3}, we get 
\begin{align}
\label{eq:I_structure_constants_proof}
\mathbb{I}^+_{{\rm PGL}_3}(\ell) \, \mathbb{I}^+_{{\rm PGL}_3}(\ell') = \underset{\ell'' \in \mathscr{A}_{\rm L}(\frak{S};\mathbb{Z})}{\textstyle \sum} c_{{\rm SL}_3}(\ell,\ell';\ell'') \, \mathbb{I}^+_{{\rm PGL}_3}(\ell'')
\end{align}
Now, since $\ell,\ell' \in \mathscr{A}_{\Delta}(\mathbb{Z}^{\redfix{T}})$, both $\mathbb{I}^+_{{\rm PGL}_3}(\ell)$ and $\mathbb{I}^+_{{\rm PGL}_3}(\ell')$, hence also their product, belong to $\mathbb{Z}[\{X_v^{\pm 1} \, | \, v\in \mathcal{V}(Q_\Delta)\}]$. So the right hand side of eq.\eqref{eq:I_structure_constants_proof} belongs to $\mathbb{Z}[\{X_v^{\pm 1} \, | \, v\in \mathcal{V}(Q_\Delta)\}]$. By Cor.\ref{cor:congruence_and_integrality_of_powers}, all $\ell'' \in \mathscr{A}_{\rm L}(\frak{S};\mathbb{Z})$ contributing to the sum in the right hand side belong to $\mathscr{A}_{\Delta}(\mathbb{Z}^{\redfix{T}})$. Then one can recognize that the resulting eq.\eqref{eq:I_structure_constants_proof} is the evaluation at the semi-field $\mathbb{R}_{>0}$ of an equality
$$
\mathbb{I}_\Delta(\ell) \, \mathbb{I}_\Delta(\ell') = \underset{\ell'' \in \mathscr{A}_{\Delta}(\mathbb{Z}^{\redfix{T}})}{\textstyle \sum} c_{{\rm SL}_3}(\ell,\ell';\ell'') \, \mathbb{I}_\Delta(\ell''),
$$
which is the desired statement in item (4). 

\vs

So, for each chosen ideal triangulation $\Delta$, the items (1) and (4) of Thm.\ref{thm:main} hold for $\mathbb{I}_\Delta$, with $\mathscr{A}_{{\rm SL}_3,\frak{S}}(\mathbb{Z}^{\redfix{T}})$ in the statements replaced by $\mathscr{A}_\Delta(\mathbb{Z}^{\redfix{T}})$, while the items (2) and (3) for $\mathbb{I}_\Delta$ hold only for this particular $\Delta$ at the moment.

\vs

Let $\Delta'$ be any other ideal triangulation. Let $\ell \in \mathscr{A}_\Delta(\mathbb{Z}^{\redfix{T}}) \subset \mathscr{A}_{\rm L}(\frak{S};\mathbb{Z})$. Then $\mathbb{I}_\Delta(\ell) \in \mathscr{O}(\mathscr{X}_{{\rm PGL}_3,\frak{S}})$, and since $\mathbb{I}_{\Delta'}(\mathscr{A}_{\Delta'}(\mathbb{Z}^{\redfix{T}}))$ is a basis of $\mathscr{O}(\mathscr{X}_{{\rm PGL}_3,\frak{S}})$ (by item (1) for $\mathbb{I}_{\Delta'}$), we have $\mathbb{I}_\Delta(\ell) = \sum_{\ell' \in \mathscr{A}_{\Delta'}(\mathbb{Z}^{\redfix{T}})} c(\ell') \mathbb{I}_{\Delta'}(\ell')$ for some $c(\ell') \in \mathbb{Q}$ which are zero for all but finitely many $\ell' \in \mathscr{A}_{\Delta'}(\mathbb{Z}^{\redfix{T}}) \subset \mathscr{A}_{\rm L}(\frak{S};\mathbb{Z})$. Evaluating at $\mathbb{R}_{>0}$ we obtain
\begin{align}
\label{eq:I_plus_ell_as_prime}
\mathbb{I}^+_{{\rm PGL}_3}(\ell) = \underset{\ell' \in \mathscr{A}_{\Delta'}(\mathbb{Z}^{\redfix{T}})}{\textstyle \sum} c(\ell') \mathbb{I}^+_{{\rm PGL}_3}(\ell').
\end{align}
Now, view all functions in eq.\eqref{eq:I_plus_ell_as_prime} as Laurent polynomials in $\{X_v^{1/3}\,|\, v\in \mathcal{V}(Q_\Delta)\}$, for $\Delta$. Since the left hand side $\mathbb{I}^+_{{\rm PGL}_3}(\ell)$ belongs to $\mathbb{Z}[\{X_v^{\pm 1}\, |\,v\in \mathcal{V}(Q_\Delta)\}]$ (because $\ell\in \mathscr{A}_\Delta(\mathbb{Z}^{\redfix{T}})$, and by item (2) for $\mathbb{I}_\Delta$), from Cor.\ref{cor:congruence_and_integrality_of_powers} for $\Delta$ we deduce that all $\ell'$'s contributing to the sum belong to $\mathscr{A}_\Delta(\mathbb{Z}^{\redfix{T}})$, hence for each of these $\ell'$ the function $\mathbb{I}^+_{{\rm PGL}_3}(\ell')$ comes from $\mathbb{I}_\Delta(\ell') \in \mathscr{O}(\mathscr{X}_{{\rm PGL}_3,\frak{S}})$. Thus, from eq.\eqref{eq:I_plus_ell_as_prime} we get
$$
\mathbb{I}_\Delta(\ell) = \underset{\ell' \in \mathscr{A}_{\Delta}(\mathbb{Z}^{\redfix{T}})}{\textstyle \sum} c(\ell') \mathbb{I}_\Delta(\ell').
$$
However, since $\mathbb{I}_\Delta$ is injective and $\mathbb{I}_\Delta(\mathscr{A}_\Delta(\mathbb{Z}^{\redfix{T}}))$ is a basis (by item (1) for $\mathbb{I}_\Delta$), it follows that the only contributing $\ell'$ in the right hand side is $\ell' =\ell$ itself, with $c(\ell')=1$. This yields:
\begin{proposition}
For each $\ell \in \mathscr{A}_\Delta(\mathbb{Z}^{\redfix{T}})$, we have $\ell \in \mathscr{A}_{\Delta'}(\mathbb{Z}^{\redfix{T}})$ and $\mathbb{I}_\Delta(\ell) = \mathbb{I}_{\Delta'}(\ell)$.
\end{proposition}
As a corollary, this proves Prop.\ref{prop:congruence_condition_is_indepdent_on_triangulation}, i.e. the sets $\mathscr{A}_\Delta(\mathbb{Z}^{\redfix{T}}) \subset \mathscr{A}_{\rm L}(\frak{S};\mathbb{Z})$ for all ideal triangulations $\Delta$ coincide with each other, as promised, and the first main theorem Thm.\ref{thm:main} holds as is written, with $\mathscr{A}_{{\rm SL}_3,\frak{S}}(\mathbb{Z}^{\redfix{T}})$ being understood as $\mathscr{A}_\Delta(\mathbb{Z}^{\redfix{T}})$ for any ideal triangulation $\Delta$.

\section{The ${\rm SL}_3$ quantum and classical trace maps and the state-sum formula}
\label{sec:SL3_trace}

In the present section we show Prop.\ref{prop:highest_term} and Prop.\ref{prop:congruence-compatibility_of_terms_of_each_basic_regular_function}\redfix{, as well as complete the proofs of Prop.\ref{prop:mutation_of_basic_semi-regular_function_at_interior_node_of_triangle} and Prop.\ref{prop:mutation_of_basic_semi-regular_function_at_edge_node_of_triangle}}, as promised. For this, we develop an ${\rm SL}_3$ version of Bonahon-Wong's ${\rm SL}_2$ quantum trace map \cite{BW}, i.e. the ${\rm SL}_3$ quantum trace. The classical version of it, which we call the ${\rm SL}_3$ classical trace, provides a tool for explicit computation of the basic semi-regular functions $\mathbb{I}^+_{{\rm PGL}_3}(\ell)$, letting us prove the above two propositions. We then explain how the ${\rm SL}_3$ quantum trace can be used for quantizing the functions $\mathbb{I}^+_{{\rm PGL}_3}(\ell)$ and $\mathbb{I}(\ell)$. Notice that, in this section, $\frak{S}$ may be a generalized marked surface having boundary, not even necessarily triangulable.

\subsection{The ${\rm SL}_3$ quantum and classical trace for stated ${\rm SL}_3$-skein algebras}
\label{subsec:SL3_trace_for_stated_A2-skein_algebra}

One goal is to study the properties of the map $\mathbb{I}^+_{{\rm PGL}_3} : \mathscr{A}_{\rm L}(\frak{S};\mathbb{Z}) \to C^\infty(\mathscr{X}^+_{{\rm PGL}_3,\frak{S}})$ defined in Def.\ref{def:translation_to_PGL3}. Crucial is the restriction to $\mathscr{A}^0_{\rm L}(\frak{S};\mathbb{Z})$, which embeds to the ${\rm SL}_3$-skein algebra $\mathcal{S}(\frak{S};\mathbb{Z})$; the image under $\mathbb{I}^+_{{\rm PGL}_3}$ of each element of $\mathcal{S}(\frak{S};\mathbb{Z})$ is a Laurent polynomial in $\{X_v^{1/3} \, | \, v\in \mathcal{V}(Q_\Delta)\}$ per each chosen $\Delta$, and we would like to investigate this Laurent polynomial. Here we will develop an {\em ${\rm SL}_3$ classical trace map}, which is a map from the ${\rm SL}_3$-skein algebra $\mathcal{S}(\frak{S};\mathbb{Z})$ to an abstract Laurent polynomial ring. As a matter of fact, we will construct its quantum version as well, the {\em ${\rm SL}_3$ quantum trace map}.  Although the main interest of the present paper is on the case when $\frak{S}$ is a punctured surface, i.e. without boundary, a full treatment of the ${\rm SL}_3$ classical/quantum trace map requires us to consider the case when $\frak{S}$ has boundary. In fact, the domain of the sought-for quantum trace map is a `stated' version of the {\em noncommutative} ${\rm SL}_3$-skein algebra. In particular, we should now consider ${\rm SL}_3$-webs in the {\em thickened surface} $\frak{S} \times (-1,1)$, which is a 3-dimensional space, and hence the crossings should now carry underpassing/overpassing information. As shall be mentioned again later, if one is interested only in the classical story, then there is no need to go three dimensions, and \redfix{one can} just work with the surface $\frak{S}$ and put $\omega^{1/2}=1$ (or $q=1$), where the computations become much simpler, as done in the first version of the present paper \cite{Kim}. In the present version, we formulate things in the 3d space in order to incorporate the quantum story.
\begin{definition}[\redfix{\cite{S05} \cite{FS} \cite{Higgins}}]
\label{def:webs_in_3d}
Let $(\Sigma,\mathcal{P})$ be a generalized marked surface, and $\frak{S} = \Sigma\setminus\mathcal{P}$. Let
$$
{\bf I} := (-1,1) 
$$
be the open interval in $\mathbb{R}$, and let the 3-manifold
$$
\frak{S} \times {\bf I}
$$
be the \ul{\em thickening} of $\frak{S}$. For a point $(x,t) \in \frak{S} \times {\bf I}$, the ${\bf I}$-coordinate $t$ is called the \ul{\em elevation} of $(x,t)$. For a subset $S$ of $\frak{S}$, we say that a point $(x,t)$ of $\frak{S} \times {\bf I}$ lies \ul{\em over} $S$ if $x\in S$. \redfix{For each boundary arc $b$ of $\frak{S}$, the corresponding boundary component $b\times {\bf I}$ of $\frak{S} \times {\bf I}$ is called a \ul{\em boundary wall}.}

\vs

An \ul{\em ${\rm SL}_3$-web $W$ in (the thickened surface) $\frak{S}\times {\bf I}$} consists of

$\bullet$ a finite subset of $\partial (\frak{S} \times {\bf I}) = (\partial \frak{S}) \times {\bf I}$, whose elements are called \ul{\em external vertices} or \ul{\em endpoints} of $W$, where we let $\partial W$ be the set of all endpoints of $W$;

$\bullet$ a finite subset of $\mathring{\frak{S}} \times {\bf I}$, whose elements are called \ul{\em internal vertices};

$\bullet$ a finite set of non-closed oriented smooth curves in $\redfix{\frak{S}} \times {\bf I}$ \redfix{whose interiors lie in $\mathring{\frak{S}} \times {\bf I}$ and that end} at points in external or internal vertices, whose elements are called \ul{\em (oriented) edges} of $W$;

$\bullet$ a finite set of closed oriented smooth curves in $\mathring{\frak{S}} \times {\bf I}$, whose elements are called \ul{\em (oriented) loops} of $W$,

\vs

subject to the following conditions:
\begin{enumerate}
\itemsep0em
\item[\rm (W5)] each external vertex $v$ is $1$-valent, i.e. exactly one edge of $W$ ends at $v$ and this edge meets $v$ once, \redfix{and $W$ meets a boundary wall transversally at an external vertex};

\item[\rm (W6)] each internal vertex $v$ is either a $3$-valent sink or a $3$-valent source, i.e. exactly three edges of $W$ end at $v$, and the orientations of them are either all toward $v$ or all outgoing from $v$.

\item[\rm (W7)] there is no self-intersection of $W$ except the 3-valent internal vertices;

\item[\rm (W8)] each of the constituent edges and loops of $W$ \redfix{is} given a \ul{\em framing}, i.e. a \redfix{continuous} choice of an element of $T_x (\frak{S}\times{\bf I}) \setminus T_x W$ \redfix{per each $x \in W$}; \redfix{}

\item[\rm (W9)] the framing at each of the \redfix{external} vertices is \ul{\em upward vertical}, i.e. is parallel to the ${\bf I}$ factor and pointing toward $1$;

\item[\rm (W10)] \redfix{for each boundary wall $b \times {\bf I}$, the endpoints of $W$ lying in $b\times {\bf I}$ } have mutually distinct elevations;

\item[\rm \redfix{(W11)}] \redfix{for each internal vertex $x$ of $W$, there is a diffeomorphism from a neighborhood of $x$ in $\frak{S} \times {\bf I}$ to $\mathbb{D} \times {\bf I}$ (where $\mathbb{D}$ is an open disc) s.t. the image of $W$ lies in $\mathbb{D} \times \{0\}$ with upward vertical framing.}
\end{enumerate}
An \ul{\em isotopy} of ${\rm SL}_3$-webs in $\frak{S} \times {\bf I}$ is an isotopy within the class of ${\rm SL}_3$-webs in $\frak{S} \times {\bf I}$.
\end{definition}
Given an ${\rm SL}_3$-web $W$ in $\frak{S}\redfix{\times {\bf I}}$, through an isotopy one can put into the situation such that, if we denote
$$
\pi :\frak{S} \times {\bf I} \to \frak{S}
$$
the natural projection, then
\begin{enumerate}
\item[\rm (P1)] the framing at every point of $W$ is upward vertical;

\item[\rm (P2)] the restriction of the projection $\pi$ to $W$ is at most $2$-to-$1$, where the point of $\pi(W)$ (as well as its preimages under $\pi$) is called a \ul{\em crossing};

\item[\rm (P3)] every crossing occurs in the interior $\redfix{\mathring{\frak{S}}}$ and is a transverse double intersection\redfix{, where no pre-image in $\mathring{\frak{S}} \times {\bf I}$ of a crossing is a 3-valent internal vertex}.
\end{enumerate}

\begin{figure}[htbp!]
\vspace{-2mm}
\begin{center}
\hspace*{-5mm}\begin{tabular}{ c | c | c }
\raisebox{-0.4\height}{
\begingroup%
  \makeatletter%
  \providecommand\color[2][]{%
    \errmessage{(Inkscape) Color is used for the text in Inkscape, but the package 'color.sty' is not loaded}%
    \renewcommand\color[2][]{}%
  }%
  \providecommand\transparent[1]{%
    \errmessage{(Inkscape) Transparency is used (non-zero) for the text in Inkscape, but the package 'transparent.sty' is not loaded}%
    \renewcommand\transparent[1]{}%
  }%
  \providecommand\rotatebox[2]{#2}%
  \newcommand*\fsize{\dimexpr\f@size pt\relax}%
  \newcommand*\lineheight[1]{\fontsize{\fsize}{#1\fsize}\selectfont}%
  \ifx\svgwidth\undefined%
    \setlength{\unitlength}{90.70866142bp}%
    \ifx\svgscale\undefined%
      \relax%
    \else%
      \setlength{\unitlength}{\unitlength * \real{\svgscale}}%
    \fi%
  \else%
    \setlength{\unitlength}{\svgwidth}%
  \fi%
  \global\let\svgwidth\undefined%
  \global\let\svgscale\undefined%
  \makeatother%
  \begin{picture}(1,0.28125)%
    \lineheight{1}%
    \setlength\tabcolsep{0pt}%
    \put(0.23728242,0.11890283){\color[rgb]{0,0,0}\makebox(0,0)[lt]{\lineheight{1.25}\smash{\begin{tabular}[t]{l}$=[3]_q {\O}=$\end{tabular}}}}%
    \put(0,0){\includegraphics[width=\unitlength,page=1]{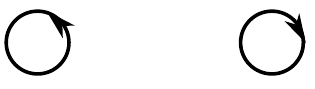}}%
  \end{picture}%
\endgroup%
} & \raisebox{-0.4\height}{
\begingroup%
  \makeatletter%
  \providecommand\color[2][]{%
    \errmessage{(Inkscape) Color is used for the text in Inkscape, but the package 'color.sty' is not loaded}%
    \renewcommand\color[2][]{}%
  }%
  \providecommand\transparent[1]{%
    \errmessage{(Inkscape) Transparency is used (non-zero) for the text in Inkscape, but the package 'transparent.sty' is not loaded}%
    \renewcommand\transparent[1]{}%
  }%
  \providecommand\rotatebox[2]{#2}%
  \newcommand*\fsize{\dimexpr\f@size pt\relax}%
  \newcommand*\lineheight[1]{\fontsize{\fsize}{#1\fsize}\selectfont}%
  \ifx\svgwidth\undefined%
    \setlength{\unitlength}{116.22047244bp}%
    \ifx\svgscale\undefined%
      \relax%
    \else%
      \setlength{\unitlength}{\unitlength * \real{\svgscale}}%
    \fi%
  \else%
    \setlength{\unitlength}{\svgwidth}%
  \fi%
  \global\let\svgwidth\undefined%
  \global\let\svgscale\undefined%
  \makeatother%
  \begin{picture}(1,0.2195122)%
    \lineheight{1}%
    \setlength\tabcolsep{0pt}%
    \put(0.40460658,0.0870799){\color[rgb]{0,0,0}\makebox(0,0)[lt]{\lineheight{1.25}\smash{\begin{tabular}[t]{l}$=-[2]_q$\end{tabular}}}}%
    \put(0,0){\includegraphics[width=\unitlength,page=1]{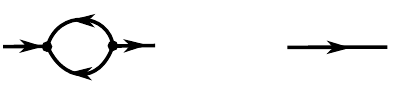}}%
  \end{picture}%
\endgroup%
} & \raisebox{-0.4\height}{
\begingroup%
  \makeatletter%
  \providecommand\color[2][]{%
    \errmessage{(Inkscape) Color is used for the text in Inkscape, but the package 'color.sty' is not loaded}%
    \renewcommand\color[2][]{}%
  }%
  \providecommand\transparent[1]{%
    \errmessage{(Inkscape) Transparency is used (non-zero) for the text in Inkscape, but the package 'transparent.sty' is not loaded}%
    \renewcommand\transparent[1]{}%
  }%
  \providecommand\rotatebox[2]{#2}%
  \newcommand*\fsize{\dimexpr\f@size pt\relax}%
  \newcommand*\lineheight[1]{\fontsize{\fsize}{#1\fsize}\selectfont}%
  \ifx\svgwidth\undefined%
    \setlength{\unitlength}{155.90551181bp}%
    \ifx\svgscale\undefined%
      \relax%
    \else%
      \setlength{\unitlength}{\unitlength * \real{\svgscale}}%
    \fi%
  \else%
    \setlength{\unitlength}{\svgwidth}%
  \fi%
  \global\let\svgwidth\undefined%
  \global\let\svgscale\undefined%
  \makeatother%
  \begin{picture}(1,0.30909091)%
    \lineheight{1}%
    \setlength\tabcolsep{0pt}%
    \put(0.3074048,0.1579046){\color[rgb]{0,0,0}\makebox(0,0)[lt]{\lineheight{1.25}\smash{\begin{tabular}[t]{l}$= $\end{tabular}}}}%
    \put(0,0){\includegraphics[width=\unitlength,page=1]{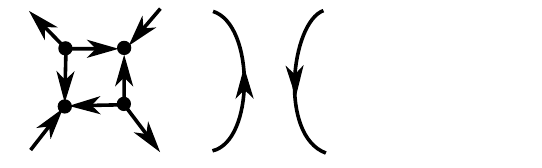}}%
    \put(0.62490477,0.1579046){\color[rgb]{0,0,0}\makebox(0,0)[lt]{\lineheight{1.25}\smash{\begin{tabular}[t]{l}$+$\end{tabular}}}}%
    \put(0,0){\includegraphics[width=\unitlength,page=2]{skein_rel7_corrected.pdf}}%
  \end{picture}%
\endgroup%
} \\
{\rm (S5)} & {\rm (S6)} & {\rm (S7)} \\  \hline
& \raisebox{-0.4\height}{
\begingroup%
  \makeatletter%
  \providecommand\color[2][]{%
    \errmessage{(Inkscape) Color is used for the text in Inkscape, but the package 'color.sty' is not loaded}%
    \renewcommand\color[2][]{}%
  }%
  \providecommand\transparent[1]{%
    \errmessage{(Inkscape) Transparency is used (non-zero) for the text in Inkscape, but the package 'transparent.sty' is not loaded}%
    \renewcommand\transparent[1]{}%
  }%
  \providecommand\rotatebox[2]{#2}%
  \newcommand*\fsize{\dimexpr\f@size pt\relax}%
  \newcommand*\lineheight[1]{\fontsize{\fsize}{#1\fsize}\selectfont}%
  \ifx\svgwidth\undefined%
    \setlength{\unitlength}{150.23622047bp}%
    \ifx\svgscale\undefined%
      \relax%
    \else%
      \setlength{\unitlength}{\unitlength * \real{\svgscale}}%
    \fi%
  \else%
    \setlength{\unitlength}{\svgwidth}%
  \fi%
  \global\let\svgwidth\undefined%
  \global\let\svgscale\undefined%
  \makeatother%
  \begin{picture}(1,0.28301887)%
    \lineheight{1}%
    \setlength\tabcolsep{0pt}%
    \put(0.23913077,0.12612741){\color[rgb]{0,0,0}\makebox(0,0)[lt]{\lineheight{1.25}\smash{\begin{tabular}[t]{l}$=q^{-2/3}$\end{tabular}}}}%
    \put(0,0){\includegraphics[width=\unitlength,page=1]{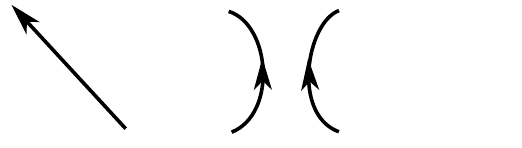}}%
    \put(0.65847036,0.12612741){\color[rgb]{0,0,0}\makebox(0,0)[lt]{\lineheight{1.25}\smash{\begin{tabular}[t]{l}$+q^{1/3}$\end{tabular}}}}%
    \put(0,0){\includegraphics[width=\unitlength,page=2]{skein_rel8_corrected.pdf}}%
  \end{picture}%
\endgroup%
} & \raisebox{-0.4\height}{
\begingroup%
  \makeatletter%
  \providecommand\color[2][]{%
    \errmessage{(Inkscape) Color is used for the text in Inkscape, but the package 'color.sty' is not loaded}%
    \renewcommand\color[2][]{}%
  }%
  \providecommand\transparent[1]{%
    \errmessage{(Inkscape) Transparency is used (non-zero) for the text in Inkscape, but the package 'transparent.sty' is not loaded}%
    \renewcommand\transparent[1]{}%
  }%
  \providecommand\rotatebox[2]{#2}%
  \newcommand*\fsize{\dimexpr\f@size pt\relax}%
  \newcommand*\lineheight[1]{\fontsize{\fsize}{#1\fsize}\selectfont}%
  \ifx\svgwidth\undefined%
    \setlength{\unitlength}{150.23622047bp}%
    \ifx\svgscale\undefined%
      \relax%
    \else%
      \setlength{\unitlength}{\unitlength * \real{\svgscale}}%
    \fi%
  \else%
    \setlength{\unitlength}{\svgwidth}%
  \fi%
  \global\let\svgwidth\undefined%
  \global\let\svgscale\undefined%
  \makeatother%
  \begin{picture}(1,0.28301887)%
    \lineheight{1}%
    \setlength\tabcolsep{0pt}%
    \put(0.23913077,0.12612741){\color[rgb]{0,0,0}\makebox(0,0)[lt]{\lineheight{1.25}\smash{\begin{tabular}[t]{l}$=q^{2/3}$\end{tabular}}}}%
    \put(0,0){\includegraphics[width=\unitlength,page=1]{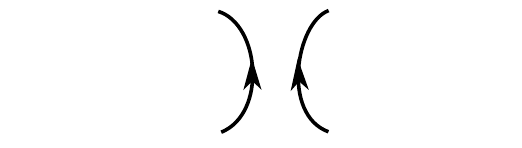}}%
    \put(0.65847036,0.12612741){\color[rgb]{0,0,0}\makebox(0,0)[lt]{\lineheight{1.25}\smash{\begin{tabular}[t]{l}$+q^{-1/3}$\end{tabular}}}}%
    \put(0,0){\includegraphics[width=\unitlength,page=2]{skein_rel9_corrected.pdf}}%
  \end{picture}%
\endgroup%
} \\ 
& {\rm (S8)} & {\rm (S9)} \\ 
\end{tabular}
\end{center}
\vspace{-3mm}
\caption{Non-commutative ${\rm SL}_3$-skein relations, drawn locally (${\O}$ means empty) in $\frak{S}$; the regions bounded by a loop, a $2$-gon, or a $4$-gon in (S5), (S6), (S7) are contractible.}
\vspace{-2mm}
\label{fig:A2-skein_relations_quantum}
\end{figure}

\begin{figure}[htbp!]
\begin{center}
\hspace*{0mm}\begin{tabular}{c|c}
\raisebox{-0.6\height}{
\begingroup%
  \makeatletter%
  \providecommand\color[2][]{%
    \errmessage{(Inkscape) Color is used for the text in Inkscape, but the package 'color.sty' is not loaded}%
    \renewcommand\color[2][]{}%
  }%
  \providecommand\transparent[1]{%
    \errmessage{(Inkscape) Transparency is used (non-zero) for the text in Inkscape, but the package 'transparent.sty' is not loaded}%
    \renewcommand\transparent[1]{}%
  }%
  \providecommand\rotatebox[2]{#2}%
  \newcommand*\fsize{\dimexpr\f@size pt\relax}%
  \newcommand*\lineheight[1]{\fontsize{\fsize}{#1\fsize}\selectfont}%
  \ifx\svgwidth\undefined%
    \setlength{\unitlength}{42.51968504bp}%
    \ifx\svgscale\undefined%
      \relax%
    \else%
      \setlength{\unitlength}{\unitlength * \real{\svgscale}}%
    \fi%
  \else%
    \setlength{\unitlength}{\svgwidth}%
  \fi%
  \global\let\svgwidth\undefined%
  \global\let\svgscale\undefined%
  \makeatother%
  \begin{picture}(1,1.16666667)%
    \lineheight{1}%
    \setlength\tabcolsep{0pt}%
    \put(0,0){\includegraphics[width=\unitlength,page=1]{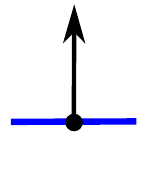}}%
    \put(0.42722329,0.10283054){\color[rgb]{0,0,0}\makebox(0,0)[lt]{\lineheight{1.25}\smash{\begin{tabular}[t]{l}$x$\end{tabular}}}}%
  \end{picture}%
\endgroup%
} \hspace{-3mm}  $=-q^{ - \frac{7}{6}}$ \hspace{-5mm} \raisebox{-0.6\height} {
\begingroup%
  \makeatletter%
  \providecommand\color[2][]{%
    \errmessage{(Inkscape) Color is used for the text in Inkscape, but the package 'color.sty' is not loaded}%
    \renewcommand\color[2][]{}%
  }%
  \providecommand\transparent[1]{%
    \errmessage{(Inkscape) Transparency is used (non-zero) for the text in Inkscape, but the package 'transparent.sty' is not loaded}%
    \renewcommand\transparent[1]{}%
  }%
  \providecommand\rotatebox[2]{#2}%
  \newcommand*\fsize{\dimexpr\f@size pt\relax}%
  \newcommand*\lineheight[1]{\fontsize{\fsize}{#1\fsize}\selectfont}%
  \ifx\svgwidth\undefined%
    \setlength{\unitlength}{42.51968504bp}%
    \ifx\svgscale\undefined%
      \relax%
    \else%
      \setlength{\unitlength}{\unitlength * \real{\svgscale}}%
    \fi%
  \else%
    \setlength{\unitlength}{\svgwidth}%
  \fi%
  \global\let\svgwidth\undefined%
  \global\let\svgscale\undefined%
  \makeatother%
  \begin{picture}(1,1.16666667)%
    \lineheight{1}%
    \setlength\tabcolsep{0pt}%
    \put(0,0){\includegraphics[width=\unitlength,page=1]{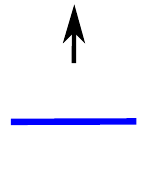}}%
    \put(0.16942518,0.13006975){\color[rgb]{0,0,0}\makebox(0,0)[lt]{\lineheight{1.25}\smash{\begin{tabular}[t]{l}$x_1$\end{tabular}}}}%
    \put(0,0){\includegraphics[width=\unitlength,page=2]{boundary_rel2_quantum.pdf}}%
    \put(0.74016697,0.12462211){\color[rgb]{0,0,0}\makebox(0,0)[lt]{\lineheight{1.25}\smash{\begin{tabular}[t]{l}$x_2$\end{tabular}}}}%
    \put(0.4330944,0.12306698){\color[rgb]{0,0,0}\makebox(0,0)[lt]{\lineheight{1.25}\smash{\begin{tabular}[t]{l}$\prec$\end{tabular}}}}%
  \end{picture}%
\endgroup%
} \hspace{0,5mm} & \raisebox{-0.6\height}{
\begingroup%
  \makeatletter%
  \providecommand\color[2][]{%
    \errmessage{(Inkscape) Color is used for the text in Inkscape, but the package 'color.sty' is not loaded}%
    \renewcommand\color[2][]{}%
  }%
  \providecommand\transparent[1]{%
    \errmessage{(Inkscape) Transparency is used (non-zero) for the text in Inkscape, but the package 'transparent.sty' is not loaded}%
    \renewcommand\transparent[1]{}%
  }%
  \providecommand\rotatebox[2]{#2}%
  \newcommand*\fsize{\dimexpr\f@size pt\relax}%
  \newcommand*\lineheight[1]{\fontsize{\fsize}{#1\fsize}\selectfont}%
  \ifx\svgwidth\undefined%
    \setlength{\unitlength}{42.51968504bp}%
    \ifx\svgscale\undefined%
      \relax%
    \else%
      \setlength{\unitlength}{\unitlength * \real{\svgscale}}%
    \fi%
  \else%
    \setlength{\unitlength}{\svgwidth}%
  \fi%
  \global\let\svgwidth\undefined%
  \global\let\svgscale\undefined%
  \makeatother%
  \begin{picture}(1,1.16666667)%
    \lineheight{1}%
    \setlength\tabcolsep{0pt}%
    \put(0,0){\includegraphics[width=\unitlength,page=1]{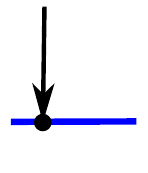}}%
    \put(0.18027884,0.10283054){\color[rgb]{0,0,0}\makebox(0,0)[lt]{\lineheight{1.25}\smash{\begin{tabular}[t]{l}$x_1$\end{tabular}}}}%
    \put(0,0){\includegraphics[width=\unitlength,page=2]{boundary_rel3_quantum.pdf}}%
    \put(0.67416761,0.10283054){\color[rgb]{0,0,0}\makebox(0,0)[lt]{\lineheight{1.25}\smash{\begin{tabular}[t]{l}$x_2$\end{tabular}}}}%
    \put(0.46250103,0.10283054){\color[rgb]{0,0,0}\makebox(0,0)[lt]{\lineheight{1.25}\smash{\begin{tabular}[t]{l}$\prec$\end{tabular}}}}%
  \end{picture}%
\endgroup%
} \hspace{-2mm} $=$ $q$ \hspace{-2mm}
 \raisebox{-0.6\height}{
\begingroup%
  \makeatletter%
  \providecommand\color[2][]{%
    \errmessage{(Inkscape) Color is used for the text in Inkscape, but the package 'color.sty' is not loaded}%
    \renewcommand\color[2][]{}%
  }%
  \providecommand\transparent[1]{%
    \errmessage{(Inkscape) Transparency is used (non-zero) for the text in Inkscape, but the package 'transparent.sty' is not loaded}%
    \renewcommand\transparent[1]{}%
  }%
  \providecommand\rotatebox[2]{#2}%
  \newcommand*\fsize{\dimexpr\f@size pt\relax}%
  \newcommand*\lineheight[1]{\fontsize{\fsize}{#1\fsize}\selectfont}%
  \ifx\svgwidth\undefined%
    \setlength{\unitlength}{42.51968504bp}%
    \ifx\svgscale\undefined%
      \relax%
    \else%
      \setlength{\unitlength}{\unitlength * \real{\svgscale}}%
    \fi%
  \else%
    \setlength{\unitlength}{\svgwidth}%
  \fi%
  \global\let\svgwidth\undefined%
  \global\let\svgscale\undefined%
  \makeatother%
  \begin{picture}(1,1.16666667)%
    \lineheight{1}%
    \setlength\tabcolsep{0pt}%
    \put(0,0){\includegraphics[width=\unitlength,page=1]{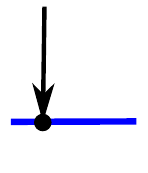}}%
    \put(0.18027884,0.10283054){\color[rgb]{0,0,0}\makebox(0,0)[lt]{\lineheight{1.25}\smash{\begin{tabular}[t]{l}$x_2$\end{tabular}}}}%
    \put(0,0){\includegraphics[width=\unitlength,page=2]{boundary_rel4_quantum.pdf}}%
    \put(0.67416761,0.10283054){\color[rgb]{0,0,0}\makebox(0,0)[lt]{\lineheight{1.25}\smash{\begin{tabular}[t]{l}$x_1$\end{tabular}}}}%
    \put(0.46250103,0.10283054){\color[rgb]{0,0,0}\makebox(0,0)[lt]{\lineheight{1.25}\smash{\begin{tabular}[t]{l}$\prec$\end{tabular}}}}%
  \end{picture}%
\endgroup%
}
 \hspace{-3mm} $+$ \hspace{-2mm}
\raisebox{-0.6\height}{
\begingroup%
  \makeatletter%
  \providecommand\color[2][]{%
    \errmessage{(Inkscape) Color is used for the text in Inkscape, but the package 'color.sty' is not loaded}%
    \renewcommand\color[2][]{}%
  }%
  \providecommand\transparent[1]{%
    \errmessage{(Inkscape) Transparency is used (non-zero) for the text in Inkscape, but the package 'transparent.sty' is not loaded}%
    \renewcommand\transparent[1]{}%
  }%
  \providecommand\rotatebox[2]{#2}%
  \newcommand*\fsize{\dimexpr\f@size pt\relax}%
  \newcommand*\lineheight[1]{\fontsize{\fsize}{#1\fsize}\selectfont}%
  \ifx\svgwidth\undefined%
    \setlength{\unitlength}{42.51968504bp}%
    \ifx\svgscale\undefined%
      \relax%
    \else%
      \setlength{\unitlength}{\unitlength * \real{\svgscale}}%
    \fi%
  \else%
    \setlength{\unitlength}{\svgwidth}%
  \fi%
  \global\let\svgwidth\undefined%
  \global\let\svgscale\undefined%
  \makeatother%
  \begin{picture}(1,1.16666667)%
    \lineheight{1}%
    \setlength\tabcolsep{0pt}%
    \put(0,0){\includegraphics[width=\unitlength,page=1]{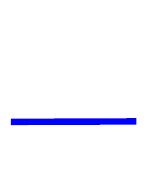}}%
    \put(0.13088999,0.10283052){\color[rgb]{0,0,0}\makebox(0,0)[lt]{\lineheight{1.25}\smash{\begin{tabular}[t]{l}$x_2$\end{tabular}}}}%
    \put(0,0){\includegraphics[width=\unitlength,page=2]{boundary_rel5_quantum.pdf}}%
    \put(0.72355657,0.10283052){\color[rgb]{0,0,0}\makebox(0,0)[lt]{\lineheight{1.25}\smash{\begin{tabular}[t]{l}$x_1$\end{tabular}}}}%
    \put(0,0){\includegraphics[width=\unitlength,page=3]{boundary_rel5_quantum.pdf}}%
    \put(0.44133449,0.10283052){\color[rgb]{0,0,0}\makebox(0,0)[lt]{\lineheight{1.25}\smash{\begin{tabular}[t]{l}$\prec$\end{tabular}}}}%
  \end{picture}%
\endgroup%
}  \\
{\rm (B1)} $s(x)=\varepsilon$, $s(x_1)=r_1(\varepsilon)$, $s(x_2)=r_2(\varepsilon)$ 
& {\rm (B2)} $s(x_1)=\varepsilon_1$, $s(x_2)= \varepsilon_2$, with  $\varepsilon_1>\varepsilon_2$  \\ \hline
\raisebox{-0.6\height}{
\begingroup%
  \makeatletter%
  \providecommand\color[2][]{%
    \errmessage{(Inkscape) Color is used for the text in Inkscape, but the package 'color.sty' is not loaded}%
    \renewcommand\color[2][]{}%
  }%
  \providecommand\transparent[1]{%
    \errmessage{(Inkscape) Transparency is used (non-zero) for the text in Inkscape, but the package 'transparent.sty' is not loaded}%
    \renewcommand\transparent[1]{}%
  }%
  \providecommand\rotatebox[2]{#2}%
  \newcommand*\fsize{\dimexpr\f@size pt\relax}%
  \newcommand*\lineheight[1]{\fontsize{\fsize}{#1\fsize}\selectfont}%
  \ifx\svgwidth\undefined%
    \setlength{\unitlength}{42.51968504bp}%
    \ifx\svgscale\undefined%
      \relax%
    \else%
      \setlength{\unitlength}{\unitlength * \real{\svgscale}}%
    \fi%
  \else%
    \setlength{\unitlength}{\svgwidth}%
  \fi%
  \global\let\svgwidth\undefined%
  \global\let\svgscale\undefined%
  \makeatother%
  \begin{picture}(1,1.16666667)%
    \lineheight{1}%
    \setlength\tabcolsep{0pt}%
    \put(0,0){\includegraphics[width=\unitlength,page=1]{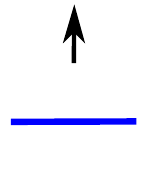}}%
    \put(0.16616776,0.10283052){\color[rgb]{0,0,0}\makebox(0,0)[lt]{\lineheight{1.25}\smash{\begin{tabular}[t]{l}$x$\end{tabular}}}}%
    \put(0,0){\includegraphics[width=\unitlength,page=2]{boundary_rel6_quantum.pdf}}%
    \put(0.72355657,0.10283052){\color[rgb]{0,0,0}\makebox(0,0)[lt]{\lineheight{1.25}\smash{\begin{tabular}[t]{l}$y$\end{tabular}}}}%
    \put(0.44133449,0.10283052){\color[rgb]{0,0,0}\makebox(0,0)[lt]{\lineheight{1.25}\smash{\begin{tabular}[t]{l}$\prec$\end{tabular}}}}%
  \end{picture}%
\endgroup%
} \hspace{-2mm} $=0$
& \raisebox{-0.6\height}{
\begingroup%
  \makeatletter%
  \providecommand\color[2][]{%
    \errmessage{(Inkscape) Color is used for the text in Inkscape, but the package 'color.sty' is not loaded}%
    \renewcommand\color[2][]{}%
  }%
  \providecommand\transparent[1]{%
    \errmessage{(Inkscape) Transparency is used (non-zero) for the text in Inkscape, but the package 'transparent.sty' is not loaded}%
    \renewcommand\transparent[1]{}%
  }%
  \providecommand\rotatebox[2]{#2}%
  \newcommand*\fsize{\dimexpr\f@size pt\relax}%
  \newcommand*\lineheight[1]{\fontsize{\fsize}{#1\fsize}\selectfont}%
  \ifx\svgwidth\undefined%
    \setlength{\unitlength}{42.51968504bp}%
    \ifx\svgscale\undefined%
      \relax%
    \else%
      \setlength{\unitlength}{\unitlength * \real{\svgscale}}%
    \fi%
  \else%
    \setlength{\unitlength}{\svgwidth}%
  \fi%
  \global\let\svgwidth\undefined%
  \global\let\svgscale\undefined%
  \makeatother%
  \begin{picture}(1,1.16666667)%
    \lineheight{1}%
    \setlength\tabcolsep{0pt}%
    \put(0,0){\includegraphics[width=\unitlength,page=1]{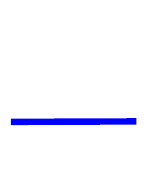}}%
    \put(0.04567213,0.10283052){\color[rgb]{0,0,0}\makebox(0,0)[lt]{\lineheight{1.25}\smash{\begin{tabular}[t]{l}$x_1$\end{tabular}}}}%
    \put(0,0){\includegraphics[width=\unitlength,page=2]{boundary_rel7_quantum.pdf}}%
    \put(0.44839,0.10283052){\color[rgb]{0,0,0}\makebox(0,0)[lt]{\lineheight{1.25}\smash{\begin{tabular}[t]{l}$x_2$\end{tabular}}}}%
    \put(0.83644529,0.10283052){\color[rgb]{0,0,0}\makebox(0,0)[lt]{\lineheight{1.25}\smash{\begin{tabular}[t]{l}$x_3$\end{tabular}}}}%
    \put(0.23672332,0.10283052){\color[rgb]{0,0,0}\makebox(0,0)[lt]{\lineheight{1.25}\smash{\begin{tabular}[t]{l}$\prec$\end{tabular}}}}%
    \put(0.66005655,0.10283052){\color[rgb]{0,0,0}\makebox(0,0)[lt]{\lineheight{1.25}\smash{\begin{tabular}[t]{l}$\prec$\end{tabular}}}}%
  \end{picture}%
\endgroup%
} \hspace{-2mm} $=$ $-q^{\frac{7}{2}}$ \hspace{-5mm} \raisebox{-0.6\height}{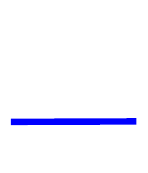} \\ 
{\rm (B3)} $s(x)=s(y)$ & {\rm (B4)} $s(x_1)=1$, $s(x_2)=2$, $s(x_3)=3$
\end{tabular}
\end{center}
\vspace{-2mm}
\caption{Boundary relations for stated ${\rm SL}_3$-skeins (horizontal blue line is boundary); the endpoints in the figure are consecutive in the elevation ordering for that boundary component (i.e. $\nexists$ other endpoint with elevation in between these)}
\vspace{-3mm}
\label{fig:stated_boundary_relations}
\end{figure}

\begin{definition}
If (P1)--(P3) holds, then we encode $W$ as the picture $\pi(W)$ in $\frak{S}$, called the \ul{\em ${\rm SL}_3$-web diagram} of $W$, where for each crossing we need to indicate the over/underpassing information (where the strand of lower elevation is indicated with broken lines), as well as the the ordering of endpoints of $\pi(W)$ lying in each boundary arc of $\frak{S}$ induced from the elevations of corresponding endpoints of $W$. If $x,y \in \pi(W)$ lie in a same boundary component of $\frak{S}$ and $x$ has higher elevation ordering than $y$, we write $x\succ y$.
\end{definition}
\begin{definition}[\cite{S05} \cite{FS} \cite{Higgins}]
\label{def:stated_A2-skein_algebra_quantum}
Let $\frak{S}$ be a generalized marked surface, not necessarily triangulable. 

\vs

$\bullet$ A \ul{\em state} of an ${\rm SL}_3$-web $W$ in $\frak{S} \times {\bf I}$ is a map $s : \partial W \to \{1,2,3\}$ which assigns a number in $\{1,2,3\}$ to each endpoint of $W$, i.e. to each external vertex of $W$. A \ul{\em stated ${\rm SL}_3$-web} in $\frak{S} \times {\bf I}$ is a pair $(W,s)$ of an ${\rm SL}_3$-web $W$ in $\frak{S} \times {\bf I}$ and a state $s$ of $W$.

\vs

$\bullet$ Let $\mathcal{R}$ be a commutative ring with unity $1$. The 
\ul{\em (non-commutative) stated ${\rm SL}_3$-skein algebra} $\mathcal{S}^\omega_{\rm s}(\frak{S};\mathcal{R})$ is the free $\mathcal{R}[\omega^{\pm1/2}]$-module with the set of all isotopy classes of stated ${\rm SL}_3$-webs in $\frak{S} \times {\bf I}$ as a free basis, mod out by the non-commutative ${\rm SL}_3$-skein relations (S5), (S6), (S7), (S8) and (S9) in Fig.\ref{fig:A2-skein_relations_quantum}, where for each positive integer $n$ we denote by
$$
[n]_q = \frac{ q^{n} - q^{-n}}{q - q^{-1}}
$$
the quantum integer, which is a Laurent polynomial in $q$ defined as
\begin{align}
\label{eq:q_and_omega}
q := \omega^9 = (\omega^{1/2})^{18};
\end{align}
so $q^{i/9}$ would mean $\omega^i$, for $i \in \frac{1}{2}\mathbb{Z}$.

\vs

$\bullet$ The \ul{\em reduced stated ${\rm SL}_3$-skein algebra} $\mathcal{S}^\omega_{\rm s}(\frak{S};\mathcal{R})_{\rm red}$ is the quotient of 
$\mathcal{S}^\omega_{\rm s}(\frak{S};\mathcal{R})$ by the \ul{\em boundary relations} in Fig.\ref{fig:stated_boundary_relations}; in the pictures, $x$ and $x_i$ are labels of endpoints, each picture is assumed to carry a respective state which is usually written as $s$, and the \ul{\em index-inversion} $(r_1(\varepsilon),r_2(\varepsilon))$ for $\varepsilon \in \{1,2,3\}$ is given by
\begin{align}
\label{eq:r1_and_r2}
(r_1(1),r_2(1)) = (1,2), \qquad
(r_1(2),r_2(2)) = (1,3), \qquad
(r_1(3),r_2(3)) = (2,3),
\end{align}
which can be though\redfix{t} of as a map from $\{1,2,3\}$ to $\{1,2,3\} \times \{1,2,3\}$.

\vs

$\bullet$ For an equivalence class of a stated ${\rm SL}_3$-web $(W,s)$ in $\frak{S} \times {\bf I}$, the corresponding element of $\mathcal{S}^\omega_{\rm s}(\frak{S};\mathcal{R})$ or $\mathcal{S}^\omega_{\rm s}(\frak{S};\mathcal{R})_{\rm red}$ is denoted by $[W,s]$ and is called a \ul{\em stated ${\rm SL}_3$-skein}.

\vs

$\bullet$ The multiplication in $\mathcal{S}^\omega_{\rm s}(\frak{S};\mathcal{R})$ and that in $\mathcal{S}^\omega_{\rm s}(\frak{S};\mathcal{R})_{\rm red}$ is given by superposition; for two stated ${\rm SL}_3$-skeins $[W,s]$ and $[W',s']$ with $W\subset \frak{S} \times (0,1)$ and $W' \subset \frak{S} \times (-1,0)$, the product $[W,s] \cdot [W',s']$ is defined as $[W\cup W', s\cup s']$, stacking the former one on top of the latter one.

\vs

$\bullet$ Define the \ul{\em (non-commutative) ${\rm SL}_3$-skein algebra} $\mathcal{S}^\omega(\frak{S};\mathcal{R})$ analogously, using ${\rm SL}_3$-webs in $\frak{S}\times {\bf I}$ (without states).
\end{definition}

Some words must be put in order. When $\frak{S}$ has no boundary, the three algebras $\mathcal{S}^\omega_{\rm s}(\frak{S};\mathcal{R})$, $\mathcal{S}^\omega_{\rm s}(\frak{S};\mathcal{R})_{\rm red}$ and $\mathcal{S}^\omega(\frak{S};\mathcal{R})$ coincide. We note that, different authors use different coefficients in the defining relations; the ones in Fig.\ref{fig:A2-skein_relations_quantum} for $\mathcal{S}^\omega_{\rm s}(\frak{S};\mathcal{R})$ are compatible with \cite{FS}, while the ones in Fig.\ref{fig:A2-skein_relations_quantum} and Fig.\ref{fig:stated_boundary_relations} for $\mathcal{S}^\omega_{\rm s}(\frak{S};\mathcal{R})_{\rm red}$ are not compatible with \cite{Higgins}. The above choice of ours turns out to be most suitable for our purposes. As we will make use of the results in \cite{Higgins}, we present the precise isomorphism between our reduced ${\rm SL}_3$-skein algebra $\mathcal{S}^\omega_{\rm s}(\frak{S};\mathcal{R})_{\rm red}$ and that studied by Higgins \cite{Higgins}, denoted by $\mathcal{S}^{SL_3}_q(\Sigma)$ there:
\begin{align}
\label{eq:isomorphism_from_ours_to_Higgins}
& \,\, \mathcal{S}^\omega_{\rm s}(\frak{S};\mathcal{R})_{\rm red} ~\longrightarrow~ \mathcal{S}^{SL_3}_q(\Sigma) \\
\nonumber
& : q^{1/18}=\omega \mapsto q^{-1}, \qquad [W,s] \mapsto \alpha_1^{n^{\rm out}(W)} \alpha_2^{n^{\rm in}(W)} \alpha_3^{n^{\rm out}_1(W,s)} \alpha_4^{n^{\rm in}_1(W,s)} [W,s'],
\end{align}
where $s'(x) = s(x) -2 \in \{-1,0,1\}$, matching the Higgins' convention on the state values $-,0,+$ (or $-1,0,1$), 
\begin{align*}
n^{\rm out}(W) & = \mbox{the number of $3$-valent internal vertices of $W$ that are sources (i.e. outgoing)}, \\
n^{\rm in}(W) & = \mbox{the number of $3$-valent internal vertices of $W$ that are sinks (i.e. incoming)}, \\
n^{\rm out}_1(W,s) & = \mbox{the number of $1$-valent external vertices $x$ of $W$ with $s(x)=1$ that are sinks}, \\
& \quad\mbox{(i.e. orientation of $W$ near $x$ is outbound, toward the boundary of $\frak{S}$)} \\
n^{\rm in}_1(W,s) & = \mbox{the number of $1$-valent external vertices $x$ of $W$ with $s(x)=1$ that are sources}, \\
& \quad\mbox{(i.e. orientation of $W$ near $x$ is inbound, toward the interior of $\frak{S}$)}
\end{align*}
and the twisting scalars are
\begin{align}
\label{eq:twisting_scalars}
\alpha_1 = q^{-5/2}, \qquad
\alpha_2 = q^{-1/2}, \qquad
\alpha_3 = -q, \qquad
\alpha_4 = -q^{-1}.
\end{align}
For example, in \cite{Higgins}, the relation (B1) is written so that the coefficient in the right hand side is $(-1)^\varepsilon q^{-\frac{1}{3} - (\varepsilon-2)}$. The relations in Fig.\ref{fig:A2-skein_relations_quantum} and Fig.\ref{fig:stated_boundary_relations} are translated from the relations in \cite{Higgins} via the isomorphism in eq.\eqref{eq:isomorphism_from_ours_to_Higgins}. As shall be seen, choosing a suitable isomorphism is a crucial step.

\vs

Throughout the present section, $\omega^{1/2}$ and $q$ will denote the quantum parameter, related to each other as in eq.\eqref{eq:q_and_omega}. The following well-known notion will become handy.
\begin{definition}
\label{def:Weyl-ordered}
For any formal variables $\wh{Z}_1,\ldots,\wh{Z}_n$ that $\omega$-commute, in the sense that $\wh{Z}_i \wh{Z}_j = \omega^{\sigma_{ij}} \wh{Z}_j \wh{Z}_i$ for some integer matrix $(\sigma_{ij})$, we define the \ul{\em Weyl-ordered} product (or monomial) as
\begin{align}
[\wh{Z}_1 \cdot \cdots \cdot \wh{Z}_n]_{\rm Weyl} := \omega^{-\frac{1}{2} \sum_{i<j} \sigma_{ij}} \wh{Z}_1 \cdot \cdots \cdot \wh{Z}_n.
\end{align}
For a non-commutative Laurent polynomial $\wh{f}$ in such variables $\wh{Z}_1^{\pm 1},\ldots,\wh{Z}_n^{\pm 1}$ with coefficients in $\mathbb{Z}[\omega^{\pm 1/2}]$, denote by $[\wh{f}]_{\rm Weyl}$ the Laurent polynomial obtained from $\wh{f}$ by replacing each appearing Laurent monomial term $\omega^? \wh{Z}_1^{k_1} \cdots \wh{Z}_n^{k_n}$ by its Weyl-ordered version $[\wh{Z}_1^{k_1}\cdots \wh{Z}_n^{k_n}]_{\rm Weyl}$. For convenience, we also define
$$
[\omega^m \wh{Z}_1 \cdot \cdots \cdot \wh{Z}_n]_{\rm Weyl}
:= [\wh{Z}_1 \cdot \cdots \cdot \wh{Z}_n]_{\rm Weyl}
$$
for $m\in \frac{1}{2}\mathbb{Z}$.

\vs

For a matrix $\wh{\bf M}$ with entries being Laurent polynomials in such variables $\wh{Z}_1^{\pm 1},\ldots,\wh{Z}_n^{\pm 1}$, denote by $[\wh{\bf M}]_{\rm Weyl}$ the matrix obtained from $\wh{\bf M}$ by replacing all entries by their Weyl-ordered versions.
\end{definition}

We now define the target ring of the sought-for quantum trace map.
\begin{definition}[Fock-Goncharov algebra; \cite{FG06} \cite{GS19} \cite{Douglas} \cite{Douglas21}]
\label{def:Fock-Goncharov_algebra_quantum}
Let $\Delta$ be an ideal triangulation of a triangulable generalized marked surface $\frak{S}$.

\vs

For each ideal triangle $t$ of $\Delta$, let $e_1,e_2,e_3$ be the sides of $t$  appearing clockwise in $\partial t$ in this order. On each side $e_\alpha$, let $v_{e_\alpha,1}$, $v_{e_\alpha,2}$ be the nodes of $Q_\Delta$ on $e_{t_\alpha}$ such that the direction $v_{e_\alpha,1} \to v_{e_\alpha,2}$ matches the clockwise orientation of $\partial t$. Let $v_t$ be the node of $Q_\Delta$ in the interior of $t$ (see the left triangle of Fig.\ref{fig:mutate_two_triangles}). Define the \ul{\em cube-root Fock-Goncharov triangle algebra} (or just \ul{\em triangle algebra} in short) $\mathcal{Z}_t^\omega$ as the free associative $\mathbb{Z}[\omega^{\pm 1/2}]$-algebra generated by $\wh{Z}_{t,v}$'s and their inverses for the seven nodes $v$ of $Q_\Delta$ appearing in $t$, with relations
$$
\wh{Z}_{t,v} \, \wh{Z}_{t,w} = \omega^{2\wh{\varepsilon}_{vw;t}}  \wh{Z}_{t,w} \wh{Z}_{t,v}, \qquad \forall v,w \in \mathcal{V}(Q_\Delta) \cap t,
$$
where
\begin{align}
\label{eq:wh_epsilon}
\wh{\varepsilon}_{vw;t} = \left\{
\begin{array}{cl}
\frac{1}{2} & \mbox{if $(v,w) = (v_{e_\alpha,1},v_{e_\alpha,2})$ for some $\alpha\in\{1,2,3\}$}, \\
-\frac{1}{2} & \mbox{if $(v,w) = (v_{e_\alpha,2},v_{e_\alpha,1})$ for some $\alpha\in\{1,2,3\}$}, \\
\varepsilon_{vw} & \mbox{otherwise}.
\end{array}
\right.
\end{align}

\vs

Consider tensor product algebra $\bigotimes_{t\in \mathcal{F}(\Delta)} \mathcal{Z}_t^\omega$, where each $\mathcal{Z}_t^\omega$ naturally embeds into, where $\mathcal{F}(\Delta)$ is the set of all ideal triangles of $\Delta$.

\vs

For each node $v$ of $Q_\Delta$, define the element $\wh{Z}_v$ of the tensor product algebra $\bigotimes_{t\in \mathcal{F}(\Delta)} \mathcal{Z}_t^\omega$ as follows:
\begin{enumerate}
\itemsep0em
\item[\rm (1)] If $v$ is an interior node $v_t$ of some triangle $t$, then $\wh{Z}_v := \wh{Z}_{t,v_t}$.

\item[\rm (2)] If $v$ is a node $v_{e_\alpha,i}$ lying in a boundary arc of $\frak{S}$, and if this node lies in triangle $t$, then $\wh{Z}_v := \wh{Z}_{t,v_{e_\alpha,i}}$.

\item[\rm (3)] If $v$ is a node lying in an interior arc of $\Delta$ so that it equals $v_{e_\alpha,1}$ of a triangle $t$ and $v_{e_\beta,2}$ of a triangle $r$, then $\wh{Z}_v := \wh{Z}_{t, v_{e_\alpha,1}} \, \wh{Z}_{r, v_{e_\beta,2}}$.

\end{enumerate}

Let $\mathcal{Z}^\omega_\Delta = \mathcal{Z}^\omega_{\Delta;\frak{S}}$ be the \ul{\em cube-root Fock-Goncharov algebra} for $\Delta$, defined as the subalgebra of $\bigotimes_{t\in \mathcal{F}(\Delta)} \mathcal{Z}_t^\omega$ generated by $\{ \wh{Z}_v^{\pm 1} \, | \, v\in \mathcal{V}(Q_\Delta)\}$.

\vs

Let
$$
\wh{X}_{t,v} := \wh{Z}_{t,v}^3, \qquad \wh{X}_v := \wh{Z}_v^3.
$$
Then the subalgebra of $\mathcal{Z}^\omega_\Delta$ generated by $\{ \wh{X}_v^{\pm 1} \, | \, v\in \mathcal{V}(Q_\Delta)\}$ is called the \ul{\em Fock-Goncharov algebra} for $\Delta$, denoted by $\mathcal{X}^q_\Delta$.
\end{definition}
In particular, note that the cube-root Fock-Goncharov algebra $\mathcal{Z}^\omega_\Delta$ for $\Delta$ satisfies
$$
\wh{Z}_v \wh{Z}_w = \omega^{2 \wh{\varepsilon}_{vw}} \wh{Z}_w \wh{Z}_v, \quad \forall v,w\in\mathcal{V}(Q_\Delta),
$$
where
\begin{align}
\nonumber
\wh{\varepsilon}_{vw} = \left\{
\begin{array}{cl}
\frac{1}{2} & \mbox{if $(v,w) = (v_{e_\alpha,1},v_{e_\alpha,2})$ for some side $e_\alpha$ of a triangle of $\Delta$ that is a boundary arc of $\Delta$}, \\
-\frac{1}{2} & \mbox{if $(v,w) = (v_{e_\alpha,2},v_{e_\alpha,1})$ for some side $e_\alpha$ of a triangle of $\Delta$ that is a boundary arc of $\Delta$}, \\
\varepsilon_{vw} & \mbox{otherwise}.
\end{array}
\right.
\end{align}
Taking the cubes, we have
$$
\wh{X}_{t,v} \wh{X}_{r,w} = q^{2 \delta_{t,r} \wh{\varepsilon}_{vw;t}} \wh{X}_{r,w} \wh{X}_{t,v}, \qquad
\wh{X}_v \wh{X}_w = q^{2 \wh{\varepsilon}_{vw}} \wh{X}_w \wh{X}_v.
$$
As the notation suggests, $\wh{X}_v$ will be the quantum counterpart of the classical coordinate function $X_v$. The $\pm \frac{1}{2}$ in the definition eq.\eqref{eq:wh_epsilon} stands for the `hidden' arrows between the nodes lying in the sides of a triangle. As informed to us by Daniel Douglas, this idea first appeared in \cite{FG06b}, and used in the quantum setting in \cite{SchSh} \cite{GS19} \cite{CS} \cite{Douglas} \cite{Douglas21}. Note that $\wh{\varepsilon}_{vw} = \varepsilon_{vw}$ if $\frak{S}$ is a punctured surface.

\vs

One more technical preliminary is the following.
\begin{lemma}[cutting process]
\label{lem:cutting_process}
Let $(\Sigma,\mathcal{P})$ be a generalized marked surface, with $\frak{S} = \Sigma\setminus\mathcal{P}$. Let $e$ be a $\mathcal{P}$-arc in $\Sigma$ (or an ideal arc in $\frak{S}$) whose interior lies in the interior of $\Sigma$. Cutting $(\Sigma,\mathcal{P})$ along $e$ yields a possibly-disconnected generalized marked surface $(\Sigma',\mathcal{P}')$, uniquely determined up to \redfix{orientation-preserving} diffeomorphism. Denote by $g : \Sigma' \to \Sigma$ a corresponding gluing map, and let ${\bf g} := g \times {\rm id} : \frak{S} \times {\bf I} \to \frak{S}' \times {\bf I}$ be the gluing map for the thickened surfaces, where $\frak{S}' = \Sigma'\setminus \mathcal{P}'$.

\vs

Let $(W,s)$ be a stated ${\rm SL}_3$-web in the thickening $\frak{S}\times {\bf I}$ of $\frak{S}$ such that $W' := {\bf g}^{-1}(W)$ is an ${\rm SL}_3$-web in the thickening $\frak{S}'\times{\bf I}$ of $\frak{S}'$; we say that $W'$ is obtained from $W$ by cutting along $e$. A state $s' : \partial W' \to \{1,2,3\}$ of $W'$ is said to be \ul{\em compatible} with $s$ if $s'(x) = s({\bf g}(x))$ for each $x \in \partial W' \cap {\bf g}^{-1}(\partial W)$ and $s'(x_1)=s'(x_2)$ for all $x_1,x_2 \in \partial W' \cap {\bf g}^{-1}(e \redfix{\times {\bf I}})$ such that ${\bf g}(x_1)={\bf g}(x_2)$.

\vs

If $\frak{S}$ is triangulable, with a chosen ideal triangulation $\Delta$, then $\Delta' := g^{-1}(\Delta)$ is an ideal triangulation of $\frak{S}'$, which is said to be obtained from $\Delta$ by cutting along $e$. Triangles of $\Delta$ are naturally in bijection with $\Delta'$, where for each triangle $t$ of $\Delta$ there is a canonical bijection from the sides of $t$ to the sides of the corresponding triangle $t'$ of $\Delta'$. The induced isomorphism $\mathcal{Z}_t^\omega \to \mathcal{Z}_{t'}^\omega$ between triangle algebras naturally induces the injection
\begin{align}
\label{eq:i_Delta_Delta_prime}
i_{\Delta,\Delta'} : \mathcal{Z}_\Delta^\omega \to \mathcal{Z}_{\Delta'}^\omega. \qed
\end{align}
\end{lemma}
For convenience, we define:
\begin{definition}
For a generalized marked surface $\frak{S}$, an ${\rm SL}_3$-web $W$ in $\frak{S} \times {\bf I}$ is called a \ul{\em 3-way web} over $\frak{S}$ if it has three external vertices, one internal vertex, no crossing, and has only one component, which consists of three edges, all meeting the internal vertex.
\end{definition}

The main object of study of the present section is the following \ul{\em ${\rm SL}_3$ quantum trace map}, forming the second main theorem of the present paper.
\begin{theorem}[the second main theorem: the ${\rm SL}_3$ quantum trace map]
\label{thm:SL3_quantum_trace_map}
There exists a family of $\mathbb{Z}[\omega^{\pm 1/2}]$-algebra homomorphisms
$$
{\rm Tr}^\omega_\Delta = {\rm Tr}^\omega_{\Delta;\frak{S}} ~:~ \mathcal{S}^\omega_{\rm s}(\frak{S};\mathbb{Z})_{\rm red} \longrightarrow \mathcal{Z}_\Delta^\omega
$$
defined for each triangulable generalized marked surface $\frak{S}$ and each ideal triangulation $\Delta$ \redfix{of $\frak{S}$}, such that:
\begin{enumerate}
\item[\rm \redfix{(QT1)}] (cutting/gluing property) Let $(W,s)$ be a stated ${\rm SL}_3$-web in the thickened surface $\frak{S}\times {\bf I}$, and $e$ \redfix{be} a constituent arc of $\Delta$ that is not a boundary arc of $\frak{S}$. Let $\frak{S}'$ be the generalized marked surface obtained from $\frak{S}$ by cutting along $e$, $\Delta'$ be the triangulation of $\frak{S}'$ obtained from $\Delta$ by cutting along $e$, and $W'$ be the ${\rm SL}_3$-web in $\frak{S}' \redfix{\times} {\bf I}$ obtained from $W$ by cutting along $e$ (Lem.\ref{lem:cutting_process}). Then
\begin{align}
i_{\Delta,\Delta'} {\rm Tr}^\omega_{\Delta;\frak{S}} ([W,s]) = \underset{s'}{\textstyle \sum} \, {\rm Tr}^\omega_{\Delta';\frak{S}'}([W',s']),
\end{align}
where the sum is over all states $s'$ of $W'$ that are compatible with $s$ in the sense as in Lem.\ref{lem:cutting_process}, and $i_{\Delta,\Delta'}$ is as in eq.\eqref{eq:i_Delta_Delta_prime}.
\item[\rm \redfix{(QT2)}] (values at triangle) Let $(W,s)$ be a stated ${\rm SL}_3$-web in the thickening $t\times {\bf I}$ of a triangle $t$, viewed as a generalized marked surface with a unique ideal triangulation $\Delta$. Denote the sides of $t$ by $e_1,e_2,e_3$, and the nodes of $Q_\Delta$ lying in $t$ by $v_{e_\alpha,1}, v_{e_\alpha,2},v_t$ (for $\alpha=1,2,3$) as in Def.\ref{def:Fock-Goncharov_algebra_quantum}. If $\wh{Z}, \wh{Z}_1,\wh{Z}_2$ are invertible elements of $\mathcal{Z}^\omega_t$, define the following $3\times 3$ matrices with entries in $\mathcal{Z}^\omega_t$:
\begin{align}
\label{eq:quantum_edge_matrix}
& \wh{\bf M}^{\rm edge}(\wh{Z}_1,\wh{Z}_2) = [ {\rm diag}(\wh{Z}_1\wh{Z}_2^2, \wh{Z}_1 \wh{Z}_2^{-1}, \wh{Z}_1^{-2}\wh{Z}_2^{-1}) ]_{\rm Weyl}, \\
\label{eq:quantum_turn_matrix}
& \wh{\bf M}^{\rm left}(\wh{Z}) = \smallmatthree{\omega^{5} \wh{Z}^2}{~ (\omega^{-1} \wh{Z}^2 + \omega^2 \wh{Z}^{-1}) ~}{\omega^{-4} \wh{Z}^{-1}}{0}{\omega^5 \wh{Z}^{-1}}{\omega^{-1} \wh{Z}^{-1}}{0}{0}{\omega^2 \wh{Z}^{-1}}, \qquad
\wh{\bf M}^{\rm right}(\wh{Z}) = \smallmatthree{\omega^{-2} \wh{Z}}{0}{0}{\omega \wh{Z}}{ \omega^{-5} \wh{Z}}{0}{\omega^4 \wh{Z}}{~(\omega^{-2} \wh{Z} + \omega \wh{Z}^{-2})~}{\omega^{-5} \wh{Z}^{-2}},
\end{align}
which can be viewed as quantum versions of the monodromy matrices (MM1)--(MM3) of \S\ref{subsec:lifting_PGL3_to_SL3}. Define also the following transposed versions
\begin{align}
\label{eq:quantum_turn_matrix_tran}
& \wh{\bf M}^{\rm left}_{\rm tran}(\wh{Z}) = \smallmatthree{\omega^{-5} \wh{Z}^2}{0}{0}{(\omega \wh{Z}^2 + \omega^{-2} \wh{Z}^{-1}) ~}{\omega^{-5} \wh{Z}^{-1}}{0}{\omega^{4} \wh{Z}^{-1}}{ \omega \wh{Z}^{-1}}{ \omega^{-2} \wh{Z}^{-1}}, \qquad
\wh{\bf M}^{\rm right}_{\rm tran}(\wh{Z}) = \smallmatthree{ \omega^{2} \wh{Z}}{~\omega^{-1} \wh{Z}~}{\omega^{-4} \wh{Z}}{0}{ \omega^{5} \wh{Z}}{~(\omega^{2} \wh{Z} + \omega^{-1} \wh{Z}^{-2})}{0}{0}{\omega^{5} \wh{Z}^{-2}}, \hspace{-5mm}
\end{align}
which may not be honest transposes of $\wh{\bf M}^{\rm left}(\wh{Z})$ and $\wh{\bf M}^{\rm right}(\wh{Z})$. For each $\alpha=1,2,3$, define
\begin{align}
\label{eq:quantum_in_and_out}
\wh{\bf M}^{\rm out}_{t,\alpha} = \wh{\bf M}^{\rm edge}(\wh{Z}_{t,v_{e_\alpha,1}}, \wh{Z}_{t,v_{e_\alpha,2}}), \qquad
\wh{\bf M}^{\rm in}_{t,\alpha} = \wh{\bf M}^{\rm edge}(\wh{Z}_{t,v_{e_\alpha,2}}, \wh{Z}_{t,v_{e_\alpha,1}}).
\end{align}

\vs

For each $k \in \{+,-\}$, $\varepsilon \in \{1,2,3\}$, and $h \in \{{\rm in}, {\rm out}\}$, define the \ul{\em fork matrix} $\wh{\bf F}^{h}_{k,\varepsilon}$ as the following $3\times 3$ matrix with entries in $\mathbb{Z}[\omega^{\pm 1/2}]$ :
\begin{align}
\label{eq:fork_matrix}
(\wh{\bf F}^{h}_{k,\varepsilon} )_{i,j} = \left\{
\begin{array}{ll}
c^h_{k,\varepsilon} & \mbox{if $(i,j)=(r_1(\varepsilon),r_2(\varepsilon))$}, \\
d^h_{k,\varepsilon} & \mbox{if $(i,j)=(r_2(\varepsilon),r_1(\varepsilon))$}, \\
0 &\mbox{otherwise}.
\end{array}
\right.
\end{align}
where $r_1(\varepsilon)$ and  $r_2(\varepsilon)$ are as in eq.\eqref{eq:r1_and_r2}, and
\begin{align*}
c^h_{+,\varepsilon} = \omega^{3/2}, \quad
d^h_{+,\varepsilon} = -\omega^{21/2}, \quad
c^h_{-,\varepsilon} = \omega^{-3/2}, \quad
d^h_{-,\varepsilon} = -\omega^{-21/2},
\end{align*}
for all $\varepsilon \in \{1,2,3\}$ and $h\in\{{\rm in}, {\rm out}\}$.

\vs

Define the \ul{\em twisted (positive) fork matrices} $\til{\bf F}^h_{+,\varepsilon}$ as

\begin{align}
\label{eq:fork_matrix_tilde}
(\til{\bf F}^{h}_{+,\varepsilon} )_{i,j} = \left\{
\begin{array}{ll}
\omega^{\frac{3}{2}p(\varepsilon)} (\wh{\bf F}^{h}_{+,\varepsilon} )_{i,j}, & \mbox{if $h={\rm out}$, $i \le j$, or if $h={\rm in}$, $i\ge j$}, \\
\omega^{-\frac{3}{2}p(\varepsilon)} (\wh{\bf F}^{h}_{+,\varepsilon} )_{i,j}, & \mbox{if $h={\rm out}$, $i \ge j$, or if $h={\rm in}$, $i\le j$},
\end{array}
\right.
\end{align}
where 
\begin{align}
\label{eq:p}
p(\varepsilon) = \left\{
\begin{array}{ll}
1 & \mbox{if $\varepsilon \in \{1,3\}$}, \\
-1 & \mbox{if $\varepsilon=2$.}
\end{array}
\right.
\end{align}

\vs

For each $\varepsilon_1,\varepsilon_2,\varepsilon_3 \in \{1,2,3\}$ and $h \in \{{\rm in},{\rm out}\}$, define
$$
\eta^h(\varepsilon_1,\varepsilon_2,\varepsilon_3) \in \mathbb{Z}
$$
as the unique integer satisfying
\begin{align}
\label{eq:eta_h}
(\wh{\bf M}^h_{t,1})_{\varepsilon_1,\varepsilon_1}
(\wh{\bf M}^h_{t,3})_{\varepsilon_3,\varepsilon_3}
(\wh{\bf M}^h_{t,2})_{\varepsilon_2,\varepsilon_2}
= \omega^{\eta^h(\varepsilon_1,\varepsilon_2,\varepsilon_3)} [(\wh{\bf M}^h_{t,1})_{\varepsilon_1,\varepsilon_1}
(\wh{\bf M}^h_{t,3})_{\varepsilon_3,\varepsilon_3}
(\wh{\bf M}^h_{t,2})_{\varepsilon_2,\varepsilon_2}]_{\rm Weyl}
\end{align}

We now define the \ul{\em outgoing $3$-way matrix} $\wh{\bf M}_t^{\rm 3out}$ and the \ul{\em incoming $3$-way matrix} $\wh{\bf M}_t^{\rm 3in}$ as $3\times 3\times 3$ matrices with entries in $\mathcal{Z}_t^\omega$ as follows. For each $\varepsilon_1,\varepsilon_2,\varepsilon_3 \in \{1,2,3\}$, the $(\varepsilon_1,\varepsilon_2,\varepsilon_3)$-th entries $(\wh{\bf M}_t^{\rm 3out})_{\varepsilon_1,\varepsilon_2,\varepsilon_3} \in \mathcal{Z}^\omega_t$ and $(\wh{\bf M}_t^{\rm 3out})_{\varepsilon_1,\varepsilon_2,\varepsilon_3} \in \mathcal{Z}^\omega_t$ are given as
\begin{align*}
(\wh{\bf M}_t^{\rm 3out})_{\varepsilon_1,\varepsilon_2,\varepsilon_3}
& = \omega^{\eta^{\rm out}(\varepsilon_1,\varepsilon_2,\varepsilon_3)} \, \underset{\varepsilon_4,\varepsilon_5}{\textstyle \sum} ( \wh{\bf M}^{\rm left}_{\rm tran}(\omega^{2g(\varepsilon_2)+2g(\varepsilon_3)} \wh{Z}_{v_t} ))_{\varepsilon_1,\varepsilon_4}
(\til{\bf F}^{\rm out}_{+,\varepsilon_3})_{\varepsilon_4,\varepsilon_5}
(\wh{\bf M}^{\rm right}(\omega^{2 g(\varepsilon_2)} \wh{Z}_{v_t}) )_{\varepsilon_5, \varepsilon_2}, \\
(\wh{\bf M}_t^{\rm 3in})_{\varepsilon_1,\varepsilon_2,\varepsilon_3}
& = \omega^{\eta^{\rm in}(\varepsilon_1,\varepsilon_2,\varepsilon_3)} \, \underset{\varepsilon_4,\varepsilon_5}{\textstyle \sum}  
( \wh{\bf M}^{\rm right}(\omega^{-2g(\varepsilon_2)-2g(\varepsilon_3)} \wh{Z}_{v_t}) )_{\varepsilon_1,\varepsilon_4}
(\til{\bf F}^{\rm in}_{+,\varepsilon_3})_{\varepsilon_4,\varepsilon_5}
(\wh{\bf M}^{\rm left}_{\rm tran}(\omega^{-2 g(\varepsilon_2)} \wh{Z}_{v_t}) )_{\varepsilon_5, \varepsilon_2},
\end{align*}
where the sums are over all $\varepsilon_4,\varepsilon_5 \in \{1,2,3\}$, and
\begin{align}
\label{eq:g}
g(\varepsilon) = \left\{
\begin{array}{cl}
-1 & \mbox{if $\varepsilon\in\{1,3\}$,} \\
2 & \mbox{if $\varepsilon=2$}.
\end{array}
\right.
\end{align}

\vs

Define the \ul{\em quantum U-turn matrices}
\begin{align}
\label{eq:quantum_U-turn_matrices}
\wh{\bf M}^{\rm U}_- 
= \smallmatthree{0}{0}{q^{-7/3}}{0}{-q^{-4/3}}{0}{q^{-1/3}}{0}{0}, \qquad
\wh{\bf M}^{\rm U}_+ 
= \smallmatthree{0}{0}{q^{7/3}}{0}{-q^{4/3}}{0}{q^{1/3}}{0}{0}.
\end{align}

\begin{enumerate}
\item[\rm \redfix{(QT2-1)}] If $W$ consists of a single left turn corner arc in $t\times {\bf I}$, with its initial point $x$ lying over the side $e_\alpha$ and the terminal point $y$ lying over the side $e_{\alpha+1}$, then
\begin{align*}
{\rm Tr}^\omega_{\Delta;t} ([W,s]) & = ( \wh{\bf M}^{\rm in}_{t,\alpha} \, \wh{\bf M}^{\rm left}(\wh{Z}_{v_t}) \, \wh{\bf M}^{\rm out}_{t,\alpha+1} )_{s(x),s(y)}.
\end{align*}

\item[\rm \redfix{(QT2-2)}] If $W$ consists of a single right turn corner arc in $t\times {\bf I}$, with its initial point $x$ lying over the side $e_{\alpha+1}$ and the terminal point $y$ lying over the side $e_{\alpha}$, then
$$
{\rm Tr}^\omega_{\Delta;t} ([W,s]) = ( \wh{\bf M}^{\rm in}_{t,\alpha+1} \, \wh{\bf M}^{\rm right}(\wh{Z}_{v_t}) \, \wh{\bf M}^{\rm out}_{t,\alpha} )_{s(x),s(y)}.
$$

\item[\rm \redfix{(QT2-3)}] Suppose $W$ consists of a single $3$-way web in $t \times {\bf I}$ with endpoints $x_1,x_2,x_3$, where $x_1,x_2$ lying over a side $e_\alpha$ and $x_3$ over a different side $e_\beta$, where $\pi(x_1)\to \pi(x_2)$ matches the clockwise orientation of $\partial t$ (where $\pi : t\times {\bf I} \to t$ is the projection), and let $\varepsilon_i := s(x_i)$ for each $i=1,2,3$. If $\varepsilon_1=\varepsilon_2$ then ${\rm Tr}^\omega_{\Delta;t} ([W,s])=0$. If $\varepsilon_1\neq \varepsilon_2$, let $\varepsilon$ be the unique element of $\{1,2,3\}$ such that $\{r_1(\varepsilon),r_2(\varepsilon)\} = \{\varepsilon_1,\varepsilon_2\}$. Let $W'$ be an ${\rm SL}_3$-web in $t \times {\bf I}$ consisting of a single corner arc connecting the thickenings of the sides $e_\alpha$ and $e_\beta$, where the endpoint on $e_\beta$ is a sink if and only if $x_3$ is a sink of $W$. Let $s'$ be the state of $W'$ that assigns $\varepsilon$ to the endpoint over $e_\alpha$ and $\varepsilon_3$ to the endpoint over $e_\beta$. Then
$$
{\rm Tr}^\omega_{\Delta;t}([W,s]) = {\bf F}^\omega_{\Delta;t}([W,s]) \, {\rm Tr}^\omega_{\Delta;t}([W',s']),
$$
where ${\rm Tr}^\omega_{\Delta;t}([W',s'])$ is given by \redfix{(QT2-1)} or \redfix{(QT2-2)}, and the value of the factor ${\bf F}^\omega_{\Delta;t}([W,s]) \in \mathbb{Z}[\omega^{\pm 1/2}]$ is defined to be $(\wh{\bf F}^h_{k,\varepsilon})_{\varepsilon_1,\varepsilon_2}$ of eq.\eqref{eq:fork_matrix}, where $h = {\rm out}$ (resp. ${\rm in}$) if $x_1,x_2,x_3$ are sinks (resp. sources) so that $W$ is outgoing (resp. incoming) 3-way web, and $k = +$ (resp. $-$) if $x_1 \succ x_2$ (resp. $x_1 \prec x_2$).

\item[\rm \redfix{(QT2-4)}] If $W$ consists of a single $3$-way web in $t\times {\bf I}$ with endpoints $x_1,x_2,x_3$ lying over the sides $e_1$, $e_2$, $e_3$ respectively, if we let $\varepsilon_\alpha := s(x_\alpha)$ for each $\alpha=1,2,3$, then
$$
{\rm Tr}^\omega_{\Delta;t} ([W,s]) = \left\{ \begin{array}{ll}
\left[ {\textstyle \prod}_{\alpha=1}^3 (\wh{\bf M}^{\rm out}_{t,\alpha})_{\varepsilon_\alpha, \varepsilon_\alpha} \right]_{\rm Weyl} \,
 (\wh{\bf M}_t^{\rm 3out})_{\varepsilon_1,\varepsilon_2,\varepsilon_3} 
& \mbox{if $W$ is an outgoing $3$-way,} \\
\left[ {\textstyle \prod}_{\alpha=1}^3 (\wh{\bf M}^{\rm in}_{t,\alpha})_{\varepsilon_\alpha, \varepsilon_\alpha} \right]_{\rm Weyl} \,
 (\wh{\bf M}_t^{\rm 3in})_{\varepsilon_1,\varepsilon_2,\varepsilon_3} 
& \mbox{if $W$ is an incoming $3$-way.}
\end{array}
\right.
$$

\item[\rm \redfix{(QT2-5)}] If $W$ consists of a single U-turn arc in $t\times{\bf I}$, i.e. $W$ consists of an oriented edge with no crossing whose two endpoints lie over one side $e_\alpha$, then
\begin{align}
\label{eq:Tr_B_U}
& \mbox{\rm case \redfix{(QT2-5-1)}:} \quad
{\rm Tr}^\omega_{\Delta;t}([ \raisebox{-0.4\height}{
\begingroup%
  \makeatletter%
  \providecommand\color[2][]{%
    \errmessage{(Inkscape) Color is used for the text in Inkscape, but the package 'color.sty' is not loaded}%
    \renewcommand\color[2][]{}%
  }%
  \providecommand\transparent[1]{%
    \errmessage{(Inkscape) Transparency is used (non-zero) for the text in Inkscape, but the package 'transparent.sty' is not loaded}%
    \renewcommand\transparent[1]{}%
  }%
  \providecommand\rotatebox[2]{#2}%
  \newcommand*\fsize{\dimexpr\f@size pt\relax}%
  \newcommand*\lineheight[1]{\fontsize{\fsize}{#1\fsize}\selectfont}%
  \ifx\svgwidth\undefined%
    \setlength{\unitlength}{42.51968504bp}%
    \ifx\svgscale\undefined%
      \relax%
    \else%
      \setlength{\unitlength}{\unitlength * \real{\svgscale}}%
    \fi%
  \else%
    \setlength{\unitlength}{\svgwidth}%
  \fi%
  \global\let\svgwidth\undefined%
  \global\let\svgscale\undefined%
  \makeatother%
  \begin{picture}(1,0.83333333)%
    \lineheight{1}%
    \setlength\tabcolsep{0pt}%
    \put(0,0){\includegraphics[width=\unitlength,page=1]{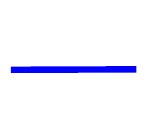}}%
    \put(0.13088999,0.12227517){\color[rgb]{0,0,0}\makebox(0,0)[lt]{\lineheight{1.25}\smash{\begin{tabular}[t]{l}$x_1$\end{tabular}}}}%
    \put(0,0){\includegraphics[width=\unitlength,page=2]{boundary_rel9_quantum.pdf}}%
    \put(0.72355658,0.12227517){\color[rgb]{0,0,0}\makebox(0,0)[lt]{\lineheight{1.25}\smash{\begin{tabular}[t]{l}$x_2$\end{tabular}}}}%
    \put(0.47661228,0.12227517){\color[rgb]{0,0,0}\makebox(0,0)[lt]{\lineheight{1.25}\smash{\begin{tabular}[t]{l}$\prec$\end{tabular}}}}%
  \end{picture}%
\endgroup%
} ])= (\wh{\bf M}^{\rm U}_+)_{s(x_1),s(x_2)}, \\
\label{eq:Tr_B_U2}
& \mbox{\rm case \redfix{(QT2-5-2)}:} \quad
{\rm Tr}^\omega_{\Delta;t}([ \raisebox{-0.4\height}{
\begingroup%
  \makeatletter%
  \providecommand\color[2][]{%
    \errmessage{(Inkscape) Color is used for the text in Inkscape, but the package 'color.sty' is not loaded}%
    \renewcommand\color[2][]{}%
  }%
  \providecommand\transparent[1]{%
    \errmessage{(Inkscape) Transparency is used (non-zero) for the text in Inkscape, but the package 'transparent.sty' is not loaded}%
    \renewcommand\transparent[1]{}%
  }%
  \providecommand\rotatebox[2]{#2}%
  \newcommand*\fsize{\dimexpr\f@size pt\relax}%
  \newcommand*\lineheight[1]{\fontsize{\fsize}{#1\fsize}\selectfont}%
  \ifx\svgwidth\undefined%
    \setlength{\unitlength}{42.51968504bp}%
    \ifx\svgscale\undefined%
      \relax%
    \else%
      \setlength{\unitlength}{\unitlength * \real{\svgscale}}%
    \fi%
  \else%
    \setlength{\unitlength}{\svgwidth}%
  \fi%
  \global\let\svgwidth\undefined%
  \global\let\svgscale\undefined%
  \makeatother%
  \begin{picture}(1,0.83333333)%
    \lineheight{1}%
    \setlength\tabcolsep{0pt}%
    \put(0,0){\includegraphics[width=\unitlength,page=1]{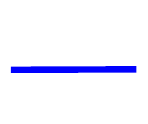}}%
    \put(0.13088999,0.12227517){\color[rgb]{0,0,0}\makebox(0,0)[lt]{\lineheight{1.25}\smash{\begin{tabular}[t]{l}$x_1$\end{tabular}}}}%
    \put(0,0){\includegraphics[width=\unitlength,page=2]{boundary_rel9_2_quantum.pdf}}%
    \put(0.72355658,0.12227517){\color[rgb]{0,0,0}\makebox(0,0)[lt]{\lineheight{1.25}\smash{\begin{tabular}[t]{l}$x_2$\end{tabular}}}}%
    \put(0.47661228,0.12227517){\color[rgb]{0,0,0}\makebox(0,0)[lt]{\lineheight{1.25}\smash{\begin{tabular}[t]{l}$\succ$\end{tabular}}}}%
  \end{picture}%
\endgroup%
} ])= (\wh{\bf M}^{\rm U}_-)_{s(x_1),s(x_2)}, \\
\label{eq:Tr_B_U3}
& \mbox{\rm case \redfix{(QT2-5-3)}:} \quad
{\rm Tr}^\omega_{\Delta;t}([ \raisebox{-0.4\height}{
\begingroup%
  \makeatletter%
  \providecommand\color[2][]{%
    \errmessage{(Inkscape) Color is used for the text in Inkscape, but the package 'color.sty' is not loaded}%
    \renewcommand\color[2][]{}%
  }%
  \providecommand\transparent[1]{%
    \errmessage{(Inkscape) Transparency is used (non-zero) for the text in Inkscape, but the package 'transparent.sty' is not loaded}%
    \renewcommand\transparent[1]{}%
  }%
  \providecommand\rotatebox[2]{#2}%
  \newcommand*\fsize{\dimexpr\f@size pt\relax}%
  \newcommand*\lineheight[1]{\fontsize{\fsize}{#1\fsize}\selectfont}%
  \ifx\svgwidth\undefined%
    \setlength{\unitlength}{42.51968504bp}%
    \ifx\svgscale\undefined%
      \relax%
    \else%
      \setlength{\unitlength}{\unitlength * \real{\svgscale}}%
    \fi%
  \else%
    \setlength{\unitlength}{\svgwidth}%
  \fi%
  \global\let\svgwidth\undefined%
  \global\let\svgscale\undefined%
  \makeatother%
  \begin{picture}(1,0.83333333)%
    \lineheight{1}%
    \setlength\tabcolsep{0pt}%
    \put(0,0){\includegraphics[width=\unitlength,page=1]{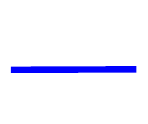}}%
    \put(0.13088999,0.12227517){\color[rgb]{0,0,0}\makebox(0,0)[lt]{\lineheight{1.25}\smash{\begin{tabular}[t]{l}$x_1$\end{tabular}}}}%
    \put(0,0){\includegraphics[width=\unitlength,page=2]{boundary_rel9_3_quantum.pdf}}%
    \put(0.72355658,0.12227517){\color[rgb]{0,0,0}\makebox(0,0)[lt]{\lineheight{1.25}\smash{\begin{tabular}[t]{l}$x_2$\end{tabular}}}}%
    \put(0.47661228,0.12227517){\color[rgb]{0,0,0}\makebox(0,0)[lt]{\lineheight{1.25}\smash{\begin{tabular}[t]{l}$\prec$\end{tabular}}}}%
  \end{picture}%
\endgroup%
}  ])= (\wh{\bf M}^{\rm U}_+)_{s(x_1),s(x_2)}, \\
\label{eq:Tr_B_U4}
& \mbox{\rm case \redfix{(QT2-5-4)}:} \quad
{\rm Tr}^\omega_{\Delta;t}([ \raisebox{-0.4\height}{
\begingroup%
  \makeatletter%
  \providecommand\color[2][]{%
    \errmessage{(Inkscape) Color is used for the text in Inkscape, but the package 'color.sty' is not loaded}%
    \renewcommand\color[2][]{}%
  }%
  \providecommand\transparent[1]{%
    \errmessage{(Inkscape) Transparency is used (non-zero) for the text in Inkscape, but the package 'transparent.sty' is not loaded}%
    \renewcommand\transparent[1]{}%
  }%
  \providecommand\rotatebox[2]{#2}%
  \newcommand*\fsize{\dimexpr\f@size pt\relax}%
  \newcommand*\lineheight[1]{\fontsize{\fsize}{#1\fsize}\selectfont}%
  \ifx\svgwidth\undefined%
    \setlength{\unitlength}{42.51968504bp}%
    \ifx\svgscale\undefined%
      \relax%
    \else%
      \setlength{\unitlength}{\unitlength * \real{\svgscale}}%
    \fi%
  \else%
    \setlength{\unitlength}{\svgwidth}%
  \fi%
  \global\let\svgwidth\undefined%
  \global\let\svgscale\undefined%
  \makeatother%
  \begin{picture}(1,0.83333333)%
    \lineheight{1}%
    \setlength\tabcolsep{0pt}%
    \put(0,0){\includegraphics[width=\unitlength,page=1]{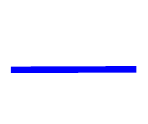}}%
    \put(0.13088999,0.12227517){\color[rgb]{0,0,0}\makebox(0,0)[lt]{\lineheight{1.25}\smash{\begin{tabular}[t]{l}$x_1$\end{tabular}}}}%
    \put(0,0){\includegraphics[width=\unitlength,page=2]{boundary_rel9_4_quantum.pdf}}%
    \put(0.72355658,0.12227517){\color[rgb]{0,0,0}\makebox(0,0)[lt]{\lineheight{1.25}\smash{\begin{tabular}[t]{l}$x_2$\end{tabular}}}}%
    \put(0.47661228,0.12227517){\color[rgb]{0,0,0}\makebox(0,0)[lt]{\lineheight{1.25}\smash{\begin{tabular}[t]{l}$\succ$\end{tabular}}}}%
  \end{picture}%
\endgroup%
}  ])= (\wh{\bf M}^{\rm U}_-)_{s(x_1),s(x_2)},
\end{align}
where the blue line is a boundary arc of $t$, and each diagram is assumed to be carrying a state $s$.

\end{enumerate}
\end{enumerate}
\end{theorem}

\begin{remark}
\label{rem:Douglas_triangle}
The values for the cases \redfix{(QT2-1)}, \redfix{(QT2-2)} and \redfix{(QT2-5)} which do not involve 3-valent vertices coincide with the ones set in \cite{Douglas} \cite{Douglas21}. The values for the remaining cases which involve 3-valent vertices are not dealt with in \cite{Douglas} \cite{Douglas21}.
\end{remark}

\begin{corollary}[the ${\rm SL}_3$ classical trace map; \cite{Kim}]
\label{cor:SL3_classical_trace_map}
Theorem \ref{thm:SL3_quantum_trace_map} holds when $\omega^{1/2}=1$.
\end{corollary}
For the classical case $\omega^{1/2}=1$, we denote
$$
{\rm Tr}_\Delta := {\rm Tr}^1_\Delta : \mathcal{S}^1_{\rm s}(\frak{S};\mathbb{Z}) \to \mathcal{Z}^1_\Delta
$$
and call it the \ul{\em ${\rm SL}_3$ classical trace map}. For convenience, we denote
$$
\mathcal{Z}_\Delta := \mathcal{Z}^1_\Delta, \qquad \mathcal{Z}_t := \mathcal{Z}^1_t,
$$
the classical cube-root Fock-Goncharov algebra and the classical cube-root triangle algebra, which are just (commutative) Laurent polynomial rings. Also, for $\omega^{1/2}=1$, there is no need to consider the 3-manifold and one can formulate solely in terms of $\frak{S}$. 
\begin{definition}
\label{def:commutative_stated_A2-skein_algebra}
Let $\frak{S}$ be a generalized marked surface. A \ul{\em stated ${\rm SL}_3$-web in $\frak{S}$} is a pair $(W,s)$ of an ${\rm SL}_3$-web $W$ in the surface $\frak{S}$ and a state $s : \partial W \to \{1,2,3\}$. Let the \ul{\em (commutative) stated ${\rm SL}_3$-skein algebra} $\mathcal{S}_{\rm s}(\frak{S};\mathbb{Z})$ be defined like the (commutative) ${\rm SL}_3$-skein algebra $\mathcal{S}(\frak{S};\mathbb{Z})$ (Def.\ref{def:A2-skein_algebra}), but modeled on stated ${\rm SL}_3$-webs in $\frak{S}$ instead of \redfix{ordinary} ${\rm SL}_3$-webs in $\frak{S}$.
\end{definition}
\begin{lemma}
\label{lem:surface_stated_skein_algebra_and_commutative_stated_skein_algebra}
There is a natural isomorphism
$$
\mathcal{S}_{\rm s}(\frak{S};\mathbb{Z}) \to \mathcal{S}_{\rm s}^1(\frak{S};\mathbb{Z}).
$$
\end{lemma}
Hence the ${\rm SL}_3$ classical trace map can be viewed as
\begin{align}
\label{eq:Tr_Delta}
{\rm Tr}_\Delta = {\rm Tr}_{\Delta;\frak{S}} : \mathcal{S}_{\rm s}(\frak{S};\mathbb{Z}) \to \mathcal{Z}_\Delta,
\end{align}
simplifying the situations a lot, because one can work only on the 2-dimensional surfaces. Indeed, in \S5 of the first version of the present paper \cite{Kim}, only Cor.\ref{cor:SL3_classical_trace_map} is proved, instead of Thm.\ref{thm:SL3_quantum_trace_map}. So, a reader who is interested only in the classical story can consult \cite{Kim} for a proof of Cor.\ref{cor:SL3_classical_trace_map} (which is Prop.5.6 there), which is much simpler than the proof of Thm.\ref{thm:SL3_quantum_trace_map}. In the present version, we establish a proof of Thm.\ref{thm:SL3_quantum_trace_map} throughout the present section. Then we will just put $\omega^{1/2} = 1$ to obtain and use the classical statements.

\vs

The normalization of the quantum left and right turn matrices and their transpose versions is justified by the following observation, which can be easily verified by straightforward computation:
\begin{lemma}[the quantum trace for a left or right turn in a triangle is Weyl-ordered]
\label{lem:quantum_turn_matrices_are_Weyl-ordered}
Each of the following four matrices (after computing the product) is Weyl-ordered, i.e. equals the Weyl-ordered version of it:
\begin{align}
\nonumber
& \wh{\bf M}^{\rm in}_{t,\alpha} \, \wh{\bf M}^{\rm left}(\wh{Z}_{v_t}) \, \wh{\bf M}^{\rm out}_{t,\alpha+1}, \qquad 
\wh{\bf M}^{\rm in}_{t,\alpha+1} \, \wh{\bf M}^{\rm right}(\wh{Z}_{v_t}) \, \wh{\bf M}^{\rm out}_{t,\alpha}, \\
\nonumber
& \wh{\bf M}^{\rm out}_{t,\alpha+1} \, \wh{\bf M}^{\rm left}_{\rm tran}(\wh{Z}_{v_t}) \, \wh{\bf M}^{\rm in}_{t,\alpha},\qquad \wh{\bf M}^{\rm out}_{t,\alpha} \, \wh{\bf M}^{\rm right}_{\rm tran}(\wh{Z}_{v_t}) \, \wh{\bf M}^{\rm in}_{t,\alpha+1}.
\end{align}
\end{lemma}
{\it Proof.} After writing down each entries, one can use the identities
\begin{align*}
\hspace{-10mm}
[\wh{Z}_{t,v_{e_\alpha,1}}^{a_1} \wh{Z}_{t,v_{e_\alpha,2}}^{a_2}]_{\rm Weyl} \wh{Z}_{v_t}^{a_3} [\wh{Z}_{t,v_{e_{\alpha+1},1}}^{a_4} \wh{Z}_{t,v_{e_{\alpha+1},2}}^{a_5}]_{\rm Weyl}
& = \omega^{(-a_1+a_2+a_4-a_5) a_3 - a_2 a_4} [ \wh{Z}_{t,v_{e_\alpha,1}}^{a_1} \wh{Z}_{t,v_{e_\alpha,2}}^{a_2} \wh{Z}_{v_t}^{a_3} \wh{Z}_{t,v_{e_{\alpha+1},1}}^{a_4} \wh{Z}_{t,v_{e_{\alpha+1},2}}^{a_5} ]_{\rm Weyl}, \\
\hspace{-10mm}
[\wh{Z}_{t,v_{e_{\alpha+1},1}}^{b_1} \wh{Z}_{t,v_{e_{\alpha+1},2}}^{b_2}]_{\rm Weyl} \wh{Z}_{v_t}^{b_3} [\wh{Z}_{t,v_{e_{\alpha},1}}^{b_4} \wh{Z}_{t,v_{e_{\alpha},2}}^{b_5}]_{\rm Weyl}
& = \omega^{(-b_1+b_2+b_4-b_5) b_3 + b_1 b_5} [ \wh{Z}_{t,v_{e_{\alpha+1},1}}^{b_1} \wh{Z}_{t,v_{e_{\alpha+1},2}}^{b_2} \wh{Z}_{v_t}^{b_3} \wh{Z}_{t,v_{e_{\alpha},1}}^{b_4} \wh{Z}_{t,v_{e_{\alpha},2}}^{b_5} ]_{\rm Weyl};
\end{align*}
then actual checking is easy, and left as exercise. \qed

\vs

A useful immediate consequence which will be convenient later is:
\begin{corollary}
\label{cor:transpose_of_quantum_turn_matrices}
One has
\begin{align*}
& (\wh{\bf M}^{\rm in}_{t,\alpha} \, \wh{\bf M}^{\rm left}(\wh{Z}_{v_t}) \, \wh{\bf M}^{\rm out}_{t,\alpha+1})^{\rm tr} = \wh{\bf M}^{\rm out}_{t,\alpha+1} \, \wh{\bf M}^{\rm left}_{\rm tran}(\wh{Z}_{v_t}) \, \wh{\bf M}^{\rm in}_{t,\alpha}, \\
& (\wh{\bf M}^{\rm in}_{t,\alpha+1} \, \wh{\bf M}^{\rm right}(\wh{Z}_{v_t}) \, \wh{\bf M}^{\rm out}_{t,\alpha})^{\rm tr} = \wh{\bf M}^{\rm out}_{t,\alpha} \, \wh{\bf M}^{\rm right}_{\rm tran}(\wh{Z}_{v_t}) \, \wh{\bf M}^{\rm in}_{t,\alpha+1},
\end{align*}
where ${\rm tr}$ denotes the transpose. \qed
\end{corollary}

The 3-way matrices $\wh{\bf M}^{\rm 3out}_t$ and $\wh{\bf M}^{\rm 3in}_t$ and hence also the triangle factor involving $3$-way webs in \redfix{(QT2-4)} may look ad hoc, and a priori seem to be depending on the labeling of sides of $t$, i.e. on the choice of which side of $t$ to call $e_3$. The following lemma shows the independence on such a choice, establishing the well-definedness of \redfix{(QT2-4)}. During the course, we compute all entries of $\wh{\bf M}^{\rm 3out}_t$ and $\wh{\bf M}^{\rm 3in}_t$.
\begin{lemma}[cyclic symmetry of $3$-way matrices]
\label{lem:cyclicity_of_3-way_matrices}
$\wh{\bf M}^{\rm 3out}_t$ and $\wh{\bf M}^{\rm 3in}_t$ have cyclic symmetries, i.e.
\begin{align*}
& (\wh{\bf M}^{\rm 3out}_t)_{\varepsilon_1,\varepsilon_2,\varepsilon_3}
= (\wh{\bf M}^{\rm 3out}_t)_{\varepsilon_2,\varepsilon_3,\varepsilon_1}
= (\wh{\bf M}^{\rm 3out}_t)_{\varepsilon_3,\varepsilon_1,\varepsilon_2},\\
& (\wh{\bf M}^{\rm 3in}_t)_{\varepsilon_1,\varepsilon_2,\varepsilon_3}
= (\wh{\bf M}^{\rm 3in}_t)_{\varepsilon_2,\varepsilon_3,\varepsilon_1}
= (\wh{\bf M}^{\rm 3in}_t)_{\varepsilon_3,\varepsilon_1,\varepsilon_2}.
\end{align*}
\end{lemma}
To ease the proof of Lem.\ref{lem:cyclicity_of_3-way_matrices}, we establish three small technical lemmas:
\begin{lemma}
\label{lem:triangle_variable_and_edge_matrices}
For each $\alpha \in \{1,2,3\}$ and $\varepsilon\in\{1,2,3\}$, one has
$$
\wh{Z}_{v_t} (\wh{\bf M}^{\rm out}_{t,\alpha})_{\varepsilon,\varepsilon}
= \omega^{2g(\varepsilon)} (\wh{\bf M}^{\rm out}_{t,\alpha})_{\varepsilon,\varepsilon} \wh{Z}_{v_t}, \qquad
\wh{Z}_{v_t} (\wh{\bf M}^{\rm in}_{t,\alpha})_{\varepsilon,\varepsilon}
= \omega^{-2g(\varepsilon)} (\wh{\bf M}^{\rm in}_{t,\alpha})_{\varepsilon,\varepsilon} \wh{Z}_{v_t},
$$
where $g(\varepsilon)$ is as defined in eq.\eqref{eq:g}, i.e. $g(1)=g(3)=-1$, $g(2)=2$.
\end{lemma}

{\it Proof of Lem.\ref{lem:triangle_variable_and_edge_matrices}.} For convenience, denote $\wh{Z}_{\alpha,i} = \wh{Z}_{t,v_{e_\alpha,i}}$, $i=1,2$. In view of eq.\eqref{eq:quantum_in_and_out}, note $(\wh{\bf M}^{\rm out}_{t,\alpha})_{1,1} = [\wh{Z}_{\alpha,1} \wh{Z}_{\alpha,2}^2]_{\rm Weyl}$, $(\wh{\bf M}^{\rm out}_{t,\alpha})_{2,2} = [\wh{Z}_{\alpha,1} \wh{Z}_{\alpha,2}^{-1}]_{\rm Weyl}$, $(\wh{\bf M}^{\rm out}_{t,\alpha})_{3,3} = [\wh{Z}_{\alpha,1}^{-2} \wh{Z}_{\alpha,2}^{-1}]_{\rm Weyl}$, while $(\wh{\bf M}^{\rm in}_{t,\alpha})_{1,1} = [\wh{Z}_{\alpha,2} \wh{Z}_{\alpha,1}^2]_{\rm Weyl}$, $(\wh{\bf M}^{\rm in}_{t,\alpha})_{2,2} = [\wh{Z}_{\alpha,2} \wh{Z}_{\alpha,1}^{-1}]_{\rm Weyl}$, $(\wh{\bf M}^{\rm in}_{t,\alpha})_{3,3} = [\wh{Z}_{\alpha,2}^{-2} \wh{Z}_{\alpha,1}^{-1}]_{\rm Weyl}$. Using $\wh{Z}_{v_t} \wh{Z}_{\alpha,1} = \omega^2 \wh{Z}_{\alpha,1} \wh{Z}_{v_t}$ and $\wh{Z}_{v_t} \wh{Z}_{\alpha,2} = \omega^{-2} \wh{Z}_{\alpha,2} \wh{Z}_{v_t}$, one can easily verify the the desired result. \qed

\vs

\begin{lemma}
\label{lem:eta_h}
The values of $\eta^h(\varepsilon_1,\varepsilon_2,\varepsilon_3)$ defined in eq.\eqref{eq:eta_h} are given as follows. For each $h\in \{{\rm in},{\rm out}\}$ and $\varepsilon_{\redfix{3}} \in \{1,2,3\}$, the matrix $\eta^h_{\varepsilon_3}$ whose $(\varepsilon_1,\varepsilon_2)$-th entry is $(\eta^h_{\varepsilon_3})_{\varepsilon_1,\varepsilon_2} = \eta^h(\varepsilon_1,\varepsilon_2,\varepsilon_3)$ is:
\begin{align*}
& \eta^{\rm out}_1 = 
\smallmatthree{2}{-1}{5}{5}{2}{-1}{-1}{-4}{-7}, \qquad
\eta^{\rm out}_2 = 
\smallmatthree{-1}{-4}{2}{2}{-1}{-4}{5}{2}{-1}, \qquad
\eta^{\rm out}_3 = 
\smallmatthree{-7}{-1}{5}{-4}{2}{-1}{-1}{5}{2} \\
& \eta^{\rm in}_1 = 
\smallmatthree{2}{5}{-1}{-1}{2}{-4}{5}{-1}{-7}, \qquad
\eta^{\rm in}_2 = 
\smallmatthree{-1}{2}{5}{-4}{-1}{2}{2}{-4}{-1}, \qquad
\eta^{\rm in}_3 = 
\smallmatthree{-7}{-4}{-1}{-1}{2}{5}{5}{-1}{2}. 
\end{align*}
\end{lemma}

{\it Proof of Lem.\ref{lem:eta_h}.} Straightforward check, e.g. using
\begin{align*}
\hspace{-10mm}
& [\wh{Z}_{t,v_{e_1,1}}^{a_1} \wh{Z}_{t,v_{e_1,2}}^{a_2}]_{\rm Weyl}
[\wh{Z}_{t,v_{e_3,1}}^{a_3} \wh{Z}_{t,v_{e_3,2}}^{a_4}]_{\rm Weyl}
[\wh{Z}_{t,v_{e_2,1}}^{a_5} \wh{Z}_{t,v_{e_2,2}}^{a_6}]_{\rm Weyl} \\
& = \omega^{a_1a_4 -a_2a_5 +a_3a_6} [ \wh{Z}_{t,v_{e_1,1}}^{a_1} \wh{Z}_{t,v_{e_1,2}}^{a_2}
\wh{Z}_{t,v_{e_3,1}}^{a_3} \wh{Z}_{t,v_{e_3,2}}^{a_4}
\wh{Z}_{t,v_{e_2,1}}^{a_5} \wh{Z}_{t,v_{e_2,2}}^{a_6}
 ]_{\rm Weyl}. \qed
\end{align*}

\begin{lemma}
\label{lem:twisted_fork_matrix_entries}
If we write the entries of the twisted (positive) fork matrices $\til{\bf F}^h_{+,\varepsilon}$ in eq.\eqref{eq:fork_matrix_tilde} as
\begin{align}
\nonumber
(\til{\bf F}^{h}_{+,\varepsilon} )_{i,j} = \left\{
\begin{array}{ll}
\til{c}^h_{+,\varepsilon} & \mbox{if $(i,j)=(r_1(\varepsilon),r_2(\varepsilon))$}, \\
\til{d}^h_{+,\varepsilon} & \mbox{if $(i,j)=(r_2(\varepsilon),r_1(\varepsilon))$}, \\
0 &\mbox{otherwise}.
\end{array}
\right.
\end{align}
then

\begin{align*}
\begin{array}{llll}
\til{c}^{\rm out}_{+,1} = \til{c}^{\rm out}_{+,3} = \omega^3, & \til{c}^{\rm out}_{+,2} = 1, & \til{c}^{\rm in}_{+,1} = \til{c}^{\rm in}_{+,3} = 1, & \til{c}^{\rm in}_{+,2} = \omega^3, \\
\til{d}^{\rm out}_{+,1} = \til{d}^{\rm out}_{+,3} = -  \omega^9, & \til{d}^{\rm out}_{+,2} = -  \omega^{12}, & \til{d}^{\rm in}_{+,1} = \til{d}^{\rm in}_{+,3} = - \omega^{12}, & \til{d}^{\rm in}_{+,2} = - \omega^9. \qed
\end{array}
\end{align*}

\end{lemma}

\vs

{\it Proof of Lem.\ref{lem:cyclicity_of_3-way_matrices}.} We first unravel the definitions of $\wh{\bf M}^{\rm 3out}_t$ and $\wh{\bf M}^{\rm 3in}_t$ a little bit. For convenience, define
\begin{align}
\label{eq:til_bf_M_3out_3in}
\til{\bf M}_t^{\rm 3out} := \omega^{-\eta^{\rm out}(\varepsilon_1,\varepsilon_2,\varepsilon_3)} \wh{\bf M}_t^{\rm 3out}, \qquad
\til{\bf M}_t^{\rm 3in} := \omega^{-\eta^{\rm in}(\varepsilon_1,\varepsilon_2,\varepsilon_3)} \wh{\bf M}_t^{\rm 3in},
\end{align}
so that by eq.\eqref{eq:eta_h}, the value ${\rm Tr}^\omega_{\Delta;t}([W,s])$ for \redfix{(QT2-4)} is
\begin{align*}
{\rm Tr}^\omega_{\Delta;t} ([W,s]_{\mbox{\tiny \redfix{(QT2-4)}}}) := \left\{ \begin{array}{ll}
(\wh{\bf M}^{\rm out}_{t,1})_{\varepsilon_1, \varepsilon_1} 
(\wh{\bf M}^{\rm out}_{t,3})_{\varepsilon_3, \varepsilon_3} 
(\wh{\bf M}^{\rm out}_{t,2})_{\varepsilon_2, \varepsilon_2} 
(\til{\bf M}_t^{\rm 3out})_{\varepsilon_1,\varepsilon_2,\varepsilon_3} 
& \mbox{if $W$ is an outgoing $3$-way,} \\
(\wh{\bf M}^{\rm in}_{t,1})_{\varepsilon_1, \varepsilon_1} 
(\wh{\bf M}^{\rm in}_{t,3})_{\varepsilon_3, \varepsilon_3} 
(\wh{\bf M}^{\rm in}_{t,2})_{\varepsilon_2, \varepsilon_2} 
(\til{\bf M}_t^{\rm 3in})_{\varepsilon_1,\varepsilon_2,\varepsilon_3}  
& \mbox{if $W$ is an incoming $3$-way.}
\end{array}
\right.
\end{align*}
Note that
\begin{align}
\nonumber
& (\wh{\bf M}^{\rm out}_{t,1})_{\varepsilon_1, \varepsilon_1} 
(\wh{\bf M}^{\rm out}_{t,3})_{\varepsilon_3, \varepsilon_3} 
(\wh{\bf M}^{\rm out}_{t,2})_{\varepsilon_2, \varepsilon_2} 
(\til{\bf M}_t^{\rm 3out})_{\varepsilon_1,\varepsilon_2,\varepsilon_3}  \\
\nonumber
& = (\wh{\bf M}^{\rm out}_{t,1})_{\varepsilon_1, \varepsilon_1} 
(\wh{\bf M}^{\rm out}_{t,3})_{\varepsilon_3, \varepsilon_3} 
\underbrace{ (\wh{\bf M}^{\rm out}_{t,2})_{\varepsilon_2, \varepsilon_2} }_{\mbox{\tiny move to right}} 
\underset{\varepsilon_4,\varepsilon_5}{\textstyle \sum} ( \wh{\bf M}^{\rm left}_{\rm tran}(\omega^{2g(\varepsilon_2)+2g(\varepsilon_3)} \wh{Z}_{v_t} ))_{\varepsilon_1,\varepsilon_4}
(\til{\bf F}^{\rm out}_{+,\varepsilon_3})_{\varepsilon_4,\varepsilon_5}
(\wh{\bf M}^{\rm right}(\omega^{2 g(\varepsilon_2)} \wh{Z}_{v_t}) )_{\varepsilon_5, \varepsilon_2}, \\
\nonumber
& \hspace*{-3mm} \overset{{\rm Lem}.\ref{lem:triangle_variable_and_edge_matrices}}{=} (\wh{\bf M}^{\rm out}_{t,1})_{\varepsilon_1, \varepsilon_1} 
(\wh{\bf M}^{\rm out}_{t,3})_{\varepsilon_3, \varepsilon_3} 
\left( \underset{\varepsilon_4,\varepsilon_5}{\textstyle \sum} ( \wh{\bf M}^{\rm left}_{\rm tran}(\omega^{2g(\varepsilon_3)} \wh{Z}_{v_t} ))_{\varepsilon_1,\varepsilon_4}
(\til{\bf F}^{\rm out}_{+,\varepsilon_3})_{\varepsilon_4,\varepsilon_5}
(\wh{\bf M}^{\rm right}(\wh{Z}_{v_t}) )_{\varepsilon_5, \varepsilon_2} \right)
({\bf M}^{\rm out}_{t,2})_{\varepsilon_2, \varepsilon_2} \\
\nonumber
& = (\wh{\bf M}^{\rm out}_{t,1})_{\varepsilon_1, \varepsilon_1} 
(\wh{\bf M}^{\rm out}_{t,3})_{\varepsilon_3, \varepsilon_3} 
( \underline{ \wh{\bf M}^{\rm left}_{\rm tran}(\omega^{2g(\varepsilon_3)} \wh{Z}_{v_t} )
\til{\bf F}^{\rm out}_{+,\varepsilon_3}
\wh{\bf M}^{\rm right}(\wh{Z}_{v_t}) } )_{\varepsilon_1, \varepsilon_2}
(\wh{\bf M}^{\rm out}_{t,2})_{\varepsilon_2, \varepsilon_2}
=: (*)^{\rm out}_{\varepsilon_1,\varepsilon_2,\varepsilon_3},
\end{align}
and by a similar computation using Lem.\ref{lem:triangle_variable_and_edge_matrices} we get
\begin{align*}
& (\wh{\bf M}^{\rm in}_{t,1})_{\varepsilon_1, \varepsilon_1} 
(\wh{\bf M}^{\rm in}_{t,3})_{\varepsilon_3, \varepsilon_3} 
(\wh{\bf M}^{\rm in}_{t,2})_{\varepsilon_2, \varepsilon_2} 
(\til{\bf M}_t^{\rm 3in})_{\varepsilon_1,\varepsilon_2,\varepsilon_3}   \\
& = (\wh{\bf M}^{\rm in}_{t,1})_{\varepsilon_1, \varepsilon_1} 
(\wh{\bf M}^{\rm in}_{t,3})_{\varepsilon_3, \varepsilon_3} 
( \underline{ \wh{\bf M}^{\rm right}(\omega^{-2g(\varepsilon_3)} \wh{Z}_{v_t} )
\til{\bf F}^{\rm in}_{+,\varepsilon_3}
\wh{\bf M}^{\rm left}_{\rm tran}(\wh{Z}_{v_t}) } )_{\varepsilon_1, \varepsilon_2}
(\wh{\bf M}^{\rm in}_{t,2})_{\varepsilon_2, \varepsilon_2} =: (*)^{\rm in}_{\varepsilon_1,\varepsilon_2,\varepsilon_3}.
\end{align*}   
Now we compute the underlined product of matrices, for each $\varepsilon_3\in\{1,2,3\}$. We just list the results here, as it is straightforward to verify.
\begin{align*}
& {\renewcommand{\arraystretch}{1.3} \begin{array}{l}
\wh{\bf M}^{\rm left}_{\rm tran}(\omega^{2 g(1)}\wh{Z}_{v_t}) \til{\bf F}^{\rm out}_{+,1} \wh{\bf M}^{\rm right}(\wh{Z}_{v_t})  \\
= \smallmatthree{\omega^{-5} (\omega^{-2} \wh{Z}_{v_t})^2}{0}{0}{\omega (\omega^{-2} \wh{Z}_{v_t})^2 + \omega^{-2} (\omega^{-2} \wh{Z}_{v_t})^{-1} ~}{\omega^{-5} (\omega^{-2} \wh{Z}_{v_t})^{-1}}{0}{\omega^{4} (\omega^{-2} \wh{Z}_{v_t})^{-1}}{\omega (\omega^{-2} \wh{Z}_{v_t})^{-1}}{\omega^{-2} (\omega^{-2} \wh{Z}_{v_t})^{-1}} 
\smallmatthree{0}{\omega^3}{0}{-\omega^9}{0}{0}{0}{0}{0} 
\smallmatthree{\omega^{-2} \wh{Z}_{v_t}}{0}{0}{\omega \wh{Z}_{v_t}}{\omega^{-5} \wh{Z}_{v_t}}{0}{\omega^4 \wh{Z}_{v_t}}{~\omega^{-2} \wh{Z}_{v_t} + \omega \wh{Z}_{v_t}^{-2}~}{\omega^{-5} \wh{Z}_{v_t}^{-2}} \\
= \smallmatthree{\omega^{-5} \wh{Z}_{v_t}^3}{\omega^{-11} \wh{Z}_{v_t}^3}{0}{\omega \wh{Z}_{v_t}^3}{~(\omega^{-5} \wh{Z}_{v_t}^3 + \omega^{-2})~}{0}{0}{\omega^4}{0}  \\
\end{array}} 
\\
& {\renewcommand{\arraystretch}{1.3} \begin{array}{ll}
\wh{\bf M}^{\rm left}_{\rm tran}(\omega^{2g(3)}\wh{Z}_{v_t}) \til{\bf F}^{\rm out}_{+,3} \wh{\bf M}^{\rm right}(\wh{Z}_{v_t})
&
\quad \wh{\bf M}^{\rm left}_{\rm tran}(\omega^{2g(2)}\wh{Z}_{v_t}) \til{\bf F}^{\rm out}_{+,2} \wh{\bf M}^{\rm right}(\wh{Z}_{v_t}) \\
= \smallmatthree{0}{0}{0}{\omega^4}{~\omega^{-2} + \omega^{1} \wh{Z}_{v_t}^{-3}~}{\omega^{-5} \wh{Z}_{v_t}^{-3} }{0}{\omega^7 \wh{Z}_{v_t}^{-3}}{\omega\wh{Z}_{v_t}^{-3}}
& 
\quad = \smallmatthree{\omega^7 \wh{Z}_{v_t}^3 }{~(\omega \wh{Z}_{v_t}^3 + \omega^4)~}{\omega^{-2}}{(\omega^{13} \wh{Z}_{v_t}^3 + \omega^{-2})}{~( \omega^7\wh{Z}_{v_t}^3 + \omega^{10} + \omega^{-5} \wh{Z}_{v_t}^{-3} + \omega^{-8} )~}{\omega^4 + \omega^{-11} \wh{Z}_{v_t}^{-3}}{0}{\omega^{-2} + \omega \wh{Z}_{v_t}^{-3}}{\omega^{-5} \wh{Z}_{v_t}^{-3}}
\end{array}} \\
& \hspace*{-2mm} {\renewcommand{\arraystretch}{1.3} \begin{array}{l}
\wh{\bf M}^{\rm right}(\omega^{-2 g(1)}\wh{Z}_{v_t}) \til{\bf F}^{\rm in}_{+,1} \wh{\bf M}^{\rm left}_{\rm tran}(\wh{Z}_{v_t})  \\
= \smallmatthree{(\omega \wh{Z}_{v_t}^3 + \omega^{-2})}{\omega^{-5}}{0}{ \omega }{ \omega^{-2}}{0}{0}{\omega }{0}
\end{array}} 
\hspace*{0mm}
{\renewcommand{\arraystretch}{1.3} \begin{array}{l}
\wh{\bf M}^{\rm right}(\omega^{-2 g(3)}\wh{Z}_{v_t}) \til{\bf F}^{\rm in}_{+,3} \wh{\bf M}^{\rm left}_{\rm tran}(\wh{Z}_{v_t})  \\
= \smallmatthree{0}{0}{0}{ \omega}{~\omega^{-2}~}{ \omega^{-5}}{0}{\omega }{ \omega^{-2} + \omega^{-5} \wh{Z}_{v_t}^{-3}}
\end{array}} 
{\renewcommand{\arraystretch}{1.3} \begin{array}{l}
\wh{\bf M}^{\rm right}(\omega^{-2 g(2)}\wh{Z}_{v_t}) \til{\bf F}^{\rm in}_{+,2} \wh{\bf M}^{\rm left}_{\rm tran}(\wh{Z}_{v_t})  \\
= \smallmatthree{\omega}{\omega^{-2}}{\omega^{-5}}{\omega^4 }{\omega}{\omega^{-2}}{0}{\omega^4}{\omega}
\end{array}}
\end{align*}

Now, in the expression for $(*)^{\rm out}_{\varepsilon_1,\varepsilon_2,\varepsilon_3}$ for each fixed $\varepsilon_3 \in \{1,2,3\}$, we will move the factor $(\wh{\bf M}^{\rm out}_{t,2})_{\varepsilon_2, \varepsilon_2}$ to the left of the underlined part. In view of Lem.\ref{lem:triangle_variable_and_edge_matrices}, for each $\varepsilon_2 \in \{1,2,3\}$ one should replace each $\wh{Z}_{v_t}$ appearing in the $\varepsilon_2$-th column of the matrix $\wh{\bf M}^{\rm left}_{\rm tran}(\omega^{2 g(\varepsilon_3)}\wh{Z}_{v_t}) \til{\bf F}^{\rm out}_{+,\varepsilon_3} \wh{\bf M}^{\rm right}(\wh{Z}_{v_t})$ by $\omega^{2g(\varepsilon_2)} \wh{Z}_{v_t}$. Then one obtains the matrix $( (\til{\bf M}_t^{\rm 3out})_{\varepsilon_1,\varepsilon_2,\varepsilon_3} )_{\varepsilon_1,\varepsilon_2\redfix{\in\{1,2,3\}}}$ for each fixed $\varepsilon_3$. In view of eq.\eqref{eq:til_bf_M_3out_3in}, multiplying these entries by $\omega^{\eta^{\rm out}(\varepsilon_1,\varepsilon_2,\varepsilon_3)}$ (whose values are given in Lem.\ref{lem:eta_h}) yields the sought-for matrix $( (\wh{\bf M}_t^{\rm 3out})_{\varepsilon_1,\varepsilon_2,\varepsilon_3} )_{\varepsilon_1,\varepsilon_2 \redfix{\in \{1,2,3\}}}$. We list the results:
\begin{align*}
& ( (\wh{\bf M}_t^{\rm 3out})_{\varepsilon_1,\varepsilon_2,1} )_{\varepsilon_1,\varepsilon_2} 
= \smallmatthree{\omega^2 \omega^{-5} (\omega^{-2}\wh{Z}_{v_t})^3}{\omega^{-1} \omega^{-11} (\omega^{4}\wh{Z}_{v_t})^3}{0}{\omega^5 \omega (\omega^{-2} \wh{Z}_{v_t})^3}{~\omega^2(\omega^{-5} (\omega^{4}\wh{Z}_{v_t})^3 + \omega^{-2})~}{0}{0}{\omega^{-4} \omega^4}{0} 
= \smallmatthree{\omega^{-9} \wh{Z}_{v_t}^3}{\wh{Z}_{v_t}^3}{0}{\wh{Z}_{v_t}^3}{~\omega^{9} \wh{Z}_{v_t}^3 + 1~}{0}{0}{1}{0}  \\
& 
{\renewcommand{\arraystretch}{1.4} \begin{array}{l}
( (\wh{\bf M}_t^{\rm 3out})_{\varepsilon_1,\varepsilon_2,3} )_{\varepsilon_1,\varepsilon_2}  \\
= \smallmatthree{0}{0}{0}{1}{~1 + \omega^{-9} \wh{Z}_{v_t}^{-3}~}{ \wh{Z}_{v_t}^{-3} }{0}{ \wh{Z}_{v_t}^{-3}}{\omega^9 Z^{-3}}
\end{array}}
\qquad\qquad
{\renewcommand{\arraystretch}{1.4} \begin{array}{l}
( (\wh{\bf M}_t^{\rm 3out})_{\varepsilon_1,\varepsilon_2,2} )_{\varepsilon_1,\varepsilon_2} \\
= \smallmatthree{\wh{Z}_{v_t}^3 }{~(\omega^{9} \wh{Z}_{v_t}^3 + 1)~}{1}{(\omega^{9} \wh{Z}_{v_t}^3 + 1)}{~(\omega^{18} Z^3+\omega^9 + \omega^{-9} + \omega^{-18} Z^{-3})~}{ 1 + \omega^{-9} \wh{Z}_{v_t}^{-3}}{0}{1 + \omega^{-9} \wh{Z}_{v_t}^{-3}}{\wh{Z}_{v_t}^{-3}}
\end{array}}
\end{align*}
By inspection on all 27 values of $(\wh{\bf M}^{\rm 3out}_t)_{\varepsilon_1,\varepsilon_2,\varepsilon_3}$, one indeed observes the cyclicity. 

\vs

Likewise, in the expression for $(*)^{\rm in}_{\varepsilon_1,\varepsilon_2,\varepsilon_3}$ for each fixed $\varepsilon_3 \in \{1,2,3\}$, we will move the factor $(\wh{\bf M}^{\rm in}_{t,2})_{\varepsilon_2, \varepsilon_2}$ to the left of the underlined part. In view of Lem.\ref{lem:triangle_variable_and_edge_matrices}, for each $\varepsilon_2 \in \{1,2,3\}$ one should replace each $\wh{Z}_{v_t}$ appearing in the $\varepsilon_2$-th column of the matrix $\wh{\bf M}^{\rm right}(\omega^{-2g(\varepsilon_3)} \wh{Z}_{v_t} )
\til{\bf F}^{\rm in}_{+,\varepsilon_3}
\wh{\bf M}^{\rm left}_{\rm tran}(\wh{Z}_{v_t})$ by $\omega^{-2g(\varepsilon_2)} \wh{Z}_{v_t}$. Then one obtains the matrix $( (\til{\bf M}_t^{\rm 3in})_{\varepsilon_1,\varepsilon_2,\varepsilon_3} )_{\varepsilon_1,\varepsilon_2\redfix{\in \{1,2,3\}}}$ for each fixed $\varepsilon_3$. In view of eq.\eqref{eq:til_bf_M_3out_3in}, multiplying these entries by $\omega^{\eta^{\rm in}(\varepsilon_1,\varepsilon_2,\varepsilon_3)}$ (whose values are given in Lem.\ref{lem:eta_h}) yields the sought-for matrix $( (\wh{\bf M}_t^{\rm 3in})_{\varepsilon_1,\varepsilon_2,\varepsilon_3} )_{\varepsilon_1,\varepsilon_2\redfix{\in\{1,2,3\}}}$. We list the results:
\begin{align*}
{\renewcommand{\arraystretch}{1.4} \begin{array}{l}
( (\wh{\bf M}_t^{\rm 3in})_{\varepsilon_1,\varepsilon_2,1} )_{\varepsilon_1,\varepsilon_2}  \\
= \smallmatthree{(\omega^9 \wh{Z}_{v_t}^3 + 1)}{1}{0}{1 }{1 }{0}{0}{1}{0}
\end{array}}
\qquad
{\renewcommand{\arraystretch}{1.4} \begin{array}{l}
( (\wh{\bf M}_t^{\rm 3in})_{\varepsilon_1,\varepsilon_2,2} )_{\varepsilon_1,\varepsilon_2}  \\
= \smallmatthree{1}{1}{1}{1}{1}{1}{0}{1}{1}
\end{array}}
\qquad
{\renewcommand{\arraystretch}{1.4} \begin{array}{l}
( (\wh{\bf M}_t^{\rm 3in})_{\varepsilon_1,\varepsilon_2,3} )_{\varepsilon_1,\varepsilon_2} \\
= \smallmatthree{0}{0}{0}{1 }{~ 1 ~}{1 }{0}{1 }{(1  +  \omega^{-9} \wh{Z}_{v_t}^{-3})}
\end{array}}
\end{align*}
By inspection, one observes the cyclicity of $(\wh{\bf M}^{\rm 3in}_t)_{\varepsilon_1,\varepsilon_2,\varepsilon_3}$. \qed  \qquad {\it End of proof of Lem.\ref{lem:cyclicity_of_3-way_matrices}.}

\vs

Two nice consequences of the above computation are:
\begin{lemma}[positivity of 3-way matrices]
\label{lem:positivity_of_3-way_matrices}
Entries of $\wh{\bf M}^{\rm 3out}_t$, $\wh{\bf M}^{\rm 3in}_t$ are Laurent polynomials in $\wh{Z}_{v_t}^3 = \wh{X}_{v_t}$ with coefficients in $\mathbb{Z}_{\ge 0}[q^{\pm1}]$ (where $q= \omega^9$). \qed
\end{lemma}

\begin{lemma}[quantum trace for a 3-way ${\rm SL}_3$-web connecting all sides of a triangle is Weyl-ordered]
\label{lem:standard_3-way_web_values_are_Weyl-ordered}
In Thm.\ref{thm:SL3_quantum_trace_map}\redfix{(QT2-5)}, ${\rm Tr}^\omega_{\Delta;t}([W,s])$ is a Weyl-ordered element of $\mathcal{Z}^\omega_t$.
\end{lemma}
{\it Proof.} Using
\begin{align*}
& \omega^a \cdot [\wh{Z}_{t,v_{e_1,1}}^{a_1} \wh{Z}_{t,v_{e_1,2}}^{a_2}
\wh{Z}_{t,v_{e_2,1}}^{a_3} \wh{Z}_{t,v_{e_2,2}}^{a_4}
\wh{Z}_{t,v_{e_3,1}}^{a_5} \wh{Z}_{t,v_{e_3,2}}^{a_6}]_{\rm Weyl} \cdot \wh{Z}_{v_t}^{a_7} \\
& = \omega^{a + (-a_1+a_2-a_3+a_4-a_5+a_6 )a_7} \cdot [\wh{Z}_{t,v_{e_1,1}}^{a_1} \wh{Z}_{t,v_{e_1,2}}^{a_2}
\wh{Z}_{t,v_{e_2,1}}^{a_3} \wh{Z}_{t,v_{e_2,2}}^{a_4}
\wh{Z}_{t,v_{e_3,1}}^{a_5} \wh{Z}_{t,v_{e_3,2}}^{a_6} \, \wh{Z}_{v_t}^{a_7} ]_{\rm Weyl},
\end{align*}
it is a straightforward check. For example, consider $(\wh{\bf M}^{\rm 3in}_t)_{3,2,1}=1$, which is the case when $a=0=a_7$; then $a + (-a_1+a_2-a_3+a_4-a_5+a_6 )a_7=0$. For another example, consider $(\wh{\bf M}^{\rm 3out}_t)_{1,1,1}=\omega^{-9} \wh{Z}^3_{v_t}$, which is the case when $a=-9$, $a_7=3$, and $a_1=1$, $a_2=2$, $a_3=1$, $a_4=2$, $a_5=1$, $a_6=2$. One verifies $a + (-a_1+a_2-a_3+a_4-a_5+a_6 )a_7=0$. Others are left as exercise. \qed

\vs

We note that what made \redfix{Lem.\ref{lem:cyclicity_of_3-way_matrices} and} Lem.\ref{lem:standard_3-way_web_values_are_Weyl-ordered} to hold is our choice of isomorphism in eq.\eqref{eq:isomorphism_from_ours_to_Higgins} and the corresponding boundary relations.

\vs

Even so, the definitions of $\wh{\bf M}^{\rm 3out}_t$ and $\wh{\bf M}^{\rm 3in}_t$ may still seem ad hoc at the moment. We will later justify them \redfix{somewhat more}, throughout the present section, especially in \S\ref{subsec:isotopy_invariance}.

\vs

We also observe another important property of ${\rm Tr}^\omega_\Delta$ about elevation reversing.
\begin{definition}
\label{def:star-structure}
Let $t$ be a triangle, and $\Delta$ be an ideal triangulation of a generalized marked surface. For the algebra $\mathcal{Z}^\omega_t$ (resp. $\mathcal{Z}^\omega_\Delta$) (Def.\ref{def:Fock-Goncharov_algebra_quantum}), define the \ul{\em $*$-structure} as the unique ring anti-homomorphism
$$
* : \mathcal{Z}^\omega_t \to \mathcal{Z}^\omega_t \quad (\mbox{resp.} ~ * : \mathcal{Z}^\omega_\Delta \to \mathcal{Z}^\omega_\Delta) \quad : u \mapsto *u,
$$
that sends $\omega^{\pm 1/2}$ to the inverse $\omega^{\mp 1/2}$ and each generator $\wh{Z}_{v,t}^{\pm 1}$ to itself $\wh{Z}_{v,t}^{\pm 1}$ (resp. $\wh{Z}_v^{\pm 1}$ to $\wh{Z}_v^{\pm 1}$).
\end{definition}
It is well known that $*$-structure leaves invariant Weyl-ordered products. Little more generally written:
\begin{definition}
\label{def:non-commutative_Laurent_monomial_and_polynomial}
For any invertible formal variables $\wh{Z}_1,\ldots,\wh{Z}_n$ in an algebra that $\omega$-commute (Def.\ref{def:Weyl-ordered}), an expression $\pm \omega^m \wh{Z}_1^{k_1} \cdots \wh{Z}_n^{k_n}$, for $k_1,\ldots,k_n \in \mathbb{Z}$ and $m \in \frac{1}{2}\mathbb{Z}$, is called an \ul{\em $\omega^{1/2}$-Laurent monomial} in $\wh{Z}_1,\ldots,\wh{Z}_n$. A $\mathbb{Z}$-linear combination of $\omega^{1/2}$-Laurent monomials is called an \ul{\em $\omega^{1/2}$-Laurent polynomial}.
\end{definition}

\begin{lemma}[Weyl-ordered monomial is invariant under $*$-map]
\label{lem:Weyl-ordered_product_is_invariant_under_star-map}
For a fixed $k_1,\ldots,k_n\in\mathbb{Z}$ in Def.\ref{def:non-commutative_Laurent_monomial_and_polynomial}, one has:
\begin{align*}
\omega^m \wh{Z}_1^{k_1} \cdots \wh{Z}_n^{k_n}
= \omega^{-m} \wh{Z}_n^{k_n} \cdots \wh{Z}_1^{k_1}
~\Leftrightarrow~
\omega^m \wh{Z}_1^{k_1} \cdots \wh{Z}_n^{k_n} = [\wh{Z}_1^{k_1} \cdots \wh{Z}_n^{k_n}]_{\rm Weyl}. \qed
\end{align*}
\end{lemma}
So, indeed $*([\wh{Z}_1^{k_1} \cdots \wh{Z}_n^{k_n}]_{\rm Weyl}) = [\wh{Z}_1^{k_1} \cdots \wh{Z}_n^{k_n}]_{\rm Weyl}$ in $\mathcal{Z}^\omega_t$ and $\mathcal{Z}^\omega_\Delta$. Meanwhile, the following is easy to observe from the defining relations of $\mathcal{S}^\omega_{\rm s}(\frak{S};\mathbb{Z})$ (Fig.\ref{fig:A2-skein_relations_quantum}).
\begin{lemma}
\label{lem:elevation-reversing_map}
Let $\frak{S}$ be a generalized marked surface. Define the \ul{\em elevation-reversing map} as the $\mathbb{Z}$-linear map
$$
{\bf r} : \mathcal{S}^\omega_{\rm s}(\frak{S};\mathbb{Z})_{\rm red} \to \mathcal{S}^\omega_{\rm s}(\frak{S};\mathbb{Z})_{\rm red}
$$
that sends $\omega^{\pm 1/2}$ to $\omega^{\mp 1/2}$ and $[W,s]$ to $[W',s']$, where $W'$ is obtained from $W$ by reversing the elevation of all points, i.e. replacing each point $(x,t) \in \frak{S} \times {\bf I}$ of $W$ by $(x,-t)$, and $s'( (x,-t))  = s(x,t)$ for each endpoint $(x,t)$ of $W$. Then ${\bf r}$ is a well-defined ring anti-homomorphism. \qed
\end{lemma}
The property of ${\rm Tr}^\omega_\Delta$ of our interest can then be written as follows, which will be used at the end of the present section:
\begin{proposition}[elevation reversing and $*$-structure are equivariant] 
\label{prop:elevation_reversing_and_star-structure}
Let $\Delta$ be an ideal triangulation of a generalized marked surface $\frak{S}$. Then
\begin{align*}
{\rm Tr}^\omega_\Delta \circ {\bf r} = * \circ {\rm Tr}^\omega_\Delta,
\end{align*}
i.e. ${\rm Tr}^\omega_\Delta({\bf r}(u)) = *({\rm Tr}^\omega_\Delta(u))$ for all $u\in \mathcal{S}^\omega_{\rm s}(\frak{S};\mathbb{Z})_{\rm red}$.
\end{proposition}
At this point, it is not easy to formulate a proof of Prop.\ref{prop:elevation_reversing_and_star-structure}; we will be able to prove it later in the present section (in \S\ref{subsec:quantum_duality_map}), after developing more machinery for ${\rm Tr}^\omega_\Delta$.

\vs

In the upcoming subsections, we shall prove the existence of the ${\rm SL}_3$ quantum and classical trace maps (i.e. prove Thm.\ref{thm:SL3_quantum_trace_map}), study the properties of the values, and relate to the map $\mathbb{I}^+_{{\rm PGL}_3}$ of our original interest in case when $\omega^{1/2}=1$.

\subsection{The biangle ${\rm SL}_3$ quantum trace}
\label{subsec:biangle_SL3_quantum_trace}

Our strategy for a proof of Thm.\ref{thm:SL3_quantum_trace_map} follows the style of Bonahon-Wong \cite{BW}. In particular, we first study the ${\rm SL}_3$ quantum trace map for a {\em biangle} $B$ (Def.\ref{def:triangulation}), i.e. a generalized marked surface diffeomorphic to a closed disc, with two marked points on the boundary, with no punctures. According to later developments by Costantino-L\^e \cite{CL} and Higgins \cite{Higgins}, it is wise to consider also a {\em monogon} $M$, i.e. a generalized marked surface diffeomorphic to a closed disc with one marked point on the boundary. Note that $B$ and $M$ are not triangulable, so they don't really fit into the setting of Thm.\ref{thm:SL3_quantum_trace_map}, hence the ${\rm SL}_3$ quantum/classical trace for them must be dealt with separately. The following is a biangle analog of Thm.\ref{thm:SL3_quantum_trace_map}. 
\begin{proposition}[the biangle ${\rm SL}_3$ quantum trace]
\label{prop:biangle_SL3_quantum_trace}
There exists a unique family of $\mathbb{Z}[\omega^{\pm1/2}]$-algebra homomorphisms
$$
{\rm Tr}^\omega_B ~:~ \mathcal{S}^\omega_{\rm s}(B;\mathbb{Z})_{\rm red} \longrightarrow \mathbb{Z}[\omega^{\pm 1/2}]
$$
defined for biangles $B$, \redfix{satisfying the following.}
\begin{enumerate}
\item[\rm (BT1)] (cutting/gluing property) Let $(W,s)$ be a stated ${\rm SL}_3$-web in $B \times {\bf I}$, and $e$ \redfix{be} an ideal arc in $B$ whose interior lies in the interior of $B$. Let $B'$ be the generalized marked surface obtained from $B$ by cutting along $e$ (Lem.\ref{lem:cutting_process}), so that $B'$ is disjoint union of two biangles $B_1$ and $B_2$. Let $W'$ be the ${\rm SL}_3$-web in $B' \times {\bf I}$ obtained from $W$ by cutting along $e$, and let $W_1$ and $W_2$ be the ${\rm SL}_3$-webs in $B_1 \times {\bf I}$ and $B_2 \times {\bf I}$ such that $W' = W_1 \cup W_2$. Then
$$
{\rm Tr}^\omega_B([W,s]) = \underset{s_1,s_2}{\textstyle \sum} {\rm Tr}^\omega_{B_1}([W_1,s_1]) \, {\rm Tr}^\omega_{B_2}([W_2,s_2]),
$$
where the sum is over all states $s_1$ and $s_2$ of $W_1$ and $W_2$ such that the state $s' := s_1 \cup s_2$ of $W' = W_1 \cup W_2$ is compatible with $s$ in the sense as in Lem.\ref{lem:cutting_process}.

\item[\rm (BT2)] (values at some elementary single-component stated ${\rm SL}_3$-webs with at most one 3-valent vertex)
\begin{enumerate}
\item[\rm (BT2-1)] If the ${\rm SL}_3$-web $W$ in $B\times {\bf I}$ consists of a single edge with no crossing connecting the thickenings of two distinct sides of $B$, and if $\varepsilon,\varepsilon'$ are the values of a state $s$ of $W$ for its two endpoints, then
$$
{\rm Tr}^\omega_B([W,s]) = \mbox{$(\varepsilon,\varepsilon')$-th entry of the $3\times 3$ identity matrix}.
$$

\item[\rm (BT2-2)] If the ${\rm SL}_3$-web $W$ in $B\times {\bf I}$ consists of a single edge with no crossing with the two endpoints lying over a single side of $B$ (i.e. is a U-turn arc), then ${\rm Tr}^\omega_B([W,s])$ is given as in eq.\eqref{eq:Tr_B_U}--\eqref{eq:Tr_B_U4} of Thm.\ref{thm:SL3_quantum_trace_map}\redfix{(QT2-5)}.

\item[\rm (BT2-3)] Suppose $W$ is a 3-way ${\rm SL}_3$-web in $B \times {\bf I}$ with endpoints $x_1,x_2,x_3$, with $x_1,x_2$ lying over one side of $B$ while $x_3$ over the other side, where $\pi(x_1)\to \pi(x_2)$ matches the clockwise orientation of the boundary $\partial B$ (where $\pi : B\times {\bf I} \to B$ is the projection), and let $\varepsilon_1,\varepsilon_2,\varepsilon_3$ be the values of a state $s$ of $W$ for the endpoints $x_1,x_2,x_3$. The value ${\rm Tr}^\omega_B([W,s])$ equals $(\wh{\bf F}^h_{k,\varepsilon_3})_{\varepsilon_1,\varepsilon_2}$ of eq.\eqref{eq:fork_matrix}, where $h = {\rm out}$ (resp. ${\rm in}$) if $x_1,x_2,x_3$ are sinks (resp. sources) so that $W$ is outgoing (resp. incoming) 3-way web, and $k = +$ (resp. $-$) if $x_1 \succ x_2$ (resp. $x_1 \prec x_2$).

\end{enumerate}
\end{enumerate}
\end{proposition}

\begin{proposition}[values of the biangle ${\rm SL}_3$ quantum trace at some more elementary ${\rm SL}_3$-webs]
\label{prop:biangle_SL3_quantum_trace_some_values}
The properties (BT1) and (BT2) of Prop.\ref{prop:biangle_SL3_quantum_trace} for the biangle ${\rm SL}_3$ quantum trace ${\rm Tr}^\omega_B$ imply:
\begin{enumerate}
\item[]

\begin{enumerate}
\item[\rm (BT2-4)] If an ${\rm SL}_3$-web $W$ in $B\times {\bf I}$ consists of two edges, each connecting the thickenings of two distinct sides of $B$, with the number of crossings being $0$ or $1$, and if $W$ is not a product of two single-edge ${\rm SL}_3$-webs of type (BT2-1), then ${\rm Tr}^\omega_B([W,s])$ is given by:
\begin{align}
\label{eq:height_exchange_and_crossing1}
& {\rm Tr}^\omega_B([\,\, \raisebox{-0.5\height}{
\begingroup%
  \makeatletter%
  \providecommand\color[2][]{%
    \errmessage{(Inkscape) Color is used for the text in Inkscape, but the package 'color.sty' is not loaded}%
    \renewcommand\color[2][]{}%
  }%
  \providecommand\transparent[1]{%
    \errmessage{(Inkscape) Transparency is used (non-zero) for the text in Inkscape, but the package 'transparent.sty' is not loaded}%
    \renewcommand\transparent[1]{}%
  }%
  \providecommand\rotatebox[2]{#2}%
  \newcommand*\fsize{\dimexpr\f@size pt\relax}%
  \newcommand*\lineheight[1]{\fontsize{\fsize}{#1\fsize}\selectfont}%
  \ifx\svgwidth\undefined%
    \setlength{\unitlength}{56.69291339bp}%
    \ifx\svgscale\undefined%
      \relax%
    \else%
      \setlength{\unitlength}{\unitlength * \real{\svgscale}}%
    \fi%
  \else%
    \setlength{\unitlength}{\svgwidth}%
  \fi%
  \global\let\svgwidth\undefined%
  \global\let\svgscale\undefined%
  \makeatother%
  \begin{picture}(1,0.875)%
    \lineheight{1}%
    \setlength\tabcolsep{0pt}%
    \put(0,0){\includegraphics[width=\unitlength,page=1]{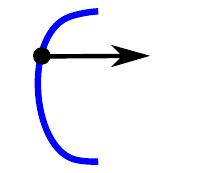}}%
    \put(0.00342324,0.57251168){\color[rgb]{0,0,0}\makebox(0,0)[lt]{\lineheight{1.25}\smash{\begin{tabular}[t]{l}$x_1$\end{tabular}}}}%
    \put(0,0){\includegraphics[width=\unitlength,page=2]{biangle_height_exchange1.pdf}}%
    \put(0.00801424,0.25622502){\color[rgb]{0,0,0}\makebox(0,0)[lt]{\lineheight{1.25}\smash{\begin{tabular}[t]{l}$x_2$\end{tabular}}}}%
    \put(0.05924823,0.51254405){\color[rgb]{0,0,0}\rotatebox{-92.070464}{\makebox(0,0)[lt]{\lineheight{1.25}\smash{\begin{tabular}[t]{l}$\prec$\end{tabular}}}}}%
    \put(0,0){\includegraphics[width=\unitlength,page=3]{biangle_height_exchange1.pdf}}%
    \put(0.85239343,0.58058091){\color[rgb]{0,0,0}\makebox(0,0)[lt]{\lineheight{1.25}\smash{\begin{tabular}[t]{l}$y_1$\end{tabular}}}}%
    \put(0.85698442,0.26429425){\color[rgb]{0,0,0}\makebox(0,0)[lt]{\lineheight{1.25}\smash{\begin{tabular}[t]{l}$y_2$\end{tabular}}}}%
    \put(0.89099145,0.50338616){\color[rgb]{0,0,0}\rotatebox{-92.070464}{\makebox(0,0)[lt]{\lineheight{1.25}\smash{\begin{tabular}[t]{l}$\succ$\end{tabular}}}}}%
  \end{picture}%
\endgroup%
} \,\, ])
= {\rm Tr}^\omega_B([\,\, \raisebox{-0.5\height}{
\begingroup%
  \makeatletter%
  \providecommand\color[2][]{%
    \errmessage{(Inkscape) Color is used for the text in Inkscape, but the package 'color.sty' is not loaded}%
    \renewcommand\color[2][]{}%
  }%
  \providecommand\transparent[1]{%
    \errmessage{(Inkscape) Transparency is used (non-zero) for the text in Inkscape, but the package 'transparent.sty' is not loaded}%
    \renewcommand\transparent[1]{}%
  }%
  \providecommand\rotatebox[2]{#2}%
  \newcommand*\fsize{\dimexpr\f@size pt\relax}%
  \newcommand*\lineheight[1]{\fontsize{\fsize}{#1\fsize}\selectfont}%
  \ifx\svgwidth\undefined%
    \setlength{\unitlength}{56.69291339bp}%
    \ifx\svgscale\undefined%
      \relax%
    \else%
      \setlength{\unitlength}{\unitlength * \real{\svgscale}}%
    \fi%
  \else%
    \setlength{\unitlength}{\svgwidth}%
  \fi%
  \global\let\svgwidth\undefined%
  \global\let\svgscale\undefined%
  \makeatother%
  \begin{picture}(1,0.875)%
    \lineheight{1}%
    \setlength\tabcolsep{0pt}%
    \put(0,0){\includegraphics[width=\unitlength,page=1]{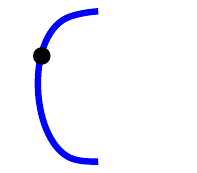}}%
    \put(-0.00171479,0.57550812){\color[rgb]{0,0,0}\makebox(0,0)[lt]{\lineheight{1.25}\smash{\begin{tabular}[t]{l}$x_2$\end{tabular}}}}%
    \put(0,0){\includegraphics[width=\unitlength,page=2]{biangle_crossing3.pdf}}%
    \put(0.00287621,0.25922147){\color[rgb]{0,0,0}\makebox(0,0)[lt]{\lineheight{1.25}\smash{\begin{tabular}[t]{l}$x_1$\end{tabular}}}}%
    \put(0.03772366,0.5151652){\color[rgb]{0,0,0}\rotatebox{-90.362304}{\makebox(0,0)[lt]{\lineheight{1.25}\smash{\begin{tabular}[t]{l}$\succ$\end{tabular}}}}}%
    \put(0,0){\includegraphics[width=\unitlength,page=3]{biangle_crossing3.pdf}}%
    \put(0.85239343,0.58058091){\color[rgb]{0,0,0}\makebox(0,0)[lt]{\lineheight{1.25}\smash{\begin{tabular}[t]{l}$y_1$\end{tabular}}}}%
    \put(0.85698442,0.26429425){\color[rgb]{0,0,0}\makebox(0,0)[lt]{\lineheight{1.25}\smash{\begin{tabular}[t]{l}$y_2$\end{tabular}}}}%
    \put(0.89099145,0.50338616){\color[rgb]{0,0,0}\rotatebox{-92.070464}{\makebox(0,0)[lt]{\lineheight{1.25}\smash{\begin{tabular}[t]{l}$\succ$\end{tabular}}}}}%
    \put(0,0){\includegraphics[width=\unitlength,page=4]{biangle_crossing3.pdf}}%
  \end{picture}%
\endgroup%
} \,\, ]) =  {\rm Tr}^\omega_B([\,\, \raisebox{-0.5\height}{
\begingroup%
  \makeatletter%
  \providecommand\color[2][]{%
    \errmessage{(Inkscape) Color is used for the text in Inkscape, but the package 'color.sty' is not loaded}%
    \renewcommand\color[2][]{}%
  }%
  \providecommand\transparent[1]{%
    \errmessage{(Inkscape) Transparency is used (non-zero) for the text in Inkscape, but the package 'transparent.sty' is not loaded}%
    \renewcommand\transparent[1]{}%
  }%
  \providecommand\rotatebox[2]{#2}%
  \newcommand*\fsize{\dimexpr\f@size pt\relax}%
  \newcommand*\lineheight[1]{\fontsize{\fsize}{#1\fsize}\selectfont}%
  \ifx\svgwidth\undefined%
    \setlength{\unitlength}{56.69291339bp}%
    \ifx\svgscale\undefined%
      \relax%
    \else%
      \setlength{\unitlength}{\unitlength * \real{\svgscale}}%
    \fi%
  \else%
    \setlength{\unitlength}{\svgwidth}%
  \fi%
  \global\let\svgwidth\undefined%
  \global\let\svgscale\undefined%
  \makeatother%
  \begin{picture}(1,0.875)%
    \lineheight{1}%
    \setlength\tabcolsep{0pt}%
    \put(0,0){\includegraphics[width=\unitlength,page=1]{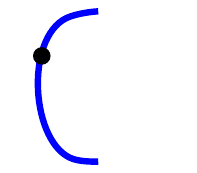}}%
    \put(-0.00171479,0.57550812){\color[rgb]{0,0,0}\makebox(0,0)[lt]{\lineheight{1.25}\smash{\begin{tabular}[t]{l}$x_1$\end{tabular}}}}%
    \put(0,0){\includegraphics[width=\unitlength,page=2]{biangle_crossing4.pdf}}%
    \put(0.00287621,0.25922147){\color[rgb]{0,0,0}\makebox(0,0)[lt]{\lineheight{1.25}\smash{\begin{tabular}[t]{l}$x_2$\end{tabular}}}}%
    \put(0.03772366,0.5151652){\color[rgb]{0,0,0}\rotatebox{-90.362304}{\makebox(0,0)[lt]{\lineheight{1.25}\smash{\begin{tabular}[t]{l}$\prec$\end{tabular}}}}}%
    \put(0,0){\includegraphics[width=\unitlength,page=3]{biangle_crossing4.pdf}}%
    \put(0.85239343,0.58058091){\color[rgb]{0,0,0}\makebox(0,0)[lt]{\lineheight{1.25}\smash{\begin{tabular}[t]{l}$y_2$\end{tabular}}}}%
    \put(0.85698442,0.26429425){\color[rgb]{0,0,0}\makebox(0,0)[lt]{\lineheight{1.25}\smash{\begin{tabular}[t]{l}$y_1$\end{tabular}}}}%
    \put(0.89099145,0.50338616){\color[rgb]{0,0,0}\rotatebox{-92.070464}{\makebox(0,0)[lt]{\lineheight{1.25}\smash{\begin{tabular}[t]{l}$\prec$\end{tabular}}}}}%
    \put(0,0){\includegraphics[width=\unitlength,page=4]{biangle_crossing4.pdf}}%
  \end{picture}%
\endgroup%
} \,\, ]) \\
\label{eq:height_exchange_and_crossing1_values}
& = 
\left\{
\begin{array}{cl}
q^{-2/3} & \mbox{if $(s(x_1),s(x_2),s(y_1),s(y_2)) \in \{(1,1,1,1),(2,2,2,2),(3,3,3,3)\}$,} \\
q^{1/3} & \mbox{if $(s(x_1),s(x_2),s(y_1),s(y_2)) \in \{(1,2,1,2),(1,3,1,3),(2,3,2,3)$,} \\
& \mbox{\hspace*{45mm} $(2,1,2,1),(3,1,3,1),(3,2,3,2)\}$,} \\
q^{-2/3} - q^{4/3} & \mbox{if $(s(x_1),s(x_2),s(y_1),s(y_2)) \in \{(1,2,2,1),(1,3,3,1),(2,3,3,2)\}$}, \\
0 & \mbox{otherwise},
\end{array}
\right.
\end{align}
\begin{align}
\label{eq:height_exchange_and_crossing3}
& {\rm Tr}^\omega_B([\,\, \raisebox{-0.5\height}{
\begingroup%
  \makeatletter%
  \providecommand\color[2][]{%
    \errmessage{(Inkscape) Color is used for the text in Inkscape, but the package 'color.sty' is not loaded}%
    \renewcommand\color[2][]{}%
  }%
  \providecommand\transparent[1]{%
    \errmessage{(Inkscape) Transparency is used (non-zero) for the text in Inkscape, but the package 'transparent.sty' is not loaded}%
    \renewcommand\transparent[1]{}%
  }%
  \providecommand\rotatebox[2]{#2}%
  \newcommand*\fsize{\dimexpr\f@size pt\relax}%
  \newcommand*\lineheight[1]{\fontsize{\fsize}{#1\fsize}\selectfont}%
  \ifx\svgwidth\undefined%
    \setlength{\unitlength}{56.69291339bp}%
    \ifx\svgscale\undefined%
      \relax%
    \else%
      \setlength{\unitlength}{\unitlength * \real{\svgscale}}%
    \fi%
  \else%
    \setlength{\unitlength}{\svgwidth}%
  \fi%
  \global\let\svgwidth\undefined%
  \global\let\svgscale\undefined%
  \makeatother%
  \begin{picture}(1,0.875)%
    \lineheight{1}%
    \setlength\tabcolsep{0pt}%
    \put(0,0){\includegraphics[width=\unitlength,page=1]{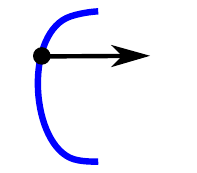}}%
    \put(0.00342324,0.57251168){\color[rgb]{0,0,0}\makebox(0,0)[lt]{\lineheight{1.25}\smash{\begin{tabular}[t]{l}$x_1$\end{tabular}}}}%
    \put(0,0){\includegraphics[width=\unitlength,page=2]{biangle_height_exchange3.pdf}}%
    \put(0.00801424,0.25622502){\color[rgb]{0,0,0}\makebox(0,0)[lt]{\lineheight{1.25}\smash{\begin{tabular}[t]{l}$x_2$\end{tabular}}}}%
    \put(0.05924823,0.51254405){\color[rgb]{0,0,0}\rotatebox{-92.070464}{\makebox(0,0)[lt]{\lineheight{1.25}\smash{\begin{tabular}[t]{l}$\prec$\end{tabular}}}}}%
    \put(0,0){\includegraphics[width=\unitlength,page=3]{biangle_height_exchange3.pdf}}%
    \put(0.85239343,0.58058091){\color[rgb]{0,0,0}\makebox(0,0)[lt]{\lineheight{1.25}\smash{\begin{tabular}[t]{l}$y_1$\end{tabular}}}}%
    \put(0.85698442,0.26429425){\color[rgb]{0,0,0}\makebox(0,0)[lt]{\lineheight{1.25}\smash{\begin{tabular}[t]{l}$y_2$\end{tabular}}}}%
    \put(0.89099145,0.50338616){\color[rgb]{0,0,0}\rotatebox{-92.070464}{\makebox(0,0)[lt]{\lineheight{1.25}\smash{\begin{tabular}[t]{l}$\succ$\end{tabular}}}}}%
  \end{picture}%
\endgroup%
} \,\, ])
= {\rm Tr}^\omega_B([\,\, \raisebox{-0.5\height}{
\begingroup%
  \makeatletter%
  \providecommand\color[2][]{%
    \errmessage{(Inkscape) Color is used for the text in Inkscape, but the package 'color.sty' is not loaded}%
    \renewcommand\color[2][]{}%
  }%
  \providecommand\transparent[1]{%
    \errmessage{(Inkscape) Transparency is used (non-zero) for the text in Inkscape, but the package 'transparent.sty' is not loaded}%
    \renewcommand\transparent[1]{}%
  }%
  \providecommand\rotatebox[2]{#2}%
  \newcommand*\fsize{\dimexpr\f@size pt\relax}%
  \newcommand*\lineheight[1]{\fontsize{\fsize}{#1\fsize}\selectfont}%
  \ifx\svgwidth\undefined%
    \setlength{\unitlength}{56.69291339bp}%
    \ifx\svgscale\undefined%
      \relax%
    \else%
      \setlength{\unitlength}{\unitlength * \real{\svgscale}}%
    \fi%
  \else%
    \setlength{\unitlength}{\svgwidth}%
  \fi%
  \global\let\svgwidth\undefined%
  \global\let\svgscale\undefined%
  \makeatother%
  \begin{picture}(1,0.875)%
    \lineheight{1}%
    \setlength\tabcolsep{0pt}%
    \put(0,0){\includegraphics[width=\unitlength,page=1]{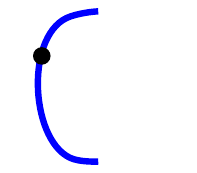}}%
    \put(-0.00171479,0.57550812){\color[rgb]{0,0,0}\makebox(0,0)[lt]{\lineheight{1.25}\smash{\begin{tabular}[t]{l}$x_2$\end{tabular}}}}%
    \put(0,0){\includegraphics[width=\unitlength,page=2]{biangle_crossing7.pdf}}%
    \put(0.00287621,0.25922147){\color[rgb]{0,0,0}\makebox(0,0)[lt]{\lineheight{1.25}\smash{\begin{tabular}[t]{l}$x_1$\end{tabular}}}}%
    \put(0.03772366,0.5151652){\color[rgb]{0,0,0}\rotatebox{-90.362304}{\makebox(0,0)[lt]{\lineheight{1.25}\smash{\begin{tabular}[t]{l}$\succ$\end{tabular}}}}}%
    \put(0,0){\includegraphics[width=\unitlength,page=3]{biangle_crossing7.pdf}}%
    \put(0.85239343,0.58058091){\color[rgb]{0,0,0}\makebox(0,0)[lt]{\lineheight{1.25}\smash{\begin{tabular}[t]{l}$y_1$\end{tabular}}}}%
    \put(0.85698442,0.26429425){\color[rgb]{0,0,0}\makebox(0,0)[lt]{\lineheight{1.25}\smash{\begin{tabular}[t]{l}$y_2$\end{tabular}}}}%
    \put(0.89099145,0.50338616){\color[rgb]{0,0,0}\rotatebox{-92.070464}{\makebox(0,0)[lt]{\lineheight{1.25}\smash{\begin{tabular}[t]{l}$\succ$\end{tabular}}}}}%
    \put(0,0){\includegraphics[width=\unitlength,page=4]{biangle_crossing7.pdf}}%
  \end{picture}%
\endgroup%
} \,\, ]) = {\rm Tr}^\omega_B([\,\, \raisebox{-0.5\height}{
\begingroup%
  \makeatletter%
  \providecommand\color[2][]{%
    \errmessage{(Inkscape) Color is used for the text in Inkscape, but the package 'color.sty' is not loaded}%
    \renewcommand\color[2][]{}%
  }%
  \providecommand\transparent[1]{%
    \errmessage{(Inkscape) Transparency is used (non-zero) for the text in Inkscape, but the package 'transparent.sty' is not loaded}%
    \renewcommand\transparent[1]{}%
  }%
  \providecommand\rotatebox[2]{#2}%
  \newcommand*\fsize{\dimexpr\f@size pt\relax}%
  \newcommand*\lineheight[1]{\fontsize{\fsize}{#1\fsize}\selectfont}%
  \ifx\svgwidth\undefined%
    \setlength{\unitlength}{56.69291339bp}%
    \ifx\svgscale\undefined%
      \relax%
    \else%
      \setlength{\unitlength}{\unitlength * \real{\svgscale}}%
    \fi%
  \else%
    \setlength{\unitlength}{\svgwidth}%
  \fi%
  \global\let\svgwidth\undefined%
  \global\let\svgscale\undefined%
  \makeatother%
  \begin{picture}(1,0.875)%
    \lineheight{1}%
    \setlength\tabcolsep{0pt}%
    \put(0,0){\includegraphics[width=\unitlength,page=1]{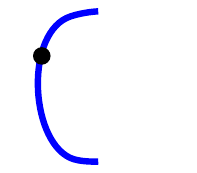}}%
    \put(-0.00171479,0.57550812){\color[rgb]{0,0,0}\makebox(0,0)[lt]{\lineheight{1.25}\smash{\begin{tabular}[t]{l}$x_1$\end{tabular}}}}%
    \put(0,0){\includegraphics[width=\unitlength,page=2]{biangle_crossing8.pdf}}%
    \put(0.00287621,0.25922147){\color[rgb]{0,0,0}\makebox(0,0)[lt]{\lineheight{1.25}\smash{\begin{tabular}[t]{l}$x_2$\end{tabular}}}}%
    \put(0.03772366,0.5151652){\color[rgb]{0,0,0}\rotatebox{-90.362304}{\makebox(0,0)[lt]{\lineheight{1.25}\smash{\begin{tabular}[t]{l}$\prec$\end{tabular}}}}}%
    \put(0,0){\includegraphics[width=\unitlength,page=3]{biangle_crossing8.pdf}}%
    \put(0.85239343,0.58058091){\color[rgb]{0,0,0}\makebox(0,0)[lt]{\lineheight{1.25}\smash{\begin{tabular}[t]{l}$y_2$\end{tabular}}}}%
    \put(0.85698442,0.26429425){\color[rgb]{0,0,0}\makebox(0,0)[lt]{\lineheight{1.25}\smash{\begin{tabular}[t]{l}$y_1$\end{tabular}}}}%
    \put(0.89099145,0.50338616){\color[rgb]{0,0,0}\rotatebox{-92.070464}{\makebox(0,0)[lt]{\lineheight{1.25}\smash{\begin{tabular}[t]{l}$\prec$\end{tabular}}}}}%
    \put(0,0){\includegraphics[width=\unitlength,page=4]{biangle_crossing8.pdf}}%
  \end{picture}%
\endgroup%
} \,\, ]) \\
\label{eq:height_exchange_and_crossing3_values}
& = 
\left\{
\begin{array}{cl}
q^{-1/3} & \mbox{if $(s(x_1),s(x_2),s(y_1),s(y_2)) \in \{(1,1,1,1),(1,2,1,2),(2,1,2,1),$} \\
& \mbox{\hspace*{45mm} $(2,3,2,3),(3,2,3,2), (3,3,3,3)\}$,} \\
q^{-1/3} - q^{5/3} & \mbox{if $(s(x_1),s(x_2),s(y_1),s(y_2)) \in \{(1,3,2,2),(2,2,3,1)\}$}, \\
q^{2/3} & \mbox{if $(s(x_1),s(x_2),s(y_1),s(y_2)) \in \{(1,3,1,3),(2,2,2,2),(3,1,3,1)\}$}, \\
q^{8/3} - q^{2/3} & \mbox{if $(s(x_1),s(x_2),s(y_1),s(y_2)) =(1,3,3,1)$}, \\
0 & \mbox{otherwise},
\end{array}
\right.
\end{align}
\begin{align}
\label{eq:height_exchange_and_crossing5}
& {\rm Tr}^\omega_B([\,\, \raisebox{-0.5\height}{
\begingroup%
  \makeatletter%
  \providecommand\color[2][]{%
    \errmessage{(Inkscape) Color is used for the text in Inkscape, but the package 'color.sty' is not loaded}%
    \renewcommand\color[2][]{}%
  }%
  \providecommand\transparent[1]{%
    \errmessage{(Inkscape) Transparency is used (non-zero) for the text in Inkscape, but the package 'transparent.sty' is not loaded}%
    \renewcommand\transparent[1]{}%
  }%
  \providecommand\rotatebox[2]{#2}%
  \newcommand*\fsize{\dimexpr\f@size pt\relax}%
  \newcommand*\lineheight[1]{\fontsize{\fsize}{#1\fsize}\selectfont}%
  \ifx\svgwidth\undefined%
    \setlength{\unitlength}{56.69291339bp}%
    \ifx\svgscale\undefined%
      \relax%
    \else%
      \setlength{\unitlength}{\unitlength * \real{\svgscale}}%
    \fi%
  \else%
    \setlength{\unitlength}{\svgwidth}%
  \fi%
  \global\let\svgwidth\undefined%
  \global\let\svgscale\undefined%
  \makeatother%
  \begin{picture}(1,0.875)%
    \lineheight{1}%
    \setlength\tabcolsep{0pt}%
    \put(0,0){\includegraphics[width=\unitlength,page=1]{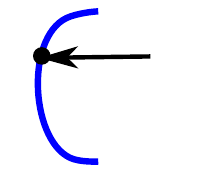}}%
    \put(0.00342324,0.57251168){\color[rgb]{0,0,0}\makebox(0,0)[lt]{\lineheight{1.25}\smash{\begin{tabular}[t]{l}$x_1$\end{tabular}}}}%
    \put(0,0){\includegraphics[width=\unitlength,page=2]{biangle_height_exchange5.pdf}}%
    \put(0.00801424,0.25622502){\color[rgb]{0,0,0}\makebox(0,0)[lt]{\lineheight{1.25}\smash{\begin{tabular}[t]{l}$x_2$\end{tabular}}}}%
    \put(0.05924823,0.51254405){\color[rgb]{0,0,0}\rotatebox{-92.070464}{\makebox(0,0)[lt]{\lineheight{1.25}\smash{\begin{tabular}[t]{l}$\prec$\end{tabular}}}}}%
    \put(0,0){\includegraphics[width=\unitlength,page=3]{biangle_height_exchange5.pdf}}%
    \put(0.85239343,0.58058091){\color[rgb]{0,0,0}\makebox(0,0)[lt]{\lineheight{1.25}\smash{\begin{tabular}[t]{l}$y_1$\end{tabular}}}}%
    \put(0.85698442,0.26429425){\color[rgb]{0,0,0}\makebox(0,0)[lt]{\lineheight{1.25}\smash{\begin{tabular}[t]{l}$y_2$\end{tabular}}}}%
    \put(0.89099145,0.50338616){\color[rgb]{0,0,0}\rotatebox{-92.070464}{\makebox(0,0)[lt]{\lineheight{1.25}\smash{\begin{tabular}[t]{l}$\succ$\end{tabular}}}}}%
  \end{picture}%
\endgroup%
} \,\, ])
= {\rm Tr}^\omega_B([\,\, \raisebox{-0.5\height}{
\begingroup%
  \makeatletter%
  \providecommand\color[2][]{%
    \errmessage{(Inkscape) Color is used for the text in Inkscape, but the package 'color.sty' is not loaded}%
    \renewcommand\color[2][]{}%
  }%
  \providecommand\transparent[1]{%
    \errmessage{(Inkscape) Transparency is used (non-zero) for the text in Inkscape, but the package 'transparent.sty' is not loaded}%
    \renewcommand\transparent[1]{}%
  }%
  \providecommand\rotatebox[2]{#2}%
  \newcommand*\fsize{\dimexpr\f@size pt\relax}%
  \newcommand*\lineheight[1]{\fontsize{\fsize}{#1\fsize}\selectfont}%
  \ifx\svgwidth\undefined%
    \setlength{\unitlength}{56.69291339bp}%
    \ifx\svgscale\undefined%
      \relax%
    \else%
      \setlength{\unitlength}{\unitlength * \real{\svgscale}}%
    \fi%
  \else%
    \setlength{\unitlength}{\svgwidth}%
  \fi%
  \global\let\svgwidth\undefined%
  \global\let\svgscale\undefined%
  \makeatother%
  \begin{picture}(1,0.875)%
    \lineheight{1}%
    \setlength\tabcolsep{0pt}%
    \put(0,0){\includegraphics[width=\unitlength,page=1]{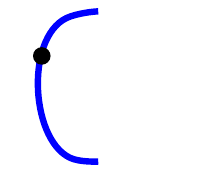}}%
    \put(-0.00171479,0.57550812){\color[rgb]{0,0,0}\makebox(0,0)[lt]{\lineheight{1.25}\smash{\begin{tabular}[t]{l}$x_2$\end{tabular}}}}%
    \put(0,0){\includegraphics[width=\unitlength,page=2]{biangle_crossing11.pdf}}%
    \put(0.00287621,0.25922147){\color[rgb]{0,0,0}\makebox(0,0)[lt]{\lineheight{1.25}\smash{\begin{tabular}[t]{l}$x_1$\end{tabular}}}}%
    \put(0.03772366,0.5151652){\color[rgb]{0,0,0}\rotatebox{-90.362304}{\makebox(0,0)[lt]{\lineheight{1.25}\smash{\begin{tabular}[t]{l}$\succ$\end{tabular}}}}}%
    \put(0,0){\includegraphics[width=\unitlength,page=3]{biangle_crossing11.pdf}}%
    \put(0.85239343,0.58058091){\color[rgb]{0,0,0}\makebox(0,0)[lt]{\lineheight{1.25}\smash{\begin{tabular}[t]{l}$y_1$\end{tabular}}}}%
    \put(0.85698442,0.26429425){\color[rgb]{0,0,0}\makebox(0,0)[lt]{\lineheight{1.25}\smash{\begin{tabular}[t]{l}$y_2$\end{tabular}}}}%
    \put(0.89099145,0.50338616){\color[rgb]{0,0,0}\rotatebox{-92.070464}{\makebox(0,0)[lt]{\lineheight{1.25}\smash{\begin{tabular}[t]{l}$\succ$\end{tabular}}}}}%
    \put(0,0){\includegraphics[width=\unitlength,page=4]{biangle_crossing11.pdf}}%
  \end{picture}%
\endgroup%
} \,\, ]) = {\rm Tr}^\omega_B([\,\, \raisebox{-0.5\height}{
\begingroup%
  \makeatletter%
  \providecommand\color[2][]{%
    \errmessage{(Inkscape) Color is used for the text in Inkscape, but the package 'color.sty' is not loaded}%
    \renewcommand\color[2][]{}%
  }%
  \providecommand\transparent[1]{%
    \errmessage{(Inkscape) Transparency is used (non-zero) for the text in Inkscape, but the package 'transparent.sty' is not loaded}%
    \renewcommand\transparent[1]{}%
  }%
  \providecommand\rotatebox[2]{#2}%
  \newcommand*\fsize{\dimexpr\f@size pt\relax}%
  \newcommand*\lineheight[1]{\fontsize{\fsize}{#1\fsize}\selectfont}%
  \ifx\svgwidth\undefined%
    \setlength{\unitlength}{56.69291339bp}%
    \ifx\svgscale\undefined%
      \relax%
    \else%
      \setlength{\unitlength}{\unitlength * \real{\svgscale}}%
    \fi%
  \else%
    \setlength{\unitlength}{\svgwidth}%
  \fi%
  \global\let\svgwidth\undefined%
  \global\let\svgscale\undefined%
  \makeatother%
  \begin{picture}(1,0.875)%
    \lineheight{1}%
    \setlength\tabcolsep{0pt}%
    \put(0,0){\includegraphics[width=\unitlength,page=1]{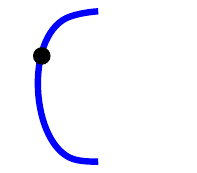}}%
    \put(-0.00171479,0.57550812){\color[rgb]{0,0,0}\makebox(0,0)[lt]{\lineheight{1.25}\smash{\begin{tabular}[t]{l}$x_1$\end{tabular}}}}%
    \put(0,0){\includegraphics[width=\unitlength,page=2]{biangle_crossing12.pdf}}%
    \put(0.00287621,0.25922147){\color[rgb]{0,0,0}\makebox(0,0)[lt]{\lineheight{1.25}\smash{\begin{tabular}[t]{l}$x_2$\end{tabular}}}}%
    \put(0.03772366,0.5151652){\color[rgb]{0,0,0}\rotatebox{-90.362304}{\makebox(0,0)[lt]{\lineheight{1.25}\smash{\begin{tabular}[t]{l}$\prec$\end{tabular}}}}}%
    \put(0,0){\includegraphics[width=\unitlength,page=3]{biangle_crossing12.pdf}}%
    \put(0.85239343,0.58058091){\color[rgb]{0,0,0}\makebox(0,0)[lt]{\lineheight{1.25}\smash{\begin{tabular}[t]{l}$y_2$\end{tabular}}}}%
    \put(0.85698442,0.26429425){\color[rgb]{0,0,0}\makebox(0,0)[lt]{\lineheight{1.25}\smash{\begin{tabular}[t]{l}$y_1$\end{tabular}}}}%
    \put(0.89099145,0.50338616){\color[rgb]{0,0,0}\rotatebox{-92.070464}{\makebox(0,0)[lt]{\lineheight{1.25}\smash{\begin{tabular}[t]{l}$\prec$\end{tabular}}}}}%
    \put(0,0){\includegraphics[width=\unitlength,page=4]{biangle_crossing12.pdf}}%
  \end{picture}%
\endgroup%
} \,\, ]) = \mbox{eq.\eqref{eq:height_exchange_and_crossing3_values}.}
\end{align}
\end{enumerate}
\end{enumerate}
\end{proposition}

\begin{remark}
\label{rem:Douglas_biangle}
The values for the cases (BT2-1), (BT2-2) and (BT2-4) which do not involve 3-valent vertices coincide with those given in \cite{Douglas} \cite{Douglas21}. We note that the cases when there are 3-valent internal vertices, i.e. (BT2-3), are not dealt with in \cite{Douglas} \cite{Douglas21}. We were informed by Daniel Douglas that he also found his values for the case (BT2-3), which are different from our values; we believe that his values can be obtained by using a different choice of twisting scalars in eq.\eqref{eq:twisting_scalars}. The values for (BT2-4) will be seen as entries of $9\times 9$ matrix later in the present subsection, which makes it easier to compare with \cite{Douglas} \cite{Douglas21}. The values for 3-way webs whose endpoints all lie over a same side of $B$ can be deduced from (BT1), (BT2-1), (BT2-2), (BT2-3), and some more relations of $\mathcal{S}^\omega_{\rm s}(B;\mathbb{Z})_{\rm red}$ to be studied (Lem.\ref{lem:3-way_boundary_relations1}--\ref{lem:3-way_boundary_relations_reversed2}).
\end{remark}

\begin{corollary}[the biangle ${\rm SL}_3$ classical trace; \cite{Kim}]
\label{cor:biangle_SL3_classical_trace}
Prop.\ref{prop:biangle_SL3_quantum_trace} holds when $\omega^{1/2}=1$. The resulting map
$$
{\rm Tr}_B = {\rm Tr}_B^1 : \mathcal{S}^1_{\rm s}(B;\mathbb{Z})_{\rm red} \cong \mathcal{S}_{\rm s}(B;\mathbb{Z})_{\rm red} \to \mathbb{Z}
$$
is called the \ul{\em biangle ${\rm SL}_3$ classical trace}.
\end{corollary}

We first prove the uniqueness part of Prop.\ref{prop:biangle_SL3_quantum_trace}, which is a lot easier than the existence part. It is useful to adapt the viewpoint of the Reshetikhin-Turaev operator invariants \cite{RT}, as done for ${\rm SL}_2$ biangle quantum trace in \cite{CL} \cite{KLS}. Namely, to each ${\rm SL}_3$-web $W$ in $B\times {\bf I}$ we associate a $\mathbb{Z}[\omega^{\pm 1/2}]$-linear operator between $\mathbb{Z}[\omega^{\pm 1/2}]$-modules, whose matrix elements are ${\rm Tr}^\omega_B([W,s])$. To give a more precise description, we should choose preferred orientations on boundary arcs of $B$. Terminology is borrowed from \cite{KLS}  \cite{CL}.
\begin{definition}
\label{def:direction}
A \ul{\em direction} of a biangle $B$ is the choice of a distinguished marked point of $B$, denoted by ${\rm dir}$. The distinguished marked point is called the \ul{\em top} marked point, while the other marked point the \ul{\em bottom} marked point. The pair $(B,{\rm dir})$ is called a \ul{\em directed biangle} and is denoted by $\vec{B}$. For a directed biangle $\vec{B}$, the \ul{\em induced orientation} on the boundary arcs of $B$ are the orientations pointing toward the distinguished marked point. The boundary arc of $B$ whose induced orientation matches the clockwise orientation (coming from the surface orientation of $B$) is called the \ul{\em left side} of $\vec{B}$, and the other boundary arc the \ul{\em right side}.
\end{definition}
We will use the notations $\mathcal{S}^\omega(\vec{B};\mathcal{R})$, $\mathcal{S}^\omega_{\rm s}(\vec{B};\mathcal{R})$ and $\mathcal{S}^\omega_{\rm s}(\vec{B};\mathcal{R})_{\rm red}$ to mean the algebras $\mathcal{S}^\omega(B;\mathcal{R})$, $\mathcal{S}^\omega_{\rm s}(B;\mathcal{R})$ and $\mathcal{S}^\omega_{\rm s}(B;\mathcal{R})_{\rm red}$ together with the information of direction on $B$.

\begin{definition}
Let $W$ be an ${\rm SL}_3$-web in $\vec{B} \times {\bf I}$. Let
$$
\mbox{$\partial_{\rm left} W$ (resp. $\partial_{\rm right} W$) }
$$
be the set of all endpoints of $W$ lying over the left side (resp. the right side) of $\vec{B}$.
\end{definition}
Consider
$$
V := \mbox{free $\mathbb{Z}[\omega^{\pm 1/2}]$-module of rank $3$, with basis $\xi_1,\xi_2,\xi_3$,}
$$
equipped with a $\mathbb{Z}[\omega^{\pm 1/2}]$-bilinear pairing
$$
\langle \,,\, \rangle : V\otimes V \to \mathbb{Z}[\omega^{\pm 1/2}], \quad \mbox{given by} \quad \langle \xi_i, \xi_j\rangle = \delta_{i,j}, ~ \forall i,j,
$$
which induces a pairing $\langle \,,\, \rangle : V^{\otimes n} \otimes V^{\otimes n} \to \mathbb{Z}[\omega^{\pm 1/2}]$ for each $n \ge 1$. For an ${\rm SL}_3$-web $W$ in $\vec{B} \times {\bf I}$, we will define a $\mathbb{Z}[\omega^{\pm 1/2}]$-linear map
$$
\rho_W : V^{\otimes \partial_{\rm right} W} \to V^{\otimes \partial_{\rm left} W},
$$
where $V^{\otimes \partial_{\rm left} W}$ means $V^{\otimes |\partial_{\rm left} W|}$ with the tensor factors enumerated by $\partial_{\rm left} W$, and likewise for $V^{\otimes \partial_{\rm right} W}$. Let $s:\partial W \to \{1,2,3\}$ be a state of $W$. Let
$$
\xi_{s;{\rm left}} := {\textstyle \bigotimes}_{x\in \partial_{\rm left}W} \xi_{s(x)} \in V^{\otimes \partial_{\rm left} W}, \qquad
\xi_{s;{\rm right}} := {\textstyle \bigotimes}_{x\in \partial_{\rm right}W} \xi_{s(x)} \in V^{\otimes \partial_{\rm right} W}.
$$
Now let $\rho_W$ be the unique $\mathbb{Z}[\omega^{\pm 1/2}]$-linear map satisfying
$$
\langle \xi_{s;{\rm left}}, \, \rho_W(\xi_{s;{\rm right}}) \rangle = {\rm Tr}^\omega_B([W,s])
$$
for all states $s$ of $W$. For this to be well-defined, one must check that the values of ${\rm Tr}^\omega_B([W,s])$ set by Prop.\ref{prop:biangle_SL3_quantum_trace} are invariant under the 180-degree rotation, which is easily verified. One way of expressing the above definition of $\rho_W$ is to say that its matrix elements are given by ${\rm Tr}^\omega_B([W,\cdot\,])$.

\begin{definition}
\label{def:composition_of_biangles}
Let $\vec{B}$ be a directed biangle. In the situation as in Prop.\ref{prop:biangle_SL3_quantum_trace}(BT1), the biangles $B_1,B_2$ naturally inherit the directions from that of $\vec{B}$. Suppose that the cutting arc $e$ is the right side of $\vec{B}_1$ and the left side of $\vec{B}_2$. 
After identifying the directed biangles $\vec{B}$, $\vec{B}_1$ and $\vec{B}_2$ with each other by direction-preserving diffeomorphisms, we say that the ${\rm SL}_3$-web $W$ in $\vec{B}\times {\bf I}$ is expressed as \ul{\em composition} of two ${\rm SL}_3$-webs $W_1$ and $W_2$ in $\vec{B}\times {\bf I}$. Write $W = W_1 \circ W_2$, allowing some ambiguity.
\end{definition}
One can easily observe that the property (BT1) of Prop.\ref{prop:biangle_SL3_quantum_trace} is encoding the matrix multiplication, or composition of operators, so that it is equivalent to
$$
\rho_{W_1 \circ W_2} = \rho_{W_1} \circ \rho_{W_2}.
$$
Property (BT2-1) says that $\rho_W = {\rm id}$ if $W$ falls into (BT2-1). Likewise, when $W$ is given by the product of single-edge ${\rm SL}_3$-webs as in (BT2-1), and if we choose an identification of $\partial_{\rm left}W$ and $\partial_{\rm right} W$  induced by the incidence relation coming from $W$, then $\rho_W = {\rm id}$. Now we are ready to consider some more elementary ${\rm SL}_3$-webs in a thickened biangle:
\begin{align}
\label{eq:some_more_elementary1}
& \raisebox{-0.5\height}{
\begingroup%
  \makeatletter%
  \providecommand\color[2][]{%
    \errmessage{(Inkscape) Color is used for the text in Inkscape, but the package 'color.sty' is not loaded}%
    \renewcommand\color[2][]{}%
  }%
  \providecommand\transparent[1]{%
    \errmessage{(Inkscape) Transparency is used (non-zero) for the text in Inkscape, but the package 'transparent.sty' is not loaded}%
    \renewcommand\transparent[1]{}%
  }%
  \providecommand\rotatebox[2]{#2}%
  \newcommand*\fsize{\dimexpr\f@size pt\relax}%
  \newcommand*\lineheight[1]{\fontsize{\fsize}{#1\fsize}\selectfont}%
  \ifx\svgwidth\undefined%
    \setlength{\unitlength}{56.69291339bp}%
    \ifx\svgscale\undefined%
      \relax%
    \else%
      \setlength{\unitlength}{\unitlength * \real{\svgscale}}%
    \fi%
  \else%
    \setlength{\unitlength}{\svgwidth}%
  \fi%
  \global\let\svgwidth\undefined%
  \global\let\svgscale\undefined%
  \makeatother%
  \begin{picture}(1,0.875)%
    \lineheight{1}%
    \setlength\tabcolsep{0pt}%
    \put(0,0){\includegraphics[width=\unitlength,page=1]{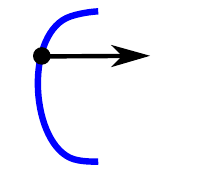}}%
    \put(0.00342324,0.57251168){\color[rgb]{0,0,0}\makebox(0,0)[lt]{\lineheight{1.25}\smash{\begin{tabular}[t]{l}$x_1$\end{tabular}}}}%
    \put(0,0){\includegraphics[width=\unitlength,page=2]{biangle_height_exchange2.pdf}}%
    \put(0.00801424,0.25622502){\color[rgb]{0,0,0}\makebox(0,0)[lt]{\lineheight{1.25}\smash{\begin{tabular}[t]{l}$x_2$\end{tabular}}}}%
    \put(0.05924823,0.51254405){\color[rgb]{0,0,0}\rotatebox{-92.070464}{\makebox(0,0)[lt]{\lineheight{1.25}\smash{\begin{tabular}[t]{l}$\succ$\end{tabular}}}}}%
    \put(0,0){\includegraphics[width=\unitlength,page=3]{biangle_height_exchange2.pdf}}%
    \put(0.85239343,0.58058091){\color[rgb]{0,0,0}\makebox(0,0)[lt]{\lineheight{1.25}\smash{\begin{tabular}[t]{l}$y_1$\end{tabular}}}}%
    \put(0.85698442,0.26429425){\color[rgb]{0,0,0}\makebox(0,0)[lt]{\lineheight{1.25}\smash{\begin{tabular}[t]{l}$y_2$\end{tabular}}}}%
    \put(0.89099145,0.50338616){\color[rgb]{0,0,0}\rotatebox{-92.070464}{\makebox(0,0)[lt]{\lineheight{1.25}\smash{\begin{tabular}[t]{l}$\prec$\end{tabular}}}}}%
  \end{picture}%
\endgroup%
} \qquad
\raisebox{-0.5\height}{
\begingroup%
  \makeatletter%
  \providecommand\color[2][]{%
    \errmessage{(Inkscape) Color is used for the text in Inkscape, but the package 'color.sty' is not loaded}%
    \renewcommand\color[2][]{}%
  }%
  \providecommand\transparent[1]{%
    \errmessage{(Inkscape) Transparency is used (non-zero) for the text in Inkscape, but the package 'transparent.sty' is not loaded}%
    \renewcommand\transparent[1]{}%
  }%
  \providecommand\rotatebox[2]{#2}%
  \newcommand*\fsize{\dimexpr\f@size pt\relax}%
  \newcommand*\lineheight[1]{\fontsize{\fsize}{#1\fsize}\selectfont}%
  \ifx\svgwidth\undefined%
    \setlength{\unitlength}{56.69291339bp}%
    \ifx\svgscale\undefined%
      \relax%
    \else%
      \setlength{\unitlength}{\unitlength * \real{\svgscale}}%
    \fi%
  \else%
    \setlength{\unitlength}{\svgwidth}%
  \fi%
  \global\let\svgwidth\undefined%
  \global\let\svgscale\undefined%
  \makeatother%
  \begin{picture}(1,0.875)%
    \lineheight{1}%
    \setlength\tabcolsep{0pt}%
    \put(0,0){\includegraphics[width=\unitlength,page=1]{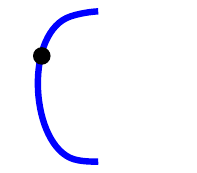}}%
    \put(-0.00171479,0.57550812){\color[rgb]{0,0,0}\makebox(0,0)[lt]{\lineheight{1.25}\smash{\begin{tabular}[t]{l}$x_2$\end{tabular}}}}%
    \put(0,0){\includegraphics[width=\unitlength,page=2]{biangle_crossing5.pdf}}%
    \put(0.00287621,0.25922147){\color[rgb]{0,0,0}\makebox(0,0)[lt]{\lineheight{1.25}\smash{\begin{tabular}[t]{l}$x_1$\end{tabular}}}}%
    \put(0.03772366,0.5151652){\color[rgb]{0,0,0}\rotatebox{-90.362304}{\makebox(0,0)[lt]{\lineheight{1.25}\smash{\begin{tabular}[t]{l}$\prec$\end{tabular}}}}}%
    \put(0,0){\includegraphics[width=\unitlength,page=3]{biangle_crossing5.pdf}}%
    \put(0.85239343,0.58058091){\color[rgb]{0,0,0}\makebox(0,0)[lt]{\lineheight{1.25}\smash{\begin{tabular}[t]{l}$y_1$\end{tabular}}}}%
    \put(0.85698442,0.26429425){\color[rgb]{0,0,0}\makebox(0,0)[lt]{\lineheight{1.25}\smash{\begin{tabular}[t]{l}$y_2$\end{tabular}}}}%
    \put(0.89099145,0.50338616){\color[rgb]{0,0,0}\rotatebox{-92.070464}{\makebox(0,0)[lt]{\lineheight{1.25}\smash{\begin{tabular}[t]{l}$\prec$\end{tabular}}}}}%
    \put(0,0){\includegraphics[width=\unitlength,page=4]{biangle_crossing5.pdf}}%
  \end{picture}%
\endgroup%
} \qquad
\raisebox{-0.5\height}{
\begingroup%
  \makeatletter%
  \providecommand\color[2][]{%
    \errmessage{(Inkscape) Color is used for the text in Inkscape, but the package 'color.sty' is not loaded}%
    \renewcommand\color[2][]{}%
  }%
  \providecommand\transparent[1]{%
    \errmessage{(Inkscape) Transparency is used (non-zero) for the text in Inkscape, but the package 'transparent.sty' is not loaded}%
    \renewcommand\transparent[1]{}%
  }%
  \providecommand\rotatebox[2]{#2}%
  \newcommand*\fsize{\dimexpr\f@size pt\relax}%
  \newcommand*\lineheight[1]{\fontsize{\fsize}{#1\fsize}\selectfont}%
  \ifx\svgwidth\undefined%
    \setlength{\unitlength}{56.69291339bp}%
    \ifx\svgscale\undefined%
      \relax%
    \else%
      \setlength{\unitlength}{\unitlength * \real{\svgscale}}%
    \fi%
  \else%
    \setlength{\unitlength}{\svgwidth}%
  \fi%
  \global\let\svgwidth\undefined%
  \global\let\svgscale\undefined%
  \makeatother%
  \begin{picture}(1,0.875)%
    \lineheight{1}%
    \setlength\tabcolsep{0pt}%
    \put(0,0){\includegraphics[width=\unitlength,page=1]{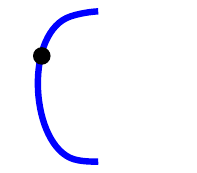}}%
    \put(-0.00171479,0.57550812){\color[rgb]{0,0,0}\makebox(0,0)[lt]{\lineheight{1.25}\smash{\begin{tabular}[t]{l}$x_1$\end{tabular}}}}%
    \put(0,0){\includegraphics[width=\unitlength,page=2]{biangle_crossing6.pdf}}%
    \put(0.00287621,0.25922147){\color[rgb]{0,0,0}\makebox(0,0)[lt]{\lineheight{1.25}\smash{\begin{tabular}[t]{l}$x_2$\end{tabular}}}}%
    \put(0.03772366,0.5151652){\color[rgb]{0,0,0}\rotatebox{-90.362304}{\makebox(0,0)[lt]{\lineheight{1.25}\smash{\begin{tabular}[t]{l}$\succ$\end{tabular}}}}}%
    \put(0,0){\includegraphics[width=\unitlength,page=3]{biangle_crossing6.pdf}}%
    \put(0.85239343,0.58058091){\color[rgb]{0,0,0}\makebox(0,0)[lt]{\lineheight{1.25}\smash{\begin{tabular}[t]{l}$y_2$\end{tabular}}}}%
    \put(0.85698442,0.26429425){\color[rgb]{0,0,0}\makebox(0,0)[lt]{\lineheight{1.25}\smash{\begin{tabular}[t]{l}$y_1$\end{tabular}}}}%
    \put(0.89099145,0.50338616){\color[rgb]{0,0,0}\rotatebox{-92.070464}{\makebox(0,0)[lt]{\lineheight{1.25}\smash{\begin{tabular}[t]{l}$\succ$\end{tabular}}}}}%
    \put(0,0){\includegraphics[width=\unitlength,page=4]{biangle_crossing6.pdf}}%
  \end{picture}%
\endgroup%
} \\
\label{eq:some_more_elementary2}
& \raisebox{-0.5\height}{
\begingroup%
  \makeatletter%
  \providecommand\color[2][]{%
    \errmessage{(Inkscape) Color is used for the text in Inkscape, but the package 'color.sty' is not loaded}%
    \renewcommand\color[2][]{}%
  }%
  \providecommand\transparent[1]{%
    \errmessage{(Inkscape) Transparency is used (non-zero) for the text in Inkscape, but the package 'transparent.sty' is not loaded}%
    \renewcommand\transparent[1]{}%
  }%
  \providecommand\rotatebox[2]{#2}%
  \newcommand*\fsize{\dimexpr\f@size pt\relax}%
  \newcommand*\lineheight[1]{\fontsize{\fsize}{#1\fsize}\selectfont}%
  \ifx\svgwidth\undefined%
    \setlength{\unitlength}{56.69291339bp}%
    \ifx\svgscale\undefined%
      \relax%
    \else%
      \setlength{\unitlength}{\unitlength * \real{\svgscale}}%
    \fi%
  \else%
    \setlength{\unitlength}{\svgwidth}%
  \fi%
  \global\let\svgwidth\undefined%
  \global\let\svgscale\undefined%
  \makeatother%
  \begin{picture}(1,0.875)%
    \lineheight{1}%
    \setlength\tabcolsep{0pt}%
    \put(0,0){\includegraphics[width=\unitlength,page=1]{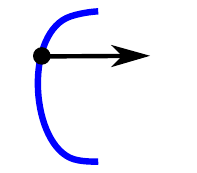}}%
    \put(0.00342324,0.57251168){\color[rgb]{0,0,0}\makebox(0,0)[lt]{\lineheight{1.25}\smash{\begin{tabular}[t]{l}$x_1$\end{tabular}}}}%
    \put(0,0){\includegraphics[width=\unitlength,page=2]{biangle_height_exchange4.pdf}}%
    \put(0.00801424,0.25622502){\color[rgb]{0,0,0}\makebox(0,0)[lt]{\lineheight{1.25}\smash{\begin{tabular}[t]{l}$x_2$\end{tabular}}}}%
    \put(0.05924823,0.51254405){\color[rgb]{0,0,0}\rotatebox{-92.070464}{\makebox(0,0)[lt]{\lineheight{1.25}\smash{\begin{tabular}[t]{l}$\succ$\end{tabular}}}}}%
    \put(0,0){\includegraphics[width=\unitlength,page=3]{biangle_height_exchange4.pdf}}%
    \put(0.85239343,0.58058091){\color[rgb]{0,0,0}\makebox(0,0)[lt]{\lineheight{1.25}\smash{\begin{tabular}[t]{l}$y_1$\end{tabular}}}}%
    \put(0.85698442,0.26429425){\color[rgb]{0,0,0}\makebox(0,0)[lt]{\lineheight{1.25}\smash{\begin{tabular}[t]{l}$y_2$\end{tabular}}}}%
    \put(0.89099145,0.50338616){\color[rgb]{0,0,0}\rotatebox{-92.070464}{\makebox(0,0)[lt]{\lineheight{1.25}\smash{\begin{tabular}[t]{l}$\prec$\end{tabular}}}}}%
  \end{picture}%
\endgroup%
}
\qquad
\raisebox{-0.5\height}{
\begingroup%
  \makeatletter%
  \providecommand\color[2][]{%
    \errmessage{(Inkscape) Color is used for the text in Inkscape, but the package 'color.sty' is not loaded}%
    \renewcommand\color[2][]{}%
  }%
  \providecommand\transparent[1]{%
    \errmessage{(Inkscape) Transparency is used (non-zero) for the text in Inkscape, but the package 'transparent.sty' is not loaded}%
    \renewcommand\transparent[1]{}%
  }%
  \providecommand\rotatebox[2]{#2}%
  \newcommand*\fsize{\dimexpr\f@size pt\relax}%
  \newcommand*\lineheight[1]{\fontsize{\fsize}{#1\fsize}\selectfont}%
  \ifx\svgwidth\undefined%
    \setlength{\unitlength}{56.69291339bp}%
    \ifx\svgscale\undefined%
      \relax%
    \else%
      \setlength{\unitlength}{\unitlength * \real{\svgscale}}%
    \fi%
  \else%
    \setlength{\unitlength}{\svgwidth}%
  \fi%
  \global\let\svgwidth\undefined%
  \global\let\svgscale\undefined%
  \makeatother%
  \begin{picture}(1,0.875)%
    \lineheight{1}%
    \setlength\tabcolsep{0pt}%
    \put(0,0){\includegraphics[width=\unitlength,page=1]{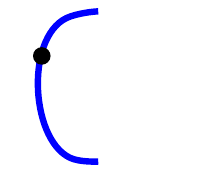}}%
    \put(-0.00171479,0.57550812){\color[rgb]{0,0,0}\makebox(0,0)[lt]{\lineheight{1.25}\smash{\begin{tabular}[t]{l}$x_2$\end{tabular}}}}%
    \put(0,0){\includegraphics[width=\unitlength,page=2]{biangle_crossing9.pdf}}%
    \put(0.00287621,0.25922147){\color[rgb]{0,0,0}\makebox(0,0)[lt]{\lineheight{1.25}\smash{\begin{tabular}[t]{l}$x_1$\end{tabular}}}}%
    \put(0.03772366,0.5151652){\color[rgb]{0,0,0}\rotatebox{-90.362304}{\makebox(0,0)[lt]{\lineheight{1.25}\smash{\begin{tabular}[t]{l}$\prec$\end{tabular}}}}}%
    \put(0,0){\includegraphics[width=\unitlength,page=3]{biangle_crossing9.pdf}}%
    \put(0.85239343,0.58058091){\color[rgb]{0,0,0}\makebox(0,0)[lt]{\lineheight{1.25}\smash{\begin{tabular}[t]{l}$y_1$\end{tabular}}}}%
    \put(0.85698442,0.26429425){\color[rgb]{0,0,0}\makebox(0,0)[lt]{\lineheight{1.25}\smash{\begin{tabular}[t]{l}$y_2$\end{tabular}}}}%
    \put(0.89099145,0.50338616){\color[rgb]{0,0,0}\rotatebox{-92.070464}{\makebox(0,0)[lt]{\lineheight{1.25}\smash{\begin{tabular}[t]{l}$\prec$\end{tabular}}}}}%
    \put(0,0){\includegraphics[width=\unitlength,page=4]{biangle_crossing9.pdf}}%
  \end{picture}%
\endgroup%
} \qquad
\raisebox{-0.5\height}{
\begingroup%
  \makeatletter%
  \providecommand\color[2][]{%
    \errmessage{(Inkscape) Color is used for the text in Inkscape, but the package 'color.sty' is not loaded}%
    \renewcommand\color[2][]{}%
  }%
  \providecommand\transparent[1]{%
    \errmessage{(Inkscape) Transparency is used (non-zero) for the text in Inkscape, but the package 'transparent.sty' is not loaded}%
    \renewcommand\transparent[1]{}%
  }%
  \providecommand\rotatebox[2]{#2}%
  \newcommand*\fsize{\dimexpr\f@size pt\relax}%
  \newcommand*\lineheight[1]{\fontsize{\fsize}{#1\fsize}\selectfont}%
  \ifx\svgwidth\undefined%
    \setlength{\unitlength}{56.69291339bp}%
    \ifx\svgscale\undefined%
      \relax%
    \else%
      \setlength{\unitlength}{\unitlength * \real{\svgscale}}%
    \fi%
  \else%
    \setlength{\unitlength}{\svgwidth}%
  \fi%
  \global\let\svgwidth\undefined%
  \global\let\svgscale\undefined%
  \makeatother%
  \begin{picture}(1,0.875)%
    \lineheight{1}%
    \setlength\tabcolsep{0pt}%
    \put(0,0){\includegraphics[width=\unitlength,page=1]{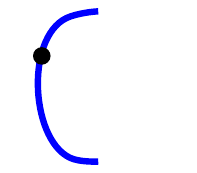}}%
    \put(-0.00171479,0.57550812){\color[rgb]{0,0,0}\makebox(0,0)[lt]{\lineheight{1.25}\smash{\begin{tabular}[t]{l}$x_1$\end{tabular}}}}%
    \put(0,0){\includegraphics[width=\unitlength,page=2]{biangle_crossing10.pdf}}%
    \put(0.00287621,0.25922147){\color[rgb]{0,0,0}\makebox(0,0)[lt]{\lineheight{1.25}\smash{\begin{tabular}[t]{l}$x_2$\end{tabular}}}}%
    \put(0.03772366,0.5151652){\color[rgb]{0,0,0}\rotatebox{-90.362304}{\makebox(0,0)[lt]{\lineheight{1.25}\smash{\begin{tabular}[t]{l}$\succ$\end{tabular}}}}}%
    \put(0,0){\includegraphics[width=\unitlength,page=3]{biangle_crossing10.pdf}}%
    \put(0.85239343,0.58058091){\color[rgb]{0,0,0}\makebox(0,0)[lt]{\lineheight{1.25}\smash{\begin{tabular}[t]{l}$y_2$\end{tabular}}}}%
    \put(0.85698442,0.26429425){\color[rgb]{0,0,0}\makebox(0,0)[lt]{\lineheight{1.25}\smash{\begin{tabular}[t]{l}$y_1$\end{tabular}}}}%
    \put(0.89099145,0.50338616){\color[rgb]{0,0,0}\rotatebox{-92.070464}{\makebox(0,0)[lt]{\lineheight{1.25}\smash{\begin{tabular}[t]{l}$\succ$\end{tabular}}}}}%
    \put(0,0){\includegraphics[width=\unitlength,page=4]{biangle_crossing10.pdf}}%
  \end{picture}%
\endgroup%
} \\
\label{eq:some_more_elementary3}
& \raisebox{-0.5\height}{
\begingroup%
  \makeatletter%
  \providecommand\color[2][]{%
    \errmessage{(Inkscape) Color is used for the text in Inkscape, but the package 'color.sty' is not loaded}%
    \renewcommand\color[2][]{}%
  }%
  \providecommand\transparent[1]{%
    \errmessage{(Inkscape) Transparency is used (non-zero) for the text in Inkscape, but the package 'transparent.sty' is not loaded}%
    \renewcommand\transparent[1]{}%
  }%
  \providecommand\rotatebox[2]{#2}%
  \newcommand*\fsize{\dimexpr\f@size pt\relax}%
  \newcommand*\lineheight[1]{\fontsize{\fsize}{#1\fsize}\selectfont}%
  \ifx\svgwidth\undefined%
    \setlength{\unitlength}{56.69291339bp}%
    \ifx\svgscale\undefined%
      \relax%
    \else%
      \setlength{\unitlength}{\unitlength * \real{\svgscale}}%
    \fi%
  \else%
    \setlength{\unitlength}{\svgwidth}%
  \fi%
  \global\let\svgwidth\undefined%
  \global\let\svgscale\undefined%
  \makeatother%
  \begin{picture}(1,0.875)%
    \lineheight{1}%
    \setlength\tabcolsep{0pt}%
    \put(0,0){\includegraphics[width=\unitlength,page=1]{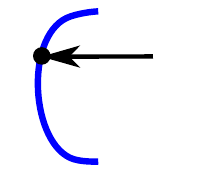}}%
    \put(0.00342324,0.57251168){\color[rgb]{0,0,0}\makebox(0,0)[lt]{\lineheight{1.25}\smash{\begin{tabular}[t]{l}$x_1$\end{tabular}}}}%
    \put(0,0){\includegraphics[width=\unitlength,page=2]{biangle_height_exchange6.pdf}}%
    \put(0.00801424,0.25622502){\color[rgb]{0,0,0}\makebox(0,0)[lt]{\lineheight{1.25}\smash{\begin{tabular}[t]{l}$x_2$\end{tabular}}}}%
    \put(0.05924823,0.51254405){\color[rgb]{0,0,0}\rotatebox{-92.070464}{\makebox(0,0)[lt]{\lineheight{1.25}\smash{\begin{tabular}[t]{l}$\succ$\end{tabular}}}}}%
    \put(0,0){\includegraphics[width=\unitlength,page=3]{biangle_height_exchange6.pdf}}%
    \put(0.85239343,0.58058091){\color[rgb]{0,0,0}\makebox(0,0)[lt]{\lineheight{1.25}\smash{\begin{tabular}[t]{l}$y_1$\end{tabular}}}}%
    \put(0.85698442,0.26429425){\color[rgb]{0,0,0}\makebox(0,0)[lt]{\lineheight{1.25}\smash{\begin{tabular}[t]{l}$y_2$\end{tabular}}}}%
    \put(0.89099145,0.50338616){\color[rgb]{0,0,0}\rotatebox{-92.070464}{\makebox(0,0)[lt]{\lineheight{1.25}\smash{\begin{tabular}[t]{l}$\prec$\end{tabular}}}}}%
  \end{picture}%
\endgroup%
} \qquad
\raisebox{-0.5\height}{
\begingroup%
  \makeatletter%
  \providecommand\color[2][]{%
    \errmessage{(Inkscape) Color is used for the text in Inkscape, but the package 'color.sty' is not loaded}%
    \renewcommand\color[2][]{}%
  }%
  \providecommand\transparent[1]{%
    \errmessage{(Inkscape) Transparency is used (non-zero) for the text in Inkscape, but the package 'transparent.sty' is not loaded}%
    \renewcommand\transparent[1]{}%
  }%
  \providecommand\rotatebox[2]{#2}%
  \newcommand*\fsize{\dimexpr\f@size pt\relax}%
  \newcommand*\lineheight[1]{\fontsize{\fsize}{#1\fsize}\selectfont}%
  \ifx\svgwidth\undefined%
    \setlength{\unitlength}{56.69291339bp}%
    \ifx\svgscale\undefined%
      \relax%
    \else%
      \setlength{\unitlength}{\unitlength * \real{\svgscale}}%
    \fi%
  \else%
    \setlength{\unitlength}{\svgwidth}%
  \fi%
  \global\let\svgwidth\undefined%
  \global\let\svgscale\undefined%
  \makeatother%
  \begin{picture}(1,0.875)%
    \lineheight{1}%
    \setlength\tabcolsep{0pt}%
    \put(0,0){\includegraphics[width=\unitlength,page=1]{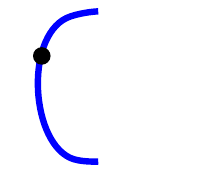}}%
    \put(-0.00171479,0.57550812){\color[rgb]{0,0,0}\makebox(0,0)[lt]{\lineheight{1.25}\smash{\begin{tabular}[t]{l}$x_2$\end{tabular}}}}%
    \put(0,0){\includegraphics[width=\unitlength,page=2]{biangle_crossing13.pdf}}%
    \put(0.00287621,0.25922147){\color[rgb]{0,0,0}\makebox(0,0)[lt]{\lineheight{1.25}\smash{\begin{tabular}[t]{l}$x_1$\end{tabular}}}}%
    \put(0.03772366,0.5151652){\color[rgb]{0,0,0}\rotatebox{-90.362304}{\makebox(0,0)[lt]{\lineheight{1.25}\smash{\begin{tabular}[t]{l}$\prec$\end{tabular}}}}}%
    \put(0,0){\includegraphics[width=\unitlength,page=3]{biangle_crossing13.pdf}}%
    \put(0.85239343,0.58058091){\color[rgb]{0,0,0}\makebox(0,0)[lt]{\lineheight{1.25}\smash{\begin{tabular}[t]{l}$y_1$\end{tabular}}}}%
    \put(0.85698442,0.26429425){\color[rgb]{0,0,0}\makebox(0,0)[lt]{\lineheight{1.25}\smash{\begin{tabular}[t]{l}$y_2$\end{tabular}}}}%
    \put(0.89099145,0.50338616){\color[rgb]{0,0,0}\rotatebox{-92.070464}{\makebox(0,0)[lt]{\lineheight{1.25}\smash{\begin{tabular}[t]{l}$\prec$\end{tabular}}}}}%
    \put(0,0){\includegraphics[width=\unitlength,page=4]{biangle_crossing13.pdf}}%
  \end{picture}%
\endgroup%
} \qquad
\raisebox{-0.5\height}{
\begingroup%
  \makeatletter%
  \providecommand\color[2][]{%
    \errmessage{(Inkscape) Color is used for the text in Inkscape, but the package 'color.sty' is not loaded}%
    \renewcommand\color[2][]{}%
  }%
  \providecommand\transparent[1]{%
    \errmessage{(Inkscape) Transparency is used (non-zero) for the text in Inkscape, but the package 'transparent.sty' is not loaded}%
    \renewcommand\transparent[1]{}%
  }%
  \providecommand\rotatebox[2]{#2}%
  \newcommand*\fsize{\dimexpr\f@size pt\relax}%
  \newcommand*\lineheight[1]{\fontsize{\fsize}{#1\fsize}\selectfont}%
  \ifx\svgwidth\undefined%
    \setlength{\unitlength}{56.69291339bp}%
    \ifx\svgscale\undefined%
      \relax%
    \else%
      \setlength{\unitlength}{\unitlength * \real{\svgscale}}%
    \fi%
  \else%
    \setlength{\unitlength}{\svgwidth}%
  \fi%
  \global\let\svgwidth\undefined%
  \global\let\svgscale\undefined%
  \makeatother%
  \begin{picture}(1,0.875)%
    \lineheight{1}%
    \setlength\tabcolsep{0pt}%
    \put(0,0){\includegraphics[width=\unitlength,page=1]{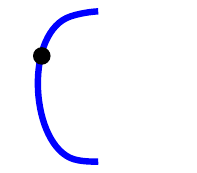}}%
    \put(-0.00171479,0.57550812){\color[rgb]{0,0,0}\makebox(0,0)[lt]{\lineheight{1.25}\smash{\begin{tabular}[t]{l}$x_1$\end{tabular}}}}%
    \put(0,0){\includegraphics[width=\unitlength,page=2]{biangle_crossing14.pdf}}%
    \put(0.00287621,0.25922147){\color[rgb]{0,0,0}\makebox(0,0)[lt]{\lineheight{1.25}\smash{\begin{tabular}[t]{l}$x_2$\end{tabular}}}}%
    \put(0.03772366,0.5151652){\color[rgb]{0,0,0}\rotatebox{-90.362304}{\makebox(0,0)[lt]{\lineheight{1.25}\smash{\begin{tabular}[t]{l}$\succ$\end{tabular}}}}}%
    \put(0,0){\includegraphics[width=\unitlength,page=3]{biangle_crossing14.pdf}}%
    \put(0.85239343,0.58058091){\color[rgb]{0,0,0}\makebox(0,0)[lt]{\lineheight{1.25}\smash{\begin{tabular}[t]{l}$y_2$\end{tabular}}}}%
    \put(0.85698442,0.26429425){\color[rgb]{0,0,0}\makebox(0,0)[lt]{\lineheight{1.25}\smash{\begin{tabular}[t]{l}$y_1$\end{tabular}}}}%
    \put(0.89099145,0.50338616){\color[rgb]{0,0,0}\rotatebox{-92.070464}{\makebox(0,0)[lt]{\lineheight{1.25}\smash{\begin{tabular}[t]{l}$\succ$\end{tabular}}}}}%
    \put(0,0){\includegraphics[width=\unitlength,page=4]{biangle_crossing14.pdf}}%
  \end{picture}%
\endgroup%
}
\end{align}
In each of eq.\eqref{eq:some_more_elementary1}--\eqref{eq:some_more_elementary3}, the three ${\rm SL}_3$-webs are isotopic. One observes that the first ${\rm SL}_3$-web in eq.\eqref{eq:some_more_elementary1}, say $W_1$, is `inverse' of the first ${\rm SL}_3$-web appearing in eq.\eqref{eq:height_exchange_and_crossing1}, say $W_2$, in the sense that $W_1 \circ W_2$ is isotopic to a product of two single-edge ${\rm SL}_3$-webs as in (BT2-1), so $\rho_{W_1} \circ \rho_{W_2} = {\rm id}$ (we are viewing the upper marked point as the top marked point, and we need some re-labeling of endpoints, to avoid confusion). Similarly, the first ${\rm SL}_3$-web in eq.\eqref{eq:some_more_elementary2} is inverse of the first ${\rm SL}_3$-web appearing in eq.\eqref{eq:height_exchange_and_crossing3}, and the first ${\rm SL}_3$-web in eq.\eqref{eq:some_more_elementary3} is inverse of the first ${\rm SL}_3$-web appearing in eq.\eqref{eq:height_exchange_and_crossing5}. Thus the values of ${\rm Tr}^\omega_B$ on the stated ${\rm SL}_3$-webs based on eq.\eqref{eq:some_more_elementary1}--\eqref{eq:some_more_elementary3} are determined by properties of Prop.\ref{prop:biangle_SL3_quantum_trace}. Now we need the following definition and observation, which will be used again later:
\begin{definition}
\label{def:elementary_A2-webs_in_biangle}
Let $\vec{B}$ be a directed biangle. An ${\rm SL}_3$-web in $\vec{B}\times {\bf I}$ is \ul{\em elementary} if it is a product of some (possibly zero) number of ${\rm SL}_3$-webs of type (BT2-1) and at most one ${\rm SL}_3$-web of type among (BT2-2)--(BT2-4) or among eq.\eqref{eq:some_more_elementary1}--\eqref{eq:some_more_elementary3}, where a crossing can occur at most once in total, which occurs at the part involving (BT2-2)--(BT2-4) or eq.\eqref{eq:some_more_elementary1}--\eqref{eq:some_more_elementary3}.
\end{definition}
\begin{lemma}
\label{lem:composition_of_elementary_biangles}
Any ${\rm SL}_3$-web in a thickened biangle $\vec{B}\times {\bf I}$ is composition of elementary ${\rm SL}_3$-webs in $\vec{B}\times {\bf I}$. \qed
\end{lemma}
By our observations above and by Prop.\ref{prop:biangle_SL3_quantum_trace_some_values} which is to be proved, the properties of Prop.\ref{prop:biangle_SL3_quantum_trace} determine the values ${\rm Tr}^\omega_B$ at elementary stated ${\rm SL}_3$-webs. Hence, by the composition property (BT1), together with Lem.\ref{lem:composition_of_elementary_biangles}, all values of ${\rm Tr}^\omega_B$ are determined by the properties of Prop.\ref{prop:biangle_SL3_quantum_trace}, finishing the proof of  the uniqueness part of Prop.\ref{prop:biangle_SL3_quantum_trace}.

\vs

There are several strategies to prove the existence part of Prop.\ref{prop:biangle_SL3_quantum_trace}. We take possibly the shortest one, by mimicking the idea of Costantino and L\^e for the ${\rm SL}_2$ case \cite{Le18} \cite{CL}, applied to the ${\rm SL}_3$ case with the help of the results of Higgins on the stated ${\rm SL}_3$-skein algebras \cite{Higgins}. Namely, it is proved in \cite{Higgins} that the reduced stated ${\rm SL}_3$-skein algebra $\mathcal{S}^\omega_{\rm s}(\vec{B};\mathbb{Z})_{\rm red}$ for a directed biangle is a Hopf algebra and is \redfix{isomorphic} to the quantum group $\mathcal{O}_q({\rm SL}_3)$. We will just use the fact that $\mathcal{S}^\omega_{\rm s}(\vec{B};\mathbb{Z})_{\rm red}$ is a bialgebra. We already know its algebra structure, i.e. the product. The coproduct comes from the cutting/gluing process.
\begin{definition}[coproduct for the stated ${\rm SL}_3$-skein algebra for a directed biangle; \cite{Higgins}]
Let $\vec{B}=(B,{\rm dir})$ be a directed biangle, and $e$ be an ideal arc in $B$ whose interior lies in the interior of $B$. Cutting $B$ along $e$ yields disjoint union of two biangles $B_1$ and $B_2$ as in Prop.\ref{prop:biangle_SL3_quantum_trace}(BT1). The direction naturally inherits to $B_1$ and $B_2$, making them directed biangles $\vec{B}_1$ and $\vec{B}_2$; assume that the left side of $\vec{B}_1$ corresponds to the left side of $\vec{B}$, and the right side of $\vec{B}_2$ to the right side of $\vec{B}$. Let $(W,s)$ be a stated ${\rm SL}_3$-web in $B\times {\bf I}$, and let $W_1$ and $W_2$ the ${\rm SL}_3$-webs in $B_1\times {\bf I}$ and $B_2\times{\bf I}$ obtained by cutting $W$ along $e$, as in Prop.\ref{prop:biangle_SL3_quantum_trace}(BT1). Define the map
\begin{align*}
& \Delta_{\vec{B},e} : \mathcal{S}^\omega_{\rm s}(\vec{B};\mathbb{Z})_{\rm red} \longrightarrow \mathcal{S}^\omega_{\rm s}(\vec{B}_1;\mathbb{Z})_{\rm red} \otimes_{\mathbb{Z}[\omega^{\pm 1/2}]} \mathcal{S}^\omega_{\rm s}(\vec{B}_2;\mathbb{Z})_{\rm red} \qquad \mbox{as} \\
& \Delta_{\vec{B},e} ( [W,s] ) := \underset{s_1,s_2}{\textstyle \sum} [W_1,s_1] \otimes [W_2,s_2],
\end{align*}
where the sum is over all states $s_1$ and $s_2$ of $W_1$ and $W_2$ such that the state $s' := s_1 \cup s_2$ of $W' = W_1 \cup W_2$ is compatible with $s$ in the sense as in Lem.\ref{lem:cutting_process}. Composing with the canonical isomorphisms $\mathcal{S}^\omega_{\rm s}(\vec{B}_i;\mathbb{Z})_{\rm red}  \cong \mathcal{S}^\omega_{\rm s}(\vec{B};\mathbb{Z})_{\rm red}$, define the map
$$
\Delta : \mathcal{S}^\omega_{\rm s}(\vec{B};\mathbb{Z})_{\rm red} \longrightarrow \mathcal{S}^\omega_{\rm s}(\vec{B};\mathbb{Z})_{\rm red} \otimes_{\mathbb{Z}[\omega^{\pm 1/2}]} \mathcal{S}^\omega_{\rm s}(\vec{B};\mathbb{Z})_{\rm red}.
$$
\end{definition}
We note that, in order for the above form of coproduct in \cite{Higgins} to translate via the isomorphism in eq.\eqref{eq:isomorphism_from_ours_to_Higgins} neatly as above, the twisting scalars in eq.\eqref{eq:twisting_scalars} must satisfy $\alpha_3\alpha_4 = 1$, which indeed holds for our choice of twisting scalars.

\vs

The counit $\epsilon : \mathcal{S}^\omega_{\rm s}(\vec{B};\mathbb{Z}) \to \mathbb{Z}[\omega^{\pm 1/2}]$ is constructed in \cite{Higgins} by first applying the `inversion' along the right side of $\vec{B}$ which has an effect of reversing the orientation on that side, then `fill in' the top marked point to obtain a stated ${\rm SL}_3$-skein over the monogon, i.e. an element of $\mathcal{S}^\omega_{\rm s}(M;\mathbb{Z})$; then use the isomorphism $\mathcal{S}^\omega_{\rm s}(M;\mathbb{Z}) \cong \mathbb{Z}[\omega^{1/2}]$ (\cite[Prop.3]{Higgins}). In the present paper, instead of recalling the detailed construction of $\epsilon$, it suffices to just know that it is a counit, and that for the stated ${\rm SL}_3$-skeins of type (BT2-1) its values are same as those under ${\rm Tr}^\omega_B$ as written in Prop.\ref{prop:biangle_SL3_quantum_trace}(BT2-1); in particular, these properties completely determine $\epsilon$.

\begin{proposition}[reduced stated ${\rm SL}_3$-skein algebra for a biangle is a bialgebra; \cite{Higgins}]
\label{prop:reduced_stated_A2-skein_algebra_for_biangle_is_a_bialgebra}
There exists a map $\epsilon : \mathcal{S}^\omega_{\rm s}(\vec{B};\mathbb{Z}) \to \mathbb{Z}[\omega^{\pm 1/2}]$ such that 
\begin{enumerate}
\item[\rm (CU1)] the product of $\mathcal{S}^\omega_{\rm s}(\vec{B};\mathbb{Z})_{\rm red}$, the unit of $\mathcal{S}^\omega_{\rm s}(\vec{B};\mathbb{Z})_{\rm red}$, the maps $\Delta$ and $\epsilon$ make $\mathcal{S}^\omega_{\rm s}(\vec{B};\mathbb{Z})_{\rm red}$ a well-defined bialgebra over $\mathbb{Z}[\omega^{\pm 1/2}]$, where $\Delta$ is the coproduct and $\epsilon$ is the counit;

\item[\rm (CU2)] If $(W,s)$ is as in Prop.\ref{prop:biangle_SL3_quantum_trace}(BT2-1), then
$$
\epsilon([W,s]) = \mbox{$(\varepsilon,\varepsilon')$-th entry of the $3\times 3$ identity matrix}.
$$
\end{enumerate}
\end{proposition}
As said above, it is not hard to see that such $\epsilon$ is unique, e.g. by a similar argument as in the proof of the uniqueness part of Prop.\ref{prop:biangle_SL3_quantum_trace}.

\vs

{\it Proof of Prop.\ref{prop:biangle_SL3_quantum_trace}.} 
We only need to show the existence. We will prove that $\epsilon$ is the sought-for map ${\rm Tr}^\omega_B$. One thing to note is that while the sought-for map ${\rm Tr}^\omega_B$ should be defined for a biangle $B$, the map $\epsilon$ is defined for a directed biangle $\vec{B}$. Given a biangle $B$, give an arbitrary direction to make it a directed biangle $\vec{B}$. As a result of our proof, we shall observe that the values of $\epsilon$ do not depend on this choice of direction on $B$. 

\vs

Since $\epsilon$ is a counit of a $\mathbb{Z}[\omega^{\pm 1/2}]$-bialgebra, it respects the product structure and hence is a $\mathbb{Z}[\omega^{\pm 1/2}]$-algebra homomorphism, and also satisfies $(\epsilon \otimes {\rm id})\circ \Delta = {\rm id} : \mathcal{S}^\omega_{\rm s}(\vec{B};\mathbb{Z})_{\rm red} \to \mathcal{S}^\omega_{\rm s}(\vec{B};\mathbb{Z})_{\rm red}$, thus $(\epsilon \otimes \epsilon) \circ \Delta = \epsilon : \mathcal{S}^\omega_{\rm s}(\vec{B};\mathbb{Z})_{\rm red} \to \mathbb{Z}[\omega^{\pm 1/2}]$ holds, which is precisely the property (BT1) for $\epsilon$. 

\vs

What remains is to check (BT2) for $\epsilon$, more precisely (BT2-1), (BT2-2), (BT2-3). Note that (BT2-1) is satisfied, by Prop.\ref{prop:reduced_stated_A2-skein_algebra_for_biangle_is_a_bialgebra}(CU2). For (BT2-2), it is convenient to make use of part of Prop.1 of \cite{Higgins}, which lists some consequences of the defining relations of the reduced stated ${\rm SL}_3$-skein algebras.
\begin{lemma}[part of {\cite[Prop.1]{Higgins}}] 
\label{lem:Higgins_Prop1}
In the version of reduced ${\rm SL}_3$-skein algebra 
$\mathcal{S}^{SL_3}_q(\Sigma)$ of \cite{Higgins} for a generalized marked surface $(\Sigma,\mathcal{P})$ and a commutative ring $\mathcal{R}$ with unity, one has the following: (blue horizontal line is boundary)
\begin{align}
\label{eq:Higgins_U}
& \raisebox{-0.6\height}{}
= - q^{-4/3} \delta_{s(x_1)+s(x_2), \, 4} \cdot \raisebox{-0.6\height}{
\begingroup%
  \makeatletter%
  \providecommand\color[2][]{%
    \errmessage{(Inkscape) Color is used for the text in Inkscape, but the package 'color.sty' is not loaded}%
    \renewcommand\color[2][]{}%
  }%
  \providecommand\transparent[1]{%
    \errmessage{(Inkscape) Transparency is used (non-zero) for the text in Inkscape, but the package 'transparent.sty' is not loaded}%
    \renewcommand\transparent[1]{}%
  }%
  \providecommand\rotatebox[2]{#2}%
  \newcommand*\fsize{\dimexpr\f@size pt\relax}%
  \newcommand*\lineheight[1]{\fontsize{\fsize}{#1\fsize}\selectfont}%
  \ifx\svgwidth\undefined%
    \setlength{\unitlength}{42.51968504bp}%
    \ifx\svgscale\undefined%
      \relax%
    \else%
      \setlength{\unitlength}{\unitlength * \real{\svgscale}}%
    \fi%
  \else%
    \setlength{\unitlength}{\svgwidth}%
  \fi%
  \global\let\svgwidth\undefined%
  \global\let\svgscale\undefined%
  \makeatother%
  \begin{picture}(1,0.83333333)%
    \lineheight{1}%
    \setlength\tabcolsep{0pt}%
    \put(0,0){\includegraphics[width=\unitlength,page=1]{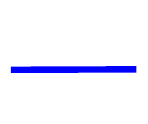}}%
  \end{picture}%
\endgroup%
}, \qquad \mbox{(in $\mathcal{S}^{SL_3}_q(\Sigma)$ of \cite{Higgins})} \\
\label{eq:Higgins_U2}
& \raisebox{-0.6\height}{}
= - q^{4/3} q^{2(s(x_2)-2)} \delta_{s(x_1)+s(x_2), \, 4} \cdot \raisebox{-0.6\height}{},  \qquad \mbox{(in $\mathcal{S}^{SL_3}_q(\Sigma)$ of \cite{Higgins})} \\
\label{eq:Higgins_U3}
& \raisebox{-0.6\height}{}
= - q^{-4/3} q^{2(s(x_1)-2)} \delta_{s(x_1)+s(x_2), \, 4} \cdot \raisebox{-0.6\height}{},  \qquad \mbox{(in $\mathcal{S}^{SL_3}_q(\Sigma)$ of \cite{Higgins})} \\
\label{eq:Higgins_U4}
& \raisebox{-0.6\height}{}
= - q^{4/3} \delta_{s(x_1)+s(x_2), \, 4} \cdot \raisebox{-0.6\height}{},  \qquad \mbox{(in $\mathcal{S}^{SL_3}_q(\Sigma)$ of \cite{Higgins})}.
\end{align}
\end{lemma}
In particular, the right hand sides are nonzero iff $(s(x_1),s(x_2)) \in \{(1,3),(2,2),(3,1)\}$. Let's apply Lem.\ref{lem:Higgins_Prop1} to (BT2-2). Suppose $[W,s] \in \mathcal{S}^\omega_{\rm s}(\vec{B};\mathbb{Z})_{\rm red}$ is as in eq.\eqref{eq:Higgins_U}, eq.\eqref{eq:Higgins_U2}, eq.\eqref{eq:Higgins_U3}, or eq.\eqref{eq:Higgins_U4}. As in Lem.\ref{lem:Higgins_Prop1}, if $(s(x_1),s(x_2))=(1,3),(2,2),(3,1)$ respectively, the image in $\mathcal{S}^{SL_3}_{q}(\Sigma)$ of this element under the isomorphism in eq.\eqref{eq:isomorphism_from_ours_to_Higgins} equals identity times $(-q^{-1})(-q^{-4/3}),-q^{-4/3},(-q)(-q^{-4/3})$ for eq.\eqref{eq:Higgins_U}, $(-q^{-1})(-q^{10/3}), -q^{4/3},(-q)(-q^{-2/3})$ for eq.\eqref{eq:Higgins_U2}, $(-q)(-q^{-10/3}), -q^{-4/3}, (-q^{-1})(-q^{2/3})$ for eq.\eqref{eq:Higgins_U3}, and $(-q)(-q^{4/3}),-q^{4/3},(-q^{-1})(-q^{4/3})$ for eq.\eqref{eq:Higgins_U4}, while for other pairs of $(s(x_1),s(x_2))$ the element is zero. Applying the inverse of the isomorphism in eq.\eqref{eq:isomorphism_from_ours_to_Higgins}, if $(s(x_1),s(x_2))=(1,3),(2,2),(3,1)$ respectively, we see that $[W,s] \in \mathcal{S}^\omega_{\rm s}(\vec{B};\mathbb{Z})_{\rm red}$ equals identity times $q^{7/3},-q^{4/3},q^{1/3}$ for eq.\eqref{eq:Higgins_U}, $q^{-7/3},-q^{-4/3},q^{-1/3}$ for eq.\eqref{eq:Higgins_U2}, $q^{7/3},-q^{4/3},q^{1/3}$ for eq.\eqref{eq:Higgins_U3}, and $q^{-7/3},-q^{-4/3},q^{-1/3}$ for eq.\eqref{eq:Higgins_U4}, while for other pairs of $(s(x_1),s(x_2))$ the element is zero. This shows (BT2-2) for $\epsilon$.

\vs

For (BT2-3) and for later use, we collect some more useful relations for the stated ${\rm SL}_3$-skein algebras (in the following five lemmas, the endpoints appearing in each figure have consecutive elevation orderings, in the sense as explained in Fig.\ref{fig:stated_boundary_relations}).
\begin{lemma}
\label{lem:3-way_boundary_relations1}
In $\mathcal{S}^\omega_{\rm s}(\frak{S};\mathcal{R})$ (hence also in $\mathcal{S}^\omega_{\rm s}(\frak{S};\mathcal{R})_{\rm red}$) for a generalized marked surface $\frak{S}$, one has:
\begin{align}
\label{eq:outgoing_fork_relations1}
\raisebox{-0.6\height} {}
& = \raisebox{-0.6\height} {
\begingroup%
  \makeatletter%
  \providecommand\color[2][]{%
    \errmessage{(Inkscape) Color is used for the text in Inkscape, but the package 'color.sty' is not loaded}%
    \renewcommand\color[2][]{}%
  }%
  \providecommand\transparent[1]{%
    \errmessage{(Inkscape) Transparency is used (non-zero) for the text in Inkscape, but the package 'transparent.sty' is not loaded}%
    \renewcommand\transparent[1]{}%
  }%
  \providecommand\rotatebox[2]{#2}%
  \newcommand*\fsize{\dimexpr\f@size pt\relax}%
  \newcommand*\lineheight[1]{\fontsize{\fsize}{#1\fsize}\selectfont}%
  \ifx\svgwidth\undefined%
    \setlength{\unitlength}{42.51968504bp}%
    \ifx\svgscale\undefined%
      \relax%
    \else%
      \setlength{\unitlength}{\unitlength * \real{\svgscale}}%
    \fi%
  \else%
    \setlength{\unitlength}{\svgwidth}%
  \fi%
  \global\let\svgwidth\undefined%
  \global\let\svgscale\undefined%
  \makeatother%
  \begin{picture}(1,1.16666667)%
    \lineheight{1}%
    \setlength\tabcolsep{0pt}%
    \put(0,0){\includegraphics[width=\unitlength,page=1]{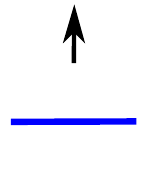}}%
    \put(0.16942518,0.13006975){\color[rgb]{0,0,0}\makebox(0,0)[lt]{\lineheight{1.25}\smash{\begin{tabular}[t]{l}$x_2$\end{tabular}}}}%
    \put(0,0){\includegraphics[width=\unitlength,page=2]{boundary_rel2_prove2_quantum.pdf}}%
    \put(0.74016697,0.12462211){\color[rgb]{0,0,0}\makebox(0,0)[lt]{\lineheight{1.25}\smash{\begin{tabular}[t]{l}$x_1$\end{tabular}}}}%
    \put(0.4330944,0.12306698){\color[rgb]{0,0,0}\makebox(0,0)[lt]{\lineheight{1.25}\smash{\begin{tabular}[t]{l}$\succ$\end{tabular}}}}%
    \put(0,0){\includegraphics[width=\unitlength,page=3]{boundary_rel2_prove2_quantum.pdf}}%
  \end{picture}%
\endgroup%
} 
\overset{{\rm (S8)}}{=} q^{-2/3} \hspace{-3mm} \raisebox{-0.6\height} {
\begingroup%
  \makeatletter%
  \providecommand\color[2][]{%
    \errmessage{(Inkscape) Color is used for the text in Inkscape, but the package 'color.sty' is not loaded}%
    \renewcommand\color[2][]{}%
  }%
  \providecommand\transparent[1]{%
    \errmessage{(Inkscape) Transparency is used (non-zero) for the text in Inkscape, but the package 'transparent.sty' is not loaded}%
    \renewcommand\transparent[1]{}%
  }%
  \providecommand\rotatebox[2]{#2}%
  \newcommand*\fsize{\dimexpr\f@size pt\relax}%
  \newcommand*\lineheight[1]{\fontsize{\fsize}{#1\fsize}\selectfont}%
  \ifx\svgwidth\undefined%
    \setlength{\unitlength}{42.51968504bp}%
    \ifx\svgscale\undefined%
      \relax%
    \else%
      \setlength{\unitlength}{\unitlength * \real{\svgscale}}%
    \fi%
  \else%
    \setlength{\unitlength}{\svgwidth}%
  \fi%
  \global\let\svgwidth\undefined%
  \global\let\svgscale\undefined%
  \makeatother%
  \begin{picture}(1,1.16666667)%
    \lineheight{1}%
    \setlength\tabcolsep{0pt}%
    \put(0,0){\includegraphics[width=\unitlength,page=1]{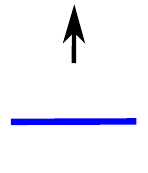}}%
    \put(0.16942518,0.13006975){\color[rgb]{0,0,0}\makebox(0,0)[lt]{\lineheight{1.25}\smash{\begin{tabular}[t]{l}$x_2$\end{tabular}}}}%
    \put(0,0){\includegraphics[width=\unitlength,page=2]{boundary_rel2_prime_quantum.pdf}}%
    \put(0.74016697,0.12462211){\color[rgb]{0,0,0}\makebox(0,0)[lt]{\lineheight{1.25}\smash{\begin{tabular}[t]{l}$x_1$\end{tabular}}}}%
    \put(0.50364996,0.12306698){\color[rgb]{0,0,0}\makebox(0,0)[lt]{\lineheight{1.25}\smash{\begin{tabular}[t]{l}$\succ$\end{tabular}}}}%
  \end{picture}%
\endgroup%
}
+ q^{1/3} \hspace{-3mm} \raisebox{-0.6\height} {
\begingroup%
  \makeatletter%
  \providecommand\color[2][]{%
    \errmessage{(Inkscape) Color is used for the text in Inkscape, but the package 'color.sty' is not loaded}%
    \renewcommand\color[2][]{}%
  }%
  \providecommand\transparent[1]{%
    \errmessage{(Inkscape) Transparency is used (non-zero) for the text in Inkscape, but the package 'transparent.sty' is not loaded}%
    \renewcommand\transparent[1]{}%
  }%
  \providecommand\rotatebox[2]{#2}%
  \newcommand*\fsize{\dimexpr\f@size pt\relax}%
  \newcommand*\lineheight[1]{\fontsize{\fsize}{#1\fsize}\selectfont}%
  \ifx\svgwidth\undefined%
    \setlength{\unitlength}{42.51968504bp}%
    \ifx\svgscale\undefined%
      \relax%
    \else%
      \setlength{\unitlength}{\unitlength * \real{\svgscale}}%
    \fi%
  \else%
    \setlength{\unitlength}{\svgwidth}%
  \fi%
  \global\let\svgwidth\undefined%
  \global\let\svgscale\undefined%
  \makeatother%
  \begin{picture}(1,1.16666667)%
    \lineheight{1}%
    \setlength\tabcolsep{0pt}%
    \put(0,0){\includegraphics[width=\unitlength,page=1]{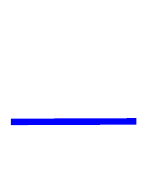}}%
    \put(0.16942518,0.13006975){\color[rgb]{0,0,0}\makebox(0,0)[lt]{\lineheight{1.25}\smash{\begin{tabular}[t]{l}$x_2$\end{tabular}}}}%
    \put(0,0){\includegraphics[width=\unitlength,page=2]{boundary_rel2_prove3_quantum.pdf}}%
    \put(0.74016697,0.12462211){\color[rgb]{0,0,0}\makebox(0,0)[lt]{\lineheight{1.25}\smash{\begin{tabular}[t]{l}$x_1$\end{tabular}}}}%
    \put(0.4330944,0.12306698){\color[rgb]{0,0,0}\makebox(0,0)[lt]{\lineheight{1.25}\smash{\begin{tabular}[t]{l}$\succ$\end{tabular}}}}%
    \put(0,0){\includegraphics[width=\unitlength,page=3]{boundary_rel2_prove3_quantum.pdf}}%
  \end{picture}%
\endgroup%
} 
\overset{{\rm (S6)}}{=} (\underbrace{ q^{-2/3} - q^{1/3} [2]_q }_{= \,\, - q^{4/3}}) \hspace{-2mm} \raisebox{-0.6\height} {}
\end{align}
\end{lemma}

\begin{lemma}
\label{lem:3-way_boundary_relations2}
In $\mathcal{S}^\omega_{\rm s}(\frak{S};\mathcal{R})_{\rm red}$ for a generalized marked surface $\frak{S}$, one has:
\begin{align}
\label{eq:outgoing_fork_relations2}
\raisebox{-0.6\height} {}
& \overset{{\rm (B2)}}{=} q \hspace{-1,5mm} \raisebox{-0.6\height} {
\begingroup%
  \makeatletter%
  \providecommand\color[2][]{%
    \errmessage{(Inkscape) Color is used for the text in Inkscape, but the package 'color.sty' is not loaded}%
    \renewcommand\color[2][]{}%
  }%
  \providecommand\transparent[1]{%
    \errmessage{(Inkscape) Transparency is used (non-zero) for the text in Inkscape, but the package 'transparent.sty' is not loaded}%
    \renewcommand\transparent[1]{}%
  }%
  \providecommand\rotatebox[2]{#2}%
  \newcommand*\fsize{\dimexpr\f@size pt\relax}%
  \newcommand*\lineheight[1]{\fontsize{\fsize}{#1\fsize}\selectfont}%
  \ifx\svgwidth\undefined%
    \setlength{\unitlength}{42.51968504bp}%
    \ifx\svgscale\undefined%
      \relax%
    \else%
      \setlength{\unitlength}{\unitlength * \real{\svgscale}}%
    \fi%
  \else%
    \setlength{\unitlength}{\svgwidth}%
  \fi%
  \global\let\svgwidth\undefined%
  \global\let\svgscale\undefined%
  \makeatother%
  \begin{picture}(1,1.16666667)%
    \lineheight{1}%
    \setlength\tabcolsep{0pt}%
    \put(0,0){\includegraphics[width=\unitlength,page=1]{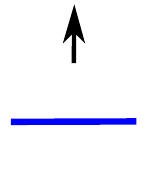}}%
    \put(0.16942518,0.13006975){\color[rgb]{0,0,0}\makebox(0,0)[lt]{\lineheight{1.25}\smash{\begin{tabular}[t]{l}$x_2$\end{tabular}}}}%
    \put(0,0){\includegraphics[width=\unitlength,page=2]{boundary_rel2_prime2_quantum.pdf}}%
    \put(0.74016697,0.12462211){\color[rgb]{0,0,0}\makebox(0,0)[lt]{\lineheight{1.25}\smash{\begin{tabular}[t]{l}$x_1$\end{tabular}}}}%
    \put(0.4330944,0.12306698){\color[rgb]{0,0,0}\makebox(0,0)[lt]{\lineheight{1.25}\smash{\begin{tabular}[t]{l}$\prec$\end{tabular}}}}%
  \end{picture}%
\endgroup%
} +  \hspace{-1,5mm} \raisebox{-0.6\height} {
\begingroup%
  \makeatletter%
  \providecommand\color[2][]{%
    \errmessage{(Inkscape) Color is used for the text in Inkscape, but the package 'color.sty' is not loaded}%
    \renewcommand\color[2][]{}%
  }%
  \providecommand\transparent[1]{%
    \errmessage{(Inkscape) Transparency is used (non-zero) for the text in Inkscape, but the package 'transparent.sty' is not loaded}%
    \renewcommand\transparent[1]{}%
  }%
  \providecommand\rotatebox[2]{#2}%
  \newcommand*\fsize{\dimexpr\f@size pt\relax}%
  \newcommand*\lineheight[1]{\fontsize{\fsize}{#1\fsize}\selectfont}%
  \ifx\svgwidth\undefined%
    \setlength{\unitlength}{42.51968504bp}%
    \ifx\svgscale\undefined%
      \relax%
    \else%
      \setlength{\unitlength}{\unitlength * \real{\svgscale}}%
    \fi%
  \else%
    \setlength{\unitlength}{\svgwidth}%
  \fi%
  \global\let\svgwidth\undefined%
  \global\let\svgscale\undefined%
  \makeatother%
  \begin{picture}(1,1.16666667)%
    \lineheight{1}%
    \setlength\tabcolsep{0pt}%
    \put(0,0){\includegraphics[width=\unitlength,page=1]{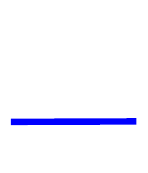}}%
    \put(0.16942518,0.13006975){\color[rgb]{0,0,0}\makebox(0,0)[lt]{\lineheight{1.25}\smash{\begin{tabular}[t]{l}$x_2$\end{tabular}}}}%
    \put(0,0){\includegraphics[width=\unitlength,page=2]{boundary_rel2_prove4_quantum.pdf}}%
    \put(0.74016697,0.12462211){\color[rgb]{0,0,0}\makebox(0,0)[lt]{\lineheight{1.25}\smash{\begin{tabular}[t]{l}$x_1$\end{tabular}}}}%
    \put(0.4330944,0.12306698){\color[rgb]{0,0,0}\makebox(0,0)[lt]{\lineheight{1.25}\smash{\begin{tabular}[t]{l}$\prec$\end{tabular}}}}%
    \put(0,0){\includegraphics[width=\unitlength,page=3]{boundary_rel2_prove4_quantum.pdf}}%
  \end{picture}%
\endgroup%
}
\overset{{\rm (S6)}}{=} (\underbrace{q - [2]_q}_{= \,\, -q^{-1} }) \hspace{-2mm}  \raisebox{-0.6\height} {} \qquad \mbox{if $s(x_1)>s(x_2)$}
\end{align}
\end{lemma}

We also state the following counterparts of the basic relations (B1)--(B4) for the cases when the strands are reversed.
\begin{lemma}
\label{lem:basic_boundary_relations_reversed}
In $\mathcal{S}^\omega_{\rm s}(\frak{S};\mathcal{R})_{\rm red}$ for a generalized marked surface $\frak{S}$, one has relations as in Fig.\ref{fig:stated_boundary_relations_reversed}.
\begin{figure}[htbp!]
\begin{center}
\hspace*{0mm}\begin{tabular}{c|c}
\raisebox{-0.6\height}{
\begingroup%
  \makeatletter%
  \providecommand\color[2][]{%
    \errmessage{(Inkscape) Color is used for the text in Inkscape, but the package 'color.sty' is not loaded}%
    \renewcommand\color[2][]{}%
  }%
  \providecommand\transparent[1]{%
    \errmessage{(Inkscape) Transparency is used (non-zero) for the text in Inkscape, but the package 'transparent.sty' is not loaded}%
    \renewcommand\transparent[1]{}%
  }%
  \providecommand\rotatebox[2]{#2}%
  \newcommand*\fsize{\dimexpr\f@size pt\relax}%
  \newcommand*\lineheight[1]{\fontsize{\fsize}{#1\fsize}\selectfont}%
  \ifx\svgwidth\undefined%
    \setlength{\unitlength}{42.51968504bp}%
    \ifx\svgscale\undefined%
      \relax%
    \else%
      \setlength{\unitlength}{\unitlength * \real{\svgscale}}%
    \fi%
  \else%
    \setlength{\unitlength}{\svgwidth}%
  \fi%
  \global\let\svgwidth\undefined%
  \global\let\svgscale\undefined%
  \makeatother%
  \begin{picture}(1,1.16666667)%
    \lineheight{1}%
    \setlength\tabcolsep{0pt}%
    \put(0,0){\includegraphics[width=\unitlength,page=1]{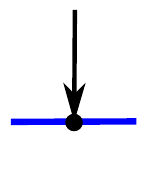}}%
    \put(0.42722329,0.10283054){\color[rgb]{0,0,0}\makebox(0,0)[lt]{\lineheight{1.25}\smash{\begin{tabular}[t]{l}$x$\end{tabular}}}}%
  \end{picture}%
\endgroup%
} \hspace{-3mm}  $=-q^{-\frac{7}{6}}$ \hspace{-5mm} \raisebox{-0.6\height} {
\begingroup%
  \makeatletter%
  \providecommand\color[2][]{%
    \errmessage{(Inkscape) Color is used for the text in Inkscape, but the package 'color.sty' is not loaded}%
    \renewcommand\color[2][]{}%
  }%
  \providecommand\transparent[1]{%
    \errmessage{(Inkscape) Transparency is used (non-zero) for the text in Inkscape, but the package 'transparent.sty' is not loaded}%
    \renewcommand\transparent[1]{}%
  }%
  \providecommand\rotatebox[2]{#2}%
  \newcommand*\fsize{\dimexpr\f@size pt\relax}%
  \newcommand*\lineheight[1]{\fontsize{\fsize}{#1\fsize}\selectfont}%
  \ifx\svgwidth\undefined%
    \setlength{\unitlength}{42.51968504bp}%
    \ifx\svgscale\undefined%
      \relax%
    \else%
      \setlength{\unitlength}{\unitlength * \real{\svgscale}}%
    \fi%
  \else%
    \setlength{\unitlength}{\svgwidth}%
  \fi%
  \global\let\svgwidth\undefined%
  \global\let\svgscale\undefined%
  \makeatother%
  \begin{picture}(1,1.16666667)%
    \lineheight{1}%
    \setlength\tabcolsep{0pt}%
    \put(0,0){\includegraphics[width=\unitlength,page=1]{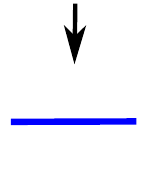}}%
    \put(0.16942518,0.13006975){\color[rgb]{0,0,0}\makebox(0,0)[lt]{\lineheight{1.25}\smash{\begin{tabular}[t]{l}$x_1$\end{tabular}}}}%
    \put(0,0){\includegraphics[width=\unitlength,page=2]{boundary_rel2_reversed_quantum.pdf}}%
    \put(0.74016697,0.12462211){\color[rgb]{0,0,0}\makebox(0,0)[lt]{\lineheight{1.25}\smash{\begin{tabular}[t]{l}$x_2$\end{tabular}}}}%
    \put(0.4330944,0.12306698){\color[rgb]{0,0,0}\makebox(0,0)[lt]{\lineheight{1.25}\smash{\begin{tabular}[t]{l}$\prec$\end{tabular}}}}%
  \end{picture}%
\endgroup%
} \hspace{0,5mm} & \raisebox{-0.6\height}{
\begingroup%
  \makeatletter%
  \providecommand\color[2][]{%
    \errmessage{(Inkscape) Color is used for the text in Inkscape, but the package 'color.sty' is not loaded}%
    \renewcommand\color[2][]{}%
  }%
  \providecommand\transparent[1]{%
    \errmessage{(Inkscape) Transparency is used (non-zero) for the text in Inkscape, but the package 'transparent.sty' is not loaded}%
    \renewcommand\transparent[1]{}%
  }%
  \providecommand\rotatebox[2]{#2}%
  \newcommand*\fsize{\dimexpr\f@size pt\relax}%
  \newcommand*\lineheight[1]{\fontsize{\fsize}{#1\fsize}\selectfont}%
  \ifx\svgwidth\undefined%
    \setlength{\unitlength}{42.51968504bp}%
    \ifx\svgscale\undefined%
      \relax%
    \else%
      \setlength{\unitlength}{\unitlength * \real{\svgscale}}%
    \fi%
  \else%
    \setlength{\unitlength}{\svgwidth}%
  \fi%
  \global\let\svgwidth\undefined%
  \global\let\svgscale\undefined%
  \makeatother%
  \begin{picture}(1,1.16666667)%
    \lineheight{1}%
    \setlength\tabcolsep{0pt}%
    \put(0,0){\includegraphics[width=\unitlength,page=1]{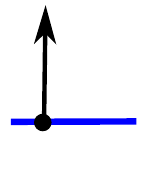}}%
    \put(0.18027884,0.10283054){\color[rgb]{0,0,0}\makebox(0,0)[lt]{\lineheight{1.25}\smash{\begin{tabular}[t]{l}$x_1$\end{tabular}}}}%
    \put(0,0){\includegraphics[width=\unitlength,page=2]{boundary_rel3_reversed_quantum.pdf}}%
    \put(0.67416761,0.10283054){\color[rgb]{0,0,0}\makebox(0,0)[lt]{\lineheight{1.25}\smash{\begin{tabular}[t]{l}$x_2$\end{tabular}}}}%
    \put(0.46250103,0.10283054){\color[rgb]{0,0,0}\makebox(0,0)[lt]{\lineheight{1.25}\smash{\begin{tabular}[t]{l}$\prec$\end{tabular}}}}%
  \end{picture}%
\endgroup%
} \hspace{-2mm} $=$ $q$ \hspace{-2mm}
 \raisebox{-0.6\height}{
\begingroup%
  \makeatletter%
  \providecommand\color[2][]{%
    \errmessage{(Inkscape) Color is used for the text in Inkscape, but the package 'color.sty' is not loaded}%
    \renewcommand\color[2][]{}%
  }%
  \providecommand\transparent[1]{%
    \errmessage{(Inkscape) Transparency is used (non-zero) for the text in Inkscape, but the package 'transparent.sty' is not loaded}%
    \renewcommand\transparent[1]{}%
  }%
  \providecommand\rotatebox[2]{#2}%
  \newcommand*\fsize{\dimexpr\f@size pt\relax}%
  \newcommand*\lineheight[1]{\fontsize{\fsize}{#1\fsize}\selectfont}%
  \ifx\svgwidth\undefined%
    \setlength{\unitlength}{42.51968504bp}%
    \ifx\svgscale\undefined%
      \relax%
    \else%
      \setlength{\unitlength}{\unitlength * \real{\svgscale}}%
    \fi%
  \else%
    \setlength{\unitlength}{\svgwidth}%
  \fi%
  \global\let\svgwidth\undefined%
  \global\let\svgscale\undefined%
  \makeatother%
  \begin{picture}(1,1.16666667)%
    \lineheight{1}%
    \setlength\tabcolsep{0pt}%
    \put(0,0){\includegraphics[width=\unitlength,page=1]{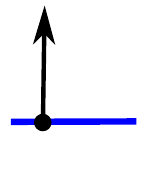}}%
    \put(0.18027884,0.10283054){\color[rgb]{0,0,0}\makebox(0,0)[lt]{\lineheight{1.25}\smash{\begin{tabular}[t]{l}$x_2$\end{tabular}}}}%
    \put(0,0){\includegraphics[width=\unitlength,page=2]{boundary_rel4_reversed_quantum.pdf}}%
    \put(0.67416761,0.10283054){\color[rgb]{0,0,0}\makebox(0,0)[lt]{\lineheight{1.25}\smash{\begin{tabular}[t]{l}$x_1$\end{tabular}}}}%
    \put(0.46250103,0.10283054){\color[rgb]{0,0,0}\makebox(0,0)[lt]{\lineheight{1.25}\smash{\begin{tabular}[t]{l}$\prec$\end{tabular}}}}%
  \end{picture}%
\endgroup%
}
 \hspace{-3mm} $+$ \hspace{-2mm}
\raisebox{-0.6\height}{
\begingroup%
  \makeatletter%
  \providecommand\color[2][]{%
    \errmessage{(Inkscape) Color is used for the text in Inkscape, but the package 'color.sty' is not loaded}%
    \renewcommand\color[2][]{}%
  }%
  \providecommand\transparent[1]{%
    \errmessage{(Inkscape) Transparency is used (non-zero) for the text in Inkscape, but the package 'transparent.sty' is not loaded}%
    \renewcommand\transparent[1]{}%
  }%
  \providecommand\rotatebox[2]{#2}%
  \newcommand*\fsize{\dimexpr\f@size pt\relax}%
  \newcommand*\lineheight[1]{\fontsize{\fsize}{#1\fsize}\selectfont}%
  \ifx\svgwidth\undefined%
    \setlength{\unitlength}{42.51968504bp}%
    \ifx\svgscale\undefined%
      \relax%
    \else%
      \setlength{\unitlength}{\unitlength * \real{\svgscale}}%
    \fi%
  \else%
    \setlength{\unitlength}{\svgwidth}%
  \fi%
  \global\let\svgwidth\undefined%
  \global\let\svgscale\undefined%
  \makeatother%
  \begin{picture}(1,1.16666667)%
    \lineheight{1}%
    \setlength\tabcolsep{0pt}%
    \put(0,0){\includegraphics[width=\unitlength,page=1]{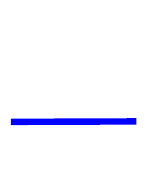}}%
    \put(0.13088999,0.10283052){\color[rgb]{0,0,0}\makebox(0,0)[lt]{\lineheight{1.25}\smash{\begin{tabular}[t]{l}$x_2$\end{tabular}}}}%
    \put(0,0){\includegraphics[width=\unitlength,page=2]{boundary_rel5_reversed_quantum.pdf}}%
    \put(0.72355657,0.10283052){\color[rgb]{0,0,0}\makebox(0,0)[lt]{\lineheight{1.25}\smash{\begin{tabular}[t]{l}$x_1$\end{tabular}}}}%
    \put(0,0){\includegraphics[width=\unitlength,page=3]{boundary_rel5_reversed_quantum.pdf}}%
    \put(0.44133449,0.10283052){\color[rgb]{0,0,0}\makebox(0,0)[lt]{\lineheight{1.25}\smash{\begin{tabular}[t]{l}$\prec$\end{tabular}}}}%
  \end{picture}%
\endgroup%
}  \\
{\rm (B1')} 
$s(x)=\varepsilon$, $s(x_1)=r_1(\varepsilon)$, $s(x_2)=r_2(\varepsilon)$ 
& {\rm (B2')} $s(x_1)=\varepsilon_1$, $s(x_2)= \varepsilon_2$, with  $\varepsilon_1>\varepsilon_2$  \\ \hline
\raisebox{-0.6\height}{
\begingroup%
  \makeatletter%
  \providecommand\color[2][]{%
    \errmessage{(Inkscape) Color is used for the text in Inkscape, but the package 'color.sty' is not loaded}%
    \renewcommand\color[2][]{}%
  }%
  \providecommand\transparent[1]{%
    \errmessage{(Inkscape) Transparency is used (non-zero) for the text in Inkscape, but the package 'transparent.sty' is not loaded}%
    \renewcommand\transparent[1]{}%
  }%
  \providecommand\rotatebox[2]{#2}%
  \newcommand*\fsize{\dimexpr\f@size pt\relax}%
  \newcommand*\lineheight[1]{\fontsize{\fsize}{#1\fsize}\selectfont}%
  \ifx\svgwidth\undefined%
    \setlength{\unitlength}{42.51968504bp}%
    \ifx\svgscale\undefined%
      \relax%
    \else%
      \setlength{\unitlength}{\unitlength * \real{\svgscale}}%
    \fi%
  \else%
    \setlength{\unitlength}{\svgwidth}%
  \fi%
  \global\let\svgwidth\undefined%
  \global\let\svgscale\undefined%
  \makeatother%
  \begin{picture}(1,1.16666667)%
    \lineheight{1}%
    \setlength\tabcolsep{0pt}%
    \put(0,0){\includegraphics[width=\unitlength,page=1]{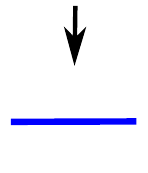}}%
    \put(0.16616776,0.10283052){\color[rgb]{0,0,0}\makebox(0,0)[lt]{\lineheight{1.25}\smash{\begin{tabular}[t]{l}$x$\end{tabular}}}}%
    \put(0,0){\includegraphics[width=\unitlength,page=2]{boundary_rel6_reversed_quantum.pdf}}%
    \put(0.72355657,0.10283052){\color[rgb]{0,0,0}\makebox(0,0)[lt]{\lineheight{1.25}\smash{\begin{tabular}[t]{l}$y$\end{tabular}}}}%
    \put(0.44133449,0.10283052){\color[rgb]{0,0,0}\makebox(0,0)[lt]{\lineheight{1.25}\smash{\begin{tabular}[t]{l}$\prec$\end{tabular}}}}%
  \end{picture}%
\endgroup%
} \hspace{-2mm} $=0$
& \raisebox{-0.6\height}{
\begingroup%
  \makeatletter%
  \providecommand\color[2][]{%
    \errmessage{(Inkscape) Color is used for the text in Inkscape, but the package 'color.sty' is not loaded}%
    \renewcommand\color[2][]{}%
  }%
  \providecommand\transparent[1]{%
    \errmessage{(Inkscape) Transparency is used (non-zero) for the text in Inkscape, but the package 'transparent.sty' is not loaded}%
    \renewcommand\transparent[1]{}%
  }%
  \providecommand\rotatebox[2]{#2}%
  \newcommand*\fsize{\dimexpr\f@size pt\relax}%
  \newcommand*\lineheight[1]{\fontsize{\fsize}{#1\fsize}\selectfont}%
  \ifx\svgwidth\undefined%
    \setlength{\unitlength}{42.51968504bp}%
    \ifx\svgscale\undefined%
      \relax%
    \else%
      \setlength{\unitlength}{\unitlength * \real{\svgscale}}%
    \fi%
  \else%
    \setlength{\unitlength}{\svgwidth}%
  \fi%
  \global\let\svgwidth\undefined%
  \global\let\svgscale\undefined%
  \makeatother%
  \begin{picture}(1,1.16666667)%
    \lineheight{1}%
    \setlength\tabcolsep{0pt}%
    \put(0,0){\includegraphics[width=\unitlength,page=1]{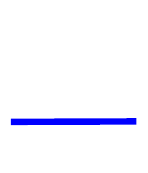}}%
    \put(0.04567213,0.10283052){\color[rgb]{0,0,0}\makebox(0,0)[lt]{\lineheight{1.25}\smash{\begin{tabular}[t]{l}$x_1$\end{tabular}}}}%
    \put(0,0){\includegraphics[width=\unitlength,page=2]{boundary_rel7_reversed_quantum.pdf}}%
    \put(0.44839,0.10283052){\color[rgb]{0,0,0}\makebox(0,0)[lt]{\lineheight{1.25}\smash{\begin{tabular}[t]{l}$x_2$\end{tabular}}}}%
    \put(0.83644529,0.10283052){\color[rgb]{0,0,0}\makebox(0,0)[lt]{\lineheight{1.25}\smash{\begin{tabular}[t]{l}$x_3$\end{tabular}}}}%
    \put(0.23672332,0.10283052){\color[rgb]{0,0,0}\makebox(0,0)[lt]{\lineheight{1.25}\smash{\begin{tabular}[t]{l}$\prec$\end{tabular}}}}%
    \put(0.66005655,0.10283052){\color[rgb]{0,0,0}\makebox(0,0)[lt]{\lineheight{1.25}\smash{\begin{tabular}[t]{l}$\prec$\end{tabular}}}}%
  \end{picture}%
\endgroup%
} \hspace{-2mm} $=$ $-q^{\frac{7}{2}}$ \hspace{-3,5mm} \raisebox{-0.6\height}{\input{boundary_rel8.pdf_tex}} \\ 
{\rm (B3')} $s(x)=s(y)$ & {\rm (B4')} $s(x_1)=1$, $s(x_2)=2$, $s(x_3)=3$
\end{tabular}
\end{center}
\vspace{-2mm}
\caption{Boundary relations for stated ${\rm SL}_3$-skeins, for reverse oriented strands}
\vspace{-2mm}
\label{fig:stated_boundary_relations_reversed}
\end{figure}

\end{lemma}
Lem.\ref{lem:basic_boundary_relations_reversed} was communicated to the author by a personal correspondence with Vijay Higgins (and translated into the above version via the isomorphism in eq.\eqref{eq:isomorphism_from_ours_to_Higgins}), and should be a straightforward exercise using the defining relations. In fact, it is our isomorphism in eq.\eqref{eq:isomorphism_from_ours_to_Higgins} that makes the coefficients in the relations in Fig.\ref{fig:stated_boundary_relations} completely symmetric with those in Fig.\ref{fig:stated_boundary_relations_reversed}. 

\vs

It is convenient to establish the orientation reversed versions of Lem.\ref{lem:3-way_boundary_relations1}--\ref{lem:3-way_boundary_relations2} too:
\begin{lemma}
\label{lem:3-way_boundary_relations_reversed1}
In $\mathcal{S}^\omega_{\rm s}(\frak{S};\mathcal{R})$ (hence also in $\mathcal{S}^\omega_{\rm s}(\frak{S};\mathcal{R})_{\rm red}$) for a generalized marked surface $\frak{S}$, one has:
\begin{align}
\label{eq:incoming_fork_relations1}
\raisebox{-0.6\height} {}
& = \raisebox{-0.6\height} {
\begingroup%
  \makeatletter%
  \providecommand\color[2][]{%
    \errmessage{(Inkscape) Color is used for the text in Inkscape, but the package 'color.sty' is not loaded}%
    \renewcommand\color[2][]{}%
  }%
  \providecommand\transparent[1]{%
    \errmessage{(Inkscape) Transparency is used (non-zero) for the text in Inkscape, but the package 'transparent.sty' is not loaded}%
    \renewcommand\transparent[1]{}%
  }%
  \providecommand\rotatebox[2]{#2}%
  \newcommand*\fsize{\dimexpr\f@size pt\relax}%
  \newcommand*\lineheight[1]{\fontsize{\fsize}{#1\fsize}\selectfont}%
  \ifx\svgwidth\undefined%
    \setlength{\unitlength}{42.51968504bp}%
    \ifx\svgscale\undefined%
      \relax%
    \else%
      \setlength{\unitlength}{\unitlength * \real{\svgscale}}%
    \fi%
  \else%
    \setlength{\unitlength}{\svgwidth}%
  \fi%
  \global\let\svgwidth\undefined%
  \global\let\svgscale\undefined%
  \makeatother%
  \begin{picture}(1,1.16666667)%
    \lineheight{1}%
    \setlength\tabcolsep{0pt}%
    \put(0,0){\includegraphics[width=\unitlength,page=1]{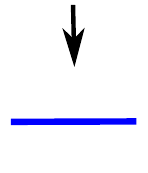}}%
    \put(0.16942518,0.13006975){\color[rgb]{0,0,0}\makebox(0,0)[lt]{\lineheight{1.25}\smash{\begin{tabular}[t]{l}$x_2$\end{tabular}}}}%
    \put(0,0){\includegraphics[width=\unitlength,page=2]{boundary_rel2_reversed_prove2_quantum.pdf}}%
    \put(0.74016697,0.12462211){\color[rgb]{0,0,0}\makebox(0,0)[lt]{\lineheight{1.25}\smash{\begin{tabular}[t]{l}$x_1$\end{tabular}}}}%
    \put(0.4330944,0.12306698){\color[rgb]{0,0,0}\makebox(0,0)[lt]{\lineheight{1.25}\smash{\begin{tabular}[t]{l}$\succ$\end{tabular}}}}%
    \put(0,0){\includegraphics[width=\unitlength,page=3]{boundary_rel2_reversed_prove2_quantum.pdf}}%
  \end{picture}%
\endgroup%
} 
\overset{{\rm (S8)}}{=} q^{-2/3} \hspace{-3mm} \raisebox{-0.6\height} {
\begingroup%
  \makeatletter%
  \providecommand\color[2][]{%
    \errmessage{(Inkscape) Color is used for the text in Inkscape, but the package 'color.sty' is not loaded}%
    \renewcommand\color[2][]{}%
  }%
  \providecommand\transparent[1]{%
    \errmessage{(Inkscape) Transparency is used (non-zero) for the text in Inkscape, but the package 'transparent.sty' is not loaded}%
    \renewcommand\transparent[1]{}%
  }%
  \providecommand\rotatebox[2]{#2}%
  \newcommand*\fsize{\dimexpr\f@size pt\relax}%
  \newcommand*\lineheight[1]{\fontsize{\fsize}{#1\fsize}\selectfont}%
  \ifx\svgwidth\undefined%
    \setlength{\unitlength}{42.51968504bp}%
    \ifx\svgscale\undefined%
      \relax%
    \else%
      \setlength{\unitlength}{\unitlength * \real{\svgscale}}%
    \fi%
  \else%
    \setlength{\unitlength}{\svgwidth}%
  \fi%
  \global\let\svgwidth\undefined%
  \global\let\svgscale\undefined%
  \makeatother%
  \begin{picture}(1,1.16666667)%
    \lineheight{1}%
    \setlength\tabcolsep{0pt}%
    \put(0,0){\includegraphics[width=\unitlength,page=1]{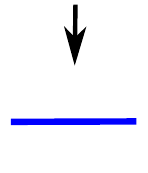}}%
    \put(0.16942518,0.13006975){\color[rgb]{0,0,0}\makebox(0,0)[lt]{\lineheight{1.25}\smash{\begin{tabular}[t]{l}$x_2$\end{tabular}}}}%
    \put(0,0){\includegraphics[width=\unitlength,page=2]{boundary_rel2_reversed_prime_quantum.pdf}}%
    \put(0.74016697,0.12462211){\color[rgb]{0,0,0}\makebox(0,0)[lt]{\lineheight{1.25}\smash{\begin{tabular}[t]{l}$x_1$\end{tabular}}}}%
    \put(0.50364996,0.12306698){\color[rgb]{0,0,0}\makebox(0,0)[lt]{\lineheight{1.25}\smash{\begin{tabular}[t]{l}$\succ$\end{tabular}}}}%
  \end{picture}%
\endgroup%
}
+ q^{1/3} \hspace{-3mm} \raisebox{-0.6\height} {
\begingroup%
  \makeatletter%
  \providecommand\color[2][]{%
    \errmessage{(Inkscape) Color is used for the text in Inkscape, but the package 'color.sty' is not loaded}%
    \renewcommand\color[2][]{}%
  }%
  \providecommand\transparent[1]{%
    \errmessage{(Inkscape) Transparency is used (non-zero) for the text in Inkscape, but the package 'transparent.sty' is not loaded}%
    \renewcommand\transparent[1]{}%
  }%
  \providecommand\rotatebox[2]{#2}%
  \newcommand*\fsize{\dimexpr\f@size pt\relax}%
  \newcommand*\lineheight[1]{\fontsize{\fsize}{#1\fsize}\selectfont}%
  \ifx\svgwidth\undefined%
    \setlength{\unitlength}{42.51968504bp}%
    \ifx\svgscale\undefined%
      \relax%
    \else%
      \setlength{\unitlength}{\unitlength * \real{\svgscale}}%
    \fi%
  \else%
    \setlength{\unitlength}{\svgwidth}%
  \fi%
  \global\let\svgwidth\undefined%
  \global\let\svgscale\undefined%
  \makeatother%
  \begin{picture}(1,1.16666667)%
    \lineheight{1}%
    \setlength\tabcolsep{0pt}%
    \put(0,0){\includegraphics[width=\unitlength,page=1]{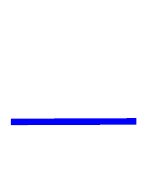}}%
    \put(0.16942518,0.13006975){\color[rgb]{0,0,0}\makebox(0,0)[lt]{\lineheight{1.25}\smash{\begin{tabular}[t]{l}$x_2$\end{tabular}}}}%
    \put(0,0){\includegraphics[width=\unitlength,page=2]{boundary_rel2_reversed_prove3_quantum.pdf}}%
    \put(0.74016697,0.12462211){\color[rgb]{0,0,0}\makebox(0,0)[lt]{\lineheight{1.25}\smash{\begin{tabular}[t]{l}$x_1$\end{tabular}}}}%
    \put(0.4330944,0.12306698){\color[rgb]{0,0,0}\makebox(0,0)[lt]{\lineheight{1.25}\smash{\begin{tabular}[t]{l}$\succ$\end{tabular}}}}%
    \put(0,0){\includegraphics[width=\unitlength,page=3]{boundary_rel2_reversed_prove3_quantum.pdf}}%
  \end{picture}%
\endgroup%
} 
\overset{{\rm (S6)}}{=} (\underbrace{ q^{-2/3} - q^{1/3} [2]_q }_{= \,\, - q^{4/3}}) \hspace{-2mm} \raisebox{-0.6\height} {}
\end{align}
\end{lemma}

\begin{lemma}
\label{lem:3-way_boundary_relations_reversed2}
In $\mathcal{S}^\omega_{\rm s}(\frak{S};\mathcal{R})_{\rm red}$ for a generalized marked surface $\frak{S}$, one has:
\begin{align}
\label{eq:incoming_fork_relations2}
\raisebox{-0.6\height} {}
& \overset{\mbox{\rm \tiny (B2')}}{=} q\hspace{-1mm} \raisebox{-0.6\height} {
\begingroup%
  \makeatletter%
  \providecommand\color[2][]{%
    \errmessage{(Inkscape) Color is used for the text in Inkscape, but the package 'color.sty' is not loaded}%
    \renewcommand\color[2][]{}%
  }%
  \providecommand\transparent[1]{%
    \errmessage{(Inkscape) Transparency is used (non-zero) for the text in Inkscape, but the package 'transparent.sty' is not loaded}%
    \renewcommand\transparent[1]{}%
  }%
  \providecommand\rotatebox[2]{#2}%
  \newcommand*\fsize{\dimexpr\f@size pt\relax}%
  \newcommand*\lineheight[1]{\fontsize{\fsize}{#1\fsize}\selectfont}%
  \ifx\svgwidth\undefined%
    \setlength{\unitlength}{42.51968504bp}%
    \ifx\svgscale\undefined%
      \relax%
    \else%
      \setlength{\unitlength}{\unitlength * \real{\svgscale}}%
    \fi%
  \else%
    \setlength{\unitlength}{\svgwidth}%
  \fi%
  \global\let\svgwidth\undefined%
  \global\let\svgscale\undefined%
  \makeatother%
  \begin{picture}(1,1.16666667)%
    \lineheight{1}%
    \setlength\tabcolsep{0pt}%
    \put(0,0){\includegraphics[width=\unitlength,page=1]{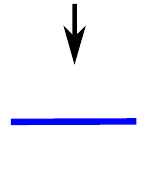}}%
    \put(0.16942518,0.13006975){\color[rgb]{0,0,0}\makebox(0,0)[lt]{\lineheight{1.25}\smash{\begin{tabular}[t]{l}$x_2$\end{tabular}}}}%
    \put(0,0){\includegraphics[width=\unitlength,page=2]{boundary_rel2_reversed_prime2_quantum.pdf}}%
    \put(0.74016697,0.12462211){\color[rgb]{0,0,0}\makebox(0,0)[lt]{\lineheight{1.25}\smash{\begin{tabular}[t]{l}$x_1$\end{tabular}}}}%
    \put(0.4330944,0.12306698){\color[rgb]{0,0,0}\makebox(0,0)[lt]{\lineheight{1.25}\smash{\begin{tabular}[t]{l}$\prec$\end{tabular}}}}%
  \end{picture}%
\endgroup%
} +  \hspace{-1mm} \raisebox{-0.6\height} {}
\overset{{\rm (S6)}}{=} (\underbrace{q - [2]_q}_{= \,\, -q^{-1} }) \hspace{-2mm}  \raisebox{-0.6\height} {} \qquad \mbox{if $s(x_1)>s(x_2)$}
\end{align}
\end{lemma}

Coming back to our proof of Prop.\ref{prop:biangle_SL3_quantum_trace}, let $(W,s)$ be as in (BT2-3), whether $x_1,x_2,x_3$ are all sinks or all sources. Suppose $x_1\succ x_2$. By relation (B3) of Fig.\ref{fig:stated_boundary_relations} and (B3') of Fig.\ref{fig:stated_boundary_relations_reversed}, we see that $[W,s] \in \mathcal{S}^\omega_{\rm s}(\vec{B};\mathbb{Z})_{\rm red}$ equals zero if $s(x_1)=s(x_2)$. In case $(s(x_1),s(x_2)) = (r_2(\varepsilon),r_1(\varepsilon))$ for $\varepsilon\in\{1,2,3\}$, by (B1) of Fig.\ref{fig:stated_boundary_relations} or (B1') of Fig.\ref{fig:stated_boundary_relations_reversed} we have $[W,s] = -q^{7/6} [W',s']$ in $\mathcal{S}^\omega_{\rm s}(\vec{B};\mathbb{Z})_{\rm red}$, where $W'$ is a single crossing-less arc connecting a point $x$ in one side and a point $x_3$ in the other side, with $s'(x)=\varepsilon$ and $s'(x_3)=s(x_3)$. Meanwhile, in case $(s(x_1),s(x_2)) = (r_1(\varepsilon),r_2(\varepsilon))$ for $\varepsilon\in\{1,2,3\}$, we see using eq.\eqref{eq:outgoing_fork_relations2} of Lem.\ref{lem:3-way_boundary_relations2} or eq.\eqref{eq:incoming_fork_relations2} of Lem.\ref{lem:3-way_boundary_relations_reversed2} that $[W,s] = -q^{-1}[W,s'']$, with $s''(x_1)=s(x_2)$, $s''(x_2)=s(x_1)$, and by (B1) of Fig.\ref{fig:stated_boundary_relations} or (B1') of Fig.\ref{fig:stated_boundary_relations_reversed} we have $[W,s''] = -q^{7/6} [W',s']$ as before. Note $[W',s']$ falls into (BT2-1), hence $\epsilon[W',s']$ equals $1$ if $s'(x)=s'(x_3)$ and equals zero otherwise. Hence $\epsilon[W,s]$ equals  $-q^{7/6}$ if $(s(x_1),s(x_2),s(x_3))=(2,1,1),(3,1,2),(3,2,3)$, equals $q^{1/6}$ if $(s(x_1),s(x_2),s(x_3))=(1,2,1),(1,3,2),(2,3,3)$, and equals zero otherwise. Now, suppose $x_1,x_2,x_3$ are sinks, but this time $x_1 \prec x_2$. By eq.\eqref{eq:outgoing_fork_relations1} of Lem.\ref{lem:3-way_boundary_relations1}, the result is $-q^{-4/3}$ times the above case with the roles of $x_1$ and $x_2$ exchanged, thus $\epsilon[W,s]$ equals $q^{-1/6}$ if $(s(x_1),s(x_2),s(x_3))=(1,2,1),(1,3,2),(2,3,3)$, equals $q^{-7/6}$ if $(s(x_1),s(x_2),s(x_3))=(2,1,1),(3,1,2),(3,2,3)$, and equals zero otherwise. So $\epsilon$ satisfies (BT2-3).

\vs

To finish the proof of the existence part of Prop.\ref{prop:biangle_SL3_quantum_trace}, observe that the values of $\epsilon$ at the stated ${\rm SL}_3$-webs appearing in (BT2-1), (BT2-2) and (BT2-3) do not depend on the choice of the direction on $B$. Hence $\epsilon : \mathcal{S}^\omega_{\rm s}(\vec{B};\mathbb{Z})_{\rm red} \to \mathbb{Z}[\omega^{\pm 1/2}]$ indeed provides a well-defined map ${\rm Tr}^\omega_B: \mathcal{S}^\omega_{\rm s}(B;\mathbb{Z})_{\rm red} \to \mathbb{Z}[\omega^{\pm 1/2}]$ with desired properties. \qed \quad {\it End of proof of Prop.\ref{prop:biangle_SL3_quantum_trace}.}

\vs

For later use, we show Prop.\ref{prop:biangle_SL3_quantum_trace_some_values}. For (BT2-4) and for later use also, it is handy to have the following:
\begin{lemma}[values of `I-webs' and `H-webs' under the biangle ${\rm SL}_3$ quantum trace]
\label{lem:value_of_I-webs_in_biangle}
Let the ${\rm SL}_3$-webs $W, W'$ in a thickened biangle $B \times {\bf I}$ be 
$$
W = \raisebox{-0.5\height}{
\begingroup%
  \makeatletter%
  \providecommand\color[2][]{%
    \errmessage{(Inkscape) Color is used for the text in Inkscape, but the package 'color.sty' is not loaded}%
    \renewcommand\color[2][]{}%
  }%
  \providecommand\transparent[1]{%
    \errmessage{(Inkscape) Transparency is used (non-zero) for the text in Inkscape, but the package 'transparent.sty' is not loaded}%
    \renewcommand\transparent[1]{}%
  }%
  \providecommand\rotatebox[2]{#2}%
  \newcommand*\fsize{\dimexpr\f@size pt\relax}%
  \newcommand*\lineheight[1]{\fontsize{\fsize}{#1\fsize}\selectfont}%
  \ifx\svgwidth\undefined%
    \setlength{\unitlength}{56.69291339bp}%
    \ifx\svgscale\undefined%
      \relax%
    \else%
      \setlength{\unitlength}{\unitlength * \real{\svgscale}}%
    \fi%
  \else%
    \setlength{\unitlength}{\svgwidth}%
  \fi%
  \global\let\svgwidth\undefined%
  \global\let\svgscale\undefined%
  \makeatother%
  \begin{picture}(1,0.875)%
    \lineheight{1}%
    \setlength\tabcolsep{0pt}%
    \put(0,0){\includegraphics[width=\unitlength,page=1]{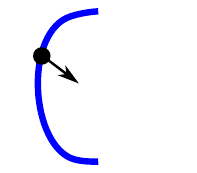}}%
    \put(-0.00171479,0.57550812){\color[rgb]{0,0,0}\makebox(0,0)[lt]{\lineheight{1.25}\smash{\begin{tabular}[t]{l}$x_1$\end{tabular}}}}%
    \put(0,0){\includegraphics[width=\unitlength,page=2]{biangle_I-web0.pdf}}%
    \put(0.00287621,0.25922147){\color[rgb]{0,0,0}\makebox(0,0)[lt]{\lineheight{1.25}\smash{\begin{tabular}[t]{l}$x_2$\end{tabular}}}}%
    \put(0,0){\includegraphics[width=\unitlength,page=3]{biangle_I-web0.pdf}}%
    \put(0.85239343,0.58058091){\color[rgb]{0,0,0}\makebox(0,0)[lt]{\lineheight{1.25}\smash{\begin{tabular}[t]{l}$y_1$\end{tabular}}}}%
    \put(0.85698442,0.26429425){\color[rgb]{0,0,0}\makebox(0,0)[lt]{\lineheight{1.25}\smash{\begin{tabular}[t]{l}$y_2$\end{tabular}}}}%
    \put(0,0){\includegraphics[width=\unitlength,page=4]{biangle_I-web0.pdf}}%
  \end{picture}%
\endgroup%
}, \qquad
W' = \raisebox{-0.5\height}{
\begingroup%
  \makeatletter%
  \providecommand\color[2][]{%
    \errmessage{(Inkscape) Color is used for the text in Inkscape, but the package 'color.sty' is not loaded}%
    \renewcommand\color[2][]{}%
  }%
  \providecommand\transparent[1]{%
    \errmessage{(Inkscape) Transparency is used (non-zero) for the text in Inkscape, but the package 'transparent.sty' is not loaded}%
    \renewcommand\transparent[1]{}%
  }%
  \providecommand\rotatebox[2]{#2}%
  \newcommand*\fsize{\dimexpr\f@size pt\relax}%
  \newcommand*\lineheight[1]{\fontsize{\fsize}{#1\fsize}\selectfont}%
  \ifx\svgwidth\undefined%
    \setlength{\unitlength}{56.69291339bp}%
    \ifx\svgscale\undefined%
      \relax%
    \else%
      \setlength{\unitlength}{\unitlength * \real{\svgscale}}%
    \fi%
  \else%
    \setlength{\unitlength}{\svgwidth}%
  \fi%
  \global\let\svgwidth\undefined%
  \global\let\svgscale\undefined%
  \makeatother%
  \begin{picture}(1,0.875)%
    \lineheight{1}%
    \setlength\tabcolsep{0pt}%
    \put(0,0){\includegraphics[width=\unitlength,page=1]{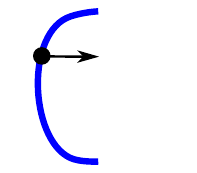}}%
    \put(-0.00171479,0.57550812){\color[rgb]{0,0,0}\makebox(0,0)[lt]{\lineheight{1.25}\smash{\begin{tabular}[t]{l}$x_1$\end{tabular}}}}%
    \put(0,0){\includegraphics[width=\unitlength,page=2]{biangle_H-web0.pdf}}%
    \put(0.00287621,0.25922147){\color[rgb]{0,0,0}\makebox(0,0)[lt]{\lineheight{1.25}\smash{\begin{tabular}[t]{l}$x_2$\end{tabular}}}}%
    \put(0,0){\includegraphics[width=\unitlength,page=3]{biangle_H-web0.pdf}}%
    \put(0.85239343,0.58058091){\color[rgb]{0,0,0}\makebox(0,0)[lt]{\lineheight{1.25}\smash{\begin{tabular}[t]{l}$y_1$\end{tabular}}}}%
    \put(0.85698442,0.26429425){\color[rgb]{0,0,0}\makebox(0,0)[lt]{\lineheight{1.25}\smash{\begin{tabular}[t]{l}$y_2$\end{tabular}}}}%
    \put(0,0){\includegraphics[width=\unitlength,page=4]{biangle_H-web0.pdf}}%
  \end{picture}%
\endgroup%
},
$$
with the elevation ordering on the endpoints is given either by $x_1\succ x_2$, $y_1\succ y_2$, or by $x_1\prec x_2$, $y_1\prec y_2$\redfix{.} 
Choose an arbitrary direction on $B$ to make it a directed biangle $\vec{B}$. 

\vs

$\bullet$ Let $s$ be a state of $W$. Then the value ${\rm Tr}^\omega_B([W,s])$ for ${\rm Tr}^\omega_B$ of Prop.\ref{prop:biangle_SL3_quantum_trace} is given as follows. Writing down the values of ${\rm Tr}^\omega_B([W,s])$ as entries of a $9\times 9$ matrix, \`a la the description right before Def.\ref{def:composition_of_biangles} by choosing $\xi_{1,1}$, $\xi_{1,2}$, $\xi_{1,3}$, $\xi_{2,1}$, $\xi_{2,2}$, $\xi_{2,3}$, $\xi_{3,1}$, $\xi_{3,2}$, $\xi_{3,3}$ as ordered bases for both the domain $V^{\otimes 2}$ (whose basis vectors are written as $\xi_{s(x_1),s(x_2)}$) and the codomain $V^{\otimes 2}$ (whose basis vectors are written as $\xi_{s(y_1),s(y_2)}$), one has
\begin{align}
\label{eq:I-web0}
{\rm Tr}^\omega_B([W,s])
= \left\{
\begin{array}{ll}
(\wh{\bf I}_+)_{(s(x_1),s(x_2)),(s(y_1),s(y_2))} & \mbox{if $x_1\succ x_2$, $y_1\succ y_2$}, \\
(\wh{\bf I}_-)_{(s(x_1),s(x_2)),(s(y_1),s(y_2))} & \mbox{if $x_1\prec x_2$, $y_1\prec y_2$},
\end{array}
\right.
\end{align}
where the matrices $\wh{\bf I}_+$ and $\wh{\bf I}_-$ are:
\begin{align*}
\wh{\bf I}_+
= \left( \begin{smallmatrix}
0 & 0 & 0 & 0 & 0 & 0 & 0 & 0 & 0  \\
0 & -q^{-1} & 0 & 1 & 0 & 0 & 0 & 0 & 0 \\
0 & 0 & -q^{-1} & 0 & 0 & 0 & 1 & 0 & 0 \\
0 & 1 & 0 & -q & 0 & 0 & 0 & 0 & 0 \\
0 & 0 & 0 & 0 & 0 & 0 & 0 & 0 & 0  \\
0 & 0 & 0 & 0 & 0 & -q^{-1} & 0 & 1 & 0 \\
0 & 0 & 1 & 0 & 0 & 0 & -q & 0 & 0 \\
0 & 0 & 0 & 0 & 0 & 1 & 0 & -q & 0 \\
0 & 0 & 0 & 0 & 0 & 0 & 0 & 0 & 0  \\
\end{smallmatrix}
\right), \quad
\wh{\bf I}_-
= \left( \begin{smallmatrix}
0 & 0 & 0 & 0 & 0 & 0 & 0 & 0 & 0  \\
0 & - q & 0 & 1 & 0 & 0 & 0 & 0 & 0 \\
0 & 0 & -q & 0 & 0 & 0 & 1 & 0 & 0 \\
0 & 1 & 0 & -q^{-1} & 0 & 0 & 0 & 0 & 0 \\
0 & 0 & 0 & 0 & 0 & 0 & 0 & 0 & 0  \\
0 & 0 & 0 & 0 & 0 & -q & 0 & 1 & 0 \\
0 & 0 & 1 & 0 & 0 & 0 & -q^{-1} & 0 & 0 \\
0 & 0 & 0 & 0 & 0 & 1 & 0 & -q^{-1} & 0 \\
0 & 0 & 0 & 0 & 0 & 0 & 0 & 0 & 0  \\
\end{smallmatrix}
\right)
\end{align*}

\vs

$\bullet$ Let $s'$ be a state of $W'$. Using the notations above, we have
\begin{align}
\label{eq:H-web0}
{\rm Tr}^\omega_B([W,s])
= \left\{
\begin{array}{ll}
(\wh{\bf H}_+)_{(s(x_1),s(x_2)),(s(y_1),s(y_2))} & \mbox{if $x_1\succ x_2$, $y_1\succ y_2$}, \\
(\wh{\bf H}_-)_{(s(x_1),s(x_2)),(s(y_1),s(y_2))} & \mbox{if $x_1\prec x_2$, $y_1\prec y_2$},
\end{array}
\right.
\end{align}
where the $9\times 9$ matrices $\wh{\bf H}_+$ and $\wh{\bf H}_-$ are:

\begin{align*}
\wh{\bf H}_+
= \left( \begin{smallmatrix}
1 & 0 & 0 & 0 & 0  & 0 & 0 & 0 & 0 \\
0 & 0 & 0 & 1 & 0 & 0 & 0 & 0 & 0 \\
0 & 0 & -q^{-1} & 0 & 1 & 0 & 0 & 0 & 0 \\
0 & 1 & 0 & 0 & 0 & 0 & 0 & 0 & 0 \\
0 & 0 & 1 & 0 & 0 & 0 & 1 & 0 & 0 \\
0 & 0 & 0 & 0 & 0 & 0 & 0 & 1 & 0 \\
0 & 0 & 0 & 0 & 1 & 0 & -q  & 0 & 0 \\
0 & 0 & 0 & 0 & 0 & 1 & 0 & 0 & 0 \\
0 & 0 & 0 & 0 & 0 & 0 & 0 & 0 & 1 
\end{smallmatrix}
\right), \quad
\wh{\bf H}_-
= \left( \begin{smallmatrix}
1 & 0 & 0 & 0 & 0  & 0 & 0 & 0 & 0 \\
0 & 0 & 0 & 1 & 0 & 0 & 0 & 0 & 0 \\
0 & 0 & -q & 0 & 1 & 0 & 0 & 0 & 0 \\
0 & 1 & 0 & 0 & 0 & 0 & 0 & 0 & 0 \\
0 & 0 & 1 & 0 & 0 & 0 & 1 & 0 & 0 \\
0 & 0 & 0 & 0 & 0 & 0 & 0 & 1 & 0 \\
0 & 0 & 0 & 0 & 1 & 0 & -q^{-1}  & 0 & 0 \\
0 & 0 & 0 & 0 & 0 & 1 & 0 & 0 & 0 \\
0 & 0 & 0 & 0 & 0 & 0 & 0 & 0 & 1 
\end{smallmatrix}
\right)
\end{align*}

\end{lemma}

{\it Proof.} To prove eq.\eqref{eq:I-web0}, cut the biangle into two biangles $B_1,B_2$ along an arc $e$ connecting the marked points that meets the ${\rm SL}_3$-web $W$ in question exactly once transversally, say at $z$, to express $W$ as composition of $W_1 = W \cap (B_1 \times {\bf I})$ and $W_2 = W \cap (B_2 \times {\bf I})$. 
\begin{align*}
W
=  \raisebox{-0.5\height}{
\begingroup%
  \makeatletter%
  \providecommand\color[2][]{%
    \errmessage{(Inkscape) Color is used for the text in Inkscape, but the package 'color.sty' is not loaded}%
    \renewcommand\color[2][]{}%
  }%
  \providecommand\transparent[1]{%
    \errmessage{(Inkscape) Transparency is used (non-zero) for the text in Inkscape, but the package 'transparent.sty' is not loaded}%
    \renewcommand\transparent[1]{}%
  }%
  \providecommand\rotatebox[2]{#2}%
  \newcommand*\fsize{\dimexpr\f@size pt\relax}%
  \newcommand*\lineheight[1]{\fontsize{\fsize}{#1\fsize}\selectfont}%
  \ifx\svgwidth\undefined%
    \setlength{\unitlength}{56.69291339bp}%
    \ifx\svgscale\undefined%
      \relax%
    \else%
      \setlength{\unitlength}{\unitlength * \real{\svgscale}}%
    \fi%
  \else%
    \setlength{\unitlength}{\svgwidth}%
  \fi%
  \global\let\svgwidth\undefined%
  \global\let\svgscale\undefined%
  \makeatother%
  \begin{picture}(1,0.875)%
    \lineheight{1}%
    \setlength\tabcolsep{0pt}%
    \put(0,0){\includegraphics[width=\unitlength,page=1]{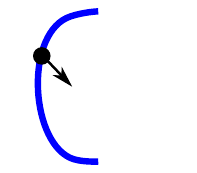}}%
    \put(-0.00171479,0.57550812){\color[rgb]{0,0,0}\makebox(0,0)[lt]{\lineheight{1.25}\smash{\begin{tabular}[t]{l}$x_1$\end{tabular}}}}%
    \put(0,0){\includegraphics[width=\unitlength,page=2]{biangle_I-web0_proof1.pdf}}%
    \put(0.00287621,0.25922147){\color[rgb]{0,0,0}\makebox(0,0)[lt]{\lineheight{1.25}\smash{\begin{tabular}[t]{l}$x_2$\end{tabular}}}}%
    \put(0,0){\includegraphics[width=\unitlength,page=3]{biangle_I-web0_proof1.pdf}}%
    \put(0.85239343,0.58058091){\color[rgb]{0,0,0}\makebox(0,0)[lt]{\lineheight{1.25}\smash{\begin{tabular}[t]{l}$y_1$\end{tabular}}}}%
    \put(0.85698442,0.26429425){\color[rgb]{0,0,0}\makebox(0,0)[lt]{\lineheight{1.25}\smash{\begin{tabular}[t]{l}$y_2$\end{tabular}}}}%
    \put(0,0){\includegraphics[width=\unitlength,page=4]{biangle_I-web0_proof1.pdf}}%
    \put(0.28192897,0.63706932){\color[rgb]{0,0,0}\makebox(0,0)[lt]{\lineheight{1.25}\smash{\begin{tabular}[t]{l}$B_1$\end{tabular}}}}%
    \put(0.52831606,0.64456034){\color[rgb]{0,0,0}\makebox(0,0)[lt]{\lineheight{1.25}\smash{\begin{tabular}[t]{l}$B_2$\end{tabular}}}}%
    \put(0,0){\includegraphics[width=\unitlength,page=5]{biangle_I-web0_proof1.pdf}}%
    \put(0.47969217,0.30896197){\color[rgb]{0,0,0}\makebox(0,0)[lt]{\lineheight{1.25}\smash{\begin{tabular}[t]{l}$z$\end{tabular}}}}%
    \put(0,0){\includegraphics[width=\unitlength,page=6]{biangle_I-web0_proof1.pdf}}%
  \end{picture}%
\endgroup%
}
\end{align*}
Choose an arbitrary direction on $B$, to make it $\vec{B}$; then $B_1,B_2$ naturally inherit directions too. By the property (BT1) of ${\rm Tr}^\omega_B : \mathcal{S}^\omega_{\rm s}(\vec{B};\mathbb{Z}) \to \mathbb{Z}[\omega^{\pm 1}]$,
\begin{align}
\label{eq:I-web_proof1}
{\rm Tr}^\omega_B([W,s]) = {\textstyle \sum}_{s_1,s_2} {\rm Tr}^\omega_B([W_1,s_1]) {\rm Tr}^\omega_B([W_2,s_2]),
\end{align}
where the sum is over all states $s_1,s_2$ compatible with $s$, in the sense of (BT1). Let $\varepsilon = s_1(z)=s_2(z)$. Then, in view of (BT2-3), ${\rm Tr}^\omega_B ([W_1,s_1]) {\rm Tr}^\omega_B ([W_2,s_2]) \neq 0$ iff $\{s_1(x_1),s_1(x_2)\} = \{s_2(y_1),s_2(y_2)\} = \{r_1(\varepsilon),r_2(\varepsilon)\}$. Conversely, ${\rm Tr}^\omega_B([W,s]) \neq 0$ iff $\{s(x_1),s(x_2)\} = \{s(y_1),s(y_2)\} = \{r_1(\varepsilon),r_2(\varepsilon)\}$ holds for some $\varepsilon\in\{1,2,3\}$; such an $\varepsilon$ is unique, if exists, hence the sum in the right hand side of eq.\eqref{eq:I-web_proof1} has only one nonzero summand, for $s_1(z)=s_2(z)=\varepsilon$. The values can be obtained by taking products of values given in (BT2-3).

\vs

To prove eq.\eqref{eq:H-web0}, cut the biangle into two biangles $B_1,B_2$ along an arc $e$ connecting the marked points that meets the ${\rm SL}_3$-web $W'$ in question exactly once transversally, say at $z_1,z_2,z_3$, to express $W'$ as composition of $W_1' = W' \cap (B_1 \times {\bf I})$ and $W_2' = W' \cap (B_2 \times {\bf I})$. 
\begin{align*}
W'
=  \raisebox{-0.5\height}{
\begingroup%
  \makeatletter%
  \providecommand\color[2][]{%
    \errmessage{(Inkscape) Color is used for the text in Inkscape, but the package 'color.sty' is not loaded}%
    \renewcommand\color[2][]{}%
  }%
  \providecommand\transparent[1]{%
    \errmessage{(Inkscape) Transparency is used (non-zero) for the text in Inkscape, but the package 'transparent.sty' is not loaded}%
    \renewcommand\transparent[1]{}%
  }%
  \providecommand\rotatebox[2]{#2}%
  \newcommand*\fsize{\dimexpr\f@size pt\relax}%
  \newcommand*\lineheight[1]{\fontsize{\fsize}{#1\fsize}\selectfont}%
  \ifx\svgwidth\undefined%
    \setlength{\unitlength}{73.7007874bp}%
    \ifx\svgscale\undefined%
      \relax%
    \else%
      \setlength{\unitlength}{\unitlength * \real{\svgscale}}%
    \fi%
  \else%
    \setlength{\unitlength}{\svgwidth}%
  \fi%
  \global\let\svgwidth\undefined%
  \global\let\svgscale\undefined%
  \makeatother%
  \begin{picture}(1,0.90384615)%
    \lineheight{1}%
    \setlength\tabcolsep{0pt}%
    \put(0,0){\includegraphics[width=\unitlength,page=1]{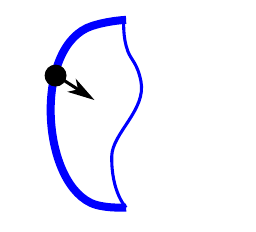}}%
    \put(0.01138983,0.5937226){\color[rgb]{0,0,0}\makebox(0,0)[lt]{\lineheight{1.25}\smash{\begin{tabular}[t]{l}$x_1$\end{tabular}}}}%
    \put(0,0){\includegraphics[width=\unitlength,page=2]{biangle_H-web0_proof1.pdf}}%
    \put(0.01580425,0.2896008){\color[rgb]{0,0,0}\makebox(0,0)[lt]{\lineheight{1.25}\smash{\begin{tabular}[t]{l}$x_2$\end{tabular}}}}%
    \put(0,0){\includegraphics[width=\unitlength,page=3]{biangle_H-web0_proof1.pdf}}%
    \put(0.83264774,0.59860027){\color[rgb]{0,0,0}\makebox(0,0)[lt]{\lineheight{1.25}\smash{\begin{tabular}[t]{l}$y_1$\end{tabular}}}}%
    \put(0.83706215,0.29447839){\color[rgb]{0,0,0}\makebox(0,0)[lt]{\lineheight{1.25}\smash{\begin{tabular}[t]{l}$y_2$\end{tabular}}}}%
    \put(0,0){\includegraphics[width=\unitlength,page=4]{biangle_H-web0_proof1.pdf}}%
    \put(0.29708948,0.66227986){\color[rgb]{0,0,0}\makebox(0,0)[lt]{\lineheight{1.25}\smash{\begin{tabular}[t]{l}$B_1$\end{tabular}}}}%
    \put(0.53968903,0.66247472){\color[rgb]{0,0,0}\makebox(0,0)[lt]{\lineheight{1.25}\smash{\begin{tabular}[t]{l}$B_2$\end{tabular}}}}%
    \put(0,0){\includegraphics[width=\unitlength,page=5]{biangle_H-web0_proof1.pdf}}%
    \put(0.38818795,0.20766551){\color[rgb]{0,0,0}\makebox(0,0)[lt]{\lineheight{1.25}\smash{\begin{tabular}[t]{l}$z_3$\end{tabular}}}}%
    \put(0,0){\includegraphics[width=\unitlength,page=6]{biangle_H-web0_proof1.pdf}}%
    \put(0.5236082,0.38163459){\color[rgb]{0,0,0}\makebox(0,0)[lt]{\lineheight{1.25}\smash{\begin{tabular}[t]{l}$z_2$\end{tabular}}}}%
    \put(0.58056282,0.48441015){\color[rgb]{0,0,0}\makebox(0,0)[lt]{\lineheight{1.25}\smash{\begin{tabular}[t]{l}$z_1$\end{tabular}}}}%
  \end{picture}%
\endgroup%
}
\end{align*}
Choose the elevation ordering on $z_1,z_2,z_3$ such that it is compatible with that on $x_1,x_2$ and on $y_1,y_2$. Namely, if $x_1\succ x_2$, $y_1\succ y_2$, then choose $z_1\succ z_2 \succ z_3$, and if $x_1 \prec x_2$, $y_1 \prec y_2$, then choose $z_1 \prec z_2 \prec z_3$. Note that each of $W_1'$ and $W_2'$ is a product of an edge component as in (BT2-1) and a 3-way component as in (BT2-3). Choose an arbitrary direction on $B$, to make it $\vec{B}$; so $B_1,B_2$ also naturally get directions. By the property (BT1) of ${\rm Tr}^\omega_B : \mathcal{S}^\omega_{\rm s}(\vec{B};\mathbb{Z}) \to \mathbb{Z}[\omega^{\pm 1}]$,
\begin{align}
\label{eq:H-web_proof1}
{\rm Tr}^\omega_B([W',s']) = {\textstyle \sum}_{s_1',s_2'} {\rm Tr}^\omega_B([W_1',s_1']) {\rm Tr}^\omega_B([W_2',s_2']),
\end{align}
where the sum is over all states $s_1',s_2'$ compatible with $s'$, in the sense of (BT1). Let $s_1',s_2'$ compatible states of $W_1',W_2'$ that has non-zero contribution to the sum in eq.\eqref{eq:H-web_proof1}. By multiplicativity of ${\rm Tr}^\omega_B$, the value at each component of $W_1',W_2'$ is nonzero. The edge components fall into (BT2-1), so $s_1'(x_2) = s_1'(z_3)$ and $s_2'(z_1) = s_2'(y_1)$. By compatibility, $s_1'(x_2) = s'(x_2)$, $s_1'(z_3) = s_2'(z_3)$, $s_2'(z_1) = s_1'(z_1)$, and $s_2'(y_1) = s'(y_1)$. Suppose  $s'(x_1)=2=s_1'(x_1)$. In view of (BT2-3) applied to the 3-way component of $W_1'$ we have $\{s_1'(z_1), s_1'(z_2)\} = \{1,3\}$. If $(s_1'(z_1),s_1'(z_2))=(1,3)$, then $1=s_1'(z_1)=s_2'(z_1)=s_2'(y_1)=s'(y_1)$, and $3 = s_1'(z_2) = s_2'(z_2)$. In view of (BT2-3) for the 3-way component of $W_2'$, it follows that either $s_2'(z_3)=1$ holds, in which case $s'(x_2)=1$, $s'(y_2)=2$, or $s_2'(z_3)=2$ holds, in which case $s'(x_2)=2$, $s'(y_2)=3$. If $(s_1'(z_1),s_1'(z_2))=(3,1)$, we have $s'(y_1)=3$ similarly as before, and from (BT2-3) for the 3-way component of $W_2'$, it follows that either $s_2'(z_3)=2$ holds, in which case $s'(x_2)=2$, $s'(y_2) = 1$, or $s_2'(z_3)=3$ holds, in which case $s'(x_2)=3$, $s'(y_2)=2$. Now suppose $s'(x_1)=1=s_1'(x_1)$. Then, in view of (BT2-3) for the 3-way component of $W_1'$ we have $\{s_1'(z_1),s_1'(z_2)\} = \{1,2\}$. By similar arguments as above, if $(s_1'(z_1),s_1'(z_2))=(1,2)$, then it follows that either $s_1'(z_3)=s_2'(z_3)=1$ holds, in which case $s'(x_2)=1$, $s'(y_1)=1$, $s'(y_2)=1$, or $s_1'(z_3)=s_2'(z_3)=3$ holds, in which case $s'(x_2)=3$, $s'(y_1)=1$, $s'(y_3)=3$. If $(s_1'(z_2),s_1'(z_2))=(2,1)$, then it follows that either $s_1'(z_3)=s_2'(z_3)=2$ holds, in which case $s'(x_2)=2$, $s'(y_1)=2$, $s'(y_2)=1$, or $s_1'(z_3)=s_2'(z_3)=3$ holds, in which case $s'(x_2)=3$, $s'(y_1)=2$, $s'(y_2)=2$. This, time suppose $s'(x_1)=3=s_1'(x_1)$. Then, in view of (BT2-3) for the 3-way component of $W_1'$ we have $\{s_1'(z_1),s_1'(z_2)\} = \{2,3\}$. By similar arguments, if $(s_1'(z_1),s_1'(z_2))=((2,3)$, then it follows that either $s_1'(z_3)=s_2'(z_3)=1$ holds, in which case $s'(x_2)=1$, $s'(y_1)=2$, $s'(y_2)=2$, or $s_1'(z_3)=s_2'(z_3)=2$ holds, in which case $s'(x_2)=2$, $s'(y_1)=2$, $s'(y_2)=3$. If $(s_1'(z_1),s_1'(z_2))=(3,2)$, then it follows that either $s_1'(z_3)=s_2'(z_3)=1$ holds, in which case $s'(x_2)=1$, $s'(y_1)=3$, $s'(y_2)=1$, or $s_1'(z_3)=s_2'(z_3)=3$ holds, in which case $s'(x_2)=3$, $s'(y_1)=3$, $s'(y_2)=3$. In these cases, note that there is only one non-zero contributing term in the sum in eq.\eqref{eq:H-web_proof1}, and the actual values can be computed from (BT2-1) and (BT2-3). \qed \qquad {\it End of proof of Lem.\ref{lem:value_of_I-webs_in_biangle}.}

\vs

{\it Proof of Prop.\ref{prop:biangle_SL3_quantum_trace_some_values}.} Let's prove (BT2-4) for ${\rm Tr}^\omega_B$. In each of eq.\eqref{eq:height_exchange_and_crossing1}, \eqref{eq:height_exchange_and_crossing3} and \eqref{eq:height_exchange_and_crossing5}, the three stated ${\rm SL}_3$-webs are identical elements of $\mathcal{S}^\omega_{\rm s}(\vec{B};\mathbb{Z})_{\rm red}$. Take eq.\eqref{eq:height_exchange_and_crossing1}, and take the middle picture. Applying the ${\rm SL}_3$-skein relation (S8) of Fig.\ref{fig:A2-skein_relations_quantum}, we get
\begin{align*}
{\rm Tr}^\omega_B ([\,\, \raisebox{-0.5\height}{} \,\, ])
= q^{-2/3} \, {\rm Tr}^\omega_B ([\,\, \raisebox{-0.5\height}{
\begingroup%
  \makeatletter%
  \providecommand\color[2][]{%
    \errmessage{(Inkscape) Color is used for the text in Inkscape, but the package 'color.sty' is not loaded}%
    \renewcommand\color[2][]{}%
  }%
  \providecommand\transparent[1]{%
    \errmessage{(Inkscape) Transparency is used (non-zero) for the text in Inkscape, but the package 'transparent.sty' is not loaded}%
    \renewcommand\transparent[1]{}%
  }%
  \providecommand\rotatebox[2]{#2}%
  \newcommand*\fsize{\dimexpr\f@size pt\relax}%
  \newcommand*\lineheight[1]{\fontsize{\fsize}{#1\fsize}\selectfont}%
  \ifx\svgwidth\undefined%
    \setlength{\unitlength}{56.69291339bp}%
    \ifx\svgscale\undefined%
      \relax%
    \else%
      \setlength{\unitlength}{\unitlength * \real{\svgscale}}%
    \fi%
  \else%
    \setlength{\unitlength}{\svgwidth}%
  \fi%
  \global\let\svgwidth\undefined%
  \global\let\svgscale\undefined%
  \makeatother%
  \begin{picture}(1,0.875)%
    \lineheight{1}%
    \setlength\tabcolsep{0pt}%
    \put(0,0){\includegraphics[width=\unitlength,page=1]{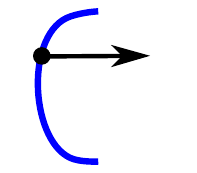}}%
    \put(0.00342324,0.57251168){\color[rgb]{0,0,0}\makebox(0,0)[lt]{\lineheight{1.25}\smash{\begin{tabular}[t]{l}$x_2$\end{tabular}}}}%
    \put(0,0){\includegraphics[width=\unitlength,page=2]{biangle_parallel1.pdf}}%
    \put(0.00801424,0.25622502){\color[rgb]{0,0,0}\makebox(0,0)[lt]{\lineheight{1.25}\smash{\begin{tabular}[t]{l}$x_1$\end{tabular}}}}%
    \put(0.05924823,0.51254405){\color[rgb]{0,0,0}\rotatebox{-92.070464}{\makebox(0,0)[lt]{\lineheight{1.25}\smash{\begin{tabular}[t]{l}$\succ$\end{tabular}}}}}%
    \put(0,0){\includegraphics[width=\unitlength,page=3]{biangle_parallel1.pdf}}%
    \put(0.85239343,0.58058091){\color[rgb]{0,0,0}\makebox(0,0)[lt]{\lineheight{1.25}\smash{\begin{tabular}[t]{l}$y_1$\end{tabular}}}}%
    \put(0.85698442,0.26429425){\color[rgb]{0,0,0}\makebox(0,0)[lt]{\lineheight{1.25}\smash{\begin{tabular}[t]{l}$y_2$\end{tabular}}}}%
    \put(0.89099145,0.50338616){\color[rgb]{0,0,0}\rotatebox{-92.070464}{\makebox(0,0)[lt]{\lineheight{1.25}\smash{\begin{tabular}[t]{l}$\succ$\end{tabular}}}}}%
  \end{picture}%
\endgroup%
} \,\, ])
+ q^{1/3} \, {\rm Tr}^\omega_B ([\,\, \raisebox{-0.5\height}{
\begingroup%
  \makeatletter%
  \providecommand\color[2][]{%
    \errmessage{(Inkscape) Color is used for the text in Inkscape, but the package 'color.sty' is not loaded}%
    \renewcommand\color[2][]{}%
  }%
  \providecommand\transparent[1]{%
    \errmessage{(Inkscape) Transparency is used (non-zero) for the text in Inkscape, but the package 'transparent.sty' is not loaded}%
    \renewcommand\transparent[1]{}%
  }%
  \providecommand\rotatebox[2]{#2}%
  \newcommand*\fsize{\dimexpr\f@size pt\relax}%
  \newcommand*\lineheight[1]{\fontsize{\fsize}{#1\fsize}\selectfont}%
  \ifx\svgwidth\undefined%
    \setlength{\unitlength}{56.69291339bp}%
    \ifx\svgscale\undefined%
      \relax%
    \else%
      \setlength{\unitlength}{\unitlength * \real{\svgscale}}%
    \fi%
  \else%
    \setlength{\unitlength}{\svgwidth}%
  \fi%
  \global\let\svgwidth\undefined%
  \global\let\svgscale\undefined%
  \makeatother%
  \begin{picture}(1,0.875)%
    \lineheight{1}%
    \setlength\tabcolsep{0pt}%
    \put(0,0){\includegraphics[width=\unitlength,page=1]{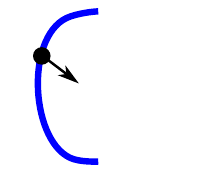}}%
    \put(-0.00171479,0.57550812){\color[rgb]{0,0,0}\makebox(0,0)[lt]{\lineheight{1.25}\smash{\begin{tabular}[t]{l}$x_2$\end{tabular}}}}%
    \put(0,0){\includegraphics[width=\unitlength,page=2]{biangle_I-web2.pdf}}%
    \put(0.00287621,0.25922147){\color[rgb]{0,0,0}\makebox(0,0)[lt]{\lineheight{1.25}\smash{\begin{tabular}[t]{l}$x_1$\end{tabular}}}}%
    \put(0,0){\includegraphics[width=\unitlength,page=3]{biangle_I-web2.pdf}}%
    \put(0.85239343,0.58058091){\color[rgb]{0,0,0}\makebox(0,0)[lt]{\lineheight{1.25}\smash{\begin{tabular}[t]{l}$y_1$\end{tabular}}}}%
    \put(0.85698442,0.26429425){\color[rgb]{0,0,0}\makebox(0,0)[lt]{\lineheight{1.25}\smash{\begin{tabular}[t]{l}$y_2$\end{tabular}}}}%
    \put(0,0){\includegraphics[width=\unitlength,page=4]{biangle_I-web2.pdf}}%
    \put(0.03252493,0.50919727){\color[rgb]{0,0,0}\rotatebox{-90.454342}{\makebox(0,0)[lt]{\lineheight{1.25}\smash{\begin{tabular}[t]{l}$\succ$\end{tabular}}}}}%
    \put(0.87905066,0.5151902){\color[rgb]{0,0,0}\rotatebox{-90.454342}{\makebox(0,0)[lt]{\lineheight{1.25}\smash{\begin{tabular}[t]{l}$\succ$\end{tabular}}}}}%
  \end{picture}%
\endgroup%
} \,\, ]).
\end{align*}
The first term in the right hand side is $q^{-2/3}$ times a product of two cases of (BT2-1), hence equals $q^{-2/3}$ if $s(x_1)=s(y_2)$, $s(x_2)=s(y_1)$ and equals zero otherwise. The value of the second term can be read from eq.\eqref{eq:I-web0} of Lem.\ref{lem:value_of_I-webs_in_biangle}, and one can verify eq.\eqref{eq:height_exchange_and_crossing1_values}; below, we write the value as $((s(x_1),s(x_2)),(s(y_1),s(y_2)))$-th entry of the matrix $\wh{\bf C}_+$, in the style as in Lem.\ref{lem:value_of_I-webs_in_biangle}:
\begin{align*}
\wh{\bf C}_+
= 
\left( \begin{smallmatrix}
q^{-2/3} & 0 & 0 & 0 & 0 & 0 & 0 & 0 & 0  \\
0 & q^{1/3} & 0 & q^{-2/3}-q^{4/3}& 0 & 0 & 0 & 0 & 0 \\
0 & 0 & q^{1/3} & 0 & 0 & 0 & q^{-2/3}-q^{4/3}  & 0 & 0 \\
0 & 0 & 0 & q^{1/3} & 0 & 0 & 0 & 0 & 0 \\
0 & 0 & 0 & 0 & q^{-2/3} & 0 & 0 & 0 & 0  \\
0 & 0 & 0 & 0 & 0 & q^{1/3} & 0 & q^{-2/3}-q^{4/3} & 0 \\
0 & 0 & 0 & 0 & 0 & 0 & q^{1/3} & 0 & 0 \\
0 & 0 & 0 & 0 & 0 & 0 & 0 & q^{1/3} & 0 \\
0 & 0 & 0 & 0 & 0 & 0 & 0 & 0 & q^{-2/3} \\
\end{smallmatrix}
\right)
\end{align*}

\vs

Now take the rightmost picture of eq.\eqref{eq:height_exchange_and_crossing3}. Applying the ${\rm SL}_3$-skein relation (S9) of Fig.\ref{fig:A2-skein_relations_quantum}, we get
\begin{align}
\label{eq:height_exchange_and_crossing3_proof}
{\rm Tr}^\omega_B ([\,\, \raisebox{-0.5\height}{} \,\, ])
= q^{2/3} \, {\rm Tr}^\omega_B ([\,\, \raisebox{-0.5\height}{
\begingroup%
  \makeatletter%
  \providecommand\color[2][]{%
    \errmessage{(Inkscape) Color is used for the text in Inkscape, but the package 'color.sty' is not loaded}%
    \renewcommand\color[2][]{}%
  }%
  \providecommand\transparent[1]{%
    \errmessage{(Inkscape) Transparency is used (non-zero) for the text in Inkscape, but the package 'transparent.sty' is not loaded}%
    \renewcommand\transparent[1]{}%
  }%
  \providecommand\rotatebox[2]{#2}%
  \newcommand*\fsize{\dimexpr\f@size pt\relax}%
  \newcommand*\lineheight[1]{\fontsize{\fsize}{#1\fsize}\selectfont}%
  \ifx\svgwidth\undefined%
    \setlength{\unitlength}{56.69291339bp}%
    \ifx\svgscale\undefined%
      \relax%
    \else%
      \setlength{\unitlength}{\unitlength * \real{\svgscale}}%
    \fi%
  \else%
    \setlength{\unitlength}{\svgwidth}%
  \fi%
  \global\let\svgwidth\undefined%
  \global\let\svgscale\undefined%
  \makeatother%
  \begin{picture}(1,0.875)%
    \lineheight{1}%
    \setlength\tabcolsep{0pt}%
    \put(0,0){\includegraphics[width=\unitlength,page=1]{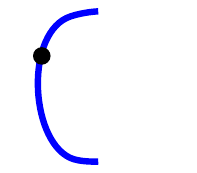}}%
    \put(-0.00171479,0.57550812){\color[rgb]{0,0,0}\makebox(0,0)[lt]{\lineheight{1.25}\smash{\begin{tabular}[t]{l}$x_1$\end{tabular}}}}%
    \put(0,0){\includegraphics[width=\unitlength,page=2]{biangle_crossing8_resolved.pdf}}%
    \put(0.00287621,0.25922147){\color[rgb]{0,0,0}\makebox(0,0)[lt]{\lineheight{1.25}\smash{\begin{tabular}[t]{l}$x_2$\end{tabular}}}}%
    \put(0.03772366,0.5151652){\color[rgb]{0,0,0}\rotatebox{-90.362304}{\makebox(0,0)[lt]{\lineheight{1.25}\smash{\begin{tabular}[t]{l}$\prec$\end{tabular}}}}}%
    \put(0,0){\includegraphics[width=\unitlength,page=3]{biangle_crossing8_resolved.pdf}}%
    \put(0.85239343,0.58058091){\color[rgb]{0,0,0}\makebox(0,0)[lt]{\lineheight{1.25}\smash{\begin{tabular}[t]{l}$y_2$\end{tabular}}}}%
    \put(0.85698442,0.26429425){\color[rgb]{0,0,0}\makebox(0,0)[lt]{\lineheight{1.25}\smash{\begin{tabular}[t]{l}$y_1$\end{tabular}}}}%
    \put(0.89099145,0.50338616){\color[rgb]{0,0,0}\rotatebox{-92.070464}{\makebox(0,0)[lt]{\lineheight{1.25}\smash{\begin{tabular}[t]{l}$\prec$\end{tabular}}}}}%
    \put(0,0){\includegraphics[width=\unitlength,page=4]{biangle_crossing8_resolved.pdf}}%
  \end{picture}%
\endgroup%
} \,\, ])
+ q^{-1/3} \, {\rm Tr}^\omega_B ([\,\, \raisebox{-0.5\height}{
\begingroup%
  \makeatletter%
  \providecommand\color[2][]{%
    \errmessage{(Inkscape) Color is used for the text in Inkscape, but the package 'color.sty' is not loaded}%
    \renewcommand\color[2][]{}%
  }%
  \providecommand\transparent[1]{%
    \errmessage{(Inkscape) Transparency is used (non-zero) for the text in Inkscape, but the package 'transparent.sty' is not loaded}%
    \renewcommand\transparent[1]{}%
  }%
  \providecommand\rotatebox[2]{#2}%
  \newcommand*\fsize{\dimexpr\f@size pt\relax}%
  \newcommand*\lineheight[1]{\fontsize{\fsize}{#1\fsize}\selectfont}%
  \ifx\svgwidth\undefined%
    \setlength{\unitlength}{56.69291339bp}%
    \ifx\svgscale\undefined%
      \relax%
    \else%
      \setlength{\unitlength}{\unitlength * \real{\svgscale}}%
    \fi%
  \else%
    \setlength{\unitlength}{\svgwidth}%
  \fi%
  \global\let\svgwidth\undefined%
  \global\let\svgscale\undefined%
  \makeatother%
  \begin{picture}(1,0.875)%
    \lineheight{1}%
    \setlength\tabcolsep{0pt}%
    \put(0,0){\includegraphics[width=\unitlength,page=1]{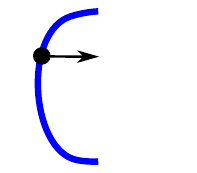}}%
    \put(-0.00171479,0.57550812){\color[rgb]{0,0,0}\makebox(0,0)[lt]{\lineheight{1.25}\smash{\begin{tabular}[t]{l}$x_1$\end{tabular}}}}%
    \put(0,0){\includegraphics[width=\unitlength,page=2]{biangle_H-web2.pdf}}%
    \put(0.00287621,0.25922147){\color[rgb]{0,0,0}\makebox(0,0)[lt]{\lineheight{1.25}\smash{\begin{tabular}[t]{l}$x_2$\end{tabular}}}}%
    \put(0,0){\includegraphics[width=\unitlength,page=3]{biangle_H-web2.pdf}}%
    \put(0.85239343,0.58058091){\color[rgb]{0,0,0}\makebox(0,0)[lt]{\lineheight{1.25}\smash{\begin{tabular}[t]{l}$y_2$\end{tabular}}}}%
    \put(0.85698442,0.26429425){\color[rgb]{0,0,0}\makebox(0,0)[lt]{\lineheight{1.25}\smash{\begin{tabular}[t]{l}$y_1$\end{tabular}}}}%
    \put(0,0){\includegraphics[width=\unitlength,page=4]{biangle_H-web2.pdf}}%
    \put(0.03227797,0.51515439){\color[rgb]{0,0,0}\rotatebox{-90.308736}{\makebox(0,0)[lt]{\lineheight{1.25}\smash{\begin{tabular}[t]{l}$\prec$\end{tabular}}}}}%
    \put(0.8832983,0.51365622){\color[rgb]{0,0,0}\rotatebox{-90.308736}{\makebox(0,0)[lt]{\lineheight{1.25}\smash{\begin{tabular}[t]{l}$\prec$\end{tabular}}}}}%
  \end{picture}%
\endgroup%
} \,\, ]).
\end{align}
The first term in the right hand side is $q^{2/3}$ times a product of two cases of (BT2-2), while the value of the second term can be read from eq.\eqref{eq:H-web0} of Lem.\ref{lem:value_of_I-webs_in_biangle}; the value is $((s(x_1),s(x_2)),(s(y_1),s(y_2)))$-th entry of the matrix $\wh{\bf C}_-$
\begin{align*}
\wh{\bf C}_- = 
\left( \begin{smallmatrix}
q^{-1/3} & 0 & 0 & 0 & 0  & 0 & 0 & 0 & 0 \\
0 & q^{-1/3} & 0 & 0 & 0 & 0 & 0 & 0 & 0 \\
0 & 0 & q^{2/3} & 0 & q^{-1/3}-q^{5/3} & 0 & q^{8/3}-q^{2/3} & 0 & 0 \\
0 & 0 & 0 & q^{-1/3} & 0 & 0 & 0 & 0 & 0 \\
0 & 0 & 0 & 0 & q^{2/3} & 0 & q^{-1/3}-q^{5/3} & 0 & 0 \\
0 & 0 & 0 & 0 & 0 & q^{-1/3} & 0 & 0 & 0 \\
0 & 0 & 0 & 0 & 0 & 0 & q^{2/3}  & 0 & 0 \\
0 & 0 & 0 & 0 & 0 & 0 & 0 & q^{-1/3} & 0 \\
0 & 0 & 0 & 0 & 0 & 0 & 0 & 0 & q^{-1/3} 
\end{smallmatrix}
\right)
\end{align*}
which one can match with eq.\eqref{eq:height_exchange_and_crossing3_values}. If we change the elevation ordering on each side of the biangles appearing in eq.\eqref{eq:height_exchange_and_crossing3_proof} and then exchange the labels of endpoints as $x_1\leftrightarrow x_2$ and $y_1 \leftrightarrow y_2$, the left hand side becomes the middle picture of eq.\eqref{eq:height_exchange_and_crossing5}; so from (BT2-2) and eq.\eqref{eq:H-web0} of Lem.\ref{lem:value_of_I-webs_in_biangle}, one gets eq.\eqref{eq:height_exchange_and_crossing5} by computation.  \qed \qquad {\it End of proof of Prop.\ref{prop:biangle_SL3_quantum_trace_some_values}.}

\vs

The values of ${\rm Tr}^\omega_B$ for the remaining elementary cases in eq.\eqref{eq:some_more_elementary1}--\eqref{eq:some_more_elementary3} can be obtained by taking `inverses' of the cases in (BT2-4), as mentioned before. Writing down the values of ${\rm Tr}^\omega_B$ for the cases in (BT2-4) as entries of a $9\times 9$ matrix, it is particularly easy to take the inverse, as the matrices $\wh{\bf C}_+$ and $\wh{\bf C}_-$ are upper triangular. For each $k=+,-$, one finds out that the inverse of $\wh{\bf C}_k$ is obtained from $\wh{\bf C}_k$ by replacing $q$ by $q^{-1}$ in the entries. That is, the entries $q^{\pm 2/3}$, $q^{\pm 1/3}$, $q^{-2/3}-q^{4/3}$, $q^{-1/3}-q^{5/3}$ and $q^{8/3}-q^{2/3}$ are replaced by $q^{\mp 2/3}$, $q^{\mp 1/3}$, $q^{2/3}-q^{-4/3}$, $q^{1/3}-q^{-5/3}$ and $q^{-8/3}-q^{-2/3}$ respectively. We will see that these results for the inverses can be obtained more conceptually (see Lem.\ref{lem:elevation_reversing_and_star-structure_for_biangles}).

\subsection{The state-sum construction}
\label{subsec:state-sum_construction}

To tackle Thm.\ref{thm:SL3_quantum_trace_map} we provide an explicit formula for computation of the value ${\rm Tr}^\omega_{\Delta;\frak{S}}([W,s])$ of the sought-for quantum trace map ${\rm Tr}^\omega_{\Delta;\frak{S}}$ for a triangulated surface. Like in Bonahon-Wong's argument \cite{BW} for the ${\rm SL}_2$ quantum trace map, we consider the split ideal triangulation $\wh{\Delta}$ for $\Delta$ (Def.\ref{def:split_ideal_triangulation}), put the complicated parts of $W$ into the biangles of $\wh{\Delta}$ by isotopy, and use the biangle ${\rm SL}_3$ quantum trace to deal with these parts. We begin by defining the $A_2$ analog of Bonahon-Wong's {\em good position} (\cite{BW}) for skeins.
\begin{definition}[good position and gool position]
\label{def:good_position}
Let $\frak{S}$ be a triangulable generalized marked surface, $\Delta$ an ideal triangulation of $\frak{S}$, and $\wh{\Delta}$ a split ideal triangulation for $\Delta$ (Def.\ref{def:split_ideal_triangulation}). An ${\rm SL}_3$-web $W$ in $\frak{S}\times {\bf I}$ is said to be in a \ul{\em good position} with respect to $\wh{\Delta}$ if it satisfies the following conditions:
\begin{enumerate}
\itemsep0em
\item[\rm (GP1)] For each triangle $\wh{t}$ of $\wh{\Delta}$ and biangle $B$ of $\wh{\Delta}$, when $\wh{t}$ and $B$ are viewed as generalized marked surfaces on their own, $W\cap (\wh{t}\times {\bf I})$ is an ${\rm SL}_3$-web in $\wh{t}\times {\bf I}$ and $W\cap (B\times {\bf I})$ is an ${\rm SL}_3$-web in $B\times {\bf I}$. 

\item[\rm (GP2)] For each triangle $\wh{t}$ of $\wh{\Delta}$, if we denote the components of the ${\rm SL}_3$-web $W\cap (\wh{t}\times {\bf I})$ in $\wh{t} \times {\bf I}$ by $W_{t,1},\ldots,W_{t,l_t}$, and denote by $I_{t,j} \subset {\bf I}$ the image of $W_{t,j}$ under the second projection $\wh{t} \times {\bf I} \to {\bf I}$, then $I_{t,1},\ldots,I_{t,l_t}$ are mutually disjoint.

\item[\rm (GP3)] For each triangle $\wh{t}$ of $\wh{\Delta}$, each component of $W\cap (\wh{t}\times {\bf I})$ has no crossing and at most one 3-valent internal vertex; if a component of $W\cap (\wh{t}\times {\bf I})$ has one 3-valent vertex, then the three endpoints of this component do not lie over a single side of $\wh{t}$.
\end{enumerate}
If furthermore the following also holds, we say $W$ is in a \ul{\em gool position}\footnote{`Gool' sounds like honey (or oyster!) in Korean.} with respect to $\wh{\Delta}$: 
\begin{enumerate}
\item[\rm (GP4)] For each triangle $\wh{t}$ of $\wh{\Delta}$, each component of the ${\rm SL}_3$-web $W\cap (\wh{t} \times {\bf I})$ over $\wh{t}$ is a corner arc (cf. Def.\ref{def:canonical_web_in_a_triangle}), i.e. is an arc that has no crossing and connects two distinct components of $\partial \wh{t} \times {\bf I}$.
\end{enumerate}
\end{definition}

By isotopy, one can push all (or almost \redfix{all}) of the 3-valent vertices into the biangles.
\begin{lemma}
\label{lem:isotopic_to_good_position}
An ${\rm SL}_3$-web $W$ in $\frak{S} \times {\bf I}$ is isotopic to an ${\rm SL}_3$-web $W'$ in a gool position with respect to $\wh{\Delta}$. \qed
\end{lemma}

\begin{definition}[state-sum trace for a gool position]
\label{def:state-sum_trace_for_gool_position}
Let $\frak{S}$, $\Delta$, and $\wh{\Delta}$ be as in Def.\ref{def:good_position}. Let $(W,s)$ be a stated ${\rm SL}_3$-web in $\frak{S}\times {\bf I}$ in a gool position with respect to $\wh{\Delta}$. 

\vs

The points of $W\cap (\wh{\Delta} \times {\bf I})$ are called \ul{\em $\wh{\Delta}$-junctures} of $W$, and a \ul{\em $\wh{\Delta}$-juncture-state} of $W$ is a map
$$
J : W\cap (\wh{\Delta} \times {\bf I}) \to \{1,2,3\}.
$$
For each ideal triangle $\wh{t}$ of $\wh{\Delta}$ corresponding to a triangle $t$ of $\Delta$, view $W\cap (\wh{t}\times {\bf I})$ as an ${\rm SL}_3$-web in $\wh{t}\times {\bf I}$, where $\wh{t}$ is viewed as a generalized marked surface on its own. Let $W_{t,1}, \ldots, W_{t,l_t}$ be components of this ${\rm SL}_3$-web $W\cap (\wh{t}\times {\bf I})$, in decreasing order of elevations. A $\wh{\Delta}$-juncture-state $J$ of $W$ restricts to a state $J_t : \partial (W\cap (\wh{t}\times {\bf I})) \to \{1,2,3\}$ for the external vertices of $W\cap (\wh{t}\times {\bf I})$, and also to a state $J_{t,j} : \partial W_{t,j} \to \{1,2,3\}$ for the external vertices of $W_{t,j}$. To each stated component $(W_{t,j}, J_{t,j})$, i.e. a pair of a component $W_{t,j}$ and a state for its external vertices, define the element
\begin{align}
\label{eq:value_of_stated_component}
\wh{\rm Tr}^\omega_t (W_{t,j},J_{t,j}) ~\in~ \mathcal{Z}^\omega_t
\end{align}
as in Thm.\ref{thm:SL3_quantum_trace_map}\redfix{(QT2-1)}--\redfix{(QT2-2)}. Define the \ul{\em triangle factor} of $W$ for the triangle $t$ (or $\wh{t}$\,) with respect to $J$ as
\begin{align}
\label{eq:triangle_factor}
\wh{\rm Tr}^\omega_t(W\cap (\wh{t}\times{\bf I}), J_t) := {\textstyle \overrightarrow{\prod}}_{j=1}^{l_t} \wh{\rm Tr}^\omega_t (W_{t,j},J_{t,j})
:= \wh{\rm Tr}^\omega_t (W_{t,1},J_{t,1})  \cdot \cdots \cdot \wh{\rm Tr}^\omega_t (W_{t,l_t},J_{t,l_t}) ~\in~ \mathcal{Z}^\omega_t
\end{align}

\vs

For each biangle $B$ of $\wh{\Delta}$, view $W\cap (B \times {\bf I})$ as an ${\rm SL}_3$-web in $B \times {\bf I}$. The $\wh{\Delta}$-juncture-state $J$ of $W$ restricts to a state $J_B : \partial (W\cap (B\times {\bf I})) \to \{1,2,3\}$. Let the \ul{\em biangle factor} of $W$ for the biangle $B$ with respect to $J$ be
$$
{\rm Tr}^\omega_B([W\cap (B\times {\bf I}), J_B])  ~\in~ \mathbb{Z}[\omega^{\pm1/2}],
$$
as given by Prop.\ref{prop:biangle_SL3_quantum_trace}.

\vs

Define the \ul{\em state-sum trace} of the stated ${\rm SL}_3$-web $(W,s)$ in a gool position with respect to $\wh{\Delta}$ as
\begin{align}
\label{eq:state-sum_formula}
\wh{\rm Tr}^\omega_\Delta(W,s) := {\textstyle \sum}_J ( {\textstyle \prod}_B {\rm Tr}^\omega_B([W\cap (B\times {\bf I}), J_B]) \, {\textstyle \bigotimes}_t \wh{\rm Tr}^\omega_t(W\cap (\wh{t} \times {\bf I}), J_t) ) ~\in~ 
{\textstyle \bigotimes}_{t\in \mathcal{F}(\Delta)} \mathcal{Z}^\omega_t
\end{align}
where the sum $\sum_J$ is over all $\wh{\Delta}$-juncture-states $J$ for $W$ that restrict to $s$ at $\partial W$, and the products $\prod_B$ and $\bigotimes_t$ are over all biangles $B$ of $\wh{\Delta}$ and triangles $t$ of $\Delta$.

\end{definition}

To use the state-sum trace as the sought-for ${\rm SL}_3$ quantum trace, one must show that the value $\wh{\rm Tr}^\omega_\Delta(W,s)$ lies in the subalgebra $\mathcal{Z}^\omega_\Delta$ of $\bigotimes_{t\in \mathcal{F}(\Delta)} \mathcal{Z}^\omega_t$ (Def.\ref{def:Fock-Goncharov_algebra_quantum}).
\begin{proposition}[balancedness of the state-sum quantum trace]
\label{prop:balancedness_of_state-sum_trace}
Let $\frak{S}$, $\Delta$, and $\wh{\Delta}$ be as in Def.\ref{def:good_position}. For a stated ${\rm SL}_3$-web $(W,s)$ in $\frak{S}\times{\bf I}$ in a gool position with respect to $\wh{\Delta}$,
$$
\wh{\rm Tr}^\omega_\Delta(W,s) ~\in~ \mathcal{Z}^\omega_\Delta ~\subset~{\textstyle \bigotimes}_{t\in \mathcal{F}(\Delta)} \mathcal{Z}^\omega_t.
$$
\end{proposition}
To prove this, we first establish the following lemma, which is interesting in its own right, and is an ${\rm SL}_3$ analog of the corresponding statement for ${\rm SL}_2$.
\begin{lemma}[charge conservation property of the biangle ${\rm SL}_3$ quantum trace]
\label{lem:charge_conservation}
Let $\vec{B}$ be a directed biangle, and $(W,s)$ be a stated ${\rm SL}_3$-web in $\vec{B}\times{\bf I}$. Let $b_{\rm left},b_{\rm right}$ be the left and the right sides of $\vec{B}$ (Def.\ref{def:direction}). For $\varepsilon \in \{1,2,3\}$, let $n^+_{{\rm left},\varepsilon}$ (resp. $n^-_{{\rm left},\varepsilon}$) be the number of endpoints $x$ of $W$ lying over $b_{\rm left}$ with $s(x)=\varepsilon$ such that $x$ is a source (resp. sink) of $W$, i.e. the orientation of $W$ near $x$ is going toward (resp. away from) the interior of $B$, or equivalently, going from left to right (resp. right to left). For $\varepsilon \in \{1,2,3\}$, let $n^+_{{\rm right},\varepsilon}$ (resp. $n^-_{{\rm right},\varepsilon}$) be the number of endpoints $x$ of $W$ lying over $b_{\rm right}$ with $s(x)=\varepsilon$ such that $x$ is a sink (resp. source) of $W$, i.e. the orientation of $W$ near $x$ is going from left to right (resp. right to left). 

\vs

For $h\in \{{\rm left},{\rm right}\}$, the \ul{\em first charge} of $(W,s)$ at the boundary arc $b_h$ is defined as
\begin{align}
\label{eq:first_charge}
\mathcal{C}^{(1)}_{h}(W,s) = n^+_{h,1} - n^+_{h,3} + n^-_{h,1} - n^-_{h,3} ~\in~ \mathbb{Z}.
\end{align}
and the \ul{\em second charge} of $(W,s)$ at boundary arc $b_h$ as
\begin{align}
\label{eq:second_charge}
\mathcal{C}^{(2)}_{h}(W,s) = n^+_{h,1} - 2n^+_{h,2} + n^+_{h,3} - n^-_{h,1} + 2n^-_{h,2} - n^-_{h,3} ~\in~ \mathbb{Z}.
\end{align}

\vs

If ${\rm Tr}^\omega_B([W,s]) \neq 0$, then
\begin{align*}
\mathcal{C}^{(1)}_{\rm left}(W,s)  = \mathcal{C}^{(1)}_{\rm right}(W,s) \quad\mbox{and}\quad
\mathcal{C}^{(2)}_{\rm left}(W,s)  = \mathcal{C}^{(2)}_{\rm right}(W,s).
\end{align*}

\end{lemma}

\begin{corollary}
\label{cor:charge_conservation}
Define
\begin{align*}
\mathcal{C}^{(3)}_h(W,s) & := \textstyle \frac{3}{2} \mathcal{C}^{(1)}_h(W,s) - \frac{1}{2} \mathcal{C}^{(2)}_h(W,s)
= n^+_{h,1} + n^+_{h,2} - 2 n^+_{h,3}
+ 2 n^-_{h,1} - n^-_{h,2} - n^-_{h,3}, \\
\mathcal{C}^{(4)}_h(W,s) & := \textstyle \frac{3}{2} \mathcal{C}^{(1)}_h(W,s) + \frac{1}{2} \mathcal{C}^{(2)}_h(W,s)
= 2 n^+_{h,1} - n^+_{h,2} - n^+_{h,3}
+ n^-_{h,1} + n^-_{h,2} - 2 n^-_{h,3}.
\end{align*}
If ${\rm Tr}^\omega_B([W,s]) \neq 0$, then
\begin{align*}
\mathcal{C}^{(3)}_{\rm left}(W,s) = \mathcal{C}^{(3)}_{\rm right}(W,s) \quad\mbox{and}\quad \mathcal{C}^{(4)}_{\rm left}(W,s) = \mathcal{C}^{(4)}_{\rm right}(W,s). \qed
\end{align*}
\end{corollary}

The first charge can be understood as
$$
\mathcal{C}^{(1)}_h(W,s) = - \underset{x\in (\partial W) \cap (b_h \times {\bf I})}{\textstyle \sum} 
(s(x)-2),
$$
i.e. minus the sum of {\em signs} for the values of the state $s$ at the endpoints of $W$ lying over $b_h$, where the sign of the state value $\varepsilon \in \{1,2,3\}$ is defined as $\varepsilon -2 \in \{-1,0,+1\}$, matching the convention of Higgins \cite{Higgins} who use the symbols $\{-,0,+\}$ as the values of states.

\vs

{\it Proof of Lem.\ref{lem:charge_conservation}.} Let $(W,s)$ be a stated ${\rm SL}_3$-web in $\vec{B} \times {\bf I}$, such that ${\rm Tr}^\omega_B([W,s])\neq 0$. Recall from Lem.\ref{lem:composition_of_elementary_biangles} that $W$ can be decomposed as composition of elementary ${\rm SL}_3$-webs (Def.\ref{def:elementary_A2-webs_in_biangle}) in thickened biangles. More precisely, there exists a finite collection of ideal arcs $e_1,\ldots,e_n$ of $B$ connecting its two marked points, dividing $\vec{B}$ into directed biangles $\vec{B}_1$,\ldots,$\vec{B}_{n+1}$, appearing in this order from the left side of $\vec{B}$ toward the right side of $\vec{B}$, so that for each $i=1,\ldots,n+1$, the ${\rm SL}_3$-web $W_i := W \cap (\vec{B}_i \times {\bf I})$ in $\vec{B}_i \times {\bf I}$ is elementary in the sense of Def.\ref{def:elementary_A2-webs_in_biangle}.

\vs

Denote the left and the right sides of $\vec{B}_i$ as $b_{{\rm left};i}$ and $b_{{\rm right};i}$. Then $b_{{\rm left};i} = e_{i-1}$ and $b_{{\rm right};i} = e_i$ for each $i=1,\ldots,n+1$, where we denote the left and the right sides $b_{\rm left}$ and $b_{\rm right}$ of $\vec{B}$ by $e_0$ and $e_{n+1}$, respectively. Let $J : W \cap ( (\partial B \cup e_1 \cup \cdots \cup e_n)\times {\bf I}) \to \{1,2,3\}$ be a juncture-state for this decomposition of $B$. By the item (BT1) of Prop.\ref{prop:biangle_SL3_quantum_trace}, we have
\begin{align}
\label{eq:Tr_B_W_s}
{\rm Tr}^\omega_B([W,s]) = {\textstyle \sum}_J {\textstyle \prod}_{i=1}^{n+1} {\rm Tr}^\omega_{B_i}([W_i, J|_{\partial W_i}])
\end{align}
where the sum is over all juncture-states $J$ restricting to $s$ at $\partial W = W\cap (\partial B \times {\bf I})$. Since ${\rm Tr}^\omega_B([W,s])\neq 0$, there exists a juncture-state $J_0$ restricting to $s$ such that the corresponding summand ${\textstyle \prod}_{i=1}^{n+1} {\rm Tr}^\omega_{B_i}([W_i, J_0|_{\partial W_i}])$ is nonzero. For this $J_0$ we therefore have ${\rm Tr}^\omega_{B_i}([W_i, J_0|_{\partial W_i}])\neq 0$ for all $i=1,\ldots,n+1$. 
\vs

For each $i=1,\ldots,n+1$, since $W_i$ is elementary, it equals the product $W_{i,1} \cdot \cdots \cdot W_{i,l_i}$ (as elements of $\mathcal{S}^\omega(\vec{B}_i;\mathbb{Z})$), where each $W_{i,j}$ falls into one of (BT2-1)--(BT2-4) of Prop.\ref{prop:biangle_SL3_quantum_trace}--\ref{prop:biangle_SL3_quantum_trace_some_values} or eq.\eqref{eq:some_more_elementary1}--\eqref{eq:some_more_elementary3}. Note ${\rm Tr}^\omega_{B_i}([W_i,J_0|_{\partial W_i}]) = \prod_{j=1}^{l_i} {\rm Tr}^\omega_{B_i}([W_{i,j}, J_0|_{\partial W_{i,j}}])$, hence ${\rm Tr}^\omega_{B_i}([W_{i,j}, J_0|_{\partial W_{i,j}}])\neq 0$ for all $j=1,\ldots,l_i$. In case $(W_{i,j}, J_0|_{\partial W_{i,j}})$ falls into (BT2-1), ${\rm Tr}^\omega_{B_i}([W_{i,j}, J_0|_{\partial W_{i,j}}])\neq 0$ iff $J_0|_{\partial W_{i,j}}$ assigns same state value $\in\{1,2,3\}$ to the two endpoints of $W_{i,j}$, while the two endpoints lie over distinct sides of $B_i$ (i.e. $e_{i-1}$ and $e_i$), hence $n^+_{e_{i-1},\varepsilon} = n^+_{e_i,\varepsilon}$ and $n^-_{e_{i-1},\varepsilon} = n^-_{e_i,\varepsilon}$ hold for all $\varepsilon\in\{1,2,3\}$, so in view of eq.\eqref{eq:first_charge} and eq.\eqref{eq:second_charge} we can observe
\begin{align}
\label{eq:first_charge_conservation_for_single_component}
\mathcal{C}^{(1)}_{\rm left}(W_{i,j},J_0|_{\partial W_{i,j}}) & = \mathcal{C}^{(1)}_{\rm right}(W_{i,j},J_0|_{\partial W_{i,j}}), \\
\label{eq:second_charge_conservation_for_single_component}
\mathcal{C}^{(2)}_{\rm left}(W_{i,j},J_0|_{\partial W_{i,j}}) & = \mathcal{C}^{(2)}_{\rm right}(W_{i,j},J_0|_{\partial W_{i,j}}).
\end{align}
In case (BT2-2), note ${\rm Tr}^\omega_{B_i}([W_{i,j}, J_0|_{\partial W_{i,j}}])\neq 0$ iff the pair of values of $J_0|_{\partial W_{i,j}}$ at the two endpoints of $W_{i,j}$ is one of $(1,3),(2,2),(3,1)$, while the two endpoints lie in a single side of $B_i$ (i.e. either $e_{i-1}$ or $e_i$), one being a source and the other a sink. So, for the one $h\in \{{\rm left},{\rm right}\}$ for which $W_{i,j}$ has no endpoints on $b_{h;i}$, manifestly $\mathcal{C}^{(1)}_h([W_{i,j}, J_0|_{\partial W_{i,j}}]) = \mathcal{C}^{(2)}_h([W_{i,j}, J_0|_{\partial W_{i,j}}])=0$. For the other $h$, one easily observes $n^+_{h,1}=n^-_{h,3}$, $n^+_{h,2}=n^-_{h,2}$ and $n^+_{h,3}=n^-_{h,1}$, hence 
 $\mathcal{C}^{(1)}_h([W_{i,j}, J_0|_{\partial W_{i,j}}]) = \mathcal{C}^{(2)}_h([W_{i,j}, J_0|_{\partial W_{i,j}}])=0$, and therefore eq.\eqref{eq:first_charge_conservation_for_single_component} and eq.\eqref{eq:second_charge_conservation_for_single_component} hold. In case (BT2-3), where one side of $B_i$ has two endpoints $x_1,x_2$ of $W_{i,j}$ and the other side of $B_i$ has one endpoint $x_3$ of $W_{i,j}$, and the three endpoints are either all sinks or all sources. By (BT2-3), note ${\rm Tr}^\omega_{B_i}([W_{i,j}, J_0|_{\partial W_{i,j}}])\neq 0$ iff $\{J_0(x_1),J_0(x_2)\} = \{r_1(J_0(x_3)), r_2(J_0(x_3))\}$. Suppose $x_1,x_2,x_3$ are sources, and $x_3$ is at $e_{i-1} = b_{{\rm left},i}$. If $J_0(x_3)=1$, then $n^+_{{\rm left},1}=n^-_{{\rm right},1}=n^-_{{\rm right},2}=1$, while the remaining $n^*_{*,*}$ are all zero; so for $(W_{i,j},J_0|_{\partial W_{i,j}})$ we have $\mathcal{C}^{(1)}_{\rm left} = 1$, $\mathcal{C}^{(1)}_{\rm right} = 1$, $\mathcal{C}^{(2)}_{\rm left} = 1$, $\mathcal{C}^{(2)}_{\rm right} = -1+2=1$. If $J_0(x_3)=2$, then $n^+_{{\rm left},2}=n^-_{{\rm right},1} = n^-_{{\rm right},3}=1$ with other $n^*_{*,*}$ being zero, so $\mathcal{C}^{(1)}_{\rm left} = 0$, $\mathcal{C}^{(1)}_{\rm right} = 1-1=0$, $\mathcal{C}^{(2)}_{\rm left} = -2$, $\mathcal{C}^{(2)}_{\rm right} = -1-1=-2$. If $J_0(x_3)=3$, then $n^+_{{\rm left},3}=n^-_{{\rm right},2} = n^-_{{\rm right},3}=1$ with other $n^*_{*,*}$ being zero, so $\mathcal{C}^{(1)}_{\rm left} = -1$, $\mathcal{C}^{(1)}_{\rm right} = -1$, $\mathcal{C}^{(2)}_{\rm left} = 1$, $\mathcal{C}^{(2)}_{\rm right} = 2-1=1$. In any case, eq.\eqref{eq:first_charge_conservation_for_single_component} and eq.\eqref{eq:second_charge_conservation_for_single_component} hold. Proof of eq.\eqref{eq:first_charge_conservation_for_single_component} and eq.\eqref{eq:second_charge_conservation_for_single_component} for the cases when $x_1,x_2,x_3$ may be sinks and $x_3$ may be at $e_i = b_{{\rm right},i}$ follows, due to the symmetry and skew-symmetry of the definition of the charges as in eq.\eqref{eq:first_charge} and eq.\eqref{eq:second_charge} under the exchange $+ \leftrightarrow -$ of the superscripts, and the symmetry of the sought-for eq.\eqref{eq:first_charge_conservation_for_single_component} and eq.\eqref{eq:second_charge_conservation_for_single_component} under $e_{i-1} \leftrightarrow e_i$ (i.e. left$\leftrightarrow$right). Now, suppose $(W_{i,j}, J_0|_{\partial W_{i,j}})$ falls into case (BT2-4), with ${\rm Tr}^\omega_{B_i}([W_{i,j}, J_0|_{\partial W_{i,j}}])\neq 0$. Take eq.\eqref{eq:height_exchange_and_crossing1}, with the upper marked point being the top marked point, so that the two component strands are both going from left to right. So $n^-_{*,*}=0$. By inspection, $n^+_{{\rm left},\varepsilon}=n^+_{{\rm right},\varepsilon}$ for all $\varepsilon\in\{1,2,3\}$. For the case when the lower marked point is the top marked point, we have $n^+_{*,*}=0$ and $n^-_{{\rm left},\varepsilon}=n^-_{{\rm right},\varepsilon}$ for all $\varepsilon\in\{1,2,3\}$. Hence eq.\eqref{eq:first_charge_conservation_for_single_component} and eq.\eqref{eq:second_charge_conservation_for_single_component} hold. Take eq.\eqref{eq:height_exchange_and_crossing3} or eq.\eqref{eq:height_exchange_and_crossing5}, with arbitrary choice of direction on the biangle. In eq.\eqref{eq:height_exchange_and_crossing3_values}, when the value of ${\rm Tr}^\omega_B$ is $q^{-1/3}$ or $q^{2/3}$, note $s(x_1)=s(y_1)$ and $s(x_2)=s(y_2)$, hence $n^\epsilon_{{\rm left},\varepsilon} = n^\epsilon_{{\rm right},\varepsilon}$ holds for all $\epsilon\in\{+,-\}$ and $\varepsilon\in\{1,2,3\}$, so eq.\eqref{eq:first_charge_conservation_for_single_component} and eq.\eqref{eq:second_charge_conservation_for_single_component} hold. For the remaining (nonzero) cases $(s(x_1),s(x_2),s(y_1),s(y_2)) \in\{(1,3,2,2),(2,2,3,1),(1,3,3,1)\}$ in eq.\eqref{eq:height_exchange_and_crossing3_values}, note for each $h\in\{{\rm left},{\rm right}\}$ that $n^+_{h,1}=n^-_{h,3}$, $n^+_{h,2}=n^-_{h,2}$ and $n^+_{h,3}=n^-_{h,1}$ hold, hence 
 $\mathcal{C}^{(1)}_h([W_{i,j}, J_0|_{\partial W_{i,j}}]) = \mathcal{C}^{(2)}_h([W_{i,j}, J_0|_{\partial W_{i,j}}])=0$, and therefore eq.\eqref{eq:first_charge_conservation_for_single_component} and eq.\eqref{eq:second_charge_conservation_for_single_component} hold. Finally, for the cases in eq.\eqref{eq:some_more_elementary1}--\eqref{eq:some_more_elementary3}, the values of ${\rm Tr}^\omega_B$ are matrix entries of the inverse matrices of the matrices for (BT2-4) which are upper triangular, as mentioned at the end of the previous subsection \S\ref{subsec:biangle_SL3_quantum_trace}; if one is just interested in when ${\rm Tr}^\omega_B$ is nonzero or not, one observes that the situation is exactly same as for (BT2-4), so eq.\eqref{eq:first_charge_conservation_for_single_component} and eq.\eqref{eq:second_charge_conservation_for_single_component} hold.

\vs

Now, summing eq.\eqref{eq:first_charge_conservation_for_single_component} and eq.\eqref{eq:second_charge_conservation_for_single_component} over $j=1,\ldots,l_i$, we get
\begin{align}
\nonumber
\mathcal{C}^{(1)}_{\rm left}(W_i,J_0|_{\partial W_i}) = \mathcal{C}^{(1)}_{\rm right}(W_i,J_0|_{\partial W_i}) \quad \mbox{and}\quad
\mathcal{C}^{(2)}_{\rm left}(W_i,J_0|_{\partial W_i}) = \mathcal{C}^{(2)}_{\rm right}(W_i,J_0|_{\partial W_i}),
\end{align}
For each $k=1,2$, observe from definition that $\mathcal{C}^{(k)}_{\rm right}(W_i,J_0|_{\partial W_i}) = \mathcal{C}^{(k)}_{\rm left}(W_{i+1},J_0|_{\partial W_{i+1}})$ holds for each $i=1,\ldots,n$, and $\mathcal{C}^{(k)}_{\rm left}(W_1,J_0|_{\partial W_1}) = \mathcal{C}^{(k)}_{\rm left}(W,s)$ and $\mathcal{C}^{(k)}_{\rm right}(W_{n+1},J_0|_{\partial W_{n+1}}) = \mathcal{C}^{(k)}_{\rm right}(W,s)$. So, by using above equalities repeatedly, one obtains $\mathcal{C}^{(k)}_{\rm left}(W,s) = \mathcal{C}^{(k)}_{\rm right}(W,s)$, as desired. \qed

\vs

{\it Proof of Prop.\ref{prop:balancedness_of_state-sum_trace}.} Let $\frak{S}$, $\Delta$, and $\wh{\Delta}$ be as in Def.\ref{def:good_position}, and $(W,s)$ be a stated ${\rm SL}_3$-web $(W,s)$ in a gool position with respect to $\wh{\Delta}$. Recall the state-sum formula for $\wh{\rm Tr}^\omega_\Delta(W,s)$ as in eq.\eqref{eq:state-sum_formula} of Def.\ref{def:state-sum_trace_for_gool_position}. Let $J$ be a $\wh{\Delta}$-juncture-state for $W$ restricting to $s$ such that the corresponding summand in eq.\eqref{eq:state-sum_formula} is nonzero; then the biangle factor ${\rm Tr}^\omega_B([W\cap (B\times{\bf I}),J_B])$ for each biangle $B$ of $\wh{\Delta}$ is nonzero. Pick any internal (i.e. non-boundary) arc $e$ of $\Delta$, and let $B$ be the corresponding biangle of $\wh{\Delta}$. Let $t,r$ be the ideal triangles of $\Delta$ having $e$ as a side, and let $\wh{t},\wh{r}$ be the corresponding triangles of $\wh{\Delta}$; say, $e\in \wh{\Delta}$ is a side of $\wh{t}$ and $e'\in\wh{\Delta}$ is a side of $\wh{r}$. Note that quiver $Q_\Delta$ has two nodes on the arc $e$, say $v_1$ and $v_2$, such that the direction $v_1 \to v_2$ matches the clockwise orientation on $\partial t$ and counterclockwise orientation on $\partial r$. We investigate the powers of $\wh{Z}_{t,v_1}$, $\wh{Z}_{t,v_2}$, $\wh{Z}_{r,v_1}$, $\wh{Z}_{r,v_2}$ in the (tensor) product of the triangle factors $\wh{\rm Tr}^\omega_t(W\cap (\wh{t}\times {\bf I}), J_t) \otimes \wh{\rm Tr}^\omega_{\redfix{r}}(W\cap (\wh{r}\times{\bf I}), J_r)$; these generators do not appear in the triangle factors for triangles other than $t,r$. 

\vs

For convenience, choose a direction of $B$ to make it a directed biangle $\vec{B}$, so that the left side $b_{\rm left}$ is a side of $\wh{t}$ and the right side $b_{\rm right}$ belongs to $\wh{r}$. For $h\in \{{\rm left},{\rm right}\}$, $\epsilon\in\{+,-\}$ and $\varepsilon\in\{1,2,3\}$, define the numbers $n^\epsilon_{h,\varepsilon}$ as in Lem.\ref{lem:charge_conservation}.
Investigating the triangle factor $\wh{\rm Tr}^\omega_t(W\cap (\wh{t}\times{\bf I}), J_t)$ as in eq.\eqref{eq:triangle_factor}, since $W$ is in a gool position, each factor as in eq.\eqref{eq:value_of_stated_component} falls into \redfix{(QT2-1)} or \redfix{(QT2-2)} of Thm.\ref{thm:SL3_quantum_trace_map}. Looking at \redfix{(QT2-1)} and \redfix{(QT2-2)}, $\wh{Z}_{t,v_1}$ or $\wh{Z}_{t,v_2}$ may appear in the entries of the matrices $\wh{\bf M}^{\rm in}_{t,e}$ and $\wh{\bf M}^{\rm out}_{t,e}$ which are all diagonal matrices, but not in the left or right turn matrices $\wh{\bf M}^{\rm left}(\wh{Z}_{t,v_t})$ and $\wh{\bf M}^{\rm right}(\wh{Z}_{t,v_t})$, or in the edge matrices $\wh{\bf M}^*_{t,*}$ for edges other than $e$. So, by investigating the edge matrices $\wh{\bf M}^*_{t,e}$, it follows that the triangle factor $\wh{\rm Tr}^\omega_t(W\cap (\wh{t}\times{\bf I}), J_t)$ equals
\begin{align*}
& \hspace*{-0mm} (\wh{Z}_{t,v_1} \wh{Z}_{t,v_2}^2)^{n^+_{{\rm left},1}} (\wh{Z}_{t,v_1} \wh{Z}_{t,v_2}^{-1})^{n^+_{{\rm left},2}}  (\wh{Z}_{t,v_1}^{-2} \wh{Z}_{t,v_2}^{-1})^{n^+_{{\rm left},3}} 
(\wh{Z}_{t,v_2} \wh{Z}_{t,v_1}^2)^{n^-_{{\rm left},1}} (\wh{Z}_{t,v_2} \wh{Z}_{t,v_1}^{-1})^{n^-_{{\rm left},2}}  (\wh{Z}_{t,v_2}^{-2} \wh{Z}_{t,v_1}^{-1})^{n^-_{{\rm left},3}}  \\
& = \wh{Z}_{t,v_1}^{n^+_{{\rm left},1} + n^+_{{\rm left},2} - 2 {\rm n}^+_{{\rm leftt},3} + 2n^-_{{\rm left},1} - n^-_{{\rm left},2} - n^-_{{\rm left},3}} \,
\wh{Z}_{t,v_2}^{ 2 n^+_{{\rm left},1} - n^+_{{\rm left},2} - n^+_{{\rm left},3} + n^-_{{\rm left},1} + n^-_{{\rm left},2} - 2n^-_{{\rm left},3} }
\end{align*}
times a Laurent polynomial in the generators of the triangle algebra $\mathcal{Z}^\omega_t$ not involving the nodes $v_1$ or $v_2$. Similarly, the triangle factor $\wh{\rm Tr}^\omega_r(W\cap (\wh{r}\times {\bf I}), J_r)$ equals
\begin{align*}
& \hspace*{-0mm} (\wh{Z}_{r,v_1} \wh{Z}_{r,v_2}^2)^{n^+_{{\rm right},1}} (\wh{Z}_{r,v_1} \wh{Z}_{r,v_2}^{-1})^{n^+_{{\rm right},2}}  (\wh{Z}_{r,v_1}^{-2} \wh{Z}_{r,v_2}^{-1})^{n^+_{{\rm right},3}} 
(\wh{Z}_{r,v_2} \wh{Z}_{r,v_1}^2)^{n^-_{{\rm right},1}} (\wh{Z}_{r,v_2} \wh{Z}_{r,v_1}^{-1})^{n^-_{{\rm right},2}}  (\wh{Z}_{r,v_2}^{-2} \wh{Z}_{r,v_1}^{-1})^{n^-_{{\rm right},3}}  \\
& = \wh{Z}_{r,v_1}^{n^+_{{\rm right},1} + n^+_{{\rm right},2} - 2 {\rm n}^+_{{\rm right},3} + 2n^-_{{\rm right},1} - n^-_{{\rm right},2} - n^-_{{\rm right},3}} \,
\wh{Z}_{r,v_2}^{ 2 n^+_{{\rm right},1} - n^+_{{\rm right},2} - n^+_{{\rm right},3} + n^-_{{\rm right},1} + n^-_{{\rm right},2} - 2n^-_{{\rm right},3} }
\end{align*}
times a Laurent polynomial in the generators of the triangle algebra $\mathcal{Z}^\omega_r$ not involving the nodes $v_1$ or $v_2$. Since ${\rm Tr}^\omega_B([W\cap (B\times {\bf I}),J_B]) \neq 0$, from Cor.\ref{cor:charge_conservation} we observe that the power of $\wh{Z}_{t,v_1}$ matches the power of $\wh{Z}_{r,v_1}$, and that the power of $\wh{Z}_{t,v_2}$ matches the power of $\wh{Z}_{r,v_2}$, hence establishing the balancedness as being asserted in the present Proposition \ref{prop:balancedness_of_state-sum_trace}. \qed

\vs

Crucial thing to show is the isotopy invariance of the state-sum trace formulated as follows, which we prove in the next subsection.
\begin{proposition}[isotopy invariance of state-sum trace for gool positions]
\label{prop:isotopy_invariance_of_state-sum_trace_for_gool_positions}
Let $\frak{S}$, $\Delta$, and $\wh{\Delta}$ be as in Def.\ref{def:good_position}. If $(W,s)$ and $(W',s')$ are isotopic stated ${\rm SL}_3$-webs in $\frak{S}\times {\bf I}$ in gool positions with respect to $\wh{\Delta}$, then $\wh{\rm Tr}^\omega_\Delta(W,s) = \wh{\rm Tr}^\omega_\Delta(W',s')$.
\end{proposition}

\vs

{\it Proof of the sought-for Thm.\ref{thm:SL3_quantum_trace_map}, assuming Prop.\ref{prop:isotopy_invariance_of_state-sum_trace_for_gool_positions}.} Let $\frak{S}$, $\Delta$, and $\wh{\Delta}$ be as in Def.\ref{def:good_position}. We will construct a map ${\rm Tr}^\omega_{\Delta;\frak{S}}$. For any stated ${\rm SL}_3$-web $(W,s)$ in $\frak{S}\times {\bf I}$, let $(W',s')$ be a stated ${\rm SL}_3$-web in a gool position with respect to $\wh{\Delta}$ and is isotopic to $(W,s)$ (it exists by Lem.\ref{lem:isotopic_to_good_position}). Define ${\rm Tr}^\omega_{\Delta;\frak{S}}([W,s])$ to be the value $\wh{\rm Tr}^\omega_\Delta(W',s')$, which lies in $\mathcal{Z}^\omega_\Delta$ according to Prop.\ref{prop:balancedness_of_state-sum_trace}:
$$
{\rm Tr}^\omega_{\Delta;\frak{S}}([W,s]) := \wh{\rm Tr}^\omega_\Delta(W',s').
$$

By Prop.\ref{prop:isotopy_invariance_of_state-sum_trace_for_gool_positions}, this value depends only on $[W,s]$ and does not depend on the choice of $(W',s')$. To show that this map ${\rm Tr}^\omega_{\Delta;\frak{S}} : \mathcal{S}^\omega_{\rm s}(\frak{S};\mathbb{Z})_{\rm red} \to \mathcal{Z}^\omega_\Delta$ is well-defined, one must show that the defining relations of the stated ${\rm SL}_3$-skein algebra $\mathcal{S}^\omega_{\rm s}(\frak{S};\mathbb{Z})_{\rm red}$ are satisfied. By pushing all the relations to biangles, one observes that this is the case. For example, take the ${\rm SL}_3$-skein relation (S8) of Fig.\ref{fig:A2-skein_relations_quantum}, so that we have $[W,s] = q^{-2/3} [W_1,s_1] + q^{1/3} [W_2,s_2]$ in $\mathcal{S}^\omega_{\rm s}(\frak{S};\mathbb{Z})_{\rm red}$, where the stated ${\rm SL}_3$-webs $(W,s)$, $(W_1,s_1)$, $(W_2,s_2)$ in $\frak{S}\times {\bf I}$ are identical except over a small disk as shown in the three figures appearing in (S8). By applying same isotopies to these three stated ${\rm SL}_3$-webs, one can push this disk to the interior of a biangle $B$ of $\wh{\Delta}$. Note that the $\wh{\Delta}$-juncture-states for these three stated ${\rm SL}_3$-webs are naturally in bijection. For each such $\wh{\Delta}$-juncture-state $J$, in eq.\eqref{eq:state-sum_formula} the only difference among the three is the biangle factor for $B$, where ${\rm Tr}^\omega_B([W\cap (B\times {\bf I}), J_B]) = q^{-2/3} {\rm Tr}^\omega_B([W_1\cap (B\times {\bf I}), J_B]) + q^{1/3} {\rm Tr}^\omega_B([W_2\cap (B\times {\bf I}), J_B])$ holds because ${\rm Tr}^\omega_B$ is a well-defined map on the stated ${\rm SL}_3$-skein algebra $\mathcal{S}^\omega_{\rm s}(B;\mathbb{Z})_{\rm red}$ (Prop.\ref{prop:biangle_SL3_quantum_trace}). Hence it follows that $\wh{\rm Tr}^\omega_\Delta(W,s) = q^{-2/3} \wh{\rm Tr}^\omega_\Delta(W_1,s_1) + q^{1/3} \wh{\rm Tr}^\omega_\Delta(W_2,s_2)$, as desired. The facts that ${\rm Tr}^\omega_{\Delta;\frak{S}}$ constructed this way is a $\mathbb{Z}[\omega^{\pm 1/2}]$-algebra homomorphism and that it satisfies the properties \redfix{(QT1)} and \redfix{(QT2-1)}--\redfix{(QT2-2)} of Thm.\ref{thm:SL3_quantum_trace_map} are built in from the very construction of the state-sum trace $\wh{\rm Tr}^\omega_\Delta$. The properties \redfix{(QT2-3)}--\redfix{(QT2-5)} would follow from the following:

\begin{proposition}[state-sum trace for a good position]
\label{prop:state-sum_trace_for_good_positions}
Let $\frak{S}$, $\Delta$, and $\wh{\Delta}$ be as in Def.\ref{def:good_position}. Let $(W,s)$ be a stated ${\rm SL}_3$-web in $\frak{S}\times {\bf I}$ in a good position with respect to $\wh{\Delta}$. Define the state-sum trace
$$
\til{\rm Tr}^\omega_\Delta(W,s) ~\in~ {\textstyle \bigotimes}_{t\in \mathcal{F}(\Delta)} \mathcal{Z}^\omega_t
$$
precisely as in eq.\eqref{eq:state-sum_formula} of Def.\ref{def:state-sum_trace_for_gool_position}, where the value $\til{\rm Tr}^\omega_t(W_{t,j},J_{j,t}) \in \mathcal{Z}^\omega_t$ (replacing eq.\eqref{eq:value_of_stated_component}) of each stated component $(W_{t,j},J_{t,j})$ for a triangle $t$ of $\Delta$ is now defined using Thm.\ref{thm:SL3_quantum_trace_map}\redfix{(QT2-1)}--\redfix{(QT2-5)}. 

\vs

If $(W',s')$ is a stated ${\rm SL}_3$-web in $\frak{S} \times {\bf I}$ in a gool position and $(W,s)$ is isotopic to $(W',s')$, then
$$
\til{\rm Tr}^\omega_\Delta(W,s) = \wh{\rm Tr}^\omega_\Delta(W',s').
$$
\end{proposition}

\subsection{Isotopy invariance of the state-sum formula}
\label{subsec:isotopy_invariance}

It remains to show Prop.\ref{prop:isotopy_invariance_of_state-sum_trace_for_gool_positions} and Prop.\ref{prop:state-sum_trace_for_good_positions} 
in order to complete our proof of Thm.\ref{thm:SL3_quantum_trace_map}, modulo Prop.\ref{prop:biangle_SL3_quantum_trace}. For both propositions, it helps to establish the following three statements first. Prop.\ref{prop:elementary_isotopy_invariance_3-way} which involves 3-valent vertices is especially important, for in \cite{CS} \cite{Douglas} \cite{Douglas21} the remaining Prop.\ref{prop:elementary_isotopy_invariance_elevation_preserving_arcs} and Prop.\ref{prop:elementary_isotopy_invariance_elevation_change} which do not involve 3-valent vertices had been considered already, but not Prop.\ref{prop:elementary_isotopy_invariance_3-way}.

\begin{figure}[htbp!]
\vspace{0mm}
\begin{center}
\raisebox{-0.5\height}{\scalebox{0.7}{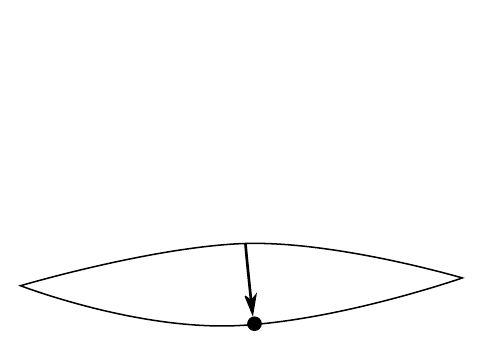}}
$\leftrightarrow$
\raisebox{-0.5\height}{\scalebox{0.7}{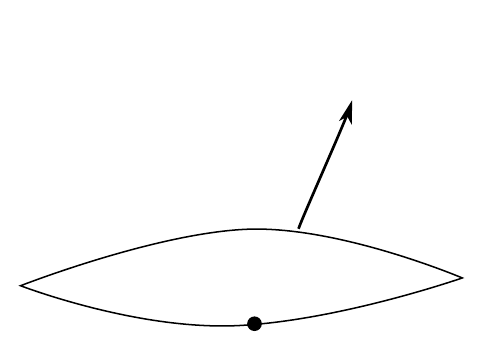}}
\end{center}
\vspace{-5mm}
\caption{Moving 3-valent vertex 1 (with $w' \succ w''$)}
\vspace{0mm}
\label{fig:isotopy_move1}
\end{figure}

\begin{figure}[htbp!]
\vspace{-2mm}
\begin{center}
\raisebox{-0.5\height}{\scalebox{0.7}{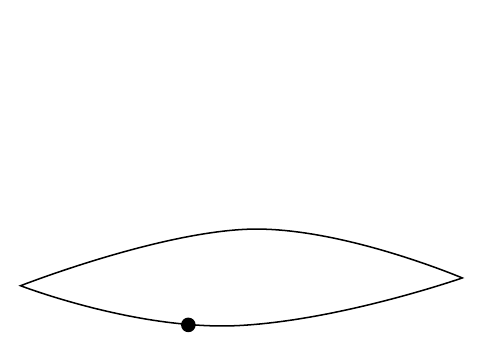}}
\hspace{-4mm} 
$\leftrightarrow$
\hspace{-3mm}
\raisebox{-0.5\height}{\scalebox{0.7}{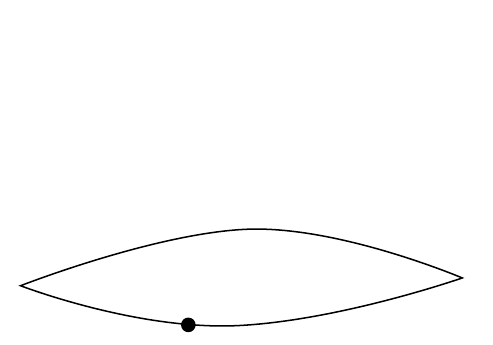}}
\qquad
\raisebox{-0.5\height}{\scalebox{0.7}{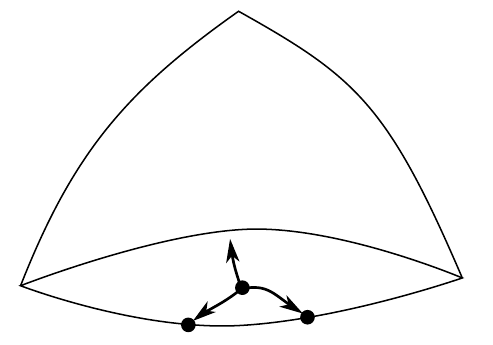}}
\hspace{-4mm} 
$\leftrightarrow$
\hspace{-3mm}
\raisebox{-0.5\height}{\scalebox{0.7}{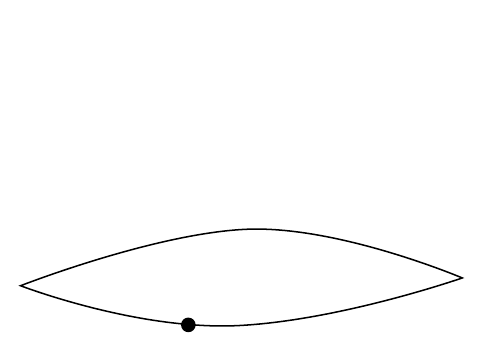}}
\end{center}
\vspace{-3mm}
\caption{Moving 3-valent vertex 2 (with $z_2\succ z_1$, $w_2\succ w_1$, or $z_2\prec z_1$, $w_2\prec w_1$)}
\vspace{0mm}
\label{fig:isotopy_move2}
\end{figure}

\begin{figure}[htbp!]
\vspace{0mm}
\begin{center}
\raisebox{-0.5\height}{\scalebox{0.7}{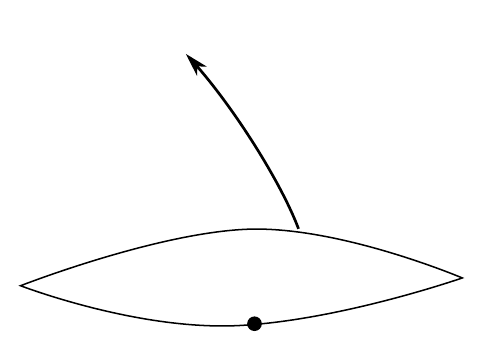}}
\hspace{-4mm} 
$\leftrightarrow$
\hspace{-3mm}
\raisebox{-0.5\height}{\scalebox{0.7}{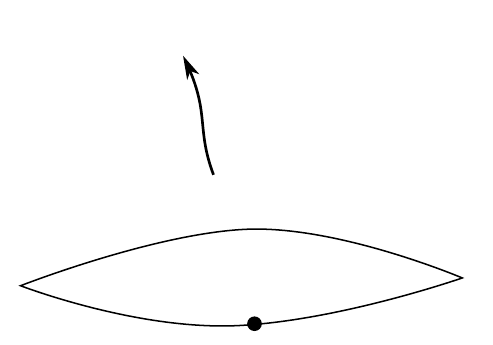}}
\qquad
\raisebox{-0.5\height}{\scalebox{0.7}{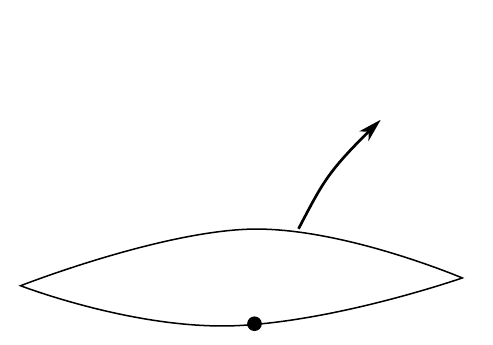}}
\hspace{-4mm} 
$\leftrightarrow$
\hspace{-3mm}
\raisebox{-0.5\height}{\scalebox{0.7}{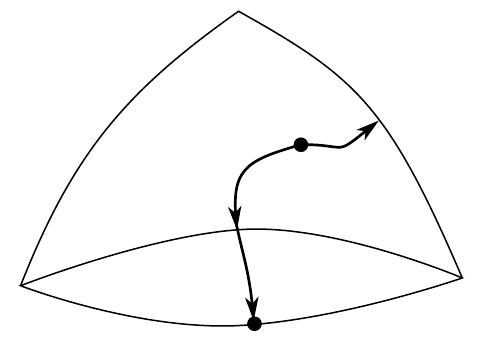}}
\end{center}
\vspace{-4mm}
\caption{Moving 3-valent vertex 3 (with $x_1\succ x_2$, $w_1\succ w_2$, $y_1\succ y_2$)}
\vspace{0mm}
\label{fig:isotopy_move3}
\end{figure}

\begin{figure}[htbp!]
\vspace{-2mm}
\begin{center}
\raisebox{-0.5\height}{\scalebox{0.7}{
\begingroup%
  \makeatletter%
  \providecommand\color[2][]{%
    \errmessage{(Inkscape) Color is used for the text in Inkscape, but the package 'color.sty' is not loaded}%
    \renewcommand\color[2][]{}%
  }%
  \providecommand\transparent[1]{%
    \errmessage{(Inkscape) Transparency is used (non-zero) for the text in Inkscape, but the package 'transparent.sty' is not loaded}%
    \renewcommand\transparent[1]{}%
  }%
  \providecommand\rotatebox[2]{#2}%
  \newcommand*\fsize{\dimexpr\f@size pt\relax}%
  \newcommand*\lineheight[1]{\fontsize{\fsize}{#1\fsize}\selectfont}%
  \ifx\svgwidth\undefined%
    \setlength{\unitlength}{141.73228346bp}%
    \ifx\svgscale\undefined%
      \relax%
    \else%
      \setlength{\unitlength}{\unitlength * \real{\svgscale}}%
    \fi%
  \else%
    \setlength{\unitlength}{\svgwidth}%
  \fi%
  \global\let\svgwidth\undefined%
  \global\let\svgscale\undefined%
  \makeatother%
  \begin{picture}(1,0.7)%
    \lineheight{1}%
    \setlength\tabcolsep{0pt}%
    \put(0,0){\includegraphics[width=\unitlength,page=1]{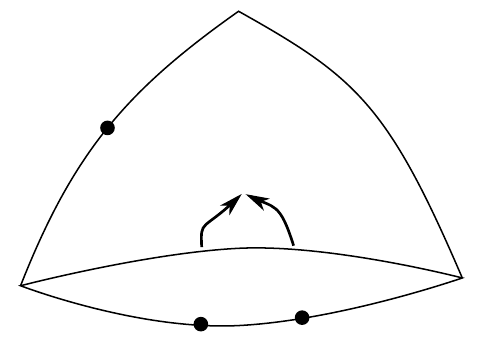}}%
    \put(0.47921241,0.2016868){\color[rgb]{0,0,0}\makebox(0,0)[lt]{\lineheight{1.25}\smash{\begin{tabular}[t]{l}$\prec$\end{tabular}}}}%
    \put(0,0){\includegraphics[width=\unitlength,page=2]{isotopy_invariance29_new.pdf}}%
  \end{picture}%
\endgroup%
}}
\hspace{-4mm} 
$\leftrightarrow$
\hspace{-3mm}
\raisebox{-0.5\height}{\scalebox{0.7}{
\begingroup%
  \makeatletter%
  \providecommand\color[2][]{%
    \errmessage{(Inkscape) Color is used for the text in Inkscape, but the package 'color.sty' is not loaded}%
    \renewcommand\color[2][]{}%
  }%
  \providecommand\transparent[1]{%
    \errmessage{(Inkscape) Transparency is used (non-zero) for the text in Inkscape, but the package 'transparent.sty' is not loaded}%
    \renewcommand\transparent[1]{}%
  }%
  \providecommand\rotatebox[2]{#2}%
  \newcommand*\fsize{\dimexpr\f@size pt\relax}%
  \newcommand*\lineheight[1]{\fontsize{\fsize}{#1\fsize}\selectfont}%
  \ifx\svgwidth\undefined%
    \setlength{\unitlength}{141.73228346bp}%
    \ifx\svgscale\undefined%
      \relax%
    \else%
      \setlength{\unitlength}{\unitlength * \real{\svgscale}}%
    \fi%
  \else%
    \setlength{\unitlength}{\svgwidth}%
  \fi%
  \global\let\svgwidth\undefined%
  \global\let\svgscale\undefined%
  \makeatother%
  \begin{picture}(1,0.7)%
    \lineheight{1}%
    \setlength\tabcolsep{0pt}%
    \put(0,0){\includegraphics[width=\unitlength,page=1]{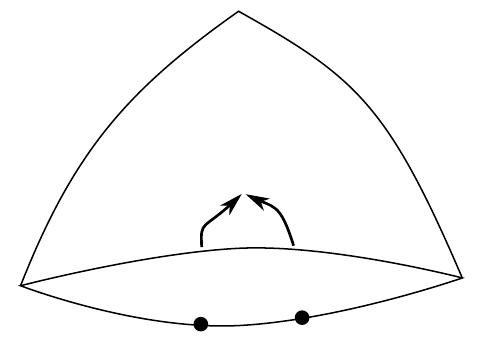}}%
    \put(0.47921241,0.2016868){\color[rgb]{0,0,0}\makebox(0,0)[lt]{\lineheight{1.25}\smash{\begin{tabular}[t]{l}$\succ$\end{tabular}}}}%
    \put(0,0){\includegraphics[width=\unitlength,page=2]{isotopy_invariance30_new.pdf}}%
  \end{picture}%
\endgroup%
}}
\qquad
\raisebox{-0.5\height}{\scalebox{0.7}{
\begingroup%
  \makeatletter%
  \providecommand\color[2][]{%
    \errmessage{(Inkscape) Color is used for the text in Inkscape, but the package 'color.sty' is not loaded}%
    \renewcommand\color[2][]{}%
  }%
  \providecommand\transparent[1]{%
    \errmessage{(Inkscape) Transparency is used (non-zero) for the text in Inkscape, but the package 'transparent.sty' is not loaded}%
    \renewcommand\transparent[1]{}%
  }%
  \providecommand\rotatebox[2]{#2}%
  \newcommand*\fsize{\dimexpr\f@size pt\relax}%
  \newcommand*\lineheight[1]{\fontsize{\fsize}{#1\fsize}\selectfont}%
  \ifx\svgwidth\undefined%
    \setlength{\unitlength}{141.73228346bp}%
    \ifx\svgscale\undefined%
      \relax%
    \else%
      \setlength{\unitlength}{\unitlength * \real{\svgscale}}%
    \fi%
  \else%
    \setlength{\unitlength}{\svgwidth}%
  \fi%
  \global\let\svgwidth\undefined%
  \global\let\svgscale\undefined%
  \makeatother%
  \begin{picture}(1,0.7)%
    \lineheight{1}%
    \setlength\tabcolsep{0pt}%
    \put(0,0){\includegraphics[width=\unitlength,page=1]{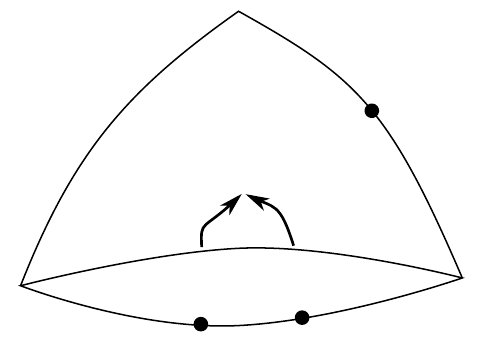}}%
    \put(0.47921241,0.2016868){\color[rgb]{0,0,0}\makebox(0,0)[lt]{\lineheight{1.25}\smash{\begin{tabular}[t]{l}$\prec$\end{tabular}}}}%
    \put(0,0){\includegraphics[width=\unitlength,page=2]{isotopy_invariance31_new.pdf}}%
  \end{picture}%
\endgroup%
}}
\hspace{-4mm} 
$\leftrightarrow$
\hspace{-3mm}
\raisebox{-0.5\height}{\scalebox{0.7}{
\begingroup%
  \makeatletter%
  \providecommand\color[2][]{%
    \errmessage{(Inkscape) Color is used for the text in Inkscape, but the package 'color.sty' is not loaded}%
    \renewcommand\color[2][]{}%
  }%
  \providecommand\transparent[1]{%
    \errmessage{(Inkscape) Transparency is used (non-zero) for the text in Inkscape, but the package 'transparent.sty' is not loaded}%
    \renewcommand\transparent[1]{}%
  }%
  \providecommand\rotatebox[2]{#2}%
  \newcommand*\fsize{\dimexpr\f@size pt\relax}%
  \newcommand*\lineheight[1]{\fontsize{\fsize}{#1\fsize}\selectfont}%
  \ifx\svgwidth\undefined%
    \setlength{\unitlength}{141.73228346bp}%
    \ifx\svgscale\undefined%
      \relax%
    \else%
      \setlength{\unitlength}{\unitlength * \real{\svgscale}}%
    \fi%
  \else%
    \setlength{\unitlength}{\svgwidth}%
  \fi%
  \global\let\svgwidth\undefined%
  \global\let\svgscale\undefined%
  \makeatother%
  \begin{picture}(1,0.7)%
    \lineheight{1}%
    \setlength\tabcolsep{0pt}%
    \put(0,0){\includegraphics[width=\unitlength,page=1]{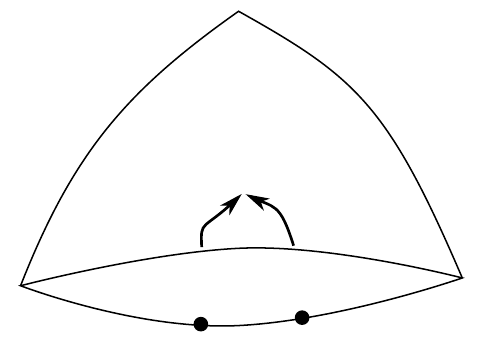}}%
    \put(0.47921241,0.2016868){\color[rgb]{0,0,0}\makebox(0,0)[lt]{\lineheight{1.25}\smash{\begin{tabular}[t]{l}$\succ$\end{tabular}}}}%
    \put(0,0){\includegraphics[width=\unitlength,page=2]{isotopy_invariance32_new.pdf}}%
  \end{picture}%
\endgroup%
}}
\end{center}
\vspace{-4mm}
\caption{Elevation exchange 1}
\vspace{0mm}
\label{fig:isotopy_move8}
\end{figure}

\begin{proposition}[isotopy invariance under the elementary moves involving 3-valent vertices]
\label{prop:elementary_isotopy_invariance_3-way}
Let $t$ be a triangle, viewed as a generalized marked surface. Let $\wh{\Delta}$ be the collection of four arcs in $t$, three of them being the boundary arcs of $t$, and the remaining one an arc connecting two marked points of $t$ whose interior lies in the interior of $t$; so $\wh{\Delta}$ divides $t$ into one triangle $\wh{t}$ and one biangle $B$. Consider the state-sum trace for stated ${\rm SL}_3$-webs in $t \times {\bf I}$ in a good position with respect to $\wh{\Delta}$ as defined in Prop.\ref{prop:state-sum_trace_for_good_positions}, which we denote by $\til{\rm Tr}^\omega_{\wh{\Delta}}$. Then $\til{\rm Tr}^\omega_{\wh{\Delta}}$ satisfies the following elementary isotopy invariance: if $(W,s)$ and $(W',s')$ are stated ${\rm SL}_3$-webs in $t\times {\bf I}$ in good positions with respect to $\wh{\Delta}$ and are related to each other by one of the moves in Fig.\ref{fig:isotopy_move1}, \ref{fig:isotopy_move2}, \ref{fig:isotopy_move3}, \ref{fig:isotopy_move8}, possibly with different possible orientations on the components, they have same values under $\til{\rm Tr}^\omega_{\wh{\Delta}}$.
\end{proposition}

\begin{figure}[htbp!]
\vspace{-2mm}
\begin{center}
\raisebox{-0.5\height}{\scalebox{0.7}{
\begingroup%
  \makeatletter%
  \providecommand\color[2][]{%
    \errmessage{(Inkscape) Color is used for the text in Inkscape, but the package 'color.sty' is not loaded}%
    \renewcommand\color[2][]{}%
  }%
  \providecommand\transparent[1]{%
    \errmessage{(Inkscape) Transparency is used (non-zero) for the text in Inkscape, but the package 'transparent.sty' is not loaded}%
    \renewcommand\transparent[1]{}%
  }%
  \providecommand\rotatebox[2]{#2}%
  \newcommand*\fsize{\dimexpr\f@size pt\relax}%
  \newcommand*\lineheight[1]{\fontsize{\fsize}{#1\fsize}\selectfont}%
  \ifx\svgwidth\undefined%
    \setlength{\unitlength}{141.73228346bp}%
    \ifx\svgscale\undefined%
      \relax%
    \else%
      \setlength{\unitlength}{\unitlength * \real{\svgscale}}%
    \fi%
  \else%
    \setlength{\unitlength}{\svgwidth}%
  \fi%
  \global\let\svgwidth\undefined%
  \global\let\svgscale\undefined%
  \makeatother%
  \begin{picture}(1,0.7)%
    \lineheight{1}%
    \setlength\tabcolsep{0pt}%
    \put(0,0){\includegraphics[width=\unitlength,page=1]{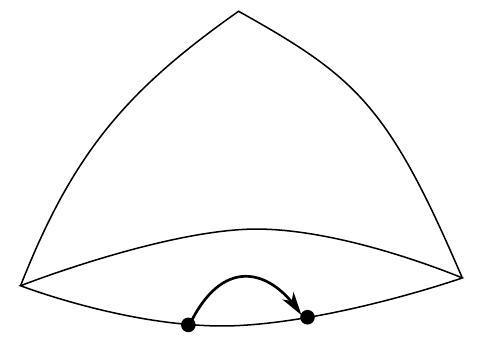}}%
    \put(0.46127484,0.50198418){\color[rgb]{0,0,0}\makebox(0,0)[lt]{\lineheight{1.25}\smash{\begin{tabular}[t]{l}$\wh{t}$\end{tabular}}}}%
    \put(0.08276521,0.36704372){\color[rgb]{0,0,0}\makebox(0,0)[lt]{\lineheight{1.25}\smash{\begin{tabular}[t]{l}$e_1$\end{tabular}}}}%
    \put(0.83239956,0.39273481){\color[rgb]{0,0,0}\makebox(0,0)[lt]{\lineheight{1.25}\smash{\begin{tabular}[t]{l}$e_2$\end{tabular}}}}%
    \put(0.21447941,0.21696364){\color[rgb]{0,0,0}\makebox(0,0)[lt]{\lineheight{1.25}\smash{\begin{tabular}[t]{l}$e_3$\end{tabular}}}}%
    \put(0.66835708,0.14990012){\color[rgb]{0,0,0}\makebox(0,0)[lt]{\lineheight{1.25}\smash{\begin{tabular}[t]{l}$B$\end{tabular}}}}%
    \put(0.21155101,0.09222444){\color[rgb]{0,0,0}\makebox(0,0)[lt]{\lineheight{1.25}\smash{\begin{tabular}[t]{l}$e_3'$\end{tabular}}}}%
    \put(0.3379864,0.08748284){\color[rgb]{0,0,0}\makebox(0,0)[lt]{\lineheight{1.25}\smash{\begin{tabular}[t]{l}$z_2$\end{tabular}}}}%
    \put(0.63594239,0.08617147){\color[rgb]{0,0,0}\makebox(0,0)[lt]{\lineheight{1.25}\smash{\begin{tabular}[t]{l}$z_1$\end{tabular}}}}%
  \end{picture}%
\endgroup%
}}
\hspace{-4mm} 
$\leftrightarrow$
\hspace{-3mm}
\raisebox{-0.5\height}{\scalebox{0.7}{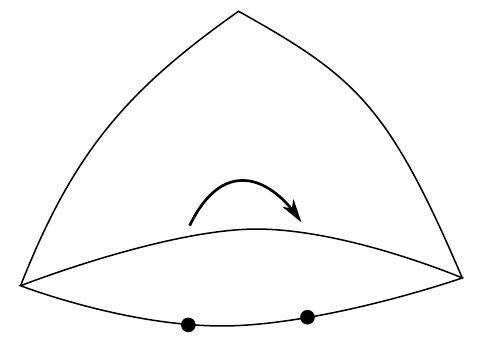}}
\qquad
\raisebox{-0.5\height}{\scalebox{0.7}{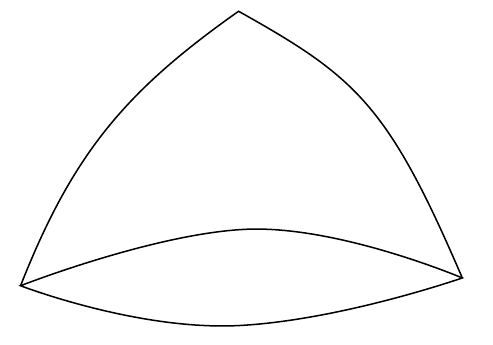}}
\hspace{-4mm} 
$\leftrightarrow$
\hspace{-3mm}
\raisebox{-0.5\height}{\scalebox{0.7}{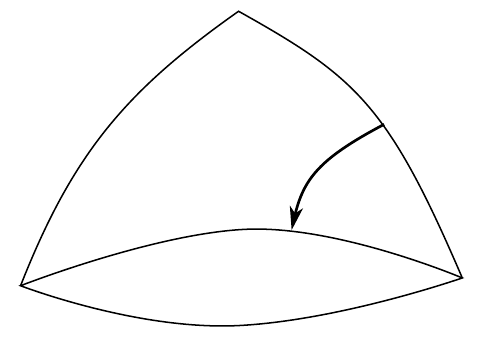}}
\end{center}
\vspace{-4mm}
\caption{Moving cup/cap 1}
\vspace{0mm}
\label{fig:isotopy_move4}
\end{figure}

\begin{figure}[htbp!]
\vspace{-2mm}
\begin{center}
\raisebox{-0.5\height}{\scalebox{0.7}{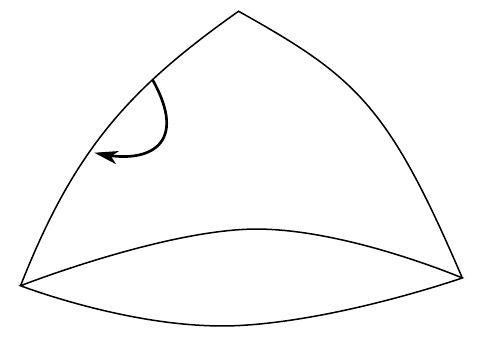}}
\hspace{-4mm} 
$\leftrightarrow$
\hspace{-3mm}
\raisebox{-0.5\height}{\scalebox{0.7}{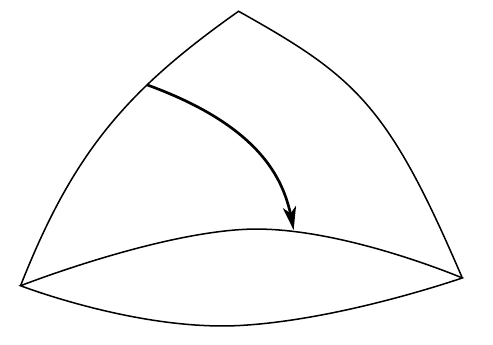}}
\qquad
\raisebox{-0.5\height}{\scalebox{0.7}{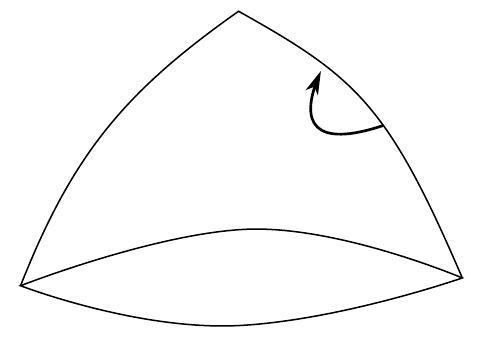}}
\hspace{-4mm} 
$\leftrightarrow$
\hspace{-3mm}
\raisebox{-0.5\height}{\scalebox{0.7}{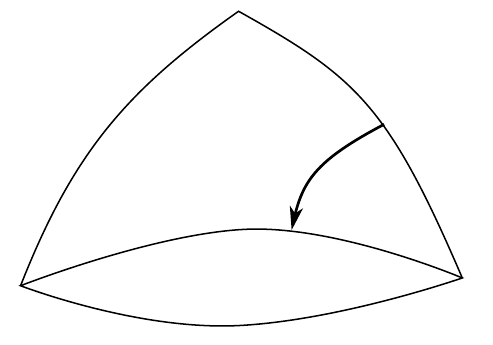}}
\end{center}
\vspace{-4mm}
\caption{Moving cup/cap 2 $\mbox{(with $x_1\succ x_2, w_1\succ w_2, y_1\succ y_2$, or $x_2\succ x_1, w_2\succ w_1, y_2\succ y_1$)}$}
\vspace{-2mm}
\label{fig:isotopy_move5}
\end{figure}

\begin{proposition}[isotopy invariance under the elevation preserving elementary moves not involving 3-valent vertices; \redfix{essentially in \cite{CS}, and} partially \redfix{in} \cite{Douglas} \cite{Douglas21}]
\label{prop:elementary_isotopy_invariance_elevation_preserving_arcs}
Analogous statement as Prop.\ref{prop:elementary_isotopy_invariance_3-way} holds for the moves in Fig.\ref{fig:isotopy_move4}, \ref{fig:isotopy_move5}, possibly with different possible orientations on the components.
\end{proposition}

\begin{figure}[htbp!]
\vspace{0mm}
\begin{center}
\raisebox{-0.5\height}{\scalebox{0.7}{
\begingroup%
  \makeatletter%
  \providecommand\color[2][]{%
    \errmessage{(Inkscape) Color is used for the text in Inkscape, but the package 'color.sty' is not loaded}%
    \renewcommand\color[2][]{}%
  }%
  \providecommand\transparent[1]{%
    \errmessage{(Inkscape) Transparency is used (non-zero) for the text in Inkscape, but the package 'transparent.sty' is not loaded}%
    \renewcommand\transparent[1]{}%
  }%
  \providecommand\rotatebox[2]{#2}%
  \newcommand*\fsize{\dimexpr\f@size pt\relax}%
  \newcommand*\lineheight[1]{\fontsize{\fsize}{#1\fsize}\selectfont}%
  \ifx\svgwidth\undefined%
    \setlength{\unitlength}{141.73228346bp}%
    \ifx\svgscale\undefined%
      \relax%
    \else%
      \setlength{\unitlength}{\unitlength * \real{\svgscale}}%
    \fi%
  \else%
    \setlength{\unitlength}{\svgwidth}%
  \fi%
  \global\let\svgwidth\undefined%
  \global\let\svgscale\undefined%
  \makeatother%
  \begin{picture}(1,0.7)%
    \lineheight{1}%
    \setlength\tabcolsep{0pt}%
    \put(0,0){\includegraphics[width=\unitlength,page=1]{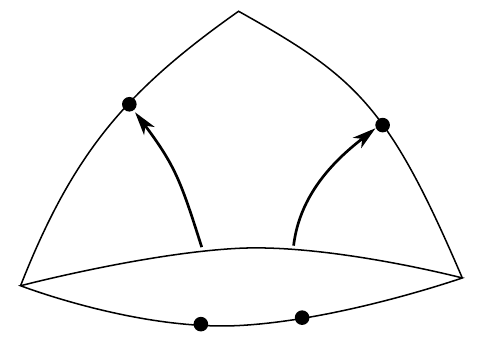}}%
    \put(0.47921241,0.2016868){\color[rgb]{0,0,0}\makebox(0,0)[lt]{\lineheight{1.25}\smash{\begin{tabular}[t]{l}$\prec$\end{tabular}}}}%
    \put(0,0){\includegraphics[width=\unitlength,page=2]{isotopy_invariance24.pdf}}%
  \end{picture}%
\endgroup%
}}
\hspace{-4mm} 
$\leftrightarrow$
\hspace{-3mm}
\raisebox{-0.5\height}{\scalebox{0.7}{
\begingroup%
  \makeatletter%
  \providecommand\color[2][]{%
    \errmessage{(Inkscape) Color is used for the text in Inkscape, but the package 'color.sty' is not loaded}%
    \renewcommand\color[2][]{}%
  }%
  \providecommand\transparent[1]{%
    \errmessage{(Inkscape) Transparency is used (non-zero) for the text in Inkscape, but the package 'transparent.sty' is not loaded}%
    \renewcommand\transparent[1]{}%
  }%
  \providecommand\rotatebox[2]{#2}%
  \newcommand*\fsize{\dimexpr\f@size pt\relax}%
  \newcommand*\lineheight[1]{\fontsize{\fsize}{#1\fsize}\selectfont}%
  \ifx\svgwidth\undefined%
    \setlength{\unitlength}{141.73228346bp}%
    \ifx\svgscale\undefined%
      \relax%
    \else%
      \setlength{\unitlength}{\unitlength * \real{\svgscale}}%
    \fi%
  \else%
    \setlength{\unitlength}{\svgwidth}%
  \fi%
  \global\let\svgwidth\undefined%
  \global\let\svgscale\undefined%
  \makeatother%
  \begin{picture}(1,0.7)%
    \lineheight{1}%
    \setlength\tabcolsep{0pt}%
    \put(0,0){\includegraphics[width=\unitlength,page=1]{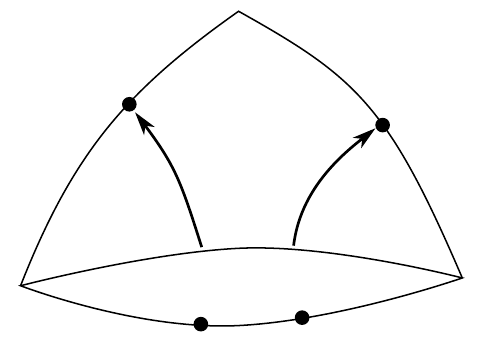}}%
    \put(0.47921241,0.2016868){\color[rgb]{0,0,0}\makebox(0,0)[lt]{\lineheight{1.25}\smash{\begin{tabular}[t]{l}$\succ$\end{tabular}}}}%
    \put(0,0){\includegraphics[width=\unitlength,page=2]{isotopy_invariance25.pdf}}%
  \end{picture}%
\endgroup%
}}
\qquad
\raisebox{-0.5\height}{\scalebox{0.7}{}}
\hspace{-4mm} 
$\leftrightarrow$
\hspace{-3mm}
\raisebox{-0.5\height}{\scalebox{0.7}{
\begingroup%
  \makeatletter%
  \providecommand\color[2][]{%
    \errmessage{(Inkscape) Color is used for the text in Inkscape, but the package 'color.sty' is not loaded}%
    \renewcommand\color[2][]{}%
  }%
  \providecommand\transparent[1]{%
    \errmessage{(Inkscape) Transparency is used (non-zero) for the text in Inkscape, but the package 'transparent.sty' is not loaded}%
    \renewcommand\transparent[1]{}%
  }%
  \providecommand\rotatebox[2]{#2}%
  \newcommand*\fsize{\dimexpr\f@size pt\relax}%
  \newcommand*\lineheight[1]{\fontsize{\fsize}{#1\fsize}\selectfont}%
  \ifx\svgwidth\undefined%
    \setlength{\unitlength}{141.73228346bp}%
    \ifx\svgscale\undefined%
      \relax%
    \else%
      \setlength{\unitlength}{\unitlength * \real{\svgscale}}%
    \fi%
  \else%
    \setlength{\unitlength}{\svgwidth}%
  \fi%
  \global\let\svgwidth\undefined%
  \global\let\svgscale\undefined%
  \makeatother%
  \begin{picture}(1,0.7)%
    \lineheight{1}%
    \setlength\tabcolsep{0pt}%
    \put(0,0){\includegraphics[width=\unitlength,page=1]{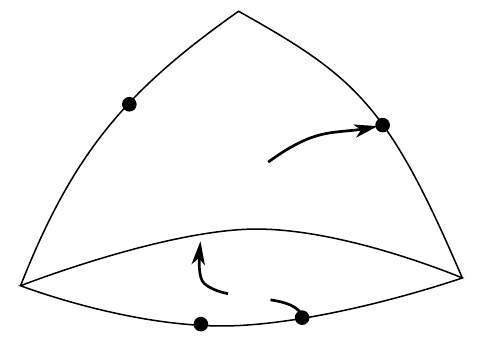}}%
    \put(0.47921241,0.2016868){\color[rgb]{0,0,0}\makebox(0,0)[lt]{\lineheight{1.25}\smash{\begin{tabular}[t]{l}$\prec$\end{tabular}}}}%
    \put(0,0){\includegraphics[width=\unitlength,page=2]{isotopy_invariance26.pdf}}%
  \end{picture}%
\endgroup%
}}
\end{center}
\vspace{-4mm}
\caption{Elevation exchange 2 / crossing}
\vspace{-2mm}
\label{fig:isotopy_move6}
\end{figure}

\begin{figure}[htbp!]
\vspace{0mm}
\begin{center}
\raisebox{-0.5\height}{\scalebox{0.7}{
\begingroup%
  \makeatletter%
  \providecommand\color[2][]{%
    \errmessage{(Inkscape) Color is used for the text in Inkscape, but the package 'color.sty' is not loaded}%
    \renewcommand\color[2][]{}%
  }%
  \providecommand\transparent[1]{%
    \errmessage{(Inkscape) Transparency is used (non-zero) for the text in Inkscape, but the package 'transparent.sty' is not loaded}%
    \renewcommand\transparent[1]{}%
  }%
  \providecommand\rotatebox[2]{#2}%
  \newcommand*\fsize{\dimexpr\f@size pt\relax}%
  \newcommand*\lineheight[1]{\fontsize{\fsize}{#1\fsize}\selectfont}%
  \ifx\svgwidth\undefined%
    \setlength{\unitlength}{141.73228346bp}%
    \ifx\svgscale\undefined%
      \relax%
    \else%
      \setlength{\unitlength}{\unitlength * \real{\svgscale}}%
    \fi%
  \else%
    \setlength{\unitlength}{\svgwidth}%
  \fi%
  \global\let\svgwidth\undefined%
  \global\let\svgscale\undefined%
  \makeatother%
  \begin{picture}(1,0.7)%
    \lineheight{1}%
    \setlength\tabcolsep{0pt}%
    \put(0,0){\includegraphics[width=\unitlength,page=1]{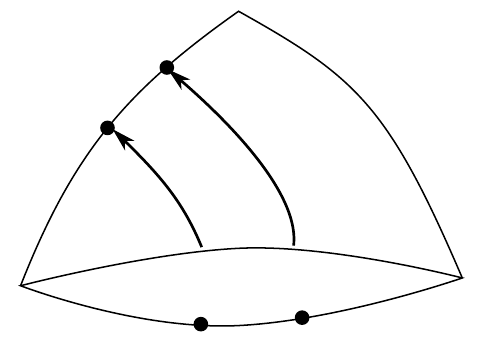}}%
    \put(0.47921241,0.2016868){\color[rgb]{0,0,0}\makebox(0,0)[lt]{\lineheight{1.25}\smash{\begin{tabular}[t]{l}$\prec$\end{tabular}}}}%
    \put(0,0){\includegraphics[width=\unitlength,page=2]{isotopy_invariance27.pdf}}%
    \put(0.39080946,0.36856361){\color[rgb]{0,0,0}\rotatebox{52.835116}{\makebox(0,0)[lt]{\lineheight{1.25}\smash{\begin{tabular}[t]{l}$\prec$\end{tabular}}}}}%
    \put(0,0){\includegraphics[width=\unitlength,page=3]{isotopy_invariance27.pdf}}%
  \end{picture}%
\endgroup%
}}
\hspace{-4mm} 
$\leftrightarrow$
\hspace{-3mm}
\raisebox{-0.5\height}{\scalebox{0.7}{
\begingroup%
  \makeatletter%
  \providecommand\color[2][]{%
    \errmessage{(Inkscape) Color is used for the text in Inkscape, but the package 'color.sty' is not loaded}%
    \renewcommand\color[2][]{}%
  }%
  \providecommand\transparent[1]{%
    \errmessage{(Inkscape) Transparency is used (non-zero) for the text in Inkscape, but the package 'transparent.sty' is not loaded}%
    \renewcommand\transparent[1]{}%
  }%
  \providecommand\rotatebox[2]{#2}%
  \newcommand*\fsize{\dimexpr\f@size pt\relax}%
  \newcommand*\lineheight[1]{\fontsize{\fsize}{#1\fsize}\selectfont}%
  \ifx\svgwidth\undefined%
    \setlength{\unitlength}{141.73228346bp}%
    \ifx\svgscale\undefined%
      \relax%
    \else%
      \setlength{\unitlength}{\unitlength * \real{\svgscale}}%
    \fi%
  \else%
    \setlength{\unitlength}{\svgwidth}%
  \fi%
  \global\let\svgwidth\undefined%
  \global\let\svgscale\undefined%
  \makeatother%
  \begin{picture}(1,0.7)%
    \lineheight{1}%
    \setlength\tabcolsep{0pt}%
    \put(0,0){\includegraphics[width=\unitlength,page=1]{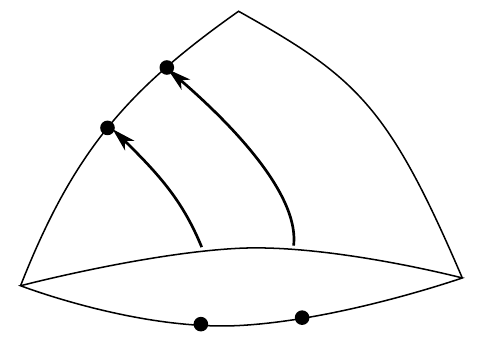}}%
    \put(0.47921241,0.2016868){\color[rgb]{0,0,0}\makebox(0,0)[lt]{\lineheight{1.25}\smash{\begin{tabular}[t]{l}$\succ$\end{tabular}}}}%
    \put(0,0){\includegraphics[width=\unitlength,page=2]{isotopy_invariance28.pdf}}%
    \put(0.39080946,0.36856361){\color[rgb]{0,0,0}\rotatebox{52.835116}{\makebox(0,0)[lt]{\lineheight{1.25}\smash{\begin{tabular}[t]{l}$\succ$\end{tabular}}}}}%
    \put(0,0){\includegraphics[width=\unitlength,page=3]{isotopy_invariance28.pdf}}%
  \end{picture}%
\endgroup%
}}
\end{center}
\vspace{-4mm}
\caption{Elevation exchange 3}
\vspace{-2mm}
\label{fig:isotopy_move7}
\end{figure}

\begin{proposition}[isotopy invariance under the elevation changing elementary moves not involving 3-valent vertices; essentially in \cite{CS}, and partially in \cite{Douglas} \cite{Douglas21}]
\label{prop:elementary_isotopy_invariance_elevation_change}
Analogous statement as Prop.\ref{prop:elementary_isotopy_invariance_3-way}, for a triangle $t$ with triangulation $\Delta$ and a split ideal triangulation $\wh{\Delta}$, holds for the moves in Fig.\ref{fig:isotopy_move6}, \ref{fig:isotopy_move7}, possibly with different possible orientations on the components.
\end{proposition}

\vs

{\it Proof of Prop.\ref{prop:elementary_isotopy_invariance_3-way}.} Label the arcs of $\wh{\Delta}$ by $e_1,e_2,e_3,e_3'$, and the names of endpoints and junctures as $x,y,z,w,w',w''$ as in the pictures. Denoting by $\Delta$ the unique triangulation of $t$, label the seven nodes of the quiver $Q_\Delta$ as in Thm.\ref{thm:SL3_quantum_trace_map}\redfix{(QT2)}. First, let $W$ and $W'$ be the ${\rm SL}_3$-webs in $t \times {\bf I}$ as in the left and the right of one of Fig.\ref{fig:isotopy_move1}, \ref{fig:isotopy_move2}, \ref{fig:isotopy_move3}. For states $s$ and $s'$, we consider
\begin{align}
\label{eq:til_Tr_W_s}
\til{\rm Tr}^\omega_{\wh{\Delta}}(W,s) & = \underset{s_1,s_2}{\textstyle \sum} {\rm Tr}^\omega_B([W\cap (B \times {\bf I}), s_1]) \, \til{\rm Tr}^\omega_t(W\cap (\wh{t} \times {\bf I}),s_2), \\
\label{eq:til_Tr_W_s_prime}
\til{\rm Tr}^\omega_{\wh{\Delta}}(W',s') & = \underset{s_1',s_2'}{\textstyle \sum} {\rm Tr}^\omega_B([W'\cap (B \times {\bf I}), s_1']) \, \til{\rm Tr}^\omega_t(W'\cap (\wh{t} \times {\bf I}),s_2'),
\end{align}
where the first sum is over all states $s_1$ and $s_2$ of $W\cap (B \times {\bf I})$ and $W\cap (\wh{t} \times {\bf I})$ that are compatible with $s$, and similarly for the second sum.

\vs

Consider the case of Fig.\ref{fig:isotopy_move1}. Let $s_1,s_2$ be states of $W\cap (B \times {\bf I})$ and $W\cap (\wh{t} \times {\bf I})$ compatible with $s$ and whose corresponding summand in eq.\eqref{eq:til_Tr_W_s} is nonzero. Note $W\cap (B \times {\bf I})$ has just one component, which is an edge connecting the thickenings of two sides of $B$, so by (BT2-1) of Prop.\ref{prop:biangle_SL3_quantum_trace}, ${\rm Tr}^\omega_B([W\cap (B \times {\bf I}),s_1]) \neq0$ implies $s_1(w)=s_1(z)$, in which case ${\rm Tr}^\omega_B([W\cap (B \times {\bf I}),s_1])=1$. By compatibility, $s_2(w)=s_1(w)$ and $s_1(z)=s(z)$, $s_2(x)=s(x)$, $s_2(y)=s(y)$. So there is only one pair of $s_1,s_2$ contributing to the sum, and hence
$$
\til{\rm Tr}^\omega_{\wh{\Delta}}(W,s) = \til{\rm Tr}^\omega_t(W\cap (\wh{t} \times {\bf I}),s_2) = {\rm Tr}^\omega_{\Delta;t}([W,s]).
$$
Denoting the values of compatible states $s_1'$ and $s_2'$ in eq.\eqref{eq:til_Tr_W_s_prime} at the internal junctures $w'$ and $w''$ as $\varepsilon_4$ and $\varepsilon_5$ respectively, we have
\begin{align*}
\hspace*{-8mm} \til{\rm Tr}^\omega_{\wh{\Delta}}(W',s') 
& = {\textstyle \sum}_{\varepsilon_4,\varepsilon_5} {\rm Tr}^\omega_B([W'\cap (B \times {\bf I}),(\varepsilon_3,\varepsilon_4,\varepsilon_5)])
{\rm Tr}^\omega_t([W'\cap (\wh{t} \times {\bf I}), (\varepsilon_1,\varepsilon_2,\varepsilon_4,\varepsilon_5)])
\end{align*}
where the sum is over all $\varepsilon_4,\varepsilon_5 \in \{1,2,3\}$, and $(\varepsilon_3,\varepsilon_4,\varepsilon_5)$ denotes the state $s'_1$ assigning these values to $z,w'\redfix{,}w''$, while $(\varepsilon_1,\varepsilon_2,\varepsilon_4,\varepsilon_5)$ denotes the state $s_2'$ assigning these values to $x,y,w',w''$.   Note $w' \succ w''$. By Prop.\ref{prop:biangle_SL3_quantum_trace}(BT2-3) we have
\begin{align}
\label{eq:isotopy_invariance_proof_biangle_factor1}
{\rm Tr}^\omega_B([W'\cap (B \times {\bf I}),(\varepsilon_3,\varepsilon_4,\varepsilon_5)])
= (\wh{\bf F}^{\rm out}_{+,\varepsilon_3})_{\varepsilon_4,\varepsilon_5}
\end{align}
Since $W'\cap (\wh{t} \times {\bf I})$ is the product of a left turn corner arc from $w'$ (over $e_3$) to $x$ (over $e_1$) and a right turn corner arc from $w''$ (over $e_3$) to $y$ (over $e_2$), with the product taken in this order, from \redfix{(QT2-1)}--\redfix{(QT2-2)} of Thm.\ref{thm:SL3_quantum_trace_map} we have
\begin{align}
\label{eq:isotopy_invariance_proof_biangle_factor2}
{\rm Tr}^\omega_t([W'\cap (\wh{t} \times {\bf I}), (\varepsilon_1,\varepsilon_2,\varepsilon_4,\varepsilon_5)])
= (\wh{\bf M}^{\rm in}_{t,3} \wh{\bf M}^{\rm left}(\wh{Z}_{v_t}) \wh{\bf M}^{\rm out}_{t,1})_{\varepsilon_4,\varepsilon_1}
(\wh{\bf M}^{\rm in}_{t,3} \wh{\bf M}^{\rm right}(\wh{Z}_{v_t}) \wh{\bf M}^{\rm out}_{t,2})_{\varepsilon_5,\varepsilon_2}.
\end{align}
By Cor.\ref{cor:transpose_of_quantum_turn_matrices} we have $(\wh{\bf M}^{\rm in}_{t,3} \wh{\bf M}^{\rm left}(\wh{Z}_{v_t}) \wh{\bf M}^{\rm out}_{t,1})_{\varepsilon_4,\varepsilon_1} = (\wh{\bf M}^{\rm out}_{t,1} \wh{\bf M}^{\rm left}_{\rm tran}(\wh{Z}_{v_t}) \wh{\bf M}^{\rm in}_{t,3})_{\varepsilon_1,\varepsilon_4}$, hence
\begin{align}
\nonumber
\til{\rm Tr}^\omega_{\wh{\Delta}}(W',s') &= {\textstyle \sum}_{\varepsilon_4,\varepsilon_5} {\rm Tr}^\omega_B([W'\cap (B \times {\bf I}),(\varepsilon_3,\varepsilon_4,\varepsilon_5)])
{\rm Tr}^\omega_t([W'\cap (\wh{t} \times {\bf I}), (\varepsilon_1,\varepsilon_2,\varepsilon_4,\varepsilon_5)]) \\
\nonumber
& = 
{\textstyle \sum}_{\varepsilon_4,\varepsilon_5}  (\wh{\bf M}^{\rm out}_{t,1} \wh{\bf M}^{\rm left}_{\rm tran}(\wh{Z}_{v_t}) \wh{\bf M}^{\rm in}_{t,3})_{\varepsilon_1,\varepsilon_4}
(\wh{\bf F}^{\rm out}_{+,\varepsilon_3})_{\varepsilon_4,\varepsilon_5}
(\wh{\bf M}^{\rm in}_{t,3} \wh{\bf M}^{\rm right}(\wh{Z}_{v_t}) \wh{\bf M}^{\rm out}_{t,2})_{\varepsilon_5,\varepsilon_2} \\
\label{eq:isotopy_invariance_proof_eq1}
& = (\wh{\bf M}^{\rm out}_{t,1} \wh{\bf M}^{\rm left}_{\rm tran}(\wh{Z}_{v_t}) \underline{ \wh{\bf M}^{\rm in}_{t,3} 
\wh{\bf F}^{\rm out}_{+,\varepsilon_3} \wh{\bf M}^{\rm in}_{t,3} } \wh{\bf M}^{\rm right}(\wh{Z}_{v_t}) \wh{\bf M}^{\rm out}_{t,2})_{\varepsilon_1,\varepsilon_2}.
\end{align}
To deal with the underlined part, we observe
\begin{lemma}
\label{lem:in_and_out_and_fork}
For each (side) $\alpha \in \{1,2,3\}$, (state value) $\varepsilon \in \{1,2,3\}$ and (sign) $k\in \{+,-\}$, one has
\begin{align*}
\wh{\bf M}^{\rm in}_{t,\alpha} 
\wh{\bf F}^{\rm out}_{k,\varepsilon} 
\wh{\bf M}^{\rm in}_{t,\alpha}
= (\wh{\bf M}^{\rm out}_{t,\alpha})_{\varepsilon,\varepsilon} \, \til{\bf F}^{\rm out}_{k,\varepsilon}, \qquad
\wh{\bf M}^{\rm out}_{t,\alpha} 
\wh{\bf F}^{\rm in}_{k,\varepsilon} 
\wh{\bf M}^{\rm out}_{t,\alpha}
= (\wh{\bf M}^{\rm in}_{t,\alpha})_{\varepsilon,\varepsilon} \, \til{\bf F}^{\rm in}_{k,\varepsilon}.
\end{align*}
\end{lemma}
{\it Proof of Lem.\ref{lem:in_and_out_and_fork}.} For each $\varepsilon_4,\varepsilon_5 \in\{1,2,3\}$ one has
$$
(\wh{\bf M}^{\rm in}_{t,\alpha} 
\wh{\bf F}^{\rm out}_{k,\varepsilon} 
\wh{\bf M}^{\rm in}_{t,\alpha})_{\varepsilon_4,\varepsilon_5} = (\wh{\bf M}^{\rm in}_{t,\alpha})_{\varepsilon_4,\varepsilon_4} 
(\wh{\bf F}^{\rm out}_{k,\varepsilon})_{\varepsilon_4,\varepsilon_5}
(\wh{\bf M}^{\rm in}_{t,\alpha})_{\varepsilon_5,\varepsilon_5}
$$
By definition of $\wh{\bf F}^{\rm out}_{k,\varepsilon}$ in eq.\eqref{eq:fork_matrix}, we have $(\wh{\bf F}^{\rm out}_{k,\varepsilon})_{\varepsilon_4,\varepsilon_5} \neq 0$ iff $\{\varepsilon_4,\varepsilon_5\} = \{r_1(\varepsilon),r_2(\varepsilon)\}$. Therefore $(\wh{\bf M}^{\rm in}_{t,\alpha} 
\wh{\bf F}^{\rm out}_{k,\varepsilon} 
\wh{\bf M}^{\rm in}_{t,\alpha})_{\varepsilon_4,\varepsilon_5} \neq 0$ iff $\{\varepsilon_4,\varepsilon_5\} = \{r_1(\varepsilon),r_2(\varepsilon)\}$ (since $\wh{\bf M}^{\rm in}_{t,\alpha}$ is diagonal), and the only nonzero entries of $\wh{\bf M}^{\rm in}_{t,\alpha} 
\wh{\bf F}^{\rm out}_{k,\varepsilon} 
\wh{\bf M}^{\rm in}_{t,\alpha}$ are
\begin{align*}
(\wh{\bf M}^{\rm in}_{t,\alpha} 
\wh{\bf F}^{\rm out}_{k,\varepsilon} 
\wh{\bf M}^{\rm in}_{t,\alpha})_{r_1(\varepsilon),r_2(\varepsilon)} & = (\wh{\bf M}^{\rm in}_{t,\alpha})_{r_1(\varepsilon),r_1(\varepsilon)} 
(\wh{\bf F}^{\rm out}_{k,\varepsilon})_{r_1(\varepsilon),r_2(\varepsilon)}
(\wh{\bf M}^{\rm in}_{t,\alpha})_{r_2(\varepsilon),r_2(\varepsilon)} \\
& = (\wh{\bf M}^{\rm in}_{t,\alpha})_{r_1(\varepsilon),r_1(\varepsilon)} 
(\wh{\bf M}^{\rm in}_{t,\alpha})_{r_2(\varepsilon),r_2(\varepsilon)}
(\wh{\bf F}^{\rm out}_{k,\varepsilon})_{r_1(\varepsilon),r_2(\varepsilon)}, \quad\mbox{and} \\
(\wh{\bf M}^{\rm in}_{t,\alpha} 
\wh{\bf F}^{\rm out}_{k,\varepsilon} 
\wh{\bf M}^{\rm in}_{t,\alpha})_{r_2(\varepsilon),r_1(\varepsilon)} 
& = (\wh{\bf M}^{\rm in}_{t,\alpha})_{r_2(\varepsilon),r_2(\varepsilon)} 
(\wh{\bf M}^{\rm in}_{t,\alpha})_{r_1(\varepsilon),r_1(\varepsilon)}
(\wh{\bf F}^{\rm out}_{k,\varepsilon})_{r_2(\varepsilon),r_1(\varepsilon)}
\end{align*}
Exchanging the roles of `${\rm in}$' and `${\rm out}$', we see that the only nonzero entries of $\wh{\bf M}^{\rm out}_{t,\alpha} 
\wh{\bf F}^{\rm in}_{k,\varepsilon} 
\wh{\bf M}^{\rm out}_{t,\alpha}$ are 
\begin{align*}
(\wh{\bf M}^{\rm out}_{t,\alpha} 
\wh{\bf F}^{\rm in}_{k,\varepsilon} 
\wh{\bf M}^{\rm out}_{t,\alpha})_{r_1(\varepsilon),r_2(\varepsilon)} 
& = (\wh{\bf M}^{\rm out}_{t,\alpha})_{r_1(\varepsilon),r_1(\varepsilon)} 
(\wh{\bf M}^{\rm out}_{t,\alpha})_{r_2(\varepsilon),r_2(\varepsilon)}
(\wh{\bf F}^{\rm in}_{k,\varepsilon})_{r_1(\varepsilon),r_2(\varepsilon)}, \quad\mbox{and} \\
(\wh{\bf M}^{\rm out}_{t,\alpha} 
\wh{\bf F}^{\rm in}_{k,\varepsilon} 
\wh{\bf M}^{\rm out}_{t,\alpha})_{r_2(\varepsilon),r_1(\varepsilon)} 
& = (\wh{\bf M}^{\rm out}_{t,\alpha})_{r_2(\varepsilon),r_2(\varepsilon)} 
(\wh{\bf M}^{\rm out}_{t,\alpha})_{r_1(\varepsilon),r_1(\varepsilon)}
(\wh{\bf F}^{\rm in}_{k,\varepsilon})_{r_2(\varepsilon),r_1(\varepsilon)}.
\end{align*}
Hence, Lem.\ref{lem:in_and_out_and_fork} follows from the following lemma, and the definition in eq.\eqref{eq:fork_matrix_tilde} of $\til{\bf F}^h_{k,\varepsilon}$.
\begin{lemma}[edge matrix inversion formula]
\label{lem:edge_matrix_inversion_formula}
For each $\alpha \in \{1,2,3\}$ and $\varepsilon\in\{1,2,3\}$, one has
\begin{align*}
&(\wh{\bf M}^{\rm in}_{t,\alpha})_{r_1(\varepsilon),r_1(\varepsilon)} 
(\wh{\bf M}^{\rm in}_{t,\alpha})_{r_2(\varepsilon),r_2(\varepsilon)}
= \omega^{\frac{3}{2} p(\varepsilon)} (\wh{\bf M}^{\rm out}_{t,\alpha})_{\varepsilon,\varepsilon}, \\
& (\wh{\bf M}^{\rm in}_{t,\alpha})_{r_2(\varepsilon),r_2(\varepsilon)} 
(\wh{\bf M}^{\rm in}_{t,\alpha})_{r_1(\varepsilon),r_1(\varepsilon)}
= \omega^{-\frac{3}{2} p(\varepsilon)} (\wh{\bf M}^{\rm out}_{t,\alpha})_{\varepsilon,\varepsilon}, \\
& (\wh{\bf M}^{\rm out}_{t,\alpha})_{r_1(\varepsilon),r_1(\varepsilon)} 
(\wh{\bf M}^{\rm out}_{t,\alpha})_{r_2(\varepsilon),r_2(\varepsilon)}
= \omega^{-\frac{3}{2} p(\varepsilon)} (\wh{\bf M}^{\rm in}_{t,\alpha})_{\varepsilon,\varepsilon},  \\
& (\wh{\bf M}^{\rm out}_{t,\alpha})_{r_2(\varepsilon),r_2(\varepsilon)} 
(\wh{\bf M}^{\rm out}_{t,\alpha})_{r_1(\varepsilon),r_1(\varepsilon)}
= \omega^{\frac{3}{2} p(\varepsilon)} (\wh{\bf M}^{\rm in}_{t,\alpha})_{\varepsilon,\varepsilon}, 
\end{align*}
where $p(\varepsilon)$ is as given in eq.\eqref{eq:p}, i.e. $p(1)=1=p(3)$, $p(2)=-1$. 
\end{lemma}

{\it Proof of Lem.\ref{lem:edge_matrix_inversion_formula}.} In view of eq.\eqref{eq:quantum_in_and_out}, if we write the three diagonal entries of $\wh{\bf M}^{\rm out}_{t,\alpha}$ as $\omega^{-1} \wh{Z}_1\wh{Z}_2^2$, $\omega^{1/2} \wh{Z}_1\wh{Z}_2^{-1}$, $\omega^{-1} \wh{Z}_1^{-2}\wh{Z}_2^{-1}$ in this order, then those of $\wh{\bf M}^{\rm in}_{t,\alpha}$ are $\omega \wh{Z}_2\wh{Z}_1^2$, $\omega^{-1/2} \wh{Z}_2 \wh{Z}_1^{-1}$, $\omega \wh{Z}_2^{-2}\wh{Z}_1^{-1}$ in this order. In view of eq.\eqref{eq:r1_and_r2}, for the first equality we check $(\omega \wh{Z}_2\wh{Z}_1^2)(\omega^{-1/2} \wh{Z}_2\wh{Z}_1^{-1}) =\omega^{3/2} (\omega^{-1}\wh{Z}_1 \wh{Z}_2^2 )$, $(\omega \wh{Z}_2\wh{Z}_1^2)(\omega\wh{Z}_2^{-2}\wh{Z}_1^{-1}) =\omega^{-3/2} (\omega^{1/2}\wh{Z}_1\wh{Z}_2^{-1})$, and $(\omega^{-1/2} \wh{Z}_2\wh{Z}_1^{-1})(\omega \wh{Z}_2^{-2}\wh{Z}_1^{-1}) = \omega^{3/2} (\omega^{-1} \wh{Z}_1^{-2}\wh{Z}_2^{-1})$. Other checks are similar. \qquad {\it End of proof of Lem.\ref{lem:edge_matrix_inversion_formula}, and Lem.\ref{lem:in_and_out_and_fork}.}

\vs

Coming back to our situation, from eq.\eqref{eq:isotopy_invariance_proof_eq1} and Lem.\ref{lem:in_and_out_and_fork} we have
\begin{align*}
\til{\rm Tr}^\omega_{\wh{\Delta}}(W',s')  = (\wh{\bf M}^{\rm out}_{t,1} \wh{\bf M}^{\rm left}_{\rm tran}(\wh{Z}_{v_t}) (\wh{\bf M}^{\rm out}_{t,3})_{\varepsilon_3,\varepsilon_3} 
\til{\bf F}^{\rm out}_{+,\varepsilon_3}  \wh{\bf M}^{\rm right}(\wh{Z}_{v_t}) \wh{\bf M}^{\rm out}_{t,2})_{\varepsilon_1,\varepsilon_2}.
\end{align*}
Notice that $(\wh{\bf M}^{\rm out}_{t,3})_{\varepsilon_3,\varepsilon_3} $ is not a matrix, but just an element of $\mathcal{Z}^\omega_t$. To change the order of product, we use Lem.\ref{lem:triangle_variable_and_edge_matrices}:
\begin{align*}
\til{\rm Tr}^\omega_{\wh{\Delta}}(W',s')  = (\wh{\bf M}^{\rm out}_{t,1} (\wh{\bf M}^{\rm out}_{t,3})_{\varepsilon_3,\varepsilon_3}
\wh{\bf M}^{\rm left}_{\rm tran}(\omega^{2g(\varepsilon_3)} \wh{Z}_{v_t})  
\til{\bf F}^{\rm out}_{+,\varepsilon_3}  \wh{\bf M}^{\rm right}(\wh{Z}_{v_t}) \wh{\bf M}^{\rm out}_{t,2})_{\varepsilon_1,\varepsilon_2}.
\end{align*}
Since $\wh{\bf M}^{\rm out}_{t,\alpha}$ are diagonal, we thus have
\begin{align*}
\redfix{\til{\rm Tr}^\omega_{\wh{\Delta}}(W',s')=}(\wh{\bf M}^{\rm out}_{t,1})_{\varepsilon_1, \varepsilon_1} 
(\wh{\bf M}^{\rm out}_{t,3})_{\varepsilon_3, \varepsilon_3} 
( \wh{\bf M}^{\rm left}_{\rm tran}(\omega^{2g(\varepsilon_3)} \wh{Z}_{v_t} )
\til{\bf F}^{\rm out}_{+,\varepsilon_3}
\wh{\bf M}^{\rm right}(\wh{Z}_{v_t})  )_{\varepsilon_1, \varepsilon_2}
(\wh{\bf M}^{\rm out}_{t,2})_{\varepsilon_2, \varepsilon_2}\redfix{,}
\end{align*}
which equals $(*)^{\rm out}_{\varepsilon_1,\varepsilon_2,\varepsilon_3}$ appearing in the proof of Lem.\ref{lem:cyclicity_of_3-way_matrices}, hence equals ${\rm Tr}^\omega_{\Delta;t}([W,s])$. So we get the desired equality
$$
\til{\rm Tr}^\omega_{\wh{\Delta}}(W',s') = {\rm Tr}^\omega_{\Delta;t}([W,s]) = \til{\rm Tr}^\omega_{\wh{\Delta}}(W,s).
$$
This is in fact partly how we came up with the values for Thm.\ref{thm:SL3_quantum_trace_map}\redfix{(QT2-4)}, which may look a bit strangely technical at a first glance. Anyhow, this finishes the proof for the case of Fig.\ref{fig:isotopy_move1}.

\vs

For the case as in Fig.\ref{fig:isotopy_move1} with the reverse orientation (with $w'\succ w''$), the proof goes similarly. In particular, the arguments \redfix{go} almost verbatim for $\til{\rm Tr}^\omega_{\wh{\Delta}}(W,s)$, while our investigation of $\til{\rm Tr}^\omega_{\wh{\Delta}}(W',s')$ should change as follows. First, eq.\eqref{eq:isotopy_invariance_proof_biangle_factor1} should be replaced by
\begin{align}
\nonumber
{\rm Tr}^\omega_B([W'\cap (B \times {\bf I}),(\varepsilon_3,\varepsilon_4,\varepsilon_5)])
= (\wh{\bf F}^{\rm in}_{+,\varepsilon_3})_{\varepsilon_4,\varepsilon_5}
\end{align}
while eq.\eqref{eq:isotopy_invariance_proof_biangle_factor2} now becomes
\begin{align}
\nonumber
{\rm Tr}^\omega_t([W'\cap (\wh{t} \times {\bf I}), (\varepsilon_1,\varepsilon_2,\varepsilon_4,\varepsilon_5)])
& = (\wh{\bf M}^{\rm in}_{t,1} \wh{\bf M}^{\rm right}(\wh{Z}_{v_t}) \wh{\bf M}^{\rm out}_{t,3})_{\varepsilon_1,\varepsilon_4}
(\wh{\bf M}^{\rm in}_{t,2} \wh{\bf M}^{\rm left}(\wh{Z}_{v_t}) \wh{\bf M}^{\rm out}_{t,3})_{\varepsilon_2,\varepsilon_5} \\
\nonumber
& = (\wh{\bf M}^{\rm in}_{t,1} \wh{\bf M}^{\rm right}(\wh{Z}_{v_t}) \wh{\bf M}^{\rm out}_{t,3})_{\varepsilon_1,\varepsilon_4}
(\wh{\bf M}^{\rm out}_{t,3} \wh{\bf M}^{\rm left}_{\rm tran}(\wh{Z}_{v_t}) \wh{\bf M}^{\rm in}_{t,2})_{\varepsilon_5,\varepsilon_2} ~ (\because\mbox{Cor.\ref{cor:transpose_of_quantum_turn_matrices}}),
\end{align}
hence
\begin{align*}
\til{\rm Tr}^\omega_{\wh{\Delta}}(W',s') & = {\textstyle \sum}_{\varepsilon_4,\varepsilon_5} {\rm Tr}^\omega_B([W'\cap (B \times {\bf I}),(\varepsilon_3,\varepsilon_4,\varepsilon_5)])
{\rm Tr}^\omega_t([W'\cap (\wh{t} \times {\bf I}), (\varepsilon_1,\varepsilon_2,\varepsilon_4,\varepsilon_5)])
 \\
& = {\textstyle \sum}_{\varepsilon_4,\varepsilon_5} (\wh{\bf M}^{\rm in}_{t,1} \wh{\bf M}^{\rm right} (\wh{Z}_{v_t}) \wh{\bf M}^{\rm out}_{t,3})_{\varepsilon_1,\varepsilon_4}
(\wh{\bf F}^{\rm in}_{+,\varepsilon_3})_{\varepsilon_4,\varepsilon_5}
(\wh{\bf M}^{\rm out}_{t,3} \wh{\bf M}^{\rm left}_{\rm tran}(\wh{Z}_{v_t}) \wh{\bf M}^{\rm in}_{t,2})_{\varepsilon_5,\varepsilon_2} \\
& = (\wh{\bf M}^{\rm in}_{t,1} \wh{\bf M}^{\rm right}(\wh{Z}_{v_t}) \underline{ \wh{\bf M}^{\rm out}_{t,3} 
\wh{\bf F}^{\rm in}_{+,\varepsilon_3} \wh{\bf M}^{\rm out}_{t,3} } \wh{\bf M}^{\rm left}_{\rm tran}(\wh{Z}_{v_t}) \wh{\bf M}^{\rm in}_{t,2})_{\varepsilon_1,\varepsilon_2} \\
& = (\wh{\bf M}^{\rm in}_{t,1} \underline{ \wh{\bf M}^{\rm right}(\wh{Z}_{v_t}) 
(\wh{\bf M}^{\rm in}_{t,3})_{\varepsilon_3,\varepsilon_3} }
\til{\bf F}^{\rm in}_{+,\varepsilon_3} 
\wh{\bf M}^{\rm left}_{\rm tran}(\wh{Z}_{v_t}) \wh{\bf M}^{\rm in}_{t,2})_{\varepsilon_1,\varepsilon_2} \quad (\because\mbox{Lem.\ref{lem:in_and_out_and_fork}}) \\
& = (\wh{\bf M}^{\rm in}_{t,1} (\wh{\bf M}^{\rm in}_{t,3})_{\varepsilon_3,\varepsilon_3} 
\wh{\bf M}^{\rm right}(\omega^{-2g(\varepsilon_3)} \wh{Z}_{v_t}) 
\til{\bf F}^{\rm in}_{+,\varepsilon_3} 
\wh{\bf M}^{\rm left}_{\rm tran}(\wh{Z}_{v_t}) \wh{\bf M}^{\rm in}_{t,2})_{\varepsilon_1,\varepsilon_2}, \quad (\because\mbox{Lem.\ref{lem:triangle_variable_and_edge_matrices}})
\end{align*}
which equals $(*)^{\rm in}_{\varepsilon_1,\varepsilon_2,\varepsilon_3}$ appearing in the proof of Lem.\ref{lem:cyclicity_of_3-way_matrices}, hence equals ${\rm Tr}^\omega_{\Delta;t}([W,s])$. Therefore we get the desired equality $\til{\rm Tr}^\omega_{\wh{\Delta}}(W',s') = {\rm Tr}^\omega_{\Delta;t}([W,s])=\til{\rm Tr}^\omega_{\wh{\Delta}}(W,s)$ for this reverse-orientation case for Fig.\ref{fig:isotopy_move1}.

\vs

Now take the case as in Fig.\ref{fig:isotopy_move2} with $z_2\succ z_1$ and $w_2\succ w_1$, possibly with all orientations reversed. Denote by $\varepsilon_1,\varepsilon_2,\varepsilon_3$ the state values of $s$ and $s'$ at endpoints $z_1,z_2,x$, respectively. Look at $W'$ first, which is on the right (i.e. the second or the fourth picture from the left in Fig.\ref{fig:isotopy_move2}). Note $W'\cap (B\times {\bf I})$ is \redfix{a} product of two edges connecting distinct sides, so in the sum in eq.\eqref{eq:til_Tr_W_s_prime}, by Prop.\ref{prop:biangle_SL3_quantum_trace}(BT2-1) the biangle factor `goes away', and we just have
\begin{align}
\label{eq:isotopy_invariance_proof_fork1}
\til{\rm Tr}^\omega_{\wh{\Delta}}(W',s') = \til{\rm Tr}^\omega_t(W' \cap (\wh{t} \times {\bf I}), (\varepsilon_1,\varepsilon_2,\varepsilon_3))
= {\rm Tr}^\omega_{\Delta;t}([W', (\varepsilon_1,\varepsilon_2,\varepsilon_3)])
\end{align}
where $(\varepsilon_1,\varepsilon_2,\varepsilon_3)$ denotes the state $s_2'$ of $W' \cap (\wh{t} \times {\bf I})$ that assigns $\varepsilon_1,\varepsilon_2,\varepsilon_3$ to $w_1,w_2,x$, and also denotes the state $s'$ of $W'$ that assigns $\varepsilon_1,\varepsilon_2,\varepsilon_3$ to $z_1,z_2,x$. On the other hand, consider $W$, which is on the left (i.e. the first or the third picture from the left in Fig.\ref{fig:isotopy_move2}). Note $W \cap (B\times {\bf I})$ consists of one 3-way ${\rm SL}_3$-web component. If the state $s_1$ of $W\cap (B\times {\bf I})$ assigns $\varepsilon, \varepsilon_1,\varepsilon_2$ to $w,z_1,z_2$, then by Prop.\ref{prop:biangle_SL3_quantum_trace}(BT2-3) we see that ${\rm Tr}^\omega_B([W\cap (B\times {\bf I}),s_1]) \neq 0$ iff $\{r_1(\varepsilon),r_2(\varepsilon)\}=\{\varepsilon_1,\varepsilon_2\}$, and the value is $(\wh{\bf F}^h_{-,\varepsilon})_{\varepsilon_1,\varepsilon_2}$ (eq.\eqref{eq:fork_matrix}), where $h \in \{{\rm in},{\rm out}\}$ indicates whether $W$ is an incoming or an outgoing $3$-way ${\rm SL}_3$-web. In particular, if there exists no $\varepsilon \in\{1,2,3\}$ with $\{r_1(\varepsilon),r_2(\varepsilon)\}=\{\varepsilon_1,\varepsilon_2\}$, then in view of eq.\eqref{eq:til_Tr_W_s} we have $\til{\rm Tr}^\omega_{\wh{\Delta}}(W,s)=0$. If there is such an $\varepsilon$, then it is unique, and in view of eq.\eqref{eq:til_Tr_W_s} we have 
\begin{align}
\label{eq:isotopy_invariance_proof_fork2}
\til{\rm Tr}^\omega_{\wh{\Delta}}(W,s)
= {\rm Tr}^\omega_B([W\cap (B\times {\bf I}),(\varepsilon_1,\varepsilon_2,\varepsilon_3)]) \, {\rm Tr}^\omega_t([W\cap (\wh{t} \times {\bf I}), (\varepsilon, \varepsilon_1)])
\end{align}
where $(\varepsilon_1,\varepsilon_2,\varepsilon_3)$ denotes the state $s_1$ of $W\cap (B\times {\bf I})$ that assigns $\varepsilon_1,\varepsilon_2,\varepsilon_3$ to $z_1,z_2,w$, and $(\varepsilon,\varepsilon_1)$ denotes the state $s_2$ of $W\cap (\wh{t} \times {\bf I})$ that assigns $\varepsilon,\varepsilon_1$ to $w,x$. Since ${\rm Tr}^\omega_B([W\cap (B\times {\bf I}),(\varepsilon_1,\varepsilon_2,\varepsilon_3)]) = (\wh{\bf F}^h_{-,\varepsilon})_{\varepsilon_1,\varepsilon_2}$, in view of Thm.\ref{thm:SL3_quantum_trace_map}\redfix{(QT2-3)} we observe that $\til{\rm Tr}^\omega_{\wh{\Delta}}(W,s)$ as in eq.\eqref{eq:isotopy_invariance_proof_fork2} equals $\til{\rm Tr}^\omega_{\wh{\Delta}}(W',s')$ as in eq.\eqref{eq:isotopy_invariance_proof_fork1}, as desired. When $z_2\prec z_1$ and $w_2\prec w_1$, just replace $\wh{\bf F}^h_{-,\varepsilon}$ by $\wh{\bf F}^h_{+,\varepsilon}$; other arguments are identical.

\vs

Next, consider the left case of Fig.\ref{fig:isotopy_move3} (with $x_1\succ x_2$ and $w_1\succ w_2$). For $W'$, note that $W' \cap (B \times {\bf I})$ consists of a single component of type (BT2-1), so the biangle factor goes away and we have
\begin{align}
\label{eq:isotopy_invariance_proof_fork3}
\til{\rm Tr}^\omega_{\wh{\Delta}}(W',s') = \til{\rm Tr}^\omega_t(W' \cap (\wh{t} \times {\bf I}), (\varepsilon_1,\varepsilon_2,\varepsilon_3))
= {\rm Tr}^\omega_{\Delta;t}([W', (\varepsilon_1,\varepsilon_2,\varepsilon_3)]),
\end{align}
where $(\varepsilon_1,\varepsilon_2,\varepsilon_3)$ denotes the state $s_2'$ of $W' \cap (\wh{t} \times {\bf I})$ that assigns $\varepsilon_1,\varepsilon_2,\varepsilon_3$ to $x_1,x_2,w$, and also denotes the state $s'$ of $W'$ that assigns $\varepsilon_1,\varepsilon_2,\varepsilon_3$ to $x_1,x_2,z$. Note $W'\cap (\wh{t} \times {\bf I})$ and $W'$ fall into Thm.\ref{thm:SL3_quantum_trace_map}\redfix{(QT2-3)}, so
\begin{align}
\label{eq:isotopy_invariance_proof_fork5}
\til{\rm Tr}^\omega_{\wh{\Delta}}(W',s')
= (\wh{\bf F}^{\rm out}_{+,\varepsilon})_{\varepsilon_1,\varepsilon_2} ( \wh{\bf M}^{\rm in}_{t,1} \wh{\bf M}^{\rm right}(\wh{Z}_{v_t}) \wh{\bf M}^{\rm out}_{t,3} )_{\varepsilon,\varepsilon_3},
\end{align}
if $\{r_1(\varepsilon),r_2(\varepsilon)\} = \{\varepsilon_1,\varepsilon_2\}$, while $\til{\rm Tr}^\omega_{\wh{\Delta}}(W',s')=0$ if there is no such $\varepsilon$. Meanwhile, for $W$, we have the state sum
\begin{align}
\label{eq:isotopy_invariance_proof_fork4}
\til{\rm Tr}^\omega_{\wh{\Delta}}(W,s)
= {\textstyle \sum}_{\varepsilon_4,\varepsilon_5} {\rm Tr}^\omega_B([W\cap (B\times {\bf I}),(\varepsilon_4,\varepsilon_5,\varepsilon_3)]) \, {\rm Tr}^\omega_t([W\cap (\wh{t} \times {\bf I}), (\varepsilon_1,\varepsilon_2,\varepsilon_4,\varepsilon_5)])
\end{align}
where $(\varepsilon_4,\varepsilon_5,\varepsilon_3)$ denotes the state of $s_2$ that assigns $\varepsilon_4,\varepsilon_5,\varepsilon_3$ to $w_1,w_2,z$, and $(\varepsilon_1,\varepsilon_2,\varepsilon_4,\varepsilon_5)$ denotes the state of $s_1$ that assigns $\varepsilon_1,\varepsilon_2,\varepsilon_4,\varepsilon_5$ to $x_1,x_2,w_1,w_2$, and the sum is over all $\varepsilon_4,\varepsilon_5 \in \{1,2,3\}$. Since $W\cap (B\times {\bf I})$ falls into Prop.\ref{prop:biangle_SL3_quantum_trace}(BT2-3), and $W\cap (\wh{t} \times {\bf I})$ is a product of two left turns as in Thm.\ref{thm:SL3_quantum_trace_map}\redfix{(QT2-1)}, we have
\begin{align*}
\til{\rm Tr}^\omega_{\wh{\Delta}}(W,s)
& = {\textstyle \sum}_{\varepsilon_4,\varepsilon_5} (\wh{\bf F}^{\rm out}_{+,\varepsilon_3})_{\varepsilon_4,\varepsilon_5} \, (\wh{\bf M}^{\rm in}_{t,3} \wh{\bf M}^{\rm left}(\wh{Z}_{v_t}) \wh{\bf M}^{\rm out}_{t,1})_{\varepsilon_4, \varepsilon_1} (\wh{\bf M}^{\rm in}_{t,3} \wh{\bf M}^{\rm left}(\wh{Z}_{v_t}) \wh{\bf M}^{\rm out}_{t,1})_{\varepsilon_5, \varepsilon_2} \\
& = {\textstyle \sum}_{\varepsilon_4,\varepsilon_5} (\wh{\bf M}^{\rm out}_{t,1} \wh{\bf M}^{\rm left}_{\rm tran}(\wh{Z}_{v_t}) \wh{\bf M}^{\rm in}_{t,3})_{\varepsilon_1, \varepsilon_4} (\wh{\bf F}^{\rm out}_{+,\varepsilon_3})_{\varepsilon_4,\varepsilon_5} (\wh{\bf M}^{\rm in}_{t,3} \wh{\bf M}^{\rm left}(\wh{Z}_{v_t}) \wh{\bf M}^{\rm out}_{t,1})_{\varepsilon_5, \varepsilon_2} \quad (\because \mbox{Cor.\ref{cor:transpose_of_quantum_turn_matrices}}) \\
& = 
(\wh{\bf M}^{\rm out}_{t,1} \wh{\bf M}^{\rm left}_{\rm tran}(\wh{Z}_{v_t}) \underline{ \wh{\bf M}^{\rm in}_{t,3}
\wh{\bf F}^{\rm out}_{+,\varepsilon_3}
\wh{\bf M}^{\rm in}_{t,3} } \wh{\bf M}^{\rm left}(\wh{Z}_{v_t}) \wh{\bf M}^{\rm out}_{t,1})_{\varepsilon_1, \varepsilon_2}  \\
& =
(\wh{\bf M}^{\rm out}_{t,1} \wh{\bf M}^{\rm left}_{\rm tran}(\wh{Z}_{v_t}) \underline{  (\wh{\bf M}^{\rm out}_{t,3})_{\varepsilon_3,\varepsilon_3} 
\til{\bf F}^{\rm out}_{+,\varepsilon_3}
\wh{\bf M}^{\rm left}(\wh{Z}_{v_t})  }\wh{\bf M}^{\rm out}_{t,1})_{\varepsilon_1, \varepsilon_2} 
\quad (\because \mbox{Lem.\ref{lem:in_and_out_and_fork}}) \\
& = 
(\wh{\bf M}^{\rm out}_{t,1}
\underbrace{  \wh{\bf M}^{\rm left}_{\rm tran}(\wh{Z}_{v_t}) 
\til{\bf F}^{\rm out}_{+,\varepsilon_3}
\wh{\bf M}^{\rm left}(\omega^{-2g(\varepsilon_3)} \wh{Z}_{v_t}) }_{=: \wh{\bf M}^{\rm lt,o,l}_{+,\varepsilon_3} } (\wh{\bf M}^{\rm out}_{t,3})_{\varepsilon_3,\varepsilon_3} 
\wh{\bf M}^{\rm out}_{t,1})_{\varepsilon_1, \varepsilon_2}. 
\quad (\because \mbox{Lem.\ref{lem:triangle_variable_and_edge_matrices}})
\end{align*}
We compute the underbraced matrix product $\wh{\bf M}^{\rm lt,o,l}_{+,\varepsilon_3}$:
\begin{align*}
& \wh{\bf M}^{\rm lt,o,l}_{+,1}
= \smallmatthree{\omega^{-5} \wh{Z}_{v_t}^2}{0}{0}{\omega \wh{Z}^2 + \omega^{-2} \wh{Z}_{v_t}^{-1} ~}{\omega^{-5} \wh{Z}_{v_t}^{-1}}{0}{\omega^{4} \wh{Z}_{v_t}^{-1}}{\omega \wh{Z}_{v_t}^{-1}}{\omega^{-2} \wh{Z}_{v_t}^{-1}}
\smallmatthree{0}{\omega^3}{0}{-\omega^9}{0}{0}{0}{0}{0}
\smallmatthree{\omega^5 (\omega^2\wh{Z}_{v_t})^2}{~ \omega^{-1} (\omega^2\wh{Z}_{v_t})^2 + \omega^{2} (\omega^2\wh{Z}_{v_t})^{-1} ~}{\omega^{-4} (\omega^2\wh{Z}_{v_t})^{-1}}{0}{\omega^5 (\omega^2\wh{Z}_{v_t})^{-1}}{\omega^{-1} (\omega^2\wh{Z}_{v_t})^{-1}}{0}{0}{\omega^2 (\omega^2\wh{Z}_{v_t})^{-1}} \\
& = \smallmatthree{\omega^{-5} \wh{Z}_{v_t}^2}{0}{0}{\omega \wh{Z}^2 + \omega^{-2} \wh{Z}_{v_t}^{-1} ~}{\omega^{-5} \wh{Z}_{v_t}^{-1}}{0}{\omega^{4} \wh{Z}_{v_t}^{-1}}{\omega \wh{Z}_{v_t}^{-1}}{\omega^{-2} \wh{Z}_{v_t}^{-1}}
\smallmatthree{0}{\omega^6 \wh{Z}_{v_t}^{-1}}{ \wh{Z}_{v_t}^{-1}}{-\omega^{18} \wh{Z}_{v_t}^2}{~(-\omega^{12} \wh{Z}_{v_t}^2 -\omega^9 \wh{Z}_{v_t}^{-1}) ~}{-\omega^{3} \wh{Z}_{v_t}^{-1}}{0}{0}{0}
= \smallmatthree{0}{\omega \wh{Z}_{v_t}}{\omega^{-5} \wh{Z}_{v_t}}{- \omega^{13} \wh{Z}_{v_t}}{0}{\omega \wh{Z}_{v_t} }{- \omega^{19} \wh{Z}_{v_t}}{ - \omega^{13} \wh{Z}_{v_t}}{0}, \\
& \wh{\bf M}^{\rm lt,o,l}_{+,2}
= \smallmatthree{0}{0}{\omega \wh{Z}_{v_t}}{0}{0}{\omega^7 \wh{Z}_{v_t} + \omega^4 \wh{Z}_{v_t}^{-2}}{-\omega^7 \wh{Z}_{v_t}}{(- \omega \wh{Z}_{v_t} - \omega^{16} \wh{Z}_{v_t}^{-2})}{0}, 
\qquad
\wh{\bf M}^{\rm lt,o,l}_{+,3} = \smallmatthree{0}{0}{0}{0}{0}{\omega^{-2} \wh{Z}_{v_t}^{-2}}{0}{-\omega^{10} \wh{Z}_{v_t}^{-2}}{0}.
\end{align*}
One can write
\begin{align*}
\til{\rm Tr}^\omega_{\wh{\Delta}}(W,s)
= ( \wh{\bf M}^{\rm out}_{t,1} \, 
\wh{\bf M}^{\rm lt,o,l}_{+,\varepsilon_3} \,
(\wh{\bf M}^{\rm out}_{t,3})_{\varepsilon_3,\varepsilon_3} \,
\wh{\bf M}^{\rm out}_{t,1})_{\varepsilon_1,\varepsilon_2} = ( \wh{\bf M}^{\rm out}_{t,1})_{\varepsilon_1,\varepsilon_1} \, 
(\wh{\bf M}^{\rm lt,o,l}_{+,\varepsilon_3})_{\varepsilon_1,\varepsilon_2} \,
(\wh{\bf M}^{\rm out}_{t,3})_{\varepsilon_3,\varepsilon_3} \,
(\wh{\bf M}^{\rm out}_{t,1})_{\varepsilon_2,\varepsilon_2},
\end{align*}
since $\wh{\bf M}^{\rm out}_{t,\alpha}$ are diagonal. By inspection, $(\wh{\bf M}^{\rm lt,o,l}_{+,\varepsilon_3})_{\varepsilon_1,\varepsilon_2} = 0$ \redfix{if} $\varepsilon_1 =\varepsilon_2$. Hence $\til{\rm Tr}^\omega_{\wh{\Delta}}(W,s)=0$ if $\varepsilon_1=\varepsilon_2$. Suppose $\varepsilon_1 \neq \varepsilon_2$, so that there is $\varepsilon$ s.t. $\{r_1(\varepsilon),r_2(\varepsilon)\} = \{\varepsilon_1,\varepsilon_2\}$. First, we change the order of the product $(\wh{\bf M}^{\rm out}_{t,3})_{\varepsilon_3,\varepsilon_3} \,
(\wh{\bf M}^{\rm out}_{t,1})_{\varepsilon_2,\varepsilon_2}$, using the following lemma, which is straightforward to verify:
\begin{lemma}
\label{lem:edge_matrices_commutation}
For any $\alpha \in \{1,2,3\}$ and $\varepsilon',\varepsilon'' \in \{1,2,3\}$, one has
\begin{align*}
(\wh{\bf M}^{\rm out}_{t,\alpha})_{\varepsilon',\varepsilon'} 
(\wh{\bf M}^{\rm out}_{t,\alpha+1})_{\varepsilon'', \varepsilon''} 
& = \omega^{2 g(\varepsilon'+1) g(\varepsilon''-1)} (\wh{\bf M}^{\rm out}_{t,\alpha+1})_{\varepsilon'', \varepsilon''}(\wh{\bf M}^{\rm out}_{t,\alpha})_{\varepsilon',\varepsilon'}, \\
(\wh{\bf M}^{\rm in}_{t,\alpha})_{\varepsilon',\varepsilon'} 
(\wh{\bf M}^{\rm in}_{t,\alpha+1})_{\varepsilon'', \varepsilon''} 
& = \omega^{2 g(\varepsilon'-1) g(\varepsilon''+1)} (\wh{\bf M}^{\rm in}_{t,\alpha+1})_{\varepsilon'', \varepsilon''}(\wh{\bf M}^{\rm in}_{t,\alpha})_{\varepsilon',\varepsilon'}
\end{align*}
where $g$ is as in eq.\eqref{eq:g}, with the argument \redfix{of $g$} being understood modulo $3$, i.e. $g(0)=-1=g(4)$. \qed
\end{lemma}
So
\begin{align*}
\til{\rm Tr}^\omega_{\wh{\Delta}}(W,s)
= \omega^{2 g(\varepsilon_3+1) g(\varepsilon_2-1)} ( \wh{\bf M}^{\rm out}_{t,1})_{\varepsilon_1,\varepsilon_1} \, 
(\wh{\bf M}^{\rm lt,o,l}_{+,\varepsilon_3})_{\varepsilon_1,\varepsilon_2} \,
(\wh{\bf M}^{\rm out}_{t,1})_{\varepsilon_2,\varepsilon_2} \,
(\wh{\bf M}^{\rm out}_{t,3})_{\varepsilon_3,\varepsilon_3}.
\end{align*}
Now, using Lem.\ref{lem:triangle_variable_and_edge_matrices}, we have
$$
(\wh{\bf M}^{\rm lt,o,l}_{+,\varepsilon_3})_{\varepsilon_1,\varepsilon_2} \,
(\wh{\bf M}^{\rm out}_{t,1})_{\varepsilon_2,\varepsilon_2}
= (\wh{\bf M}^{\rm out}_{t,1})_{\varepsilon_2,\varepsilon_2} (\wh{\bf M}^{\rm lt,o,l, col}_{+,\varepsilon_3})_{\varepsilon_1,\varepsilon_2} \,
$$
where $\wh{\bf M}^{\rm lt,o,l, col}_{+,\varepsilon_3}$ is obtained from $\wh{\bf M}^{\rm lt,o,l}_{+,\varepsilon_3}$ by replacing $\wh{Z}_{v_t}$ appearing in each $\varepsilon_2$-th column with $\omega^{2g(\varepsilon_2)} \wh{Z}_{v_t}$:
\begin{align*}
\hspace{0mm}
\begin{array}{lll}
\wh{\bf M}^{\rm lt,o,l,col}_{+,1} &
\wh{\bf M}^{\rm lt,o,l,col}_{+,2} 
& 
\wh{\bf M}^{\rm lt,o,l,col}_{+,3}  \\
= \smallmatthree{0}{\omega^5 \wh{Z}_{v_t}}{\omega^{-7} \wh{Z}_{v_t}}{- \omega^{11} \wh{Z}_{v_t}}{0}{\omega^{-1} \wh{Z}_{v_t} }{- \omega^{17} \wh{Z}_{v_t}}{ - \omega^{17} \wh{Z}_{v_t}}{0}
& 
= 
\smallmatthree{0}{0}{\omega^{-1} \wh{Z}_{v_t}}{0}{0}{\omega^5 \wh{Z}_{v_t} + \omega^8 \wh{Z}_{v_t}^{-2}}{-\omega^5 \wh{Z}_{v_t}}{(- \omega^5 \wh{Z}_{v_t} - \omega^{8} \wh{Z}_{v_t}^{-2})}{0}
& 
= \smallmatthree{0}{0}{0}{0}{0}{\omega^{2} \wh{Z}_{v_t}^{-2}}{0}{-\omega^{2} \wh{Z}_{v_t}^{-2}}{0}
\end{array}
\end{align*}
Note
\begin{align*}
\til{\rm Tr}^\omega_{\wh{\Delta}}(W,s)
& = \omega^{2 g(\varepsilon_3+1) g(\varepsilon_2-1)} \underline{ ( \wh{\bf M}^{\rm out}_{t,1})_{\varepsilon_1,\varepsilon_1} \, 
(\wh{\bf M}^{\rm out}_{t,1})_{\varepsilon_2,\varepsilon_2} } \,
(\wh{\bf M}^{\rm lt,o,l,col}_{+,\varepsilon_3})_{\varepsilon_1,\varepsilon_2} \,
(\wh{\bf M}^{\rm out}_{t,3})_{\varepsilon_3,\varepsilon_3} \\
& = 
\omega^{2 g(\varepsilon_3+1) g(\varepsilon_2-1)} \omega^{\frac{3}{2}p(\varepsilon) {\rm sgn}(\varepsilon_1-\varepsilon_2)} ( \wh{\bf M}^{\rm in}_{t,1})_{\varepsilon,\varepsilon} \,
(\wh{\bf M}^{\rm lt,o,l,col}_{+,\varepsilon_3})_{\varepsilon_1,\varepsilon_2} \,
(\wh{\bf M}^{\rm out}_{t,3})_{\varepsilon_3,\varepsilon_3}.
\qquad (\because\mbox{Lem.\ref{lem:edge_matrix_inversion_formula}})
\end{align*}

From eq.\eqref{eq:isotopy_invariance_proof_fork5} we have
$$
\til{\rm Tr}^\omega_{\wh{\Delta}}(W',s') 
= (\wh{\bf F}^{\rm out}_{+,\varepsilon})_{\varepsilon_1,\varepsilon_2} ( \wh{\bf M}^{\rm in}_{t,1})_{\varepsilon,\varepsilon} (\wh{\bf M}^{\rm right}(\wh{Z}_{v_t}) )_{\varepsilon,\varepsilon_3}
(\wh{\bf M}^{\rm out}_{t,3} )_{\varepsilon,\varepsilon_3},
$$
and hence it suffices to show
$$
\omega^{2 g(\varepsilon_3+1) g(\varepsilon_2-1)} \omega^{\frac{3}{2}p(\varepsilon) {\rm sgn}(\varepsilon_1-\varepsilon_2)}
(\wh{\bf M}^{\rm lt,o,l,col}_{+,\varepsilon_3})_{\varepsilon_1,\varepsilon_2}
= (\wh{\bf F}^{\rm out}_{+,\varepsilon})_{\varepsilon_1,\varepsilon_2} (\wh{\bf M}^{\rm right}(\wh{Z}_{v_t}))_{\varepsilon,\varepsilon_3},
$$
or equivalently, in view of the definition in eq.\eqref{eq:fork_matrix_tilde} of $\til{\bf F}^h_{k,\varepsilon}$, to show
$$
\omega^{2 g(\varepsilon_3+1) g(\varepsilon_2-1)}
(\wh{\bf M}^{\rm lt,o,l,col}_{+,\varepsilon_3})_{\varepsilon_1,\varepsilon_2}
= (\til{\bf F}^{\rm out}_{+,\varepsilon})_{\varepsilon_1,\varepsilon_2} (\wh{\bf M}^{\rm right}(\wh{Z}_{v_t}))_{\varepsilon,\varepsilon_3},
$$
For each $\varepsilon_3\in \{1,2,3\}$, define $\wh{\bf M}^{\rm r}_{+,\varepsilon_3}$ as the $3\times 3$ matrix whose $(\varepsilon_1,\varepsilon_2)$-th entry equals zero if $\varepsilon_1=\varepsilon_2$, and equals $\omega^{-2g(\varepsilon_3+1)g(\varepsilon_2-1)} (\wh{\bf M}^{\rm right}(\wh{Z}_{v_t}))_{\varepsilon,\varepsilon_3}$ if $\varepsilon_1\neq \varepsilon_2$, where $\varepsilon$ is determined by $\{r_1(\varepsilon),r_2(\varepsilon)\} = \{\varepsilon_1,\varepsilon_2\}$; that is, $\wh{\bf M}^{\rm r}_{+,\varepsilon_3}$ is made from entries of the $\varepsilon_3$-th column of $\wh{\bf M}^{\rm right}(\wh{Z}_{v_t})$ in a certain way. We compute these matrices:
\begin{align*}
\begin{array}{lll}
\wh{\bf M}^{\rm r}_{+,1} 
&
\wh{\bf M}^{\rm r}_{+,2}
&
\wh{\bf M}^{\rm r}_{+,3} \\
= 
\smallmatthree{0}{\omega^{2} \wh{Z}_{v_t}}{\omega^{-7} \wh{Z}_{v_t}}{\omega^{2} \wh{Z}_{v_t}}{0}{\omega^{-4} \wh{Z}_{v_t}}{\omega^5 \wh{Z}_{v_t}}{\omega^8 \wh{Z}_{v_t}}{0}
&
 = 
\smallmatthree{0}{0}{\omega^{-1} \wh{Z}_{v_t}}{0}{0}{\omega^{2} \wh{Z}_{v_t} + \omega^5 \wh{Z}_{v_t}^{-2}}{\omega^{-7} \wh{Z}_{v_t}}{\omega^{-4} \wh{Z}_{v_t} + \omega^{-1} \wh{Z}_{v_t}^{-2}}{0}
&
 = 
\smallmatthree{0}{0}{0}{0}{0}{\omega^{-1} \wh{Z}_{v_t}^{-2}}{0}{\omega^{-7} \wh{Z}_{v_t}^{-2}}{0}.
\end{array}
\end{align*}
Now, to each $(\varepsilon_1,\varepsilon_2)$-th entry, we multiply $(\til{\bf F}^{\rm out}_{+,\varepsilon})_{\varepsilon_1,\varepsilon_2}$. For convenience, if we define $\til{\bf F}^{\rm out}_+ := \sum_{\alpha=1}^3 \til{\bf F}^{\rm out}_{+,\varepsilon}$, we have $(\til{\bf F}^{\rm out}_{+,\varepsilon})_{\varepsilon_1,\varepsilon_2} = (\til{\bf F}^{\rm out}_+)_{\varepsilon_1,\varepsilon_2}$. Note
\begin{align}
\label{eq:bf_F_out_tilde}
\til{\bf F}^{\rm out}_+ = \til{\bf F}^{\rm out}_{+,1} + \til{\bf F}^{\rm out}_{+,2} + \til{\bf F}^{\rm out}_{+,3} = \smallmatthree{0}{\omega^3}{1}{-\omega^9}{0}{\omega^3}{-\omega^{12}}{-\omega^9}{0}.
\end{align}
It is now easy to check that $(\wh{\bf M}^{\rm lt,o,l,col}_{+,\varepsilon_3})_{\varepsilon_1,\varepsilon_2} = (\til{\bf F}^{\rm out}_+)_{\varepsilon_1,\varepsilon_2} (\wh{\bf M}^{\rm r}_{+,\varepsilon_3})_{\varepsilon_1,\varepsilon_2}$ holds, as desired. This finishes the proof of $\til{\rm Tr}^\omega_{\wh{\Delta}}(W',s') =\til{\rm Tr}^\omega_{\wh{\Delta}}(W,s)$ for the left case of Fig.\ref{fig:isotopy_move3}.

\vs

Consider now the right case of Fig.\ref{fig:isotopy_move3} (with $y_1\succ y_2$, $w_1 \succ w_2$). Similarly as in the left case of Fig.\ref{fig:isotopy_move3} (involving eq.\eqref{eq:isotopy_invariance_proof_fork3}, eq.\eqref{eq:isotopy_invariance_proof_fork5}, Prop.\ref{prop:biangle_SL3_quantum_trace}(BT2-1), Thm.\ref{thm:SL3_quantum_trace_map}\redfix{(QT2-3))}, if we denote the state values of $s$ and $s'$ at $y_1,y_2,z$ by $\varepsilon_1,\varepsilon_2,\varepsilon_3$, we have
\begin{align*}
\til{\rm Tr}^\omega_{\wh{\Delta}}(W',s') 
= (\wh{\bf F}^{\rm out}_{+,\varepsilon})_{\varepsilon_1,\varepsilon_2} ( \wh{\bf M}^{\rm in}_{t,2} \wh{\bf M}^{\rm left}(\wh{Z}_{v_t}) \wh{\bf M}^{\rm out}_{t,3} )_{\varepsilon,\varepsilon_3}
= (\wh{\bf F}^{\rm out}_{+,\varepsilon})_{\varepsilon_1,\varepsilon_2} ( \wh{\bf M}^{\rm in}_{t,2} )_{\varepsilon,\varepsilon} \,
(\wh{\bf M}^{\rm left}(\wh{Z}_{v_t}) )_{\varepsilon,\varepsilon_3} \,
(\wh{\bf M}^{\rm out}_{t,3} )_{\varepsilon_3,\varepsilon_3},
\end{align*}
if $\{r_1(\varepsilon),r_2(\varepsilon)\} = \{\varepsilon_1,\varepsilon_2\}$, while $\til{\rm Tr}^\omega_{\wh{\Delta}}(W',s')=0$ if there is no such $\varepsilon$. Arguments for $(W,s)$ go similarly, and \redfix{yield} the following, if we denote by $\varepsilon_4,\varepsilon_5$ the juncture-state values for $w_1,w_2$:
\begin{align*}
\til{\rm Tr}^\omega_{\wh{\Delta}}(W,s)
& = {\textstyle \sum}_{\varepsilon_4,\varepsilon_5} (\wh{\bf F}^{\rm out}_{+,\varepsilon_3})_{\varepsilon_4,\varepsilon_5} \, (\wh{\bf M}^{\rm in}_{t,3} \wh{\bf M}^{\rm right}(\wh{Z}_{v_t}) \wh{\bf M}^{\rm out}_{t,2})_{\varepsilon_4, \varepsilon_1} (\wh{\bf M}^{\rm in}_{t,3} \wh{\bf M}^{\rm right}(\wh{Z}_{v_t}) \wh{\bf M}^{\rm out}_{t,2})_{\varepsilon_5, \varepsilon_2} \\
& = {\textstyle \sum}_{\varepsilon_4,\varepsilon_5} (\wh{\bf M}^{\rm out}_{t,2} \wh{\bf M}^{\rm right}_{\rm tran}(\wh{Z}_{v_t}) \wh{\bf M}^{\rm in}_{t,3})_{\varepsilon_1, \varepsilon_4} (\wh{\bf F}^{\rm out}_{+,\varepsilon_3})_{\varepsilon_4,\varepsilon_5} (\wh{\bf M}^{\rm in}_{t,3} \wh{\bf M}^{\rm right}(\wh{Z}_{v_t}) \wh{\bf M}^{\rm out}_{t,2})_{\varepsilon_5, \varepsilon_2} \quad (\because \mbox{Cor.\ref{cor:transpose_of_quantum_turn_matrices}}) \\
& = 
(\wh{\bf M}^{\rm out}_{t,2} \wh{\bf M}^{\rm right}_{\rm tran}(\wh{Z}_{v_t}) \underline{ \wh{\bf M}^{\rm in}_{t,3}
\wh{\bf F}^{\rm out}_{+,\varepsilon_3}
\wh{\bf M}^{\rm in}_{t,3} } \wh{\bf M}^{\rm right}(\wh{Z}_{v_t}) \wh{\bf M}^{\rm out}_{t,2})_{\varepsilon_1, \varepsilon_2}  \\
& =
(\wh{\bf M}^{\rm out}_{t,2} \wh{\bf M}^{\rm right}_{\rm tran}(\wh{Z}_{v_t}) \underline{  (\wh{\bf M}^{\rm out}_{t,3})_{\varepsilon_3,\varepsilon_3} 
\til{\bf F}^{\rm out}_{+,\varepsilon_3}
\wh{\bf M}^{\rm right}(\wh{Z}_{v_t})  }\wh{\bf M}^{\rm out}_{t,2})_{\varepsilon_1, \varepsilon_2} 
\quad (\because \mbox{Lem.\ref{lem:in_and_out_and_fork}}) \\
& = 
(\wh{\bf M}^{\rm out}_{t,2}
\underbrace{  \wh{\bf M}^{\rm right}_{\rm tran}(\wh{Z}_{v_t}) 
\til{\bf F}^{\rm out}_{+,\varepsilon_3}
\wh{\bf M}^{\rm right}(\omega^{-2g(\varepsilon_3)} \wh{Z}_{v_t}) }_{=: \wh{\bf M}^{\rm rt,o,r}_{+,\varepsilon_3} } (\wh{\bf M}^{\rm out}_{t,3})_{\varepsilon_3,\varepsilon_3} 
\wh{\bf M}^{\rm out}_{t,2})_{\varepsilon_1, \varepsilon_2}. 
\quad (\because \mbox{Lem.\ref{lem:triangle_variable_and_edge_matrices}})
\end{align*}
We compute the underbraced matrix product $\wh{\bf M}^{\rm rt,o,r}_{+,\varepsilon_3}$:
\begin{align*}
\hspace{0mm}
\begin{array}{lll}
\wh{\bf M}^{\rm rt,o,r}_{+,1}
&
\wh{\bf M}^{\rm rt,o,r}_{+,2}
&
\wh{\bf M}^{\rm rt,o,r}_{+,3}
\\
=  \smallmatthree{0}{\omega^2 \wh{Z}_{v_t}^2}{0}{-\omega^{14} \wh{Z}_{v_t}^2}{0}{0}{0}{0}{0}
&
=  \smallmatthree{0}{\omega^{-4} \wh{Z}_{v_t}^2 + \omega^{11} \wh{Z}_{v_t}^{-1}}{\omega^5 \wh{Z}_{v_t}^{-1}}{-\omega^8 \wh{Z}_{v_t}^2 - \omega^5 \wh{Z}_{v_t}^{-1}}{0}{0}{-\omega^{11} \wh{Z}_{v_t}^{-1}}{0}{0}
&
=  \smallmatthree{0}{\omega^{-1} \wh{Z}_{v_t}^{-1}}{\omega^{-7} \wh{Z}_{v_t}^{-1}}{- \omega^{11} \wh{Z}_{v_t}^{-1}}{0}{\omega^{-1} \wh{Z}_{v_t}^{-1}}{- \omega^{17} \wh{Z}_{v_t}^{-1}}{-\omega^{11} \wh{Z}_{v_t}^{-1}}{0}
\end{array}
\end{align*}
Then
\begin{align*}
\til{\rm Tr}^\omega_{\wh{\Delta}}(W,s)
& = ( \wh{\bf M}^{\rm out}_{t,2} \,
\wh{\bf M}^{\rm rt,o,r}_{+,\varepsilon_3} \,
(\wh{\bf M}^{\rm out}_{t,3})_{\varepsilon_3,\varepsilon_3} \,
\wh{\bf M}^{\rm out}_{t,2} )_{\varepsilon_1,\varepsilon_2} \\
& = ( \wh{\bf M}^{\rm out}_{t,2} )_{\varepsilon_1,\varepsilon_1} \,
( \wh{\bf M}^{\rm rt,o,r}_{+,\varepsilon_3} )_{\varepsilon_1,\varepsilon_2} \,
\underline{ (\wh{\bf M}^{\rm out}_{t,3})_{\varepsilon_3,\varepsilon_3} \,
(\wh{\bf M}^{\rm out}_{t,2} )_{\varepsilon_2,\varepsilon_2} } \quad (\because\mbox{$\wh{\bf M}^{\rm out}_{t,\alpha}$ are diagonal}) \\
& = \omega^{-2g(\varepsilon_2+1)g(\varepsilon_3-1)} ( \wh{\bf M}^{\rm out}_{t,2} )_{\varepsilon_1,\varepsilon_1} \,
( \wh{\bf M}^{\rm rt,o,r}_{+,\varepsilon_3} )_{\varepsilon_1,\varepsilon_2} \,
(\wh{\bf M}^{\rm out}_{t,2} )_{\varepsilon_2,\varepsilon_2} 
(\wh{\bf M}^{\rm out}_{t,3})_{\varepsilon_3,\varepsilon_3} \,
 \quad (\because \mbox{Lem.\ref{lem:edge_matrices_commutation}})
\end{align*}
By inspection, if $\varepsilon_1=\varepsilon_2$, then $( \wh{\bf M}^{\rm rt,o,r}_{+,\varepsilon_3} )_{\varepsilon_1,\varepsilon_2}=0$, hence $\til{\rm Tr}^\omega_{\wh{\Delta}}(W,s)=0$.
Using Lem.\ref{lem:triangle_variable_and_edge_matrices}, we have
$$
(\wh{\bf M}^{\rm rt,o,r}_{+,\varepsilon_3})_{\varepsilon_1,\varepsilon_2} \,
(\wh{\bf M}^{\rm out}_{t,2})_{\varepsilon_2,\varepsilon_2}
= (\wh{\bf M}^{\rm out}_{t,2})_{\varepsilon_2,\varepsilon_2} (\wh{\bf M}^{\rm rt,o,r, col}_{+,\varepsilon_3})_{\varepsilon_1,\varepsilon_2} \,
$$
where $\wh{\bf M}^{\rm rt,o,r, col}_{+,\varepsilon_3}$ is obtained from $\wh{\bf M}^{\rm rt,o,r}_{+,\varepsilon_3}$ by replacing $\wh{Z}_{v_t}$ appearing in each $\varepsilon_2$-th column with $\omega^{2g(\varepsilon_2)} \wh{Z}_{v_t}$:
\begin{align*}
\hspace{0mm}
\begin{array}{lll}
\wh{\bf M}^{\rm rt,o,r,col}_{+,1}
&
\wh{\bf M}^{\rm rt,o,r,col}_{+,2}
&
\wh{\bf M}^{\rm rt,o,r,col}_{+,3}
\\
= \smallmatthree{0}{\omega^{10} \wh{Z}_{v_t}^2}{0}{-\omega^{10} \wh{Z}_{v_t}^2}{0}{0}{0}{0}{0}
&
= \smallmatthree{0}{\omega^{4} \wh{Z}_{v_t}^2 + \omega^{7} \wh{Z}_{v_t}^{-1}}{\omega^7 \wh{Z}_{v_t}^{-1}}{-\omega^4 \wh{Z}_{v_t}^2 - \omega^7 \wh{Z}_{v_t}^{-1}}{0}{0}{-\omega^{13} \wh{Z}_{v_t}^{-1}}{0}{0}
&
= \smallmatthree{0}{\omega^{-5} \wh{Z}_{v_t}^{-1}}{\omega^{-5} \wh{Z}_{v_t}^{-1}}{- \omega^{13} \wh{Z}_{v_t}^{-1}}{0}{\omega \wh{Z}_{v_t}^{-1}}{- \omega^{19} \wh{Z}_{v_t}^{-1}}{-\omega^{7} \wh{Z}_{v_t}^{-1}}{0}
\end{array}
\end{align*}
For $\varepsilon_1\neq \varepsilon_2$, if $\{r_1(\varepsilon),r_2(\varepsilon)\} = \{\varepsilon_1,\varepsilon_2\}$, we have
\begin{align*}
\til{\rm Tr}^\omega_{\wh{\Delta}}(W,s)
& = \omega^{-2g(\varepsilon_2+1)g(\varepsilon_3-1)} ( \wh{\bf M}^{\rm out}_{t,2} )_{\varepsilon_1,\varepsilon_1} \,
(\wh{\bf M}^{\rm out}_{t,2} )_{\varepsilon_2,\varepsilon_2}  \,( \wh{\bf M}^{\rm rt,o,r,col}_{+,\varepsilon_3} )_{\varepsilon_1,\varepsilon_2} \,
(\wh{\bf M}^{\rm out}_{t,3})_{\varepsilon_3,\varepsilon_3} \\
& = 
\omega^{-2g(\varepsilon_2+1)g(\varepsilon_3-1)} \omega^{\frac{3}{2}p(\varepsilon) {\rm sgn}(\varepsilon_1-\varepsilon_2)} ( \wh{\bf M}^{\rm in}_{t,2} )_{\varepsilon,\varepsilon}
( \wh{\bf M}^{\rm rt,o,r,col}_{+,\varepsilon_3} )_{\varepsilon_1,\varepsilon_2} \,
(\wh{\bf M}^{\rm out}_{t,3})_{\varepsilon_3,\varepsilon_3},
\qquad (\because\mbox{Lem.\ref{lem:edge_matrix_inversion_formula}})
\end{align*}
hence now it suffices to show
\begin{align*}
\omega^{-2g(\varepsilon_2+1)g(\varepsilon_3-1)} \,
\omega^{\frac{3}{2}p(\varepsilon) {\rm sgn}(\varepsilon_1-\varepsilon_2)} ( \wh{\bf M}^{\rm rt,o,r,col}_{+,\varepsilon_3} )_{\varepsilon_1,\varepsilon_2}
= (\wh{\bf F}^{\rm out}_{+,\varepsilon})_{\varepsilon_1,\varepsilon_2} \,
(\wh{\bf M}^{\rm left}(\wh{Z}_{v_t}) )_{\varepsilon,\varepsilon_3},
\end{align*}
or equivalently (from eq.\eqref{eq:fork_matrix_tilde}),
\begin{align*}
\omega^{-2g(\varepsilon_2+1)g(\varepsilon_3-1)} \, ( \wh{\bf M}^{\rm rt,o,r,col}_{+,\varepsilon_3} )_{\varepsilon_1,\varepsilon_2}
= (\til{\bf F}^{\rm out}_{+,\varepsilon})_{\varepsilon_1,\varepsilon_2} \,
(\wh{\bf M}^{\rm left}(\wh{Z}_{v_t}) )_{\varepsilon,\varepsilon_3},
\end{align*}
For each $\varepsilon_3 \in \{1,2,3\}$, define $\wh{\bf M}^{\rm l}_{+,\varepsilon_3}$ as the $3 \times 3$ matrix whose $(\varepsilon_1,\varepsilon_2)$-th entry equals zero if $\varepsilon_1=\varepsilon_2$, and equals $\omega^{2g(\varepsilon_2+1)g(\varepsilon_3-1)} (\wh{\bf M}^{\rm left}(\wh{Z}_{v_t}) )_{\varepsilon,\varepsilon_3}$ if $\varepsilon_1\neq \varepsilon_2$, where $\varepsilon$ is determined by $\{r_1(\varepsilon),r_2(\varepsilon)\} = \{\varepsilon_1,\varepsilon_2\}$. We compute them:
\begin{align*}
\begin{array}{lll}
\wh{\bf M}^{\rm l}_{+,1}
& 
\wh{\bf M}^{\rm l}_{+,2}
&
\wh{\bf M}^{\rm l}_{+,3}
\\
= \smallmatthree{0}{\omega^7 \wh{Z}_{v_t}^2}{0}{\omega \wh{Z}_{v_t}^2}{0}{0}{0}{0}{0}
&
= \smallmatthree{0}{\omega \wh{Z}_{v_t}^2 + \omega^4 \wh{Z}_{v_t}^{-1}}{\omega^7 \wh{Z}_{v_t}^{-1}}{\omega^{-5} \wh{Z}_{v_t}^2 + \omega^{-2} \wh{Z}_{v_t}^{-1}}{0}{0}{\omega \wh{Z}_{v_t}^{-1}}{0}{0}
&
= \smallmatthree{0}{\omega^{-8} \wh{Z}_{v_t}^{-1}}{\omega^{-5} \wh{Z}_{v_t}^{-1}}{\omega^{4} \wh{Z}_{v_t}^{-1}}{0}{\omega^{-2} \wh{Z}_{v_t}^{-1}}{\omega^{7} \wh{Z}_{v_t}^{-1}}{\omega^{-2}\wh{Z}_{v_t}^{-1}}{0}
\end{array}
\end{align*}
Multiplying each $(\varepsilon_1,\varepsilon_2)$-th entry by $(\til{\bf F}_+)_{\varepsilon_1,\varepsilon_2}$ (eq.\eqref{eq:bf_F_out_tilde}), one can easily check that $(\wh{\bf M}^{\rm rt,o,r,col}_{+,\varepsilon_3})_{\varepsilon_1,\varepsilon_2} = (\til{\bf F}^{\rm out}_+)_{\varepsilon_1,\varepsilon_2} (\wh{\bf M}^{\rm l}_{+,\varepsilon_3})_{\varepsilon_1,\varepsilon_2}$ holds, as desired. This finishes the proof of $\til{\rm Tr}^\omega_{\wh{\Delta}}(W',s') =\til{\rm Tr}^\omega_{\wh{\Delta}}(W,s)$ for the left case of Fig.\ref{fig:isotopy_move3}.

\vs

Now, take the left case of Fig.\ref{fig:isotopy_move3} (with $x_1\succ x_2$ and $w_1\succ w_2$), with orientations reversed. Using similar arguments and notations, we get
\begin{align}
\nonumber
\til{\rm Tr}^\omega_{\wh{\Delta}}(W',s')
= (\wh{\bf F}^{\rm in}_{+,\varepsilon})_{\varepsilon_1,\varepsilon_2} ( \wh{\bf M}^{\rm in}_{t,3} \wh{\bf M}^{\rm left}(\wh{Z}_{v_t}) \wh{\bf M}^{\rm out}_{t,1} )_{\varepsilon_3,\varepsilon} 
& \stackrel{{\rm Cor.}\ref{cor:transpose_of_quantum_turn_matrices}}{=}
(\wh{\bf F}^{\rm in}_{+,\varepsilon})_{\varepsilon_1,\varepsilon_2} ( \wh{\bf M}^{\rm out}_{t,1}  \wh{\bf M}^{\rm left}_{\rm tran} (\wh{Z}_{v_t}) \wh{\bf M}^{\rm in}_{t,3} )_{\varepsilon,\varepsilon_3}  \\
\nonumber
& = (\wh{\bf F}^{\rm in}_{+,\varepsilon})_{\varepsilon_1,\varepsilon_2}  \, (\wh{\bf M}^{\rm out}_{t,1} )_{\varepsilon,\varepsilon} \, (\wh{\bf M}^{\rm left}_{\rm tran}(\wh{Z}_{v_t}) )_{\varepsilon,\varepsilon_3} \, ( \wh{\bf M}^{\rm in}_{t,3} )_{\varepsilon_3,\varepsilon_3},
\end{align}
while
\begin{align*}
\til{\rm Tr}^\omega_{\wh{\Delta}}(W,s)
& = {\textstyle \sum}_{\varepsilon_4,\varepsilon_5} (\wh{\bf F}^{\rm in}_{+,\varepsilon_3})_{\varepsilon_4,\varepsilon_5} \, 
(\wh{\bf M}^{\rm in}_{t,1} \wh{\bf M}^{\rm right}(\wh{Z}_{v_t}) \wh{\bf M}^{\rm out}_{t,3})_{\varepsilon_1, \varepsilon_4} 
(\wh{\bf M}^{\rm in}_{t,1} \wh{\bf M}^{\rm right}(\wh{Z}_{v_t}) \wh{\bf M}^{\rm out}_{t,3})_{\varepsilon_2, \varepsilon_5} \\
& = {\textstyle \sum}_{\varepsilon_4,\varepsilon_5} 
(\wh{\bf M}^{\rm in}_{t,1} \wh{\bf M}^{\rm right}(\wh{Z}_{v_t}) \wh{\bf M}^{\rm out}_{t,3})_{\varepsilon_1, \varepsilon_4}
(\wh{\bf F}^{\rm in}_{+,\varepsilon_3})_{\varepsilon_4,\varepsilon_5} 
(\wh{\bf M}^{\rm out}_{t,3} \wh{\bf M}^{\rm right}_{\rm tran}(\wh{Z}_{v_t}) \wh{\bf M}^{\rm in}_{t,1})_{\varepsilon_5, \varepsilon_2} \quad (\because \mbox{Cor.\ref{cor:transpose_of_quantum_turn_matrices}}) \\
& = 
(\wh{\bf M}^{\rm in}_{t,1} \wh{\bf M}^{\rm right}(\wh{Z}_{v_t}) 
\underline{ \wh{\bf M}^{\rm out}_{t,3}
\wh{\bf F}^{\rm in}_{+,\varepsilon_3}
\wh{\bf M}^{\rm out}_{t,3} }
\wh{\bf M}^{\rm right}_{\rm tran}(\wh{Z}_{v_t}) \wh{\bf M}^{\rm in}_{t,1} )_{\varepsilon_1, \varepsilon_2}  \\
& =
(\wh{\bf M}^{\rm in}_{t,1} \wh{\bf M}^{\rm right}(\wh{Z}_{v_t}) 
\underline{ (\wh{\bf M}^{\rm in}_{t,3})_{\varepsilon_3,\varepsilon_3}
\til{\bf F}^{\rm in}_{+,\varepsilon_3} 
\wh{\bf M}^{\rm right}_{\rm tran}(\wh{Z}_{v_t}) } \wh{\bf M}^{\rm in}_{t,1} )_{\varepsilon_1, \varepsilon_2}  
\quad (\because \mbox{Lem.\ref{lem:in_and_out_and_fork}}) \\
& = 
( \wh{\bf M}^{\rm in}_{t,1} 
\underbrace{ \wh{\bf M}^{\rm right}(\wh{Z}_{v_t}) 
\til{\bf F}^{\rm in}_{+,\varepsilon_3} 
\wh{\bf M}^{\rm right}_{\rm tran}(\omega^{2g(\varepsilon_3)} \wh{Z}_{v_t}) }_{=:\wh{\bf M}^{\rm r,i,rt}_{+,\varepsilon_3}}
(\wh{\bf M}^{\rm in}_{t,3})_{\varepsilon_3,\varepsilon_3}
\wh{\bf M}^{\rm in}_{t,1}  )_{\varepsilon_1, \varepsilon_2}. 
\quad (\because \mbox{Lem.\ref{lem:triangle_variable_and_edge_matrices}})
\end{align*}
Then
\begin{align*}
\til{\rm Tr}^\omega_{\wh{\Delta}}(W,s)
& = ( \wh{\bf M}^{\rm in}_{t,1} 
\wh{\bf M}^{\rm r,i,rt}_{+,\varepsilon_3}
(\wh{\bf M}^{\rm in}_{t,3})_{\varepsilon_3,\varepsilon_3}
\wh{\bf M}^{\rm in}_{t,1}  )_{\varepsilon_1, \varepsilon_2} \\
& = ( \wh{\bf M}^{\rm in}_{t,1} )_{\varepsilon_1,\varepsilon_1}
(\wh{\bf M}^{\rm r,i,rt}_{+,\varepsilon_3} )_{\varepsilon_1,\varepsilon_2}
\underline{ (\wh{\bf M}^{\rm in}_{t,3})_{\varepsilon_3,\varepsilon_3}
(\wh{\bf M}^{\rm in}_{t,1}  )_{\varepsilon_2, \varepsilon_2} } \\
& = \omega^{2g(\varepsilon_3-1)g(\varepsilon_2+1)} ( \wh{\bf M}^{\rm in}_{t,1} )_{\varepsilon_1,\varepsilon_1}
\underline{ (\wh{\bf M}^{\rm r,i,rt}_{+,\varepsilon_3} )_{\varepsilon_1,\varepsilon_2}
(\wh{\bf M}^{\rm in}_{t,1}  )_{\varepsilon_2, \varepsilon_2} }
(\wh{\bf M}^{\rm in}_{t,3})_{\varepsilon_3,\varepsilon_3} \quad (\because \mbox{Lem.\ref{lem:edge_matrices_commutation}}) \\
& = \omega^{2g(\varepsilon_3-1)g(\varepsilon_2+1)} ( \wh{\bf M}^{\rm in}_{t,1} )_{\varepsilon_1,\varepsilon_1}
(\wh{\bf M}^{\rm in}_{t,1}  )_{\varepsilon_2, \varepsilon_2}
(\wh{\bf M}^{\rm r,i,rt,col}_{+,\varepsilon_3} )_{\varepsilon_1,\varepsilon_2}
(\wh{\bf M}^{\rm in}_{t,3})_{\varepsilon_3,\varepsilon_3}, \quad (\because \mbox{Lem.\ref{lem:triangle_variable_and_edge_matrices}})
\end{align*}
where $\wh{\bf M}^{\rm r,i,rt, col}_{+,\varepsilon_3}$ is obtained from $\wh{\bf M}^{\rm r,i,rt}_{+,\varepsilon_3}$ by replacing $\wh{Z}_{v_t}$ appearing in each $\varepsilon_2$-th column with $\omega^{-2g(\varepsilon_2)} \wh{Z}_{v_t}$. By Lem.\ref{lem:edge_matrix_inversion_formula} and eq.\eqref{eq:fork_matrix_tilde}, it then suffices to show
\begin{align*}
\omega^{2g(\varepsilon_3-1)g(\varepsilon_2+1)}
(\wh{\bf M}^{\rm r,i,rt,col}_{+,\varepsilon_3} )_{\varepsilon_1,\varepsilon_2}
= (\til{\bf F}^{\rm in}_{+,\varepsilon})_{\varepsilon_1,\varepsilon_2} (\wh{\bf M}^{\rm left}_{\rm tran}(\wh{Z}_{v_t}) )_{\varepsilon,\varepsilon_3}
\end{align*}
for each $\varepsilon_3$ and whenever $\{r_1(\varepsilon),r_2(\varepsilon)\} = \{\varepsilon_1,\varepsilon_2\}$. The actual computation is similar as before, and is left as an exercise.

\vs

Take the right case of Fig.\ref{fig:isotopy_move3} (with $y_1\succ y_2$ and $w_1\succ w_2$), with orientations reversed. Using similar arguments and notations, we get
\begin{align}
\nonumber
\til{\rm Tr}^\omega_{\wh{\Delta}}(W',s')
= (\wh{\bf F}^{\rm in}_{+,\varepsilon})_{\varepsilon_1,\varepsilon_2} ( \wh{\bf M}^{\rm in}_{t,3} \wh{\bf M}^{\rm right}(\wh{Z}_{v_t}) \wh{\bf M}^{\rm out}_{t,2} )_{\varepsilon_3,\varepsilon}
& \stackrel{{\rm Cor.}\ref{cor:transpose_of_quantum_turn_matrices}}{=} 
(\wh{\bf F}^{\rm in}_{+,\varepsilon})_{\varepsilon_1,\varepsilon_2} ( \wh{\bf M}^{\rm out}_{t,2} \wh{\bf M}^{\rm right}_{\rm tran}(\wh{Z}_{v_t}) \wh{\bf M}^{\rm in}_{t,3} )_{\varepsilon,\varepsilon_3}
\\
\nonumber
&= (\wh{\bf F}^{\rm in}_{+,\varepsilon})_{\varepsilon_1,\varepsilon_2}  \, (\wh{\bf M}^{\rm out}_{t,2} )_{\varepsilon,\varepsilon} (\wh{\bf M}^{\rm right}_{\rm tran}(\wh{Z}_{v_t}) )_{\varepsilon,\varepsilon_3} \, ( \wh{\bf M}^{\rm in}_{t,3} )_{\varepsilon_3,\varepsilon_3},
\end{align}
while
\begin{align*}
\til{\rm Tr}^\omega_{\wh{\Delta}}(W,s)
& = {\textstyle \sum}_{\varepsilon_4,\varepsilon_5} (\wh{\bf F}^{\rm in}_{+,\varepsilon_3})_{\varepsilon_4,\varepsilon_5} \, 
(\wh{\bf M}^{\rm in}_{t,2} \wh{\bf M}^{\rm left}(\wh{Z}_{v_t}) \wh{\bf M}^{\rm out}_{t,3})_{\varepsilon_1, \varepsilon_4} 
(\wh{\bf M}^{\rm in}_{t,2} \wh{\bf M}^{\rm left}(\wh{Z}_{v_t}) \wh{\bf M}^{\rm out}_{t,3})_{\varepsilon_2, \varepsilon_5} \\
& = {\textstyle \sum}_{\varepsilon_4,\varepsilon_5} 
(\wh{\bf M}^{\rm in}_{t,2} \wh{\bf M}^{\rm left}(\wh{Z}_{v_t}) \wh{\bf M}^{\rm out}_{t,3})_{\varepsilon_1, \varepsilon_4}
(\wh{\bf F}^{\rm in}_{+,\varepsilon_3})_{\varepsilon_4,\varepsilon_5} 
(\wh{\bf M}^{\rm out}_{t,3} \wh{\bf M}^{\rm left}_{\rm tran}(\wh{Z}_{v_t}) \wh{\bf M}^{\rm in}_{t,2})_{\varepsilon_5, \varepsilon_2} \quad (\because \mbox{Cor.\ref{cor:transpose_of_quantum_turn_matrices}}) \\
& = 
(\wh{\bf M}^{\rm in}_{t,2} \wh{\bf M}^{\rm left}(\wh{Z}_{v_t}) 
\underline{ \wh{\bf M}^{\rm out}_{t,3}
\wh{\bf F}^{\rm in}_{+,\varepsilon_3}
\wh{\bf M}^{\rm out}_{t,3} }
\wh{\bf M}^{\rm left}_{\rm tran}(\wh{Z}_{v_t}) \wh{\bf M}^{\rm in}_{t,2} )_{\varepsilon_1, \varepsilon_2}  \\
& =
(\wh{\bf M}^{\rm in}_{t,2} \wh{\bf M}^{\rm left}(\wh{Z}_{v_t}) 
\underline{ (\wh{\bf M}^{\rm in}_{t,3})_{\varepsilon_3,\varepsilon_3}
\til{\bf F}^{\rm in}_{+,\varepsilon_3} 
\wh{\bf M}^{\rm left}_{\rm tran}(\wh{Z}_{v_t}) } \wh{\bf M}^{\rm in}_{t,2} )_{\varepsilon_1, \varepsilon_2}  
\quad (\because \mbox{Lem.\ref{lem:in_and_out_and_fork}}) \\
& = 
( \wh{\bf M}^{\rm in}_{t,2} 
\underbrace{ \wh{\bf M}^{\rm left}(\wh{Z}_{v_t}) 
\til{\bf F}^{\rm in}_{+,\varepsilon_3} 
\wh{\bf M}^{\rm left}_{\rm tran}(\omega^{2g(\varepsilon_3)} \wh{Z}_{v_t}) }_{=:\wh{\bf M}^{\rm l,i,lt}_{+,\varepsilon_3}}
(\wh{\bf M}^{\rm in}_{t,3})_{\varepsilon_3,\varepsilon_3}
\wh{\bf M}^{\rm in}_{t,2}  )_{\varepsilon_1, \varepsilon_2}. 
\quad (\because \mbox{Lem.\ref{lem:triangle_variable_and_edge_matrices}})
\end{align*}
Then
\begin{align*}
\til{\rm Tr}^\omega_{\wh{\Delta}}(W,s)
& = ( \wh{\bf M}^{\rm in}_{t,2} 
\wh{\bf M}^{\rm l,i,lt}_{+,\varepsilon_3}
(\wh{\bf M}^{\rm in}_{t,3})_{\varepsilon_3,\varepsilon_3}
\wh{\bf M}^{\rm in}_{t,2}  )_{\varepsilon_1, \varepsilon_2} \\
& = ( \wh{\bf M}^{\rm in}_{t,2} )_{\varepsilon_1,\varepsilon_1}
(\wh{\bf M}^{\rm l,i,lt}_{+,\varepsilon_3} )_{\varepsilon_1,\varepsilon_2}
\underline{ (\wh{\bf M}^{\rm in}_{t,3})_{\varepsilon_3,\varepsilon_3}
(\wh{\bf M}^{\rm in}_{t,2}  )_{\varepsilon_2, \varepsilon_2} } \\
& = \omega^{-2g(\varepsilon_2-1)g(\varepsilon_3+1)} ( \wh{\bf M}^{\rm in}_{t,2} )_{\varepsilon_1,\varepsilon_1}
\underline{ (\wh{\bf M}^{\rm l,i,lt}_{+,\varepsilon_3} )_{\varepsilon_1,\varepsilon_2}
(\wh{\bf M}^{\rm in}_{t,2}  )_{\varepsilon_2, \varepsilon_2} }
(\wh{\bf M}^{\rm in}_{t,3})_{\varepsilon_3,\varepsilon_3} \quad (\because \mbox{Lem.\ref{lem:edge_matrices_commutation}}) \\
& = \omega^{-2g(\varepsilon_2-1)g(\varepsilon_3+1)} ( \wh{\bf M}^{\rm in}_{t,2} )_{\varepsilon_1,\varepsilon_1}
(\wh{\bf M}^{\rm in}_{t,2}  )_{\varepsilon_2, \varepsilon_2}
(\wh{\bf M}^{\rm l,i,lt,col}_{+,\varepsilon_3} )_{\varepsilon_1,\varepsilon_2}
(\wh{\bf M}^{\rm in}_{t,3})_{\varepsilon_3,\varepsilon_3}, \quad (\because \mbox{Lem.\ref{lem:triangle_variable_and_edge_matrices}})
\end{align*}
where $\wh{\bf M}^{\rm l,i,lt, col}_{+,\varepsilon_3}$ is obtained from $\wh{\bf M}^{\rm l,i,lt}_{+,\varepsilon_3}$ by replacing $\wh{Z}_{v_t}$ appearing in each $\varepsilon_2$-th column with $\omega^{-2g(\varepsilon_2)} \wh{Z}_{v_t}$. By Lem.\ref{lem:edge_matrix_inversion_formula} and eq.\eqref{eq:fork_matrix_tilde}, it then suffices to show
\begin{align*}
\omega^{-2g(\varepsilon_2-1)g(\varepsilon_3+1)}
(\wh{\bf M}^{\rm l,i,lt,col}_{+,\varepsilon_3} )_{\varepsilon_1,\varepsilon_2}
= (\til{\bf F}^{\rm in}_{+,\varepsilon})_{\varepsilon_1,\varepsilon_2} (\wh{\bf M}^{\rm right}_{\rm tran}(\wh{Z}_{v_t}) )_{\varepsilon,\varepsilon_3}
\end{align*}
for each $\varepsilon_3$ and whenever $\{r_1(\varepsilon),r_2(\varepsilon)\} = \{\varepsilon_1,\varepsilon_2\}$. The actual computation is similar as before, and is left as an exercise. 

\vs

Finally, the isotopy invariance for a move in Fig.\ref{fig:isotopy_move8} can be obtained as a consequence of that for a move in Fig.\ref{fig:isotopy_move2}. Split the biangle of $B$ in Fig.\ref{fig:isotopy_move8} into two (by cutting it along an ideal arc connecting the two marked points), and pull the 3-valent vertex into the biangle adjacent to the triangle using the move in Fig.\ref{fig:isotopy_move2}; then we use isotopy invariance of ${\rm Tr}^\omega_B$. One can also directly prove it, without resorting to the isotopy invariance of ${\rm Tr}^\omega_B$ (i.e. to Prop.\ref{prop:biangle_SL3_quantum_trace}). \qed \quad {\it End of proof of Prop.\ref{prop:elementary_isotopy_invariance_3-way}.}

\vs

The remaining cases, Prop.\ref{prop:elementary_isotopy_invariance_elevation_preserving_arcs} and Prop.\ref{prop:elementary_isotopy_invariance_elevation_change}, are objects of study in \cite{Douglas} \cite{Douglas21}, where, as said in \cite{Douglas21}, \redfix{we believe that these two propositions are} essentially proved in \cite[\redfix{Thms. 2.5, 2.6, 3.1}]{CS}\redfix{; although \cite{CS} is written in a slightly different language, it should be straightforward to verify this, but we do not do so here}. \redfix{Note that} some moves are checked in \cite{Douglas} \cite{Douglas21} with the aid of a computer calculation. We can use the results of \cite{Douglas} \cite{Douglas21}, thanks to Remarks \ref{rem:Douglas_triangle} and \ref{rem:Douglas_biangle}. In fact, one needs to be careful, as the conventions of \cite{Douglas} \cite{Douglas21} are a bit different from ours; namely, the order of the superposition product in the definition of ${\rm SL}_3$-skein algebras implicitly used in \cite{Douglas} \cite{Douglas21} is opposite to ours. So, to really match with \cite{Douglas} \cite{Douglas21}, one should apply the elevation-reversing map ${\bf r}$ of Lem.\ref{lem:elevation-reversing_map}. The effect of applying ${\bf r}$ on the biangle ${\rm SL}_3$ quantum trace can be expressed as the following biangle version of the equivariance statement, Prop.\ref{prop:elevation_reversing_and_star-structure}:
\begin{lemma}[equivariance under elevation-reversing and $*$-map for biangles]
\label{lem:elevation_reversing_and_star-structure_for_biangles}
For a biangle $B$,
\begin{align}
\label{eq:biangle_equivariance}
{\rm Tr}^\omega_B \circ {\bf r} = * \circ {\rm Tr}^\omega_B
\end{align}
holds, where ${\bf r}$ is as in Lem.\ref{lem:elevation-reversing_map}, and $*$ on the right hand side is $* : \mathbb{Z}[\omega^{\pm 1/2}] \to \mathbb{Z}[\omega^{\pm 1/2}]$ is the ring (anti-)homomorphism sending $\omega^{\pm 1/2}$ to $\omega^{\mp 1/2}$.
\end{lemma}
Indeed, by the item (BT1) of Prop.\ref{prop:biangle_SL3_quantum_trace}, it suffices to show eq.\eqref{eq:biangle_equivariance} applied to the cases in (BT2), for which one can easily verify eq.\eqref{eq:biangle_equivariance} by inspection; we note that it is our choice of  isomorphism in eq.\eqref{eq:isomorphism_from_ours_to_Higgins} and the corresponding boundary relations that made this equivariance to hold in biangles. Anyways, by observing that in \cite{Douglas} \cite{Douglas21} the quiver $Q_\Delta$ is drawn with the orientations opposite to ours, so that $q$ and $\omega$ of \cite{Douglas} \cite{Douglas21} correspond to our $q^{-1}$ and $\omega^{-1}$, one can finally match our setting with that of \cite{Douglas} \cite{Douglas21}. Daniel Douglas informed us that he also checked that the other relations hold, using a computer. For completeness, here we present how to show Prop.\ref{prop:elementary_isotopy_invariance_elevation_preserving_arcs} by hand. Prop.\ref{prop:elementary_isotopy_invariance_elevation_change} can be directly checked, as in \cite{Douglas} \cite{Douglas21}, and we expect that it can also be proved using Prop.\ref{prop:elementary_isotopy_invariance_3-way} (perhaps together with Fig.\ref{fig:isotopy_move1} with $w' \prec w''$ which can be proved in a similar manner as for $w' \succ w''$) and Prop.\ref{prop:elementary_isotopy_invariance_elevation_preserving_arcs}. Note that, for the purpose of the ${\rm SL}_3$ classical trace, we do not need Prop.\ref{prop:elementary_isotopy_invariance_elevation_change}.

\vs

{\it Proof of Prop.\ref{prop:elementary_isotopy_invariance_elevation_preserving_arcs}.} Take the left case of Fig.\ref{fig:isotopy_move4}, possibly with all orientations reversed. Denote by $\varepsilon_1,\varepsilon_2$ the state values of $s$ and $s'$ at $z_1,z_2$ respectively. Look at $W$ on the left (i.e. the first picture from the left in Fig.\ref{fig:isotopy_move4}). There is nothing in the triangle, so we have $\til{\rm Tr}^\omega_{\wh{\Delta}}(W,s) = {\rm Tr}^\omega_B([W\cap (B\times {\bf I}),(z_1\mapsto \varepsilon_1,z_2 \mapsto \varepsilon_2)])$. Look at $W'$, which is on the right (i.e. the second picture from the left in Fig.\ref{fig:isotopy_move4}). Suppose first that $z_1\succ z_2$, $w_1\succ w_2$, or that $z_1\prec z_2$, $w_1\prec w_2$. Then $W'\cap (B\times {\bf I})$ is the product of two edges connecting distinct sides, so in the sum in eq.\eqref{eq:til_Tr_W_s_prime}, by Prop.\ref{prop:biangle_SL3_quantum_trace}(BT2-1) the biangle factor `goes away', and we just have $\til{\rm Tr}^\omega_{\wh{\Delta}}(W',s') = \til{\rm Tr}^\omega_t(W' \cap (\wh{t}\times {\bf I}), (w_1 \mapsto \varepsilon_1, w_2 \mapsto \varepsilon_2))$. In view of the values stipulated by Thm.\ref{thm:SL3_quantum_trace_map}\redfix{(QT2-6)} and Prop.\ref{prop:biangle_SL3_quantum_trace}(BT2-2), we have $\til{\rm Tr}^\omega_{\wh{\Delta}}(W,s)=\til{\rm Tr}^\omega_{\wh{\Delta}}(W',s')$, as desired. Now suppose that $z_1\succ z_2$, $w_1 \prec w_2$, or that $z_1 \prec z_2$, $w_1 \succ w_2$. Split $B$ into two biangles $B_1,B_2$ by an ideal arc $e$ such that $e\times {\bf I}$ meets each of the two edges of $W' \cap (B\times {\bf I})$ exactly once, say at $u_1,u_2$, so that $z_1,u_1,w_1$ are on a same edge. Isotope $W'\cap (B\times {\bf I})$ in $B\times {\bf I}$ so that the elevation ordering for $u_1,u_2$ is same as that for $w_1,w_2$. Then, by the isotopy invariance just proved, one can pull the U-turn part of $W'$ living in $\wh{t} \times {\bf I}$ into the biangle $B_2$ that is adjacent to $\wh{t}$. Then one can use the isotopy invariance of ${\rm Tr}^\omega_B$, to prove the sought-for isotopy invariance for the move of the left picture of Fig.\ref{fig:isotopy_move4}. This proof depends on Cor.\ref{prop:biangle_SL3_quantum_trace}; one can also prove it directly, without using Cor.\ref{prop:biangle_SL3_quantum_trace}.

\vs

Take the right case of Fig.\ref{fig:isotopy_move4}, with $w_1\prec w_2$. Denote by $\varepsilon_1,\varepsilon_2$ the state values of $s$ and $s'$ at $x,y$ respectively. For $W$ (i.e. the third picture from the left of Fig.\ref{fig:isotopy_move4}), there is nothing in the biangle, so $\til{\rm Tr}^\omega_{\wh{\Delta}}(W,s) = \til{\rm Tr}^\omega_t(W\cap (\wh{t}\times {\bf I}),(x\mapsto \varepsilon_1,y\mapsto \varepsilon_2))$, which equals
$$
( \wh{\bf M}^{\rm in}_{t,2} \, \wh{\bf M}^{\rm right}(\wh{Z}_{v_t}) \, \wh{\bf M}^{\rm out}_{t,1} )_{\varepsilon_2,\varepsilon_1},
$$
in view of Thm.\ref{thm:SL3_quantum_trace_map}\redfix{(QT2-2)}. For $W'$ (i.e. the fourth picture from the left of Fig.\ref{fig:isotopy_move4}), $W'\cap (B\times {\bf I})$ consists of one U-turn component, hence the value under ${\rm Tr}^\omega_B$ is governed by Prop.\ref{prop:biangle_SL3_quantum_trace}(BT2-2), especially eq.\eqref{eq:Tr_B_U2}. Meanwhile, $W'\cap (\wh{t}\times {\bf I})$ consists of two left turn corner arcs, so the values under $\til{\rm Tr}^\omega_t$ are given by Thm.\ref{thm:SL3_quantum_trace_map}\redfix{(QT2-1)}. In the state sum in eq.\eqref{eq:til_Tr_W_s_prime}, denoting by $\varepsilon_4,\varepsilon_5$ the states assigned to the junctures $w_1,w_2$, we get
\begin{align*}
\til{\rm Tr}^\omega_{\wh{\Delta}}(W',s')
& = {\textstyle \sum}_{\varepsilon_4,\varepsilon_5} (\wh{\bf M}^{\rm in}_{t,2} \, \wh{\bf M}^{\rm left}(\wh{Z}_{v_t}) \, \wh{\bf M}^{\rm out}_{t,3})_{\varepsilon_2,\varepsilon_5}
(\wh{\bf M}^{\rm U}_-)_{\varepsilon_5,\varepsilon_4} 
( \wh{\bf M}^{\rm in}_{t,3} \, \wh {\bf M}^{\rm left}(\wh{Z}_{v_t}) \, \wh{\bf M}^{\rm out}_{t,1} )_{\varepsilon_4,\varepsilon_1} \\
& = ( \wh{\bf M}^{\rm in}_{t,2} \, \wh{\bf M}^{\rm left}(\wh{Z}_{v_t}) \, \wh{\bf M}^{\rm out}_{t,3} \,\, \wh{\bf M}^{\rm U}_- \,\, \wh{\bf M}^{\rm in}_{t,3} \, \wh {\bf M}^{\rm left}(\wh{Z}_{v_t}) \, \wh{\bf M}^{\rm out}_{t,1})_{\varepsilon_2,\varepsilon_1},
\end{align*}
hence the problem of showing $\til{\rm Tr}^\omega_{\wh{\Delta}}(W,s) = \til{\rm Tr}^\omega_{\wh{\Delta}}(W',s')$ boils down to showing the matrix identity
\begin{align}
\label{eq:matrix_identity_to_verify1}
\wh{\bf M}^{\rm in}_{t,2} \, \wh{\bf M}^{\rm right}(\wh{Z}_{v_t}) \, \wh{\bf M}^{\rm out}_{t,1}
= 
\wh{\bf M}^{\rm in}_{t,2} \, \wh{\bf M}^{\rm left}(\wh{Z}_{v_t}) \, \wh{\bf M}^{\rm out}_{t,3} \,\, \wh{\bf M}^{\rm U}_- \,\, \wh{\bf M}^{\rm in}_{t,3} \, \wh {\bf M}^{\rm left}(\wh{Z}_{v_t}) \, \wh{\bf M}^{\rm out}_{t,1}.
\end{align}
Before showing this, note that since all remaining cases of Fig.\ref{fig:isotopy_move4} and Fig.\ref{fig:isotopy_move5} are such that each of $W$ and $W'$ is a single oriented curve, the problem for each of them also boils down to checking an identity of products of quantum monodromy matrices for segments. We first collect all such matrix identities to check, and show them altogether. First, still for the right case of Fig.\ref{fig:isotopy_move4} but with opposite orientations, with $w_1\succ w_2$, we should show
\begin{align}
\label{eq:matrix_identity_to_verify2}
\wh{\bf M}^{\rm in}_{t,1} \,
\wh{\bf M}^{\rm left}(\wh{Z}_{v_t}) \,
\wh{\bf M}^{\rm out}_{t,2}
= 
\wh{\bf M}^{\rm in}_{t,1} \,
\wh{\bf M}^{\rm right}(\wh{Z}_{v_t}) \,
\wh{\bf M}^{\rm out}_{t,3} 
\,\,
\wh{\bf M}^{\rm U}_+
\,\,
\wh{\bf M}^{\rm in}_{t,3} \,
\wh{\bf M}^{\rm right}(\wh{Z}_{v_t}) \,
\wh{\bf M}^{\rm out}_{t,2}.
\end{align}
For the left case of Fig.\ref{fig:isotopy_move5} with $x_1\prec x_2$, $w_1\prec w_2$, and for the right case of Fig.\ref{fig:isotopy_move5} with reversed orientations with $y_1\succ y_2$, $w_1\succ w_2$, we should show
\begin{align}
\label{eq:matrix_identity_to_verify3}
\wh{\bf M}^{\rm U}_-
= \wh{\bf M}^{\rm in}_{t,1} \,
\wh{\bf M}^{\rm right}(\wh{Z}_{v_t}) \,
\wh{\bf M}^{\rm out}_{t,3} \,\, 
\wh{\bf M}^{\rm U}_- \,\,
\wh{\bf M}^{\rm in}_{t,3} \,
\wh{\bf M}^{\rm left}(\wh{Z}_{v_t}) \,
\wh{\bf M}^{\rm out}_{t,1}, \\
\label{eq:matrix_identity_to_verify4}
\wh{\bf M}^{\rm U}_-
= \wh{\bf M}^{\rm in}_{t,2} \,
\wh{\bf M}^{\rm left}(\wh{Z}_{v_t}) \,
\wh{\bf M}^{\rm out}_{t,3} \,\, 
\wh{\bf M}^{\rm U}_- \,\, 
\wh{\bf M}^{\rm in}_{t,3} \,
\wh{\bf M}^{\rm right}(\wh{Z}_{v_t}) \,
\wh{\bf M}^{\rm out}_{t,2}.
\end{align}
The orientation reversed versions of these two, with $x_1\succ x_2$, $w_1\succ w_2$, $y_1\succ y_2$, yield identical matrix identities as themselves, with each $\wh{\bf M}^{\rm U}_-$ replaced by $(\wh{\bf M}^{\rm U}_+)^{\rm tr}$, where ${\rm tr}$ means the transpose:
\begin{align}
\label{eq:matrix_identity_to_verify5}
(\wh{\bf M}^{\rm U}_+)^{\rm tr}
= \wh{\bf M}^{\rm in}_{t,1} \,
\wh{\bf M}^{\rm right}(\wh{Z}_{v_t}) \,
\wh{\bf M}^{\rm out}_{t,3} \,\, 
(\wh{\bf M}^{\rm U}_+)^{\rm tr} \,\,
\wh{\bf M}^{\rm in}_{t,3} \,
\wh{\bf M}^{\rm left}(\wh{Z}_{v_t}) \,
\wh{\bf M}^{\rm out}_{t,1}, \\
\label{eq:matrix_identity_to_verify6}
(\wh{\bf M}^{\rm U}_+)^{\rm tr}
= \wh{\bf M}^{\rm in}_{t,2} \,
\wh{\bf M}^{\rm left}(\wh{Z}_{v_t}) \,
\wh{\bf M}^{\rm out}_{t,3} \,\, 
(\wh{\bf M}^{\rm U}_+)^{\rm tr} \,\, 
\wh{\bf M}^{\rm in}_{t,3} \,
\wh{\bf M}^{\rm right}(\wh{Z}_{v_t}) \,
\wh{\bf M}^{\rm out}_{t,2}.
\end{align}
To prove these matrix identities, we observe:
\begin{lemma}[compatibility relations among quantum monodromy matrices]
\label{lem:compatibility_relations_among_quantum_monodromy_matrices}
For each $\alpha \in \{1,2,3\}$ and $k\in \{+,-\}$, one has
\begin{align}
\label{eq:quantum_compatiblity1}
& \wh{\bf M}^{\rm out}_{t,\alpha} \, 
\wh{\bf M}^{\rm U}_k
\wh{\bf M}^{\rm in}_{t,\alpha} = \wh{\bf M}^{\rm U}_k = \wh{\bf M}^{\rm in}_{t,\alpha} \, \wh{\bf M}^{\rm U}_k \wh{\bf M}^{\rm out}_{t,\alpha}, \\
\label{eq:quantum_compatiblity2}
& \wh{\bf M}^{\rm left}(\wh{Z}) \, 
\wh{\bf M}^{\rm U}_- \, 
\wh{\bf M}^{\rm left}(\wh{Z}) = \wh{\bf M}^{\rm right}(\wh{Z}), \qquad
\wh{\bf M}^{\rm right}(\wh{Z}) \, \wh{\bf M}^{\rm U}_+ \, \wh{\bf M}^{\rm right}(\wh{Z}) = \wh{\bf M}^{\rm left}(\wh{Z}), \\
\label{eq:quantum_compatiblity3}
& \wh{\bf M}^{\rm left}(\wh{Z}) \, \wh{\bf M}^{\rm U}_- \, \wh{\bf M}^{\rm right}(\wh{Z}) = \wh{\bf M}^{\rm U}_- = \wh{\bf M}^{\rm right}(\wh{Z}) \, \wh{\bf M}^{\rm U}_- \, \wh{\bf M}^{\rm left}(\wh{Z}), \\
\label{eq:quantum_compatiblity4}
& \wh{\bf M}^{\rm left}(\wh{Z}) \, (\wh{\bf M}^{\rm U}_+)^{\rm tr} \, \wh{\bf M}^{\rm right}(\wh{Z}) = (\wh{\bf M}^{\rm U}_+)^{\rm tr} = \wh{\bf M}^{\rm right}(\wh{Z}) \, (\wh{\bf M}^{\rm U}_+)^{\rm tr} \, \wh{\bf M}^{\rm left}(\wh{Z})
\end{align}
\end{lemma}
{\it Proof of Lem.\ref{lem:compatibility_relations_among_quantum_monodromy_matrices}.} It is a straightforward exercise to check; see \cite[Thm.9.2]{FG06} for a classical version. Note that, for eq.\eqref{eq:quantum_compatiblity1}, the (proof of the) classical version of almost yields a proof of the above quantum version, which is easy to verify by hand. The remaining matrix identities are not hard to check by hand either. For example, the first equality in eq.\eqref{eq:quantum_compatiblity2} can be checked as,
\begin{align*}
& \wh{\bf M}^{\rm left}(\wh{Z}) \, 
\wh{\bf M}^{\rm U}_- \, 
\wh{\bf M}^{\rm left}(\wh{Z}) \\
& = \smallmatthree{\omega^5 \wh{Z}^2}{~\omega^{-1} \wh{Z}^2 + \omega^2 \wh{Z}^{-1} ~}{\omega^{-4} \wh{Z}^{-1}}{0}{\omega^5 \wh{Z}^{-1}}{\omega^{-1} \wh{Z}^{-1}}{0}{0}{\omega^2 \wh{Z}^{-1}}
\smallmatthree{0}{0}{q^{-7/3}}{0}{-q^{-4/3}}{0}{q^{-1/3}}{0}{0}
\smallmatthree{\omega^5 \wh{Z}^2}{~\omega^{-1} \wh{Z}^2 + \omega^2 \wh{Z}^{-1} ~}{\omega^{-4} \wh{Z}^{-1}}{0}{\omega^5 \wh{Z}^{-1}}{\omega^{-1} \wh{Z}^{-1}}{0}{0}{\omega^2 \wh{Z}^{-1}} \\
& = \smallmatthree{\omega^5 \wh{Z}^2}{~\omega^{-1} \wh{Z}^2 + \omega^2 \wh{Z}^{-1} ~}{\omega^{-4} \wh{Z}^{-1}}{0}{\omega^5 \wh{Z}^{-1}}{\omega^{-1} \wh{Z}^{-1}}{0}{0}{\omega^2 \wh{Z}^{-1}}
\smallmatthree{0}{0}{q^{-7/3}\omega^2 \wh{Z}^{-1}}{0}{-q^{-4/3}\omega^5 \wh{Z}^{-1}}{-q^{-4/3}\omega^{-1} \wh{Z}^{-1}}{q^{-1/3}\omega^5 \wh{Z}^2}{~q^{-1/3}\omega^{-1} \wh{Z}^2 + q^{-1/3}\omega^2 \wh{Z}^{-1} ~}{q^{-1/3}\omega^{-4} \wh{Z}^{-1}} \\
& = \smallmatthree{\omega^{-2} \wh{Z}}{0}{0}{\omega\wh{Z}}{\omega^{-5} \wh{Z}}{0}{\omega^4 \wh{Z}}{\omega^{-2} \wh{Z} + \omega \wh{Z}^{-2}}{\omega^{-5} \wh{Z}^{-2}} = \wh{\bf M}^{\rm right}(\wh{Z}),
\end{align*}
where we used $q = \omega^{9}$, and the second equality in eq.\eqref{eq:quantum_compatiblity4} can be checked as
\begin{align*}
& \wh{\bf M}^{\rm right}(\wh{Z}) \, (\wh{\bf M}^{\rm U}_+)^{\rm tr} \, \wh{\bf M}^{\rm left}(\wh{Z}) \\
& = \smallmatthree{\omega^{-2} \wh{Z}}{0}{0}{\omega\wh{Z}}{\omega^{-5} \wh{Z}}{0}{\omega^4 \wh{Z}}{\omega^{-2} \wh{Z} + \omega \wh{Z}^{-2}}{\omega^{-5} \wh{Z}^{-2}} 
\smallmatthree{0}{0}{q^{1/3}}{0}{-q^{4/3}}{0}{q^{7/3}}{0}{0} 
\smallmatthree{\omega^5 \wh{Z}^2}{~\omega^{-1} \wh{Z}^2 + \omega^2 \wh{Z}^{-1} ~}{\omega^{-4} \wh{Z}^{-1}}{0}{\omega^5 \wh{Z}^{-1}}{\omega^{-1} \wh{Z}^{-1}}{0}{0}{\omega^2 \wh{Z}^{-1}} \\
& = \smallmatthree{\omega^{-2} \wh{Z}}{0}{0}{\omega\wh{Z}}{\omega^{-5} \wh{Z}}{0}{\omega^4 \wh{Z}}{\omega^{-2} \wh{Z} + \omega \wh{Z}^{-2}}{\omega^{-5} \wh{Z}^{-2}} 
\smallmatthree{0}{0}{q^{1/3}\omega^2 \wh{Z}^{-1}}{0}{-q^{4/3}\omega^5 \wh{Z}^{-1}}{-q^{4/3}\omega^{-1} \wh{Z}^{-1}}{q^{7/3}\omega^5 \wh{Z}^2}{~(q^{7/3}\omega^{-1} \wh{Z}^2 + q^{7/3}\omega^2 \wh{Z}^{-1}) ~}{q^{7/3}\omega^{-4} \wh{Z}^{-1}} \\
& = \smallmatthree{0}{0}{q^{1/3}}{0}{-q^{4/3}}{0}{q^{7/3}}{0}{0} = (\wh{\bf M}^{\rm U}_+)^{\rm tr}.
\end{align*}
Remaining ones are left as exercises. {\it End of proof of Lem.\ref{lem:compatibility_relations_among_quantum_monodromy_matrices}.}

\vs

Using Lem.\ref{lem:compatibility_relations_among_quantum_monodromy_matrices}, it is easy to see that the sought-for equalities eq.\eqref{eq:matrix_identity_to_verify1},\eqref{eq:matrix_identity_to_verify2},\eqref{eq:matrix_identity_to_verify3},\eqref{eq:matrix_identity_to_verify4}, \eqref{eq:matrix_identity_to_verify5} and \eqref{eq:matrix_identity_to_verify6} are satisfied. 

\vs

Finally, the right case for Fig.\ref{fig:isotopy_move4} with $w_1 \prec w_2$ can either be checked directly (in the classical setting, it is same as the case $w_1 \succ w_2$), or follows from the case $w_1\succ w_2$ together with the next Prop.\ref{prop:elementary_isotopy_invariance_elevation_change} \qed \quad {\it End of proof of Prop.\ref{prop:elementary_isotopy_invariance_elevation_preserving_arcs}.}

\vs

\begin{figure}[htbp!]
\vspace{-3mm}
\begin{center}
\raisebox{-0.5\height}{\scalebox{0.5}{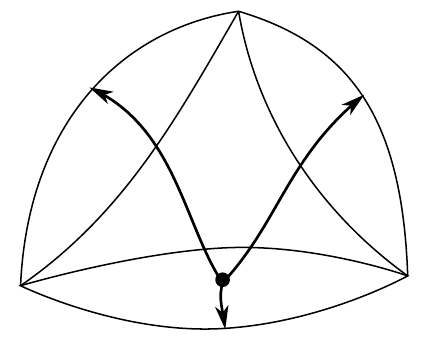}}
\hspace{-4mm} 
$\leftrightarrow$
\hspace{-3mm}
\raisebox{-0.5\height}{\scalebox{0.5}{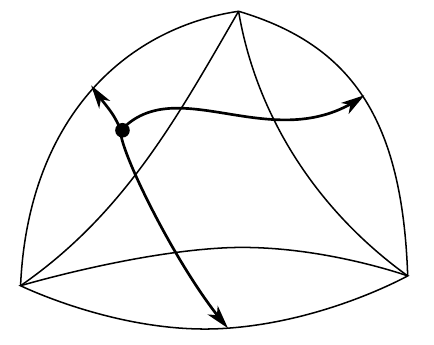}}
\quad
\raisebox{-0.5\height}{\scalebox{0.5}{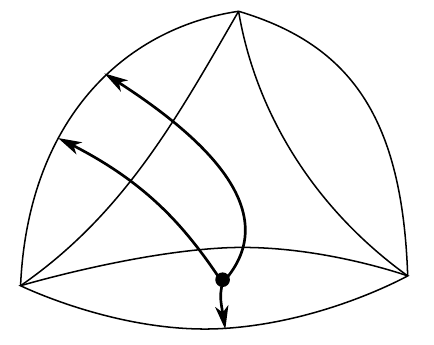}}
\hspace{-4mm} 
$\leftrightarrow$
\hspace{-3mm}
\raisebox{-0.5\height}{\scalebox{0.5}{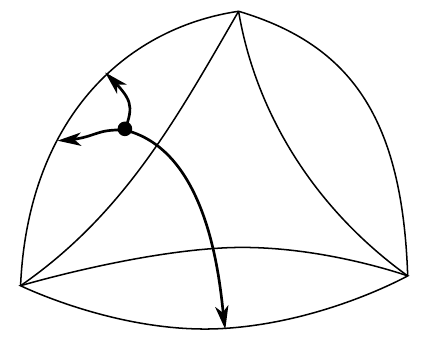}}
\quad
\raisebox{-0.5\height}{\scalebox{0.5}{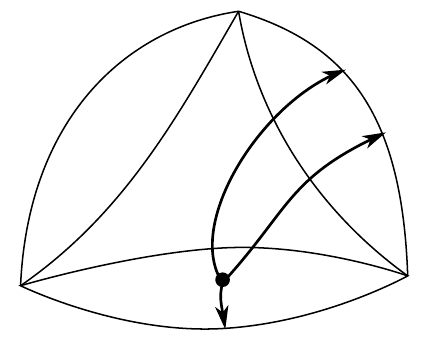}}
\hspace{-4mm} 
$\leftrightarrow$
\hspace{-3mm}
\raisebox{-0.5\height}{\scalebox{0.5}{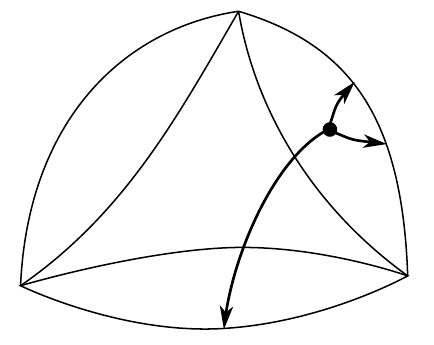}}
\end{center}
\vspace{-3mm}
\caption{Gool moves for 3-ways}
\vspace{-2mm}
\label{fig:gool_move2}
\end{figure}

\begin{figure}[htbp!]
\vspace{-1mm}
\begin{center}
\raisebox{-0.5\height}{\scalebox{0.5}{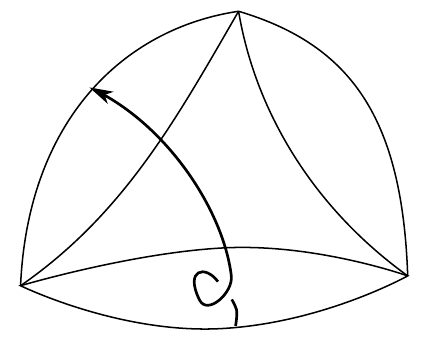}}
\hspace{-4mm} 
$\leftrightarrow$
\hspace{-3mm}
\raisebox{-0.5\height}{\scalebox{0.5}{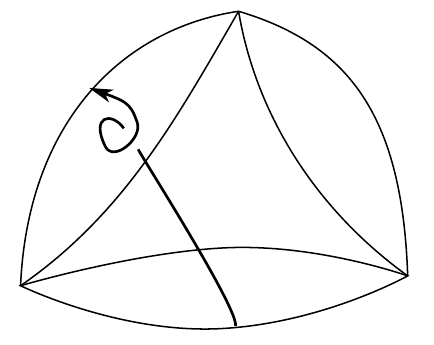}}
\quad
\raisebox{-0.5\height}{\scalebox{0.5}{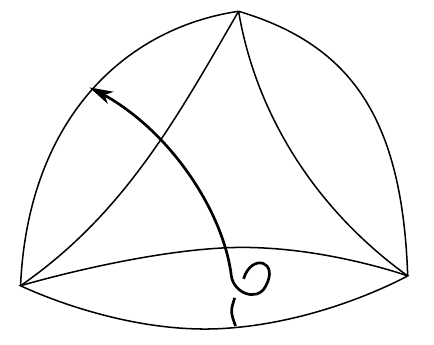}}
\hspace{-4mm} 
$\leftrightarrow$
\hspace{-3mm}
\raisebox{-0.5\height}{\scalebox{0.5}{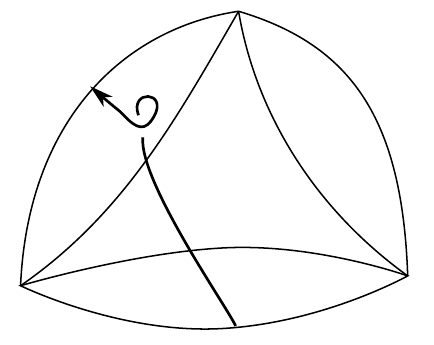}}
\end{center}
\vspace{-4mm}
\caption{Gool moves for kinks}
\vspace{-2mm}
\label{fig:gool_move_kink}
\end{figure}

We now observe the following topological lemma, whose proof can be obtained in the style of Lem.24 of \cite{BW}.
\begin{lemma}[moves between gool positions]
\label{lem:moves_between_gool_positions}
Let $\frak{S}$, $\Delta$, and $\wh{\Delta}$ be as in Def.\ref{def:good_position}. Let $W$ and $W'$ be ${\rm SL}_3$-webs in $\frak{S} \times {\bf I}$ in gool positions with respect to $\wh{\Delta}$, such that $W$ is isotopic to $W'$ as ${\rm SL}_3$-webs in $\frak{S}\times {\bf I}$. Then $W$ can be connected to $W'$ by a sequence of ${\rm SL}_3$-webs $W=W_1,W_2,\ldots,W_n=W'$ in gool positions with respect to $\wh{\Delta}$, such that for each $i=1,\ldots,n-1$, $W_i$ is related to $W_{i+1}$ either by an isotopy within the class of ${\rm SL}_3$-webs in gool positions with respect to $\wh{\Delta}$ or by one of the moves in Fig.\ref{fig:isotopy_move4}, \ref{fig:isotopy_move5}, \ref{fig:isotopy_move6}, \ref{fig:isotopy_move7}, \ref{fig:gool_move2}, \ref{fig:gool_move_kink}, possibly with different possible orientations on the components. \qed
\end{lemma}

In fact, the statement of Lem.\ref{lem:moves_between_gool_positions} should be more refined. Namely, in the definition of the above moves, each biangle in the pictures may not precisely be one entire biangle in the split ideal triangulation $\wh{\Delta}$. Before applying the move as depicted in the picture, one may have to divide a biangle of $\wh{\Delta}$ into several biangles, by introducing some ideal arcs in this biangle connecting the two marked points. This will yield a generalized version $\wh{\Delta}'$ of split ideal triangulation, which may have more than one biangles per each edge of $\Delta$. We require that we draw the new arcs so that the ${\rm SL}_3$-web in question is still transverse to the thickening of the edges of $\wh{\Delta}'$. Then apply the moves as in the above pictures, for the part of ${\rm SL}_3$-web living over the union of a triangle of $\wh{\Delta}'$ and its three neighboring biangles. For example, a move like \raisebox{-0.5\height}{\scalebox{0.4}{
\begingroup%
  \makeatletter%
  \providecommand\color[2][]{%
    \errmessage{(Inkscape) Color is used for the text in Inkscape, but the package 'color.sty' is not loaded}%
    \renewcommand\color[2][]{}%
  }%
  \providecommand\transparent[1]{%
    \errmessage{(Inkscape) Transparency is used (non-zero) for the text in Inkscape, but the package 'transparent.sty' is not loaded}%
    \renewcommand\transparent[1]{}%
  }%
  \providecommand\rotatebox[2]{#2}%
  \newcommand*\fsize{\dimexpr\f@size pt\relax}%
  \newcommand*\lineheight[1]{\fontsize{\fsize}{#1\fsize}\selectfont}%
  \ifx\svgwidth\undefined%
    \setlength{\unitlength}{124.72440945bp}%
    \ifx\svgscale\undefined%
      \relax%
    \else%
      \setlength{\unitlength}{\unitlength * \real{\svgscale}}%
    \fi%
  \else%
    \setlength{\unitlength}{\svgwidth}%
  \fi%
  \global\let\svgwidth\undefined%
  \global\let\svgscale\undefined%
  \makeatother%
  \begin{picture}(1,0.79545455)%
    \lineheight{1}%
    \setlength\tabcolsep{0pt}%
    \put(0,0){\includegraphics[width=\unitlength,page=1]{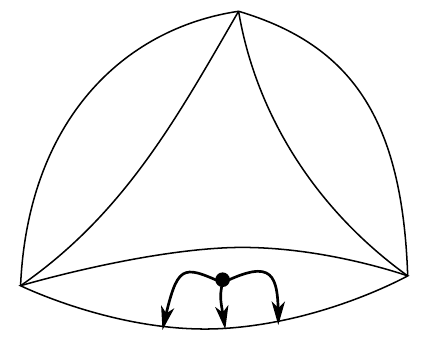}}%
  \end{picture}%
\endgroup%
}}
\hspace{-3mm} 
$\leftrightarrow$
\hspace{-4mm}
\raisebox{-0.5\height}{\scalebox{0.4}{
\begingroup%
  \makeatletter%
  \providecommand\color[2][]{%
    \errmessage{(Inkscape) Color is used for the text in Inkscape, but the package 'color.sty' is not loaded}%
    \renewcommand\color[2][]{}%
  }%
  \providecommand\transparent[1]{%
    \errmessage{(Inkscape) Transparency is used (non-zero) for the text in Inkscape, but the package 'transparent.sty' is not loaded}%
    \renewcommand\transparent[1]{}%
  }%
  \providecommand\rotatebox[2]{#2}%
  \newcommand*\fsize{\dimexpr\f@size pt\relax}%
  \newcommand*\lineheight[1]{\fontsize{\fsize}{#1\fsize}\selectfont}%
  \ifx\svgwidth\undefined%
    \setlength{\unitlength}{124.72440945bp}%
    \ifx\svgscale\undefined%
      \relax%
    \else%
      \setlength{\unitlength}{\unitlength * \real{\svgscale}}%
    \fi%
  \else%
    \setlength{\unitlength}{\svgwidth}%
  \fi%
  \global\let\svgwidth\undefined%
  \global\let\svgscale\undefined%
  \makeatother%
  \begin{picture}(1,0.79545455)%
    \lineheight{1}%
    \setlength\tabcolsep{0pt}%
    \put(0,0){\includegraphics[width=\unitlength,page=1]{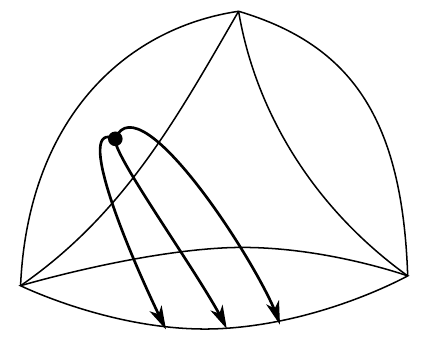}}%
  \end{picture}%
\endgroup%
}} is redundant, as one can show that this can be obtained as composition of the above moves, applied in the sense just described.

\vs

{\it Proof of Prop.\ref{prop:isotopy_invariance_of_state-sum_trace_for_gool_positions}.} In view of Lem.\ref{lem:moves_between_gool_positions}, it suffices to show that $\wh{\rm Tr}^\omega_\Delta(W,s) = \wh{\rm Tr}^\omega_\Delta(W',s')$ in case $(W,s)$ and $(W',s')$ are related by one of those moves, in the sense just described above, using generalized split ideal triangulation $\wh{\Delta}'$. In each union of one triangle and three neighboring (thin) biangles, one observes that these moves, except for the ones in Fig.\ref{fig:gool_move_kink}, can be obtained as compositions of the moves dealt with in Propositions \ref{prop:elementary_isotopy_invariance_3-way}, \ref
{prop:elementary_isotopy_invariance_elevation_preserving_arcs} and \ref{prop:elementary_isotopy_invariance_elevation_change}. Now, write $\wh{\rm Tr}^\omega_\Delta(W,s)$ and $\wh{\rm Tr}^\omega_\Delta(W',s')$ using a state-sum formula adapted to $\wh{\Delta}'$, instead of $\wh{\Delta}$; so we may have more numbers of junctures, and more numbers of biangle factors. By Prop.\ref{prop:biangle_SL3_quantum_trace}(BT1), the new state-sum formulas give same answers as before which used $\wh{\Delta}$. Now, in these new state-sum expressions, the parts involving the above mentioned union of a triangle and three biangles have equal values for $(W,s)$ and $(W',s')$, due to Propositions \ref{prop:elementary_isotopy_invariance_3-way}, \ref
{prop:elementary_isotopy_invariance_elevation_preserving_arcs} and \ref{prop:elementary_isotopy_invariance_elevation_change}. Thus $\wh{\rm Tr}^\omega_\Delta(W,s) = \wh{\rm Tr}^\omega_\Delta(W',s')$. The case of Fig.\ref{fig:gool_move_kink} follows from the well-definedness of ${\rm Tr}^\omega_B$, and the relation
$$
\raisebox{-0.5\height}{\scalebox{1.0}{
\begingroup%
  \makeatletter%
  \providecommand\color[2][]{%
    \errmessage{(Inkscape) Color is used for the text in Inkscape, but the package 'color.sty' is not loaded}%
    \renewcommand\color[2][]{}%
  }%
  \providecommand\transparent[1]{%
    \errmessage{(Inkscape) Transparency is used (non-zero) for the text in Inkscape, but the package 'transparent.sty' is not loaded}%
    \renewcommand\transparent[1]{}%
  }%
  \providecommand\rotatebox[2]{#2}%
  \newcommand*\fsize{\dimexpr\f@size pt\relax}%
  \newcommand*\lineheight[1]{\fontsize{\fsize}{#1\fsize}\selectfont}%
  \ifx\svgwidth\undefined%
    \setlength{\unitlength}{121.88976378bp}%
    \ifx\svgscale\undefined%
      \relax%
    \else%
      \setlength{\unitlength}{\unitlength * \real{\svgscale}}%
    \fi%
  \else%
    \setlength{\unitlength}{\svgwidth}%
  \fi%
  \global\let\svgwidth\undefined%
  \global\let\svgscale\undefined%
  \makeatother%
  \begin{picture}(1,0.34883721)%
    \lineheight{1}%
    \setlength\tabcolsep{0pt}%
    \put(0.02400617,0.15545937){\color[rgb]{0,0,0}\makebox(0,0)[lt]{\lineheight{1.25}\smash{\begin{tabular}[t]{l}$q^{8/3}$\end{tabular}}}}%
    \put(0,0){\includegraphics[width=\unitlength,page=1]{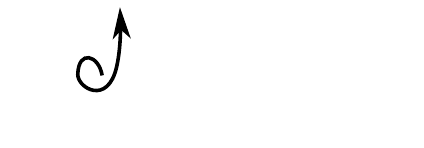}}%
    \put(0.35627355,0.15545937){\color[rgb]{0,0,0}\makebox(0,0)[lt]{\lineheight{1.25}\smash{\begin{tabular}[t]{l}$=$\end{tabular}}}}%
    \put(0.565479,0.15545937){\color[rgb]{0,0,0}\makebox(0,0)[lt]{\lineheight{1.25}\smash{\begin{tabular}[t]{l}$=q^{-8/3}$\end{tabular}}}}%
    \put(0,0){\includegraphics[width=\unitlength,page=2]{skein_kink_rel1.pdf}}%
  \end{picture}%
\endgroup%
}}
$$
in $\mathcal{S}^\omega_{\rm s}(B;\mathbb{Z})_{\rm red}$, which is a consequence of the defining relations of the ${\rm SL}_3$-skein algebras in Fig.\ref{fig:A2-skein_relations_quantum} (see e.g. \cite[Prop.1]{Higgins}). \qed

\vs

{\it Proof of Prop.\ref{prop:state-sum_trace_for_good_positions}.} Let $\frak{S}$, $\Delta$, and $\wh{\Delta}$ be as in Def.\ref{def:good_position}. Let $(W,s)$ be a stated ${\rm SL}_3$-web in $\frak{S} \times {\bf I}$ in a good position with respect to $\wh{\Delta}$. Let $t$ be a triangle of $\Delta$, and $\wh{t}$ be the corresponding triangle of $\wh{\Delta}$, such that the ${\rm SL}_3$-web $W\cap (\wh{t}\times {\bf I})$ contains a component that is either a U-turn arc or a 3-way web. One can push the whole of such U-turn arc or the 3-valent vertex to a neighboring biangle $B$, by an isotopy. Before pushing, one can divide $B$ into two biangles $B_1$ and $B_2$ by considering an ideal arc $e$ in $B$, such that $B_1$ is adjacent to $\wh{t}$ and $W\cap (B_1 \times {\bf I})$ consists only of the components of the form as in Prop.\ref{prop:biangle_SL3_quantum_trace}(BT2-1). Then, push the U-turn arc or the 3-valent vertex living in $\wh{t}$ into the biangle $B_1$. Then, by Propositions \ref{prop:elementary_isotopy_invariance_3-way} and \ref
{prop:elementary_isotopy_invariance_elevation_preserving_arcs}, the value of $(W\cap ( (\wh{t} \cup B_1) \times {\bf I}), s|_{\partial(W\cap (\wh{t} \cup B_1))})$ under the state-sum trace $\til{\rm Tr}^\omega_{\wh{t} \cup B_1}$ is unchanged by such a pushing. Meanwhile, one can observe that the state-sum trace $\til{\rm Tr}^\omega_\Delta$ as defined in Prop.\ref{prop:state-sum_trace_for_good_positions} using the split ideal triangulation $\wh{\Delta}$ equals the new state-sum trace defined also as in Prop.\ref{prop:state-sum_trace_for_good_positions} (i.e. as in eq.\eqref{eq:state-sum_formula} of Def.\ref{def:state-sum_trace_for_gool_position}) but this time for the finer decomposition $\wh{\Delta} \cup\{e\}$ with one more number of biangles, using Prop.\ref{prop:biangle_SL3_quantum_trace}(BT1). This shows that the value under $\til{\rm Tr}^\omega_\Delta$ does not change after such a pushing. By a finite number of such pushing moves, one can put $(W,s)$ into an ${\rm SL}_3$-web $(W',s')$ in a gool position with respect to $\wh{\Delta}$; so $\til{\rm Tr}^\omega_\Delta(W,s) = \til{\rm Tr}^\omega_\Delta(W',s')$. Meanwhile, we have $\til{\rm Tr}^\omega_\Delta(W',s') = \wh{\rm Tr}^\omega_\Delta(W',s')$ by construction. Hence $\til{\rm Tr}^\omega_\Delta(W,s) = \wh{\rm Tr}^\omega_\Delta(W',s')$. \qed

\vs

Consequently, Thm.\ref{thm:SL3_quantum_trace_map} is finally proved, modulo a complete proof of Prop.\ref{prop:elementary_isotopy_invariance_elevation_change}, via the argument at the end of \S\ref{subsec:state-sum_construction}. We state a computationally useful corollary:
\begin{corollary}[the state-sum formula for the ${\rm SL}_3$ quantum trace]
Let $\frak{S}$, $\Delta$, and $\wh{\Delta}$ be as in Def.\ref{def:good_position}. For any stated ${\rm SL}_3$-web $(W,s)$ in $\frak{S}\times {\bf I}$ in a good position with respect to $\wh{\Delta}$, one has
$$
{\rm Tr}^\omega_\Delta([W,s]) = \til{\rm Tr}^\omega_\Delta(W,s). \qed
$$
\end{corollary}

\subsection{Congruence of terms, and the highest term}

Recall that one major motivation for our study of the ${\rm SL}_3$ quantum trace map ${\rm Tr}^\omega_\Delta$ was to prove Prop.\ref{prop:highest_term} and Prop.\ref{prop:congruence-compatibility_of_terms_of_each_basic_regular_function}, which are on the highest term and the congruence of terms of the basic semi-regular function $\mathbb{I}^+_{{\rm PGL}_3}(\ell) \in C^\infty(\mathscr{X}^+_{{\rm PGL}_3})$ for each $\ell \in \mathscr{A}_{\rm L}(\frak{S};\mathbb{Z})$. For that purpose, we only need to deal with the ${\rm SL}_3$ classical trace ${\rm Tr}_\Delta$ (eq.\eqref{eq:Tr_Delta}), i.e. ${\rm Tr}^1_\Delta$ when $\omega^{1/2}=1$. For convenience, the generator variables $\wh{Z}_v$ and $\wh{Z}_{t,v}$ for $\mathcal{Z}^1_\Delta$ and $\mathcal{Z}^1_t$ \redfix{may} be denoted by $Z_v$ and $Z_{t,v}$ respectively, without the hats. Also the quantum monodromy matrices $\wh{\bf M}^{\rm out}_{t,\alpha}$, $\wh{\bf M}^{\rm out}_{t,\alpha}$, $\wh{\bf M}^{\rm left}$, $\wh{\bf M}^{\rm right}$ will be used without hats, to emphasize the classical case. We also use the classical version of the ${\rm SL}_3$ biangle trace ${\rm Tr}^\omega_B$, as
constructed in Cor.\ref{cor:biangle_SL3_classical_trace}. In the present subsection, which is written in terms of the surface $\frak{S}$ instead of the 3d space $\frak{S} \times {\bf I}$ (using Lem.\ref{lem:surface_stated_skein_algebra_and_commutative_stated_skein_algebra}), we establish the counterparts of Prop.\ref{prop:highest_term} and Prop.\ref{prop:congruence-compatibility_of_terms_of_each_basic_regular_function} for the ${\rm SL}_3$ classical trace map ${\rm Tr}_\Delta$.  
\begin{proposition}[congruence of terms of the ${\rm SL}_3$ classical trace]
\label{prop:congruence_of_terms_of_SL3_classical_trace}
Let $\Delta$ be an ideal triangulation of a triangulable generalized marked surface $\frak{S}$, and let $(W,s)$ be a stated ${\rm SL}_3$-web in $\frak{S}$. Then the value of the ${\rm SL}_3$ classical trace map ${\rm Tr}_\Delta([W,s])\in \mathcal{Z}_\Delta$ can be written as a Laurent polynomial in the generators $\{Z_v \, |\, v\in \mathcal{V}(Q_\Delta)\}$ of $\mathcal{Z}_\Delta$ with integer coefficients so that all appearing Laurent monomials are congruent to each other in the following sense: for any two Laurent monomials $\prod_v Z_v^{\alpha_v}$ and $\prod_v Z_v^{\beta_v}$ appearing in this Laurent polynomial (with $(\alpha_v)_v,(\beta_v)_v \in \mathbb{Z}^{\mathcal{V}(Q_\Delta)}$), we have $\alpha_v - \beta_v \in 3\mathbb{Z}$ for all $v\in \mathcal{V}(Q_\Delta)$.
\end{proposition}

{\it Proof.} We may assume that $(W,s)$ is in a gool position with respect to a split ideal triangulation $\wh{\Delta}$ for $\Delta$, by applying an isotopy. Let's use the state-sum formula for ${\rm Tr}_\Delta([W,s])$ as in eq.\eqref{eq:state-sum_formula} in Def.\ref{def:state-sum_trace_for_gool_position}. Note that the biangle factors ${\rm Tr}_B([W\cap B, J_B])$ are integers. The triangle factor $\wh{\rm Tr}_t(W\cap \wh{t}, J_t)$ is a product of factors $\wh{\rm Tr}_t (W_{t,j},J_{t,j})$, whose values are, as described in \redfix{(QT2-1)}--\redfix{(QT2-2)} of Thm.\ref{thm:SL3_quantum_trace_map} (or, Prop.\ref{cor:SL3_classical_trace_map}), entries of certain products of matrices. By inspection, all the nonzero entries are congruent to each other. \qed

\vs

\begin{proposition}[the highest term of the ${\rm SL}_3$ classical trace]
\label{prop:highest_term_of_SL3_classical_trace}
Let $\Delta$ be an ideal triangulation of a triangulable generalized marked surface $\frak{S}$, and $\wh{\Delta}$ be a split ideal triangulation for $\Delta$. Let $W$ be a (reduced) non-elliptic ${\rm SL}_3$-web in $\frak{S}$ in a canonical position with respect to $\wh{\Delta}$ (Def.\ref{def:canonical_wrt_split_ideal_triangulation}) that has no endpoints. View $W$ as an ${\rm SL}_3$-lamination by giving the weight $1$; let ${\rm a}_v(W) \in \frac{1}{3}\mathbb{Z}$, $v\in \mathcal{V}(Q_\Delta)$, be the tropical coordinates defined in Def.\ref{def:tropical_coordinates}. Then ${\rm Tr}_\Delta([W,{\O}]) \in \mathcal{Z}_\Delta$ can be written as a Laurent polynomial in $\{Z_v \, | \, v\in \mathcal{V}(Q_\Delta)\}$ with integer coefficients so that $\prod_v Z_v^{3{\rm a}_v(W)}$ is the unique Laurent monomial among all appearing Laurent monomials that has higher partial ordering than any other appearing Laurent monomials, and this highest Laurent monomial has coefficient $1$.
\end{proposition}

\begin{figure}[htbp!]
\vspace{-5mm}
\begin{center}
\raisebox{-0.5\height}{\scalebox{0.8}{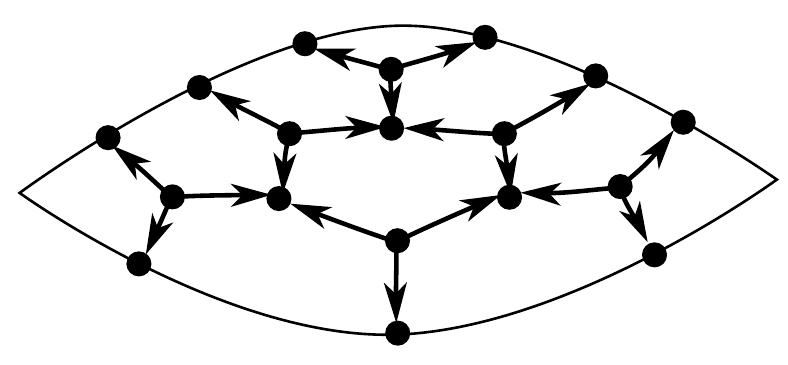}}
$\to$
\raisebox{-0.5\height}{\scalebox{0.8}{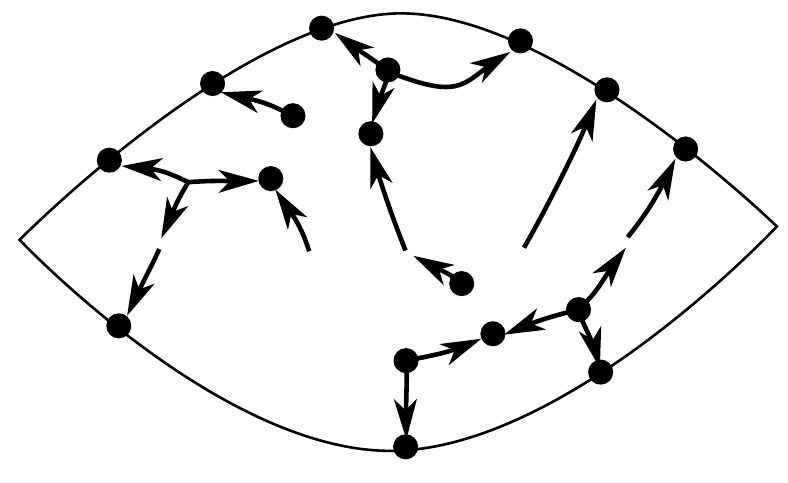}}
\end{center}
\vspace{-5mm}
\caption{Pyramid $P_d$ in a biangle $B$ and its decomposition (for $d=3$)}
\vspace{-2mm}
\label{fig:pyramid_in_biangle}
\end{figure}

This very important proposition is proved in several steps. The core lies in the following treatment of single-component canonical ${\rm SL}_3$-webs in a triangle (Def.\ref{def:canonical_web_in_a_triangle}).
\begin{proposition}[the highest term of the ${\rm SL}_3$ classical trace for a triangle]
\label{prop:highest_term_of_SL3_classical_trace_for_a_triangle}
Let $t$ be a triangle, viewed as a generalized marked surface. Let $W$ be a canonical ${\rm SL}_3$-web in $t$ (Def.\ref{def:canonical_web_in_a_triangle}). Let $\Delta$ be the unique triangulation of $t$, so that $Q_\Delta$ has seven nodes. For each $v\in \mathcal{V}(Q_\Delta)$, let ${\rm a}_v(W) \in \frac{1}{3}\mathbb{Z}$ be the tropical coordinate of $W$ as defined in Def.\ref{def:tropical_coordinates}, when $W$ is viewed as an ${\rm SL}_3$-lamination in $t$ with weight $1$. Denote by ${\bf 1}_W$ the state of $W$ assigning $1 \in \{1,2,3\}$ to all endpoints of $W$. Then
\begin{enumerate}
\item[\rm (HT1)] ${\rm Tr}_\Delta([W,{\bf 1}_W]) \in \mathcal{Z}_\Delta = \mathcal{Z}_{\wh{t}}$ can be written as a Laurent polynomial in $\{Z_v \, | \, v\in \mathcal{V}(Q_\Delta)\}$ with integer coefficients so that $\prod_{v\in \mathcal{V}(Q_\Delta)} Z_v^{3{\rm a}_v(W)}$ is the unique Laurent monomial of the highest partial ordering.

\item[\rm (HT2)] For any other state $s$ of $W$, ${\rm Tr}_\Delta([W,s]) \in \mathcal{Z}_\Delta = \mathcal{Z}_t$ can be written as a Laurent polynomial in $\{Z_v \, | \, v\in \mathcal{V}(Q_\Delta)\}$ so that each appearing Laurent monomial has strictly lower partial ordering than $\prod_{v\in \mathcal{V}(Q_\Delta)} Z_v^{3{\rm a}_v(W)}$.
\end{enumerate}
\end{proposition}
The hard case is when $W$ involves a pyramid (Def.\ref{def:canonical_web_in_a_triangle}). To compute ${\rm Tr}_\Delta([W,s])$, we will push all 3-valent vertices into the biangle for one of the sides of the triangle. Thus we find it convenient to build some lemmas for ${\rm SL}_3$-webs in a biangle.
\begin{lemma}[the ${\rm SL}_3$ classical trace for a pyramid $P_d$ in a biangle]
\label{lem:P_d}
Let $d\in \mathbb{Z}$, $d\neq 0$, and let $P_d$ be the ${\rm SL}_3$-web in a biangle $B$ as in the left picture of Fig.\ref{fig:pyramid_in_biangle} (there, an example is drawn for $d=3$), called a  \ul{\em degree $d$ pyramid in a biangle} $B$. One way of constructing $P_d$ is from a degree $d$ pyramid $H_d$ in a triangle as in Def.\ref{def:canonical_web_in_a_triangle} by removing (i.e. forgetting, or `filling in') one marked point of the triangle to turn it into a biangle. Label the endpoints of $P_d$ as $x_1,x_2,\ldots,x_{|d|}$, $y_1,y_2,\ldots,y_{|d|}$, $z_{|d|},\ldots,z_2,z_1$, appearing clockwise in this order on $\partial B$, so that $x_*,y_*$ lie on one side and $z_*$ lie on the other side. Suppose $s$ is a state of $P_d$ assigning $1$ to all $x_i$'s and $z_i$'s. Then ${\rm Tr}_B([P_d,s]) \neq 0$ if and only if $s$ assigns $2$ to all $y_i$'s, and for that $s$ we have ${\rm Tr}_B([P_d,s])=1$.
\end{lemma}

{\it Proof of Lem.\ref{lem:P_d}.} We will show the statement for $d>0$. The proof for $d<0$ is completely symmetric. We use induction on $d$. The base case is $d=1$. One notes that $P_1$ is a $3$-way web in $B$ falling into the case of (BT2-3) of Prop.\ref{prop:biangle_SL3_quantum_trace}, and it is easy to verify the desired statement. Now, let $d\ge 2$, and suppose we \redfix{have shown that} the statement holds for $P_{d-1}$. Observe that $P_d$ `contains' $P_{d-1}$ in its lower right corner (which is easier to see for $H_d$ and $H_{d-1}$ in a triangle), so that there exists an ideal arc $e$ in $B$ (drawn as a dotted line in the left picture of Fig.\ref{fig:pyramid_in_biangle}) connecting the two marked points of $B$, cutting $B$ into two biangles $B_1$ and $B_2$, such that the non-elliptic ${\rm SL}_3$-web $P_d \cap B_1$ in $B_1$ consists of one edge connecting the two sides of $B_1$ (having a red dot as one endpoint in the right picture of Fig.\ref{fig:pyramid_in_biangle}) and a degree $d-1$ pyramid $P_{d-1}$ in $B_1$ (having blue dots as some endpoints in the right picture of Fig.\ref{fig:pyramid_in_biangle}). Let's label the junctures of $P_d$ at $e$, i.e. the elements of $P_d \cap e$, as $w_1,w_2,\ldots,w_d$, $u_2,\ldots,u_d$, as in the right picture of Fig.\ref{fig:pyramid_in_biangle}. By the cutting property (BT1) of Prop.\ref{prop:biangle_SL3_quantum_trace} (or, of Cor.\ref{cor:biangle_SL3_classical_trace}) one has
\begin{align}
\label{eq:Tr_P_d}
{\rm Tr}_B([P_d,s]) = \underset{s_1,s_2}{\textstyle \sum} {\rm Tr}_{B_1}([P_d\cap B_1,s_1]) \, {\rm Tr}_{B_2}([P_d\cap B_2,s_2])
\end{align}
where the sum is over all states $s_1,s_2$ of the ${\rm SL}_3$-webs $P_d\cap B_1$ in $B_1$ and $P_d\cap B_2$ in $B_2$ compatible with $s$, in the sense as in Prop.\ref{prop:biangle_SL3_quantum_trace}(BT1). In particular, any such $s_1$ assigns $1$ to all $z_1,z_2,\ldots,z_d$ and any such $s_2$ assigns $1$ to all $x_1,\ldots,x_d$. If ${\rm Tr}_{B_1}([P_d\cap B_1,s_1])\neq 0$, then the value under ${\rm Tr}_{B_1}$ of each of the two components of $(P_d \cap B_1,s_1)$ must be nonzero, by multiplicativity of ${\rm Tr}_{B_1}$. The edge component, which connects the endpoints $z_1$ and $w_1$, falls into the case Prop.\ref{prop:biangle_SL3_quantum_trace}(BT2-1), hence it has nonzero ${\rm Tr}_{B_1}$ value iff $s_1$ assigns same value to $z_1$ and $w_1$, so $s_1(w_1)=1$.

\vs

We now investigate the ${\rm SL}_3$-web $P_d \cap B_2$ in $B_2$. It consists of $d$ components, where $d-1$ of them are edges connecting $u_i$ in $e$ and $y_i$ in the other side of $B_2$ (with $i=2,3,\ldots,d$); see the right picture of Fig.\ref{fig:pyramid_in_biangle}. Denote the remaining component as the ${\rm SL}_3$-web $K_d$ in a biangle $B_2$. Its endpoints on one side are $x_1,x_2,\ldots,x_d, y_1$ appearing in this order along clockwise orientation on $\partial B_2$, and the endpoints on the other side are $w_d,w_{d-1},\ldots, w_2,w_1$ appearing in this order along clockwise orientation on $\partial B_2$. Note $x_1,\ldots,x_d,y_1,w_1$ are sinks, $w_2,\ldots,w_d$ are sources, and there are $2d-1$ internal $3$-valent vertices. We prove:
\begin{lemma}
\label{lem:K_d}
Suppose $s_2$ is a state for $K_d$ assigning $1$ to all $x_1,\ldots,x_d$ and $w_1$. Then ${\rm Tr}_{B_2}([K_d,s_2])\neq 0$ if and only if $s_2$ assigns $1$ to all $w_2,\ldots,w_d$ and $2$ to $y_1$. For this $s_2$, we have ${\rm Tr}_{B_2}([K_d,s_2])=1$.
\end{lemma}

\begin{figure}[htbp!]
\vspace{-4mm}
\begin{center}
\raisebox{-0.5\height}{\scalebox{0.8}{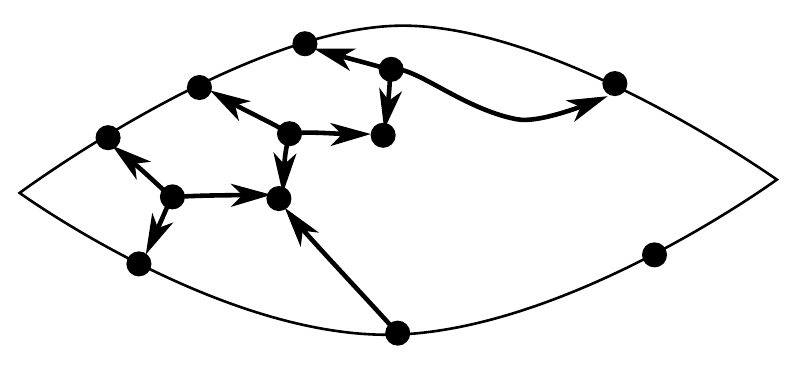}}
$\to$
\raisebox{-0.5\height}{\scalebox{0.8}{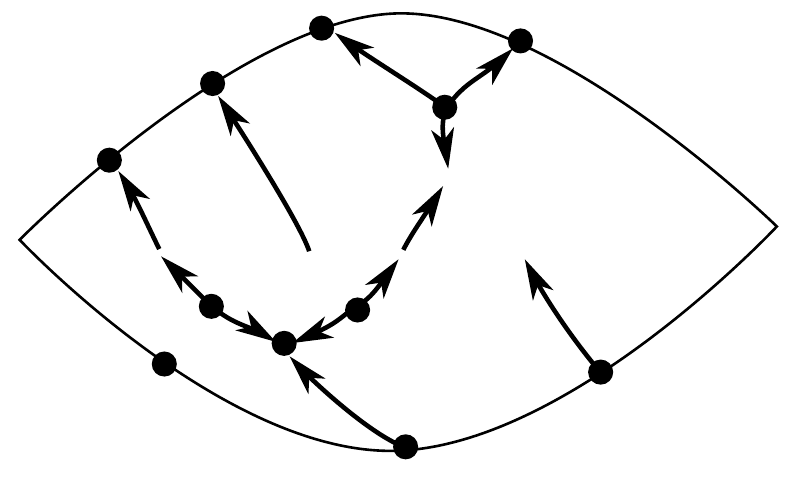}}
\end{center}
\vspace{-5mm}
\caption{${\rm SL}_3$-web $K_d$ in a biangle $B_2$ and its decomposition (for $d=3$)}
\vspace{-2mm}
\label{fig:K_d_in_biangle}
\end{figure}

{\it Proof of Lem.\ref{lem:K_d}.} We use induction on $d$. For the base case  $K_1$, note that the endpoints on one side of $B_2$ are $x_1,y_1$, while there is only one endpoint $z_1$ in the other side. There is only one internal $3$-valent vertex, so $K_1$ is just a 3-way web falling into (BT2-3) of Prop.\ref{prop:biangle_SL3_quantum_trace}, so the statement of Lem.\ref{lem:K_d} holds. Let $d\ge 2$, and suppose Lem.\ref{lem:K_d} holds for $K_{d-1}$. Observes that $K_d$ `contains' $K_{d-1}$ in its lower left corner, so that there exists an ideal arc $e'$ in $B_2$ (drawn as a dotted line in the left picture of Fig.\ref{fig:K_d_in_biangle}) cutting $B_2$ into biangles $B_3$ and $B_4$, such that the ${\rm SL}_3$-web $K_d \cap B_3$ in $B_3$ consists of $K_{d-1}$ in $B_3$ (having red dots and the blue dot as vertices in the right picture of Fig.\ref{fig:K_d_in_biangle}) and an edge connecting the two sides of $B_3$ (having purple dot as a vertex in the right picture of Fig.\ref{fig:K_d_in_biangle}). Label the junctures of $K_d$ at $e'$ as $r_1,\ldots,r_{d-1},r,r'$ as in Fig.\ref{fig:K_d_in_biangle}. By the cutting property (BT1) of Prop.\ref{prop:biangle_SL3_quantum_trace} one has
\begin{align}
\label{eq:Tr_K_d}
{\rm Tr}_{B_2} ([K_d,s_2]) = \underset{s_3,s_4}{\textstyle \sum} {\rm Tr}_{B_3}([K_d\cap B_3,s_3]) \, {\rm Tr}_{B_4}([K_d\cap B_4,s_4])
\end{align}
where the sum is over all states $s_3,s_4$ of the ${\rm SL}_3$-webs $K_d\cap B_3$ in $B_3$ and $K_d\cap B_4$ in $B_4$ compatible with $s_2$, in the sense as in Prop.\ref{prop:biangle_SL3_quantum_trace}(BT1). In particular, any such $s_4$ assigns $1$ to all $x_1,\ldots,x_d$, and any such $s_3$ assigns $1$ to $w_1$. 

\vs

Note $K_d \cap B_4$ has $d$ components, where $d-1$ of them are edges connecting $r_i$ in $e'$ and $x_i$ in the other side of $B_4$ (with $i=1,\ldots,d-1$). The remaining component can be called an {\em \redfix{\rm I}-web}, denoted by $I$. If ${\rm Tr}_{B_4}([K_d\cap B_4,s_4])\neq 0$, then the value under ${\rm Tr}_{B_4}$ of each of the $d$ components of $(K_d \cap B_4,s_4)$ must be nonzero, by multiplicativity of ${\rm Tr}_{B_4}$. The edge component, which connects the endpoints $r_i$ and $x_i$, falls into the case Prop.\ref{prop:biangle_SL3_quantum_trace}(BT2-1), hence it has nonzero ${\rm Tr}_{B_4}$ value iff $s_4$ assigns same value to $r_i$ and $x_i$, so $s_4(r_i)=1$ for all $i=1,\ldots,d-1$. For the remaining I-web $I$, we apply Lem.\ref{lem:value_of_I-webs_in_biangle}. It follows that, under the condition $s_4(x_d)=1$, we have ${\rm Tr}_{B_4} ([I,s_4|_{\partial I}])\neq 0$ iff $(s_4(y_1),s_4(r),s_4(r'))$ is one of $(2,2,1)$, $(3,3,1)$, $(2,1,2)$ or $(3,1,3)$. In the former two cases the value of ${\rm Tr}_{B_4} ([I,s_4|_{\partial I}])$ is $1$, while in the latter two cases this value is $-1$.

\vs

Let $s_4$ be as above, and let $s_3$ be some compatible state of $K_d\cap B_3$ such that ${\rm Tr}_{B_3}([K_d\cap B_3,s_3])\neq 0$. That is, so far we are requiring that $s_3,s_4$ be compatible with $s_2$, and that ${\rm Tr}_{B_4}([K_d\cap B_4,s_4])\neq 0$, ${\rm Tr}_{B_3}([K_d\cap B_3,s_3])\neq 0$. By multiplicativity of ${\rm Tr}_{B_3}$, the value under ${\rm Tr}_{B_3}$ of each of the $d$ components of $(K_d \cap B_3,s_3)$ must be nonzero. By compatibility, we have $s_3(r_i)=s_4(r_i)=1$ for all $i=1,\ldots,d-1$; we also had $s_3(w_1)=1$. Hence, the induction hypothesis applies for the $K_{d-1}$ component of $(K_d \cap B_3,s_3)$; so the value under ${\rm Tr}_{B_3}$ of this $K_{d-1}$ component is nonzero iff $s_3$ assigns $1$ to all $w_2,\ldots,w_{d-1}$ and $2$ to $r$, and in this case, the value is $1$. So $s_3(r)=2$. By compatibility, $s_4(r) = s_3(r)=2$, hence by the above observation on the I-web $I$, it must be $s_4(r') = 1$, and ${\rm Tr}_{B_4}([I,s_4|_{\partial I}])=1$. Again by compatibility, $s_3(r')=s_4(r')=1$. The edge component of $(K_d \cap B_3,s_3)$, which connects $r'$ and $w_d$, falls into Prop.\ref{prop:biangle_SL3_quantum_trace}(BT2-1), hence it has nonzero ${\rm Tr}_{B_3}$ value iff $s_3$ assigns same value to $r'$ and $w_d$, hence it follows $s_3(w_d)=1$.

\vs

To summarize, the unique pair of states  $s_3,s_4$ whose corresponding summand in eq.\eqref{eq:Tr_K_d} is nonzero assign the values $1$ to $r_1,\ldots,r_{d-1}$, $r'$, $w_2,\ldots, w_d$, and the value $2$ to $y_1,r$. For this choice of states, the summand is $1$.  This finishes proof of Lem.\ref{lem:K_d}. \qed

\vs

We go back to proof of Lem.\ref{lem:P_d}, investigating the sum in eq.\eqref{eq:Tr_P_d}. Let $s_1,s_2$ be states of $P_d\cap B_1$ and $P_d\cap B_2$ compatible with $s$, and whose corresponding summand of eq.\eqref{eq:Tr_P_d} is nonzero. Recall that we already know $s_1$ assigns value $1$ to $z_1,z_2,\ldots,z_d,w_1$, and $s_2$ assigns $1$ to $x_1,\ldots,x_d$. By multiplicativity of ${\rm Tr}_{B_2}$, it follows that the value of ${\rm Tr}_{B_2}$ at the $K_d$ component is nonzero. Since $s_2$ assigns $1$ to $x_1,\ldots,x_d,w_1$, Lem.\ref{lem:K_d} that we just showed applies, and so the value of ${\rm Tr}_{B_2}$ at this component is nonzero iff $s_2$ assigns $1$ to all $w_2,\ldots,w_d$ and $2$ to $y_1$, in which case the value is $1$. By compatibility, $s_1$ assigns $1$ to all $w_2,\ldots,w_d$. Since $s_1$ also assigns $1$ to $z_2,\ldots,z_d$, the induction hypothesis (of our proof of Lem.\ref{lem:P_d}) applies to the $P_{d-1}$ component of $P_d \cap B_1$, hence the value of ${\rm Tr}_{B_1}$ at this component is nonzero iff $s_1$ assigns $2$ to all $u_2,\ldots,u_d$, in which case the value is $1$. By compatibility, $s_2$ assigns $2$ to all $u_2,\ldots,u_d$. Each edge component of $(P_d\cap B_2,s_2)$, connecting $u_i$ and $y_i$, falls into Prop.\ref{prop:biangle_SL3_quantum_trace}(BT2-1), hence it has nonzero ${\rm Tr}_{B_2}$ value iff $s_2$ assigns same value to $u_i$ and $y_i$, so $s_2(y_i)=2$ for all $i=2,\ldots,d$. To summarize, there is only one pair of $s_1,s_2$ contributing to the sum in eq.\eqref{eq:Tr_P_d}, which assign $2$ to all $y_1,\ldots,y_d$, and the corresponding summand in the sum in eq.\eqref{eq:Tr_P_d} is $1$. This finishes the proof of Lem.\ref{lem:P_d}. \qed 

\vs

We now prove Prop.\ref{prop:highest_term_of_SL3_classical_trace_for_a_triangle}.

\vs

{\it Proof of Prop.\ref{prop:highest_term_of_SL3_classical_trace_for_a_triangle}.} By the fact that ${\rm Tr}_\Delta$ is a ring homomorphism and from the additivity of the tropical coordinates ${\rm a}_v$ (Lem.\ref{lem:additivity_of_coordinates}), it suffices to prove the statement for each canonical ${\rm SL}_3$-web $W$ having a single component. Denote by $e_1,e_2,e_3$ the sides of the triangle $\wh{t}$ appearing clockwise this order along $\partial \wh{t}$. Denote the nodes of $Q_\Delta$ by $v_{e_\alpha,1}$, $v_{e_\alpha,2}$, $v_t$ as in Def.\ref{def:Fock-Goncharov_algebra_quantum}.

\vs

Suppose that $W$ is a left turn corner arc in $t$, as in Thm.\ref{thm:SL3_quantum_trace_map}\redfix{(QT2-1)}; so 
$$
{\rm Tr}_\Delta ([W,s]) = ({\bf M}^{\rm in}_{t,\alpha})_{s(x),s(x)} \, ({\bf M}^{\rm left}(Z_{v_t}))_{s(x),s(y)} \, ({\bf M}^{\rm out}_{t,\alpha+1})_{s(y),s(y)},
$$
where $\alpha,x,y$ are as in Thm.\ref{thm:SL3_quantum_trace_map}\redfix{(QT2-1)}. We used the fact that ${\bf M}^{\rm in}_{t,\alpha}$ and ${\bf M}^{\rm out}_{t,\alpha+1}$ are diagonal matrices. So, ${\bf M}^{\rm in}_{t,\alpha}$ involves variables $Z_{v_{e_\alpha,1}}$, $Z_{v_{e_\alpha,2}}$ but no others, ${\bf M}^{\rm left}(Z_{v_t})$ involves $Z_{v_t}$ but no others, and ${\bf M}^{\rm out}_{t,\alpha+1}$ involves variables $Z_{v_{e_{\alpha+1},1}}$, $Z_{v_{e_{\alpha+1},2}}$ but no others. In view of eq.\eqref{eq:quantum_edge_matrix}, eq.\eqref{eq:quantum_turn_matrix} and eq.\eqref{eq:quantum_in_and_out},
$$
{\rm Tr}_\Delta([W,{\bf 1}_W]) = 
Z_{v_{e_\alpha,2}} Z_{v_{e_\alpha,1}}^2 Z_{v_t}^2 
Z_{v_{e_{\alpha+1},1}} Z_{v_{e_{\alpha+1},2}}^2.
$$
In view of eq.\eqref{eq:tropical_coordinates_for_a_left_turn_W}, item (HT1) is satisfied. By inspection of the monodromy matrices ${\bf M}^{\rm in}_{t,\alpha}$, ${\bf M}^{\rm out}_{t,\alpha+1}$ and ${\bf M}^{\rm left}(Z_{v_t})$, it follows that this Laurent monomial indeed has higher or equal partial ordering than any other Laurent monomials appearing in ${\rm Tr}_\Delta([W,s])$. Also, if $s(x)\neq 1$, then $({\bf M}^{\rm in}_{t,\alpha})_{s(x),s(x)}$ has strictly lower partial ordering than $({\bf M}^{\rm in}_{t,\alpha})_{1,1} = Z_{v_{e_\alpha,2}} Z_{v_{e_\alpha,1}}^2$. If $s(y) \neq 1$, then $({\bf M}^{\rm out}_{t,\alpha+1})_{s(y),s(y)}$ has strictly lower partial ordering than $({\bf M}^{\rm out}_{t,\alpha+1})_{1,1} = Z_{v_{e_{\alpha+1},1}} Z_{v_{e_{\alpha+1},2}}^2$. Thus (HT2) is satisfied.

\vs

When $W$ is a right turn corner arc as in Thm.\ref{thm:SL3_quantum_trace_map}\redfix{(QT2-2)}, the proof goes completely parallel. We just have to check (HT1) precisely. Indeed, 
$$
{\rm Tr}_\Delta([W,{\bf 1}_W]) = 
Z_{v_{e_{\alpha+1},2}} Z_{v_{e_{\alpha+1},1}}^2
Z_{v_t}
Z_{v_{e_\alpha,1}} Z_{v_{e_\alpha,2}}^2,
$$
hence (HT1) is satisfied, in view of eq.\eqref{eq:tropical_coordinates_for_a_right_turn_W}.

\vs

\begin{figure}[htbp!]
\vspace{-2mm}
\begin{center}
\raisebox{-0.5\height}{\scalebox{0.8}{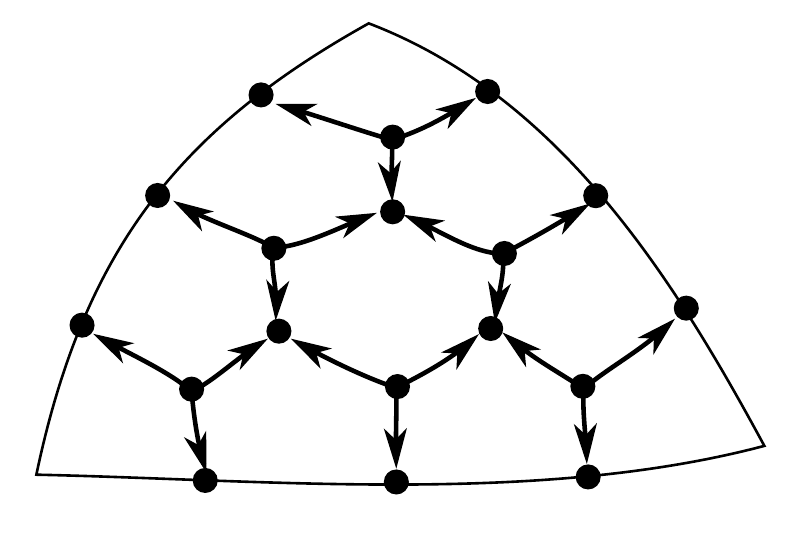}}
$\to$
\raisebox{-0.5\height}{\scalebox{0.8}{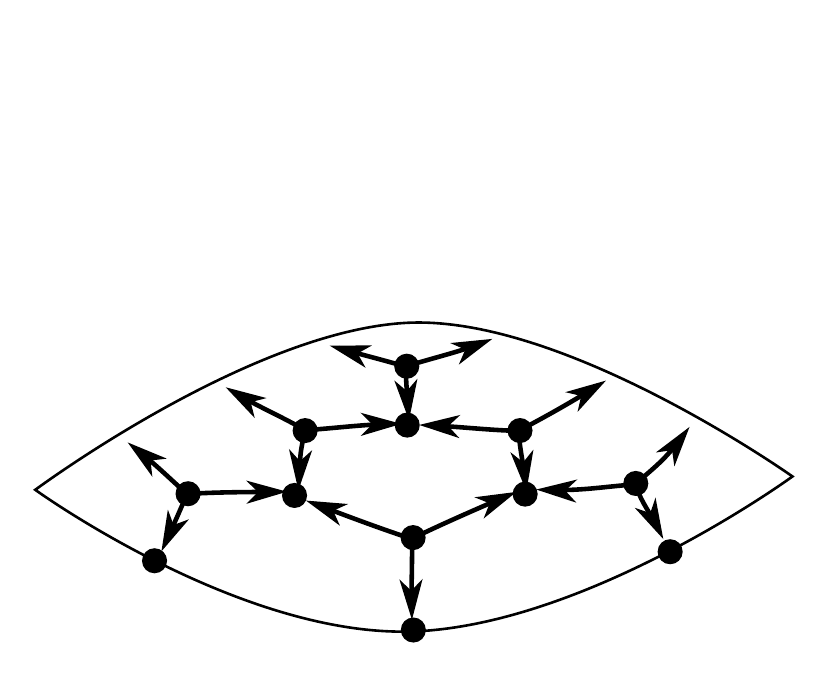}}
\end{center}
\vspace{-7mm}
\caption{Pushing a pyramid $H_d$ from triangle to biangle (for $d=3$)}
\vspace{-2mm}
\label{fig:pyramid_in_triangle_pushing}
\end{figure}

Now suppose that $W$ is a degree $d$ pyramid $H_d$ for some nonzero $d\in \mathbb{Z}$. We first present how to deal with $d>0$ in detail; the case $d<0$ is completely parallel, and we give a note the end of this proof of (HT1) about a difference between $d>0$ and $d<0$. To compute ${\rm Tr}_\Delta([W,s])$, we decompose $t$ into one triangle $\wh{t}$ and one biangle $B$, as done in Prop.\ref{prop:elementary_isotopy_invariance_3-way}. Let's say that the biangle is attached at the side $e_3$ of $t$. Let $e_1,e_2,e_3$ be sides of $\wh{t}$, and let $e_3'$ be the other side of $B$. Push all 3-valent vertices of $W$ to the biangle $B$ to form a stated ${\rm SL}_3$-web $(W',s')$ in $t$ as shown in Fig.\ref{fig:pyramid_in_triangle_pushing}, and apply the state-sum formula in eq.\eqref{eq:state-sum_formula} in Def.\ref{def:state-sum_trace_for_gool_position} to define $\wh{\rm Tr}_\Delta(W',s') \in \mathcal{Z}_\Delta$. By the isotopy invariance of the state-sum formula that we proved, we know ${\rm Tr}_\Delta([W,s]) = \wh{\rm Tr}_\Delta(W',s')$ (if one wants to write down a proof explicitly, one may want to consider a genuine split ideal triangulation of $\Delta$, which has three biangles). 

\vs

Denote the endpoints of $W$ by $x_1,\ldots,x_d$, $y_1,\ldots,y_d$, $z_1,\ldots,z_d$, appearing clockwise in this order along $\partial t$, where $x_i$'s are in $e_1$, $y_i$'s in $e_2$, and $z_i$'s in $e_3$; see Fig.\ref{fig:pyramid_in_triangle_pushing}. Inherit these labelings to $W'$, so that $z_1,\ldots,z_d$ lie in the outermost side $e_3'$ of the biangle $B$. Denote the junctures of $W'$ at the common arc $e_3$ of $\wh{t}$ and $B$ as $w_1,\ldots,w_d$, $u_1,\ldots,u_d$ as in Fig.\ref{fig:pyramid_in_triangle_pushing}. Then
\begin{align}
\label{eq:Tr_W_s}
{\rm Tr}_\Delta([W,s]) = \underset{s_1,s_2}{\textstyle \sum} {\rm Tr}_B([W'\cap B,s_1]) \, {\rm Tr}_{\wh{t}}([W'\cap \wh{t},s_2])
\end{align}
where the sum is over all states $s_1,s_2$ of the ${\rm SL}_3$-webs $W'\cap B$ in $B$ and $W'\cap \wh{t}$ in $\wh{t}$ compatible with $s$, in the sense as in Prop.\ref{prop:biangle_SL3_quantum_trace}(BT1). Note that $W'\cap \wh{t}$ has $2d$ components, which are corner arcs. Denote the component connecting $w_i$ and $x_i$ by $W_i^{\rm left}$, and the component connecting $u_i$ to $y_i$ by $W_i^{\rm right}$, for $i=1,\ldots,d$. So
\begin{align}
\label{eq:Tr_W_prime}
{\rm Tr}_{\wh{t}}([W'\cap \wh{t},s_2])
= {\textstyle \prod}_{i=1}^d {\rm Tr}_{\wh{t}}([W_i^{\rm left}, s_2]) {\rm Tr}_{\wh{t}}([W_i^{\rm right}, s_2]),
\end{align}
where the $s_2$'s appearing in the right hand side mean appropriate restrictions. By Thm.\ref{thm:SL3_quantum_trace_map}\redfix{(QT2-1)} and \redfix{(QT2-2)} we have
\begin{align}
\label{eq:Tr_W_i_left}
{\rm Tr}_{\wh{t}}([W_i^{\rm left}, s_2]) & = ({\bf M}^{\rm in}_{t,3})_{s_2(w_i),s_2(w_i)} \, ({\bf M}^{\rm left}(Z_{v_t}))_{s_2(w_i),s_2(x_i)} \, ({\bf M}^{\rm out}_{t,1})_{s_2(x_i),s_2(x_i)}, \\
\label{eq:Tr_W_i_right}
{\rm Tr}_{\wh{t}}([W_i^{\rm right}, s_2]) & = ({\bf M}^{\rm in}_{t,3})_{s_2(u_i),s_2(u_i)} \, ({\bf M}^{\rm right}(Z_{v_t}))_{s_2(u_i),s_2(y_i)} \, ({\bf M}^{\rm out}_{t,2})_{s_2(y_i),s_2(y_i)}.
\end{align}

\vs

We first compute ${\rm Tr}_\Delta([W,{\bf 1}_W])$. Let $s_1,s_2$ be states of $W'\cap B$ and $W'\cap \wh{t}$ compatible with ${\bf 1}_W$ and such that the corresponding summand in eq.\eqref{eq:Tr_W_s} is nonzero. In particular, $s_1$ assigns $1$ to $z_1,\ldots,z_d$, and $s_2$ assigns $1$ to $x_1,\ldots,x_d,y_1,\ldots,y_d$. By multiplicativity of ${\rm Tr}_{\wh{t}}$, the value under ${\rm Tr}_{\wh{t}}$ of each component of $(W'\cap \wh{t}, s_2)$, i.e. eq.\eqref{eq:Tr_W_i_left} and eq.\eqref{eq:Tr_W_i_right}, must be nonzero. Since $s_2(x_i)=1$ and ${\bf M}^{\rm left}(Z_{v_t})$ is upper triangular, it follows $s_2(w_i)=1$, for all $i=1,\ldots,d$. By compatibility, $s_1(w_i)=1$ for all $i=1,\ldots,d$. Since $W' \cap B$ is the ${\rm SL}_3$-web $P_d$, where $s_1$ assigns $1$ to $w_1,\ldots,w_d,z_1,\ldots,z_d$, Lem.\ref{lem:P_d} applies. So ${\rm Tr}_B([W'\cap B,s_1]) \neq 0$ implies $s_1(u_i)=2$ for all $i=1,\ldots,d$, in which case ${\rm Tr}_B([W'\cap B,s_1])=1$. By compatibility, $s_2(u_i)=2$ for all $i=1,\ldots,d$. So there is a unique such pair of states $s_1,s_2$. For this pair of states, we have
\begin{align}
\label{eq:HT1_proof1}
{\rm Tr}_{\wh{t}}([W_i^{\rm left},s_2]) = (Z_{v_{e_3,2}} Z_{v_{e_3,1}}^2) Z_{v_t}^2 (Z_{v_{e_1,1}} Z_{v_{e_1,2}}^2), ~
{\rm Tr}_{\wh{t}}([W_i^{\rm right},s_2]) = (Z_{v_{e_3,2}} Z_{v_{e_3,1}}^{-1}) Z_{v_t} (Z_{v_{e_2,1}} Z_{v_{e_2,2}}^2) \hspace{-10mm}\,
\end{align}
in view of eq.\eqref{eq:quantum_edge_matrix} and \eqref{eq:quantum_in_and_out}. Thus (for $d>0$)
\begin{align*}
{\rm Tr}_\Delta([W,{\bf 1}_W]) & = 
1\cdot \left( (Z_{v_{e_3,2}} Z_{v_{e_3,1}}^2) Z_{v_t}^{2} (Z_{v_{e_1,1}} Z_{v_{e_1,2}}^2) 
(Z_{v_{e_3,2}} Z_{v_{e_3,1}}^{-1}) Z_{v_t} (Z_{v_{e_2,1}} Z_{v_{e_2,2}}^2) \right)^d \\
& = Z_{v_t}^{3d} \, {\textstyle \prod}_{\alpha=1}^3 Z_{v_{e_\alpha,1}}^d Z_{v_{e_\alpha,2}}^{2d}.
\end{align*}
One can easily verify that these powers are indeed $3$ times the tropical coordinates of the degree $d$ pyramid $W=H_d$ (eq.\eqref{eq:tropical_coordinates_for_Hd_positive}), hence (HT1) is satisfied. In fact, this last step becomes slightly different for the case when $d<0$. Reversing all the arrows, left turns become right turns and vice versa, and this time using the fact that ${\bf M}^{\rm right}(Z_{v_t})$ is lower triangular, one observes that for the only contributing state-summand for ${\rm Tr}_\Delta([W,{\bf 1}_W])$, the right turns have state $1$ at the both endpoints lying on edge $e_1$ and $e_3$, while the left turns have states $1$ at the initial endpoints lying on edge $e_2$ and $2$ at the terminal endpoint lying on edge $e_3$. Hence eq.\eqref{eq:HT1_proof1} should become
$$
{\rm Tr}_{\wh{t}}([W^{\rm right}_i,s_2]) = (Z_{v_{e_1,2}} Z_{v_{e_1,1}}^2) Z_{v_t} (Z_{v_{e_3,1}} Z_{v_{e_3,2}}^2),
~
{\rm Tr}_{\wh{t}}([W^{\rm left}_i,s_2]) = (Z_{v_{e_2,2}} Z_{v_{e_2,1}}^2) (Z_{v_t}^2 + Z_{v_t}^{-1}) (Z_{v_{e_3,1}} Z_{v_{e_3,2}}^{-1})
$$
in view of eq.\eqref{eq:quantum_edge_matrix} and \eqref{eq:quantum_in_and_out}. Thus (for $d<0$)
\begin{align*}
{\rm Tr}_\Delta([W,{\bf 1}_W]) & = 
1\cdot \left( (Z_{v_{e_1,2}} Z_{v_{e_1,1}}^2) Z_{v_t}  (Z_{v_{e_3,1}} Z_{v_{e_3,2}}^2) 
(Z_{v_{e_2,2}} Z_{v_{e_2,1}}^{-1}) (Z_{v_t}^2+Z_{v_t}^{-1}) (Z_{v_{e_3,1}} Z_{v_{e_3,2}}^{-1}) \right)^{-d} \\
& = (Z_{v_t}^3+1)^{-d} \, {\textstyle \prod}_{\alpha=1}^3 Z_{v_{e_\alpha,1}}^{-2d} Z_{v_{e_\alpha,2}}^{-d}.
\end{align*}
One can verify that in the unique highest term $Z_{v_t}^{-3d} \, {\textstyle \prod}_{\alpha=1}^3 Z_{v_{e_\alpha,1}}^{-2d} Z_{v_{e_\alpha,2}}^{-d}$ for ${\rm Tr}_\Delta([W,{\bf 1}_W])$, the powers of generators are indeed $3$ times the tropical coordinates of the degree $d$ pyramid $W=H_d$ (eq.\eqref{eq:tropical_coordinates_for_Hd_positive}), hence (HT1) is satisfied.

\vs

Now let's prove (HT2). Let $s$ be any state of $W$, and consider eq.\eqref{eq:Tr_W_s}. Let $s_1,s_2$ be any pairs of states of $W'\cap B$ and $W'\cap \wh{t}$ compatible with $s$. The biangle factor ${\rm Tr}_B([W'\cap B,s_1])$ is an integer, so it does not involve any generator of $\mathcal{Z}_t$. The triangle factor is given by eq.\eqref{eq:Tr_W_prime}, \eqref{eq:Tr_W_i_left} and \eqref{eq:Tr_W_i_right}. The only places where the variable $Z_{v_t}$ appears are ${\bf M}^{\rm left}(Z_{v_t})$ of eq.\eqref{eq:Tr_W_i_left} and ${\bf M}^{\rm right}(Z_{v_t})$ of eq.\eqref{eq:Tr_W_i_right}. In view of eq.\eqref{eq:quantum_turn_matrix}, the highest power of $Z_{v_t}$ in ${\bf M}^{\rm left}(Z_{v_t})$ is $Z_{v_t}^2$ and that in ${\bf M}^{\rm right}(Z_{v_t})$ is $Z_{v_t}$. Hence, the maximum possible power of $Z_{v_t}$ that can appear in a summand of the sum expression for ${\rm Tr}_\Delta([W,s])$ in eq.\eqref{eq:Tr_W_s} is $(Z_{v_t}^2)^d (Z_{v_t})^d = Z_{v_t}^{3d}$, which is the power of $Z_{v_t}$ that does appear in ${\rm Tr}_\Delta([W,{\bf 1}_W])$. Now, among eq.\eqref{eq:Tr_W_i_left} and \eqref{eq:Tr_W_i_right}, the only place where the variables $Z_{v_{e_1,1}}$ or $Z_{v_{e_1,2}}$ is ${\bf M}^{\rm out}_{t,1}$ of eq.\eqref{eq:Tr_W_i_left}. In view of eq.\eqref{eq:quantum_edge_matrix} and eq.\eqref{eq:quantum_in_and_out}, the highest monomial in $Z_{v_{e_1,1}}$ and $Z_{v_{e_1,2}}$ appearing in ${\bf M}^{\rm out}_{t,1}$ is the $(1,1)$-entry $({\bf M}^{\rm out}_{t,1})_{1,1} = Z_{v_{e_1,1}} Z_{v_{e_1,2}}^2$. So the highest possible Laurent monomial in $Z_{v_{e_1,1}}$ and $Z_{v_{e_1,2}}$ that can appear in the sum expression for ${\rm Tr}_\Delta([W,s])$ is $(Z_{v_{e_1,1}} Z_{v_{e_1,2}}^2)^d$, which may happen when $s$ assigns values $1$ to all $x_1,\ldots,x_d$, e.g. in ${\rm Tr}_\Delta([W,{\bf 1}_W])$. If $s$ does not assign $1$ to some $x_i$, then one notes that $({\bf M}^{\rm out}_{t,1})_{s_2(x_i),s_2(x_i)} = ({\bf M}^{\rm out}_{t,1})_{s(x_i),s(x_i)}$ in eq.\eqref{eq:Tr_W_i_left} is either $Z_{v_{e_1,1}} Z_{v_{e_1,2}}^{-1}$ or $Z_{v_{e_1,1}}^{-2} Z_{v_{e_1,2}}^{-1}$, hence is of strictly lower partial order than $({\bf M}^{\rm out}_{t,1})_{1,1}=Z_{v_{e_1,1}} Z_{v_{e_1,2}}^2$. 

\vs

Now, go back to the beginning, before we split $t$ into the triangle $\wh{t}$ and a biangle $B$ at the side $e_3$. This time, decompose $t$ into a triangle and a biangle where the biangle is at a different side than $e_3$. Apply the same arguments as we have seen so far, which is possible because the ${\rm SL}_3$-web $W=H_d$ has cyclic symmetry. Then we obtain similar results about the variables $Z_{v_{e_2,1}}$ and $Z_{v_{e_2,2}}$ lying in the side $e_2$, and also for the variables $Z_{v_{e_3,1}}$ and $Z_{v_{e_3,2}}$. So, for each $i=1,2,3$, the highest possible Laurent monomial in $Z_{v_{e_\alpha,1}}$ and $Z_{v_{e_\alpha,2}}$ that can appear in the sum expression for ${\rm Tr}_\Delta([W,s])$ is $(Z_{v_{e_\alpha,1}} Z_{v_{e_\alpha,2}}^2)^d$, which may happen when $s$ assigns values $1$ to all endpoints of $W$ lying in $e_\alpha$. If $s$ does not assign $1$ to some endpoint in $e_\alpha$, then one notes that the Laurent monomials in $Z_{v_{e_\alpha,1}}$ and $Z_{v_{e_\alpha,2}}$ appearing in the summands of eq.\eqref{eq:Tr_W_s} have strictly lower partial order than $(Z_{v_{e_\alpha,1}} Z_{v_{e_\alpha,2}}^2)^d$. This completes the proof of (HT2). \qed \qquad {\it End of proof of Prop.\ref{prop:highest_term_of_SL3_classical_trace_for_a_triangle}.}

\vs

Before proceeding to a proof of the general highest-term statement, Prop.\ref{prop:highest_term_of_SL3_classical_trace}, we need one more easy lemma about biangles.
\begin{lemma}
\label{lem:crossbar_web_value}
Let $W$ be a crossbar ${\rm SL}_3$-web in a biangle $B$ (Def.\ref{def:crossbar}), and let ${\bf 1}_W$ be the state of $W$ assigning the value $1\in \{1,2,3\}$ to all the endpoints of $W$. Then ${\rm Tr}_B([W,{\bf 1}_W])=1$.
\end{lemma}

{\it Proof of Lem.\ref{lem:crossbar_web_value}.} Note that $W$ can be decomposed as composition of crossbar webs having exactly one crossbar, i.e. contains exactly two internal $3$-valent vertices. That is, there is a finite collection of ideal arcs $e_1,\ldots,e_n$ dividing $B$ into biangles $B_1,\ldots,B_{n+1}$, as in the proof of Lem.\ref{lem:charge_conservation}, so that  $W_i = W\cap B_i$ is a crossbar web in $B_i$ with two internal 3-valent vertices. We have the state-sum formula eq.\eqref{eq:Tr_B_W_s}, with $s = {\bf 1}_W$. Denote the endpoints of $W_i$ lying in one side of $B_i$ by $x_1,\ldots,x_n$, and the endpoints of $W_i$ lying on the other side by $y_1,\ldots,y_n$, so that $x_j$ is connected to $y_j$ by an edge of $W_i$ for all $j=1,\ldots,n$ except for some two adjacent $j$'s. 

\vs

Assume that $s_i$ is a state of $W_i$ assigning $1$ to all $x_1,\ldots,x_n$. Assume ${\rm Tr}_{B_i}([W_i,s_i])\neq 0$. By multiplicativity of ${\rm Tr}_{B_i}$, the value under ${\rm Tr}_{B_i}$ of each component of $(W_i,s_i)$ must be nonzero. The edge component, connecting $x_k$ and $y_k$, falls into Prop.\ref{prop:biangle_SL3_quantum_trace}(BT2-1), hence it follows $s_i(y_k)=1$. For the single-crossbar component $(L,s_i|_{\partial L})$, connecting $x_j,x_{j+1}$ and $y_j,y_{j+1}$, Lem.\ref{lem:value_of_I-webs_in_biangle} applies, telling us that ${\rm Tr}_{B_i} ([L,s_i|_{\partial L}])\neq 0$ iff $s$ assigns $1$ to $y_j$ and $y_{j+1}$, in which case ${\rm Tr}_{B_i} ([L,s_i|_{\partial L}])=1$.

\vs

We go back to the state-sum formula eq.\eqref{eq:Tr_B_W_s} with $s = {\bf 1}_W$. Let $J$ be any juncture-state of $W$ compatible with $s = {\bf 1}_W$ whose corresponding summand in the sum in eq.\eqref{eq:Tr_B_W_s} is nonzero. Then ${\rm Tr}_{B_i}(W_i, J|_{\partial W_i}) \neq 0$ for all $i=1,\ldots,n+1$. Look at the first biangle $B_1$ whose one side equals one side of $B$; so $J|_{\partial W_1}$ assigns $1$ to all endpoints lying on this side of $B_1$. Applying the above observation we made, ${\rm Tr}_{B_i}(W_i, J|_{\partial W_i}) \neq 0$ iff $J|_{\partial W_i}$ assigns $1$ to all endpoints of $W_i$, in which case ${\rm Tr}_{B_i}(W_i, J|_{\partial W_i})=1$. Then we go to biangle $B_2$, where now we know $J|_{\partial W_2}$ assigns $1$ to all endpoints lying in one side. Apply the above observation. Repeating this till the end, we deduce that $J$ must assign $1$ to all junctures, in which case the corresponding summand in eq.\eqref{eq:Tr_B_W_s} is $1$. So ${\rm Tr}_B([W,{\bf 1}_W])=1$ as desired. \qed

\vs

We finally can provide a proof of Prop.\ref{prop:highest_term_of_SL3_classical_trace}.

\vs

{\it Proof of Prop.\ref{prop:highest_term_of_SL3_classical_trace}.} Let $W$ be a reduced non-elliptic ${\rm SL}_3$-web in a generalized marked surface $\frak{S}$ in a canonical position with respect to a split ideal triangulation $\wh{\Delta}$ of $\frak{S}$. We use the state-sum formula eq.\eqref{eq:state-sum_formula} of Def.\ref{def:state-sum_trace_for_gool_position}. \redfix{Note that in Def.\ref{def:state-sum_trace_for_gool_position}, we required that the ${\rm SL}_3$-web should be in a gool position. However, now using Thm.\ref{thm:SL3_quantum_trace_map}(QT1) and the results of \S\ref{subsec:isotopy_invariance}, one can easily deduce that the state-sum formula as in eq.\eqref{eq:state-sum_formula} actually holds for any ${\rm SL}_3$-web $W$ in $\frak{S}\times {\bf I}$ provided that, for each triangle $\wh{t}$ and each biangle $B$ of $\wh{\Delta}$, $W\cap (\wh{t} \times {\bf I})$ and $W\cap (B \times {\bf I})$ are well-defined ${\rm SL}_3$-webs in $\wh{t} \times {\bf I}$ and $B\times {\bf I}$. Of course, the case of ${\rm SL}_3$ classical trace map is easier, because there is no need to consider the elevations in ${\bf I}$, and $W$ is an ${\rm SL}_3$-web living in the surface $\frak{S}$.} Let $J_0$ be the $\wh{\Delta}$-juncture-state of $W$ that assigns the value $1\in\{1,2,3\}$ to all junctures. Since each $W\cap B$ is a crossbar ${\rm SL}_3$-web in a biangle $B$, by Lem.\ref{lem:crossbar_web_value} the biangle factor ${\rm Tr}_B([W\cap B,(J_0)_B])$ equals $1$. For each triangle $t$ of $\Delta$, by Prop.\ref{prop:highest_term_of_SL3_classical_trace_for_a_triangle}(HT1) the triangle factor $\wh{\rm Tr}_t(W\cap \wh{t}, (J_0)_t)$ equals $\prod_{v\in \mathcal{V}(Q_\Delta) \cap t} Z_{t;v}^{3{\rm a}_v(W \cap \wh{t}\,\,)} \in \mathcal{Z}_t$\redfix{; here $\wh{\rm Tr}_t$ means $\wh{\rm Tr}^\omega_t$ as in eq.\eqref{eq:triangle_factor}, with $\omega^{1/2}=1$}. For the node $v_t$ of $Q_\Delta$ lying in the interior of $t$, we have $Z_{t;v_t}^{3{\rm a}_{v_t}(W\cap \wh{t}\,\,)} = Z_{v_t}^{3{\rm a}_{v_t}(W)} \in \mathcal{Z}_\Delta$. Now let $v$ be a node of $Q_\Delta$ lying in an internal arc of $\Delta$, say a common side of triangles $t$ and $r$. Let $B$ be the biangle in between the triangles $\wh{t}$ and $\wh{r}$ of $\wh{\Delta}$. By the well-definedness of the tropical coordinates at arcs of $\Delta$, note that ${\rm a}_v(W\cap \wh{t}\,\,) = {\rm a}_v(W\cap \wh{r}) = {\rm a}_v(W)$. Then note $Z_{t;v}^{3{\rm a}_v(W\cap \wh{t}\,\,)} Z_{r;v}^{3{\rm a}_v(W\cap \wh{r})} = Z_v^{3{\rm a}_v(W)} \in \mathcal{Z}_\Delta$. Therefore, the summand of eq.\eqref{eq:state-sum_formula} corresponding to $J_0$ exactly equals $\prod_{v\in \mathcal{V}(Q_\Delta)} Z_v^{3{\rm a}_v(W)}$. Now, let $J$ be any $\wh{\Delta}$-juncture-state of $W$ different from $J_0$. Note ${\rm Tr}_B([W\cap B,J_B])$ is an integer, hence does not involve any generator of $\mathcal{Z}_\Delta$. By Prop.\ref{prop:highest_term_of_SL3_classical_trace_for_a_triangle}(HT2), for each triangle $t$ of $\Delta$, the triangle factor $\wh{\rm Tr}_t(W\cap \wh{t}, J_t) \in \mathcal{Z}_t$ only involves Laurent monomials $\mathcal{Z}_t$ having lower or equal partial order than $\wh{\rm Tr}_t(W\cap \wh{t}, (J_0)_t)$ which is a single monomial in $\mathcal{Z}_t$. Also by Prop.\ref{prop:highest_term_of_SL3_classical_trace_for_a_triangle}(HT2), there exists a triangle $t$ such that the triangle factor $\wh{\rm Tr}_t(W\cap \wh{t}, J_t) \in \mathcal{Z}_t$ only involves Laurent monomials of $\mathcal{Z}_t$ having strictly lower partial order than $\wh{\rm Tr}_t(W\cap \wh{t}, (J_0)_t)$. This finishes the proof of Prop.\ref{prop:highest_term_of_SL3_classical_trace}. \qed

\subsection{The relationship with the basic semi-regular functions}

In order to prove Prop.\ref{prop:highest_term} and Prop.\ref{prop:congruence-compatibility_of_terms_of_each_basic_regular_function}, we should translate the results from the previous subsection about the ${\rm SL}_3$ classical (state-sum) trace ${\rm Tr}_\Delta$ into those of basic semi-regular functions $\mathbb{I}^+_{{\rm PGL}_3}(\ell) \in C^\infty(\mathscr{X}^+_{{\rm PGL}_3,\frak{S}})$.

\begin{definition}
\label{def:iota_Delta}
For each ideal triangulation $\Delta$ of a triangulable punctured surface $\frak{S}$, define 
$$
\iota_\Delta : \mathcal{Z}_\Delta \to C^\infty(\mathscr{X}^+_{{\rm PGL}_3, \frak{S}})
$$
as the unique ring homomorphism sending $Z_v^{\pm 1} \in \mathcal{Z}_\Delta$ to $X_v^{\pm 1/3} \in C^\infty(\mathscr{X}^+_{{\rm PGL}_3,\frak{S}})$, $\forall v\in \mathcal{V}(Q_\Delta)$.
\end{definition}

\begin{proposition}[${\rm SL}_3$ classical trace and basic semi-regular function]
\label{prop:SL3_classical_trace_and_I_PGL3}
Let $\frak{S}$ be a triangulable punctured surface, $\Delta$ an ideal triangulation of $\frak{S}$, and $\wh{\Delta}$ a split ideal triangulation of $\Delta$. Let $\ell \in \mathscr{A}_{\rm L}(\frak{S};\mathbb{Z})$ be an ${\rm SL}_3$-lamination in $\frak{S}$ that can be represented as a (reduced) non-elliptic ${\rm SL}_3$-web  $W$ in $\frak{S}$ such that
\begin{enumerate}
\itemsep0em
\item[\rm (E1)] $W$ contains no peripheral loops, 

\item[\rm (E2)] All weights (of components of $W$) are $1$.
\end{enumerate}
Then
\begin{align}
\label{eq:SL3_classical_trace_and_I_PGL3}
\mathbb{I}^+_{{\rm PGL}_3}(\ell) = \iota_\Delta {\rm Tr}_{\Delta;\frak{S}}([W;{\O}])
\end{align}
\end{proposition}

{\it Proof of Prop.\ref{prop:SL3_classical_trace_and_I_PGL3}.} First, not precisely being fit to the current situation, assume that $W$ is an ${\rm SL}_3$-web consisting of a single oriented \redfix{non-contractible} loop, say $\gamma$, which is not necessarily simple. By applying an isotopy, we may assume that $\gamma$ meets $\wh{\Delta}$ transversally in a minimal possible number of points. We apply the construction in \S\ref{subsec:lifting_PGL3_to_SL3} of the monodromy matrix for $\gamma$. The $\wh{\Delta}$-junctures of $\gamma$, i.e. the points of $\wh{\Delta} \cap \gamma$, divide $\gamma$ into segments $\gamma_1,\ldots,\gamma_N$, so that $\gamma$ is the concatenation $\gamma = \gamma_1.\gamma_2.\cdots.\gamma_N$. The $\gamma_i$ in a triangle of $\wh{\Delta}$ work as a triangle segment, and $\gamma_i$ in a biangle of $\wh{\Delta}$ work as a juncture segment, so that $f^+_\gamma = {\rm tr}({\bf M}_{\gamma_1} \cdots {\bf M}_{\gamma_N})$. Also, the sequence $\gamma_1,\ldots,\gamma_N$ alternates between triangle segments and juncture segments. Let $\vec{\alpha} = (\alpha_1,\alpha_2,\ldots,\alpha_N) \in \{1,2,3\}^N$. Denote by $\wh{\bf M}_{\gamma_i}$ the matrix defined by the same formula as ${\bf M}_{\gamma_i}$, where the entries are thought of as elements of $\mathcal{Z}_\Delta$, so that ${\bf M}_{\gamma_i} = \iota_\Delta \wh{\bf M}_{\gamma_i}$, where $\iota_\Delta$ applied to a matrix means $\iota_\Delta$ applied to each entry. Denoting by $(\wh{\bf M}_{\gamma_i})_{\alpha_i,\alpha_{i+1}}$ the $(\alpha_i,\alpha_{i+1})$-th entry of $\wh{\bf M}_{\gamma_i}$, as usual, we have
$$
f^+_\gamma = \iota_\Delta {\rm tr}(\wh{\bf M}_{\gamma_1}\cdots \wh{\bf M}_{\gamma_N})
= \iota_\Delta {\textstyle \sum}_{\vec{\alpha} \in \{1,2,3\}^N} {\textstyle \prod}_{i=1}^N (\wh{\bf M}_{\gamma_i})_{\alpha_i,\alpha_{i+1}}
$$
where $\alpha_{N+1}:=\alpha_1$. View $\alpha_i$ as being associated to the $\wh{\Delta}$-juncture of $\gamma$ that is the initial point of $\gamma_i$ (which is the terminal point of $\gamma_{i-1}$; let $\gamma_0:=\gamma_N$). So $\vec{\alpha}$ can be viewed as a $\wh{\Delta}$-juncture-state $J=J^{\vec{\alpha}}$ of $\gamma$, and the above sum is over all $\wh{\Delta}$-juncture-states $J$. For each juncture segment $\gamma_i$, note $\wh{\bf M}_{\gamma_i}$ is diagonal, so $(\wh{\bf M}_{\gamma_i})_{\alpha_i,\alpha_{i+1}}=0$ unless $\alpha_i = \alpha_{i+1}$. Hence only the $\wh{\Delta}$-juncture-states $J^{\vec{\alpha}}$ that are {\em biangle-coherent} may contribute to the above sum, where we say that $\wh{\Delta}$-juncture-state is biangle-coherent if it assigns the same value to the two endpoints of each segment of $\gamma$ living in a biangle. For each juncture segment $\gamma_i$ as in Fig.\ref{fig:juncture_segment}, we have $\wh{\bf M}_{\gamma_i}={\rm diag}(Z_1 Z_2^2, Z_1 Z_2^{-1}, Z_1^{-2} Z_2^{-1})$. Suppose the initial and terminal points of $\gamma_i$ \redfix{live} in triangles $\wh{t}$ and $\wh{r}$ of $\wh{\Delta}$ corresponding to triangles $t$ and $r$ of $\Delta$. Define
$$
\wh{\bf M}_{\gamma_i}^{\rm ini} := {\rm diag}(Z_{t,1} Z_{t,2}^2, Z_{t,1} Z_{t,2}^{-1}, Z_{t,1}^{-2} Z_{t,2}^{-1}), \quad
\wh{\bf M}_{\gamma_i}^{\rm ter} := {\rm diag}(Z_{r,1} Z_{r,2}^2, Z_{r,1} Z_{r,2}^{-1}, Z_{r,1}^{-2} Z_{r,2}^{-1}),
$$
so that $\wh{\bf M}_{\gamma_i} = \wh{\bf M}_{\gamma_i}^{\rm ini} \wh{\bf M}_{\gamma_i}^{\rm ter}$ and $(\wh{\bf M}_{\gamma_i})_{\alpha_i,\alpha_i} = (\wh{\bf M}_{\gamma_i}^{\rm ini})_{\alpha_i,\alpha_i} (\wh{\bf M}_{\gamma_i}^{\rm ter})_{\alpha_i,\alpha_i}$. Meanwhile, for each triangle segment $\gamma_j$, living in triangle $\wh{t}$ (or $t$), one observes from Def.\ref{def:state-sum_trace_for_gool_position} (eq.\eqref{eq:value_of_stated_component}) that
$$
\wh{\rm Tr}_t(\gamma_j,  (J^{\vec{\alpha}})|_{\partial\gamma_j}) 
= ( \wh{\bf M}_{\gamma_{j-1}}^{\rm ter} \wh{\bf M}_{\gamma_j} \wh{\bf M}_{\gamma_{j+1}}^{\rm ini} )_{\alpha_j, \alpha_{j+1}} ~\in~ \mathcal{Z}_t.
$$
Now, assuming that $\gamma_1$ is a triangle segment (so that $\gamma_N$ is a juncture segment), for each biangle-coherent $\wh{\Delta}$-juncture-state $J^{\vec{\alpha}}$ (so that $\alpha_2=\alpha_3$, $\alpha_4=\alpha_5$, \ldots, $\alpha_{N-2}=\alpha_{N-1}$, $\alpha_N=\alpha_1$), observe
\begin{align*}
{\textstyle \prod}_{i=1}^N (\wh{\bf M}_{\gamma_i})_{\alpha_i,\alpha_{i+1}} & = (\wh{\bf M}_{\gamma_1})_{\alpha_1,\alpha_2} (\wh{\bf M}_{\gamma_2})_{\alpha_2,\alpha_3} (\wh{\bf M}_{\gamma_3})_{\alpha_3,\alpha_4} \cdots (\wh{\bf M}_{\gamma_N})_{\alpha_N,\alpha_1} \\
& = (\wh{\bf M}_{\gamma_1})_{\alpha_1,\alpha_2} (\wh{\bf M}_{\gamma_2}^{\rm ini} \, \wh{\bf M}_{\gamma_2}^{\rm ter})_{\alpha_2,\alpha_3} (\wh{\bf M}_{\gamma_3})_{\alpha_3,\alpha_4}  \cdots (\wh{\bf M}_{\gamma_N}^{\rm ini} \, \wh{\bf M}_{\gamma_N}^{\rm ter})_{\alpha_N,\alpha_1} \\
& = (\wh{\bf M}_{\gamma_1})_{\alpha_1,\alpha_2} (\wh{\bf M}_{\gamma_2}^{\rm ini})_{\alpha_2,\alpha_3} \, (\wh{\bf M}_{\gamma_2}^{\rm ter})_{\alpha_2,\alpha_3} (\wh{\bf M}_{\gamma_3})_{\alpha_3,\alpha_4}  \cdots (\wh{\bf M}_{\gamma_N}^{\rm ini})_{\alpha_N,\alpha_1} \, (\wh{\bf M}_{\gamma_N}^{\rm ter})_{\alpha_N,\alpha_1} \\
& = {\textstyle \prod}_{j=1,3,5,\ldots,N-1} (\wh{\bf M}_{\redfix{\gamma_{j-1}}}^{\rm ter})_{\alpha_{j-1},\alpha_j} (\wh{\bf M}_{\gamma_j})_{\alpha_j,\alpha_{j+1}} (\wh{\bf M}_{\redfix{\gamma_{j+1}}}^{\rm ini})_{\alpha_{j+1},\alpha_{j+2}} \\
& = {\textstyle \prod}_{j=1,3,5,\ldots,N-1} (\wh{\bf M}_{\redfix{\gamma_{j-1}}}^{\rm ter} \wh{\bf M}_{\gamma_j} \wh{\bf M}_{\redfix{\gamma_{j+1}}}^{\rm ini})_{\alpha_j,\alpha_{j+1}} \\
& = {\textstyle \prod}_t {\textstyle \prod}_{{\rm triangle~segments}~\gamma_j~{\rm in~}t} \, \wh{\rm Tr}_t(\gamma_j, (J^{\vec{\alpha}})_{\partial \gamma_j}) \\
& = {\textstyle \prod}_t \wh{\rm Tr}_t(\gamma\cap \wh{t}, (J^{\vec{\alpha}})|_{\partial (\gamma \cap \wh{t})} ).
\end{align*}
Meanwhile, consider the state-sum formula in eq.\eqref{eq:state-sum_formula}, which we can apply because $W=\gamma$ is in a gool position with respect to $\wh{\Delta}$ (Def.\ref{def:good_position}). The ${\rm SL}_3$-web $W\cap B$ in each biangle $B$ consists of edge components of type as in Prop.\ref{prop:biangle_SL3_quantum_trace}(BT2-1), hence it follows that each $\wh{\Delta}$-juncture-state $J$ whose corresponding summand in the sum in eq.\eqref{eq:state-sum_formula} is nonzero is biangle-coherent, in the sense as defined above, and for such $J$'s the biangle factors ${\rm Tr}_B([W\cap B,J_B])$ are $1$. So it follows that
$$
f^+_\gamma = \iota_\Delta {\textstyle \sum}_{\vec{\alpha}} {\textstyle \prod}_t \wh{\rm Tr}_t(\gamma\cap \wh{t}, (J^{\vec{\alpha}})|_{\partial (\gamma \cap \wh{t})} )
= \iota_\Delta {\textstyle \sum}_J {\textstyle \prod}_t \wh{\rm Tr}_t(W \cap \wh{t}, J_t) = \iota_\Delta {\rm Tr}_{\Delta;\frak{S}}([W,{\O}]),
$$
where the middle equality holds because both are sums over biangle-coherent juncture-states. 

\vs

Now, coming back \redfix{to} the original situation of the problem, let $W$ be a reduced non-elliptic ${\rm SL}_3$-web in $\frak{S}$ representing $\ell \in \mathscr{A}_{\rm L}(\frak{S};\mathbb{Z})$, satisfying the conditions (E1) and (E2) as in the statement of the present proposition. \redfix{Recall from Def.\ref{def:translation_to_PGL3} that $\mathbb{I}^+_{{\rm PGL}_3}(\ell) = \Psi^* (\mathbb{I}_{{\rm SL}_3}(\ell)(\mathbb{R}))$, and from eq.\eqref{eq:same_for_non-peripheral_loop} that $\mathbb{I}_{{\rm SL}_3}(\ell) = F^*\mathbb{I}^0_{{\rm SL}_3}(\ell)$. In view of Prop.\ref{prop:skein_algebra_isomorphism} and Cor.\ref{cor:A2-bangles_basis_for_L_SL3}, note that $\mathbb{I}^0_{{\rm SL}_3}(\ell) = \Phi(W)$. Hence we have}
\begin{align}
\nonumber
\redfix{\mathbb{I}^+_{{\rm PGL}_3}(\ell) = \Psi^*(F^*(\Phi(W))(\mathbb{R})),}
\end{align}
where $W$ is viewed as an element of the ${\rm SL}_3$-skein algebra $\mathcal{S}(\frak{S};\mathbb{Z})$, which is naturally isomorphic to the reduced stated ${\rm SL}_3$-skein algebra $\mathcal{S}_{\rm s}(\frak{S};\mathbb{Z})_{\rm red}$, as $\partial \frak{S} = {\O}$. Thus we should prove the equality
$$
\Psi^*(F^*(\Phi(W))(\mathbb{R}))
= \iota_\Delta {\rm Tr}_{\Delta;\frak{S}}([W,{\O}])
$$
for all ${\rm SL}_3$-webs $W$ satisfying conditions (E1) and (E2). Note that all maps $\Psi^*,F^*,\redfix{\Phi}$, evaluation at $\mathbb{R}$, $\iota_\Delta$, and ${\rm Tr}_{\Delta;\frak{S}}$ are ring homomorphisms, and that the map $\Phi$ defined on $\mathcal{S}(\frak{S};\mathbb{Q})$ and the map ${\rm Tr}_{\Delta;\frak{S}}$ defined on $\mathcal{S}_{\rm s}(\frak{S};\mathbb{Z})_{\rm red}$ respect the defining ${\rm SL}_3$-skein relations. It is known that $\mathcal{S}(\frak{S};\mathbb{Q})$ is generated by oriented \redfix{non-contractible} loops, so it suffices \redfix{to} show the above equality when $W$ is an oriented \redfix{non-contractible} loop $\gamma$. In this case, the left-hand-side $\Psi^*(F^*(\Phi(W))(\mathbb{R}))$ equals the trace-of-monodromy $f^+_\gamma$ (Def.\ref{def:f_plus_gamma}), by construction. And we showed $f^+_\gamma = \iota_\Delta {\rm Tr}_{\Delta;\frak{S}}([W,{\O}])$ above. \qed

\vs

Before proceeding, we state one immediate but non-trivial consequence of the proof of Prop.\ref{prop:SL3_classical_trace_and_I_PGL3}:
\begin{corollary}[the ${\rm SL}_3$ classical trace is independent on triangulations]
\label{cor:state-sum_is_flip_compatible}
Let $W$ be an ${\rm SL}_3$-web in a triangulable punctured surface $\frak{S}$, without external vertices. Let $\Delta,\Delta'$ be ideal triangulations of $\frak{S}$. Then
$$
\iota_\Delta {\rm Tr}_\Delta([W,{\O}]) = \iota_{\Delta'} {\rm Tr}_{\Delta'}([W,{\O}]). \qed
$$
\end{corollary}
\begin{remark}
In order to relate ${\rm Tr}_\Delta(W;{\O})\in \mathcal{Z}_\Delta$ and ${\rm Tr}_{\Delta'}(W;{\O}) \in \mathcal{Z}_{\Delta'}$ directly, one first needs to come up with a coordinate change isomorphism between (\redfix{the} fields of fractions of) some `balanced' subalgebras of the algebras $\mathcal{Z}_\Delta$ and $\mathcal{Z}_{\Delta'}$; compare with the ${\rm SL}_2$ case studied in \cite{BW} \cite{Hiatt}. \redfix{This is done in a follow-up paper \cite{Kim21} to the present one, in the form of the following statement, which used to be Conjecture 5.72 in a previous version of the present paper (ver3).}
\end{remark}
\begin{proposition}[the compatibility of the ${\rm SL}_3$ quantum trace under changes of triangulations; \cite{Kim21}]
\label{prop:quantum_coordinate_change}
\redfix{A quantum version of Cor.\ref{cor:state-sum_is_flip_compatible} holds in the following sense. For a triangulable generalized marked surface $\frak{S}$, for two ideal triangulations $\Delta$ and $\Delta'$, let
$$
\Phi^q_{\Delta\Delta'}~:~{\rm Frac}(\mathcal{X}^q_{\Delta'}) \to {\rm Frac}(\mathcal{X}^q_\Delta)
$$
be the quantum coordinate change map, obtained as composition of the quantum mutations $\mu^q_k$ (e.g. from \cite{FG09}) corresponding to the sequence of classical mutations $\mu_k$ relating the cluster $\mathscr{X}$-seeds for $\Delta$ and $\Delta'$. For each $\Delta$, define the \ul{\em balanced subalgebra} $\wh{\mathcal{Z}}^\omega_\Delta$ of $\mathcal{Z}^\omega_\Delta$ as the subalgebra spanned by the monomials $[\prod_{v\in\mathcal{V}(Q_\Delta)}\wh{X}_v^{a_v}]_{\rm Weyl} = [\prod_{v\in\mathcal{V}(Q_\Delta)}\wh{Z}_v^{3a_v}]_{\rm Weyl}$ with $(a_v)_v \in (\frac{1}{3}\mathbb{Z})^{\mathcal{V}(Q_\Delta)}$ being balanced in the sense of Prop.\ref{prop:tropical_coordinate_is_well-defined}. Then there exist balanced quantum coordinate change maps
$$
\Theta^\omega_{\Delta\Delta'} ~:~ {\rm Frac}(\wh{\mathcal{Z}}^\omega_{\Delta'}) \to {\rm Frac}(\wh{\mathcal{Z}}^\omega_\Delta)
$$
extending $\Phi^q_{\Delta\Delta'}$, recovering the classical formula as $\omega^{1/2}\to 1$, and satisfying the consistency relations $\Theta^\omega_{\Delta\Delta''} = \Theta^\omega_{\Delta\Delta'}\Theta^\omega_{\Delta'\Delta''}$. The ${\rm SL}_3$ quantum traces are compatible under these balanced quantum coordinate change maps:}
$$
\redfix{{\rm Tr}^\omega_\Delta = \Theta^\omega_{\Delta\Delta'} \circ {\rm Tr}^\omega_{\Delta'}.}
$$
\end{proposition}
\redfix{We refer the readers to \cite{Kim21} for more details for the above proposition.}

\vs

After a long journey, we finally prove the following.

\vs

{\it Proof of Prop.\ref{prop:highest_term} and Prop.\ref{prop:congruence-compatibility_of_terms_of_each_basic_regular_function}.} Let $\frak{S}$ be a triangulable punctured surface, $\Delta$ be an ideal triangulation of $\frak{S}$, and let $\ell \in \mathscr{A}_{\rm L}(\frak{S};\mathbb{Z})$. One can write $\ell = \ell_1 \cup \ell_2$ as disjoint union, where $\ell_1$ consists only of peripheral loops, and $\ell_2$ has no peripheral loop. Recall eq.\eqref{eq:I_plus_peripheral}, which says 
\begin{align}
\nonumber
\mathbb{I}^+_{{\rm PGL}_3}(\ell_1) = \underset{v \in \mathcal{V}(Q_\Delta)}{\textstyle \prod} X_v^{{\rm a}_v(\ell_1)}.
\end{align}
Meanwhile, $\ell_2$ can be represented as an ${\rm SL}_3$-web $W_2$ satisfying (E1) and (E2) of Prop.\ref{prop:SL3_classical_trace_and_I_PGL3}, hence eq.\eqref{eq:SL3_classical_trace_and_I_PGL3} holds for $\ell_2$:
\begin{align}
\nonumber
\mathbb{I}^+_{{\rm PGL}_3}(\ell_2) = \iota_\Delta {\rm Tr}_{\Delta;\frak{S}}([W_2;{\O}])
\end{align}
Since ${\rm Tr}_{\Delta;\frak{S}}([W_2;{\O}]) \in \mathcal{Z}^1_\Delta$ (Prop.\ref{prop:balancedness_of_state-sum_trace}), and in view of Def.\ref{def:iota_Delta}, it follows that $\mathbb{I}^+_{{\rm PGL}_3}(\ell_2)$ can be written as a Laurent polynomial in $\{X_v^{1/3} \,|\, v\in \mathcal{V}(Q_\Delta)\}$ with integer coefficients. By Prop.\ref{prop:highest_term_of_SL3_classical_trace}, such a Laurent polynomial expression can be chosen so that there is a unique highest order Laurent monomial, which is $\prod_{v\in \mathcal{V}(Q_\Delta)} X_v^{{\rm a}_v(\ell_2)}$ and is of coefficient $1$. And by Prop.\ref{prop:congruence_of_terms_of_SL3_classical_trace}, such a Laurent polynomial expression can be chosen so that other Laurent monomials appearing in this expression are $\prod_{v\in \mathcal{V}(Q_\Delta)} X_v^{{\rm a}_v(\ell_2)}$ times some integer powers of $X_v$'s. By (partial) multiplicativity of $\mathbb{I}^+_{{\rm PGL}_3}$ (see Lem.\ref{lem:product_factors}), 
we have $\mathbb{I}^+_{{\rm PGL}_3}(\ell) = \mathbb{I}^+_{{\rm PGL}_3}(\ell_1)\mathbb{I}^+_{{\rm PGL}_3}(\ell_2)$, hence it follows that $\mathbb{I}^+_{{\rm PGL}_3}(\ell)$ can be written as a Laurent polynomial $\{X_v^{1/3} \,|\, v\in \mathcal{V}(Q_\Delta)\}$ with integer coefficients, so that $\prod_v X_v^{{\rm a}_v(\ell_1)} \prod_v X_v^{{\rm a}_v(\ell_2)} = \prod_v X_v^{{\rm a}_v(\ell)}$ ($\because$ Lem.\ref{lem:additivity_of_coordinates}) is the unique highest order term with coefficient $1$, while the other terms are $\prod_v X_v^{{\rm a}_v(\ell)}$ times some integer powers of $X_v$'s. \qed

\subsection{\redfix{On the effect of a single mutation}}
\label{subsec:on_the_effect_of_a_single_mutation}

\redfix{Now, in the proof of the first main theorem, Thm.\ref{thm:main}, given in the previous section, what remain to be proved are Prop.\ref{prop:mutation_of_basic_semi-regular_function_at_interior_node_of_triangle} and Prop.\ref{prop:mutation_of_basic_semi-regular_function_at_edge_node_of_triangle}. These are on the effect of the mutation at a single node on a basic semi-regular function $\mathbb{I}^+_{{\rm PGL}_3}(\ell)$, $\ell \in \mathscr{A}_{\rm L}(\frak{S};\mathbb{Z})$, when $\frak{S}$ is a triangulable punctured surface, and they together yield Cor.\ref{cor:mutation_of_congruent_ell} which in turn is a crucial step in the proof of Thm.\ref{thm:main}. In \S\ref{subsec:mutation_of_basic_regular_functions} we proved these two propositions in cases when $\ell$ can be represented by an ${\rm SL}_3$-web without any $3$-valent vertices, i.e. when $\ell$ can be represented by loops. In view of the sought-for statements, by using the additivity of the tropical coordinates of ${\rm SL}_3$-laminations (Lem.\ref{lem:additivity_of_coordinates}), and also the results proved in \S\ref{subsec:mutation_of_basic_regular_functions} for the case when $\ell$ consists only of peripheral loops, now we can just deal with the case when $\ell$ does not contain any peripheral loop. Then $\ell$ can be represented by an ${\rm SL}_3$-web with non-negative weights, hence can be represented by a (reduced) non-elliptic ${\rm SL}_3$-web $W$ in $\frak{S}$ with all weights being $1$. Thus, by Prop.\ref{prop:SL3_classical_trace_and_I_PGL3}, eq.\eqref{eq:SL3_classical_trace_and_I_PGL3} holds: $\mathbb{I}^+_{{\rm PGL}_3}(\ell) = \iota_\Delta {\rm Tr}_{\Delta;\frak{S}}([W;{\O}])$, when $\Delta$ is an ideal triangulation of $\frak{S}$. So, we will try to mutate the expression ${\rm Tr}_{\Delta;\frak{S}}([W;{\O}]) \in \mathcal{Z}^\omega_\Delta$ (or its image under $\iota_\Delta$, to be more precise) at a single node of the quiver $Q_\Delta$. As done in \S\ref{subsec:mutation_of_basic_regular_functions}, we will refer to Fig.\ref{fig:mutate_two_triangles} for the labels of the sides and nodes that are relevant.}

\vs

\begin{figure}[htbp!]
\vspace{-2mm}
\begin{center}
\hspace{0mm} \raisebox{-0.5\height}{\scalebox{0.8}{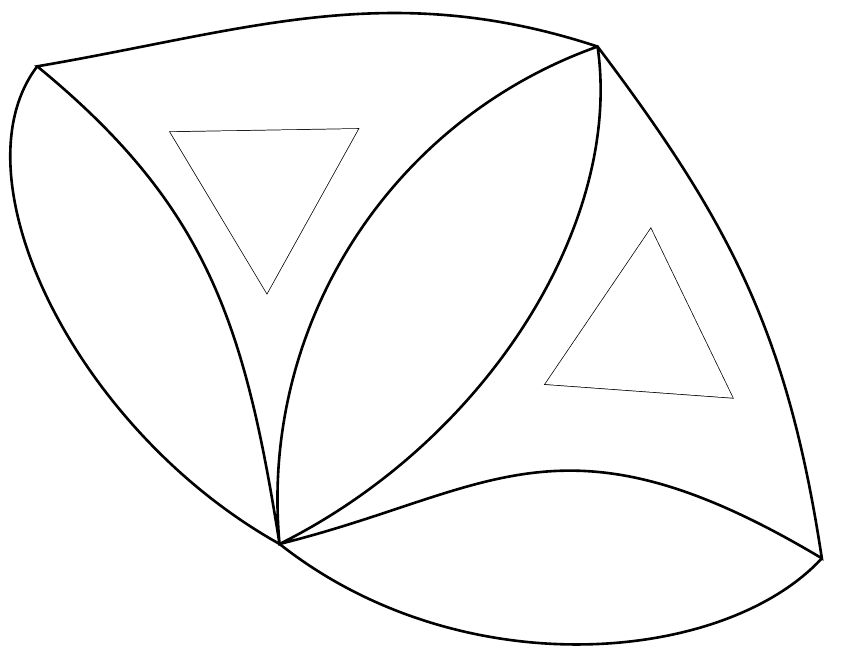}}
$\to$
\raisebox{-0.5\height}{\scalebox{0.8}{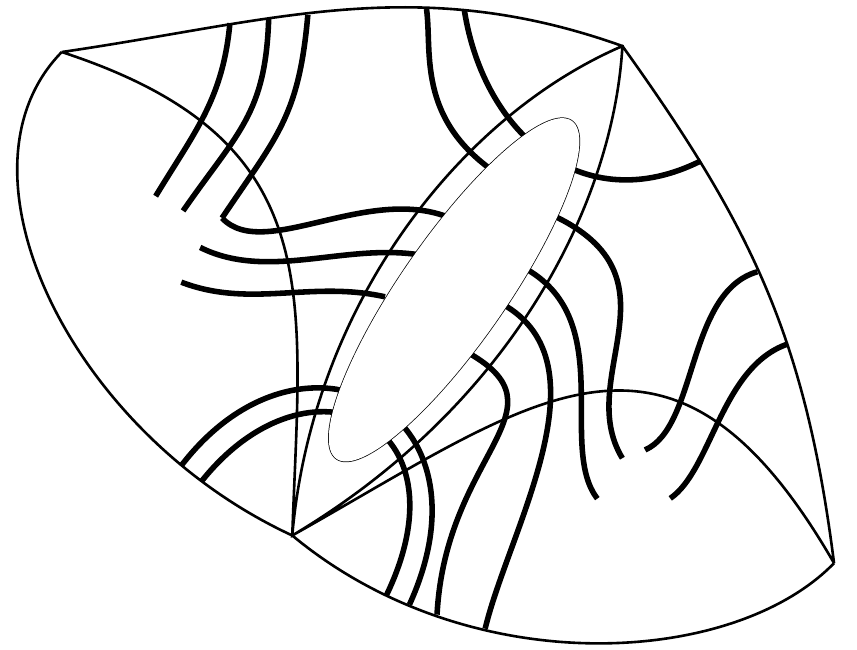}} 
\end{center}
\vspace{-4mm}
\caption{Pushing pyramids from quadrilateral to neighboring biangles}
\vspace{-2mm}
\label{fig:pyramid_in_quadrilateral_pushing}
\end{figure}

\redfix{Let's begin with Prop.\ref{prop:mutation_of_basic_semi-regular_function_at_interior_node_of_triangle}, where we are mutating at the interior node $v_t$ in an ideal triangle $t$ of $\Delta$, which is the triangle on the left of the quadrilateral in Fig.\ref{fig:mutate_two_triangles}. In particular, the side names are $e_1,e_2,e_3$ in the clockwise order, and the nodes on the sides are $v_{e_i,1}$, $v_{e_i,2}$ on each $e_i$. We make use of the state-sum formula as in eq.\eqref{eq:state-sum_formula}, with $\omega^{1/2}=1$ and the elevations in ${\bf I}$ ignored:
\begin{align}
\nonumber
{\rm Tr}_\Delta([W,s]) = {\textstyle \sum}_J ( {\textstyle \prod}_B {\rm Tr}_B([W\cap B, J_B]) \, {\textstyle \bigotimes}_t {\rm Tr}_{\wh{t}}([W\cap \wh{t}, J_t]) ) ~\in~ \mathcal{Z}_\Delta ~\in~ 
{\textstyle \bigotimes}_{t\in \mathcal{F}(\Delta)} \mathcal{Z}_t,
\end{align}
where the target ring $\mathcal{Z}_{\wh{t}}$ of the ${\rm SL}_3$ classical trace ${\rm Tr}_{\wh{t}}$ for the triangle $\wh{t}$ is naturally identified with $\mathcal{Z}_t$. Recall that the sum is over all $\wh{\Delta}$-juncture-states $J : W\cap \wh{\Delta} \to \{1,2,3\}$ of $W$ that restrict to $s : \partial W \to \{1,2,3\}$, $\prod_B$ is over all biangles $B$ of $\wh{\Delta}$, $\prod_t$ is over all triangles $t$ of $\Delta$ (or triangles $\wh{t}$ of $\wh{\Delta}$), while $J_B = J|_{\partial (W\cap B)}$ and $J_t = J|_{\partial (W\cap \wh{t})}$. As mentioned in the previous subsection, this state-sum formula works for any ${\rm SL}_3$-web $W$ in $\frak{S}$ not just for one that is in a gool position, provided that for each triangle $\wh{t}$ and a biangle $B$ of $\wh{\Delta}$, $W\cap \wh{t}$ and $W\cap B$ are well-defined ${\rm SL}_3$-webs in $\wh{t}$ and $B$. Therefore, once we put $W$ by isotopy into a canonical position with respect to $\wh{\Delta}$ (as in Def.\ref{def:canonical_wrt_split_ideal_triangulation}), the formula in eq.\eqref{eq:state-sum_formula} works. The left picture in Fig.\ref{fig:pyramid_in_quadrilateral_pushing} shows part of $W$ in a canonical position with respect to $\wh{\Delta}$, where the triangle $\wh{t}$ of $\wh{\Delta}$ corresponds to the triangle $t$ of $\Delta$. For our purpose at the moment, one can just focus on $\wh{t}$ and the two neighboring biangles (one can also draw one more biangle on top of $\wh{t}$, but we omitted it), and ignore the triangle $\wh{r}$ and the biangle below $\wh{r}$. In Fig.\ref{fig:pyramid_in_quadrilateral_pushing} we omitted indicating the orientations of edges and loops of $W$. Note that $W \cap \wh{t}$ is a canonical ${\rm SL}_3$-web in $\wh{t}$ (Def.\ref{def:canonical_web_in_a_triangle}), hence is a disjoint union of some number of corner arcs, and a degree $d_1$ pyramid $H_{d_1}$ for some $d_1\in \mathbb{Z}$.}

\vs

\redfix{When we mutate at the node $v_t$, only the cluster $\mathscr{X}$-variables for the nodes of $Q_\Delta$ lying in $t$ may transform, and other variables stay the same. Note that the node $v_t$ lies in $t$ but not in any other triangle, while each node $v$ lying in a side edge of $t$ also lies in the corresponding side edge of a neighboring triangle, say $r$. In fact, in every (nonzero) term of ${\rm Tr}_\Delta([W,s]) \in {\textstyle \bigotimes}_{t\in \mathcal{F}(\Delta)} \mathcal{Z}_t$, we showed in Prop.\ref{prop:balancedness_of_state-sum_trace} that the power of the variable $Z_{t,v} \in \mathcal{Z}_t$ associated to such a node $v$ for triangle $t$ equals that of the variable $Z_{r,v} \in \mathcal{Z}_r$ for $r$; recall that $Z_{t,v} Z_{r,v}$ is defined as the variable $Z_v \in \mathcal{Z}_\Delta$. So, to study the effect of mutation at $v_t$ on ${\rm Tr}_\Delta([W,s])$, it suffices to study the effect on the factor ${\rm Tr}_{\wh{t}}([W\cap \wh{t},J_t])$ for the triangle $\wh{t}$, and indeed that is what is done in the partial proof of Prop.\ref{prop:mutation_of_basic_semi-regular_function_at_interior_node_of_triangle} given in \S\ref{subsec:mutation_of_basic_regular_functions}, in case when $W\cap \wh{t}$ consists only of corner arcs. So now we study the case when $W\cap \wh{t}$ is any general canonical ${\rm SL}_3$-web in the triangle $\wh{t}$, e.g. as presented in Fig.\ref{fig:pyramid_in_quadrilateral_pushing}. In particular, it contains a pyramid web $H_{d_1}$, with $d_1 \in \mathbb{Z}$. We divide $\wh{t}$ into a triangle $\til{t}$ and a biangle $B_1$ as in the right picture of Fig.\ref{fig:pyramid_in_quadrilateral_pushing}, by adding one more ideal arc parallel to one side of $\wh{t}$ (it corresponds to the side $e_2$ in Fig.\ref{fig:mutate_two_triangles}), and push the pyramid web $H_{d_1}$ by isotopy into the biangle $B_1$ so that it becomes the web $P_{d_1}$ which was dealt with in Lem.\ref{lem:P_d}. By the results of \S\ref{subsec:isotopy_invariance} one has the state-sum formula}
$$
\redfix{{\rm Tr}_{\wh{t}}([W\cap \wh{t}, J_t]) = \sum_{J_{\til{t}}, J_{B_1}} {\rm Tr}_{B_1}([W\cap B_1, J_B]) \, {\rm Tr}_{\til{t}}([W\cap \til{t},J_{\til{t}}]),}
$$
\redfix{summed over all states $J_B$ and $J_{\til{t}}$ compatible with $J_t$ in an appropriate sense, where the target ring $\mathcal{Z}_{\til{t}}$ of ${\rm Tr}_{\til{t}}$ is naturally identified with $\mathcal{Z}_{\wh{t}} = \mathcal{Z}_t$. Note now that $W\cap \til{t}$ consists only of corner arcs, hence the results of the partial proof of Prop.\ref{prop:mutation_of_basic_semi-regular_function_at_interior_node_of_triangle} obtained in \S\ref{subsec:mutation_of_basic_regular_functions} apply. Therefore, by applying the mutation at $v_t$ to ${\rm Tr}_{\til{t}}([W\cap \til{t},J_{\til{t}}]) \in \mathcal{Z}_t$ we get}
$$
\redfix{{\rm Tr}_{\til{t}}([W\cap \til{t},J_{\til{t}}]) ~\in~ {X'}_{\hspace{-1,5mm}v_t}^{- \til{\rm a}_{v_t}(W\cap \til{t}\,)} \, ({\textstyle \prod}_{v\in (\mathcal{V}(Q') \cap t)\setminus\{v_t\}} {X'}_{\hspace{-1,5mm}v}^{\til{\rm a}_v(W\cap \til{t}\,)}) \cdot \mathbb{Z}[\{{X'}_{\hspace{-1,5mm}v}^{\pm 1} \, | \, v\in (\mathcal{V}(Q') \cap t) \}],}
$$
\redfix{where, $\til{\rm a}_v(W\cap \til{t}\,)$ stands for the tropical coordinate at $v$ of $W\cap \til{t}$ viewed as an ${\rm SL}_3$-lamination (with all weights $1$) in $\til{t}$. Here, the nodes $v$ of the quiver $Q'$ (obtained from $Q_\Delta$ by mutating at $v_t$) living in $t$ are naturally identified with the nodes $v$ of $Q_{\til{\Delta}}$, where $\til{\Delta}$ is the unique ideal triangulation of $\til{t}$ viewed as a generalized marked surface. On the other hand, for each node $v$ in $\mathcal{V}(Q')\cap t = \mathcal{V}(Q_\Delta) \cap t$ we have ${\rm a}_v(\ell) = {\rm a}_v(W) = \wh{\rm a}_v(W\cap \wh{t}\, )$, where $\wh{\rm a}_v(W\cap \wh{t}\, )$ is the tropical coordinate at $v$ of $W\cap \wh{t}$ viewed as an ${\rm SL}_3$-lamination in $\wh{t}$. We claim that, for each node $v$ lying in $t$, we have
$$
\wh{\rm a}_v(W\cap \wh{t}\,) \equiv \til{\rm a}_v(W\cap \til{t}\, )  \quad {\rm modulo} ~ \mathbb{Z},
$$
which would yield the desired result of Prop.\ref{prop:mutation_of_basic_semi-regular_function_at_interior_node_of_triangle} for ${\rm Tr}_\Delta([W,{\O}])$, i.e. for $\mathbb{I}^+_{{\rm PGL}_3}(\ell)$.} \redfix{By the additivity of tropical coordinates (Lem.\ref{lem:additivity_of_coordinates}), and since the corner arcs of $W\cap \wh{t}$ stay the same in $W\cap \til{t}$ (after isotopy), it suffices to check the above equality only for the case when $W\cap \wh{t}$ is the pyramid web $H_{d_1}$ for some $d_1\in \mathbb{Z}$. Then, as can be read from Fig.\ref{fig:DS_coordinates}, if $d_1>0$, we have $\wh{\rm a}_{v_t}(W\cap \wh{t}\,) = d_1$, $\wh{\rm a}_{v_{e_\alpha,1}}(W\cap \wh{t}\,) = d_1/3$, $\wh{\rm a}_{v_{e_\alpha,2}}(W\cap \wh{t}\,) = 2d_1/3$, for $\alpha=1,2,3$, while $\til{a}_{v_t}(W\cap \til{t}\,) = d_1$, $\til{\rm a}_{v_{e_\alpha,1}}(W\cap \til{t}\,) = d_1/3$, $\til{\rm a}_{v_{e_\alpha,2}}(W\cap \til{t}\,) = 2d_1/3$, for $\alpha=1,3$, and $\til{\rm a}_{v_{e_2,1}}(W\cap \til{t}\,) = 4d_1/3$, $\til{\rm a}_{v_{e_2,2}}(W\cap \til{t}\,) = 2d_1/3$. The only discrepancy is at the node $v_{e_2,1}$, with the difference being $d_1$ which is an integer. Similarly, when $d_1<0$, the only discrepancy is at the node $v_{e_2,2}$, with $\wh{\rm a}_{v_{e_2,2}}(W\cap \wh{t}\,) = d_1/3$ and $\til{\rm a}_{v_{e_2,2}}(W\cap \til{t}\,) = 4d_1/3$, the difference being $d_1 \in \mathbb{Z}$, as desired. This concludes the proof of Prop.\ref{prop:mutation_of_basic_semi-regular_function_at_interior_node_of_triangle}.}

\vs

\redfix{We proceed to the proof of Prop.\ref{prop:mutation_of_basic_semi-regular_function_at_edge_node_of_triangle}, where we are mutating at a node $v_0$ lying on an ideal arc of $\Delta$; under the notation of Fig.\ref{fig:mutate_two_triangles}, we mutate at the node $v_{e_1,1}$. We assume that $W$ is in a canonical position with respect to $\wh{\Delta}$, and use the state-sum formula for ${\rm Tr}_\Delta([W,{\O}])$ as used in the above proof of Prop.\ref{prop:mutation_of_basic_semi-regular_function_at_interior_node_of_triangle}. This time, the mutation may transform the variables for nodes living in the triangles $t$ and $r$, but not others. So, in the state-sum formula, we only need to care about the factors for the triangles $\wh{t}$ and $\wh{r}$, and the biangle $B_0$ in between them (see Fig.\ref{fig:pyramid_in_quadrilateral_pushing}), i.e. focus on the factor
$$
{\rm Tr}_{B_0}([W\cap B_0, J_{B_0}]) \, {\rm Tr}_{\wh{t}}([W\cap \wh{t}, J_t]) \, {\rm Tr}_{\wh{r}}([W\cap \wh{r}, J_r]) ~\in ~ \mathcal{Z}_t \otimes \mathcal{Z}_r,
$$
that contributes a nonzero term in the original state-sum formula for ${\rm Tr}_\Delta([W,{\O}])$; one could also package these factors appearing in ${\rm Tr}_\Delta([W,{\O}])$ into a sum
$$
\sum_{J_{B_0 \cup \wh{t} \cup \wh{r}}} {\rm Tr}_{B_0\cup \wh{t} \cup \wh{r}}( [W\cap (B_0 \cup \wh{t} \cup \wh{r}\,), J_{B_0 \cup \wh{t} \cup \wh{r}}]),
$$
in an appropriate sense.} \redfix{Indeed, the effect of the mutation at $v_{e_1,1}$ on this expression is what is investigated in the partial proof of 
Prop.\ref{prop:mutation_of_basic_semi-regular_function_at_edge_node_of_triangle} in \S\ref{subsec:mutation_of_basic_regular_functions}, in case when $W\cap (B_0 \cup \wh{t} \cup \wh{r}\,)$ has no $3$-valent vertex.}

\vs

\redfix{Using similar ideas we employed in the above proof of Prop.\ref{prop:mutation_of_basic_semi-regular_function_at_interior_node_of_triangle}, we will push by isotopy the $3$-valent vertices of $W$ to biangles adjacent to $\wh{t}$ and $\wh{r}$ other than $B_0$, so that $W \cap (B_0 \cup \wh{t} \cup \wh{r}\,)$ has no $3$-valent vertex, and in the meantime, verify that the tropical coordinates for the original $W$ (or $\ell$) at the nodes living in $t$ and $r$ differ from the new coordinates for the isotoped $W$ at those nodes by integers. Similarly as before, we divide the triangle $\wh{t}$ into the triangle $\til{t}$ and a biangle $B_1$ (located at the side edge $e_2$), and the triangle $\wh{r}$ into the triangle $\til{r}$ and a biangle $B_2$ (at the side edge $e_6$), and push the pyramid web $H_{d_1}$ of $W\cap \wh{t}$ by isotopy to $B_1$ and the pyramid web $H_{d_2}$ of $W\cap \wh{r}$ by isotopy to $B_2$, as shown in Fig.\ref{fig:pyramid_in_quadrilateral_pushing}; one can view the situation so that $W$ stays the same (i.e. no isotoping is being done) and we chose the additional arcs forming the new biangles $B_1$ and $B_2$ as such. Then, if one considers the ${\rm SL}_3$-web $W\cap \til{t}$ in $\til{t}$ and the ${\rm SL}_3$-web $W\cap \til{r}$ in $\til{r}$, viewed as ${\rm SL}_3$-laminations in $\til{t}$ and $\til{r}$, the tropical coordinates of them differ from the those of the original ${\rm SL}_3$-web $W$ (i.e. of $W\cap \wh{t}$ and of $W\cap \wh{r}$) at the nodes living in $t$ and $r$ by integers, as seen in the above proof of Prop.\ref{prop:mutation_of_basic_semi-regular_function_at_interior_node_of_triangle}. In the meantime, by the results of \S\ref{subsec:isotopy_invariance} one has the state-sum formula}
\begin{align*}
& \sum_{J_{B_0 \cup \wh{t} \cup \wh{r}}} {\rm Tr}_{B_0\cup \wh{t} \cup \wh{r}}( [W\cap (B_0 \cup \wh{t} \cup \wh{r}\,), J_{B_0 \cup \wh{t} \cup \wh{r}}]) \\
& = \sum_{J_{B_1},J_{B_2}, J_{B_0\cup \til{t} \cup \til{r}}} {\rm Tr}_{B_1}([W\cap B_1, J_{B_1}]) \, {\rm Tr}_{B_2}([W\cap B_2, J_{B_2}]) \, {\rm Tr}_{B_0\cup \til{t} \cup \til{r}}( [W\cap (B_0 \cup \til{t} \cup \til{r}\,), J_{B_0 \cup \til{t} \cup \til{r}}]);
\end{align*}
\redfix{hence we could just focus on the factor ${\rm Tr}_{B_0\cup \til{t} \cup \til{r}}( [W\cap (B_0 \cup \til{t} \cup \til{r}\,), J_{B_0 \cup \til{t} \cup \til{r}}])$.}

\vs

\begin{figure}[htbp!]
\vspace{-2mm}
\begin{center}
\hspace{0mm} \raisebox{-0.5\height}{\scalebox{0.9}{
\begingroup%
  \makeatletter%
  \providecommand\color[2][]{%
    \errmessage{(Inkscape) Color is used for the text in Inkscape, but the package 'color.sty' is not loaded}%
    \renewcommand\color[2][]{}%
  }%
  \providecommand\transparent[1]{%
    \errmessage{(Inkscape) Transparency is used (non-zero) for the text in Inkscape, but the package 'transparent.sty' is not loaded}%
    \renewcommand\transparent[1]{}%
  }%
  \providecommand\rotatebox[2]{#2}%
  \newcommand*\fsize{\dimexpr\f@size pt\relax}%
  \newcommand*\lineheight[1]{\fontsize{\fsize}{#1\fsize}\selectfont}%
  \ifx\svgwidth\undefined%
    \setlength{\unitlength}{175.7480315bp}%
    \ifx\svgscale\undefined%
      \relax%
    \else%
      \setlength{\unitlength}{\unitlength * \real{\svgscale}}%
    \fi%
  \else%
    \setlength{\unitlength}{\svgwidth}%
  \fi%
  \global\let\svgwidth\undefined%
  \global\let\svgscale\undefined%
  \makeatother%
  \begin{picture}(1,0.88709677)%
    \lineheight{1}%
    \setlength\tabcolsep{0pt}%
    \put(0,0){\includegraphics[width=\unitlength,page=1]{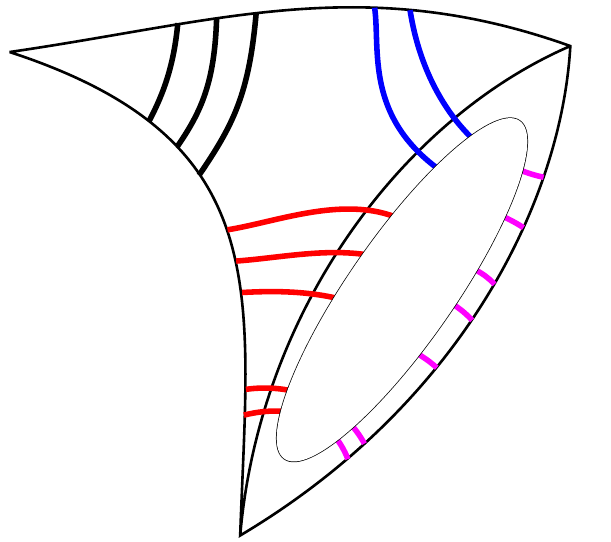}}%
    \put(0.64149781,0.39374985){\makebox(0,0)[lt]{\lineheight{1.25}\smash{\begin{tabular}[t]{l}$L_0$\end{tabular}}}}%
    \put(0,0){\includegraphics[width=\unitlength,page=2]{crossbar_web_pushing1.pdf}}%
    \put(0.50198643,0.67339946){\makebox(0,0)[lt]{\lineheight{1.25}\smash{\begin{tabular}[t]{l}$\til{t}$\end{tabular}}}}%
    \put(0.83928756,0.71438695){\makebox(0,0)[lt]{\lineheight{1.25}\smash{\begin{tabular}[t]{l}$B_0$\end{tabular}}}}%
  \end{picture}%
\endgroup%
}} \hspace{5mm}
$\to$  \hspace{5mm} 
\raisebox{-0.5\height}{\scalebox{0.9}{
\begingroup%
  \makeatletter%
  \providecommand\color[2][]{%
    \errmessage{(Inkscape) Color is used for the text in Inkscape, but the package 'color.sty' is not loaded}%
    \renewcommand\color[2][]{}%
  }%
  \providecommand\transparent[1]{%
    \errmessage{(Inkscape) Transparency is used (non-zero) for the text in Inkscape, but the package 'transparent.sty' is not loaded}%
    \renewcommand\transparent[1]{}%
  }%
  \providecommand\rotatebox[2]{#2}%
  \newcommand*\fsize{\dimexpr\f@size pt\relax}%
  \newcommand*\lineheight[1]{\fontsize{\fsize}{#1\fsize}\selectfont}%
  \ifx\svgwidth\undefined%
    \setlength{\unitlength}{175.7480315bp}%
    \ifx\svgscale\undefined%
      \relax%
    \else%
      \setlength{\unitlength}{\unitlength * \real{\svgscale}}%
    \fi%
  \else%
    \setlength{\unitlength}{\svgwidth}%
  \fi%
  \global\let\svgwidth\undefined%
  \global\let\svgscale\undefined%
  \makeatother%
  \begin{picture}(1,0.88709677)%
    \lineheight{1}%
    \setlength\tabcolsep{0pt}%
    \put(0,0){\includegraphics[width=\unitlength,page=1]{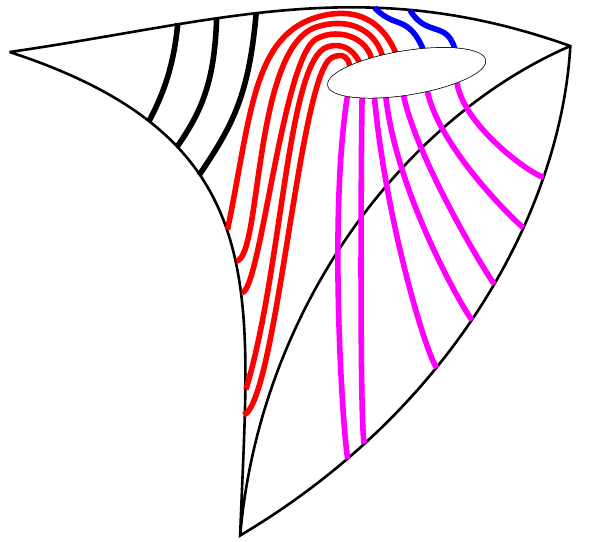}}%
    \put(0.64385953,0.75216797){\rotatebox{9.23432342}{\makebox(0,0)[lt]{\lineheight{1.25}\smash{\begin{tabular}[t]{l}$L_0$\end{tabular}}}}}%
    \put(0,0){\includegraphics[width=\unitlength,page=2]{crossbar_web_pushing2.pdf}}%
    \put(0.50679888,0.54916607){\makebox(0,0)[lt]{\lineheight{1.25}\smash{\begin{tabular}[t]{l}$\check{t}$\end{tabular}}}}%
    \put(0.82576963,0.69279075){\makebox(0,0)[lt]{\lineheight{1.25}\smash{\begin{tabular}[t]{l}$B_0$\end{tabular}}}}%
    \put(0,0){\includegraphics[width=\unitlength,page=3]{crossbar_web_pushing2.pdf}}%
    \put(0.79928442,0.79628298){\makebox(0,0)[lt]{\lineheight{1.25}\smash{\begin{tabular}[t]{l}$B_3$\end{tabular}}}}%
  \end{picture}%
\endgroup%
}} 
\end{center}
\vspace{-4mm}
\caption{Pushing a crossbar web from a biangle to another biangle}
\vspace{-2mm}
\label{fig:crossbar_web_pushing}
\end{figure}

\redfix{We need one more isotoping, because $W\cap (B_0 \cup \til{t} \cup \til{r}\,)$ can still have $3$-valent vertices in the biangle $B_0$. So, we fatten one of the sides of $\til{t}$ and $\til{r}$, other than the ones corresponding to $e_1$, to a biangle, and push the $3$-valent vertices of $W\cup B_0$ into that biangle. In Fig.\ref{fig:crossbar_web_pushing}, we choose to divide $\til{t}$ into $\check{t}$ and a biangle $B_3$ (at the side $e_3$) and push the `entire' ${\rm SL}_3$-web $W\cap B_0$ into $B_3$ by isotopy. By a similar state-sum formula resulting from \S\ref{subsec:isotopy_invariance}}
$$
{\rm Tr}_{B_0\cup \til{t} \cup \til{r}}( [W\cap (B_0 \cup \til{t} \cup \til{r}\,), J_{B_0 \cup \til{t} \cup \til{r}}])
= \sum_{J_{B_3}, J_{B_0\cup \check{t} \cup \til{r}}} {\rm Tr}_{B_3}([W\cap B_3, J_{B_3}]) \, {\rm Tr}_{B_0\cup \check{t} \cup \til{r}}( [W\cap (B_0 \cup \check{t} \cup \til{r}\,), J_{B_0 \cup \check{t} \cup \til{r}}])
$$
\redfix{we could focus on ${\rm Tr}_{B_0\cup \check{t} \cup \til{r}}( [W\cap (B_0 \cup \check{t} \cup \til{r}\,), J_{B_0 \cup \check{t} \cup \til{r}}])$. Now the ${\rm SL}_3$-web $W\cap (B_0 \cup \check{t} \cup \til{r}\,)$ living in $B_0 \cup \check{t} \cup \til{r}$ does not have any $3$-valent vertex, so the results of the partial proof of Prop.\ref{prop:mutation_of_basic_semi-regular_function_at_edge_node_of_triangle} obtained in \S\ref{subsec:mutation_of_basic_regular_functions} applies. Therefore, by applying the mutation at $v_0 := v_{e_1,1}$ to ${\rm Tr}_{B_0\cup \check{t} \cup \til{r}}( [W\cap (B_0 \cup \check{t} \cup \til{r}\,), J_{B_0 \cup \check{t} \cup \til{r}}])$ we get}
\begin{align*}
& {\rm Tr}_{B_0\cup \check{t} \cup \til{r}}( [W\cap (B_0 \cup \check{t} \cup \til{r}\,), J_{B_0 \cup \check{t} \cup \til{r}}]) \\
& \in~ {X''}_{\hspace{-2mm}v_0}^{- \check{\rm a}_{v_0}(\check{W}) + \sum_{v\in \mathcal{V}(Q_\Delta)} [\varepsilon_{v_0 v}]_+ \check{\rm a}_v(\check{W})} \, (\underset{v\in (\mathcal{V}(Q'') \cap (t\cup r))\setminus\{v_0\}}{\textstyle \prod} {X''}_{\hspace{-2mm}v}^{\check{\rm a}_v(\check{W})}) \cdot \mathbb{Z}[\{{X''}_{\hspace{-2mm}v}^{\pm 1} \, | \, v\in \mathcal{V}(Q'')\cap (t\cup r)\}],
\end{align*}
\redfix{with appropriate identifications in the style of the above proof of Prop.\ref{prop:mutation_of_basic_semi-regular_function_at_interior_node_of_triangle}, where $\check{W} := W\cap (B_0 \cup \check{t} \cup \til{r}\,)$ is viewed as an ${\rm SL}_3$-lamination in the surface $B_0 \cup \check{t} \cup \til{r}$, for which the tropical coordinates are denoted by $\check{\rm a}_v(\check{W})$; here $X''_*$ stands for the cluster $\mathscr{X}$-variables resulting after applying the mutation at $v_0$ to the cluster $\mathscr{X}$-chart associated to the triangulation $\Delta$, and $Q''$ the quiver obtained after the mutation at $v_0$ applied to $Q_\Delta$. Now, all that remains is to show 
$$
{\rm a}_v(\ell) \equiv \check{\rm a}_v(\check{W}) \quad {\rm modulo} ~ \mathbb{Z}
$$
for all nodes $v\in \mathcal{V}(Q'') \cap (t\cup r)$.}

\vs

\redfix{We have seen that the original coordinates ${\rm a}_v(\ell)$ coincide up to integers with the coordinates of the ${\rm SL}_3$-lamination $W\cap \til{t}$ (with weights $1$) in $\til{t}$ and the ${\rm SL}_3$-lamination $W\cap \til{r}$ in $\til{r}$, where these $W$ denote the original $W$ before we isotoped the $3$-valent vertices the into $B_3$. We should then just check whether the coordinates for $W\cap \til{t}$ for the nodes $v$ in $t$ coincide up to integers with $\check{a}_v(\check{W})$; i.e. in Fig.\ref{fig:crossbar_web_pushing}, compare the tropical coordinates for the ${\rm SL}_3$-lamination in $\til{t}$ in the left picture and those for the ${\rm SL}_3$-lamination in $\check{t}$ in the right picture. This can be easily checked for each corner arc of $W\cap \til{t}$ in the left picture, explained as follows. The black corner arcs of $W\cap \til{t}$ in the left picture of Fig.\ref{fig:crossbar_web_pushing} stay the same in $\check{W}$ in the right picture, so they contribute the same amount to the tropical coordinates. Since the web $W\cap B_0$ (i.e. $L_0$ in the picture) is a crossbar web in a biangle, by Lem.\ref{lem:minimal_crossbar_A2-web_determined_by_signature}(MC2), the total number of red or blue arcs going into $B_0$ from the left (resp. out of $B_0$ to the left) is same as the number of purple arcs going out of $B_0$ to the right (resp. into $B_0$ from the right). This gives a non-unique bijection from the set of red or blue arcs to the set of purple arcs (in fact, a crossbar web $W\cap B_0$ can be turned into a wiring diagram, so there is a natural preferred bijection, but we don't need it here).}

\vs

\redfix{Consider a red arc of $W\cap \til{t}$ in the left picture of Fig\ref{fig:crossbar_web_pushing} which connects the sides corresponding to $e_2$ and $e_1$ of Fig.\ref{fig:mutate_two_triangles}. If it is oriented from $e_2$ toward $e_1$, then, in view of Fig.\ref{fig:DS_coordinates}, the contributions from this arc to the tropical coordinates for $W\cap \til{t}$ are ${\rm a}_{v_{e_1,1}} = 1/3$, ${\rm a}_{v_{e_1,2}}=2/3$, ${\rm a}_{v_{e_2,1}}=2/3$, ${\rm a}_{v_{e_2,2}}=1/3$, ${\rm a}_{v_{e_3,1}}={\rm a}_{v_{e_3,2}}=0$, ${\rm a}_{v_t} = 1/3$. After isotoping, the sum of contributions of the isotoped red arc in $\check{t}$ (note that this isotoped red arc goes from $e_2$ to $e_3$) and the purple arc in $\check{t}$ corresponding to this red arc (note that this purple arc goes from $e_3$ to $e_1$) to the coordinates of $\check{W}$ are $\til{\rm a}_{v_{e_1,1}} = 0+1/3$, $\til{\rm a}_{v_{e_1,2}} = 0+2/3$, $\til{\rm a}_{v_{e_2,1}} = 2/3+0$, $\til{\rm a}_{v_{e_2,2}} = 1/3+0$, $\til{\rm a}_{v_{e_3,1}} = 1/3+2/3$, $\til{\rm a}_{v_{e_3,2}} = 2/3+1/3$, $\til{\rm a}_{v_t} = 2/3+2/3$. By inspection, the difference between the old and the new coordinates are integers, as desired. For a red arc of $W\cap \til{t}$ going from $e_1$ to $e_2$, the contributions to previous coordinates are ${\rm a}_{v_{e_1,1}} = 2/3$, ${\rm a}_{v_{e_1,2}}=1/3$, ${\rm a}_{v_{e_2,1}}=1/3$, ${\rm a}_{v_{e_2,2}}=2/3$, ${\rm a}_{v_{e_3,1}}={\rm a}_{v_{e_3,2}}=0$, ${\rm a}_{v_t} = 2/3$, while after isotoping, the contributions to new coordinates by the isotoped red arc (from $e_3$ to $e_2$) and the corresponding purple arc (from $e_1$ to $e_3$) are $\til{\rm a}_{v_{e_1,1}} = 0+2/3$, $\til{\rm a}_{v_{e_1,2}} = 0+1/3$, $\til{\rm a}_{v_{e_2,1}} = 1/3+0$, $\til{\rm a}_{v_{e_2,2}} = 2/3+0$, $\til{\rm a}_{v_{e_3,1}} = 2/3+1/3$, $\til{\rm a}_{v_{e_3,2}} = 1/3+2/3$, $\til{\rm a}_{v_t} = 1/3+1/3$. So the difference of the old and the new coordinates are integers. For a blue arc of $W\cap \til{t}$, after isotoping, we consider the purple arc of $\check{W}$ in $\check{t}$ corresponding to this blue arc. Since the blue arc and the corresponding purple arc are same type of corner arcs in triangles, their contribution to the tropical coordinates are exactly the same. This finishes the proof of Prop.\ref{prop:mutation_of_basic_semi-regular_function_at_edge_node_of_triangle}.}

\vs

At last, this justifies our proof of the first main theorem, Thm.\ref{thm:main}, given in the previous section.

\subsection{The quantum ${\rm SL}_3$-${\rm PGL}_3$ duality map}
\label{subsec:quantum_duality_map}

Making use of the ${\rm SL}_3$ quantum trace map ${\rm Tr}^\omega_\Delta$ constructed in the second main theorem, Thm.\ref{thm:SL3_quantum_trace_map}, we propose a quantum version of the duality map $\mathbb{I} : \mathscr{A}_{{\rm SL}_3,\frak{S}}(\mathbb{Z}^{\redfix{T}}) \to \mathscr{O}(\mathscr{X}_{{\rm PGL}_3,\frak{S}})$ of Thm.\ref{thm:main}, as an analog of the ${\rm SL}_2$-${\rm PGL}_2$ quantum duality map constructed in \cite{AK}.
\begin{definition}
\label{def:quantum_duality_map_omega}
Let $\frak{S}$ be a triangulable punctured surface, and $\Delta$ be an ideal triangulation of $\frak{S}$. Let $\mathcal{Z}^\omega_\Delta$ be the cube-root Fock-Goncharov algebra for $\Delta$ (Def.\ref{def:Fock-Goncharov_algebra_quantum}). Define the \ul{\em ${\rm SL}_3$-${\rm PGL}_3$ quantum duality map}
$$
\wh{\mathbb{I}}^\omega_\Delta : \mathscr{A}_{\rm L}(\frak{S};\mathbb{Z}) \to \mathcal{Z}^\omega_\Delta
$$
as follows. Let $\ell \in \mathscr{A}_{\rm L}(\frak{S};\mathbb{Z})$. Represent $\ell$ as disjoint union $\ell = \ell_1 \cup \cdots \cup \ell_n$ (Def.\ref{def:disjoint}) of single-component ${\rm SL}_3$-laminations $\ell_1,\ldots,\ell_n$, whose underlying non-elliptic ${\rm SL}_3$-webs in $\frak{S}$ are mutually non-isotopic, and each $\ell_i$ that is not a peripheral loop has weight $1$. Define $\wh{\mathbb{I}}^\omega_\Delta(\ell_i)$ as:
\begin{enumerate}
\itemsep0em
\item[\rm (Q1)] If $\ell_i$ consists of a peripheral  loop, then 
$$
\wh{\mathbb{I}}^\omega_\Delta(\ell_i) := [ {\textstyle \prod}_{v\in \mathcal{V}(Q_\Delta)} \wh{Z}_v^{{\rm a}_v(\ell_i)} ]_{\rm Weyl}
$$

\item[\rm (Q2)] Otherwise, if the underlying ${\rm SL}_3$-web in $\frak{S}$ for $\ell_i$ is $W_i$, let $\wh{W}_i$ be the ${\rm SL}_3$-web in $\frak{S} \times {\bf I}$ obtained by embedding $W_i$ at a constant elevation surface $\frak{S} \times \{c\}$, equipped with the upward vertical framing. Then,
$$
\wh{\mathbb{I}}^\omega_\Delta(\ell_i) := {\rm Tr}^\omega_\Delta([\wh{W}_i,{\O}])
$$
\end{enumerate}
Define
$$
\wh{\mathbb{I}}^\omega_\Delta(\ell) := \wh{\mathbb{I}}^\omega_\Delta(\ell_1) \cdot \cdots \cdot\, \wh{\mathbb{I}}^\omega_\Delta(\ell_n).
$$
By convention, set $\wh{\mathbb{I}}^\omega_\Delta({\O}):=1$.
\end{definition}
A basic observation:
\begin{lemma}
Let $\Delta$ be a triangulation of a punctured surface $\frak{S}$. If $\ell_1,\ell_2 \in \mathscr{A}_{\rm L}(\frak{S};\mathbb{Z})$ are disjoint, then $\wh{\mathbb{I}}^\omega_\Delta(\ell_1)\wh{\mathbb{I}}^\omega_\Delta(\ell_2) = \wh{\mathbb{I}}^\omega_\Delta(\ell_2) \wh{\mathbb{I}}^\omega_\Delta(\ell_1)$.
\end{lemma}
{\it Proof.} It suffices to prove this when each of $\ell_1,\ell_2$ is represented by a single-component ${\rm SL}_3$-web in $\frak{S}$. If both are non-peripheral, then the commutativity follows from the product structure of $\mathcal{S}^\omega_{\rm s}(\frak{S};\mathbb{Z})_{\rm red} \cong \mathcal{S}^\omega(\frak{S};\mathbb{Z})$ and the multiplicativity of ${\rm Tr}^\omega_\Delta$. When one of $\ell_1,\ell_2$ is a peripheral loop, it is a straightforward exercise (cf. \cite[Lem.3.9]{AK}). \qed

\begin{corollary}
\label{Cor:product_to_sum_quantum}
Let $\Delta$ be an ideal triangulation of a punctured surface $\frak{S}$. If $\ell,\ell' \in \mathscr{A}_{\rm L}(\frak{S};\mathbb{Z})$,
\begin{align}
\label{eq:product_to_sum_quantum1}
\wh{\mathbb{I}}^\omega_\Delta(\ell) \, \wh{\mathbb{I}}^\omega_\Delta(\ell') = \underset{\ell'' \in \mathscr{A}_{\rm L}(\frak{S};\mathbb{Z})}{\textstyle \sum}  \wh{c}^{\,\omega}(\ell,\ell';\ell'') \, \wh{\mathbb{I}}^\omega_\Delta(\ell''),
\end{align}
for some $\wh{c}^{\,\omega}(\ell,\ell';\ell'') \in \mathbb{Z}[\omega^{\pm 3}]$ that does not depend on $\Delta$, where the sum is a finite sum.
\end{corollary}

{\it Proof.} If $\ell$ or $\ell'$ consists only of peripheral loops, then $\wh{\mathbb{I}}^\omega_\Delta(\ell) \wh{\mathbb{I}}^\omega_\Delta(\ell') = \wh{\mathbb{I}}^\omega_\Delta(\ell \cup \ell')$, and we are done. Suppose each $\ell$ and $\ell'$ is represented by a single-component ${\rm SL}_3$-web in $\frak{S}$ with weight $1$, say $W$ and $W'$, whose constant-elevation lifts in $\mathcal{S}^\omega(\frak{S};\mathbb{Z})$ are denoted by $[\wh{W}]$ and $[\wh{W}']$. By using ${\rm SL}_3$-skein relations in Fig.\ref{fig:A2-skein_relations_quantum}, one obtains
\begin{align}
\label{eq:product_to_sum_quantum_duality_proof1}
[\wh{W}] \cdot [\wh{W}'] = {\textstyle \sum}_{[\wh{W}'']} \wh{c}^{\,\omega}(\ell,\ell';[\wh{W}'']) [\wh{W}''], 
\end{align}
which is a finite sum, with $\wh{c}^{\,\omega}(\ell,\ell' ; [\wh{W}'']) \in \mathbb{Z}[\omega^{\pm 3}]$, where each $\wh{W}''$ is an ${\rm SL}_3$-web in $\frak{S} \times {\bf I}$ without crossing. Thus $\wh{W}''$ is isotopic to a constant-elevation lift of an ${\rm SL}_3$-web $W''$ in $\frak{S}$. By using the relations in Fig.\ref{fig:A2-skein_relations_quantum}, one could assume that $W''$ has no internal 2-gon or 4-gon, so that $W''$ is non-elliptic, which is automatically reduced because $\frak{S}$ is without boundary. Thus $W''$ with weight $1$ forms an ${\rm SL}_3$-lamination $\ell'' \in \mathscr{A}_{\rm L}(\frak{S};\mathbb{Z})$. Writing $\wh{c}^{\,\omega}(\ell,\ell';[\wh{W}''])$ as $\wh{c}^{\,\omega}(\ell,\ell';\ell'')$, and applying ${\rm Tr}^\omega_\Delta$ to eq.\eqref{eq:product_to_sum_quantum_duality_proof1}, we get the desired result. 

\vs

For general $\ell,\ell'$, apply the above observation repeatedly, for each of their components. \qed

\vs

We now establish the quantum versions of Prop.\ref{prop:congruence_of_terms_of_SL3_classical_trace} and Prop.\ref{prop:highest_term_of_SL3_classical_trace_for_a_triangle}
\begin{proposition}[congruence of the Laurent monomial degrees of terms of the ${\rm SL}_3$ quantum trace]
\label{prop:congruence_of_terms_of_SL3_quantum_trace1}
Let $\Delta$ be an ideal triangulation of a triangulable generalized marked surface $\frak{S}$. For any stated ${\rm SL}_3$-web $(W,s)$ in $\frak{S} \times {\bf I}$, ${\rm Tr}^\omega_\Delta([W,s]) \in \mathcal{Z}^\omega_\Delta$ can be written as a $\omega^{1/2}$-Laurent polynomial (Def.\ref{def:non-commutative_Laurent_monomial_and_polynomial}) in the generators $\{\wh{Z}_v \, |\, v\in \mathcal{V}(Q_\Delta)\}$ of $\mathcal{Z}^\omega_\Delta$ so that all appearing $\omega^{1/2}$-Laurent monomials are congruent to each other: for any two $\omega^{1/2}$-Laurent monomials $\epsilon \, \omega^m [\prod_v \wh{Z}_v^{\alpha_v}]_{\rm Weyl}$ and $\epsilon' \omega^{m'} [\prod_v \wh{Z}_v^{\beta_v}]_{\rm Weyl}$ appearing, with $(\alpha_v)_v,(\beta_v)_v \in \mathbb{Z}^{\mathcal{V}(Q_\Delta)}$, $\epsilon,\epsilon'\in \{+1,-1\}$, $m,m' \in \frac{1}{2}\mathbb{Z}$, we have $\alpha_v - \beta_v \in 3\mathbb{Z}$ for all $v\in \mathcal{V}(Q_\Delta)$.
\end{proposition}

{\it Proof.} Works almost verbatim as in Prop.\ref{prop:congruence_of_terms_of_SL3_classical_trace}. \qed

\vs

A quantum version of Prop.\ref{prop:highest_term_of_SL3_classical_trace_for_a_triangle}: 

\begin{definition}
On the set of all $\omega^{1/2}$-Laurent monomials in $\{\wh{Z}_v \, | \, v\in \mathcal{V}(Q_\Delta)\}$ (Def.\ref{def:non-commutative_Laurent_monomial_and_polynomial}), define the preorder as follows: for $(\alpha_v)_{v\in \mathcal{V}(Q_\Delta)}, (\beta_v)_{v\in \mathcal{V}(Q_\Delta)} \in \mathbb{Z}^{\mathcal{V}(Q_\Delta)}$, $\epsilon,\epsilon' \in \{+1,-1\}$, $m,m'\in \frac{1}{2}\mathbb{Z}$,
$$
\epsilon \, \omega^m \textstyle \prod_v \wh{Z}_v^{\alpha_v} \succ \epsilon' \omega^{m'} \prod_v \wh{Z}_v^{\beta_v} \qquad \overset{\rm def.}{\Longleftrightarrow} \qquad \alpha_v \ge \beta_v, ~ \forall v\in \mathcal{V}(Q_\Delta).
$$
By convention, the zero monomial is set to be of the lowest preorder, i.e. $\epsilon \, \omega^m \prod_v \wh{Z}_v^{\alpha_v} \succ 0$.
\end{definition}

\begin{proposition}[the highest term of the ${\rm SL}_3$ quantum trace for a triangle]
\label{prop:highest_term_of_SL3_quantum_trace_for_a_triangle}
Let $t$ be a triangle, viewed as a generalized marked surface. Let $W$ be a canonical ${\rm SL}_3$-web in $t$ (Def.\ref{def:canonical_web_in_a_triangle}). Let $\Delta$ be the unique triangulation of $t$, so that $Q_\Delta$ has seven nodes. For each $v\in \mathcal{V}(Q_\Delta)$, let ${\rm a}_v(W) \in \frac{1}{3}\mathbb{Z}$ be the tropical coordinate of $W$ as defined in Def.\ref{def:tropical_coordinates}, when $W$ is viewed as an ${\rm SL}_3$-lamination in $t$ with weight $1$. There exists an ${\rm SL}_3$-web $\wh{W}$ in $t \times {\bf I}$ that projects to $W$, such that if we denote by ${\bf 1}_{\wh{W}}$ the state of $\wh{W}$ assigning $1 \in \{1,2,3\}$ to all endpoints of $\wh{W}$, then the following holds:
\begin{enumerate}
\item[\rm (QHT1)]  ${\rm Tr}^\omega_\Delta([\wh{W},{\bf 1}_{\wh{W}}]) \in \mathcal{Z}^\omega_\Delta = \mathcal{Z}^\omega_{\wh{t}}$ can be written as a $\omega^{1/2}$-Laurent polynomial in $\{ \wh{Z}_v \, | \, v\in \mathcal{V}(Q_\Delta)\}$ so that the unique $\omega^{1/2}$-Laurent monomial of the highest preorder is $\omega^m [\prod_{v\in \mathcal{V}(Q_\Delta)} \wh{Z}_v^{3{\rm a}_v(W)}]_{\rm Weyl}$ for some $m\in \frac{1}{2}\mathbb{Z}$.

\item[\rm (QHT2)] For any other state $s$ of $\wh{W}$, ${\rm Tr}^\omega_\Delta([\wh{W},s]) \in \mathcal{Z}^\omega_\Delta = \mathcal{Z}^\omega_t$ can be written as a $\omega^{1/2}$-Laurent polynomial in $\{\wh{Z}_v \, | \, v\in \mathcal{V}(Q_\Delta)\}$ so that each appearing Laurent monomial has strictly lower preorder than $[\prod_v \wh{Z}_v^{3{\rm a}_v(W)}]_{\rm Weyl}$.
\end{enumerate}
\end{proposition}

{\it Proof.} If we choose $\wh{W}$ appropriately, we claim that the proof of Prop.\ref{prop:highest_term_of_SL3_classical_trace_for_a_triangle} works almost verbatim, except that we cannot pin down $m$ in (QHT1). Let $W_1,\ldots,W_n$ be the components of $W$. We choose $\wh{W}$ such that the corresponding components $\wh{W}_1,\ldots,\wh{W}_n$ are located at mutually disjoint elevations (like in (GP2) of Def.\ref{def:good_position}), e.g. so that $\wh{W} = \wh{W}_1 \cdot \wh{W}_2 \cdot \cdots \cdot \wh{W}_n$. Then it suffices to show (QHT1) and (QHT2) for each component. So one can assume from the beginning that $W$ has only one component. When $W$ is a left turn or a right turn arc, then (QHT1) and (QHT2) are easy to observe from the values of ${\rm Tr}^\omega_t$ written in \redfix{(QT2-1)} and \redfix{(QT2-2)} of Thm.\ref{thm:SL3_quantum_trace_map}. Suppose now $W$ is a degree $d$ pyramid $H_d$ for a nonzero $d\in \mathbb{Z}$, as in the proof of Prop.\ref{prop:highest_term_of_SL3_classical_trace_for_a_triangle}. Call the three sides of $t$ as $e_1,e_2,e_3$, and the endpoints of $W$ as $x_1,\ldots,x_{|d|}$ (living in $e_1$), $y_1,\ldots,y_{|d|}$ (living in $e_2$), and $z_1,\ldots,z_{|d|}$ (living in $e_3$), exactly as in Fig.\ref{fig:pyramid_in_triangle_pushing}. Denote by the same labels the corresponding endpoints of $\wh{W}$. Note that $\wh{W}$ is completely determined by the choice of elevation orderings of endpoints living in each side of $t$. Choose $\wh{W}$ so that we have the elevation orderings as $x_1\succ x_2\succ \cdots \succ x_{|d|}$ over $e_1$, $y_1\succ \cdots \succ y_{|d|}$ over $e_2$, and $z_{|d|} \succ \cdots \succ z_1$ over $e_3$; so, an endpoint located more toward the `left' of Fig.\ref{fig:pyramid_in_triangle_pushing} has higher elevation. Then, apply the proof of Prop.\ref{prop:highest_term_of_SL3_classical_trace_for_a_triangle}, keeping mind this strategy of \redfix{choosing higher elevation} for a juncture or endpoint located more toward the left of the figures. More specifically, in the right picture of Fig.\ref{fig:K_d_in_biangle} we choose the elevations for the inner junctures so that $r_1\succ r_2\succ \cdots \succ r_{|d|-1} \succ r \succ r'$, and in the right picture of Fig.\ref{fig:pyramid_in_triangle_pushing} we choose the elevations for the inner junctures so that $w_1\succ \cdots \succ w_{|d|} \succ u_1 \succ \cdots u_{|d|}$. Then indeed, for each ${\rm SL}_3$-web appearing during the proof, the components are at mutually disjoint elevations, hence decomposes as product of the components as assumed in the proof of Prop.\ref{prop:highest_term_of_SL3_classical_trace_for_a_triangle}, so that a similar proof works, to yield the desired result, because the analysis of signs and non-zero-ness of the entries (or values) in Prop.\ref{prop:biangle_SL3_quantum_trace}(BT2-1), (BT2-3) and Lem.\ref{lem:value_of_I-webs_in_biangle} (as well as Thm.\ref{thm:SL3_quantum_trace_map}\redfix{(QT2-1)}--\redfix{(QT2-2))} is the same for the quantum setting and the classical setting, by inspection. \qed

\vs

In the investigation of the highest term of the ${\rm SL}_3$ quantum trace over the entire surface $\frak{S}$, what plays a crucial role is Prop.\ref{prop:elevation_reversing_and_star-structure}, which is the equivariance under the elevation reversing and the $*$-structure. We give a proof only now, using the state-sum formula.

\vs

{\it Proof of Prop.\ref{prop:elevation_reversing_and_star-structure}.} Let $\Delta$ be an ideal triangulation of a generalized marked surface $\frak{S}$, and let $[W,s] \in\mathcal{S}^\omega_{\rm s}(\frak{S};\mathbb{Z})_{\rm red}$. Let $\wh{\Delta}$ be a split ideal triangulation for $\Delta$, and put $(W,s)$ into a gool position with respect to $\wh{\Delta}$ through an isotopy (Lem.\ref{lem:isotopic_to_good_position}). Consider the state-sum trace $\wh{\rm Tr}^\omega_\Delta(W,s)$, as in eq.\eqref{eq:state-sum_formula} of Def.\ref{def:state-sum_trace_for_gool_position}):
\begin{align}
\nonumber
\wh{\rm Tr}^\omega_\Delta(W,s) = {\textstyle \sum}_J ( {\textstyle \prod}_B {\rm Tr}^\omega_B([W\cap (B\times {\bf I}), J_B]) \, {\textstyle \bigotimes}_t \wh{\rm Tr}^\omega_t(W\cap (\wh{t} \times {\bf I}), J_t) ) ~\in~ 
{\textstyle \bigotimes}_{t\in \mathcal{F}(\Delta)} \mathcal{Z}^\omega_t
\end{align}
The equivariance for the biangle factors is shown in Lem.\ref{lem:elevation_reversing_and_star-structure_for_biangles}. It remains to study the effect of elevation reversing on the triangle factors. Note that the elevation reversing ${\bf r}$ reverses the elevation order of the components of $W\cap (\wh{t} \times {\bf I})$, hence reverses the product order of factors of $\wh{\rm Tr}^\omega_t(W\cap (\wh{t} \times {\bf I}), J_t) )$. 

\vs

So it suffices to show the equivariance for a single left or right turn arc over $\wh{t}$. For a stated ${\rm SL}_3$-skein $[W,s] \in \mathcal{S}^\omega_{\rm s}(\wh{t} \,;\mathbb{Z})_{\rm red}$ that is a single left or right turn over $\wh{t}$, the elevation reversal yields a same stated ${\rm SL}_3$-skein, i.e. ${\bf r}[W,s]  = [W,s]$, while we know from Lem.\ref{lem:quantum_turn_matrices_are_Weyl-ordered} that ${\rm Tr}^\omega_t([W,s])$ is Weyl-ordered, hence $* ({\rm Tr}^\omega_t([W,s])) = {\rm Tr}^\omega_t([W,s])$ by Lem.\ref{lem:Weyl-ordered_product_is_invariant_under_star-map}. Thus indeed the equivariance holds for a single left or right turn over $\wh{t}$. To summarize, we have shown $\wh{\rm Tr}^\omega_\Delta({\bf r}(W,s)) = *(\wh{\rm Tr}^\omega_\Delta(W,s))$, thus ${\rm Tr}^\omega_\Delta({\bf r}[W,s]) = *({\rm Tr}^\omega_\Delta([W,s]))$, as desired. \qed

\vs

One last step before the quantum highest term statement is the quantum version of Lem.\ref{lem:crossbar_web_value}:
\begin{lemma}
\label{lem:crossbar_web_value_quantum}
Let $W$ be a crossbar ${\rm SL}_3$-web in a biangle $B$ (Def.\ref{def:crossbar}), $\wh{W}$ be any ${\rm SL}_3$-web in $B\times {\bf I}$ that projects to $W$, and ${\bf 1}_{\wh{W}}$ be the state of $\wh{W}$ assigning the value $1\in \{1,2,3\}$ to all the endpoints of $\wh{W}$. Then ${\rm Tr}_B([\wh{W},{\bf 1}_W])= \omega^m$ for some $m\in \frac{1}{2}\mathbb{Z}$.
\end{lemma}

{\it Proof.} Works almost verbatim as in Lem.\ref{lem:crossbar_web_value}. Due to the presence of the elevation orderings, in each elementary ${\rm SL}_3$-web over a biangle $B_i$, one may see some height-exchange ${\rm SL}_3$-web or a single-crossing ${\rm SL}_3$-webs as in Prop.\ref{prop:biangle_SL3_quantum_trace_some_values}(BT2-4) or eq.\eqref{eq:some_more_elementary1}--\eqref{eq:some_more_elementary3}, in addition to the single-crossbar (i.e. the `H-web' $W'$ of Lem.\ref{lem:value_of_I-webs_in_biangle}). Still, looking at the values for these cases, the induction proof of Lem.\ref{lem:crossbar_web_value} works the same way. \qed

\vs

We finally obtain a quantum version of Prop.\ref{prop:highest_term_of_SL3_classical_trace}, the highest term statement.
\begin{proposition}[the highest term of the ${\rm SL}_3$ quantum trace]
\label{prop:highest_term_of_SL3_quantum_trace}
Let $\Delta$ be an ideal triangulation of a triangulable generalized marked surface $\frak{S}$, and $\wh{\Delta}$ be a split ideal triangulation for $\Delta$. Let $W$ be a (reduced) non-elliptic ${\rm SL}_3$-web in $\frak{S}$ in a canonical position with respect to $\wh{\Delta}$ (Def.\ref{def:canonical_wrt_split_ideal_triangulation}) that has no endpoints. View $W$ as an ${\rm SL}_3$-lamination by giving the weight $1$; let ${\rm a}_v(W) \in \frac{1}{3}\mathbb{Z}$, $v\in \mathcal{V}(Q_\Delta)$, be the tropical coordinates defined in Def.\ref{def:tropical_coordinates}. Let $\wh{W}$ be an ${\rm SL}_3$-web in $\frak{S} \times {\bf I}$ that projects to $W$. 
Then ${\rm Tr}^\omega_\Delta([\wh{W},{\O}]) \in \mathcal{Z}^\omega_\Delta$ can be written as a $\omega^{1/2}$-Laurent polynomial in $\{\wh{Z}_v \, | \, v\in \mathcal{V}(Q_\Delta)\}$ so that $[\prod_v \redfix{\wh{Z}}_v^{3{\rm a}_v(W)}]_{\rm Weyl}$\redfix{$=[\prod_v \wh{X}_v^{{\rm a}_v(W)}]_{\rm Weyl}$} is the unique $\omega^{1/2}$-Laurent monomial of the highest preorder.
\end{proposition}

{\it Proof.} Almost verbatim proof as that of Prop.\ref{prop:highest_term_of_SL3_classical_trace}, if we are a bit careful when choosing the elevations of the junctures. Namely, we should choose (the elevations of the points of) $\wh{W}$ so that for each triangle $\wh{t}$ of $\wh{\Delta}$,  $\wh{W} \cap (\wh{t} \times {\bf I})$ is as in the proof of Prop.\ref{prop:highest_term_of_SL3_quantum_trace_for_a_triangle}; this is certainly possible (e.g. by applying a `vertical isotopy' to any chosen $\wh{W}$), and under such a choice,  Prop.\ref{prop:highest_term_of_SL3_quantum_trace_for_a_triangle} holds, which together with Lem.\ref{lem:crossbar_web_value_quantum} makes a similar proof as in Prop.\ref{prop:highest_term_of_SL3_classical_trace} to work. And note that ${\rm Tr}^\omega_\Delta([\wh{W},{\O}])$ depends only on $W$, due to the isotopy invariance of ${\rm Tr}^\omega_\Delta$. As a result, we get that the unique highest term of ${\rm Tr}^\omega_\Delta([\wh{W},{\O}])$ is $\omega^m [\prod_v \redfix{\wh{Z}}_v^{3{\rm a}_v(W)}]_{\rm Weyl}$ for some $m\in \frac{1}{2} \mathbb{Z}$. Since $[\wh{W},{\O}]$ equals its elevation reversed version ${\bf r}[\wh{W},{\O}]$, it follows from Prop.\ref{prop:elevation_reversing_and_star-structure} that $*({\rm Tr}^\omega_\Delta([\wh{W},{\O}])) = {\rm Tr}^\omega_\Delta([\wh{W},{\O}])$. By definition of the $*$-map (Def.\ref{def:star-structure}), note that the highest term of $*({\rm Tr}^\omega_\Delta([\wh{W},{\O}]))$ equals the image under the $*$-map of the highest term of ${\rm Tr}^\omega_\Delta([\wh{W},{\O}])$. Thus it follows that the highest term of ${\rm Tr}^\omega_\Delta([\wh{W},{\O}])$ is $*$-invariant. Hence by Lem.\ref{lem:Weyl-ordered_product_is_invariant_under_star-map} it follows that $m=1$. \qed

\vs

Before proceeding, one obtains the following useful corollary or each punctured surface $\frak{S}$ (answering a question posed to \redfix{the author} by Vijay Higgins), using an argument similar to the proof of Cor.\ref{cor:congruence_and_integrality_of_powers}, keeping in mind that the non-elliptic ${\rm SL}_3$-webs form a basis of $\mathcal{S}^\omega_{\rm s}(\frak{S};\mathbb{Z})_{\rm red} \cong \mathcal{S}^\omega(\frak{S};\mathbb{Z})$ (\cite{SW} \cite[Thm.2]{FS}).
\begin{corollary}
For a triangulable punctured surface $\frak{S}$ and an ideal triangulation $\Delta$ of $\frak{S}$, the ${\rm SL}_3$ quantum trace map ${\rm Tr}^\omega_\Delta$ is injective. \qed
\end{corollary}
\redfix{Combining Prop.\ref{prop:highest_term_of_SL3_quantum_trace} and Prop.\ref{prop:congruence_of_terms_of_SL3_quantum_trace1} one obtains:}
\begin{corollary}[congruence of the Laurent monomial degrees of terms of the ${\rm SL}_3$ quantum trace \redfix{values at non-elliptic ${\rm SL}_3$-webs}]
\label{cor:congruence_of_terms_of_SL3_quantum_trace2}
Let $\Delta$ be an ideal triangulation of a triangulable generalized marked surface $\frak{S}$. Let $W$ and $\wh{W}$ be as in Prop.\ref{prop:highest_term_of_SL3_quantum_trace}. Then ${\rm Tr}^\omega_\Delta([\wh{W},{\O}]) \in \mathcal{Z}^\omega_\Delta$ can be written as a $\omega^{1/2}$-Laurent polynomial in $\{\wh{Z}_v \, |\, v\in \mathcal{V}(Q_\Delta)\}$ so that each appearing $\omega^{1/2}$-Laurent monomial is of the form $\pm \redfix{\omega}^m [\prod_v \wh{X}_v^{\alpha_v}\redfix{\wh{Z}}_v^{3{\rm a}_v(W)}]_{\rm Weyl}$ for some $(\alpha_v)_v\in \mathbb{Z}^{\mathcal{V}(Q_\Delta)}$ and $m \in \frac{1}{2}\mathbb{Z}$. 
\end{corollary}
\redfix{We conjecture that $m \in \frac{1}{2}\mathbb{Z}$ appearing in Cor.\ref{cor:congruence_of_terms_of_SL3_quantum_trace2} belongs to $9\mathbb{Z}$, so that $\omega^m$ is an integer power of $q$. We leave this as a future research problem. Perhaps, a good approach would be to interpret this power $m$ as a writhe of a lift of $W$ in a certain branched 3-fold cover of $\frak{S}$, i.e. to compare our ${\rm SL}_3$ quantum trace with the quantum holonomy construction in \cite{Gabella}; see \cite{KLS}, where such a comparison is made rigorous between the Bonahon-Wong ${\rm SL}_2$ quantum trace \cite{BW} and the quantum holonomy of \cite{Gabella}.}

\vs

Combining all above, we obtain a statement for the quantum duality map.
\begin{theorem}
\label{thm:quantum_duality_map}
Let $\frak{S}$ be a triangulable punctured surface, and $\Delta$ be an ideal triangulation of $\frak{S}$. Denote by $\mathcal{X}^q_\Delta$ the Fock-Goncharov algebra, defined as the free associative $\mathbb{Z}[\redfix{q^{\pm 1/18}}]$-algebra generated by $\{ \wh{X}_v^{\pm 1} \, | \, v\in \mathcal{V}(Q_\Delta)\}$ mod out by the relations
$$
\wh{X}_v \wh{X}_w = q^{2\varepsilon_{vw}} \wh{X}_w \wh{X}_v, \qquad \forall v,w\in \mathcal{V}(Q_\Delta).
$$
Restricting $\wh{\mathbb{I}}^\omega_\Delta$ (Def.\ref{def:quantum_duality_map_omega}) to $\mathscr{A}_{{\rm SL}_3,\frak{S}}(\mathbb{Z}^{\redfix{T}}) \subset \mathscr{A}_{\rm L}(\frak{S};\mathbb{Z})$, one obtains a \ul{\em quantum ${\rm SL}_3$-${\rm PGL}_3$ duality map}
$$
\mathbb{I}^q_\Delta : \mathscr{A}_{{\rm SL}_3,\frak{S}}(\mathbb{Z}^{\redfix{T}}) \to \mathcal{X}^q_\Delta,
$$
satisfying
\begin{enumerate}
\item[\rm (1)] When $q^{1/18}=1$, this map $\mathbb{I}^q_\Delta$ recovers the classical duality map $\mathbb{I}$ of Thm.\ref{thm:main}.

\item[\rm (2)] For $\ell \in \mathscr{A}_{{\rm SL}_3,\frak{S}}(\mathbb{Z}^{\redfix{T}})$, $*(\mathbb{I}^q_\Delta(\ell)) = \mathbb{I}^q_\Delta(\ell)$.

\item[\rm (3)] For $\ell \in \mathscr{A}_{{\rm SL}_3,\frak{S}}(\mathbb{Z}^{\redfix{T}})$, the unique highest Laurent monomial of $\mathbb{I}^q_\Delta(\ell)$ is $[\prod_{v\in \mathcal{V}(Q_\Delta)} \wh{X}_v^{{\rm a}_v(\ell)}]_{\rm Weyl}$.

\item[\rm (4)] If $\ell \in \mathscr{A}_{{\rm SL}_3,\frak{S}}(\mathbb{Z}^{\redfix{T}})$ consists only of peripheral loops, then $\mathbb{I}^q_\Delta(\ell)= [\prod_{v\in \mathcal{V}(Q_\Delta)} \wh{X}_v^{{\rm a}_v(\ell)}]_{\rm Weyl}$.

\item[\rm (5)] For any $\ell,\ell' \in \mathscr{A}_{{\rm SL}_3,\frak{S}}(\mathbb{Z}^{\redfix{T}})$, we have
\begin{align}
\nonumber
\mathbb{I}^q_\Delta(\ell) \, \mathbb{I}^q_\Delta(\ell') = \underset{\ell'' \in \mathscr{A}_{{\rm SL}_3,\frak{S}}(\mathbb{Z}^{\redfix{T}})}{\textstyle \sum}  c^q(\ell,\ell';\ell'') \, \mathbb{I}^q_\Delta(\ell'')
\end{align}
where $c^q(\ell,\ell';\ell'')\in \mathbb{Z}[\redfix{q^{\pm 1/3}}]$ are zero for all but finitely many $\ell''$, and do not depend on $\Delta$.

\item[\rm (6)] If $\ell,\ell' \in \mathscr{A}_{{\rm SL}_3,\frak{S}}(\mathbb{Z}^{\redfix{T}})$ and if $\ell$ consists only of peripheral loops, then $\mathbb{I}^q_\Delta(\ell) \mathbb{I}^q_\Delta(\ell') = \mathbb{I}^q_\Delta(\ell \cup \ell') = \mathbb{I}^q_\Delta(\ell') \mathbb{I}^q_\Delta(\ell)$.

\end{enumerate}
\end{theorem}

{\it Proof.} We do not provide all details for a proof. A general strategy is to first show the corresponding statements for $\wh{\mathbb{I}}^\omega_\Delta$, by deducing them from those for 
the ${\rm SL}_3$ quantum trace ${\rm Tr}^\omega_\Delta$. For example, for $\ell \in \mathscr{A}_{{\rm SL}_3,\frak{S}}(\mathbb{Z}^{\redfix{T}})$, the fact that $\mathbb{I}^q_\Delta(\ell)$ lies in $\mathcal{X}^q_\Delta$ can be deduced from a version of \redfix{Cor.\ref{cor:congruence_of_terms_of_SL3_quantum_trace2}}. The only part that needs a care is the item (5). For $\ell,\ell' \subset\mathscr{A}_{{\rm SL}_3,\frak{S}}(\mathbb{Z}^{\redfix{T}}) \subset \mathscr{A}_{\rm L}(\frak{S};\mathbb{Z})$, what we do have is eq.\eqref{eq:product_to_sum_quantum1} of Cor.\ref{Cor:product_to_sum_quantum}. From a version of the statement (3) on the highest term for $\wh{\mathbb{I}}^\omega_\Delta(\ell)$, $\wh{\mathbb{I}}^\omega_\Delta(\ell')$, $\wh{\mathbb{I}}^\omega_\Delta(\ell'')$ (for $\ell'' \in \mathscr{A}_{\rm L}(\frak{S};\mathbb{Z})$), a version of \redfix{Cor.\ref{cor:congruence_of_terms_of_SL3_quantum_trace2}} for these, and an observation that the highest term of the product of $\wh{\mathbb{I}}^\omega_\Delta(\ell)$ and $\wh{\mathbb{I}}^\omega_\Delta(\ell')$ equals the product of the highest terms of them which in turn equals $[\prod_v \wh{X}_v^{{\rm a}_v(\ell)}]_{\rm Weyl} [\prod_v \wh{X}_v^{{\rm a}_v(\ell')}]_{\rm Weyl} = q^m [\prod_v \wh{X}_v^{{\rm a}_v(\ell)+{\rm a}_v(\ell')}]_{\rm Weyl}$ for some integer $m$, by applying a similar argument as in our proof in \S\ref{subsec:basis} of Cor.\ref{cor:congruence_and_integrality_of_powers}, one can show that in the right hand side of eq.\eqref{eq:product_to_sum_quantum1} (for the case when $\ell,\ell'$ in the left hand side belong to $\mathscr{A}_{{\rm SL}_3,\frak{S}}(\mathbb{Z}^{\redfix{T}})$), the $\ell''$ that have nonzero contributions all belong to $\mathscr{A}_{{\rm SL}_3,\frak{S}}(\mathbb{Z}^{\redfix{T}})$, and the coefficients $\wh{c}^{\,\omega}(\ell,\ell';\ell'')$ belong to $\mathbb{Z}[\redfix{q^{\pm 1/18}}]$. In particular, for the unique $\ell''_0$ s.t. $[\prod_v \wh{X}_v^{{\rm a}_v(\ell''_0)}]_{\rm Weyl}$ equals the highest term of the right hand side of eq.\eqref{eq:product_to_sum_quantum1}, we have $\wh{c}^{\,\omega}(\ell,\ell';\ell''_0) = q^m$. These coefficients $\wh{c}^{\,\omega}(\ell,\ell';\ell'')$ \redfix{in fact} belong to $\redfix{\mathbb{Z}[\omega^{\pm 3}]=}\mathbb{Z}[q^{\pm 1/3}]$, from Cor.\ref{Cor:product_to_sum_quantum}. \qed

\begin{conjecture}
\redfix{Cor.\ref{cor:congruence_of_terms_of_SL3_quantum_trace2}} holds with $m\in \redfix{9}\mathbb{Z}$, and hence the $\mathbb{Z}[\redfix{q^{\pm 1/18}}]$ in the statement of Thm.\ref{thm:quantum_duality_map} can be replaced by $\mathbb{Z}[q^{\pm 1}]$, i.e. $\mathbb{I}^q_\Delta(\ell)$ for $\ell \in \mathscr{A}_{{\rm SL}_3,\frak{S}}(\mathbb{Z}^{\redfix{T}})$ is a non-commutative Laurent polynomial in $\{\wh{X}^{\pm 1}_v \, | \, v\in \mathcal{V}(Q_\Delta)\}$ with coefficients in $\mathbb{Z}[q^{\pm 1}]$.
\end{conjecture}

\begin{proposition}[\cite{Kim21}]
\label{prop:compatibility_of_I_q}
\redfix{$\mathbb{I}^q_\Delta$ and $\mathbb{I}^q_{\Delta'}$ for different triangulations are related by the quantum coordinate change map between ${\rm Frac}(\mathcal{X}^q_\Delta)$ and ${\rm Frac}(\mathcal{X}^q_{\Delta'})$ as mentioned in Prop.\ref{prop:quantum_coordinate_change}:}
$$
\redfix{\mathbb{I}^q_\Delta = \Phi^q_{\Delta\Delta'} \circ \mathbb{I}^q_{\Delta'}}
$$
\end{proposition}
\redfix{As said, Prop.\ref{prop:quantum_coordinate_change} is shown in \cite{Kim21}; in addition, a version of Prop.\ref{prop:quantum_coordinate_change} for a peripheral loop (i.e. just about the highest term) is also developed in \cite{Kim21}, which together with Prop.\ref{prop:quantum_coordinate_change} implies the above Prop.\ref{prop:compatibility_of_I_q}, which used to be Conjecture 5.85 in a previous version of the present paper (ver3). With the help of Prop.\ref{prop:compatibility_of_I_q},} we obtain a deformation quantization map $\mathscr{O}(\mathscr{X}_{{\rm PGL}_3,\frak{S}}) 
\to \mathscr{O}^q_{\redfix{\rm tri}}(\mathscr{X}_{{\rm PGL}_3,\frak{S}})$ of the moduli space $\mathscr{X}_{{\rm PGL}_3,\frak{S}}$, as discussed in \S\ref{subsec:intro_quantum}. \redfix{A natural conjecture is whether our construction yields a deformation quantization map from $\mathscr{O}(\mathscr{X}_{{\rm PGL}_3,\frak{S}})$ to $\mathscr{O}^q(\mathscr{X}_{{\rm PGL}_3,\frak{S}})$ which is the ring of quantum universally Laurent elements (for all quantum cluster seeds), not just to $\mathscr{O}^q_{\rm tri}(\mathscr{X}_{{\rm PGL}_3,\frak{S}})$ which is the ring of quantum universally Laurent elements only for ideal triangulations. That is:}

\redfix{\begin{conjecture}
$\forall \ell \in \mathscr{A}_{{\rm SL}_3,\frak{S}}(\mathbb{Z}^T)$, the element $\mathbb{I}^q_\Delta(\ell)$ is (quantum) Laurent in every quantum cluster seed, that is, it stays Laurent after applying an arbitrary sequence of quantum mutations.
\end{conjecture}}

\section{Conjectures}
\label{sec:conjectures}

We state some naturally arising conjectures and questions, besides those which already appeared in the text.
\begin{question}
If a rational function $f$ on $\mathscr{X}_{{\rm PGL}_3,\frak{S}}$ is regular on all the cluster $\mathscr{X}$-charts for ideal triangulations of $\frak{S}$, then is it regular on all cluster $\mathscr{X}$-charts, hence is a regular function on $\mathscr{X}_{{\rm PGL}_3\redfix{,\frak{S}}}$? That is, does universally Laurent for triangulations imply universally Laurent for all cluster $\mathscr{X}$-charts?
\end{question}
It might be convenient to have an affirmative answer to the above, but maybe it is more natural for us to consider more general class of ideal triangulations, called tagged ideal triangulations (in particular incorporating the self-folded triangles). For these, the construction of Fock-Goncharov $\mathscr{X}$-coordinates must be modified; see \cite[\S10.7]{FG06} for a discussion, \cite{AB} for the ${\rm SL}_2$ case\redfix{, and \cite{FP} for a higher rank version of tagging}; see also \cite{JK}. Next, we consider:

\begin{conjecture}[Laurent coefficient positivity]
\label{conj:Laurent_coefficient_positivity}
For each $\ell \in \mathscr{A}_{{\rm SL}_3,\frak{S}}(\mathbb{Z}^{\redfix{T}})$, the basic regular function $\mathbb{I}(\ell) \in \mathscr{O}(\mathscr{X}_{{\rm PGL}_3,\frak{S}})$ can be written, for any cluster $\mathscr{X}$-chart, as a Laurent polynomial in the generators with {\em non-negative} integer coefficients.
\end{conjecture}
We have a partial result, due to our state-sum formula:
\begin{proposition}
Conjecture \ref{conj:Laurent_coefficient_positivity} holds for a cluster $\mathscr{X}$-chart for an ideal triangulation $\Delta$ of a punctured surface $\frak{S}$, for each $\ell$ that can be represented by a non-elliptic ${\rm SL}_3$-web in a canonical position with respect to $\wh{\Delta}$ such that there is at most one internal 3-valent vertex in each triangle of $\wh{\Delta}$ (i.e. degree of the pyramid in each triangle is in $\{-1,0,1\}$) and no internal 3-valent vertex in biangles of $\wh{\Delta}$.
\end{proposition}
To prove the full version for cluster $\mathscr{X}$-chart for all ideal triangulations, one must \redfix{try to} show \redfix{for example} that the values of the ${\rm SL}_3$ classical trace for all pyramids $H_d$ in a triangle are Laurent polynomials with non-negative coefficients. For that, one should use results established in \S\ref{sec:SL3_trace}\redfix{; however, unfortunately, such positivity does not seem to hold in a triangle (already for $d=2$), so one might have to look for some other idea}. Once Conjecture \ref{conj:Laurent_coefficient_positivity} is settled, then one can try to check whether $\mathbb{I}(\ell)$ are {\em extremal}, i.e. cannot be a sum of two universally positive Laurent functions (as predicted in Conjecture 12.3 of \cite{FG06}). 

\vs

Another perhaps more important kind of positivity is the following:
\begin{question}[structure constant positivity]
\label{question:structure_constant_positivity}
Does our $A_2$-bangles basis of $\mathscr{O}(\mathscr{X}_{{\rm PGL}_3,\frak{S}})$ have non-negative structure constants? Namely, are the structure constants $c(\ell,\ell';\ell'')$ in eq.\eqref{eq:I_structure_constants} of Thm.\ref{thm:main} are non-negative?
\end{question}
One can ask similar question for the basis of $\mathscr{O}(\mathscr{X}_{{\rm SL}_3,\frak{S}})$ in Def.\ref{def:I_SL3} or the basis of $\mathscr{O}(\mathscr{L}_{{\rm SL}_3,\frak{S}})$ as in Cor.\ref{cor:A2-bangles_basis_for_L_SL3} and Cor.\ref{cor:structure_constants_SL3_0}. Then the question is related to a similar question for the basis of the ${\rm SL}_3$-skein algebra $\mathcal{S}(\frak{S};\mathbb{Z})$, consisting of non-elliptic ${\rm SL}_3$-webs (\cite{SW}). Note that such a positivity holds true for Fock-Goncharov's basis of $\mathscr{O}(\mathscr{X}_{{\rm PGL}_2,\frak{S}})$, and a core part of the proof relies on the corresponding positivity of a certain basis of the ${\rm A}_1$-type skein algebra (i.e. the usual Kauffman bracket skein algebra), proved by Dylan Thurston \cite{DThurston}. Notice that one important aspect of this statement and proof for the $A_1$ type (or ${\rm SL}_2$) is that the positive basis is not a bangles basis, but is a bracelets basis (see Def.\ref{def:bangles_and_bracelets}). So we propose a new basis that is an $A_2$-analog of the bracelets basis of Kauffman bracket skein algebra.
\begin{definition}
Define a map
$$
\mathbb{I}^{\rm br}_{{\rm SL}_3} : \mathscr{A}_{\rm L}(\frak{S};\mathbb{Z}) \to \mathscr{O}(\mathscr{X}_{{\rm SL}_3,\frak{S}})
$$
as in Def.\ref{def:I_SL3} with the following modification. Let $\ell \in \mathscr{A}_{\rm L}(\frak{S};\mathbb{Z})$, and let $\ell = \ell_1 \cup \cdots \cup \ell_n$ be the disjoint union of single-component ${\rm SL}_3$-webs that are mutually non-isotopic. Define $\mathbb{I}^{\rm br}_{{\rm SL}_3}(\ell_i)$ as $\mathbb{I}_{{\rm SL}_3}(\ell_i)$ unless $\ell_i$ consists of a non-peripheral loop $\gamma_i$, say with weight $k_i \in \mathbb{Z}_{>0}$, in which case we define 
$$
\mathbb{I}_{{\rm SL}_3}(\ell_i) := F^* \Phi([W_{\gamma_i}^{(k_i)}]),
$$
where $W_{\gamma_i}^{(k_i)}$ is as in Def.\ref{def:bangles_and_bracelets}, that is, as the trace of monodromy along $\gamma_i^{k_i}$, the $k_i$-time-going-along-$\gamma_i$. Define $\mathbb{I}_{{\rm SL}_3}^{\rm br}(\ell) := \mathbb{I}_{{\rm SL}_3}(\ell_1) \,\cdot  \cdots \cdot \, \mathbb{I}_{{\rm SL}_3}(\ell_n)$.
\end{definition}
\begin{conjecture}
$\mathbb{I}^{\rm br}_{{\rm SL}_3}$ is injective, and its image forms a basis of $\mathscr{O}(\mathscr{X}_{{\rm SL}_3,\frak{S}})$.
\end{conjecture}
It is not hard to prove this conjecture by showing an analogous statement for the ${\rm SL}_3$-skein algebra, by observing that the `base change' transformation between $\mathbb{I}_{{\rm SL}_3}$ and $\mathbb{I}^{\rm br}_{{\rm SL}_3}$ is `upper triangular'. Let's denote the resulting basis $\mathbb{I}^{\rm br}_{{\rm SL}_3}(\mathscr{A}_{\rm L}(\frak{S};\mathbb{Z}))$ of $\mathscr{O}(\mathscr{X}_{{\rm SL}_3,\frak{S}})$ by \ul{\em $A_2$-bracelets basis} of $\mathscr{O}(\mathscr{X}_{{\rm SL}_3,\frak{S}})$.
\begin{conjecture}
By mimicking the present paper's construction $\mathbb{I}_{{\rm SL}_3} \leadsto \mathbb{I}^+_{{\rm PGL}_3} \leadsto \mathbb{I}$ to $\mathbb{I}^{\rm br}_{{\rm SL}_3}$, one can obtain a map
$$
\mathbb{I}^{\rm br} : \mathscr{A}_{{\rm SL}_3,\frak{S}}(\mathbb{Z}^{\redfix{T}}) \to \mathscr{O}(\mathscr{X}_{{\rm PGL}_3,\frak{S}}),
$$
which is injective and whose image forms a basis of $\mathscr{O}(\mathscr{X}_{{\rm PGL}_3,\frak{S}})$.
\end{conjecture}
We call this conjectural basis $\mathbb{I}^{\rm br}(\mathscr{A}_{{\rm SL}_3,\frak{S}}(\mathbb{Z}^{\redfix{T}}))$ the \ul{\em $A_2$-bracelets basis} of $\mathscr{O}(\mathscr{X}_{{\rm PGL}_3,\frak{S}})$, by a slight abuse of notation. We then formulate the positivity conjecture.
\begin{conjecture}
$\mathbb{I}^{\rm br}$ (or $\mathbb{I}$) satisfies the Laurent coefficient positivity as in Conjecture \ref{conj:Laurent_coefficient_positivity}.
\end{conjecture}
\begin{conjecture}
$\mathbb{I}^{\rm br}$ (or $\mathbb{I}$)  satisfies the structure constant positivity as in Question \ref{question:structure_constant_positivity}.
\end{conjecture}
Meanwhile, the work of Gross-Hacking-Keel-Kontsevich \cite{GHKK} yields a duality map
$$
\mathbb{I}^{\rm GHKK} : \mathscr{A}_{ |Q_\Delta| }(\mathbb{Z}^{\redfix{T}}) \to \mathscr{O}( \mathscr{X}_{|Q_\Delta|})
$$
where $\mathscr{A}_{ |Q_\Delta| }$ and $\mathscr{X}_{|Q_\Delta|}$ denote the cluster $\mathscr{A}$- and $\mathscr{X}$- varieties for the quiver mutation equivalence class $|Q_\Delta|$ of the exchange matrix of the $3$-triangulation quiver $Q_\Delta$ for a triangulation $\Delta$ of a generalized marked surface $\frak{S}$. The definition is quite combinatorial and computationally much involved (and abstract), and the existence is guaranteed by the result of Goncharov-Shen \cite{GS18}. \redfix{In fact, as of now, a crucial ingredient called a `consistent scattering diagram' is only known to exist, but no explicit enough construction is known. This makes $\mathbb{I}^{\rm GHKK}$ not completely constructive in some sense, even for the simplest surfaces like once-punctured torus. In any case,} of course, a natural question is whether their duality map equals ours, which is much more geometric and intuitive. 
\begin{conjecture}
\label{conj:GHKK}
$\mathbb{I}^{\rm GHKK}$ coincides with $\mathbb{I}^{\rm br}$ (or with $\mathbb{I}$).
\end{conjecture}
In fact, even for ${\rm SL}_2$-${\rm PGL}_2$ (or $A_1$), Gross-Hacking-Keel-Kontsevich's duality map is not known to coincide with Fock-Goncharov's \cite{FG06}, except for couple of small surfaces; \redfix{recently, Mandel and Qin announced that they proved that the two ${\rm SL}_2$-${\rm PGL}_2$ duality maps coincide, in an upcoming paper \cite{MQ}}. For ${\rm SL}_3$-${\rm PGL}_3$ (or $A_2$), one may get some evidence by computing some examples (but computation of $\mathbb{I}^{\rm GHKK}$ is almost impossible at the moment!). \redfix{When tackling Conjecture \ref{conj:GHKK}, perhaps it will help to investigate the lowest terms of $\mathbb{I}_\Delta(\ell)$ and $\mathbb{I}^q_\Delta(\ell)$; we claim that, in case $\ell$ has no peripheral loops, similar arguments of the present paper would show that they are $\prod_{v\in\mathcal{V}(Q_\Delta)} X_v^{-{\rm a}_v(\ol{\ell})}$ and $[\prod_{v\in\mathcal{V}(Q_\Delta)} \wh{X}_v^{-{\rm a}_v(\ol{\ell})}]_{\rm Weyl}$, where $\ol{\ell}$ is the ${\rm SL}_3$-lamination obtained from $\ell$ by reversing the orientations of all components. A more thorough investigation of our duality map (or its slight modification) using the lowest terms will be done in an upcoming joint paper with Linhui Shen \cite{KS}.}

\vs

\redfix{Note that $\mathbb{I}$ and $\mathbb{I}^{\rm br}$ yield different bases that share many properties, and it's not clear at the moment which one is more canonical. An answer to Conjecture \ref{conj:GHKK} would settle such a canonicity. One may think that, in view of the ${\rm SL}_2$-${\rm PGL}_2$ case, the $A_2$-bracelets basis coming from $\mathbb{I}^{\rm br}$ has a better chance. However, in fact, the only difference between these two bases is on the treatment of loops, and as for the ${\rm SL}_3$-laminations containing $3$-valent vertices, both bases are built on the non-elliptic ${\rm SL}_3$-webs. It might be the case that the non-elliptic ${\rm SL}_3$-webs are convenient only for a topological viewpoint, and are not canonical objects in terms of cluster varieties, e.g. for the purpose of Conjecture \ref{conj:GHKK}, and for example it might be a better idea to consider the version obtained by applying Fomin-Pylyavskyy's `arborization' process \cite{FP16} to these webs, to lessen the number of 3-valent vertices. Or there might even be another choice. So maybe it's better to refer to the $A_2$-bangles basis of the present paper as the $A_2$-{\it non-elliptic} basis.}

\vspace{1mm}

Another natural direction of research is on the quantization of our duality map $\mathbb{I}$. A quantum duality map was left as merely a conjecture in the first arXiv version of the present paper \cite[Conj.6.11]{Kim}. \redfix{Since the second} version, we \redfix{have} developed the ${\rm SL}_3$ quantum trace ${\rm Tr}^\omega_\Delta$ in \S\ref{sec:SL3_trace}, and proposed a quantum duality map $\mathbb{I}^q_\Delta$ in \S\ref{subsec:quantum_duality_map}\redfix{, whose naturality is proved in the subsequent work \cite{Kim21}}. One can think of comparing with other people's constructions. On the one hand, we are informed by Thang L\^e that he is working with T. Yu on the ${\rm SL}_n$ quantum trace \cite{Lcoming}, where the method of construction is different from the one of the present paper. On the other hand, a quantum version of Gross-Hacking-Keel-Kontsevich's duality map has been constructed by Davison-Mandel \cite{DM}\redfix{, which however is as non-constructive as its classical counterpart $\mathbb{I}^{\rm GHKK}$.}

\vspace{-0,7mm}

\begin{conjecture}
\, \quad \,
\begin{enumerate}
\item[\rm (1)] Davison-Mandel's quantum duality map \cite{DM} essentially coincides with our $\mathbb{I}^q_\Delta$ (or bracelets version of $\mathbb{I}^q_\Delta$).

\item[\rm (2)] L\^e\redfix{-Yu}'s ${\rm SL}_n$ quantum trace map \cite{Lcoming} essentially coincides with our ${\rm Tr}^\omega_\Delta$ when $n=3$.
\end{enumerate}
\end{conjecture}
Conjectures for the classical map $\mathbb{I}$ should have quantum counterparts, such as:
\begin{conjecture}
\label{conj:quantum_Laurent_positivity}
For each $\ell \in\mathscr{A}_{{\rm SL}_3,\frak{S}}(\mathbb{Z}^{\redfix{T}})$, $\mathbb{I}^q_\Delta(\ell)$ is a (non-commutative) Laurent polynomial in $\{\wh{X}^{\pm 1}_v \, | \, v\in \mathcal{V}(Q_\Delta)\}$ with coefficients in $\mathbb{Z}_{\ge 0}[\redfix{q^{\pm 1/18}}]$ \redfix{(emphasis on $\mathbb{Z}_{\ge 0}$)}.
\end{conjecture}
\begin{conjecture}
The coefficients $\redfix{c^q}{(\ell,\ell';\ell'')}$ in the item (5) of Thm.\ref{thm:quantum_duality_map} belong to $\mathbb{Z}_{\ge 0}[\redfix{q^{\pm 1/3}}]$ \redfix{(emphasis on $\mathbb{Z}_{\ge 0}$)}.
\end{conjecture}
\redfix{We thank Tsukasa Ishibashi for pointing out that the result of a previous joint work of the author with So Young Cho, Hyuna Kim and Doeun Oh \cite{CKKO} implies Conjecture \ref{conj:quantum_Laurent_positivity} when $\ell$ consists only of loops, i.e. has no $3$-valent vertex.}

\vs

Another possible further topic is on the $A_2$-type laminations for the tropical integer points of the cluster $\mathscr{X}$-variety. Note that for the $A_1$-type theory, or the ${\rm SL}_2$-${\rm PGL}_2$ theory, there are dual notions of $\mathscr{A}$-laminations and $\mathscr{X}$-laminations \cite{FG06}. Our ${\rm SL}_3$-lamination is an $A_2$-analog of Fock-Goncharov's $\mathscr{A}$-laminations, so there should also be the $A_2$-type counterpart for Fock-Goncharov's $\mathscr{X}$-laminations. \redfix{We proposed some preliminary steps for such new type of laminations in the previous versions (arXiv ver. 2 and 3) of the present paper, which we removed in the present version, as they are more or less fully realized in a recent work of Ishibashi and Kano \cite{IK}, which should be a major step toward the ${\rm PGL}_3$-${\rm SL}_3$ duality map $\mathscr{X}_{{\rm PGL}_3,\frak{S}}(\mathbb{Z}^t) \to \mathscr{O}(\mathscr{A}_{{\rm SL}_3,\frak{S}})$ (or its appropriate modification) and its quantum version. See also \cite{IY}.}

\vs

Some more topics are: to extend the duality map $\mathbb{I}$ to generalized marked surfaces \redfix{(with boundary)} with $\mathscr{X}_{{\rm PGL}_3,\frak{S}}$ replaced by $\mathscr{P}_{{\rm PGL}_3,\frak{S}}$ of \cite{GS19} perhaps by using the constructions of \S\ref{sec:SL3_trace} (one may want to consult \cite{IO}) \redfix{(this is in progress with L. Shen in \cite{KS})}, to consider the tropical rational and tropical real points of $\mathscr{A}_{{\rm SL}_3,\frak{S}}$ and $\mathscr{X}_{{\rm PGL}_3,\frak{S}}$ and find geometric meanings, to interpret the set $\mathscr{A}^+_{{\rm SL}_3,\frak{S}}(\mathbb{Z}^{t})$ of \cite{GS15} in terms of our ${\rm SL}_3$-laminations, to compare with results in the physics literature \cite{Xie} \cite{CGT} \cite{GLM} \cite{GMN13} \cite{Gabella} \cite{KLS} \redfix{, and to use the ${\rm SL}_3$ quantum trace ${\rm Tr}^\omega_\Delta$ to develop a (finite dimensional) representation theory of ${\rm SL}_3$ skein algebras as an ${\rm SL}_3$ analog of Bonahon-Wong's series of works on the ${\rm SL}_2$ counterpart, (see e.g. the first one \cite{BW16}), which might also lead to 3d topological quantum field theories or 2d conformal field theories.}

\vs

Lastly, an obvious way to explore is the higher \redfix{rank} generalization to ${\rm SL}_m$-${\rm PGL}_m$. \redfix{In the previous versions of the present paper (up to arXiv ver. 3) we devoted a whole subsection to some suggestion toward this direction, based on the observations of Xie \cite{Xie}, and the constructions in \cite{S01} \cite{MOY} \cite{CKM}. In the present version we removed it. It seems that several groups of mathematicians are working on it, so that soon there will be progress in the literature. For example, we refer the readers to \cite{Lcoming}.}


\begin{thebibliography}{GHKK18}

\bibitem[A19]{A19} D.G.L. Allegretti, {\it Categorified canonical bases and framed BPS states}, Sel. Math. New Ser. {\bf 25}, 69 (2019). \quad arXiv:1806.10394

\bibitem[AB20]{AB} D.G.L. Allegretti and T. Bridgeland, {\it The monodromy of meromorphic projective structures}, T. Am. Math. Soc. {\bf 373} (2020), 6321--6367. \quad arXiv:1802.02505

\bibitem[AK17]{AK} D.G.L. Allegretti and H. Kim, {\it A duality map for quantum cluster varieties from surfaces}, Adv. Math. {\bf 306} (2017), 1164--1208. \quad arXiv:1509.01567


\bibitem[BW11]{BW} F. Bonahon and \redfix{H}. Wong, {\it Quantum traces for representations of surface groups in ${\rm SL}_2(\mathbb{C})$}, Geom. Topol. {\bf 15} (2011), 1569--1615. \quad arXiv:1003.5250

\redfix{\bibitem[BW16]{BW16} F. Bonahon and H. Wong, {\it Representations of the Kauffman bracket skein algebra I: invariants and miraculous cancellations}, Invent. Math. {\bf 204} (2016) 195--243. \quad arXiv:1206.1638}


\bibitem[CKM14]{CKM} S. Cautis, J. Kamnitzer, and S. Morrison, {\it Webs and quantum skew Howe duality}, Math. Ann. {\bf 360} (2014), 351--390. \quad arXiv:1210.6437

\bibitem[CS22]{CS} L.O. Che\redfix{k}hov and M. Shapiro, \redfix{{\it Log-canonical coordinates for symplectic groupoid and cluster algebras}, Int. Math. Res. Notices (2022), rnac101.} \quad arXiv:2003.07499

\redfix{\bibitem[CKKO20]{CKKO} S. Cho, H. Kim, H. Kim and D. Oh, {\it Laurent positivity of quantized canonical bases for quantum cluster varieties from surfaces}, Commun. Math. Phys. {\bf 373} (2020), 655-705. \quad arXiv:1710.06217}


\bibitem[CGT15]{CGT} I. Coman, M. Gabella, and J. Teschner, {\it Line operators in theories of class $\mathcal{S}$, quantized moduli space of flat connections, and Toda field theory}, J. High Energy Phys. {\bf 10} (2015), 143. \quad arXiv:1505.05898

\bibitem[CL19]{CL} F. Costantino and T.T.Q. L\^e, {\it Stated skein algebras of surfaces}, to appear in J. Eur. Math. Soc. arXiv:1907.11400

\bibitem[DM21]{DM} B. Davison and T. Mandel, {\it Strong positivity for quantum theta bases of quantum cluster algebras}, \redfix{Invent. Math. {\bf 226}(3) (2021), 725--843.} \quad arXiv:1910.12915

\bibitem[D20]{Douglas} D.C. Douglas, {\it Classical and quantum traces coming from ${\rm SL}_n(\mathbb{C})$ and ${\rm U}_q(\frak{sl}_n)$}, Ph.D. thesis, 2020.

\bibitem[D21]{Douglas21} D.C. Douglas, {\it Quantum traces for ${\rm SL}_n(\mathbb{C})$: the case $n=3$}, arXiv:2101.06817

\bibitem[DS20a]{DS1} D.C. Douglas and Z. Sun, {\it Tropical Fock-Goncharov coordinates for ${\rm SL}_3$-webs on surfaces I: construction}, arXiv:2011.01768

\bibitem[DS20b]{DS2} D.C. Douglas and Z. Sun, {\it Tropical Fock-Goncharov coordinates for ${\rm SL}_3$-webs on surfaces II: naturality}, arXiv:2012.14202

\bibitem[FG06a]{FG06} V.V. Fock and A.B. Goncharov, {\it Moduli spaces of local systems and higher Teichm\"uller theory}, Publ. Math. Inst. Hautes \'{E}tudes Sci. {\bf 103} (2006), 1--211. arXiv:math/0311149v4

\bibitem[FG06b]{FG06b} V.V. Fock and A.B. Goncharov, ``Cluster $\mathscr{X}$-varieties, amalgamation, and Poisson-Lie groups, in {\em Algebraic geometry and number theory}, Progr. Math. {\bf 253}, Birkh\"auser, 2006, 27--68.


\bibitem[FG07]{FG07} V.V. Fock and A.B. Goncharov, Dual Teichm\"uller and lamination spaces, in {\it Handbook of Teichm\"uller theory. Vol. I}, IRMA Lect. Math. Theor. Phys. {\bf 11}, Eur. Math. Soc., Z\"urich, 2007, pp. 647--684. \quad arXiv:math/0510312

\bibitem[FG09]{FG09} V. V. Fock and A. B. Goncharov, {\it The quantum dilogarithm and representations of the quantum cluster varieties}, Invent. Math. {\bf 175}(2) (2009) 223--286.

\redfix{\bibitem[FP16]{FP16} S. Fomin and P. Pylyavskyy, {\it Tensor diagrams and cluster algebras}, Adv. Math. {\bf 300} (2016), 717--787.\quad arXiv:1210.1888} 

\bibitem[FZ07]{FZ07} S. Fomin and A. Zelevinsky, {\it Cluster algebras IV: Coefficients}, Compositio Math. {\bf 143} (2007), 112--164. \quad arXiv:math/0602259

\redfix{\bibitem[FP21]{FP} C. Fraser and P. Pylyavskyy, {\it Tensor diagrams and cluster combinatorics at punctures}, arXiv:2107.13069}

\bibitem[FS22]{FS} C. Frohman and A. Sikora, {\it ${\rm SU}(3)$-skein algebras and webs on surfaces}, \redfix{Math. Z. {\bf 300} (2022), 33--56}. arXiv:2002.08151



\bibitem[G17]{Gabella} M. Gabella, {\it Quantum Holonomies from Spectral Networks and Framed BPS States}, Commun. Math. Phys. {\bf 351}(2) (2017), 563--598. \quad arXiv:1603.05258

\redfix{\bibitem[GMN13]{GMN13} D. Gaiotto, G.W. Moore, and A. Neitzke, {\it Spectral networks}, Ann. Henri Poincar\'{e} {\bf 14}(7) (2013), 1643--1731. \quad arXiv:1204.4824 \quad MR3115984}


\bibitem[GLM15]{GLM} D. Galakhov, P. Longhi, and G.W. Moore, {\it Spectral networks with spin}, Commun. Math. Phys. {\bf 340}(1) (2015), 171-232. \quad arXiv:1408.0207

\bibitem[GS15]{GS15} A.B. Goncharov and L. Shen, {\it Geometry of canonical bases and mirror symmetry}, Invent. Math. {\bf 202} (2015), 487--633. \quad arXiv:1309.5922

\bibitem[GS18]{GS18} A.B. Goncharov and L. Shen, {\it Donaldson-Thomas transformations of moduli spaces of ${\rm G}$-local systems}, Adv. Math. {\bf 327} (2018), 225--348. \quad arXiv:1602.06479

\bibitem[GS19]{GS19} A.B. Goncharov L. Shen, {\it Quantum geometry of moduli spaces of local systems and representation theory}, arXiv:1904.10491

\bibitem[GHK15]{GHK} M. Gross, P. Hacking, and S. Keel, {\it Birational geometry of cluster algebras}, Algebraic Geometry {\bf 2}(2) (2015), 137--175. \quad arXiv:1309.2573

\bibitem[GHKK18]{GHKK} M. Gross, P. Hacking, S. Keel, and M. Kontsevich, {\it Canonical bases for cluster algebras}, J. Am. Math. Soc. {\bf 31}\redfix{(2)} (2018), 497--608. \quad arXiv:1411.1394

\bibitem[H10]{Hiatt} C. Hiatt, {\it Quantum traces in quantum Teichm\"uller theory}, Algebr. Geom. Topol. {\bf 10}(3) (2010), 1245--1283. \quad arXiv:0809.5118

\bibitem[H20]{Higgins} V. Higgins, {\it Triangular decomposition of ${\rm SL}_3$ skein algebras}, arXiv:2008.09419

\redfix{\bibitem[IK22]{IK} T. Ishibashi and S. Kano, {\it Unbounded $\frak{sl}_3$-laminations and their shear coordinates}, arXiv:2204.08947}

\bibitem[IO20]{IO} T. Ishibashi and H. Oya, {\it Wilson lines and their Laurent positivity}, arXiv:2011.14260

\bibitem[IY21]{IY} T. Ishibashi and W. Yuasa, {\it Skein and cluster algebras of marked surfaces without punctures for $\frak{sl}_3$}, arXiv:2101.00643

\redfix{\bibitem[JK]{JK} S.-J. Jung and H. Kim, {\it On self-folded triangulations in ${\rm SL}_n$ cluster varieties for surfaces}, in preparation.}

\bibitem[K20]{Kim} H. Kim, {\it $A_2$-laminations as basis for ${\rm PGL}_3$ cluster variety for surface}, arXiv:2011.14765v1 (previous version of the present paper)

\redfix{\bibitem[K21]{Kim21} H. Kim, {\it Naturality of the ${\rm SL}_3$ quantum trace maps for surfaces}, arXiv:2104.06286}

\bibitem[KLS18]{KLS} H. Kim, T.T.Q. L\^{e}, and M. Son, {\it ${\rm SL}_2$ quantum trace in quantum Teichm\"uller theory via writhe}, arXiv:1812.11618v2

\redfix{\bibitem[KS]{KS} H. Kim and L. Shen, in preparation.}

\redfix{\bibitem[KT99]{KT99} A. Knutson and T. Tao, {\it The honeycomb model of $GL_n(\mathbb{C})$ tensor products I: proof of the saturation conjecture}, J. Am. Math. Soc. {\bf 12}(4) (1999), 1055-1090.}

\bibitem[K96]{Kuperberg} G. Kuperberg, {\it Spiders for Rank $2$ Lie Algebras}, Commun. Math. Phys. {\bf 180} (1996), 109--151. \quad arXiv:q-alg/9712003

\bibitem[L16]{Ian} I. Le, {\it Higher laminations and affine buildings}, Geom. Topol. {\bf 20} (2016), 1673--1735. \quad arXiv:1209.0812 

\bibitem[L17]{Le17} T.T.Q. L\^{e}, {\it Quantum Teichm\"uller spaces and quantum trace map}, J. Inst. Math. Jussieu (2017), 1--43. \quad arXiv:1511.06054

\bibitem[L18]{Le18} T.T.Q. L\^{e}, {\it Triangular decomposition of skein algebras}, Quantum Topol. {\bf 9} (2018), 591--632. \quad arXiv:1609.04987

\bibitem[LY]{Lcoming} T.T.Q. L\^{e} \redfix{and T. Yu}, {\it Quantum trace for ${\rm SL}_n$-skein algebras}, in preparation.

\redfix{\bibitem[MQ]{MQ} T. Mandel and F. Qin, {\it Bracelet bases are theta bases}, in preparation.}

\bibitem[MOY98]{MOY} H. Murakami, T. Ohtsuki, and S. Yamada, {\it Homfly polynomial via an invariant of colored plane graphs}, Enseign. Math. (2) {\bf 44} (1998), 325--360.

\bibitem[P76]{P76} C. Procesi, {\it The invariant theory of $n\times n$ matrices}, Adv. Math. {\bf 19} (1976), 306--381.

\bibitem[RT90]{RT} N.Y. Reshetikhin and V.G. Turaev. {\it 
Ribbon graphs and their invariants derived from quantum groups}, Comm\redfix{un}. Math. Phys. {\bf 127}(1) (1990), 1--26. 


\bibitem[SSh19]{SchSh} G. Schrader and A. Shapiro, {\it A cluster realization of $U_q(\frak{sl}_n)$ from quantum character varieties}, Invent. Math. {\bf 216} (2019), 799-846.


\bibitem[S01]{S01} A.S. Sikora, {\it ${\rm SL}_n$-character varieties as spaces of graphs}, \redfix{T. Am. Math. Soc.} {\bf 35} no.7 (2001), 2773--2804. \quad arXiv:math/9806016

\bibitem[S05]{S05} A.S. Sikora, {\it Skein theory for $SU(n)$-quantum invariants}, \redfix{Alg. Geom. Topol.} {\bf 5} (2005), 865--897. \quad arXiv:math/0407299

\bibitem[SW07]{SW} A.S. Sikora and B.W. Westbury, {\it Confluence theory of graphs}, Alg. Geom. Topol. {\bf 7} (2007), 439--478. \quad arXiv:math/0609832

\bibitem[S20]{Shen} L. Shen, {\it Duals of semisimple Poisson-Lie groups and cluster theory of moduli spaces of $G$-local systems}, \redfix{Int. Math. Res. Notices (2021), rnab094}. \quad arXiv:2003.07901

\bibitem[T14]{DThurston} D.P. Thurston, {\it Positive bases for surface skein algebras}, Proc. Natl. Acad. Sci. {\bf 111}(27) (2014), 9725--9732. \quad arXiv:1310.1959

\bibitem[X13]{Xie} D. Xie, {\it Higher laminations, webs and $\mathcal{N}=2$ line operators}\redfix{,} arXiv:1304.2390
\end{thebibliography}
\end{document}